 \documentclass{amsbook}
\usepackage{lmodern}
\numberwithin{equation}{subsection}

 \makeindex
\usepackage{amsmath}
\usepackage{amsfonts}
 \usepackage{amssymb}
 \usepackage{textcomp}
\usepackage{mathrsfs}
\usepackage{minitoc}
\usepackage{latexsym}
\usepackage[pdf]{diagrams} 
\usepackage{verbatim}
\usepackage{graphicx}
\usepackage{nomencl}
\usepackage{colortbl}
\usepackage{hyperref}
\renewcommand{\marginpar}[2][]{}

\DeclareFontEncoding{OT2}{}{} %
  \newcommand{\textcyr}[1]{%
    {\fontencoding{OT2}\fontfamily{wncyr}\fontseries{m}\fontshape{n}%
     \selectfont #1}}
\newcommand{\Sha}{{\mbox{\textcyr{Sh}}}}

\def\PPP{\mathbf{P}}
\def\hyp{\mathrm{hyp}}
\def\Tamagawa{\mathrm{Tamagawa}}
\def\kij{i\kern-0.03em{j}}
\def\JL{J\kern-0.03em{L}}
\def\comp{\#Y}

\def\Sengun{\c{S}eng\"{u}n}
\def\hangingchad{phantom class}
\def\ev{\mathrm{ev}}
\def\LL{\mathscr{L}}
\def\LLL{\widetilde{\LL}}
\newcommand{\NN}{\mathbf{N}}
\newcommand{\coclosed}{\mathrm{coclosed}}
\reversemarginpar
\newcommand{\Gdual}{\widehat{G}}
\newcommand{\height}{\mathrm{Ht}}

\newcommand{\tr}{\mathrm{tr}}
\newcommand{\torsnew}{h_{\mathrm{tors}}^\new}
\newcommand{\diagA}{\mathbf{A}}
\newcommand{\torsnewE}{h_{\mathrm{tors}}^{E,\new}}
\newcommand{\iso}{\iota}
\newcommand{\torscong}{h_{\mathrm{cong}}^{\new}}
\newcommand{\Sigmainvertq}{\Sigma[\frac{1}{{\q}}]}
\newcommand{\KSigmainvertq}{K_{\Sigma}[\frac{1}{{\q}}]}
\newcommand{\logdetR}{\log \det {}^*}
\newcommand{\Kmax}{K_{\mathrm{max}}}
\newcommand{\num}{\mathrm{num}}

\newcommand{\escapeclause}{except possibly at orbifold primes $($which 
are at most $3)$}
\newcommand{\emb}{\mathrm{emb}}
 \newcommand{\bdy}{\mathrm{bdy}}

\newcommand{\trace}{\mathrm{tr}}
\newcommand{\alg}{\mathrm{alg}}
\newcommand{\detprime}{\det{}'}

\newcommand{\qnew}{\q-\mathrm{new}}
\newcommand{\disc}{\mathrm{disc}}

\newcommand{\Tmax}{T_{\mathrm{max}}}

\newcommand{\bm}{\mathrm{BM}}

\newcommand{\PU}{\mathrm{PU}}
\newcommand{\varF}{\upsilon}

\newcommand{\one}{E}

\def\Cat{\mathscr{C}}

\newcommand{\area}{\mathrm{area}}

\def\eisideal{\mathfrak{I}}

\def\heart{\xi}
\def\betavar{\gamma}

 \def\lif{\mathrm{lif}}
\def\cl{\mathrm{lif}}
\def\divi{\mathrm{div}}
\def\tors{\mathrm{tors}}
\def\baromega{\bar{\omega}}
\def\tf{\mathrm{tf}}
\newcommand{\varchi}{\mathcal{X}}

\def\etale{\'{e}tale }
\def\CL{\mathrm{Cl}}
\def\Vbar{\overline{V}}
\def\Lbar{\overline{L}}
\def\Def{\mathrm{Def}}
\def\Adele{\mathbb{A}}
\def\Afinite{\mathbb{A}_{\kern-0.05em{f}}}
\def\testring{R}
\def\Kmax{ \mathrm{PU}_2 \times \PGL_2(\widehat{\OL})}
\def\Sigmar{\Sigma \kern+0.1em{\rf}}

\def\PP{{\frakn}}
\def\Dfl{D_{\mathrm{fl}}}
\def\Dord{D_{\mathrm{ord}}}
\def\Ind{\mathrm{Ind}}
\newcommand{\carrot}{\wedge}

\newcommand{\cusp}{\mathrm{cusp}}

\newcommand{\Eis}{\mathrm{Eis}}
\def\Hfl{H^1_{\mathrm{fl}}}
\def\Hord{H^1_{\mathrm{ord}}}

\def\Hflv{H^1_{\mathrm{fl},v}}
\def\Hordv{H^1_{\mathrm{ord},v}}

\def\weightk{\underline{k}}
\newcommand{\Sigmainf}{\Sigma_{\infty}}

\newcommand{\Tree}{\mathcal{T}}
\newcommand{\Doubletree}{\mathcal{T}^{\pm}}
\newcommand{\hrel}{h_{\mathrm{rel}}}

\def\KM{K_{\Sigma}}
\def\KMO{K_{\Sigma,0}}
\def\KN{K_{\Sigma,1}}
\def\YO{Y(\KM)}
\def\YN{Y(\KN)}

\def\YOq{Y(K_{\Sigma/\q})}

\def\WT{\widetilde{\T}}
\def\End{\mathrm{End}}
\def\Tnew{\T^{\mathrm{new}}}
\def\frakn{\mathfrak{n}}
\def\new{\mathrm{new}}

\newcommand{\finite}{\kern-0.05em{f}}

\def\Qbar{\overline{\Q}}
\newcommand{\gk}{\mathfrak{g}/\mathfrak{k}}

\def\rec{\mathrm{rec}}

\def\rf{\mathfrak{r}}

\def\Fbar{\overline{\F}}

\def\Gal{\mathrm{Gal}}
\def\eps{\epsilon}
\def\Sym{\mathrm{Sym}}

\def\GammaP{\Gamma^{\kern+0.07em{p}}}
\def\GammaP{\Gamma}
\def\AN{A}
\newcommand{\M}{\mathscr{M}}

\def\Sym{\mathrm{Sym}}
\def\Atimes{\kern+0.1em{\A^{\kern-0.1em{\times}}}}

\def\vol{\mathrm{vol}}
\def\OL{\mathscr{O}}
\def\T{\mathbf{T}}
\def\m{\mathfrak{m}}
\def\A{\mathbb{A}}
\newcommand{\indexK}{[\PGL_2(\widehat{\OO}):K]}
\def\G{\mathbb{G}}
\def\analT{\tau_{\mathrm{an}}} %
\def\an{\mathrm{an}}
\def\RT{\mathrm{RT}}

\def\a{\mathfrak{a}}
\def\b{\mathfrak{b}}

\def\Gammaab{\Gamma_{\kern-0.1em{\a\b}}}

\def\Mp{\M_p}
\def\rhobar{\overline{\rho}}
\def\Frob{\mathrm{Frob}}

\def\PGL{\mathrm{PGL}}

\def\Ht{\tilde{H}}

\def\coker{\mathrm{coker}}
\def\Cl{\mathrm{Cl}}
\def\q{{\mathfrak{q}}}
\def\PSL{\mathrm{PSL}}
\def\SL{\mathrm{SL}}

\def\con{\mathrm{cong}}
\def\Hom{\mathrm{Hom}}
\def\F{\mathbf{F}}
\def\Q{\mathbf{Q}}
\def\GL{\mathrm{GL}}
\def\H{\mathcal{H}}
\def\p{{\mathfrak{p}}}
\def\Z{\mathbf{Z}}
\def\C{\mathbf{C}}

\newarrow{Equals}{=}{=}{=}{=}{=}
\newarrow{Dotsto}{.}{.}{.}{.}{>}
\newarrow{Dotto}{.}{.}{.}{.}{.}
\newarrow{Doppo}{}{.}{}{.}{}

\def\new{\mathrm{new}}
\def\old{\mathrm{old}}
\def\Tnew{\T^{\new}}

\def\cF{\mathcal{F}}

\newtheorem{conj}[subsection]{Conjecture} 
\newtheorem{question}[subsection]{Question} 
\newtheorem{df}[subsection]{Definition}

\newtheorem*{df*}{Definition}
\newtheorem*{lemma*}{Lemma}
\newtheorem*{theorem*}{Theorem}
\newtheorem*{theoremA*}{Theorem~A}
\newtheorem*{theoremB*}{Theorem~B}
\newtheorem*{theoremC*}{Theorem C}
\newtheorem*{theoremD*}{Theorem~A$^{\dagger}$}

\newtheorem{theorem}[subsection]{Theorem}
\newtheorem{lemma}[subsection]{Lemma}

\newtheorem{example}[subsection]{Example}
\newtheorem{corollary}[subsection]{Corollary}
\newtheorem{definition}[subsection]{Definition}

\newtheorem{remarkable}[subsection]{Remark}
\newtheorem*{assumption*}{Assumption}

\theoremstyle{remark}
\newtheorem{remark}[subsection]{Remark}
\newtheorem*{remark*}{Remark}

\def\con{\mathrm{cong}}
\def\R{\mathbf{R}}

\newcommand{\B}{\mathbf{B}}
 \newcommand{\Ad}{\mathrm{Ad}}
 \newcommand{\Ext}{\mathrm{Ext}}

\newcommand{\adele}{\mathbb{A}}
\newcommand{\ab}{\mathrm{ab}}
\newcommand{\localsystem}{\mathscr{M}}
\newcommand{\reg}{\mathrm{reg}}
\newcommand{\pibar}{\overline{\pi}}
\newcommand{\data}{\mathscr{D}}

\newcommand{\valuefield}{\C}

\renewcommand{\a}{\mathfrak{a}}
\newcommand{\ellcurve}{C}

\newcommand{\order}{\mathscr{O}}
\newcommand{\N}{\mathbf{N}}

 \renewcommand{\H}{\mathbf{H}}

\newcommand{\im}{\mathrm{im}}

\newcommand{\n}{{\frakn}}

\newcommand{\Spec}{\mathrm{Spec}}

\newcommand{\Norm}{\mathrm{N}}
\newcommand{\OO}{\mathscr{O}}

\DeclareFontFamily{OT1}{rsfs}{}

\newcommand{\CM}[2]{\ensuremath{H^1({#2}, {#1})}}
\newcommand{\HM}[2]{\ensuremath{H_1({#2}, {#1})}}

\newcommand{\constWp}{\Q_p/\Z_p}
\newcommand{\constMp}{\Z_p}

\DeclareFontShape{OT1}{rsfs}{n}{it}{<-> rsfs10}{}
\DeclareMathAlphabet{\mathscr}{OT1}{rsfs}{n}{it}

\newcommand{\SO}{\mathrm{SO}}

\begin{document}

\frontmatter	
\title{A torsion  Jacquet--Langlands correspondence}
\author{Frank Calegari and Akshay Venkatesh}

 \begin{abstract}
  We study  torsion in the homology of arithmetic groups and give evidence
that  it plays a role in the Langlands program.  We prove, among other results,  a numerical form of a Jacquet--Langlands correspondence in the torsion setting. 
 \end{abstract}

\maketitle
\tableofcontents

\textit{Acknowledgments}
\addcontentsline{toc}{section}{Acknowledgments}

The first author (F.C.) would like to thank Matthew Emerton for many conversations regarding
possible 
integral formulations of reciprocity, and to thank Nathan Dunfield; the original calculations
suggesting that a torsion Jacquet--Langlands theorem might be true were done
during the process of writing our joint paper~\cite{CD}, and the data thus
produced proved very useful for  suggesting some of the phenomena we have
studied in this book. He would also like to thank Kevin Hutchinson for some helpful
remarks concerning a theorem of Suslin.
F.C.  was supported during the preparation of this
paper by a  Sloan fellowship and by the National Science Foundation 
(CAREER Grant DMS-0846285).

\medskip

The second author (A.V.) would like to express his gratitude to Nicolas Bergeron,
with whom he wrote the paper~\cite{AB}, for many fruitful discussions and thoughts
about the analytic behavior of analytic torsion; to Laurent Clozel, who suggested to him
the importance of investigating torsion in cohomology; to Kartik Prasanna, 
for helpful discussions and references concerning his work;   and to Avner Ash for encouraging words.  
He was supported during the preparing of this book by a Sloan fellowship, an 
National Science Foundation grant,
and a David and Lucile Packard Foundation fellowship; he gratefully thanks these organizations for their support. 

\medskip

 The mathematical debt this book owes to the work of Cheeger \cite{Cheeger} and M\"{u}ller \cite{Muller} should be clear: the main result about comparison of torsion homology
 makes essential use of their theorem. 

 \medskip

Some of the original ideas of this book were conceived during the conference ``Explicit Methods in Number Theory''
in Oberwolfach during the summer of 2007. Various parts of the manuscript were written while the authors were resident at the following institutions:
New York University, Stanford University, Northwestern University,  the Institute for Advanced Study,  the University of Sydney,
and the Brothers~K coffee shop; we thank all of these institutions for their hospitality.

\medskip
We thank Romyar Sharifi  and  Krzysztof Klosin for helpful comments, especially on the material of Chapter~\ref{chapter:ch3}.  We thank Aurel Page for providing a presentation of 
$\displaystyle{\PGL_2 \left(\Z\left[\frac{1 + \sqrt{-491}}{2}\right] \right)}$
used in~\S~\ref{section:K2examples}.

\mainmatter

%
%

%
%
%
%

%

 \chapter{Introduction}
 \label{chapter:ch1}
 \section{Introduction}

The main goal of this manuscript is to substantiate the
claim  that torsion in the homology of arithmetic groups  plays a significant
 role in the arithmetic Langlands program. 
 
 \medskip
 
The Langlands program conjecturally relates homology
 of arithmetic groups   to Galois representations;
 an emerging extension of this program predicts that this relationship exists  not only
 for the characteristic zero homology but also for the torsion.
One consequence of these conjectures is that the integral  homology
of arithmetic groups for different inner forms of the same group will be related in a
 non-obvious way, and this conjectural correspondence can be studied
even without discussing Galois representations.

 \medskip
 
An interesting class of arithmetic groups is  provided by arithmetic
  Fuchsian subgroups $\Gamma \leqslant \PGL_2(\C) \cong \SL_2(\C)/\{\pm I \}$,
  \nomenclature{$\Gamma$}{A Fuchsian group}%
   in particular, $\PGL_2(\OL_F)$ and its congruence subgroups,
   where  $F$ is an imaginary quadratic field. For such groups, it has been
   observed (both numerically and theoretically) that $H_1(\Gamma, \mathbf{Z}) = \Gamma^{\mathrm{ab}}$ has a large torsion
   subgroup~\cite{AB,Sengun}.  
 This paper studies and verifies several predictions about the homology of these groups
 (and their corresponding hyperbolic $3$-manifolds).  
 
 \medskip

In particular, we establish a number of results showing that torsion behaves according
to the predictions of the Langlands program: \begin{enumerate} \label{firstpagephen}
  \item Relationships between $H_{1, \tors}(\Gamma)$ and $H_{1, \tors}(\Gamma')$,
for certain {\em incommensurable} groups $\Gamma, \Gamma'$  (``Jacquet--Langlands correspondence,''  Chapter 7). 
 \item Relationships between $H_{1, \tors}(\Gamma)$ and $H_{1, \tors}(\tilde{\Gamma})$,
  where $\tilde{\Gamma}$ is a congruence subgroup of $\Gamma$ (``level raising,''
  see Chapter 4). 
\item $H_1(\Gamma,\F_p) \neq 0$ when the 
Galois cohomology group $H^1(\mathrm{Gal}(\overline{F}/F),\mu_p^{-1}) $
admits everywhere unramified classes,
and
$\Gamma$ is contained in a congruence subgroup $\Gamma_0(\p)$
of $\PGL_2(\OL_F)$ (see  Theorem~\ref{theorem410}, 
and \S \ref{section:Eisenstein}). 
  \end{enumerate}
    
We also study several other phenomena related to torsion: 
 \begin{enumerate}
\item[4.] Regulators for arithmetic manifolds (Chapter 5), 
whose behavior influences the behavior of torsion in many of our numerical examples; 
\item[5.]  Phenomena related to restricting an {\em even} Galois representation over $\Q$
to an imaginary quadratic field (Chapter 8) 
  \item[6.] A large collection of numerical examples (Chapter 9);
\end{enumerate}

 Our main  tool to analyze~(1) is {\em analytic torsion}. 
To better understand~(1), and also to study~(2) and (3), we will use other techniques, for example: 
 the congruence subgroup property, Waldspurger's formulas, and Tate's theorem on $K_2$.   
 Similarly the analysis of (4)--(6) uses various {\em ad hoc} tools.  
   \medskip

Although the final scope is much broader, this book   grew out of an attempt to understand phenomenon~(1).   In the next subsection~\S~\ref{section:introtheorems} we will  introduce the paper by describing
a theorem of type~(1), and how its generalizations are related to~(2) and~(3).

\medskip
An apology for length is in order. As mentioned,  this book began solely as an investigation
of a possible Jacquet--Langlands correspondence for torsion classes. However,  along the way,
we took many detours to explore related phenomena, some of which was inspired by the data computed for the first author
by Nathan Dunfield.
In view of the almost complete lack of rigorous understanding of torsion  for (non-Hermitian) locally symmetric spaces,   we have included many of these results, even when what we can prove is rather modest. Moreover, the results of this book involve several areas of mathematics, and we have therefore 
tried to be somewhat self-contained, often including proofs of known results, as well
as precise formulations of folklore conjectures.

  \subsection{}  \label{section:introtheorems} 

   Let $F/\Q$ be an imaginary quadratic field of odd class number\footnote{The assumption on the class number is made in this section (and this section only) 
purely for exposition --- it implies that   the relevant adelic quotient has only one connected component, so that the manifold $Y_0(\n)$ considered below
is indeed the ``correct'' object to consider.},
   and ${\frakn}$ an ideal of $F$; let ${\p}, {\q}$ be prime ideals of $\OO_F$
  that do not divide ${\frakn}$. 
Let $\Gamma_0({\frakn})$ be the subgroup of $\PGL_2(\OO_F)$
corresponding to matrices $ \left( \begin{array}{cc} a & b \\ c & d \end{array} \right)$
where ${\frakn}|c$.

  \medskip

  Let $D$ be the unique quaternion algebra over $F$ ramified at $\p, \q$, 
 $\OO_D$ a maximal order in $D$. Let $\Gamma_0'$
 be the image of $\OO_D^{\times}/\OO_F^{\times}$ in $\PGL_2(\C)$; 
 it is a co-compact lattice. 
  There is  a canonical conjugacy class of surjections
 $\Gamma_0' \twoheadrightarrow \PGL_2(\OO_F/{\frakn})$, and thus
 one may define the analogue $\Gamma_0'({\frakn})$ of the subgroup
 $\Gamma_0({\frakn})$ in this context as the preimage in $\Gamma_0'$ of the upper triangular subgroup of $\PGL_2(\OO_F/{\frakn})$. 
 It can be alternately described as $\mathscr{O}_{D, {\frakn}}^{\times}/\OO_F^{\times}$, 
 where $\mathscr{O}_{D, {\frakn}}$ is an Eichler order of level ${\frakn}$.

 \medskip
 
 Let $Y_0({\frakn})$ be the  finite volume, hyperbolic $3$-manifold ${\H}^3/\Gamma_0({\frakn})$. 
 Note that $Y_0({\frakn})$ is {\em non-compact}.
 Let $Y_0'({\frakn})$ be the
 quotient of ${\H}^3$ by $\Gamma_0'({\frakn})$;  it is a
 {\em compact} hyperbolic $3$-manifold.  
 
 \medskip
 
   We refer to \cite{MaclachlanReid} for background and further details
  on these types of construction, especially Chapter 6 for orders in quaternion algebras. 
 
 \medskip
 
  The homology group of both $Y_0({\frakn})$ and $Y_0'({\frakn})$
is equipped with a natural action of the {\em Hecke algebra} (see~\S~\ref{ss:HeckeDef} for definitions).  A theorem of Jacquet--Langlands
implies that there is a {\em Hecke-equivariant injection}  
\begin{equation} \label{jl} \dim H_1( Y_0'({\frakn}), \C) \hookrightarrow  \dim H_1(Y_0({\frakn} \p \q), \C).\end{equation}
The proof of~\eqref{jl} is based on the fact that the length spectra of the two manifolds
are closely related.

When considering \emph{integral} homology, it will be convenient to replace
the group $H_1(Y,\Z)$ by a slightly different group $H^{E^*}_{1}(Y,\Z)$ (the
\emph{dual-essential homology}, see definition~\ref{df:essentialhomology}
and~\S~\ref{ss:Noncompactcongruence}). 
The difference between the two groups is easy to understand and roughly accounted for
by classes which arise from congruence quotients of the corresponding arithmetic group
(cf.  \S \ref{s:congess}).

\medskip

We present several results which are examples of what
can be proven by our methods:

 \begin{theoremA*}   \label{theorem:unnamedtheorem1}  (Proved in
~\S~\ref{section:relations}.)  Suppose that  $\dim H_1(Y_0({\frakn} \p \q), \C)= 0$.
 Then the order of the finite group $H^{E^*}_1(Y_0'({\frakn}), \Z)$ divides the order of the finite group 
 $H^{E^*}_1(Y_0({\frakn} \p \q), \Z)^{\p\q-\new}$, \escapeclause. 
 \end{theoremA*}
 
 The superscript $\p\q-\new$ means  ``new at $\p$ and $\q$'':   we take the quotient
 of the homology group by classes pulled-back from levels ${\frakn}\p$ and ${\frakn} \q$.

 The proof uses the (usual) Jacquet--Langlands correspondence, 
the   Cheeger--M{\"u}ller theorem \cite{Cheeger,Muller} and
 the congruence subgroup property for $S$-arithmetic $\SL_2$. 
In order to apply the Cheeger--M{\"u}ller theorem to {\em non-compact} hyperbolic manifolds, 
we are forced to address a certain number of technical issues,  although we
are able to take a shortcut especially  adapted to our situation.

 \begin{remark}
 Our results apply under substantially weaker hypotheses; the strong vanishing assumption
 is for simplicity only. 
   
 However, even this rather strong vanishing hypothesis is frequently satisfied. 
 Note, for example, that the cusps of $Y_0(\mathfrak{m})$   do not
 contribute to its complex cohomology  (Lemma~\ref{vanishingboundaryhomology}): their cross-section is 
 a quotient of a torus by negation.  Similarly, $\SL_2(\OO_F)$
 very often has base-change homology, but this often
 does not extend to $\PGL_2(\OO_F)$, cf. Remark \ref{remark:ktrivial}, second point.

 The constraint on $2$ and $3$ also arises from orbifold issues:
 for $3$  this is not a serious constraint and could presumably be removed with more careful analysis.
  However, we avoid the prime $2$ at various other points
  in the text  for a number of other technical reasons.
 \end{remark}

\begin{remark} 
There is (conditional on the congruence subgroup property for certain division algebras) a corresponding result in a situation where $Y_0, Y_0'$ both arise from quaternion algebras. In that case, one can {\em very quickly} see, 
using the Cheeger--M{\"u}ller theorem and the classical Jacquet-Langlands 
correspondence, that there should be {\em some} relationship between 
between the sizes of torsion homology of $Y_0$ and $Y_0'$. 
The reader may wish to examine this very short and  simple argument (see Chapter 7, up to the
first paragraph of the proof of Theorem~\ref{TJL1}). 

Nonetheless, this book focuses on the case where one of $Y_0$ and $Y_0'$ is split, thus
forcing us to handle, as mentioned above, a host of analytic complications. 
One reason is that, 
in order to sensibly ``interpret'' Theorem~\ref{TJL1} -- e.g., in order to massage
it into an arithmetically suggestive form like Theorem ~A or (below) Theorem~B,      one needs the congruence subgroup property
for the $S$-units -- and, as mentioned, this is not known for a quaternion algebra.  
But also we find and study many interesting phenomenon, both analytic and arithmetic,  peculiar to the split case (e.g., the study of Eisenstein deformations in Chapter 9, or the relationship to $K$-theory provided by Theorem
\ref{theorem:K2popularversion}). 
\end{remark}

\medskip

It is natural to ask for a more precise version of Theorem~A, replacing
divisibility by an equality.     As we shall see, this is {\em false}  --
and false for  an interesting reason related to the nontriviality of $K_2(\OO_F)$. For example we show:

\begin{theoremD*}   \label{theorem:theoremD}
(Proved in
~\S~\ref{section:relations}; see also following sections.)
Notation and assumptions as in Theorem $A$, suppose that ${\frakn}$ is the trivial ideal.
Then the ratio of orders 
$$\frac{|H^{E^*}_1(Y_0( \p \q), \Z)^{\p\q-\new}|} {|H^{E^*}_1(Y_0'(1), \Z)|}$$
 is divisible by the order of $K_2(\OO)$, up to   primes dividing $6 \gcd(N\p-1, N\q-1)$. 
\end{theoremD*}
Although our computational evidence is somewhat limited, we have
no reason to suppose that the divisibility in this theorem is not actually an
\emph{equality}.
We refer to Theorem~\ref{theorem:K2popularversion}  for a related result,~\S~\ref{section:Eisenstein}
for a discussion of the ``Galois side'' of the story, 
and~\S~\ref{section:K2examples} for  a numerical example.
To further understand the relationship of $K$-theory and the torsion Langlands program
seems a very interesting task.   %

\medskip

\medskip

If $H_1(Y'_0({\frakn}), \C) \neq 0$, then we wish to compare
two infinite groups. In this context, a host of new complications, including
issues related to the level raising phenomena (i.e. (2)   described on page \pageref{firstpagephen}), arise. In this case, instead
of controlling the ratios of the torsion subgroups between the split
and non split side, we relate this ratio to a ratio of periods arising
from automorphic forms $\pi$ and their  Jacquet--Langlands
correspondent $\pi^{\JL}$. This period ratio is related to work of Prasanna for $\GL(2)/\Q$
(see~\ref{section:prasanna}).

\medskip
In this context, a sample of what we prove is:

 \begin{theoremB*}  \label{theorem:unnamedtheorem2} (Proved in
~\S~\ref{section:relations}.)  Suppose that  $l$ is a prime $\notin \{2,3\}$  such that:  
 \begin{enumerate}
 \item $H_1(Y_0({\frakn} \q),\Z_l) = 0$,
 \item $H_1(Y_0({\frakn} \p),\Z_l)$ and $H_1(Y_0({\frakn} \p \q),\Z_l)$ are $l$--torsion free,
 \item $H_1(Y'_0({\frakn}),\Q_l) = 0$ --- equivalently: $H_1(Y'_0({\frakn}), \C)=0$.
 \end{enumerate}
Then $l$ divides $H_1^{E^*}(Y'_0({\frakn}),\Z)$ if and only if $l$ divides $\Delta$, where
 $(N(\p) - 1) \Delta$ is the determinant of 
 $T_{\q}^2 - (N(\q)+1)^2$ acting on $H_1(Y_0(\frakn  \p), \Q)$, and $T_{\q}$ %
 denotes the $\q$-Hecke operator.
 \end{theoremB*}

Thus the naive analog of Theorem A -- simply asking for a divisibility 
between the sizes of torsion parts -- is {\em false} here:
We have extra torsion on the compact $Y_0'({\frakn})$,  at primes dividing $\Delta$. By ``extra,''  we mean that there is no corresponding torsion for $Y_0(\mathfrak{n p q})$. 
How should we interpret this extra torsion, i.e. how to interpret the number $\Delta$ in terms
of modular forms on the split manifold $Y_0$?

\medskip

Suppose for simplicity that $\ell$ does not divide $N(\p)-1$. Then
  $\ell$ divides $\Delta$ exactly when $T_{\q}^2 - (N(\q)+1)^2$
has a nontrivial kernel on $H_1(Y_0({\frakn}{\p}), \F_{\ell})$. 
In view of this, an interpretation of $\Delta$ is given by the following ``level raising theorem.''
We state in a somewhat imprecise form for the moment. 
(To compare with the previous theorem, replace $\mathfrak{m}$
by $\mathfrak{n p}$). 

\begin{theoremC*}  \label{theorem:unamiC} (Proved in~\S~\ref{sec:llsplit}; see also~\S~\ref{sec:lr}.)
Suppose $[c]  \in H_1(Y_0(\mathfrak{m}),\F_{l})$ %
is a non-Eisenstein (see Definition~\ref{df:Eisenstein}) Hecke eigenclass, in the kernel of $T^2_{\q} -   (\Norm(\q) + 1)^2$.

Then there exists $[\tilde{c}] \in H_1(Y_0(\mathfrak{m} \q), \F_{l})$ %
with the same generalized Hecke eigenvalues as $c$ at primes not dividing $\q$, and 
not in the image of the pullback degeneracy map
$H_1(Y_0(\mathfrak{m}))^2 \rightarrow H_1(Y_0(\mathfrak{m} \q))$. 
\end{theoremC*}

\medskip
To understand these Theorems (and their relationship) better,   let us examine the situation over $\Q$ and 
explain analogues of Theorems B and C in that setting.

\medskip
  Let $N$ be an integer and $p,q$ primes,
and denote by $\Gamma_0(N), \Gamma_0(Np)$ etc. the corresponding
congruence subgroups of $\SL_2(\Z)$; denote by $\Gamma_0^*(N)$
the units in the level $N$ Eichler order inside the quaternion algebra over $\Q$
ramified at $p, q$. 
Define $\Delta$ analogously, replacing $H_1(Y_0(\mathfrak{n p}), \C)$
by the homology $H_1(\Gamma_0(N p), \C)$ and $T_{\q}$
by the usual Hecke operator $T_q$. 

\medskip

Then, if $l$ divides $\Delta$,
there exists a Hecke eigenclass $f \in H_1(\Gamma_0(Np), \mathbf{F}_{l})$
which is annihilated by $T_q^2 - (q +1)^2$. Under this assumption there exists:
\begin{itemize} 
\item[(b)]
a Hecke eigenform $f^* \in H_1(\Gamma_0^*(N), \F_{l})$

\item[(c)]
a ``new at $q$'' Hecke eigenform
$\tilde{f} \in H_1(\Gamma_0(N pq), \F_{l})$  %
 \end{itemize}
 both with the same (mod-$l$) Hecke eigenvalues as $f$.

\medskip
Indeed, (c) follows from Ribet's level-raising theorem, and (b) then follows
using (c) and Jacquet--Langlands.  In terms of representation theory this phenomena can be described thus: $f$
gives rise to a principal series representation $\pi$ of $\GL_2(\Q_q)$
on an $\F_{l}$-vector space.  The condition that $T^2_{\q} - (N(\q) + 1)^2 $
kills $f$ forces this representation to be {\em reducible}:
it contains 
a twist of the Steinberg representation.

\medskip

Our Theorems B and C can be seen as analogues   of (b) and (c) respectively. In all cases one
produces homology at primes $l$ dividing $\Delta$ on a quaternion algebra or at higher level.
There are significant differences: most notably, in case (b) the class $g$ lifts to characteristic zero
whereas in case (B) it is a pure torsion class.  More distressing is that Theorem~B, as well
as Theorem~A and all our related results, gives no information about Hecke eigenvalues. 
\medskip

\section{An Example}   \label{section:example} 

We present now a single example to give the reader some  numerical sense of the phenomenon being studied; we shall discuss the relationship between the homology of two specific 
hyperbolic $3$-manifolds.   Our example here is not an example of the situation described   in \S~\ref{section:introtheorems}  -- here,  both manifolds are compact -- but it does fit in the more general framework
we consider in Chapter 7.   (Also, Chapter 9 studies in detail many numerical examples
in the setting of \S~\ref{section:introtheorems}.)

Let $F = \Q(\theta)$ where $\theta^3 - \theta + 1 = 0$. 
If $k$ is a rational integer and there exists a unique prime ideal of norm $k$ 
in the ring of integers $\OO_F$ of $F$, we denote that  ideal by ${\p}_k$. 
Similarly, if there is a unique ideal of norm $k$, we denote that ideal by
${\frakn}_k$.

For any ideal $\n$ of $\OO_F$ that is relatively prime to $\p_5$ (resp. $\p_7$),
we define  compact arithmetic hyperbolic $3$-manifolds $W_0(\n)$ 
(respectively $M_0(\n)$) as the quotients of $\H^3$ by the norm one elements of the  level $\n$ Eichler order in the quaternion  algebras over $F$ ramified at the real place and $\p_5$ (resp. $\p_7$).    We give a more geometric discussion of $W_0(1)$ and $M_0(1)$ 
in~\S~\ref{Examplegeometricdesc}; for example $W_0(1)$ is related to the Weeks manifold. 

Specialized to this example, the results and conjectures of this book suggest:
\begin{quote} {\em The first homology
groups of  $ W_0(\p_7\cdot \n)$  and $M_0(\p_5 \cdot \n)$
are closely related  --- although not isomorphic.  } \end{quote}

For example, we compare the
  integral
  homology of  the manifolds
  $W_0({\p}_{7} \cdot {\frakn}_{5049} )$ and $M_0({\p}_5 \cdot {\frakn}_{5049} )$. In order to present the results more succinctly,
  we compute with coefficients in $\Z_{S}:=\Z[\frac{1}{S}]$, where $S$ is divisible by all prime numbers
  $< 40$.
  {\small
  $$\begin{aligned}
  & H_1(W_0({\p}_{7} \cdot {\frakn}_{5049}))  \\
  = 
 &  \ \  (\Z/43)^4 \oplus (\Z/61)^2    \oplus  (\Z/127)
  \oplus (\Z/139)^2 \oplus (\Z/181) \oplus (\Z_S)^{81} \oplus (\Z/67)^2  \\
   & H_1(M_0({\p}_5 \cdot {\frakn}_{5049})) \\
   = 
   &  \ \  (\Z/43)^4 \oplus (\Z/61)^2    \oplus   (\Z/127)
  \oplus (\Z/139)^2 \oplus (\Z/181) \oplus (\Z_S)^{113} \end{aligned}$$
  }

  The surprising similarity between these two groups turns out to be a general phenomenon,
  and this is the main concern of the present paper.

  The reader familiar with the Jacquet--Langlands theorem for $\GL(2)/\Q$ may note a pleasant feature of our situation:  The resemblance between the two groups above is manifest 
  without any computations with Hecke operators:  one sees torsion elements of the same prime order in both groups.  Over $\Q$, the Jacquet--Langlands correspondence implies (among other things) a relationship between the
cohomology of modular curves and Shimura curves at the level of Hecke modules;
but if  one forgets the Hecke structure it asserts nothing more than
an isomorphism between two free finite rank abelian groups.

\subsection{Geometric description of \texorpdfstring{$W_0(1)$}{W_0(1)} and \texorpdfstring{$M_0(1)$}{M_0(1)}.} \label{Examplegeometricdesc}

Although it is not strictly necessary for our purposes, we discuss some of the geometry of the manifolds from the prior example:

 Let $W$ be the Weeks manifold. It is a closed arithmetic hyperbolic
manifold of  volume
$$\mathrm{Vol}(W) = \frac{3 \cdot 23 \sqrt{23} \cdot \zeta_F(2) }{4 \pi^4} = 0.9427073628869 \ldots$$
where $F = \Q(\theta)$ and $\theta^3 - \theta + 1 = 0$. The Weeks manifold is the smallest volume
arithmetic  manifold~\cite{CFJR}, as well as the smallest volume orientable hyperbolic
manifold~\cite{Gabai}.  It may be obtained by $(5,-2)$ and $(5,-1)$ Dehn surgery on the Whitehead link.

Let $M$ be the arithmetic manifold $m017(-3,2)$ from
the  Hodgson--Weeks census.  %
$M$ can be obtained by
 $(7,-2)$ and $(7,-2)$ Dehn surgery on the Whitehead link (see~\cite{Gabai}, p.27).

  \begin{figure}[!ht]
\begin{center}
  \includegraphics[width=90mm]{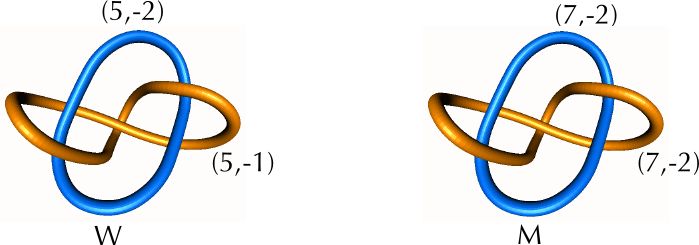}
\end{center}
\caption{The manifolds $W$ and $M$ as surgeries on the Whitehead link complement}
\end{figure}

The manifolds $M$ and $W$ are 
\emph{not} commensurable.
Nonetheless, they are closely related: 
 the invariant trace field of both $M$ and $W$ is $F$; moreover, 
$$ \vol(M) = 2 \  \vol(W), \ \ \mathrm{CS}(M) =   2 \  \mathrm{CS}(W)
 + \frac{1}{12}  \mod \  \frac{1}{2} \cdot \Z,$$
where $\mathrm{CS}$ denotes the Chern--Simons invariant
which depends (up to sign) on a choice
of orientation of the manifold.

The manifolds $W_0(1), M_0(1)$ of the prior example can be described thus:

\begin{itemize}
\item[$W_0(1)$:]
  The manifold $W$ has symmetry group $D_{12}$, the dihedral group of order $12$.
  The Weeks manifold therefore admits a unique quotient $W_0(1)$ of degree three.
  \item[$M_0(1)$:]
  The manifold $M$
has symmetry group $D_8$, and admits a unique orbifold quotient $M_0(1)$ corresponding to the
kernel of the projection $D_8 \rightarrow (\Z/2)^2$.
\end{itemize}
 
As for the manifold  
$W_0({\p}_7 \cdot {\frakn}_{5049})$,
it is a
covering of $W_0(1)$ 
with fiber given by  $\mathbf{P}^1(\OO_F/{\p}_7 \cdot {\frakn}_{5049})$ and monodromy group
  $\SL_2(\OO_F/{\p}_7 \cdot {\frakn}_{5049})$. For $M_0({\p}_5 \cdot {\frakn}_{5049})$ the same
  holds with ${\p}_7 \cdot {\frakn}_{5049}$ replaced by ${\p}_5 \cdot {\frakn}_{5049}$.  

\section{A guide to reading this book}
\label{section:guide}

We now describe the contents of this book (this is not, however, a complete table of contents, simply a guide to some of the important subsections).

\medskip

We hope that at least parts of this book will be accessible
to three audiences --- those interested in hyperbolic geometry, those interested in arithmetic
Langlands, and those interested in analytic torsion. One place to start (besides the introduction)
is Chapter~\ref{part:examples}, which contains some detailed examples.

\begin{enumerate}
\item Chapter~\ref{chapter:ch1} is the introduction.
\begin{itemize}
\item[\S~\ref{section:introtheorems}] Some sample theorems.
\item[\S~\ref{section:example}] An Example: two Dehn surgeries on the Whitehead link.
\item[\S~\ref{section:guide}] A guide to the paper, which you are now reading.
\end{itemize}
\item Chapter~\ref{chapter:ch2} gives a general discussion of the context of the paper -- in particular, extending the Langlands program to torsion classes, and how results such as Theorem~A would fit into that extension. \begin{itemize}
\item[\S~\ref{sec:background}] Background on automorphic forms and ``reciprocity over $\Z$.''
\item[\S~\ref{sec:gl2inner}] Conjectures about reciprocity related to torsion classes for $\GL_2$;
in particular, we formulate an $R=\T$ conjecture that guides many of our computations and conjectures. 
\end{itemize}

\item Chapter 3 gives notation:
\begin{itemize}
\item[\S~\ref{sec:notn0}] Summary of the most important notation; the reader could read this
and refer to the other parts of this Chapter only as needed. 

\item[\S~\ref{sec:notn2}] Notation concerning modular forms, 
and related notions such as Eisenstein classes and congruence homology.\end{itemize}

\item Chapter~\ref{chapter:ch3}, {\em Raising the Level: newforms and oldforms} compares spaces of modular forms at different levels. 

\begin{itemize}
\item[\S~\ref{sec:ihara}]  Ihara's lemma over imaginary quadratic fields. Similar results have also been proven by Klosin. 
\item[\S~\ref{sec:lloneprime}] This section gives certain ``dual'' results to Ihara's lemma. 
\item[\S~\ref{sec:lr}] We prove level raising, an analogue  of a result of Ribet over $\mathbf{Q}$. 

\item[\S~\ref{sec:clss}] We study more carefully how level raising is related
to cohomology of $S$-arithmetic groups. 

In particular 
we are able to produce torsion classes related to $K_2$. This latter result matches with a Galois-theoretic analysis undertaken in~\S~\ref{section:Eisenstein}.

\end{itemize} 
 \item Chapter~\ref{chapter:ch4}, {\em Analytic torsion and regulators}, discusses analytic torsion and the closely 
connected notion of regulator. The regulator measures the volume of the first homology group,
with respect to the metric defined by harmonic forms; it has a substantial influence on the behavior of the torsion homology. 
\begin{itemize}
\item[\S~\ref{sec:CMthm}]   An exposition
of the Cheeger--M{\"u}ller theorem in the present context.   
\item[\S~\ref{regcusp}]  The regulator
and its relation to periods and Faltings heights. 
The main content is an old theorem of Waldspurger. 
\item[\S~\ref{regLproof}] Proof of Theorem~\ref{regulatorL} from the previous section. 
\item[\S~\ref{sec:regchangelevel}] Comparison of regulator at different levels:   the change
of regulator is related to level raising.   \item[\S~\ref{reglowlevel}]  Comparisons of regulators between  a  Jacquet--Langlands 
pair; 
our conjecture is that this is governed by level-{\em lowering} congruences, and we give heuristic reasons for this.   
\end{itemize}

\item Chapter~\ref{chapter:ch5}, {\em The split case},  analyzes various features of the noncompact case. 
In particular, it gives details required to extend the Cheeger--M{\"u}ller theorem in the non-compact case. 

\begin{itemize}
\item[\S~\ref{section:split}] Notation and basic properties of the manifolds under consideration. 
\item[\S~\ref{sec:rtatnc}] Definitions of analytic and Reidemeister torsion in the split case. 
The difficulty is that, for a non-compact Riemannian manifold $M$, cohomology classes
need not be representable by {\em square integrable} harmonic forms.

\item[\S~\ref{sec:Eis-series}] A review of the theory of Eisenstein series, or that part of it which is necessary. This contains, in particular, the definition of the regulator in the non-compact case,
which is used at several prior points in the paper. 

\item[\S \ref{sec:eisintegrality}] Bounds on Eisenstein torsion, and the Eisenstein part of the regulator, in terms of $L$-values.

\item[\S~\ref{analysissec}] This contains the core part of the analysis, and contains gives the proofs.  
The key point is the analysis of small eigenvalues of the Laplacian on the ``truncation''
of a non-compact hyperbolic $3$-manifold.  This analysis nicely matches with the trace formula on the non-compact manifold. 
\end{itemize}

\item Chapter~\ref{chapter:ch6}, {\em  Jacquet--Langlands}, studies ``Jacquet--Langlands pairs'' of hyperbolic $3$-manifolds and in particular proves Theorem~A, Theorem~$A^{\dagger}$, and Theorem~B
quoted in~\S~\ref{section:introtheorems}.
\begin{itemize}

 \item[\S \ref{ss:JLclassic}] Recollections on the classical Jacquet--Langlands correspondence.
\item[\S \ref{sec:newnewnew}]  Some background on the notion of newform. 

\item[\S~\ref{sec:coh}] We state and prove the theorem on comparison of homology
in a very crude form. 
\item[\S~\ref{section:refined}] We show that
that certain volume factors occurring in the comparison theorem correspond exactly to congruence homology. 

\item[ \S~\ref{section:essential}] We introduce the notion of {\em essential homology} and {\em dual-essential homology}: two variants of homology which (in different ways)
``cut out'' congruence homology.  

\item[\S~\ref{section:relations}] We begin to convert the prior Theorems into 
actual comparison theorems between orders of  newforms. We consider
in this section simple  cases in which we can control as many of the auxiliary phenomena (level lowering, level raising, $K$-theoretic classes) as possible.

In particular, we prove Theorems A , $A^{\dagger}$, and B quoted earlier in the introduction. 

\item[\S~\ref{section:nonny} and~\S~\ref{sec:nonny2}] We attempt to generalize
the results of \S~\ref{section:relations} as far as possible.  This uses the spectral sequence of Chapter~\ref{chapter:ch3}; the strength of our results
is unfortunately rather limited because of our poor knowledge about homology of $S$-arithmetic groups.     The results here are summarized in \S \ref{section:finalsummary};
we advise the reader to look at \S \ref{section:finalsummary} first. 
\end{itemize}

\item Chapter~\ref{chapter:ch7}, {\em  Eisenstein deformations, and even Galois representations}: We present some circumstantial evidence for  the conjectures
presented in Chapter 2. 

\begin{itemize}
\item[\S~\ref{involutions}] We analyze the effect of an orientation-reversing involution on homology; in combination with the alternating pairing on torsion homology, it (roughly) splits the torsion homology into two dual pieces.  We call this phenomenona ``doubling.''
\item[\S~\ref{section:Eisenstein}] We analyze the deformation ring of a residually reducible representation and show that this matches with the doubling phenomenon. 
\item[\S~\ref{section:phantomclasses}] We discuss the notion of phantom classes of level $\q$, and how classes from $K_2$ give rise to Eisenstein homology
at level $\q \rf$ for any auxiliary  prime $\rf$, but not necessarily at level $\q$.
\item[\S~\ref{section:even}] We restrict an even residual Galois representation over $\Q$ to an imaginary quadratic field. We show that the resulting deformation ring  again
shows a behavior matching with doubling. 

\item[\S~~\ref{ss:Serre}] We  discuss finiteness theorems
for even Galois representations over $\Q$ that would follow
from modularity over imaginary quadratic fields, and discuss a (as yet unsuccesful)
attempt to locate a  even Galois representation of tame level $1$ with big image.  \end{itemize}

 \item Chapter~\ref{chapter:ch8}, {\em Examples}, presents several numerical examples
 studying torsion cohomology and the Jacquet--Langlands correspondence -- mostly 
in the case of the field $\Q(\sqrt{-2})$. In particular:
 \begin{itemize}
 \item[\S~\ref{sec:mans}] The manifolds under consideration.
 \item[\S~\ref{sec:nocharzero}] The simplest examples of JL correspondence,
 where there are no forms in characteristic zero.
 \item[\S~\ref{sec:oldcharzero}] Examples where there are only oldforms in characteristic zero.
 \item[\S~\ref{sec:newcharzero}] These are the most complicated examples, where
 there exist newforms in characteristic zero; they exhibit a combination of level lowering and level raising that affects torsion.
 \item[\S~\ref{sec:Eisenstein7}] Examples with Eisenstein classes, illustrating
 some of the results of~\S~\ref{sec:clss}.
 \item[\S~\ref{section:phantomtwo}] We compute the primes $\q$ with $N(\q) \le 617$ where we expect (and do) find 
 phantom classes of residue characteristic $p = 3$ and $5$.
 \item[\S~\ref{section:K2examples} ] A numerical example of Theorem~\ref{theorem:K2popularversion} relating
 $K_2(\OL_F)$ to cohomology, in the case $F = \Q(\sqrt{-491})$. 
 \item[\S~\ref{section:tables}] Tables of data illustrating the correspondence, as computed for us by Nathan Dunfield.
 \end{itemize}
 \item Chapter~\ref{chapter:ch9}, {\em Concluding remarks} discusses some potential generalizations
 of our results and conjectures (\S~\ref{higherweights}) as well as some interesting questions
 (\S~\ref{sec:qns}).

\end{enumerate}

\chapter{Some Background and Motivation}

 In this chapter, we assume familiarity with the basic vocabulary of the Langlands program.

Our goal is to formulate precisely conjectures
relating homology to Galois representations --- at least, in the context of inner forms of $\GL_2$ over number fields --- with emphasis on writing them in a generality that applies to torsion homology.  

We then explain why these conjectures suggest that a numerical Jacquet--Langlands correspondence, relating the size of torsion homology groups between two arithmetic manifolds, 
should hold (Lemma~\ref{lem:ConjecturesImplyJL}).

\label{chapter:ch2}

\section{Reciprocity over \texorpdfstring{$\Z$}{Z}}
\label{sec:background}

\subsection{Arithmetic Quotients}  \label{section:cohomological}
We recall here briefly the construction of arithmetic quotients of symmetric space, and the connection to automorphic representations.
 
Let $\G$ be a connected linear reductive
algebraic group over $\Q$, and write $G_{\infty} = \G( \R)$.
Let $A$ be a maximal $\Q$-split torus in the center of $\G$,
and let $A^0_{\infty} \subset G_{\infty}$ denote the connected component of the $\R$-points of $A$.
Furthermore, let $K_{\infty}$ denote a maximal compact of $G_{\infty}$ with
connected component $K^{0}_{\infty}$.

The quotient $ G_{\infty}/K_{\infty}$ carries a $G_{\infty}$-invariant Riemannian metric,
with respect to which it is a Riemannian symmetric space.  
For our purposes it will be slightly better to work with the (possibly disconnected) quotient 
$$S := G_{\infty}/A^0_{\infty} K^0_{\infty},$$
 because
 $G_{\infty}$ always acts on $S$ in an orientation-preserving fashion,   and $S$ does not have a factor of zero curvature.  (Later in the document we will work almost exclusively
 in a case where this coincides with $G_{\infty}/K_{\infty}$, anyway). 
 
For example:

\begin{itemize}
\item If $G_{\infty} = \GL_2(\R)$, then $S$ is the union of the upper- and lower- half plane, i.e. $S = \mathbf{C}-\mathbf{R}$  with invariant  metric $\displaystyle{\frac{|dz|^2}{\mathrm{Im}(z)^2}}$. 
\item If $G_{\infty} = \mathrm{SO}(n,1)(\R)$, then $S$ is isometric to hyperbolic $n$-space $\H^n$.
\item If $G_{\infty} = \GL_2(\C)$, then $S = \H^3$. 
\end{itemize}

The symmetric space $\H^3$ occurs in both
 contexts because of the exceptional isomorphism
 $\mathfrak{sl}_2(\C) \simeq \mathfrak{so}(3,1)(\R)$.

Now let $\Adele$ be the adele ring of $\Q$,  and
$\Afinite$ the finite  adeles.  
For any compact open subgroup $K$ of
$\G(\Afinite)$, we may define an ``arithmetic manifold'' (or rather ``arithmetic orbifold'')
$Y(K)$ as follows:
$$Y(K):=  \G(F) \backslash (S \times \G(\Afinite))/K = \G(F) \backslash \G(\Adele)/A^0_{\infty} K^0_{\infty} K.$$ The orbifold $Y(K)$ may or may
not be compact, and may be disconnected; in all cases, it has finite volume with respect to the measure defined by the metric structure.

Since the components of $S$ are contractible, the connected components
of $Y(K)$ will be $K(\pi,1)$-spaces, and the cohomology of each component of $Y(K)$ is canonically isomorphic to
the group cohomology of the fundamental group.

 \subsection{Cohomology of arithmetic quotients and reciprocity.}
 \label{section:matsu}

Assume for now that $Y(K)$ is compact.
 Matsushima's formula~\cite{BW} states that the cohomology
 $H^i(Y(K),\C)$ decomposes as a direct sum of the  $(\mathfrak{g},K)$-cohomology:
 $$  H^i(Y(K),\C) = \bigoplus m(\pi) H^i(\mathfrak{g},K;\pi_{\infty}).$$
 Here the right hand side sum is taken over automorphic representations $\pi = \pi_{\infty} \otimes \pi_{\finite}$
 for $\G(\adele)=G_{\infty} \times \G(\Afinite)$, and $m(\pi)$ is the dimension of the $K$-fixed subspace in $\pi_{\finite}$.
 The effective content of this equation is that ``cohomology can be represented by harmonic forms.''

 Appropriate variations of these results apply equally with coefficient systems;
 in the non-compact case the work of Franke~\cite{Fr} provides a
 (somewhat more involved) substitute.

 \subsection{Reciprocity in characteristic zero.}
 
 Let us suppose $\pi$ is a representation that ``contributes to cohomology,''
 that is to say, $m(\pi) \neq 0$ and $H^i(\mathfrak{g},K;\pi_{\infty}) \neq 0$
 for suitable $K$.  Assume moreover that $\G$ is {\em simply connected}.
 Let $\widehat{G}$ denote the dual group of $\G$, considered as a reductive algebraic group over
 $\Qbar$, and let ${}^L G$ be  the semi-direct
 product of $\widehat{G}$ with $\Gal(\Qbar/\Q)$, the action
 of $\Gal(\Qbar/\Q)$ being the standard one (see e.g. \cite[\S 2.1]{BG}). 
  Then one conjectures (see, for example, Conjecture~3.2.1 of~\cite{BG} for
  a more throrough and precise treatment)
the existence (for each $p$) of a continuous irreducible Galois representation
$$\rho = \rec(\pi): \Gal(\overline{\Q}/\Q) \rightarrow {}^{L}{G}(\Qbar_p),$$
such that
\begin{enumerate}
\item For every prime $\ell \neq p$, the representation $\rho|G_{\ell}$ is associated to $\pi_\ell$
via the local Langlands correspondence; 
\item $\rho$ is unramified outside finitely many primes; 
\item If $c$ is a complex conjugation, then $\rho(c)$ is {\em odd} (see \cite[Proposition 6.1]{AB}
for precise definition;  if $\G$ is split, it means that the trace of $\rho(c)$
in the adjoint representation should be minimal among all involutions).

\item
$\rho|G_{F_{p}}$ is de Rham for any place $v|p$,  with prescribed Hodge-Tate weights:
they correspond to the conjugacy class of cocharacter $\mathbb{G}_m \rightarrow \widehat{G}$
 that are dual to the half-sum of positive roots for $G$.  
 \end{enumerate}
For a more precise formulation (including what is meant by ``de Rham'' for
targets besides $\GL_n$) the reader should consult~\cite{BG}.

\subsection{Cohomological reciprocity over \texorpdfstring{$\Z$}{Z}}  \label{Z:cores}

Starting with observations at least as far back as Grunewald~\cite{Grunewald}
in 1972 (continuing with  further work of Grunewald, Helling and Mennicke~\cite{Germans}, 
Ash~\cite{Ash}, Figueiredo~\cite{Fig}, and many others),
it has become apparent that something much more general  should be true. 

Let us suppose now that $\G$ is split simply connected, with complex dual group $\Gdual(\C)$. 
Recall that for almost every prime $\ell$ the ``Hecke algebra'' $\mathscr{H}_{\ell}$ 
  acts by correspondences on $Y(K)$; 
   moreover~\cite{Gross} $$\mathscr{H}_{\ell} \otimes \Z[\ell^{-1}] \cong \mathrm{Rep}(\widehat{G}(\C)) \otimes \Z[\ell^{-1}].$$     Accordingly, if  $\sigma$ is any representation of $\widehat{G}(\C)$, let $T_{\sigma}(\ell)$
be the corresponding element of $\mathscr{H}_{\ell} \otimes \Z[\ell^{-1}]$.  For example,
if $\G = \SL_2$, then $\widehat{G}(\C) = \PGL_2(\C)$; the element $T_{\mathrm{Ad}}$
associated to the adjoint representation of $\widehat{G}(\C)$ is the Hecke operator usually denoted $T_{\ell^2}$, of degree $\ell^2+\ell+1$. 

We may then rephrase the classical reciprocity conjecture in the less precise form:

\begin{itemize}
\item To every {\em Hecke eigenclass} 
$\alpha \in H^i(Y(K), \Qbar_p)$, there exists a matching Galois representation
$\rho_{\alpha}: \Gal(\Qbar/\Q) \rightarrow\Gdual(\Qbar_p)$, so  conditions  
(2) --- (4) above hold, and for every $\sigma$ and almost every $\ell$, 
$$ T_{\sigma}(\ell) \alpha = \mathrm{trace} (\sigma \circ \rho_{\alpha} (\mathrm{Frob}_{\ell}) ).$$  
 Conversely,  every 
Galois representation satisfying conditions (2), (3), (4) arises thus (for some $i$). 
\end{itemize}

One difficulty with this conjecture is that it only (conjecturally) explains Galois
representations in characteristic zero. Thus it is of interest to replace the coefficients $\Qbar_p$ of
cohomology $H^i(Y(K), \Qbar_p)$ with torsion coefficients. 
For cleanness of formulation it is simplest to replace
$\Qbar_p$ above by $\Fbar_p$, or, more generally, a finite  Artinian local $W(\Fbar_p)$-algebra
$A$. The following is the ``torsion'' counterpart to the conjecture above:

\begin{itemize}
\item   Let $A$ be a finite length local $W(\Fbar_p)$-algebra.
To every {\em Hecke eigenclass} 
$\alpha$ in the group $H^i(Y(K), A)$, there exists a matching Galois representation
$\rho_{\alpha}$ with image in  $ \widehat{G}(A)$ such that conditions (2), (4) above hold, 
and, for every $\sigma$ and almost every $\ell$, 
$$ T_{\sigma}(\ell) \alpha = \mathrm{trace} (\sigma \circ \rho_{\alpha} (\mathrm{Frob}_{\ell}) ).$$ %
 Conversely,  every 
Galois representation satisfying conditions (2),   (4) arises thus  (for some $i$).
 \end{itemize}

Note that we have deemed requirement $(3)$ to be empty.
By the universal coefficient theorem,  the cohomology with coefficients in $A$ 
(for any $A$) is determined by the cohomology with $\Z$-coefficients,
and this motivates a study of the integral
cohomology groups  
 $H^{\bullet}(Y(K),\Z)$.   We therefore refer to these conjectures together as ``reciprocity over $\Z$.''
 Since   $H^{\bullet}(Y(K),\Z)$ may possess torsion,
reciprocity over $\Z$ is {\em not} a consequence of the first conjecture alone.
  
 The majority of the numerical evidence for this conjecture --- particularly in the cases where it involves
mod-$p$ torsion classes  and $A$ is a finite field ---  is due to A. Ash and his collaborators;
see \cite{Ash}.
 \medskip

 Since the cohomology groups of $Y(K)$ are  finitely generated,
 any eigenclass $\alpha \in H^i(Y(K),\Qbar_p)$ will actually come from
 $H^i(Y(K),L)$ for some finite extension $L/\Q_p$ (a similar remark applies to eigenclasses
 in $W(\Fbar_p)$-algebras). However, the conjectures are more naturally stated over
 $\Qbar_p$ than over $L$,
 since even if $\alpha \in H^i(Y(K),L)$, and
 even if there is a natural choice of dual group $\widehat{G}$ over $\Q$,
 it may be the case that $\rho_{\alpha}$ may only
 be conjugated to lie 
 in $\widehat{G}(L')$ for some nontrivial extension $L'/L$.
 One can already see this for representations of finite groups, for example,
 two dimensional representations of the quaternion group have traces in $\Q_p$
 for all $p$, but
 can not be conjugated into $\GL_2(\Q_p)$ unless $p \equiv 1 \mod 4$.

 \smallskip
 
There is some flexibility as to  whether we  work with integral \emph{homology} or integral 
\emph{cohomology}. It will also be useful at some points to work with $\Q/\Z$-coefficients.
The universal coefficient theorem relates these groups, and there are various ways in which
to phrase the conjectures. For example, if $Y=Y(K)$ happens to be a compact $3$-manifold, the torsion classes that don't arise from characteristic zero classes
live in $H_1(Y,\Z)$, $H^2(Y,\Z)$,  $H^1(Y,\Q/\Z)$ and $H_2(Y,\Q/\Z)$.
For the purposes of this book, it will be most convenient (aesthetically) to work with $H_1(Y,\Z)$, and hence we phrase
our conjectures in terms of homology.

\smallskip

One may also ask which {\em degree} the homology classes of interest live in.
This is a complicated question; the paper~\cite{CE1} 
suggest that the Hecke algebra  completed at a maximal ideal corresponding to a representation with large image 
 acts faithfully in degree
$ \frac{1}{2}(\dim(G_{\infty}/K_{\infty}) - \mathrm{rank}(G_\infty) + \mathrm{rank}(K_{\infty}))$.

\subsection{Discussion}

 \label{sophistry}

One may well ask if it is worth the effort to study reciprocity integrally when
many of the naturally occuring Galois representations of interest
(those attached to motives, for example) arise in characteristic zero.
On an even more elementary level, one might ask whether
 reciprocity over $\Z$ is significantly different from reciprocity over $\Q$.
 Let us explain why the answer to the second question is ``yes'', and explain
 why this provides several justifications for studying torsion classes.

  The main difference between torsion classes and characteristic zero
 classes occurs when $G_{\infty}$ does \emph{not} admit discrete series
 (although the study of torsion in the cohomology of Shimura varieties is
 also very interesting!).
 Whenever $G_{\infty}$ does
\emph{not} admit discrete series, the set of automorphic Galois
representations will not be Zariski
dense in the deformation space of  all (odd) $p$-adic Galois representations into
${}^L G$. When $\G = \GL(2)$ over an imaginary quadratic field $F$, this follows
from Theorem~7.1 of~\cite{CM}, and Ramakrishna's arguments are presumably sufficiently
general to establish a similar theorem in the general setting where $G_{\infty}$ does not have
discrete series.
Thus, if one believes that the collection of \emph{all} Galois
representations is an object of interest, the representations obtained
either geometrically (via the cohomology of a variety) or automorphically
(via a  classical automorphic form)
 do not suffice to study them.

An interesting example is obtained by taking an {\em even Galois representation}
$\rhobar: \mathrm{Gal}(\overline{\Q}/\Q) \rightarrow \GL_2(\mathbf{F}_p)$. 
While there is no known ``modular'' object over $\Q$ that is associated to $\rhobar$,
the restriction $\rhobar|G_F$ to the Galois group of an imaginary quadratic field $F$ should be,
by ``reciprocity over $\Z$'', associated to an integral cohomology class for a suitable Bianchi group (cf. Chapter 9); 
and there is no reason to expect that this class lifts to characteristic zero.
(One knows in many cases that there are no geometric lifts of $\rhobar$ over $\Q$ ---
see~\cite{CalEven2}).

On the other hand, suppose that one is only interested in geometric Galois representations.
One obstruction to proving modularity results for such representations
 is apparent failure of the Taylor--Wiles
method when $G_{\infty}$ does not have discrete series.
A recent approach to circumventing these obstructions can be found in the work of the first author
and David Geraghty~\cite{CalG}. One of the main theorems of~\cite{CalG} is to prove minimal modularity
lifting theorems for imaginary quadratic fields 
\emph{contingent} on the existence of certain Galois representations.
 It is \emph{crucial} in the approach of \emph{ibid}. that one work
with the \emph{all} the cohomology classes over $\Z$, not merely those which lift to characteristic zero.  
In light of this,  it is important to study  the nature of torsion classes even
if the goal is ultimately only to study motives.

Having decided that the (integral) cohomology of arithmetic groups should play
an important role in the study of reciprocity and Galois representations, it is natural
to ask what the role of  Langlands principle of \emph{functoriality} will be.
Reciprocity compels us to believe the existence of (some form of) functoriality on the level
of the cohomology of arithmetic groups.
One of the central goals of this book is {\em to give some evidence towards arithmetic functoriality
in a context in which it is not a consequence of a classical functoriality for
automorphic forms.  }

\section{Inner forms of \texorpdfstring{$\GL(2)$}{GL(2)}: conjectures } \label{sec:gl2inner}

We now formulate more precise conjectures in
the case when $\G$ is an inner form of $\PGL(2)$ or $\GL(2)$ over a number field $F$. 
In most of this book we shall be working with $\PGL(2)$ rather than $\GL(2)$, but at certain points (e.g.~\S~\ref{section:even}) it will be convenient to have the conjectures formulated for $\GL(2)$. 
 
In particular, we formulate conjectures concerning $R=\T$, multiplicity one, and Jacquet--Langlands. 
 We will not strive for any semblance
of generality: we concentrate on two  cases corresponding to the trivial local system, and we shall only consider the semistable case, for ease of notation as much
as anything else.
The originality of these conjectures (if any) consists of 
conjecturally identifying the action of the entire Hecke
algebra $\T$ of endomorphisms
 on cohomology in terms of Galois representations, not merely $\T \otimes \Q$ (following
Langlands, Clozel, and others), nor the finite field quotients $\T/\m$ (generazing Serre's
conjecture --- for the specific
case of imaginary quadratic fields, this is considered
in work of Figueiredo~\cite{Fig}, Torrey~\cite{Torrey}, and also discussed
in unpublished notes of Clozel~\cite{ClozelA,ClozelB}).

\subsection{Universal Deformation Rings} \label{section:universal} 

Let $F$ be a number field and let $\Sigma$  denote a finite set of primes in $\OL_F$

Let $k$ be a finite field of characteristic $p$, and let
$$\rhobar: G_{F} \rightarrow \GL_2(k)$$
be a continuous absolutely irreducible odd Galois representation
(we also consider reducible representations in~\S~\ref{section:Eisenstein}).
Recall that oddness in this context means that $\rhobar(c_v)$ has determinant $-1$
for any complex conjugation $c_v$ corresponding to a real place $v$ of $F$; no 
 condition  is imposed at the complex places.

Let $\eps : G_F \rightarrow \Z_p^{\times}$ denote the cyclotomic character, and for any place $v$, let $I_v \subset G_v$
denote inertia and decomposition subgroups of $G_F$, respectively.  %

Let us suppose for all $v|p$ the following is satisfied:
\begin{enumerate}
\item $\det(\rhobar)= \eps \phi$, where $\phi |G_v$ is trivial for $v|p$.  If $\G =\PGL_2$ we require $\phi \equiv 1$. 
\item  For $v|p$ and $v \notin \Sigma$, either:
\begin{enumerate}
\item $\rhobar|_{G_v}$ is finite flat,
\item $\rhobar|_{G_v}$ is ordinary, but not finite flat.
\end{enumerate}
\item For $v \in \Sigma$, $\rhobar|_{G_v} \sim \left(\begin{matrix}
 \eps \chi &  *  \\ 0 & \chi \end{matrix}\right)$, where $\chi$ is an unramified character
 of $G_{v}$.
 \item For $v \notin \Sigma \cup \{v | p\}  $, $\rhobar|G_v$ is unramified.
 \end{enumerate}
Associated to $\rhobar$, we can, in the usual way (see~\cite{Fermat,Mazur,Wiles})
define a deformation problem that associates to any local  Artinian ring $(A,\m,k)$
deformations $\rho$ of $\rhobar$ such that:
\begin{enumerate}
\item 
$\det(\rho) = \eps \cdot \widetilde{\phi}$, where $\phi$ is the Teichmuller lift of  $\phi$.
\item For $v|p$ and $v \notin \Sigma$, either:
\begin{enumerate}
\item $\rho|G_v$ is 
finite flat,
\item $\rho|G_v$ is ordinary, and $\rhobar|G_v$ is not finite flat.
\end{enumerate}
\item For $v \in \Sigma$,  $\rho|G_v \sim  \left(\begin{matrix}
 \eps \widetilde{\chi} &  *  \\ 0 & \widetilde{\chi}  \end{matrix}\right)$,  where $\widetilde{\chi}$
 is the Teichmuller lift of $\chi$.
 \item If $v$ is a real place of $F$, then $\rho|G_v$ is odd. (This is only a nontrivial condition
 if $\mathrm{char}(k) = 2$).
  \item If $v \notin \Sigma \cup \{v|p\} $, $\rho|G_v$ is unramified.
 \end{enumerate}
We obtain in this way a universal deformation ring $R_{\Sigma}$.
One may determine the na\"{\i}ve ``expected dimension''
of these algebras (more precisely, their relative dimension over $\Z_p$) using the
global Euler characteristic formula. This  expected dimension is
$-r_2$, where $r_2$ is the number of complex places of $\OL_F$. Equivalently, one expects
that these rings are finite over $\Z_p$, and finite if $F$ has at least
one complex place. This computation can not be taken too seriously, however,
as it is inconsistent with the existence of elliptic curves over any field $F$.

\subsection{Cohomology} \label{section:ar1}

Let  $\G/F$ be an inner form of $\GL(2)$
corresponding to a quaternion algebra ramified at some set of primes 
$S \subset \Sigma \cup \{v| \infty\}$. Let $K_{\Sigma} = \prod_{v} K_{\Sigma, v}$ denote an open compact subgroup
of $\G(\Afinite)$  such that:
\begin{enumerate}
\item  If $v \in \Sigma$ but $v \notin S$,  $K_{\Sigma,v} =  \left\{ g \in \GL_2(\OL_{v}) \ | \
g \equiv 
\left(\begin{matrix} * & * \\ 0 & * \end{matrix} \right) \mod \pi_{v} \right\}$.
\item If $v \in S$, then $K_{\Sigma, v}$ consists of the image
of $B_v^{\times}$ inside $\G(F_v)$, where $B_v$ is a maximal order of the underlying quaternion algebra. 
\item If $v \notin \Sigma$, then 
$K_{\Sigma,v} = \G(\OL_v)$.
\end{enumerate}
Moreover, if $\q$ is a finite place in $\Sigma$ that is not in $S$ (so that situation $1$ above
occurs), let $K_{\Sigma/\q}$ be defined as $K_{\Sigma}$ except that $K_{\Sigma/\q}$
is maximal at $\q$.

\medskip

Let $Y(K_{\Sigma})$ be the locally symmetric space associated with $\mathbb{G}$.
If $\mathbb{G}$ is an inner form of $\PGL_2$, we are identifying $K_{\Sigma}$ with its image in 
$\PGL_2(\Afinite)$; each component of $Y$ is then uniformized by a product of hyperbolic $2$- and $3$-spaces. If $\mathbb{G}$ is an inner form of $\GL_2$, 
then $Y(K_{\Sigma})$ is a torus bundle over the corresponding $\PGL_2$-arithmetic manifold. 

 Let $q$ be the number of archimedean places for which $\G \times_F F_v$ is split.

\begin{df} Let $\Psi^{\vee}: H_{q}(Y(K_{\Sigma/\q}),\Z)^2 \rightarrow H_{q}(Y(K_{\Sigma}),\Z)$ be the
transfer map (see~\S~\ref{section:levelmaps}). 
Then the space of new-at-$\q$ forms of level $K_{\Sigma}$ is defined to be  the
group $H_{q}(Y(K_{\Sigma}),\Z)_{\qnew}:= \coker(\Psi^{\vee})$.
\end{df}

We note that there are other candidate definitions for the space of newforms ``over $\Z$''; the definition
above appears to work best in our applications. 

There is an obvious extension of this definition to the space of forms that are new at any
set of primes in $\Sigma \setminus S$. Let us refer to the space of \emph{newforms} as the
quotient   $H_{q}(Y(K_{\Sigma}),\Z)^{\new}$  by the image of the transfer map
for \emph{all} $\q$ in $\Sigma \setminus S$.

If $Y(K_{\Sigma})$ and $Y(K_{\Sigma/\q})$ are the corresponding arithmetic quotients, then
the integral homology and cohomology groups  $H_{\bullet}(Y,\Z)$ and $H^{\bullet}(Y,\Z)$ admit an action by Hecke endomorphisms (see~\S~\ref{ss:HeckeDef}  for the
definition of the ring of Hecke operators).
We may now define the corresponding Hecke algebras.

\begin{df}

Let $\T_{\Sigma}$ --- respectively the {\em new Hecke algebra} $\Tnew_{\Sigma}$ --- 
denote the subring of $\End \ H_{q}(Y(\Sigma),\Z)$  $($respectively  $\End \ H_{q}(Y(K_{\Sigma}),\Z)^{\new})$
generated by Hecke endomorphisms.  \label{df:hecke}
 \end{df}

It is easy to see that $\Tnew_{\Sigma}$ is a quotient of $\T_{\Sigma}$.
Of course, if $\Sigma = S$, then $\Tnew_{\Sigma} = \T_{\Sigma}$. 

We say that a maximal ideal $\m$ of $\T_{\Sigma}$ is Eisenstein \index{Eisenstein}
if the image of $\T_{\q}$  in the field $\T_{\Sigma}/\m$ is given by the sum of
two Hecke characters evaluated at $\q$,  for all but finitely many $\q$
(see Definition~\ref{df:Eisenstein} for a more precise discussion of 
the various possible definitions of what it means to be Eisenstein).

\label{section:recipe2}
\begin{conj} [$R=\T$] Let $\m$ be a maximal ideal of  \label{conj:serre}
$\Tnew_{\Sigma}$, with $k = \T_{\Sigma}/\m$ a field of characteristic $p$,
 and suppose that $\m$ not Eisenstein.
 Then:
 \begin{enumerate}
 \item There exists a Galois representation
 $$\rho_{\m}: G_F \rightarrow \GL_2(\Tnew_{\Sigma,\m})$$
 satisfying the usual compatibility between Hecke eigenvalues of the characteristic
 polynomial of Frobenius.
\item The induced map %
 $R_{\Sigma} \rightarrow \Tnew_{\Sigma,\m}$ is an isomorphism.
 \end{enumerate}
\label{conj:reciprocity}
\end{conj}

\begin{remarkable} \emph{ In a recent preprint, the first author and Geraghty have proved a minimal modularity lifting theorem
for imaginary quadratic fields $F$
under the assumption that there exists a map from $R_{\Sigma}$ to $\T_{\m}$ (see~\cite{CalG}). This provides
strong evidence for this conjecture. Moreover, in the minimal case, one also deduces (under the conjectural
existence of Galois representations) that the Hecke algebra is free as an $R_{\Sigma}$-module, and is multiplicity
one whenever it is infinite (see Conjecture~\ref{conj:mult}). \label{remark:geraghty}
}
\end{remarkable}

\subsection{Jacquet--Langlands correspondence} \label{ss:jlc}
\label{section:discusscongruence}

\index{Jacquet--Langlands pair}

Let $F/\Q$
be a field with one complex place, 
let  $\G/F$ and $\G'/F$ be inner forms of $\PGL(2)/F$ that  are ramified at the set of finite primes $S$ and $S'$ respectively and ramified at all real infinite places of $F$.
 Let $Y(K_{\Sigma})$ be the arithmetic manifold
constructed as in~\S~\ref{section:ar1}, and denote by $Y'(K_{\Sigma})$ the analogous
construction for $\G'$. 
We refer to a pair of such manifolds as a {\em Jacquet--Langlands pair.}

It follows from the Jacquet--Langlands correspondence that
there is a Hecke-equivariant isomorphism
$$H_1(Y(K_{\Sigma}),\C)^{\new} \cong H_1(Y'(K_{\Sigma}),\C)^{\new}.$$

One is tempted to conjecture the following 
$$|H_1(Y(K_{\Sigma}),\Z)^{\new}| =^{?}  |H_1(Y'(K_{\Sigma}),\Z)^{\new}|,$$
and, even more ambitiously, a corresponding isomorphism of Hecke modules.
We shall see, however, that this is \emph{not} the case.

Consider, first, an example in the context of the upper half-plane: Let $N \in \Z$ be prime, $N \equiv 1 \mod p$,  and consider
the group $\Gamma_0(N)  \subset \PSL_2(\Z)$. Then there is a map
$\Gamma_0(N) \rightarrow \Gamma_0(N)/\Gamma_1(N) \rightarrow \Z/p\Z$,
giving rise   to a class in $H^1(\Gamma_0(N),\Z/p\Z)$.
On the other hand, we see that this class capitulates
when we pull back to a congruence cover (compare the annihilation
of the Shimura subgroup under the map of Jacobians $J_0(N) \rightarrow J_1(N)$).

Similar situations arise in the context of the $Y(K_{\Sigma})$. 
We denote the cohomology that arises in this way \emph{congruence}, because it
vanishes upon restriction to a congruence subgroup.  See~\ref{s:congess} for a precise definition. 
 A simple argument (Lemma~\ref{cong-is-Eis}) shows that all congruence classes are necessarily Eisenstein.

For $3$-manifolds, it may be the case that congruence homology
does not lift to characteristic zero.  
Let us return to the examples of the first section. We observe that the Weeks
manifold $W$ satisfies $H_1(W,\Z) = \Z/5 \oplus \Z/5$, whereas
$M$ satisfies  $H_1(M,\Z) = \Z/7 \oplus \Z/7$.
 In the case of the Weeks manifold, the fundamental
group injects (and has Zariski dense image) into the unique index $6$ subgroup of
the norm one units of the unique quaternion algebra over $\Q_5$. In turn, this subgroup
is easily seen to admit a surjection onto $\Z/5 \oplus \Z/5$,  which ``explains''
the homology of $W$, that is, all the homology of $W$ is congruence.
A similar phenomenon (over $\Q_7$) explains all the homology of $M$.
Since these classes are invisible in the adelic limit (i.e., they become
trivial upon increasing the level), it is not at all clear how they interact
with the classical theory of automorphic forms.

We define the \emph{essential}
homology, denoted $H^E_1$, by excising congruence homology (see Definition~\ref{df:essentialhomology}) Although this definition is perhaps is ad hoc, our first impression was that the
following equality should hold:
$$|H^E_1(Y(K_{\Sigma}),\Z)^{\new}| =^{?} |H^E_1(Y'(K_{\Sigma}),\Z)^{\new}|.$$
However, even this conjecture is too strong --- we shall see that ``corrections''
to this equality may occur arising from $K$-theoretic classes.
(In particular, see Theorem~$A^{\dagger}$ from the introduction,
and~\S~\ref{section:K2examples}, which gives an explicit example
where these two spaces are finite and have different orders. One should also compare
Lemma~\ref{lemma:differentnew}, which gives another example
where these two spaces are finite and have different orders, and
\S~\ref{section:phantomclasses} for some theoretical musings on these issues. 
We leave open the possibility that this equality might always hold when $D$
and $D'$ are both non-split, equivalently, $Y$ and $Y'$ are both compact, although there
is no reason to suppose it should be true in that case either.)
In the examples in which these groups have different orders, the disparity arises from
certain unusual Eisenstein classes. Thus, one may salvage a plausible conjecture by localizing the situation
at a maximal ideal of the Hecke algebra.

\begin{conj} If $\m$ is not  Eisenstein,
there is an equality \label{conj:equality}
$$|H_1(Y(K_{\Sigma}),\Z)^{\new}_{\m}| = |H_1(Y'(K_{\Sigma}),\Z)^{\new}_{\m}|.$$
\end{conj}

  Let us   try to explain why a conjecture such as~\ref{conj:equality}
 might be expected. 
Let  $\rhobar: G_F \rightarrow \GL_2(k)$ be an
absolutely irreducible Galois
representation giving rise to a universal deformation ring $R_{\Sigma}$ as in 
\S~\ref{section:universal}. Let $Y(K_{\Sigma})$ be the arithmetic manifold
constructed as in~\S~\ref{section:ar1}, and denote by $Y'(K_{\Sigma})$ the analogous
construction for $\G'$. If Conjecture~\ref{conj:reciprocity} holds, then we obtain an
identification of $R_{\Sigma}$ with the ring of Hecke endomorphisms acting on both
 $H_1(Y(K_{\Sigma}),\Z)^{\new}_{\m}$ and $H_1(Y'(K_{\Sigma}),\Z)^{\new}_{\m}$.

What is the structure of these groups as Hecke modules? One might first guess that they
are isomorphic. However, this hope is dashed even for $\GL(2)/\Q$, as evidenced by the failure
(in certain situations)
of multiplicity one for non-split quaternion algebras as shown by Ribet~\cite{Ribet}. In the context
of the $p$-adic Langlands program, this failure is not an accident, but rather a consequence of
 local-global compatibility, and should be governed purely by local properties of $\rhobar$.
 A similar remarks apply to multiplicites for $p = 2$ on the split side when $\rhobar$ is unramified
 at $2$ and $\rhobar(\Frob_2)$ is a scalar. In this spirit, we propose the following:

\begin{conj}[multiplicity one] \label{conj:multone} Let 
 $\rhobar: G_F \rightarrow \GL_2(k)$ be an  absolutely irreducible Galois representation, 
and let $\m$ be the associated maximal ideal.
Suppose that:
\begin{enumerate}
\item $p$ divides neither $2$ nor $d_F$.
\item If $\q \in S$, then either:
\begin{enumerate}
\item  $\rhobar|G_{\q}$ is ramified at $\q$, or,
\item $\rhobar|G_{\q}$ is unramified at $\q$, and
$\rhobar(\Frob_{\q})$ is not a scalar.
\end{enumerate}
\end{enumerate} \label{conj:mult}
Then  $H_1(Y(K_{\Sigma}),\Z)^{\new}_{\m}$ is free of rank one over $\Tnew_{\m}$.
\end{conj}

One might argue that condition $1$ is somewhat timid, since (as for $\GL(2)/\Q$) when
$p = 2$ and $2 \nmid d_F$, one should also have multiplicity one under some mild conditions on
$\rhobar|G_v$ for $v|2$, and  when $p|d_F$ and $F/\Q$ is not too ramified at $p$, multiplicity one should also hold.
(See also Remark~\ref{remark:geraghty}.)
However, our only real data along these lines is for $p = 2$ and
$F = \Q(\sqrt{-2})$, when the failure of multiplicity one is frequent.
A trivial consequence of Conjectures~\ref{conj:mult} and~\ref{conj:reciprocity} is the following.
\begin{lemma} \label{lem:ConjecturesImplyJL} Assume Conjectures~\ref{conj:mult} and~\ref{conj:reciprocity}.
Suppose that:
\begin{enumerate}
\item $p$ divides neither $2$ nor $d_F$.
\item If $\q \in S$ or $S'$, then either:
\begin{enumerate}
\item  $\rhobar|G_{\q}$ is ramified at $\q$, or,
\item $\rhobar|G_{\q}$ is unramified at $\q$, and
$\rhobar(\Frob_{\q})$ is not a scalar.
\end{enumerate}
\end{enumerate} \label{lemma:equi}
Then there is an isomorphism 
 $$H_1(Y(K_{\Sigma}),\Z)^{\new}_{\m}\simeq H_1(Y'(K_{\Sigma}),\Z)^{\new}_{\m}.$$
\end{lemma}

Let us summarize where we are: %
Later in the paper, we will prove theorems relating the size of $H_1(Y(K_{\Sigma}), \Z)$ and $H_1(Y'(K_{\Sigma}, \Z))$. What we have shown
up to here is that,  at least if one completes at non-Eisenstein maximal ideals, such a result {\em is to be expected}
if one admits the multiplicity one conjecture Conjectures~\ref{conj:mult}
and the reciprocity conjecture Conjecture~\ref{conj:reciprocity} noted above. In this way, 
we   view the results of this book as giving evidence for these conjectures.

   It is an interesting  question to formulate a  conjecture on  how the $\m$-adic integral
homology should be related when multiplicity one (i.e. Conjecture~\ref{conj:mult}) fails. Asking that both sides
are \emph{isogenous} is equivalent to the rational  Jacquet--Langlands (a theorem!)
since any torsion can be absorbed by an isogeny. It should at least be the case that
one space of newforms  is zero if and only if the other space is.

\subsection{Auxiliary level structure; other local systems} \label{ls1} 
One can generalize the foregoing conjectures to add level structures at auxiliary primes. 
That is, we may add level structure at some set $N$ of primes (not containing any
primes in $\Sigma$ or dividing $p$).  

Similarly, one can consider not simply homology with constant coefficients, 
but with values in a local system. 

We will not explicitly do so in this document; the modifications are routine
and do not seem to add any new content to our theorems. (This is perhaps not exactly true as regards nontrivial local systems: In some respects our theorems would become better. See
\S~\ref{chapter:ch9}.)

\chapter{Notation} 
\label{sec:notn1}

 The most important notation necessary for browsing this book --- 
 the manifolds $Y(K)$ we consider, and the notation for their homology groups --- is contained in 
 Sections~\S~\ref{sec:notn0}. 
 
 The remainder of the Chapter contains more specialized notation as well as recalling various background results. 
 
 The reader might consult these other sections as necessary when reading the text. 
\section{A summary of important notation} \label{sec:notn0}
For $A$ an abelian group, we denote by  $A_{\tf}$ the  maximal \index{$\tf$} \index{$\divi$} 
\index{$\tors$}
torsion-free quotient of $A$ and by $A_{\divi}$ the maximal divisible subgroup of $A$. 
We shall use the notation $A_{\tors}$ in two contexts:
if $A$ is finitely generated over $\Z$ or $\Z_p$, it will mean the
torsion subgroup of $A$; on the other hand, if $ \Hom(A,\Q/\Z)$ or $\Hom(A,\Q_p/\Z_p)$ 
is finitely generated over $\Z_p$, it will mean the maximal torsion {\em quotient} of $A$. 
These notations are compatible when both apply at once.

Let $F$ be a number field with 
one complex place and $r_1 = [F:\Q]-2$ real places.   We denote the ring of integers of $F$ by
  $\order$ and set $F_{\infty} := F \otimes \R \cong \C \times \R^{r_1}$. 
  We fix, once and for all, a complex embedding $F \hookrightarrow \C$ 
  (this is unique up to complex conjugation).
  The adele ring of $F$ will be denoted $\Adele$ (occasionally by $\Adele_F$ when we wish to emphasize the base field $F$)
   and the ring of {\em finite}
adeles $\Afinite$.

Let $D$ be a quaternion algebra over $F$
ramified at a set $S$ of places that contains all real places. 
(By abuse of notation we occasionally use $S$
to denote the set of finite places at which $D$ ramifies;
this abuse should be clear from context.) 

Let $\G$ be the algebraic group over $\Q$ defined as $$\G = \mathrm{Res}_{F/\Q} \GL_1(D)/\mathbb{G}_m.$$
 Let $G_{\infty} = \G(\R)$ and
let $K_{\infty}$ be a maximal compact of $G_{\infty}$. The associated symmetric space $G_{\infty}/K_{\infty}$ is isometric to the hyperbolic $3$-space $\H^3$; we fix such an  isometric identification.

For $\Sigma \supset S$ a subset of the finite places of $F$, we define
the associated arithmetic quotient, a (possibly disconnected) hyperbolic three  orbifold
of finite volume: \index{$Y(\Sigma)$} \index{$Y(K_{\Sigma})$} \index{$Y_0({\frakn})$}
$$
\begin{aligned}
 \label{yKdef}  Y(K_{\Sigma})   & \ \ \   \left( \mbox{alternate notation:  $Y(\Sigma)$ or $ Y_0({\frakn})$, 
 where ${\frakn} = \prod_{{\p} \in \Sigma} {\p}$.}\right) \\
 :=   & \  \G(F) \backslash \G(\Adele)  \slash K_{\infty}   K_{\Sigma}   \\     
  = &  \   \G(F) \backslash (\mathbf{H}^3 \times \G(\Afinite)) / K_{\Sigma},
  \end{aligned}
  $$
where $K_{\Sigma} = \prod K_v$, 
where, fixing for every $v \notin S$ an identification of $D_v^{\times}$ with $\GL_2(F_v)$,
and so also of $\G(F_v)$ with $\PGL_2(F_v)$, we have set
 $$ \displaystyle{K_{v} =  \begin{cases} \PGL_2(\OL_v), & v \notin \Sigma;  \\ \
  g \in \PGL_2(\OL_{v}) \ | \
g \equiv 
\left(\begin{matrix} * & * \\ 0 & * \end{matrix} \right) \mod \pi_{v}, & v \in \Sigma - S; 
\\ \mbox{ the image of $B_v^{\times}$ in $\G(F_v)$} & \mbox{otherwise,}
  \end{cases}},
$$
where $B_v$ is a maximal order of $D_v$.

Then $Y(\Sigma)$ is an analog
of the modular curve $Y_0(N)$ with level $N$ the product of all primes in $\Sigma$.
 Note that there is an abuse of language in all of our notation $Y(\Sigma)$, $Y(K_{\Sigma})$, 
$Y_0({\frakn})$: 
we do not explicitly include in this notation the set of places $S$ defining the quaternion algebra $D$. {\em When we write $Y(\Sigma)$, there is always an implicit 
choice of a fixed  subset $S \subset \Sigma$.}

Almost all of our theorems carry over to general level structures. Accordingly, 
if $K \subset \G(\Afinite)$ is an open compact subgroup, we denote by $Y(K)$
the corresponding manifold $\G(F) \backslash  (\mathbf{H}^3 \times \G(\Afinite)) / K_{\Sigma}$;
of course, when $K = K_{\Sigma}$, this recovers the manifold above.

It will often be convenient in this book to use the phrase ``hyperbolic three manifold''
to include disconnected three-manifolds, each of whose components are hyperbolic. 
Similarly, in  order to avoid having to add
 a parenthetical remark ``(or possibly orbifold)'' after every use of the word manifold, 
 we shall often elide the distinction and  use the word manifold.  \index{orbifold prime}
 By ``orbifold prime'' for the orbifold $M$ we shall mean a prime that divides the order of one of the point  stabilizer groups. When we take cohomology of such an $M$, then,
 we always mean in the sense of orbifolds, that is to say: if $\tilde{M} \rightarrow M$ is a
 $G$-covering by a genuine manifold $\tilde{M}$, as will always exist in our settings,  the cohomology of $M$
 is understood to be the cohomology of $M \times E/G$, where $E$ is a space on which $G$ acts freely.

Our primary interest in this book is in the homology (or cohomology)
of the arithmetic manifolds $Y(K_{\Sigma})$ just defined. 
Since the expression $H_1(Y(K_{\Sigma}),\Z)$ is somewhat cumbersome, we shall,
hopefully without confusion, write:
$$H_1(\Sigma,\Z) := H_1(Y(K_{\Sigma}),\Z).$$
Since $\Sigma$ is only a finite set  of primes, it is hoped that no ambiguity will result, since
the left hand side has no other intrinsic meaning.
In the usual way (\S~\ref{ss:HeckeDef}) one defines the action of a Hecke operator 
$$T_{\p} : H_1(\Sigma, \Z) \longrightarrow H_1(\Sigma, \Z)$$
for every $\p \notin S$.  These commute for different $\p$.  When $\p \in \Sigma - S$,
this operator is what is often called the ``$U$-operator.'' 

If $\q \in \Sigma$ is a prime that does not lie in $S$, then we write
$\Sigma = \Sigma/\q \cup \{\q\}$, and hence:
$$H_1(\Sigma/\q,\Z) := H_1(Y(K_{\Sigma/\q}),\Z).$$

 Suppose that $\q \in \Sigma$ is not in $S$. There are two natural degeneracy maps
$\YO \rightarrow \YOq$ (see~\S~\ref{ss:HeckeDef});  these give rise to a pair of maps: 
$$\Phi:H^1(\Sigma /\q,\Q/\Z)^2 \rightarrow H^1(\Sigma,\Q/\Z),$$
$$\Phi^{\vee}:  H^1(\Sigma,\Q/\Z) \rightarrow H^1(\Sigma/\q,\Q/\Z)^2,$$
where the first map $\Phi$ is the pullback, and the latter map  is the transfer. In the same way, we obtain maps:
$$\Psi: H_1(\Sigma,\Z) \rightarrow H^1(\Sigma/\q,\Z)^2,$$
$$\Psi^{\vee}: H_1(\Sigma/\q,\Z)^2 \rightarrow H_1(\Sigma,\Z).$$
The map $\Psi$ can be obtained from $\Phi$ (and similarly
$\Psi^{\vee}$ from $\Phi^{\vee}$)  by applying the functor $\Hom(-,\Q/\Z)$.

\subsection{Additional notation in the split case} \label{ss:splitnotn}

When $\G = \PGL_2$, a situation analyzed in detail in Chapter 
\ref{chapter:ch5}, we will use some further specialized notation. In this case
$F$ is imaginary quadratic.  We put
\index{$\B$ (Borel subgroup)}
$$ \B = \left( \begin{array}{cc} * & * \\ 0 & * \end{array} \right) \subset \PGL_2,$$
a Borel subgroup, and denote by $\NN$ its unipotent radical \index{$\alpha: \B \rightarrow \mathbf{G}_m$.} 
$ \left( \begin{array}{cc} 1 & * \\ 0 & 1 \end{array} \right)$;
we denote by $\diagA$ the diagonal torus.  We denote by $\alpha:  \B \rightarrow \mathbf{G}_m$ the positive root:
\begin{equation} \label{alphadef} \alpha:  \left( \begin{array}{cc} x & * \\ 0 & 1 \end{array} \right)  \mapsto x. \end{equation}

We also define $\PU_2 \subset \PGL_2(\C)$ to be the stabilizer of the
standard Hermitian form $|x_1|^2 + |x_2|^2$.  Then
$ \PU_2 \times \prod_{v} \PGL_2(\OO_v)$, the product being taken over finite $v$,  is a maximal compact subgroup of $\PGL_2(\Adele)$.

Chapter~\ref{chapter:ch5} uses a significant amount of further notation, and we do not summarize it here since it is localized to that Chapter. 

  \newpage

\section{Fields and adeles} 

\subsection{The number field \texorpdfstring{$F$}{F}} \label{ss:notn:F} 
Let $F$ be a number field with 
one complex place\footnote{Most of the notation we give makes sense without the assumption that $F$ has one complex place.} and $r_1 = [F:\Q]-2$ real places.   We denote the ring of integers of $F$ by
  $\order$ and set $F_{\infty} := F \otimes \R \cong \C \times \R^{r_1}$. 
  We fix, once and for all, a complex embedding $F \hookrightarrow \C$ 
  (this is unique up to complex conjugation).

  \medskip
  
The adele ring of $F$ will be denoted $\Adele$, and the ring of {\em finite}
adeles $\Afinite$. If $T$ is any set of places, we denote by $\Adele^{(T)}$
the adele ring omitting places in $T$, i.e. the restricted direct product of $F_v$ for $v \notin T$. 
Similarly --- if $T$ consists only of finite places ---  we write $\Afinite^{(T)}$ for the finite adele ring omitting places in $T$. 
We denote by $\OO_{\adele}$  or $\widehat{\OO}$ the closure of $\OO_F$ in $\Afinite$; it is 
isomorphic to the profinite completion of the ring $\OO_F$.

For $v$ any place of $F$, and $x \in F_v$, we denote by $|x|_v$ the normalized absolute value on $F_v$, i.e.,
the effect of multiplication by $x$ on a Haar measure. In the case of $\C$ this differs from the usual absolute value:
$|z|_{\C} = |z|^2$ for $z \in \C$, where we will denote by $|z|$ the ``usual'' absolute value $|x+iy| = \sqrt{x^2+y^2}$. 
For $x \in \adele$ we write, as usual, $|x| = \prod_{v} |x_v|_v$. 

We denote by $\zeta_F$ the zeta-function of the field $F$,
and by $\zeta_v$ its local factor at the finite place $v$ of $F$.  We denote by $\xi_F$ the completed $\zeta_F$, i.e. including $\Gamma$-factors at archimedean places.

When we deal with more general $L$-functions (associated to modular forms), 
we follow   the convention that $L$ refers to the ``finite part'' of an $L$-function, i.e., excluding all archimedean factors, and $\Lambda$ refers to the completed $L$-function. 

We denote by $\Cl(\OO_F)$ (or simply $\Cl(F)$) the class group of $\OO_F$,
and we denote its order by $h_F$. \index{$h_F$}  \index{$\Cl$}

We denote by $w_F$ the the number of roots of unity in $F$,  \index{$w_F^{(2)}$}
i.e. (denoting by $\mu$ the set of all roots of unity in an algebraic closure, as a Galois module),
the number of Galois-invariant elements of $\mu$.

 We denote by $w_F^{(2)}$
the number of Galois-invariant elements of $\mu \otimes \mu$, equivalently,
the greatest common divisor of $\mathrm{N}({\q})^2-1$, where
${\q}$ ranges through all sufficiently large prime ideals of $F$. 
One easily sees that if $F$ is imaginary quadratic, then 
$w_F$ divides $12$ and 
$w_F^{(2)}$ divides $48$. 
In particular it is divisible only at most by $2$ and $3$.

We denote by $w_H$ \index{$w_H$} the number of roots of unity in the Hilbert class field of $F$. 
If $F$ has a real place, then $w_H = 2$; otherwise, $w_H$ is divisible at most by $2$ and $3$.

\subsection{Quaternion Algebras} 

Let $D$ be a quaternion algebra over $F$
ramified at a set $S$ of places that contains all real places. 
(By abuse of notation we occasionally use $S$
to denote the set of finite places at which $D$ ramifies;
this abuse should be clear from context.)

Recall that $D$ is determined up to isomorphism by the set
  $S$ of ramified places; moreover, the association of $D$ to $S$
  defines a bijection between isomorphism classes of quaternion algebras, 
  and subsets of places with even cardinality. 
Let $\psi$ be the standard anti-involution\footnote{
It is uniquely characterized
 by the property
that $\left( (t-x)(t-\psi(x) \right)^2 $ belongs to $F[t]$ and is the characteristic polynomial
of multiplication by $x \in D$ on $D$. }
 of $D$; it preserves \label{section:involutions}
the center $F$, and 
the reduced norm of $D$ is given by $x \mapsto
x \cdot \psi(x)$. The additive identity $0 \in D$ is the only element of norm zero.

Finally, let $\G$ be the algebraic group over $F$ defined as $$\G = \GL_1(D)/\mathbb{G}_m.$$ %

When $v \notin S$, the group $\mathbb{G} \times_F F_v$ is isomorphic
to $\PGL_2$ (over $F_v$). We fix once and for all an isomorphism 
\begin{equation} \label{isodef} \iso_v: \mathbb{G} \times_F F_v \stackrel{\sim}{\rightarrow} \PGL_2. \end{equation}
This will be used only as a convenience to identify subgroups. 

 Let $B/\OL_F$ be a maximal order in  $D$.
 Let $B^{\times}$ denote the invertible elements of $B$; they are precisely
 the elements of $B$ whose reduced norm lies in $\OL^{\times}_F$.
 Let $B^1 \subset B^{\times}$ denote the elements of $B$ of reduced norm one.

 Later we shall compare two quaternion algebras $D, D'$; we attach to $D'$
data $ S', \G', A'$.

\subsection{Quaternion Algebras over Local Fields}
\label{quatinfo}

Let $F_v/\Q_p$ denote the completion of $F$ with ring
of integers $\OO_{F,v}$, uniformizer $\pi_v$, and residue field $\OO_{F,v}/\pi_v \OO_{F,v}= k_v$.

If $D_v$ is not isomorphic to the algebra of $2 \times 2$ matrices, it 
  can be represented
by the symbol
$\left( \frac{\pi,\alpha}{F_v}\right)$, i.e. 
$$ F_v[i,j] / \{i^2 = \pi_v, j^2 = \alpha, ij=-ji \}.$$
where $\alpha \in \OO_{F,v}$ is such that the image of $\alpha$ in $k$ is a quadratic non-residue.
Then $B_v:=\OO_{F,v}[i,j,ij]$ is a maximal order and all such are conjugate. 

\medskip

The maximal order $B_v$ admits a unique maximal  bi-ideal $\m_v$, which is explicitly given
by $(i)$ with the presentation above. Then $\ell :=B_v/\m_v$
is a quadratic extension of $k$. 
  The image of $B^{\times}_v$ in $B_v/\m_v$ is $l^{\times}$, whereas
the image of  the norm one elements $B^1_v$ in $B_v/\m_v$ is given by those elements in $l^{\times}$ whose norm to $k$ is trivial.   Denote this subgroup
by $l^1$.

\section{The hyperbolic \texorpdfstring{$3$}{3}-manifolds} \label{sec:YKdef}

\subsection{Level Structure} \label{ls}

We now fix notations for the level structures that we study. 
They all correspond to   compact open subgroups of $G(\Afinite)$ 
of ``product type,''  that is to say, of the form $\prod_v K_{v}$.

\medskip

We have fixed an isomorphism of $\G(F_v)$
with $\PGL_2(F_v)$ whenever $v \notin S$ (see~\eqref{isodef}) We  use it in what follows to identify subgroups of $\PGL_2(F_v)$
with subgroups of $\G(F_v)$.

For any finite place, 
set $K_{0,v} \subset \G(\OL_v)$
to be the image in $\PGL_2$ of $$\displaystyle{K_{0,v} = \left\{ g \in \ GL_2(\OL_{v}) \ | \
g \equiv 
\left(\begin{matrix} * & * \\ 0 & * \end{matrix} \right) \mod \pi_{v} \right\}}$$
if $\G$ is split and the image of $B_v^{\times} \rightarrow \G(F_v)$ otherwise.   

 Define $K_{1,v}$ to be  the image in $\PGL_2$ of  $$\displaystyle{K_{1,v} = \left\{ g \in \GL_2(\OL_{v}) \ | \
g \equiv 
\left(\begin{matrix} 1 & * \\ 0 & * \end{matrix} \right) \mod \pi_{v} \right\}}$$
if $v$ is split and the image of $1+\mathfrak{m}_v$ inside $\G(F_v)$ otherwise. 

Thus, $K_{1,v}$ is normal in $K_{0,v}$, and the quotient is isomorphic
to $k_v^{\times}$
 or $l_v^{\times}/k_v^{\times}$, according to whether $\G$ is split or not; here $l_v$ is the quadratic extension of the residue field. 

Let $\Sigma$ denote a finite set of primes containing every element of $S$.
We then set
$$K_{\Sigma,0} := \prod_{v \in \Sigma} K_{0,v}  \prod_{v \notin~\Sigma} \G(\OL_v).$$
and similarly $K_{\Sigma, 1}$. 

\begin{quote}
{\em We write simply $K_{\Sigma}$ for $K_{\Sigma, 0}$. }
 \end{quote}
 
In particular, the quotient
$K_{\Sigma}/K_{\Sigma,1}$ is abelian, and it has size
$$\prod_{\q \in S} (N(\q) +1) \prod_{\q \in \Sigma -S} (N(\q) - 1).$$

In all instances, let  $K^{(1)}_v \subset K_v$ denote the subgroup of norm one elements,
that is, the elements either of determinant one if $\G(\OL_v) = \GL_2(\OL_v)$ or
norm one if $\G(\OL_v) = B^{\times}$.

\subsection{Groups and Lie algebras} 
Let $G_{\infty} = \G(\R)$ and
let $K_{\infty}$ be a maximal compact of $G_{\infty}$. 
We may suppose that $K_{\infty}$ is carried, under 
the morphism $G_{\infty} \rightarrow \PGL_2(\C) $
derived from~\eqref{isodef}, to the standard subgroup $\PU_2$
corresponding to the stabilizer of the Hermitian form $|z_1|^2+|z_2|^2$. 
  \index{$\mathfrak{g}$} \index{$\mathfrak{k}$}
 
 We denote by $\mathfrak{g}$ the Lie algebra of $G_{\infty}$
 and by $\mathfrak{k}$ the Lie algebra of $K_{\infty}$.

\subsection{\texorpdfstring{$\YO$}{Y} and \texorpdfstring{$\YN$}{Y}}  \label{YMdef} 

 If $K$ is a compact open subgroup of $\G(\Afinite)$, we define
the associated arithmetic quotient\footnote{Note that, for general groups $\mathbb{G}$, it might
 be more suitable, depending on the application, to replace $K_{\infty}$ in this definition by $K_{\infty}^0 A_{\infty}^0$ (so that $Y(K)$ has a finite invariant measure
and has an orientation preserved by $G_{\infty}$. } 
$$
\begin{aligned}
 \label{yKdeftwo} Y(K)   :=   & \  \G(F) \backslash \G(\Adele)  \slash K_{\infty}   K   \\     
  = &  \   \G(F) \backslash (G_{\infty}/K_{\infty} \times \G(\Afinite)) / K,
  \end{aligned}
  $$
 The manifold $Y(K)$ is a (possibly disconnected) hyperbolic three manifold
of finite volume, since $G_{\infty}/K_{\infty}$ is isometric to $\mathbf{H}$ (see~\S~\ref{H3} for that identification).

 As before, let $\Sigma$ be a finite set of finite places containing all places in $S$. 
 Define the arithmetic quotients $\YO$ and $\YN$ as:
$$\YO:=\G(F) \backslash \G(\Adele)/K _{\infty}  \KMO,$$
$$\YN:=\G(F) \backslash \G(\Adele)/K _{\infty}   \KN.$$
 
 By analogy with usual notation, we sometimes write $Y_0({\frakn})$
 and $Y_1({\frakn})$ for $\YO$ or $\YN$, where ${\frakn} = \prod_{\p \in \Sigma} \p$.

 \begin{remarkable} \label{conncompremark} {\em 
We may express this more explicitly:
\begin{equation} \label{A def} \YO= \coprod_{\AN}
Y_0(\Sigma,\a) = 
\coprod_{\AN} \Gamma_0(\Sigma,\a) \backslash \H^3, \end{equation} 
$$\YN= \coprod_{\AN}
Y_1(\Sigma,\a) = 
\coprod_{\AN} \GammaP_1(\Sigma,\a) \backslash \H^3,$$
where:
\begin{enumerate}
\item The indexing set $\AN$ is equal to  (the finite set) $F^{\times} \backslash \Atimes \slash
\kern+0.1em{\OO^{\times}_{\Adele}} {\A^{\times 2}}$, which is the quotient of the class group of $F$ by squares. 
It is important to notice that $\AN$ 
depends only on $F$, and not on $D$ or $\Sigma$.

More generally, the  component set for $Y(K)$ (for general $K$) may be identified with a quotient 
of some ray class group
of $\OL_F$.

\item Each $\Gamma_i(\Sigma,\a)$ is a subgroup of $\G(F)$,
realized as a finite volume discrete arithmetic subgroup of
$\PGL_2(\C)$ via the unique complex
embedding $F \rightarrow \C$ and the fixed isomorphism of $\G(\C)$ with $\PGL_2(\C)$. 
Moreover, the corresponding group will be co-compact providing that $\G \neq \PGL_2$.
Explicitly,
let $\alpha \in \Atimes$ represent the class $\a\in \AN$; then 
\begin{equation} \label{gammasigmaadef} \Gamma_i(\Sigma,\a) \simeq 
\left\{ \gamma \in \G(F) \subset \G(\A) \ \left| \ 
\left( \begin{matrix} \alpha & 0 \\ 0 & 1 \end{matrix}
\right)^{-1} \gamma\left( \begin{matrix} \alpha & 0 \\ 0 & 1 \end{matrix}
\right) \in K_{\Sigma,i} \right\} \right. \end{equation} 
\end{enumerate}}
\end{remarkable}

\begin{example}
\emph{
 We specialize our notation to various classical situations in order to orient the reader:
 }
 
 \emph{
 Let $\G = \GL(2)/\Q$, $S  =\varnothing$, $\Sigma = \{p,q\}$, and $\Gamma = \PSL_2(\Z)$.
 Suppose instead that  $\G' = \GL_1(D)/\Q$, where $D/\Q$ is ramified at $p$
and $q$ (and so indefinite at $\infty$),  $S = \{p,q\}$, $\Sigma = \{p,q\}$.
Denote the corresponding manifolds by $Y$ and $Y'$. }

\emph{ 
Then
$\YO$ can be identified with $Y_0(pq)$, and 
$\YN$ can be identified with $Y_1(pq)$, where 
$Y_0(n)$ and $Y_1(n)$ are the modular curves of level $n$ parameterizing
(respectively) elliptic curves with a  cyclic subgroup of order $n$, 
and elliptic curves with a point of exact order $n$.  
} 

\emph{
On the other hand,  $ Y'(\Sigma)$ can be identified with a Shimura variety, namely, 
the variety parameterizing ``fake elliptic curves'' (principally polarized abelian surfaces $A$
together with an embedding of $\OO_D$ into $\End(A)$  satisfying some extra conditions). 
On the other hand, $Y'(K_{\Sigma,1})$ can be identified with the same moduli 
space together with a level structure: if ${\p}, {\q}$
are the maximal bi-ideals of $\OO_D$ corresponding to $p, q$,
then the level structure is given by specifying a generator of $K$,
where $A[{\p} {\q}] \simeq K \oplus K$ under the idempotents
$\displaystyle{e = \left(\begin{matrix} 1 & 0 \\ 0 & 0 \end{matrix}\right)}$ and $1-e$
of $\OO_D/{\p}{\q}$, having identified
the latter with
a matrix algebra,
cf.~\cite{BuzzardThesis}). 
}
\end{example}

\subsection{Orbifold primes for \texorpdfstring{$Y(K_{\Sigma})$}{Y}} 

The $Y(K_{\Sigma})$ are orbifolds. As such, to each point there is canonically associated
an isotropy group. 
It is not difficult to see that the only primes dividing the order of an orbifold isotropy group
are at most the primes dividing $w_F^{(2)}$ (see~\S~\ref{ss:notn:F}).

In particular, the only possible orbifold primes in the split case, or indeed if $F$ is quadratic imaginary,  are $2$ and $3$.
At certain points in the text we will invert $w_F^{(2)}$; this is simply to avoid such complications. 
 
 We will sometimes say ``$A=B$ up to orbifold primes.'' By this we mean that the ratio $A/B$
 is a rational number whose numerator and denominator are divisible only by orbifold primes.

 \section{Homology, cohomology, and spaces of modular foms} \label{sec:notn2} 
 
 Recall (\S~\ref{sec:notn0}) that we denote  $H_1(Y(K_{\Sigma}),\Z)$ by $H_1(\Sigma,\Z)$.
 These (co)homology groups may also be interpreted via group (co)homology.
Suppose that $\Gamma \subset \PGL_2(\C)$ is a discrete subgroup; then there
 are functorial isomorphisms
$$H^i(\Gamma \backslash \H) \simeq H^i(\Gamma, A),$$

There is a similar identification
on group homology.
 From the fact that $\YO$ is a disconnected sum of finitely many  $3$-manifolds
 which are of the form $\Gamma_0(M,\a) \backslash \H^3$,
we deduce %
that each $\Gamma_0(M,\a)$ has cohomological
dimension at most $3$ (away from orbifold primes, i.e. on modules on which $w_F^{(2)}$ is invertible).

 \subsection{Duality and the linking form}
 \label{section:linking}
 
 We discuss the compact case; for the noncompact case
 see~\S~\ref{section:linking2}.
 
 For a compact $3$-manifold $M$, $H_0(M,\Z) = H_3(M,\Z) = \Z$, and
 Poincar\'{e} duality defines isomorphisms $H_1(M,\Z) \simeq H^2(M,\Z)$ and
 $H_2(M,\Z) \simeq H^1(M.\Z)$. From the universal coefficient theorem, we further
 deduce that 
 $$H^1(M,\Q/\Z) \simeq \Hom(H_1(M,\Z),\Q/\Z)).$$
 For an abelian group $A$, let $A^{\vee} = \Hom(A,\Q/\Z)$. \index{$A^{\vee}$} 
 If $A$ is finite, then $A^{\vee \vee} = A$.
 The exact sequence
 $0 \rightarrow \Z \rightarrow \Q \rightarrow \Q/\Z \rightarrow 0$
 induces a long exact sequence in homology; esp.
\begin{equation} H^1(\Q) \stackrel{\alpha}{\rightarrow}
 H^1(\Q/\Z) \stackrel{\delta}{\rightarrow} H_1(\Z)^{\vee} \stackrel{\beta}{\rightarrow} H^2(\Q) \end{equation} 
 $\delta$ induces an isomorphism of $\mathrm{coker}(\alpha)$ 
 with $\ker(\beta)$, i.e. an isomorphism of the ``torsion parts'' of $H^1(\Q/\Z)$
 and $H^1(\Q/\Z)^{\vee}$; the pairing $(x,y) \mapsto (x, \delta y)$ or $(x, \delta^{-1} y)$
 therefore induces a {\em perfect} pairing on these torsion groups.  
 
 \medskip
 
 This is the
`` linking form''; it can also be considered as a perfect pairing 
 $$H_1(M, \Z)_{\tors} \times H_1(M, \Z)_{\tors} \rightarrow \Q/\Z.$$
 obtained from intersection-product together with the connecting homomorphism
 $H_2(\Q/\Z)_{\tors}  \stackrel{\sim}{\rightarrow} H_1(\Z)_{\tors}$. 
 
 \subsubsection{Duality for orbifolds}  \label{sss:dualityorbifold} 
 
 For compact {\em orbifolds}, these conclusions still hold
 so long as we localize away from the order of any torsion primes: If $N$
 is such that any isotropy group has order dividing $N$,
 then Poincar{\'e} duality holds with $\Z[\frac{1}{N}]$ coefficients, 
and the foregoing goes through replacing $\Z$  by $\Z[\frac{1}{N}]$ 
and $\Q/\Z \cong \bigoplus_p \Q_p/\Z_p$ by $\bigoplus_{(p, N)=1} \Q_p/\Z_p$.

  Over $\Z$ the following still holds: The pairing
 $$H_1(Y(K_{\Sigma}), \Z)_{\tf} \times H_2(Y(K_{\Sigma}), \Z)_{\tf} \rightarrow \Z$$
(obtained from the corresponding pairing with $\Q$ coefficients, rescaled so the image lies in $\Z$) is perfect away from orbifold primes, i.e. the discriminant of the corresponding Gram matrix
is divisible only by orbifold primes.

 \subsection{Atkin-Lehner involutions.} \label{subsec:ali}
 
 The manifold $Y(K_{\Sigma})$ has a canonical action of the group $(\Z/2\Z)^{\Sigma}$, generated by the so-called Atkin-Lehner involutions. Indeed, for every $v \in \Sigma$
 for which $\G$ is split 
 the element 
 
\begin{equation} \label{wvdef} w_v := \left( \begin{array}{cc} 0 & 1 \\ \varpi_v & 0 \end{array}\right)  \in \G(F_v) \end{equation}
 normalizes $K_{0,v}$;  for $v \in \Sigma$ for which $\G$ is not split
 we take $w_v$ to be any element of $\G(F_v)$ not in $K_{0,v}$.

The group generated by the $w_v$ generates
 a subgroup of $N(K_{\Sigma})/K_{\Sigma}$ isomorphic to $(\Z/2\Z)^{\Sigma}$.

 The element $w_v$ also normalizes $K_{1, v}$ and it acts by inversion
 on the abelian group $K_v/K_{1,v}$ (which is isomorphic to $k_v^{\times}$
 or $l_v^{\times}/k_v^{\times}$, according to whether $\G$ is split or not). 
 
If $M$ is a space of modular forms with an action of an (understood) Atkin-Lehner involution $w$, we write $M^{\pm}$ for the $+$ and $-$ eigenspaces of $w$ on $M$. An example of particular importance is the cohomology of $Y_0(\q)$ for $\q$ prime; thus, for example, we write
$$H_1(Y_0(\q), \Q)^{-} := \{ z \in H_1(Y_0(\q), \Q): w_{\q} z = -z\}.$$
 
 (We do not use this notation when it is ambiguous which Atkin-Lehner involution is being referred to.) 

\subsection{Hecke Operators} \label{ss:HeckeDef}
 
For any $g \in \G(\Afinite)$ we may consider the diagram 
\begin{equation} \label{Heckediagram}  \YO  \leftarrow Y(g K_{\Sigma}  g^{-1} \cap K_{\Sigma}) \rightarrow  Y(K_{\Sigma}) \end{equation} 
that arise from, respectively, the identity on $\G(\Afinite)$ and right multiplication by $g$
on $\G(\Afinite)$.

\medskip

Now fix any prime ${\p} \notin S$ and let $g_{\p} \in \G(F_{\p})$ be the preimage of 
$ \left(\begin{matrix} \alpha & 0 \\ 0 & 1 \end{matrix} \right)$
under the fixed isomorphism of~\eqref{isodef}, and let $g \in \G(\Afinite)$ be the element
with component $g_{\p}$ at the prime $\p$ and all other components trivial. 
In this case~\eqref{Heckediagram} produces two maps $Y(\Sigma \cup \{\q\}) \rightarrow Y(\Sigma)$, 
often referred to as the two degeneracy maps. 
The corresponding
operator on homology or cohomology of $Y(\Sigma)$, obtained by pullback followed by pushforward, is the $\p$th Hecke operator, denoted $T_{\p}$.

 \medskip
 
These Hecke operators  preserve $H^{\bullet}(\YO,\Z)$, but
not the cohomology of the connected components. Indeed, the action on the component
group is via the determinant map on $\G(\Afinite)$ and the natural action of
$\Afinite^{\times}$ on $\AN =  F^{\times} \backslash \Afinite^{\times}/U.$

We now define the Hecke algebra (cf. definition~\ref{df:hecke}).
\begin{df}
Let $\T_{\Sigma}$  \index{$\T$} denote the subring of $\End \ H_{q}(Y(K_{\Sigma}),\Z)$
generated by Hecke endomorphisms 
$T_{\q}$ for all primes $\q$ not dividing $\Sigma$. 
\end{df}

\begin{df}  If $f$ is an eigenform for the Hecke algebra, 
we denote by $a(f, \q)$ the eigenvalue of $T_{\q}$ on $f$. 
\end{df}

 \subsection{The abstract Hecke algebra} \label{CompletionConvention}
 
 It is sometimes convenient to use instead the {\em abstract} Hecke algebra.   
 
By this we mean: For any finite set of places
$\Omega$, let $\mathscr{T}_{\Omega} := \Z[\mathcal{T}_{\q}: \q \notin \Omega]$ 
be the free commutative algebra on {\em formal} symbols $\mathcal{T}_{\q}$. 
Then there is an obvious surjection $\mathscr{T}_{\Omega}\twoheadrightarrow \T_{\Omega}$.
We will very occasionally denote $\mathscr{T}_{\Omega}$ simply {\em by} $\T_{\Omega}$;
when this is so we clarify in advance.

A typical situation where this is of use is when we are working with   two levels $\Sigma, \Sigma'$ simultaneously:

\medskip Suppose given a maximal ideal $\mathfrak{m}$ of $\T_{\Sigma}$.
 We will  on occasion refer to ``the completion of the cohomology $H_q(Y(\Sigma'))$ at $\mathfrak{m}$.''
By convention this means -- unless otherwise specified -- the following:  
Any ideal $\mathfrak{m}$ of $\T_{\Sigma}$ gives rise to a maximal ideal $\mathfrak{m}'$ of
$\mathscr{T}_{\Sigma \cup \Sigma'}$ via preimage under the maps $$ \mathscr{T}_{\Sigma \cup \Sigma'} \rightarrow  \T_{\Sigma},$$
and by $H_q(\Sigma')_{\mathfrak{m}}$ we mean, {\em by convention},
the completion of $H_q(\Sigma')$ at the maximal ideal
$\mathfrak{m}'$ of $\mathscr{T}_{\Sigma \cup \Sigma'}$. 

\subsubsection{Localization versus Completion}

  The Hecke algebras $\T \otimes \Z_p$ are always
semi-local rings, that is, the completions $\T_{\m}$ are finite $\Z_p$-algebras.
It follows that the \emph{completion} $H_q(\Z_p)_{\m}$ of $H_q(\Z_p)$ is equal
to the \emph{localization} $H_q(\Z_p)_{(\m)}$.

Similarly, whenever $H_q(\Z)$ is a finite
group, completion may be identified with localization. 

When we are dealing with $\Q_p/\Z_p$-coefficients, 
the localization of $H_q(\Q_p/\Z_p)$ at $\m$ is isomorphic
not to the completion but rather to the $\m^{\infty}$-torsion:
$$ \varinjlim H_q(\Q_p/\Z_p)[\mathfrak{m}^n]  = H_q(\Q_p/\Z_p)[\m^{\infty}] \stackrel{\sim}{\longrightarrow} H_q(\Q_p/\Z_p)_{(\m)}.$$
By a slight abuse of terminology, we will denote these spaces
also by $H_q(\Q_p/\Z_p)_{\m}$.

In general, the
completion $\T_{\m}$ is different from $\T_{(\m)}$ if $\T$ is infinite. We shall only
ever be concerned with modules over the first algebra.

  \subsection{Degeneracy maps}   \label{section:levelmaps}
Suppose that $\q \in \Sigma$ is not in $S$. The two natural degeneracy maps
$\YO \rightarrow \YOq$ noted in the previous section give rise to a pair of maps: 
$$\Phi:H^1(\Sigma /\q,\Q/\Z)^2 \rightarrow H^1(\Sigma,\Q/\Z),$$
$$\Phi^{\vee}:  H^1(\Sigma,\Q/\Z) \rightarrow H^1(\Sigma/\q,\Q/\Z)^2,$$
where the first map is induced from the degeneracy maps whilst the second
is the transfer homomorphism induced from the first. In the same way, we obtain maps:
$$\Psi: H_1(\Sigma,\Z) \rightarrow H_1(\Sigma/\q,\Z)^2,$$
$$\Psi^{\vee}: H_1(\Sigma/\q,\Z)^2 \rightarrow H_1(\Sigma,\Z).$$
The map $\Psi$ can be obtained from $\Phi$ (and similarly
$\Psi^{\vee}$ from $\Phi^{\vee}$)  by applying the functor $\Hom(-,\Q/\Z)$.

By taking the difference (resp. sum) of the maps $\Phi^{\vee}$,  we obtain
$$H_1(\Sigma, \Z)^- \rightarrow H_1(\Sigma/\q, \Z),   \ \ H_1(\Sigma, \Z)^+ \rightarrow H_1(\Sigma/\q, \Z)$$
where $+$ and $-$ refer to eigenspaces for the Atkin-Lehner involution $w_{\q}$ corresponding to the prime $\q$. 

When we write ``degeneracy map'' from either $H_1(\Sigma, \Z)^-$
or the companion group $H_1(\Sigma, \Z)^+$ to $H_1(\Sigma/\q, \Z)$, we always have in mind
these differenced or summed versions of $\Phi^{\vee}$. Similar remarks apply 
for $\Phi, \Psi, \Psi^{\vee}$.

\begin{lemma} 
The composite maps \label{lemma:composite}
$\Phi^{\vee} \circ \Phi$ and $\Psi  \circ \Psi^{\vee}$ are equal to
$$\left( \begin{matrix}   (N(\q) + 1) & T_{\q} \\
T_{\q} &  (N(\q) + 1) \end{matrix} \right),$$
where  $T_{\q}$ is the
Hecke operator. \end{lemma}

\begin{remark}
Were we dealing with $\GL_2$ rather than $\PGL_2$ (or the full group of units
of a quaternion algebra, and not its quotient by center) this would be replaced by
$$\left( \begin{matrix} \langle \q \rangle (N(\q) + 1) & T_{\q} \\
T_{\q} &  (N(\q) + 1) \end{matrix} \right),$$
where $\langle \q \rangle$ and $T_{\q}$ are the diamond and
Hecke operators respectively. In many of our results
where the expression $T_q^2 - (N(\q)+1)^2$ appears, 
it can be generalized to the case of $\GL_2$
simply by replacing it by $T_q^2 - \langle \q \rangle (N(\q)+1)^2$. 
\end{remark}

 \section{Normalization of metric and measures}  \label{metricmeasurenotn}

 \subsection{The metric on hyperbolic \texorpdfstring{$3$}{3}-space}  \label{H3}
 
 By $\H^3$ we denote hyperbolic $3$-space , i.e. triples $(x_1, x_2, y) \in \R^2 \times \R_{>0}$. We shall sometimes identify it with
 $\C \times \R_{>0}$ via $(x_1, x_2, y) \mapsto (x_1 + i x_2, y)$. 
We equip it with the metric of constant curvature  $$  g:=\frac{dx_1^2 +dx_2^2 +dy^2}{y^2}.$$

 There is a metric-preserving action of $\PGL_2(\C)$ on  $\H^3$, in which 
 the subgroup $\PU_2$ preserves the point $(0,0,1)$,
the action of the element $\left(\begin{array}{cc} 1 & z \\ 0 & 1 \end{array}\right)$ is given by
 translation in the $x_1+ix_2$ variable, and the action of  
 $\left(\begin{array}{cc} a & 0 \\ 0 & 1 \end{array}\right)$
 is given by $(z, y) \in \C \times \R_{>0} \mapsto (az, |a|y).$ 
 
 \medskip
This normalizes a hyperbolic metric on every manifold $Y(K)$.

\medskip
 
 On the other hand, 
the action of $\PGL_2(\C)$ on $\H^3$ gives rise to another canonical metric  $g_{K}$, 
namely, that induced by the Killing form on $\PGL_2(\C)$ considered as a real Lie group. 
(For each point $x \in \H^3$,
the tangent space $T_x \H^3$ is identified with a certain quotient of the Lie algebra,
on which the Killing form induces a metric.) 
Then\footnote{If $h = \left(  \begin{array}{cc} 1 & 0 \\ 0 & -1 \end{array} \right)  \in \mathfrak{sl}_{2,\C}$, 
then $\langle h, h \rangle = 16$.  But the image of $h$ in the tangent space
at $(0,0,1)$ has hyperbolic length $2$. 
}
$g_K = 2 g,$ so that the Laplacian for $g$ corresponds to four times  the Casimir. 

  Let $\mathfrak{pgl}_2 \supset \mathfrak{pu}_2$
 denote the Lie algebra of $\PGL_2(\C) \supset \mathrm{PU}_2$;
 the quotient ${\p} := \mathfrak{pgl}_2/\mathfrak{pu_2}$
 it is identified naturally with the tangent space to ${\H}^3$
 at the $\mathrm{PU}_2$-fixed point.  The tangent bundle
 to ${\H}^3$ is then naturally identified with the quotient
 $\PGL_2(\C) \times {\p} / \mathrm{PU}_2$,
 where the $\mathrm{PU}_2$ acts via the adjoint action on~${\p}$. 
 This identification is $\PGL_2(\C)$-equivariant, where $\PGL_2(\C)$
 acts on $\PGL_2(\C)$ by left-translation and trivially on~${\p}$.

\subsection{Measures} \label{subsec:measures}

There is a volume form on $Y(K)$ specified by the Riemannian metric. When we speak
of the volume of a component of $Y(K)$, as in~\S~\ref{notnvol}, we shall always 
mean with reference to that measure, unless otherwise stated.

We fix for later reference a measure on $\adele$: 
For each finite place $v$, equip $F_v$ with the additive Haar measure  $dx$ that assigns mass $1$ 
to the maximal compact subring; for each archimedean place $v$, 
we assign $F_v$ the usual Lebesgue measure (for $v$ complex, we mean
the measure that assigns the set $\{z \in F_v: |z| \leq 1\}$ the mass $\pi$). 

We equip $F_v^{\times}$ with the measure $\frac{dx}{|x|} \cdot \zeta_{F_v}(1)$ (this assigns
measure $1$ to the maximal compact subgroup). 

In~\S~\ref{sec:data} we will depart slightly from these normalizations for compatibility 
with the work of Waldspurger, to which~\S~\ref{sec:data} is really a corollary. Similarly
we will depart slightly in~\S~\ref{Eismod}.

\subsection{Volumes} \label{notnvol}

We shall use the following convention:
\begin{quote} If $X$ is --- possibly disconnected --- hyperbolic $3$-manifold, then $\vol(X)$ denotes the {\em product} of the volumes of the connected components of $X$.
\end{quote}

\begin{theorem}[Volumes; Borel~\cite{Borel}] Let $S$ be the set of primes that ramify in $D$.
Let $\Sigma$ be a set of primes containing all elements of $S$.  Then \label{theorem:borel}
$$\vol(Y_0(\Sigma,\a))
= \frac{ \mu 
 \cdot \zeta_{F}(2) \cdot |d_F|^{3/2}}{2^{m}  (4 \pi^2)^{[F:\Q] - 1}} \cdot   \prod_{\q \in S} (N \q - 1) 
 \prod_{\q \in \Sigma \setminus S}
(N \q + 1),$$
for some  $m \in \mathbf{N}$ and some integer $\mu$ depending only on $F$. 
\end{theorem}

 According to our convention, then, the volume of $Y_0(\Sigma)$
 is given by the right hand side of the above formula, {\em raised to the power}
 of the number of components,
which is equal to the order of $\Cl(\OO_F)/ 2 \Cl(\OO_F)$ if
$\Cl(\OO_F)$ is the class group of $F$.
(See~\eqref{gammasigmaadef}).

\section{\texorpdfstring{$S$}{S}-arithmetic groups} \label{Sarithmetic}
We shall need at a few points the $S$-arithmetic analogue  of the spaces
$Y(K)$. Roughly, these 
bear the same relation to $\G(\OO[{\q}^{-1}])$ as 
the space $Y(K)$ bears to $\G(\OO)$.

Let $\q$ be a prime that at which $\G$ is split. 
Recall that $\G(F_{\q}) \cong \PGL_2(F_{\q})$ acts on an infinite $(\Norm(\q)+1)$-valent tree $\mathcal{T}_q$, 
the Bruhat--Tits building.   For the purpose of this book a ``tree'' is a topological space, not a combinatorial object: it is a CW-complex 
obtained by gluing intervals representing edges to vertices in the usual fashion.
This can be understood as an analogue  of the action of $G_{\infty}$
on $ G_{\infty}/K_{\infty}$. In particular $\mathcal{T}_q$ is contractible and the vertex and edge stabilizers are maximal compact subgroups.

We shall also deal with a ``doubled'' tree that places the role, relative to $\PGL_2(F_{\q})$, of the ``upper and lower half-planes'' for $\PGL_2(\R)$. Namely, let $\Doubletree_{\q}$
be the union of two copies of $\mathcal{T}_{\q}$, and the $\PGL_2(F_{\q})$-action is given by
the product of the usual action on $\mathcal{T}_{\q}$ and the determinant-valuation
action of $\PGL_2(F_{\q})$ on $\Z/2\Z$ (in which $g$ acts nontrivially if and only if
$\det g$ is of odd valuation).

We set, for an arbitrary open compact subgroup $K \subset \G(\Afinite)$, 
$$   Y ( K[\frac{1}{\q}]) =  \G(F) \backslash \left(\H^3 \times \Tree_{\q} \times  \G(\Adele^{(\infty, \q)})/K^{(\q)}\right). $$
$$   Y ( K[\frac{1}{\q}])^{\pm} =  \G(F) \backslash \left(\H^3 \times \Doubletree_{\q}  \times  \G(\Adele^{(\infty, \q)})/K^{(\q)}\right)$$
where $K^{(\q)}$ is the projection of $K$ to $\G(\Adele^{(\infty, \q)})$. 
\index{ $Y ( K[\frac{1}{\q}])$} \index{$Y ( K[\frac{1}{\q}])^{\pm}$}

   We shall sometimes also denote this space $Y(K[\frac{1}{{\q}}])$   simply by $Y(\frac{1}{{\q}})$ when the choice of $K$  is implicit. We will often denote it
   by $Y(\Sigma[\frac{1}{{\q}}])$ when $K = K_{\Sigma}$. 
    Evidently, $Y(\frac{1}{\q}) $ is the quotient of $Y^{\pm}(\frac{1}{\q})$
by a natural involution (that switches the two components of $\Doubletree_{\q}$).

    We abbreviate their cohomology in a similar fashion to that previously discussed, 
    including superscripts $\pm$ when necessary to distinguish between the above spaces, e.g.:
\begin{eqnarray} \label{coabbrev} 
   H^*(\Sigmainvertq, \Z)^{\pm} &:=& H^*( Y(K_{\Sigma}[\frac{1}{\q}])^{\pm},  \Z ).\end{eqnarray}
   \index{$ H^*(\Sigmainvertq, \Z)^{\pm}$}
Inverting the prime $2$, the  fixed space of $H^*(\Sigmainvertq, -)^{\pm}$ under the
 action induced by the mentioned involution is $H^*(\Sigmainvertq, -)$.

The two fixed points for the maximal compact $K_{\q}$ on the tree $\Doubletree_{\q}$
gives rise to two inclusions $Y(\Sigma/\q)$ into $Y(\Sigma[1/\q])^{\pm}$;
these two inclusions collapse to the same map in the quotient $Y(\Sigma[1/\q])$. 
 In fact, topologically, 
\begin{equation} \label{glueq} \mbox{   $Y (\Sigma[\frac{1}{\q}])^{\pm}$ is $ Y(\Sigma) \times [0,1]$ glued to two copies of $Y(\Sigma/\q)$;} \end{equation} 
this structure being obtained by dividing $\mathcal{T}_{\q}$ into $1$-cells and $0$-cells; 
the gluing maps are the two degeneracy maps $Y(\Sigma)  \rightarrow Y(\Sigma/\q)$;
one obtains a corresponding description of $Y(\Sigma[\frac{1}{\q}])$ after quotienting by the involutions.

\begin{example}   \label{pgl2exampleSarithmetic} Suppose that $\G  = \PGL(2)/F$, that $\Cl(\OL_F)$ is odd, and that $K = \PGL_2(\OL_{\A})$.
Then the spaces $Y(K)$, $Y^+(K[\frac{1}{\q}])$, and $Y (K[\frac{1}{\q}])$
are all connected $K(\pi, 1)$-spaces, and
\begin{eqnarray*} \pi_1( Y(K)) &\cong& \PGL_2(\OO); \\
\pi_1( Y(K[\frac{1}{\q}]) ) &\cong& \PGL_2(\OO[\frac{1}{\q}]); \\
\pi_1\left( Y (K[\frac{1}{\q}]^{\pm} \right) & \cong & \PGL_2(\OO[\frac{1}{\q}])^{(\ev)}, 
\end{eqnarray*}
where $\PGL_2^{(\ev)} \subset \PGL_2(\OL[\frac{1}{\q}])$ consists of those elements whose determinant
has even valuation at $\q$. 
In particular, the homology of each of these spaces is identified
with the group homology of the right-hand groups. 

\end{example}

In general we may write 
$$  \displaystyle{
 Y(K_{\Sigma}[\frac{1}{{\q}}])^{\pm} =
\coprod_{\AN/\q} \Gamma^{(\q)}_0(\Sigma,\a) \backslash  \left(  \mathbf{H}^3 \times \Doubletree_{\q} \right)}
$$
with the notation of~\eqref{A def}; 
here the group $\Gamma^{(\q)}_0(\Sigma, \a)$ is defined similarly
to $\Gamma_0(\Sigma, \a)$, but replacing $K_{\Sigma, 0}$
by $K_{\Sigma,0} \cdot \G(F_{\q})$;   and $\AN_{\q}$
is now the quotient of $\AN$ by the class of $\q$.

Denote by $ \Gamma_0^{(\q)}(\Sigma, \a)^{(\ev)}$
those elements of $ \Gamma_0^{(\q)}(\Sigma, \a)$ whose determinant has even valuation at $\q$. 
Note that $\Gamma_0^{(\q)}  \neq  \Gamma_0^{(\q)}(\Sigma, \a)^{(\ev)}$
if and only if  
$\Gamma_0^{(\q)}(\Sigma, \a)$ contains an element that switches
the two trees in $\Doubletree_{\q}$; this is so (for any $\a$)
if and only if $\q$ is a square in the ideal class group. \footnote{We are asking for the existence 
of an element of $\G(F)$
that belongs to the set
 $
\left( \begin{matrix} \alpha & 0 \\ 0 & 1 \end{matrix}
\right) K_{\Sigma, v} \left( \begin{matrix} \alpha & 0 \\ 0 & 1 \end{matrix}
\right) ^{-1}  $ for every finite $v$ besides $\q$, and whose
determinant has odd valuation at $\q$ itself.  One now applies the strong approximation theorem.}
 
 \medskip

Therefore $Y(\Sigma[\frac{1}{\q}])^{\pm}$ is homeomorphic to
\begin{equation} \label{y1qid}  \coprod_{A/\q} \Gamma_0^{(\q)}(\Sigma, \a)^{(\ev)} \backslash  (\H^3 \times \Tree_{\q}),   \mbox{ $\q$ a square}, \end{equation} 
and {\em two copies of the same}, if $\q$ is not a square. 
Note that in all cases it has exactly the same number of connecting components
as $Y(\Sigma)$.

\section{Congruence homology}  \label{s:congess}

 One class of elements in cohomology that
should be considered ``trivial'' are those classes 
arising from congruence covers (cf. the discussion
of~\S~\ref{section:discusscongruence}). 

\subsubsection{Congruence homology for arithmetic groups} 
\index{Congruence Homology}
Before the definition for $Y(K)$, let us give the corresponding definition for arithmetic groups: 

If $\Gamma$ is any arithmetic group, it admits a map 
to its congruence completion  $\widehat{\Gamma}$:
the completion of $\Gamma$ for the topology defined by congruence subgroups.
This yields a natural surjection
$$H_1(\Gamma, \Z) \longrightarrow H_1(\widehat{\Gamma}, \Z)$$
which we call the ``congruence quotient'' of homology, or, by a slight abuse of notation, 
the ``congruence homology'' $H_{1, \con}(\Gamma, \Z)$.

For example,   if $\Gamma = \Gamma_0(N) \subset \SL_2(\Z)$, 
then the morphism $$H_1(\Gamma_0(N), \Z) \longrightarrow (\Z/N\Z)^{\times}$$
 arising from $\left( \begin{array} {cc} a & b \\ c & d \end{array} \right) \mapsto a$,
 factors through the congruence homology.

\subsubsection{Congruence homology for \texorpdfstring{$Y(K)$}{Y(K)}}  
\label{subsubsection:reducetogroup}
Recall that $Y(K)= \coprod \Gamma_i \backslash \H^3$
for various arithmetic subgroups $\Gamma_i$ (cf.~\S~\ref{conncompremark}). 
For any coefficient group $A$,  understood to have trivial $\Gamma_i$-action, we define $H_{1, \con}(Y(K),  A)$ as the quotient of $H_1(Y(K),A)$ defined by 
the map: 
\begin{equation} \label{drm} 
 \bigoplus H_1(\Gamma_i,  A)
\twoheadrightarrow \bigoplus H_1(\widehat{\Gamma}_i,  A)  ,  \end{equation} 
Note that, at least for $H_1$, this map is indeed surjective.

\subsubsection{Reformulation}

There is a convenient adelic formalism to work with the right-hand side of~\eqref{drm}.
Although it has no essential content it makes the statements more compact: 

Let us denote by $Y(K)^{\carrot}$ the quotient
$$Y(K)^{\carrot} :=  \overline{\G(F)} \backslash \G(\Afinite) / K.$$
Here, $\overline{\G(F)}$ is the  the closure of $\G(F)$ in $\G(\Afinite)$;  by the strong approximation theorem 
 the closure of $\G(F)$ contains the $\Afinite$-points of the derived group $[\G, \G]$, and
 therefore $\overline{\G(F)}$ is a {\em normal} subgroup of $\G(\Afinite)$, with abelian quotient;
{\em in particular, $\G(\Afinite)$   acts on  $Y(K)^{\carrot}$ by ``left multiplication.''}

  As a set, $Y(K)^{\carrot}$ is finite; it is identified with the set of components of $Y(K)$. 
We will be regarding it, however, as a  groupoid: the objects are given by $ \G(\Afinite)/K$
  and the morphisms are given by left multiplication by $\overline{\G(F)}$.  
  Then there are finitely many isomorphism
  classes of objects -- we refer to these isomorphism classes as the ``underlying set'' of $Y(K)^{\carrot}$ --  and each has a profinite isotropy group (i.e., group of self-isomorphisms).  
The underlying set of $Y(K)^{\carrot}$ can thus be identified with components of $Y(K)$, and the isotropy group
  is exactly  the congruence completion of the $\pi_1$ of the
 corresponding component of $Y(K)$.  
 
  Thus, informally speaking, we regard $Y(K)^{\carrot}$ as a ``profinite orbifold'': a finite
  set of points, to each of which is associated a profinite isotropy group. \index{profinite orbifold}
  When we speak of its homology, we mean the homology of the classifying space of this groupoid --
  i.e., more concretely,   $\bigoplus_x H_*(J_x)$, the sum being taken
over a set of representatives for isomorphism classes, and $J_x$
being the isotropy group of $x$.  We will often identify $Y(K)^{\carrot}$ with its classifying space;
for instance, there is a natural continuous morphism $Y(K) \rightarrow Y(K)^{\carrot}$.

 In particular, $\G(\Afinite)$ acts on the homology of $Y(K)^{\carrot}$; 
also  the  natural $Y(K) \rightarrow Y(K)^{\carrot}$
 induces a map on homology, and we may reinterpret the prior definition:
  
  For an abelian group $A$, 
   $H_{1, \con}(Y(K),A)$   is the quotient of $H_1(Y(K))$ defined by
   $$H_{1}(Y(K), A) \twoheadrightarrow H_1(Y(K)^{\carrot}, A),$$
   and dually we define $H^1_{\con}(Y(K), A)$ as the image of
   $$H^1(Y(K)^{\carrot}, A) \hookrightarrow H^1(Y(K), A).$$ 
If $[c] \in  H^1(Y(K),-)$  lies in the image of
$H^1_{\con}(Y(K),-)$, then we say that
$[c]$ is congruence.

Now the stabilizer of every point in $Y(K)^{\carrot}$ is the same: 
it is the set of all elements of $K$ whose determinant (reduced norm, if $\G$ corresponds to a quaternion algebra) has the same square class (in $\Afinite^{\times}/(\Afinite^{\times})^2$)
as an element of $F^{\times}$. Calling this subgroup $K^+$, we may always identify
 $H^1_{\con}(Y(K),A)$
 with $$ \mbox{ functions: }  \mbox{ underlying set of }Y(K)^{\carrot} \rightarrow H^1(K^+, A)$$
 Moreover, the action
 of $\G(\Afinite)$ corresponds simply to permuting $Y(K)^{\carrot}$ considered as a set. 
 We discuss computing the cohomology of $K^+$ in~\S~\ref{CongruenceHomologyComputation}.

 \begin{remarkable} {\em 
 We may alternately characterize the congruence quotient of homology in the following way:
 It is the largest quotient of $H_1(Y(K),  A)$ in which the homology of  every ``sufficiently deep'' covering $Y(K')$ vanishes:
 $$H_{1,\con}(Y(K), A) = 
\coker \left( \lim_{\underset{K'}{\leftarrow}} H_1(Y(K') , A) \rightarrow H_1(Y(K), A) \right).$$
$$H^1_{\con}(Y(K), A) =
\ker \left(H^1(Y(K), A)  \rightarrow \lim_{\underset{K'}{\rightarrow}} H^1(Y(K'), A) \right).$$
}
\end{remarkable}

 \subsubsection{Variations: \texorpdfstring{$\q$}{q}-congruence, \texorpdfstring{$S$}{S}-arithmetic } \label{Sarithmeticcongruence} 
Fix a prime $\q$.  For simplicity in what follows we suppose that 
$\G$ is split at $\q$, and 
$K = K_{\q}
\times K^{(\q)}$, where $K^{(\q)} \subset \G(\Afinite^{(\q)})$ and $K_{\q}$ corresponds to $\PGL_2(\OO_{\q})$ under the fixed isomorphism of \eqref{isodef}.

Note that, when we regard $Y(K)^{\carrot}$ as a groupoid, all the morphisms sets are
subsets of $\overline{\G(F)}$, in particular,   subsets of $\G(\Afinite)$.
Thus we may obtain a new groupoid by projecting these morphism-sets to any quotient of $\G(\Afinite)$. In particular: 
 
\begin{itemize}
\item[(a)] (``Project at $\q$:'') Denote by $Y(K)^{\carrot}_{\q}$ the profinite-orbifold obtained by
 projecting each set of morphisms to
to $\G(F_{\q})$.   \item[(b)] (``Project away from $\q$:'') Denote by  
$Y(K)^{(\q), \carrot}$  or by $Y(K[\frac{1}{\q}])^{\carrot, \pm} $  -- the reason for the second notation will be explained below -- the profinite-orbifold
by projecting each set of morphisms to  $\G(\Afinite^{(\q)})$. 
\end{itemize}

Thus, informally, both of these have the same underlying set as $Y(K)$, but 
in the case of $Y(K)^{\carrot}_{\q}$ the isotropy groups are the projection of $K^+$
to $\G(F_{\q})$, whereas in the case of $Y(K)^{(\q), \carrot}$ the isotropy groups
are the projection of $K^+$ to $\G(\Afinite^{(\q)})$.

These two constructions actually arise naturally: 
\medskip

Firstly, one can imitate all the prior discussion but replacing all congruence subgroups
 by simply ``congruence at $\q$'' subgroups.  For example, we would replace
 the congruence completions in~\S~\ref{subsubsection:reducetogroup}
by the completions for the ``congruence at $\q$'' topology, i.e., the closures
inside $\G(F_{\q})$. 
This leads to notions of $\q$-congruence cohomology and homology. 
We may then describe the
the $\q$-congruence cohomology as the image of
 $$H^1(Y(K)^{\carrot}_{\q}, A) \rightarrow H^1(Y(K), A)$$
 induced by the natural map from $Y(K)$ to (the classifying space of) $Y(K)^{\carrot}$. 
Similarly, the $\q$-congruence homology is the quotient defined by
$H_1(Y(K)) \rightarrow H_1(Y(K)^{\carrot}_{\q})$.  

\medskip

As for $Y(K[\frac{1}{\q}])^{\carrot, \pm} $, it arises  in the  $S$-arithmetic case.  In fact,
one can alternately define it as the quotient (interpreted similarly to before)\footnote{
Indeed, we have a natural map $\G(\Afinite)/K \rightarrow  \left(\{\pm 1\} \times  \G(\Afinite^{(\q)})  \right)/K^{(\q)}$,
induced by the determinant of valuation at $\q$; and a compatible natural morphism
$\overline{\G(F)} \rightarrow \overline{\G(F)}'$.  Thus one obtains a morphism
of quotient groupoids. This morphism induces a bijection on underlying sets and
an isomorphism on isotropy groups: The isotropy group 
corresponding to $g \in \G(\Afinite)$ ``upstairs'' 
consists of elements of $g^{-1} K g$ whose determinant
has the same square class as an element of $F^{\times}$;
downstairs, the corresponding isotropy group  consists
of elements of $g^{-1} K^{(\q)} g$ whose determinant has the same
square class as an element of $F^{\times}$ which also has even valuation at $\q$,
and one projects onto the other.} 
$$\overline{\G(F)}' \backslash \left(\{\pm 1\} \times  \G(\Afinite^{(\q)})  \right)/K^{(\q)}$$
where $K^{(\q)}$ acts trivially on the $\pm 1$ factor,
$\G(F)$ acts on $\pm 1$ by valuation of determinant, 
and $\overline{\G(F)}'$ is the closure of $\G(F)$ acting on
$\pm 1 \times \G(\Afinite^{(\q)}$, rather than the whole closure
of $\G(F)$ inside $\G(\Afinite)$. 
 
\medskip

Then the evident map $\H^3 \times \Tree_{\q}^{\pm} \rightarrow \{\pm 1\}$ induces 
\begin{equation} \label{curryfavor}  Y(K[\frac{1}{\q}])^{\pm} \rightarrow Y(K[\frac{1}{\q}])^{\carrot, \pm} \end{equation} 
 and the   congruence cohomology of $Y(K[\frac{1}{\q}]^{\pm})$ (again, defined as in~\S~\ref{subsubsection:reducetogroup}, i.e.
replacing each group in~\eqref{y1qid} with its congruence completion) can be alternately described as the image on cohomology of the induced map of \eqref{curryfavor}. Again, in \eqref{curryfavor}, the right-hand side really means the classifying space
of the groupoid.

 \medskip
  Note also that $\G(\Afinite)$ still acts on $H_1(Y(K)^{(\q), \carrot})$, 
  just as in the case of $H_1(Y(K)^{\carrot})$, using again
  the ``valuation of determinant'' to determine its action on $\pm 1$.

  \subsubsection{Degeneracy maps and Hecke operators}
 \label{HeckeCongruence}

An inclusion $K_1 \subset K_2$ induces push-forward maps
$H_*(Y(K_1)^{\carrot}) \rightarrow H_*(Y(K_2)^{\carrot})$
and transfer maps in the reverse direction. 
The 
 {\em 
degeneracy maps} (see
\S~\ref{section:levelmaps}) and also Hecke operators  act on the congruence quotient of homology.

 In fact, we can compute this Hecke action:
For almost every  prime $\q$, we have an automorphism $[\q]$
  of $H_{1, \con}$ or $H^1_{\con}$,  given by the action of any element of $\G(\Afinite)$ whose reduced norm
  is an idele corresponding to $\q$. (This is well-defined away from finitely many primes). 
With this notation,
the action of $T_{\q}$ on congruence 
homology or congruence cohomology  is given by multiplication by \begin{equation} \label{cong-is-Eis}
[\q] (1 +\Norm(\q)).\end{equation}

 \medskip
 Let us prove this. 
 Let $\q$ be a prime and let $g_{\q}$
 as in \S \ref{ss:HeckeDef}, i.e. an element of $\G(F_{\q}) \subset \G(\Afinite)$ supported at $\q$
 corresponding to the $\q$-Hecke operator.

  We have a commutative diagram, where the 
  horizontal rows are defined as in~\eqref{Heckediagram}
\begin{equation} \label{HDX}
 \begin{diagram}
 H_1(Y(K)^{\carrot} ) & \rTo & 
H_1(Y(K \cap g_{\q} K g_{\q}^{-1})^{\carrot} ) & \rTo^{\phi}&
H_1(Y(K)^{\carrot}) \\
\dTo & & \dTo & & \dTo \\
 H_1(Y(K)^{(\q), \carrot} ) & \rTo & 
H_1(Y(K \cap g_{\q} K g_{\q}^{-1})^{(\q), \carrot} ) & \rTo^{\phi}&
H_1(Y(K)^{(\q), \carrot}). 
\end{diagram}
\end{equation} 

The maps $\phi$  are induced by the maps ``multiplication by $g_{\q}$.''

Recall that the bottom row is obtained, roughly speaking, by projecting
$K$ ``away from $\q$'' to $\G(\Afinite^{(\q)})$. In particular,
$Y(K \cap x K x^{-1})^{(\q), \carrot} = Y(K)^{(\q), \carrot}$. 
From this we easily see that the bottom row is simply given by 
$[\q] \cdot (\Norm(\q)+1)$, 
where $[\q]$ is as above.
Finally, the map $Y(K)^{\carrot} \rightarrow Y(K)^{(\q), \carrot}$
induces an isomorphism for $H_1$ for all but finitely many $\q$, as we shall prove below.  

\medskip
Indeed,  similar analysis shows that
  the two degeneracy maps
  $H_{1, \con}(\Sigma, \Z) \rightarrow H_{1, \con}(\Sigma/\q, \Z)$
  differ exactly by the action of $[\q]$ on the target.

 \subsubsection{Liftable congruence homology} \label{LiftableCongruenceHomology}

We define $h_{\cl}(\Sigma)$  (the subscript stands for ``liftable,'' see below for discussion) as the order of the cokernel
of $$H_1(\Sigma,\Z)_{\tors} \rightarrow H_1(\Sigma, \Z)_{\con}.$$ \index{$h_{\cl}(\Sigma)$}

In general, the image of $H_1(\Sigma, \Z)_{\tors}$ in congruence homology
need not be stable under the action of $\G(\Afinite)$.  This is basically the fact that different connected components of $Y(K)$
might have different homology, although their congruence homology is identical. Because of that fact,  we will  later need
the following variant of the definition: For any ideal class
$[\mathfrak{r}]$, we define  \index{$h_{\cl}(\Sigma;\mathfrak{r})$}
  $h_{\cl}(\Sigma; \mathfrak{r})$ 
  as the quotient of $H_1(\Sigma, \Z)_{\con}$ by the span
  of  $[\mathfrak{r}]^i H_1(\Sigma, \Z)_{\tors}$ for all $i$, i.e. the cokernel of
  $$ \bigoplus_{i} H_1(\Sigma, \Z)_{\tors} \stackrel{\bigoplus [\mathfrak{r}]^i}{\longrightarrow} H_1(\Sigma, \Z)_{\con}.$$ 
  Note that $[\mathfrak{r}]^2$ acts trivially, so $i = 0$ and $1$ suffice. 
  
  \medskip
  
  Since the map $H_1(\Sigma,\Z) \rightarrow H_1(\Sigma, \Z)_{\con}$ (and
its analogue  for $\Z_p$) is surjective,
the cokernel of the map  $H_1(\Sigma,\Z)_{\tors} \rightarrow H_1(\Sigma, \Z)_{\con}$ consists of 
the congruence classes which can only be accounted for by characteristic zero classes, that is,
they \emph{lift} to a class of infinite order. The $p$-part of the order $h_{\cl}(\Sigma)$ of
this cokernel may
 also be described as the order of
$H^1(\Sigma,\Q_p/\Z_p)_{\mathrm{div}} \cap
H^1_{\con}(\Sigma, \Q_p/\Z_p)$.

\subsection{Computation} \label{CongruenceHomologyComputation}

Under mild assumptions on $K$, it is simple to actually explicitly compute congruence cohomology and homology.  We do so in the present section (see in particular
\eqref{congcompute}). \medskip

Write $\PSL_2(F_v)$ for the image of $\SL_2(F_v)$ in $\PGL_2(F_v)$;
the quotient of $\PGL_2(F_v)$ by $\PSL_2(F_v)$
is then isomorphic to the group of square classes in $F_v$.
For the non-split quaternion algebra $D_v/F_v$, the image of $D^{1}_{v}$ in
$D^{\times}_v/F^{\times}_v$  also has index of exponent $2$.
Recall the definition of $K_v^{(1)}$: 
it is the intersection of $K_v$ with $\PSL_2(F_v)$ or the image of $D^1_v$, according
to whether $\G$ splits or not at $v$. 

\medskip

\subsubsection{\texorpdfstring{$p$}{p}-convenient subgroups} \label{sss:pcs} 
The inclusion $K_v^{(1)} \rightarrow K_v$ often induces an isomorphism on
$H_1(-,\Q_p/\Z_p)$. This property will be very useful:

\begin{definition} \label{convenient} $K_v$ is {\em $p$-convenient}
if the inclusion $K_v^{(1)} \rightarrow K_v$ induces
an isomorphism on $H_1(-, \Q_p/\Z_p)$. 
  A subgroup $K \subset \G(\Afinite)$ of the form
$\prod_{v}  K_v$ is $p$-convenient if $K_v$ is $p$-convenient for all odd $p$.
\end{definition}
\medskip

Now suppose $K = \prod K_v$, and let $K^+$
be the isotropy group of any point of $Y(K)^{\vee}$. 
The natural inclusion $K^+ \rightarrow K$ fits inside
$$\prod K_v^{(1)} \hookrightarrow K^+ \hookrightarrow \prod K_v.$$
If $K_v$ is $p$-convenient, 
 then these maps induce an isomorphism
\begin{equation} \label{congcompute} H^1(K^+, \Q_p/\Z_p) \stackrel{\sim}{\leftarrow} \prod_{v} H^1(K_v, \Q_p/\Z_p).\end{equation} 
In particular, the congruence cohomology can be computed locally from this formula and does not depend on the choice of connected component.
 Similarly, under the same assumption on $K$, the $\q$-congruence cohomology of $Y(K)$
 (respectively the congruence cohomology of $Y(K[\frac{1}{\q}])^{\pm}$)
 consists of one copy of $H^1(K_v, \Q_p/\Z_p)$
 (respectively $\prod_{v \neq \q} H^1(K_v, \Q_p/\Z_p)$) for each connected component of $Y(K)$.

Under these conditions we have an isomorphism:
\begin{equation} \label{disjointness} H^1_{\q-\con}(\Sigma, \Q_p/\Z_p) \oplus H^{1}_{\con}(\Sigma/\q, \Q_p/\Z_p)  \stackrel{\sim}{\rightarrow} H^1_{\con}(\Sigma, \Q_p/\Z_p), \end{equation} 
 where we may take the second map to be either of the degeneracy maps from level $\Sigma/\q$
 to $\Sigma$. 
 \subsubsection{}
 We now turn to showing that many of the level structures of interest are indeed $p$-convenient
 and explicitly computing their cohomologies:

 \begin{lemma} \label{lemma:congcomp} Assume that $p > 2$.  We have the following:
 \begin{enumerate}
\item[(i)] If either $p > 3$ or $\Norm({v}) = \# \OO_v/{v} > 3$, then
 $$H^1(\PGL_2(\OO_{{v}}), \Q_p/\Z_p)  = H^1(\PSL_2(\OO_{{v}}), \Q_p/\Z_p)= 0.$$
 \item[(ii)] If $p = 3$ and $\OO_{{v}}/{v} = \F_3$, then $H^1(\PSL_2(\OO_{{v}}), \Q_p/\Z_p) = 0$, but
$$H^1(\PGL_2(\OO_{{v}}), \Q_p/\Z_p) \simeq H^1(\PGL_2(\F_3),\Q_3/\Z_3)
= H^1(A_4,\Q_3/\Z_3) = \Z/3\Z.$$
\item[(iii)] If $K_{{v}} = \G(\OO_{{v}})$ and $\G$ is \emph{non-split} at ${v}$, then
$$H^1(K^{(1)}_{{v}},\Q_p/\Z_p) \stackrel{\sim}{\leftarrow}  H^1(K_{{v}},\Q_p/\Z_p) =  \Z_p/(\Norm({v})+1) \Z_p.$$
\item[(iv)] If $K_{0,v} \subset \PGL_2(\OO_{v})$ is of $\Gamma_0(v)$-type,
$$H^1(K^{(1)}_{0,{{v}}},\Q_p/\Z_p)  \stackrel{\sim}{\leftarrow}   H^1(K_{0,{{v}}},\Q_p/\Z_p) = \Z_p/(\Norm({{v}})-1) \Z_p.$$
\item[(v)] If the image of the determinant map $K_{v} \subset \PGL_2(\OO_v) \rightarrow
\OO^{\times}_v/\OO^{\times 2}_v$ is trivial, then $H^1(K^{(1)}_{v},\Q_p/\Z_p) \stackrel{\sim}{\leftarrow}  H^1(K_v,\Q_p/\Z_p)$.
 \end{enumerate} 
 In particular, if $K_v$ is satisfies one of $(i)$, $(iii)$, $(iv)$, or $(v)$, then it is $p$-convenient
 in the language of Definition~\ref{convenient}.
 \end{lemma}
 
 \begin{proof}
 By inflation-restriction, there is an exact sequence:
 $$H^1(\PSL_2(k_{{v}}),\Q_p/\Z_p)
 \hookrightarrow  H^1(\PSL_2(\OO_{{v}}), \Q_p/\Z_p)
 \rightarrow H^1(K({{v}}),\Q_p/\Z_p)^{\PSL_2(k_{{v}})},$$
 where $K({{v}})$ is the kernel of $\PSL_2(\OO_{{v}}) \rightarrow \PSL_2(k_{{v}})$.
 If $\# k_{{v}} > 3$, then $\PSL_2(k_{{v}})$ is simple and non-abelian, and hence the first group vanishes.
 If $\# k_{{v}} = 3$ and $p > 3$, then $\PSL_2(k_v)$ is a group of order prime to $p$ and the cohomology
 still vanishes.
 The group $K({{v}})$ is pro-${{v}}$, and so $H^1(K({{v}}),\Q_p/\Z_p)$ vanishes unless ${{v}}|p$.
 By Nakayama's lemma, it suffices to prove that 
 $$H^1(K({{v}}),\Z/p\Z)^{\PSL_2(k_{{v}})} =  H^1(K({{v}}),\Q_p/\Z_p)[p]^{\PSL_2(k_{{v}})}$$
  vanishes (the group $K(v)$ is pro-finite and so has no divisible quotients).
 Yet $H^1(K({{v}}),\F_p)$ is the adjoint representation of $\PSL_2(k_{{v}})$, which is irreducible 
 providing that $\# k_{{v}} \ge 3$. 
 The same argument applies to $\PGL_2(\OO_{{v}})$, with one modification: the group $\PGL_2(k_{{v}})$
 is the extension of a cyclic group of order $2$ by a simple non-abelian group, and thus
 $H^1(\PGL_2(k_{{v}}),\Q_p/\Z_p) = 0$ providing that $\# k_v > 3$ or $p > 3$.
 When $N({v}) = 3$ and $p = 3$, the same argument applies except that
 $\PGL_2(\F_3) = S_4$ and $\PGL_2(\F_3) = A_4$; the first group has  abelianization $\Z/2\Z$ whereas
 the second has abelianization $\Z/3\Z$. This suffices to prove $(i)$ and $(ii)$.

 Now suppose that $\G$ is non-split at ${v}$. Recall from~\ref{quatinfo} that we have exact sequences:
 $$0 \rightarrow L({v}) \rightarrow B^{\times}_{{v}}/\OO^{\times}_{{v}} \rightarrow l^{\times}/k^{\times} \rightarrow 0,$$
 $$0 \rightarrow K({v}) \rightarrow B^{1}_{{v}} \rightarrow l^1 \rightarrow 0,$$
 where $k = \OO_{{v}}/{v}$, $l/k$ is a quadratic extension, and $l^1 \subset l$ are the elements of norm $1$.
 (Here $K({v})$ and $L({v})$ are defined as the kernels of the corresponding reduction maps.)
  We certainly have 
  $$H^1(l^1,\Q_p/\Z_p) = H^1(l^{\times}/k^{\times},\Q_p/\Z_p) = \Z_p/(N({v}) + 1) \Z_p,$$
  and the natural map $l^1 \rightarrow l^{\times}/k^{\times}$ (which is multiplication by $2$)
  induced from the map $B^{1}_v \rightarrow B^{\times}_v/\OO^{\times}_v$ induces an isomorphism
  on cohomology. To prove the result, it suffices to prove that $H^1(K({v}),\F_p)$ has no $l^{\times}$ invariants.
  We may certainly assume that ${v}$ has residue characteristic $p$. In particular, $K(v)$ and $L(v)$
 are naturally isomorphic, since (by Hensel's Lemma)  every element in $\OO^{\times}_q$ which
  is $1 \mod {v}$ is a square.
Then, given the explicit decription of
  $B^{1}_{{v}}$ in section~\ref{quatinfo}, the Frattini quotient of $K({v})$ is
  given explicitly by $(1+ \m_{{v}})/(1 + \m^2_{{v}}) \simeq (\Z/p\Z)^2$, where
  $\m_{{v}}$ is generated by $i$.  This quotient is given explicitly by elements of the form
  $$\eta = 1 + (a + b j) i = 1 + a \cdot i - b \cdot \kij \mod 1 + \m^2_v.$$ %
  Yet $j \eta j^{-1} = \eta$ implies that  $a = b = 0$, and thus $\eta$ is trivial. This proves $(iii)$.
  
  \medskip
  
  Now suppose that $K = K_{v,0}$ is of type $\Gamma_0(v)$. There are sequences:
    $$0 \rightarrow L(v) \rightarrow K_{v,0} \rightarrow B(\PGL_2(\OO_v/v)) \rightarrow 0,$$
  $$0 \rightarrow K(v) \rightarrow K^{1}_v \rightarrow B(\SL_2(\OO_v/v)) \rightarrow 0.$$
  Here $B$ denotes the corresponding Borel subgroup.
  As above, the natural map $K^{1}_v \rightarrow K_{v,0}$ induces an isomorphism
  $$\Hom(k^{\times},\Q_p/\Z_p) \simeq H^1(B(\SL_2(\OO_v)/v) \rightarrow
  H^1(B(\PGL_2(\OO_v)) \simeq \Hom(k^{\times},\Q_p/\Z_p)$$
  on cohomology, which is induced explicitly by the squaring map..Once again we may
  assume that the residue characteristic of $v$ is $p$, and hence (since $p \ne 2$)
  there are natural isomorphisms $K(v) \simeq L(v)$.
  The action of the Borel on the adjoint representation $H^1(K(v),\F_p)$ is a non-trivial
  extension of three characters of which only the middle character is trivial; in particular
  $H^1(K(v),\F_p)^B$ is trivial, proving $(iv)$.
  
If the assumption of part $(v)$ holds, then the map $K^{(1)}_v \rightarrow K_v$ is an isomorphism,
and the claim is trivial.
 \end{proof}

\subsubsection{Variations for \texorpdfstring{$p=3$}{p=3}} \label{p3Var}
Our main numerical example concerns the congruence subgroups $\Gamma_0(\n)$ of $\PGL_2(\OO)$ for
$\OO = \Z[\sqrt{-2}]$, which admits primes $v|3$ such that $\OO/v = \F_3$.
We make here some further remarks about this case to show that it retains
good properties even in a commonly occuring case which is not $p$-convenient: 

Suppose that:
\begin{enumerate}
 \item $K_v$ is $p$-convenient for all $v$ not dividing $3$.
 \item If $v|3$, then either $K_v$  is $p$-convenient or $K_v = \PGL_2(\OO_v)$ and the following hold:
 \begin{enumerate}
 \item The square class of $-1$ belongs to $\det(K_v)$
for every finite $v$;
\item If $\lambda \in F^{\times}$ has even valuation  at all finite primes,
then the class of $\lambda \in  F^{\times}_v/F^{\times 2}_3$
coincides with the class of $+1$ or $-1$ for every $v|3$.
\end{enumerate}
\end{enumerate}
 In that case there exists an element
$\varepsilon \in K$ whose determinant has the same square class
as $-1$ at every $v$. Denote by $K_3^{(\pm 1)}$ those elements
of $K_3$ with determinant $\pm 1$. Then  we have inclusions:
$$  \langle \varepsilon \rangle  \prod_{v} K_v^{(1)}  \hookrightarrow K^+ \hookrightarrow 
\prod K_v,$$
The cohomology of the former group, with $\Q_3/\Z_3$-coefficients, is identified with
 $\prod H^1(K^{(1)}, \Q_p/\Z_p)^{\varepsilon},$
 where the superscript denotes ``$\varepsilon$-fixed.'' But (by assumption) $\varepsilon$ acts trivially on $H^1(K_v^{(1)})$ for $v$ not dividing $3$; and
 its fixed space for $v|3$ is simply
 $$\bigoplus_{v|3} H^1(\PSL_2(\OO_v),\Q_p/\Z_p)^{\eps} = \bigoplus_{v|3} H^1(\PGL_2(\OO_v),\Q_p/\Z_p) = 0.$$
 Hence we also obtain an isomorphism:
  \begin{equation} H^1(K, \Q_p/\Z_p) \stackrel{\sim}{\longleftarrow}
 \prod H^1(K_v, \Q_p/\Z_p)
 \label{congcomputetwo}\end{equation} 
 in this case (i.e., under assumptions (1) and (2) above) as well.

\section{Eisenstein classes} \index{Eisenstein} \label{section:eiz}
\begin{df}[Eisenstein classes] \label{eisdef}  A maximal ideal $\m$ of $\T_{\Sigma}$ is \emph{Eisenstein} if \label{df:Eisenstein}
$T_{\q} \equiv \psi_1(\q) + \psi_2(\q) \mod \m$ for all but finitely many $\q$, where
$\psi_1$ and $\psi_2$
are characters of the adelic class group
$$\psi_i : F^{\times} \backslash \Afinite^{\times} \rightarrow (\T_{\Sigma}/\m)^{\times}.$$
\end{df}
 
 Associated to $\m$ is a natural reducible semisimple Galois representation
 $\rhobar_{\m}$, via class field theory.

 A natural source of Eisenstein classes, for instance, is congruence homology, 
 as we saw in Remark~\ref{cong-is-Eis}. 
  This implies, then,  that $H^1_E$ and $H^1$  (and $H^E_1$ and $H_1$) coincide after localization
 at any non-Eisenstein ideal $\mathfrak{m}$ (notation as in the prior section).

More generally, there are several variants of the Eisenstein definition, 
less general than  the definition just presented, but useful in certain specific contexts: 
    \begin{itemize}
 \item[D0:] $\m$ is Eisenstein if $T_{\q}  \mbox { mod } \m= \psi_1(\q) + \psi_2(\q)$ for
 characters $\psi_i$ of $F^{\times} \backslash \Afinite^{\times}$.
 \item[D1:] $\m$ is cyclotomic-Eisenstein if $T_{\q} = 1 + N(\q) \mod \m$ for   all but finitely many {\em principal} prime ideals $\q$. 
  \item[D2:] $\m$ is congruence-Eisenstein in $\T_{\Sigma}$
  if $\m$ has support in $H_{1,\con}(\Sigma,\Z)$.
  \item[D3:]  (only makes sense for $\G$ split --) $\m$ is cusp-Eisenstein in $\T_{\Sigma}$ if it ``comes from the cusps'',
  in the sense that there does \emph{not}
  exist a Hecke equivariant splitting:
  $$H_1(\Sigma,\Z)_{\m} = \mathrm{im}(H_1(\partial \Sigma,\Z)_{\m})
  \oplus \ker(H^{\bm}_1(\Sigma,\Z)_{\m} \rightarrow H_0(\partial \Sigma,\Z)_{\m}).$$
(Here we use the notation of~\S~\ref{section:split}, which defines the right-hand groups
when $\G= PGL_2$). 
\end{itemize}
We may think of $D0$ and $D1$ as  Galois-theoretic definitions, $D2$ as a homological definition, and
$D3$ as an automorphic definition.
\S \ref{HeckeCongruence} 
  shows that
   $D2 \Rightarrow D1$, and it is easy to see that $D3 \Rightarrow D0$,
 and $D1 \Rightarrow D0$:

 It is classes of the type $D3$ that are measured by special values of $L$-functions
 (see~\S~\ref{sec:eisintegrality}). 
 Since ideals of type $D3$ can only exist in the split case, there are many examples of
 classes that satisfy $D2$ but not $D3$.
 In~\S~\ref{subsection:pathologies}, we observe a class (\hangingchad)
 that is $D1$ but not $D2$ (see Remark~\ref{remark:hangingchad}).
 Since the order of congruence homology groups is controlled by congruence conditions on the
 level, and since prime divisors of $L$-values are not bounded in such a manner, one
  expects to find  classes that satisfy
 $D3$ but not $D2$.
We summarize this discussion by the following graph describing the partial
ordering:
{\small
$$
\begin{diagram}
D2 & &  & \\
\dTo & &  &  D3 \\
D1 & &  \ldTo(3,3) &   \\
\dTo &  & &  \\
D0 & & & \\
\end{diagram}
$$
}

Throughout the text, when we say ``Eisenstein'', we shall mean Eisenstein of type $D0$ unless we specify otherwise.

\section{Automorphic representations.  Cohomological representations.}
   \label{autrepnotns}

For us {\em automorphic representation} $\pi$
is an irreducible representation 
of the Hecke algebra  of $\G(\adele)$
that is equipped with an embedding $\pi \hookrightarrow C^{\infty}(\G(\adele)/\G(F))$.

Thus $\pi$ factors as a restricted tensor product $\bigotimes_v \pi_v$,
where $\pi_v$ is an irreducible smooth representation of $\G(F_v)$
for $v$ finite, and, for $v$ archimedean, is an irreducible Harish-Chandra module
for $\G(F_v)$.  We often write $$\pi_{\infty} = \otimes_{v |\infty} \pi_v,
\pi_f = \otimes_{v \ finite} \ \ \pi_v.$$
Thus $\pi_{\infty}$ is naturally a $(\mathfrak{g}, K_{\infty})$-module. 

We say $\pi$ is {\em cohomological} if  $\pi_{\infty}$ has nonvanishing
$(\mathfrak{g}, K_{\infty})$-cohomology for every infinite place $v$, i.e., 
in the category of $(\mathfrak{g}, K_{\infty})$-modules
the group $$H^i(\mathfrak{g}, K_{\infty}; \pi_{\infty}) := \mathrm{Ext}^i(\mathrm{trivial}, \pi_{\infty})$$ does not vanish for some $i$.   In the present case, this group, if nonzero,
is isomorphic as a vector space to $\Hom_{K_{\infty}}( \wedge^i \mathfrak{g}/\mathfrak{k}, \pi_{\infty})$.

\medskip

We have Matsushima's formula:
 $$  H^i(Y(K),\C) =  \bigoplus m(\pi) H^i(\mathfrak{g},K_{\infty};\pi_{\infty}).$$
where $m(\pi)$ is the dimension of $K$-invariants on $\pi_f$.
(For an explicit map from the right-hand side for $i=1$, to $1$-forms
on $Y(K)$, 
realizing the isomorphism above, see~\eqref{omegadef}.)

\medskip

{\em Warning:} Note that it is usual to define cohomological to mean that  
it is cohomological after twisting by some finite-dimensional representation.
In that case, $\pi$ contributes to the cohomology of some
nontrivial local system on $Y(K)$. 
For our purposes in this book, it will be most convenient to regard cohomological as meaning with reference to the {\em trivial} local system.

\section{Newforms and the level raising/level lowering complexes} \label{Sec:Newforms}

\subsection{Newforms}

\begin{df}[\label{newformdef} Newforms.]
The space of newforms $H_1(\Sigma,\Z)^{\new}$ and
$H^1(\Sigma,\Q/\Z)^{\new}$ are defined respectively as: 
$$ \mathrm{coker}\left(   \bigoplus_{\q \in \Sigma \setminus S} H_1(\Sigma/\q,\Z)^2 \stackrel{\Psi^{\vee}}{\rightarrow} H_1(\Sigma,\Z) \right)$$ 
$$\mathrm{ker} \left( H^1(\Sigma,\Q/\Z) \stackrel{\Phi^{\vee}}{\rightarrow}  \bigoplus_{\Sigma \setminus S}
H^1(\Sigma/\q,\Q/\Z)^2 \right) $$
\end{df}

If $\q \in \Sigma$, we shall sometimes use the corresponding notion of $\q$-new: 
 $$H_1(\Sigma, \Z)^{\qnew} =  
 \mathrm{coker}  \left( H_1(\Sigma/\q,\Z)^2 \longrightarrow H_1(\Sigma, \Z)\right).$$ 
The corresponding notions with $\Q$ or $\C$ coefficients are obtained by tensoring. 
In particular, $H_1(\Sigma, \C)^{\qnew} = 0$ exactly when
$\dim H_1(\Sigma, \C) = 2 \dim H_1(\Sigma/\q, \C)$. 

Similarly, given some subset $T \subset \Sigma$, we may define in an analogous
way the notion of ``$T$-new'' as the cokernel of all degeneracy maps
from $\Sigma / {\q}$ to $\Sigma$, where ${\q} \in T$. 

\subsection{Level-raising and level lowering complexes.}  \label{sec:LRC}

The maps used in Definition~\ref{newformdef} extends 
in a natural way to a complex, 
 obtained by
composing the maps $\Psi^{\vee}$ and $\Phi^{\vee}$ in the obvious way, alternating
signs appropriately:
$$0 \rightarrow H_1(S,\Z)^{2^d} \rightarrow \ldots \rightarrow \bigoplus_{\Sigma \setminus S}
H_1(\Sigma/\q,\Z)^2 \rightarrow H_1(\Sigma,\Z),$$
$$H^1(\Sigma,\Q/\Z) \rightarrow  \bigoplus_{\Sigma \setminus S}
H^1(\Sigma/\q,\Q/\Z)^2  \rightarrow \ldots \rightarrow  H^1(S,\Q/\Z)^{2^d} \rightarrow 0,$$  
where $d = |\Sigma \setminus S|$. These are referred to as --- respectively --- the level raising
and level lowering complexes.  The level raising complex is exact everywhere but the last term after tensoring with $\Q$; the level lowering complex is similarly
 exact everywhere but the first term after tensoring with $\Q$.
 \medskip

 These complexes arise naturally in our paper because, in our analysis of analytic torsion, 
 we are led naturally to an alternating ratio of orders of torsion groups that exactly correspond to the complexes above. {\em But they also arise naturally in the analysis of cohomology of $S$-arithmetic groups}, as we explain in Chapter 4.    We are not able to use this coincidence in as definitive way as we should like,  owing to our lack of fine knowledge about the cohomology of $S$-arithmetic groups.

 \medskip
 
We will see Chapter~\ref{chapter:ch3}, using the relationship to the cohomology of $S$-arithmetic groups, that these complexes are {\em not} in general
 exact over $\Z$.    But has one at least the following result:
 \begin{quote}
 For every ${\q} \in \Sigma \backslash S$, there exists a natural self-map of degree $-1$ $H_{{\p}}$
of either complex so that
$d H_{{\q}} + H_{{\q}} d =  O_{{\q}}$,
where $O_{{\q}}$ is defined as $\Phi_{\q}^{\vee} \circ \Phi_\q, \Phi_\q \circ \Phi_\q^{\vee},
\Psi_\q \circ \Psi_\q^{\vee} ,\Psi_\q^{\vee }\circ \Psi_\q$, according to what makes sense on each direct summand. 
\end{quote}
For example, if we consider the level lowering complex
on the torsion-free quotients $H_{1,\tf}$, 
the only  homology arises at primes $\ell$ for which there 
exists $f \in H_1(S, \Z/\ell)$ annihilated by {\em all}
$T_{{\q}}^2 -   (\Norm(\q)+1)^2$, for ${\q} \in \Sigma-S$.

\chapter{Raising the Level: newforms and oldforms }
\label{chapter:ch3} \label{CHAPTER:CH3}

In this chapter, we give a treatment of several matters related to comparing spaces of modular forms at different levels.

\medskip

Among the main results are  Theorem~\ref{theorem:ribet} (level-raising)
and  Theorem~\ref{theorem:K2popularversion} (relationship between $K_2$ and spaces of modular forms;
see also~\S~\ref{section:Eisenstein} for the ``Galois side'' of this story, and \S \ref{section:K2examples}
 for numerical examples). 

\medskip

Our results in this chapter are relevant both in that 
they show that torsion homology behaves as the Langlands predicts (e.g. level-raising from~\S~\ref{sec:lr},
or the results of~\S~\ref{sec:clss}), and in that they are important
in our attempt to understand better 
the Jacquet-Langlands correspondence for torsion (see Chapter~\ref{chapter:ch6}).

\section{Ihara's lemma}  See also~\cite{Klosin}. We give a self-contained
 treatment.  \label{sec:ihara}
  \subsection{Remarks on the congruence subgroup property} 
 \label{CSPremarks}
 
 Let $T$ be a finite set of places of $F$. 
 
Recall that, for $\Gamma \leqslant \G(F)$ a $T$-arithmetic group, the {\em congruence kernel}
is the kernel of the map $\Gamma^* \rightarrow \widehat{\Gamma}$,  
where $\Gamma^*$ and $\widehat{\Gamma}$ are the completions of $\Gamma$ for the topologies
defined by all finite index subgroups  and congruence subgroups respectively. 
This congruence kernel depends only on $\G$, the field $F$, and the set $T$; these being fixed, it is independent of choice of $\Gamma$.   Indeed, it coincides with the kernel
of the map $\G(F)^* \rightarrow \widehat{G(F)}$, where these are the completions of $G(F)$
for the topologies defined by finite index or congruence subgroups of $\Gamma$, respectively;
these topologies on $\G(F)$ are independent of $\Gamma$. 

We refer to this as {\em the congruence kernel for $\G$ over $\OO_F[T^{-1}]$.}

To say that the congruence kernel has prime-to-$p$ order implies, in particular, the following statement: Any normal subgroup $\Gamma_1 \leqslant \Gamma$
of $p$-power index is, in fact, congruence. 
(Indeed, a homomorphism from $\Gamma$ to the finite $p$-group $P:=\Gamma/\Gamma_1$
extends to $\Gamma^* \rightarrow P$, and then must factor through $\widehat{\Gamma}$.)
 
In particular this implies that the natural map $\Gamma \rightarrow \widehat{\Gamma}$
from $\Gamma$ to its congruence completion induces 
an isomorphism $H^1(\Gamma, M) \rightarrow H^1(\widehat{\Gamma}, M)$, whenever $M$
is a $p$-torsion module with trivial action.

 \subsection{Statement of Ihara's lemma}
\begin{lemma}[Ihara's Lemma]  \label{iharalemma} 
Suppose that $\q$ is a finite prime in $\Sigma$, and  \label{lemma:ihara}
\begin{equation} \mbox{
 The congruence kernel for $\G$ over $\OO[\frac{1}{\q}]$ has prime-to-$p$ order}
\label{conj:congruence} \end{equation}
Suppose that one of the following conditions is satisfied:
 \begin{enumerate}
 \item $p > 3$,
 \item $N(\q) > 3$,
 \item $p = 3$, $N(\q) = 3$, and
\begin{enumerate}
\item The square class of $-1$ belongs to $\det(K_v)$
for every finite $v$;
\item If $\lambda \in F^{\times}$ has even valuation  at all finite primes,
then the class of $\lambda \in  F^{\times}_v/F^{\times 2}_3$
coincides with the class of $+1$ or $-1$ for every $v|3$.
\end{enumerate}
\end{enumerate}
Then:  \begin{enumerate}
\item The kernel of the level raising map $\Phi : H^1(\Sigma/\q, \Q_p/\Z_p)^2 \rightarrow H^1(\Sigma,\Q_p/\Z_p)$
is isomorphic  to the congruence cohomology $H^1_{\con}(\Sigma/\q,\Q_p/\Z_p)$, 
embedded in  $H^1(\Sigma/\q)^2$ via the twisted-diagonal embedding:
$$H^1_{\con}(\Sigma/\q, \Q_p/\Z_p) \stackrel{x \mapsto (-x, [\q] x)} \hookrightarrow
 H^1_{\con} (\Sigma/\q)  \oplus  H^1_{\con} (\Sigma/\q)  \hookrightarrow H^1(\Sigma/\q)^2,$$
where $[\q]$ is as defined in~\eqref{cong-is-Eis}.

\item The cokernel of the level lowering map
$\Psi$ on $H_1(\Sigma,\Z_p) \rightarrow H_1(\Sigma/\q, \Z_p)^2$ is $H_{1,\con}(\Sigma/\q,\Z_p)$ considered
as a quotient of $H_1(\Sigma/\q, \Z_p)^2$ via the dual of the twisted-diagonal embedding. 

\end{enumerate}
 \end{lemma}
  By a theorem of Serre \cite[p 499, Corollaire 2] {Serre} and an easy argument to pass between $\PGL_2$ and $\SL_2$, the condition above
  (concerning the congruence kernel)  holds for  $\G = \PGL_2$, 
when $p$ does 
  not divide the number of roots of unity $w_F = |\mu_F|$ in $F$.

The proof uses the cohomology of $S$-arithmetic groups.  Speaking roughly,
the proof goes as follows:  
$\G(\OO[{\q}^{-1}])$ is (almost) an amalgam of two copies of $\G(\OO)$
along a congruence subgroup $\Gamma_0(\q)$; the homology exact sequence associated
to an amalgam will yield the result. Practically, it is more useful to phrase the results
in terms of the $S$-arithmetic spaces introduced in~\S~\ref{Sarithmetic}.

 \begin{lemma} \label{lyndon-lemma} There are short exact sequences 
$$
\begin{diagram}
0 & \rTo &  H^1(  \Sigmainvertq), \Q/\Z)^{\pm} & \rTo & H^1(\Sigma/\q,\Q/\Z)^2 & \rTo^{\quad \Phi \quad} & 
H^1(\Sigma,\Q/\Z).\end{diagram}$$
and dually 
\begin{equation} \label{homology-lyndon} 
\begin{diagram}
 H_2( \Sigmainvertq, \Z)^{\pm} & \rTo &  H_1(\Sigma,  \Z) & \rTo^{\quad \Psi \quad} & H_1(\Sigma/\q, \Z)^2 & \rTo & 
H_1(\Sigmainvertq,\Z)^{\pm} &\rTo& 0. \end{diagram} \end{equation} 
\end{lemma}

\begin{proof} Apply a Mayer-Vietoris sequence to the description of $Y(\KSigmainvertq)$ --
  given in~\eqref{glueq}.  \end{proof} 
This can be translated into more group-theoretic terms, and we give that translation
in~\S~\ref{section:grouptheoretic}.

The action of the natural involution on the spaces in~\eqref{glueq}
allows us to split the sequence into ``positive'' and ``negative'' parts.
Warning: The involution on $Y(\Sigma[\frac{1}{\q}])$
induces an involution on each group of the sequence, but the induced 
map of commutative diagrams does not commute; there are sign factors, because
the sign of the connecting homomorphism in Mayer--Vietoris depends
implicitly on the choice of {\em order} of the covering sets. 

For example, the part of the sequence that computes
the positive eigenspace on $H_*(\Sigma[\frac{1}{\q}], \Z[\frac{1}{2}])$ is:
{\small
\begin{equation}
\begin{diagram}  \label{homology-lyndontwo} 
 H_2(\Sigmainvertq, \Z[\frac{1}{2}]) & \rTo &  H_1(\Sigma,  \Z[\frac{1}{2}])^{-} & \rTo^{\quad \bar{\Psi} \quad} & H_1(\Sigma/\q, \Z[\frac{1}{2}]) & \rTo & 
H_1(\Sigmainvertq,\Z[\frac{1}{2}]) \rightarrow 0, \end{diagram} \end{equation} 
}
where $H_1(\Sigma,  \Z[\frac{1}{2}])^{-} $ is the $-$ eigenspace
for the Atkin-Lehner involution $w_{\q}$, and   $\bar{\Psi}$ here denotes the map $H_1(\Sigma) \rightarrow H_1(\Sigma/\q)$
that is obtained as the difference of the two push-forward maps, or, equivalently,  the composite $H_1(\Sigma)  \longrightarrow H_1(\Sigma)^2 \stackrel{\Psi}{\rightarrow} H_1(\Sigma/\q)$, the first map being $x \mapsto (x, -x)$.

\begin{proof}  (of Ihara's lemma~\ref{iharalemma}). 
We describe the argument in cohomology, the argument in homology being dual.  

Return to the sequence
$$
\begin{diagram}
0 & \rTo &  H^1(  \Sigmainvertq, \Q_p/\Z_p)^{\pm} & \rTo & H^1(\Sigma/\q,\Q_p/\Z_p)^2 & \rTo^{\quad \Phi \quad} & 
H^1(\Sigma,\Q_p/\Z_p).\end{diagram}.$$ 

Now 
$$H^1(\Sigmainvertq, \Q_p/\Z_p)^{\pm} = H^1_{\con}(\Sigmainvertq, \Q_p/\Z_p)^{\pm},$$
this follows from the analogue of~\eqref{drm}, and the (assumed) fact that, for each of the arithmetic groups 
$X =  \Gamma_0^{(\q)}(\Sigma, \a)^{(\ev)}$  the map $X \rightarrow \hat{X}$
to its congruence completion
induces
an isomorphism on $H^1$ with $p$-torsion coefficients. 
Thus we have a diagram

$$
\begin{diagram}
   H^1_{\con}(  \Sigmainvertq, \Q_p/\Z_p)^{\pm} & \rTo & H^1_{\con}(\Sigma/\q,\Q_p/\Z_p)^2   \\\dTo^{\sim} &  & \dTo  &&   \\
  H^1(  \Sigmainvertq, \Q_p/\Z_p)^{\pm} & \rTo & H^1(\Sigma/\q,\Q_p/\Z_p)^2 & \rTo^{\quad \Phi \quad} & 
H^1(\Sigma,\Q_p/\Z_p) \end{diagram}$$
We claim that the upper horizontal arrow
has image the  twisted-diagonally embedded congruence cohomology.
Indeed, the map
$Y(K)^{\carrot} \rightarrow Y(K[\frac{1}{\q}])^{\carrot, \pm} $
induces a bijection on the underlying points, as discussed
in~\S~\ref{Sarithmeticcongruence}, 
and the computations of~\eqref{CongruenceHomologyComputation}
give the desired result. 
\end{proof}

  \subsection{Amalgams} \label{amalgams}
We shall later  (see~\S~\ref{edgemapistamesymbol})
need to compute explicitly certain connecting maps, and for this a discussion
of the homology sequence associated to an amalgamated free product will be useful.
This will also allow us to give a purely group-theoretic proof of Ihara's lemma.

Let $G$ be a group, and
let $A$ and $B$ be two groups (not necessarily distinct) together with
 fixed embeddings $G \rightarrow A$ and $G \rightarrow B$ respectively.
 Recall that the \emph{amalgamated free product} $A *_{G} B$ of $A$ and $B$ along $G$
 is the quotient of the free product $A * B$ obtained by identifying the images
 of $G$ with each other.  For short we say simply that $A *_G B$ is the ``free amalgam'' of $A$ and $B$ over $G$.

Suppose a group $X$ acts on a regular tree $\mathcal{T}$, transitively on edges
but preserving the natural bipartition of the vertices. Let 
$e$ be an edge, with stabilizer $G$, and let $A, B$
be the stabilizers of the two end vertices of $e$. Then 
there are clearly embeddings $G \hookrightarrow A, G \hookrightarrow B$, and moreover
\begin{equation}\label{serretree} \mbox{
The canonical map   $A *_G B \rightarrow X$ is an isomorphism.} \end{equation} This
fact~\eqref{serretree}  is (for example) a consequence of the Seifert--Van Kampen theorem
applied to the space $\mathcal{T} \times E /G$; here $E$
is a contractible space on which $G$ acts properly discontinuously.

Here is an example of particular interest. Let $\q$ be a place of $F$, so that $\PSL_2(F_{\q})$
 acts on the tree $\mathcal{T}_{\q}$. The remarks above give
 $$ A *_G B \stackrel{\sim}{\rightarrow} \PSL_2(F_{\q})$$
 where $A$ is the image of $\SL_2(\OL_{\q})$, 
 $B = a_{\pi} A a_{\pi}^{-1}$ where $a_{\pi} = \left( \begin{array}{cc} \pi & 0 \\ 0 & 1\end{array} \right)$ and $G$ is equal to $A \cap B$.   There is no corresponding decomposition
 of $\PGL_2(F_{\q})$ as it does not preserve a bipartition of the vertices of the tree. 
  
 We return to the general case. 
 Associated to $A$, $B$, $G$, and $A *_G B$ are the
 following long exact sequences, due to Lyndon %
 (see
 Serre~\cite{Serre}, p.169)):
 
 {\small
\begin{equation} \label{lyndon}  \rightarrow H_{i+1}( A *_G B,M) \rightarrow H_i(G,M) \rightarrow H_i(A,M) \oplus H_i(B,M)
 \rightarrow H_i(A *_G B,M) \rightarrow \ldots$$
  $$ \rightarrow H^{i}( A *_G B,M) \rightarrow H^i(A,M) \oplus H^i(B,M)
 \rightarrow H^i(G,M) \rightarrow H^{i+1}( A *_G B,M) \rightarrow \ldots \end{equation} 
 }
 which can also be derived by applying a Mayer--Vietoris sequence
 to $\mathcal{T} \times E/G$.

For later reference we detail how to construct the connecting homomorphism $H_2(X, \Z) \rightarrow
H_1(G, \Z)$, or rather its dual
\begin{equation} \label{moo2} H^1(G, \Q/\Z) \longrightarrow H^2(X, \Q/\Z)\end{equation}  

\begin{lemma} Take $\kappa \in H^1(G, \Q/\Z)$; its image in $H^2(X, \Q/\Z)$
is represented by the central extension of $X$   $$\tilde{X}_{\kappa} =  \frac{ \left( A * B   \times \Q/\Z \right) }{N},$$
where $N$ is the normal closure of $G$ in $A*B$,  embedded in $ A * B \times \Q/\Z$ via the graph
of the unique extension $\tilde{\kappa}: N \rightarrow \Q/\Z$. 
\end{lemma} 
Note that $\tilde{X}_{\kappa}$ is indeed, via the natural map, a central extension of $X$,  which is isomorphic to the quotient of $A * B$ by the {\em normal subgroup} $N$
generated by the image of $G\hookrightarrow A*B$.     Note also
that there are canonical splittings $  A \rightarrow \tilde{X}_{\kappa}$
and similarly for $B$, descending from the natural inclusion $A \hookrightarrow A*B$
and $B \hookrightarrow A*B$.

We leave the proof to the reader, but we explain why $\kappa$ extends uniquely to $N$: 
Compare
the above exact sequence~\eqref{lyndon} in low degree to the exact sequence for the cohomology of a group quotient:
{\small 
 $$
\begin{diagram}
H^1(X, \Q/\Z) &\rTo &  H^1(A*B, \Q/\Z)  &  \rTo &   H^1(N, \Q/\Z)^{X}  \\
\dTo^{\sim} && \dTo_{\sim}  & &  \dTo   \\
H^1(X, \Q/\Z) & \rTo & H^1(A,\Q/\Z) \oplus H^1(B, \Q/\Z)   &  \rTo &    H^1(G, \Q/\Z)   \end{diagram} $$
$$
\begin{diagram}
&  \rTo & H^2(X, \Q/\Z)  &\rTo& H^2(A*B, \Q/\Z) \\
 &  &  \dTo_{=} && \dTo_{\sim}  \\
&  \rTo & H^2(X, \Q/\Z) & \rTo& H^2(A, \Q/\Z) \oplus H^2(B, \Q/\Z) 
 \end{diagram} $$ }
 
 This is {\em commutative} (the only nontrivial point to be checked is the third square) from which we conclude
 that the restriction map 
 $$H^1(N, \Q/\Z)^X \rightarrow H^1(G, \Q/\Z)$$ 
 is an isomorphism, i.e.,
 every homomorphism $\kappa: G\rightarrow \Q/\Z$ extends uniquely to a homomorphism $\tilde{\kappa}: N \rightarrow \Q/\Z$ which is also $X$-invariant.

 \subsection{Group-theoretic proof of 
  Lemma~\ref{lyndon-lemma} } \label{section:grouptheoretic}
  We give, as mentioned after the proof of this Lemma, a group-theoretic proof. 
   This method of proof follows work of Ribet~\cite{Ribeticm}.
 It  will result from exact sequence~\eqref{lyndon} applied to the fundamental groups
of $Y(\KSigmainvertq)$, $Y(K_{\Sigma/\q})$, and $\YO$ (playing the roles of
$A *_G B$, $A = B$, and $G$ respectively); to be precise, we need to take into account the disconnectedness of these spaces.

 \medskip
 
 We follow the notation after Example~\ref{pgl2exampleSarithmetic};
 in particular,
$$\displaystyle{Y(\KSigmainvertq)^{\pm} = 
\coprod_{\AN_{\q}} \Gamma^{(\q)}_0(\Sigma,\a) \backslash  ( G_{\infty}/K_{\infty} \times \Doubletree_{\q}) }.$$
The action of $\Gamma^{(\q)}_0(\Sigma, \a)$ on $\Doubletree_{\q}$ induces an isomorphism:

\begin{equation} \label{rsamalg} \Gamma^{(\q)}_0(\Sigma, \a)^{(\ev)} \cong \Gamma_0(\Sigma/\q, \a) *_{\Gamma_0(\Sigma, \a)} \Gamma_0(\Sigma/\q, \a \q), \end{equation}

 We take the direct sum of the Lyndon homology sequence~\eqref{lyndon} over $\a \in A_{\q}$, i.e.
{\small 
\begin{equation} \label{lyndonihara}   H_i {\Gamma_0(\Sigma, \a)}  \rightarrow \bigoplus_{\a \in A} \left( H_i  \Gamma_0(\Sigma/\q, \a)   \oplus 
H_i \Gamma_0(\Sigma/\q, \a \q) \right)
\rightarrow 
\bigoplus_{\a \in A} H_i( \Gamma^{(\q)}_0(\Sigma, \a)^{(\ev)} ) 
\rightarrow  \dots  \end{equation}  }

For brevity, we have omitted the coefficient group. 

 As we have seen in 
~\eqref{y1qid}, the final  group is identified with $H_i(\Sigma[\frac{1}{\q}])^{\pm}$, 
 whereas both $\bigoplus_{\a \in A}   H_i  \Gamma_0(\Sigma/\q, \a)  $
 and $\bigoplus_A H_i \Gamma_0(\Sigma/\q, \a \q)  $ are identified with
 $H_i(\Sigma/\q)$. 

This gives an exact sequence as in   Lemma~\ref{lyndon-lemma};
we omit  the verification that the maps appearing are indeed $\Phi$ and $\Psi$.

\section{No newforms in characteristic zero.} \label{sec:lloneprime}

In this section, we study consequences of Ihara's lemma for {\em dual} maps
$\Phi^{\vee}, \Psi^{\vee}$ under the following assumption:
\begin{quote} There are no newforms in characteristic zero, \end{quote}
that is to say $H_1(\Sigma, \C)^{\qnew} = 0$, so that $H_1(\Sigma/\q, \C)^2 \cong H_1(\Sigma, \C)$
via the natural maps in either direction.   

The advantage of this condition is as follows: If $H_1(\Sigma, \C)$ were literally zero,
and $Y(\Sigma)$ is compact, then
the torsion in $H_1$ is self-dual, and so we may ``dualize'' the statement of Ihara's lemma if one has only torsion classes to worry about.  It turns out that the $\q$-new part being zero
is enough.

We make our task more complicated by not localizing away from Eisenstein primes.
This is because we wish to contol certain numerical factors in this case, although it is irrelevant
to the applications of~\S~\ref{section:ribet}.

\begin{lemma} \label{regulator-compare-1} \label{theorem:regcomp}
Suppose $H_1(\Sigma, \C)^{\qnew} = 0$. Suppose that $p > 2$ satisfies
the conditions of Ihara's Lemma~\ref{iharalemma}.
Then cokernel of the map
$\Psi^{\vee}_{\tf}:  H_1(\Sigma/\q,\Z_p)^2_{\tf} \rightarrow H_1(\Sigma,\Z_p)_{\tf}$
has order, up a $p$-adic unit, given by  
$$ \frac{1}{h_{\cl}(\Sigma/\q;\q)}  %
\det \left( T_{\q}^2 - (1 + N(\q))^2 \big|  H_1(\Sigma/\q, \C) \right) ,$$ 

\end{lemma}

Recall that the notation $h_{\cl}$ has been defined in~\S~\ref{LiftableCongruenceHomology}:
it roughly speaking measures congruence classes that lift to characteristic zero.
\begin{remarkable}
{\em 
\begin{enumerate}

\item
If one assumes that the congruence subgroup property holds
for $\G$ over $\OL[1/\q]$ (as one expects), then an analogous
statement holds if $\Z_p$ is replaced by $\Z$, at least up to powers
of $2$ and $3$. 

\item The determinant appearing above is none other than (the $p$-power part of)
$\prod_{f}  (a(f,\q)^2 -   (1 + N(\q))^2)$
the product ranging over a basis of Hecke eigenforms 
for $H_1(\Sigma/\q,\C)$, and
 $a(f,\q)$  denotes the   eigenvalue of
$T_{\q}$ on $f$.   
\end{enumerate} } 
\end{remarkable}
\begin{proof} By Ihara's lemma, the cokernel
of  $\Psi: H_1(\Sigma,\Z_p) \rightarrow  H_1(\Sigma/\q,\Z_p)^2.$
is   congruence homology.  
Consider the following diagram, with exact rows and columns: 
$$
\begin{diagram}
&& (\ker \Psi)_{\tors} & \rTo &  \ker \Psi & \rTo & 0  & &  \\
&&  \dTo & &  \dTo &   & \dTo & &  \\ 
    0 & \rTo &  H_1(\Sigma, \Z_p)_{\tors} & \rTo_{\Psi} &  H_1(\Sigma , \Z_p) & \rTo & H_1(\Sigma, \Z_p)_{\tf}& \rTo & 0  \\
&& \dTo_{\Psi_{\tors}} & & \dTo_{\Psi} &   & \dTo_{\Psi_{\tf}}  & & \\
 0& \rTo & H_1(\Sigma/\q, \Z_p)^2_{\tors} & \rTo_{\Psi} &H_1(\Sigma/\q, \Z_p)^2 & \rTo & H_1(\Sigma/\q, \Z_p)^2_{\tf}  &
 \rTo  & 0  \\
&& \dOnto & & \dOnto &   &  \dOnto & &  \\
&&\coker(\Psi_{\tors}) & \rTo
      & \coker(\Psi) & \rTo &  \coker(\Psi_{\tf})  & \rTo &  0 \\
\end{diagram}
$$
We see that
the order of $\coker(\Psi_{\tf})$  is
equal to  the order of the cokernel of the map
from $\coker(\Psi_{\tors})$ to  $\coker(\Psi)$,  or
equivalently, the order of the cokernel of $H_1(\Sigma/\q, \Z_p)^2_{\tors} \rightarrow
H_1(\Sigma/\q, \Z_p)_{\con}$. This map factors as
$$H_1(\Sigma/\q,\Z_p)^2_{\tors}  \rightarrow H_1(\Sigma/\q, \Z_p)^2_{\con}
\rightarrow  H_1(\Sigma/\q, \Z_p)_{\con},$$
where the last map is given by $(x,y) \mapsto (x + [{\q}] y)$. We deduce that the cokernel has order $h_{\cl}(\Sigma/\q;  {\q})$. 
\medskip

Consider the sequence:
$$H_1(\Sigma/\q,\Z)^2_{\tf}
\stackrel{\Psi^{\vee}_{\tf}}{\rightarrow} H_1(\Sigma,\Z)_{\tf} \stackrel{\Psi_{\tf}}{\rightarrow}  H_1(\Sigma/\q,\Z)_{\tf}^2.$$
Since these groups are torsion free, there
is an equality
$|\coker (\Psi^{\vee}_{\tf} \circ \Psi_{\tf} )| = |\coker(\Psi_{\tf})| \cdot |\coker(\Psi^{\vee}_{\tf})|$.
On the other hand, by 
 Lemma~\ref{lemma:composite}, we therefore deduce that
$$|\coker(\Psi_{\tf} \circ \Psi^{\vee}_{\tf})|  = 
\det \left( T_{\q}^2 - (1 + N(\q))^2 \big|  H_1(\Sigma/\q, \C) \right)  .$$
Since we have already computed the order of $\coker(\Psi_{\tf})$, this allows
us to  determine $\coker(\Psi^{\vee}_{\tf})$.
\end{proof}

\begin{theorem} \label{wotsisname} 
Suppose $H_1(\Sigma, \C)^{\qnew} = 0$.
Let $\m$ be an non-Eisenstein maximal\footnote{As the proof will make clear, it is even enough that $\m$
not be cyclotomic-Eisenstein, or in the nonsplit case, that  $\m$ not be congruence-Eisenstein.} ideal of $\T_{\Sigma}$, whose residue characteristic $p$
satisfies the conditions of Ihara's lemma.  Suppose that $\ker(\Phi^{\vee}_{\mathfrak{m}})$
is finite.

 Then 
$$\Phi^{\vee}_{\mathfrak{m}}:  H^1(\Sigma, \Q_p/\Z_p)_{\mathfrak{m}} \rightarrow H^1 (\Sigma/\q, \Q_p/\Z_p)^2_{\mathfrak{m}}$$
is surjective. 
Equivalently, by duality,  the map
$$\Psi^{\vee}: H_1(\Sigma/\q,\Z_p)^2_{\mathfrak{m}} \rightarrow H_1(\Sigma,\Z_p)_{\mathfrak{m}}$$
 is injective.  
 \label{theorem:surjective}
\end{theorem}

It should be noted that the assumption 
is rare over $\Q$, but rather common over imaginary quadratic fields,
especially after localizing at a maximal ideal $\m$. 

\begin{remark}   It is too much to ask that the map $H^1(\Sigma, \Q_p/\Z_p) \rightarrow H^1(\Sigma/\q, \Q_p/\Z_p)^2$
be surjective; this will not be so simply because of the behavior
of congruence homology at $p$.  The theorem shows that this if we 
localize away from Eisenstein ideals (thus killing the congruence homology)
there is no remaining obstruction. 
In~\S~\ref{section:essential}, we introduce refined versions of cohomology and cohomology in which
this congruence homology is excised. For an appropriate version of this cohomology,
the corresponding map of Theorem~\ref{theorem:surjective} \emph{will} be surjective; see
Theorem~\ref{essentialgood}.
 \end{remark}

\begin{proof} 
We give the proof for $\G$ nonsplit; see~\S~\ref{theoremsurjectivesplitproof}
  for the split case.

The argument will be the same whether we first tensor with $\T_{\Sigma,\m}$ or not, hence,
we omit $\m$ from the notation.
Consider the following {\em commutative} %
diagram:
$$
\begin{diagram}
0& & \ker(\Phi^{\vee}_{\mathfrak{m}}) &   & \ker(\Psi_{\mathfrak{m}})  & & 0  \\
\dTo & &  \dTo &  & \dTo && \dTo \\ 
H^1(\Sigma, \Q_p)_{\mathfrak{m}} & \rTo& H^1(\Sigma, \Q_p/\Z_p)_{\m} & \rTo_{\delta_{\Sigma}} & H_1(\Sigma, \Z_p)_{\m}  &\rTo &  H_1(\Sigma, \Q_p)_{\m} \\ 
\dTo_{\Phi^{\vee}_{\Q, \m}} && \dTo_{\Phi^{\vee}_{\m}} & & \dTo_{\Psi_{\m}} &  &  \dTo  \\
H^1(\Sigma/\q, \Q_p)_{\mathfrak{m}}^2 &\rTo & H^1(\Sigma/\q, \Q_p/\Z_p)_{\m}^2& \rTo_{\delta_{\Sigma/\q}} & H_1(\Sigma/\q, \Z_p)_{\m}^2 & \rTo & H_1(\Sigma/\q, \Q_p)_{\m}^2  \\\dOnto & & \dOnto &  &    \dTo  && \dTo \\
0 && \coker(\Phi^{\vee}_{\m}) &
  & 0&&0   \\
\end{diagram}
$$
Here every vertical column is exact ( that $\Psi_{\mathfrak{m}}$
is surjective follows from Ihara's lemma) and the middle rows are exact.

 Since $\ker(\Phi^{\vee}_{\mathfrak{m}})$ is finite, so is $\ker(\Psi_{\mathfrak{m}})$.

\medskip

Consider an element $[c]$ in $H^1(\Sigma/\q, \Q_p/\Z_p)_{\m}^2$. All elements in this group are
torsion, and hence $m[c] = 0$ for some positive integer $m$.
Lift $\delta_{\Sigma/\q}([c])$ to a class $[b] \in \HM{\constMp}{\Sigma}$. 
Since
$\Psi(m[b]) = m \delta_{\Sigma/\q}([c]) = \delta_{\Sigma/\q}(m[c]) = 0$, it follows that $m[b] \in \ker(\Psi)$. Since
$\ker(\Psi)$ is finite it follows that $m[b]$ has finite order, and hence $[b]$ has finite
order. Since $H_1(\Sigma,\Q_p)_{\m}$ is torsion free,  $\delta_{\Sigma}$ surjects onto the torsion classes
in $H_1(\YO,\Z_p)$, and $[b]$ lifts to a class  $[a] \in \CM{\constWp}{{\Sigma}}$:
$$ \delta_{\Sigma} ( [a] ) = [b].$$

Now $\Phi^{\vee}_{\mathfrak{m}}([a]) - [c]$ lies in the kernel of $\delta_{\Sigma/\q}$, i.e.
to the image of $H^1(\Sigma/\q, \Q_p)^2_{\mathfrak{m}}$;
since $\Phi^{\vee}_{\Q,\m}$ is an isomorphism
we conclude that we may modify $[a]$ so that its image is exactly $[c]$. 
 \end{proof}

\section{Level raising} \label{sec:lr}
Suppose that $\q \notin \Sigma$, and fix  a pair of levels 
$\Sigma$ and $\Sigma' = \Sigma \cup \{\q\}$. Associated
to these levels are various Hecke algebras which we now compare:
\medskip

The  Hecke algebra $\T_{\Sigma}$ at level $\Sigma$ contains $T_{\q}$,
which is not in $\T_{\Sigma'}$. We also have the Hecke algebra $\T^{\new}_{\Sigma'}$  
which is the image of $T_{\Sigma'}$ inside
the endomorphism ring of the the space of $\q$-new forms of level $\Sigma'$.

\medskip
We have the following diagrams of rings:
$$\begin{diagram}
\T^{\an}_{\Sigma} & \lOnto & \T_{\Sigma'} \\
\dInto  &   & \dOnto \\
\T_{\Sigma}  &  & \T^{\new}_{\Sigma'} \end{diagram}$$
Since $\T_{\Sigma}$ is finite over $\Z$, any
maximal ideal $\m_{\Sigma}$  of $\T_{\Sigma}$ has finite residue field. It follows that
any such maximal ideal gives rise to
a maximal ideal in $\T^{\an}_{\Sigma}$, and hence a maximal
ideal of $\T_{\Sigma'}$. 
Let us call these ideals $\m_{\Sigma}$, $\m^{\an}_{\Sigma}$, and 
$\m_{\Sigma'}$ respectively. 
  In these terms, the problem of \emph{level raising} can be phrased:
\begin{quote} When does $\m_{\Sigma'}$ give rise to a maximal ideal of  
$\T^{\new}_{\Sigma'}$? \end{quote}
 i.e., when is $\m_{\Sigma'}$ 
 obtained by pulling back a maximal ideal of $\T^{\new}_{\Sigma'}$?
 
We reformulate this in slightly more down-to-earth terms. 
The Hecke algebra $\T_{\Sigma'}$ acts faithfully (by definition)
on $H^1(\Sigma',\Q_p/\Z_p)$, and the quotient $\T^{\new}_{\Sigma'}$
acts faithfully on the subspace $H^1(\Sigma',\Q_p/\Z_p)^{\new}$. 
The ideal $\m_{\Sigma'}$ of $\T_{\Sigma'}$ gives rise to a corresponding
ideal of $\T^{\new}_{\Sigma'}$ if and only if the module
$H^1(\Sigma',\Q_p/\Z_p)^{\q-\new}$ has support at $\m_{\Sigma'}$ --
that is to say, if 
$H^1(\Sigma',\Q_p/\Z_p)^{\q-\new}[\m_{\Sigma'}] \ne 0.$
More explicitly, given a mod-$p$ eigenform of level $\Sigma$,
when  is there a mod-$p$ newform which has the same Hecke eigenvalues for $T_{\p}$ 
for all $\p$ (not dividing $\Sigma'$) as the original eigenform?

We give an answer to this question below, which is the 
 imaginary quadratic analogue of a theorem of Ribet~\cite{Ribeticm,Ribet100}. \label{section:ribet}

\begin{theorem} Assume
 the congruence kernel for $\G$ over $\OO[\frac{1}{\q}]$ has prime-to-$p$ order.
 Let
 $\m_{\Sigma} \subset \T_{\Sigma}$ be a non-Eisenstein maximal ideal;
as above,   $\m_{\Sigma}$ gives rise to a maximal ideal $\m:=\m_{\Sigma'}$ of
$\T_{\Sigma'}$. 

If $T^2_{\q} - (1 + N(\q))^2 \in \m_{\Sigma}$
then $\m$ gives rise to a maximal ideal of $\T^{\new}_{\Sigma'}$, in the sense described above.   In particular, there exist  non-zero eigenclasses $[c] \in H^1(\Sigma',\Q_p/\Z_p)^{\new}$ which
are annihilated by $\m$. If, furthermore, 
 $H^1(\Sigma,\Q_p/\Z_p)_{\m}$ is finite, then the order of
$H^1(\Sigma',\Q_p/\Z_p)_{\m}$ is strictly larger than
 $H^1(\Sigma,\Q_p/\Z_p)^2_{\m}$. 
 \label{theorem:ribet}
\end{theorem}

\begin{remarkable}
\emph{
 In the classical setting, this amounts to the assertion that, if the eigenvalue
$a_p$ of a level $\Gamma_0(N)$ modular form satisfies $a^2_p \equiv (p+1)^2$ modulo $\ell$, 
then there exists a congruent newform of level $\Gamma_0(N \ell)$. In our setting, the same result is true --- although
the newform (or the oldform, for that matter) {\em may not lift to characteristic $0$}. It turns out, in fact, that this provides a surprising
efficient way of \emph{finding} lifts of torsion classes to characteristic zero by adding $\Gamma_0(\q)$-structure
(see part~\ref{part:examples} for more discussion on this point).
}
\end{remarkable}

\begin{remarkable}
\emph{
On the level of Galois representations, the condition on Hecke eigenvalues above exactly
predicts that the residual representation $\rhobar_{\m}$ at the place
$v$ corresponding to $\q$ satisfies
$\rhobar |G_v \sim \displaystyle{\left( \begin{matrix} \eps \chi & * \\ 0 & \chi \end{matrix} \right)}$.
This, in turn, is exactly the condition that is required to 
 \emph{define} the universal deformation ring
$R^{\new}_{\Sigma'}$ recording deformations  of $\rhobar$ that are ``new of level $\Gamma_0(\q)$'' at $\q$. Thus this
theorem could be thought of in the following way: $\Tnew_{\Sigma',\m}$
is non-zero exactly when $R^{\new}_{\Sigma'}$ is non-zero. This is, of course, a  weaker version
of the claim that $R^{\new}_{\Sigma'} = \Tnew_{\Sigma',\m}$.
}
\end{remarkable}

\begin{proof} (of Theorem~\ref{theorem:ribet}). 
The Hecke algebra $\T_{\Sigma'}$ acts on all the relevant modules, and thus it makes sense to
consider the $\m$-torsion of any such module where $\m = \m_{\Sigma'}$.
Notation as in the statement, we must show that the level lowering map:
$$\Phi^{\vee}[\m]: H^1(\Sigma',\Q_p/\Z_p)[\m] \rightarrow H^1(\Sigma,\Q_p/\Z_p)[\m]^{2}$$
has a nontrivial kernel. It suffices to show the same for $(\Phi^{\vee} \circ \Phi)[\m]$, because  by Ihara's lemma, 
  the map $\Phi[\m]$ is injective (since we have completed at a non-Eisenstein ideal, we need not worry about congruence
cohomology and essential cohomology).   But $$ \Phi^{\vee} \circ \Phi = \left( \begin{matrix}   (N(\q) + 1) & T_{\q} \\
T_{\q} &  (N(\q) + 1) \end{matrix} \right).$$
 Since the ``determinant'' $ T^2_{\q} -  (1 + N(\q))^2 $
acts trivially on $H^1(\Sigma, \Q_p/\Z_p)[\m_{\Sigma}]$, 
it has  nontrivial kernel on the (possibly larger) space
 $H^1(\Sigma,\Q_p/\Z_p)^2[\m]$
and we have proved the first assertion. 
\medskip

In order to see that $H^1(\Sigma',\Q_p/\Z_p)_{\m}$ is strictly bigger than
$H^1(\Sigma,\Q_p/\Z_p)^2_{\m}$, we may assume that $\ker(\Phi^{\vee}_{\m})$ is finite.
By Theorem~\ref{wotsisname}, it follows  that $\Phi^{\vee}_{\m}$ is surjective, and hence
$$|H^1(\Sigma',\Q_p/\Z_p)_{\m}| = |\ker(\Phi^{\vee}_{\m})| \cdot |H^1(\Sigma,\Q_p/\Z_p)^2|.$$
Since we just proved that $\ker(\Phi^{\vee}_{\m})$ is non-trivial, the result follows.
(For an example of what happens when $H^1(\Sigma,\Q_p/\Z_p)_{\m}$ is infinite, see the remark at the end of the proof of this theorem.)

\end{proof}  

\begin{remarkable} \label{remark:steinbergnew} {\em 
There is a slightly larger Hecke algebra $\WT_{\Sigma'}$ which
contains $\T_{\Sigma'}$ together with the Hecke operator $U_{\q}$.
Although $U_{\q}$ does not act on $H^1(\Sigma,\Q_p/\Z_p)$, it \emph{does} act naturally
on the image of $H^1(\Sigma',\Q_p/\Z_p)$ in  $H^1(\Sigma,\Q_p/\Z_p)^2$ (this is surjective
away from Eisenstein primes). Moreover, on this space, the action of $U_{\q}$ satisfies the
equation
$$U^2_{\q} - T_{\q} U_{\q} + N(\q) = 0.$$
Under the level raising hypothesis, 
 it follows that an $\m$ in $\T_{\Sigma}$ also gives
rise to an $\m$ in $\WT^{\new}_{\Sigma'}$.
Conversely, suppose that $[c] \in H^1(\Sigma',\F_p) \subset H^1(\Sigma',\Q_p/\Z_p)$ generates a
Steinberg representation at $\q$; then, by a standard calculation, $U_{\q}$ acts on $[c]$ by $\pm 1$,
and thus (from the above relation) the corresponding ideal $\m$ of $\WT_{\Sigma'}$ may arise
via level raising only when
$1 \pm T_{\q} + N(\q) \in \m$, or equivalently when $T^2_{\q} - (1 + N(\q))^2 \in \m$.
} \end{remarkable} 

\begin{remarkable} \label{nonewcohrem} {\em 
 Statement (ii) implies that  
 $H^1(\Sigma', \Q_p/\Z_p)_{\mathfrak{m}}$ is not isomorphic 
 to $H^1(\Sigma, \Q_p/\Z_p)_{\mathfrak{m}}^2$;
 equivalently, that $H_1(\Sigma', \Z_p)_{\mathfrak{m}}$
 is not isomorphic to $H_1(\Sigma, \Z_p)_{\mathfrak{m}}^2$.
 This is {\em false} without the finiteness assumption, as we now explain:
 
 \medskip
 
Fix a level $\Sigma$ and put $Y=Y(K_{\Sigma})$. 
Let $\m$ be a non-Eisenstein maximal ideal  such that $\T/\m = \mathbf{F}_{\ell}$,
the finite field of size $\ell$; denote by $\Z_{\ell}$ 
the Witt vectors of this field.
Suppose that $H_1(Y,\Z)_{\m} \simeq \Z_{\ell}$. Suppose that
$\p$  is a prime such that
$$T_{\p} \equiv 1 + N(\p) \mod \pi, \ \ N(\p) \not\equiv -1 \mod \pi$$
At level $\Sigma \cup \{\p\}$,  let $\m$ denote the maximal ideal
of the corresponding ring $\T$ on which $U_{\p}$ and $U_{\q}$ act by $+1$.
By level raising, we know that $H_1(\Sigma \cup \p,\Z)^{\new}_{\m} \neq 0 $; nonetheless,
 is perfectly consistent 
that there is an isomorphism
$$H_1(\Sigma \cup \p,\Z)_{\m}^{-} \simeq \Z_{\ell}.$$
 
In fact, all that the the prior proof shows is that, in this case, the cokernel of the {\em transfer} homomorphism
$$H_1(\Sigma ,\Z)_{\m} \rightarrow H_1(\Sigma \cup \{\p\},\Z)_{\m}^{-}$$
is $\Z_{\ell}/\eta_{\p}$, where $\eta_{\p} = T_{\p} - 1 - N(\p)$.
 For an explicit example of this, see Example~\ref{example:data3}.

 However, it seems likely that this situation will not happen
 with level raising at {\em multiple} primes.  We analyze this 
 later in Theorem~\ref{twoprimes}. 
 } \end{remarkable} 
  
\section[The spectral sequence \dots]{The spectral sequence computing the cohomology of \texorpdfstring{$S$}{S}-arithmetic groups} \label{sec:sss}

{\em In this section~\S~\ref{sec:sss} we allow $F$ to be a general number field, i.e.,
we relax the assumption that $F$ has a unique complex place.  Most of the notation established in Chapter~\ref{sec:notn1} still applies;
we comment on adjustments when necessary. Put $R = \Z[\frac{1}{w_F^{(2)}}]$, where $w_F$ is the number of roots of unity in $F$,
and $w_F^{(2)}$ is as in~\S~\ref{ss:notn:F}. }

\medskip

We shall analyze the cohomology of analogues   of 
$Y(\KSigmainvertq)$ where one inverts not merely $\q$ but an arbitrary (possibly infinite) set $T$ of finite primes.    Most of our applications of this are conditional
on conjectures about cohomology of $S$-arithmetic groups, but for
an example of an unconditional statement see Theorem~\ref{theorem:K2popularversion}. 

\medskip

Since we work over a general field $F$ let us briefly recall our notation:
$\G$ is an inner form of $\PGL_2$, ramified at a set of places $S$. 
By a slight abuse of notation, use the same letter for the set of {\em finite} places in $S$. 
We fix a set of finite places $\Sigma$ containing $ S$ -- the ``level''
and a further set $T$ of finite places disjoint from $\Sigma$ -- the ``primes to be inverted.'' 
For $K$ an open compact subgroup of $\G(\Afinite)$, we $K^{(T)}$ be
the projection of $ K$ to $\G(\Afinite^{(T)})$.   
Finally put $G_{\infty} = \G(F \otimes_{\Q} \R)$ and $K_{\infty}$ a maximal compact
subgroup of $G_{\infty}$.

\medskip

Let $\mathscr{B}_T$ be the product of the Bruhat--Tits buildings of $\PGL_2(F_v)$, for $v \in T$;
we regard each building as a contractible simplicial complex, and so $\mathscr{B}_T$ is a contractible square complex. 
In particular, $\mathscr{B}_T$ has a natural filtration:
$$\mathscr{B}_T^0 \subset \mathscr{B}_T^1 \subset \mathscr{B}_T^2 \subset \dots$$
where $\mathscr{B}_T^{(j)}$ comprises the union of cells of dimension $\leq j$.

Consider the quotient 
$$Y(K_{\Sigma}[\frac{1}{T}])  :=  \G(F) \backslash \left( G_{\infty}/K_{\infty} \times \mathscr{B}_T \times \G(\Afinite^{(T)})    / K_{\Sigma}^{(T)} \right).$$
This is compatible with the notation of~\S~\ref{Sarithmetic}, in the case when $F$ is imaginary quadratic. As in that case, we use the abbreviations $Y(\Sigma[1/T])$ and similar notation for its cohomology. 

\medskip

This has a natural filtration by  spaces $Y_T^j$ defined by replacing $\mathscr{B}_T$
with $\mathscr{B}_T^j$. 
The space $Y_T^{j} - Y_T^{j-1}$ is seen to be a smooth manifold of dimension
$\dim(Y_{\{\infty\}})+ j$.  

\medskip
\label{epsilondef}
Let $\epsilon: T \rightarrow \{\pm 1\}$ be a choice of sign for every place $v \in T$. 
Associated to $\epsilon$ there is a natural character  $\chi_{\epsilon}  : \G(F) \rightarrow \{\pm 1\}$,
namely $\prod_{v \in T : \epsilon(v) = -1} \chi_v$; 
here $\chi_v$ is the
``parity of the valuation of determinant'', 
obtained via the natural maps $$\G(F) \stackrel{\det}{\longrightarrow} F^{\times}/(F^{\times})^2 \rightarrow \prod_{v} F_v^{\times}/(F_v^{\times})^2 \  \stackrel{v}{\longrightarrow} \pm 1,$$
where the final map is the parity of the valuation. 

Correspondingly, we obtain a {\em sheaf of $\testring$-modules}, denoted $\mathcal{F}_{\epsilon}$,  on the space $Y(K_{\Sigma}[\frac{1}{T}])$. Namely, the total space of the local system $\mathcal{F}_{\epsilon}$
corresponds to the quotient  of 
$$\left( G_{\infty}/K_{\infty} \times \mathscr{B}_T \times \G(\Afinite^{(T )}) / K^T \right) \times R$$
 by the action of $\G(F)$: the natural action on the first factor,
and the action via $\chi_\epsilon$ on the second factor. 

\medskip

Write $\Sigma = S \cup T$,  and $R = \Z[\frac{1}{ w_F^{(2)}}]$.

\begin{theorem}  \label{SarithmeticSS} There exists an $E_1$ homology spectral sequence 
 abutting to the homology groups $H_*(Y(K_{\Sigma}[\frac{1}{T}]), \mathcal{F}_{\epsilon})$, where
$$E^{1}_{p,q} = \bigoplus_{V \subset T, |V| = p} H_q(S \cup V, \testring)^{\overline{\epsilon}},$$
 the superscript $\overline{\epsilon}$ denotes  
the eigenspace where each Atkin-Lehner involution $w_{\p}$  (see~\S~\ref{subsec:ali} for definition) acts by $-\epsilon_{\p}$, for $\p \in V$. 
Up to signs, the differential 
$$d_1: H_q(S \cup V \cup \{{\q}\})^{\overline{\epsilon}} \rightarrow 
H_q(S \cup V)^{\overline{\epsilon}},$$
given by the difference (resp. sum) of the two degeneracy maps,
according to whether $\epsilon({\q}) = 1$ (resp. $-1$). 
\end{theorem}

Here, the ``difference of the two degeneracy maps'' means
$$H_q(S \cup V \cup \{\q\}, \Z)  \stackrel{ (\mathrm{id}, -\mathrm{id})}{ \longrightarrow} H_q(S \cup V \cup  \{\q\}, \Z)^2
\stackrel{\Psi}{\longrightarrow} H_q(S \cup V, \Z).$$

For example, when $\epsilon_{\p} = 1$ for every $\p$, 
  the $E_1$ term of this sequence looks like:
 \begin{equation}  \label{SSconvenientref}
\begin{diagram}
H_2(S,R)  & \lTo & \bigoplus H_2(S \cup \{\p\},\testring)^{-} & \lTo & 
\bigoplus  H_2(S \cup \{\p,\q\},\testring)^{--}  \\
 H_1(S,\testring)  & \lTo & \bigoplus  H_1(S \cup \p,\testring)^{-}  & \lTo & 
\bigoplus H_1(S \cup \{\p,\q\} ,\testring)^{--} \\
\testring \qquad & \longleftarrow &  0 \qquad  & \longleftarrow &  \qquad 0  \\
\end{diagram}
\end{equation} 

\medskip 
\medskip
\noindent
where, for instance, $  H_1(S \cup \{\p,\q\} ,\testring)^{--}$
is the subspace of $H_1((S \cup \{\p,\q\})$ where $w_{\p}, w_{\q}$ {\em both}
act by $-1$.   This sequence is converging to the cohomology of $Y(K_{\Sigma}[\frac{1}{T}])$;
if the class number of $F$ is odd, this is simply the group cohomology of $\PGL_2(\OO_F[\frac{1}{T}])$.

One could likely prove this theorem simply by iteratively applying the argument
of  Lemma~\ref{lyndon-lemma}, thus successively passing to larger and larger $S$-arithmetic groups,  but we prefer to give the   more general approach. See also the paper \cite{Moss}.

 \medskip

\proof The filtration $Y^i_T$ gives rise to a homology spectral sequence
  (see,
for example, the first chaper of~\cite{HatcherSpectral}):
$$ E^1_{pq} = H_{q}(Y_T^p, Y_T^{p-1}; \mathcal{F}_{\epsilon}) \implies H_{p+q}(Y(K_{\Sigma}[\frac{1}{T}]); \mathcal{F}_{\epsilon}).$$

The space $\mathring{Y}_T^p = Y_T^p - Y_T^{p-1}$ is a smooth manifold diffeomorphic to
$$\bigcup_{V \subset T, |V| = p} (0,1)^{V} \times Y(S \cup V)/ \langle w_{\p} \rangle,$$ 
the quotient of $Y(S \cup V) \times (0,1)^p$ by the group $W_V  = \langle w_{\p} \rangle_{\p \in V}$ 
of Atkin-Lehner involutions $w_{\p}$ for $\p \in V$; here, each $w_{\p} \ \ (\p \in V)$ acts on $(0,1)$ via $x \mapsto 1-x$;  the restriction of the sheaf $\mathcal{F}_{\epsilon}$ is the local system corresponding to the map $w_{\p} \rightarrow -\eps_{\p}$ for each
$\p \in V$.    

On the other hand, $Y_T^p$ may be identified with
$$ \left( \bigcup_{V \subset T, |V| \leq p} [0,1]^{V} \times Y(S \cup V)/   W_{V}  \right)/ \sim,$$ 
where the equivalence relation $\sim$  is generated by the rule
\begin{equation} \label{explictequiv}  (t \in [0,1]^{V \cup \q} \times y) \sim \begin{cases}  (t|_{V} \in [0,1]^{V} \times \bar{y}) ,&  t(\q) = 0 \\ 
(t|_{V} \in [0,1]^{V} \times \bar{y}') , &  t(\q) = 1; \end{cases} \end{equation} 
and $y \mapsto \bar{y}, \bar{y}'$ are the two degeneracy maps $Y(S \cup V \cup  \q) \rightarrow Y(S \cup V )$.

We have a canonical (Thom) isomorphism
$$H_{j-s}(M) \stackrel{\sim}{\longrightarrow} H_j( [0,1]^s \times M, \partial [0,1]^s \times M)$$
for any topological space $M$; here we have relative homology on the right,
and $\partial [0,1]^s$ denotes the boundary of $[0,1]^s$.  This isomorphism
sends a chain $c$ on $M$ to the relative chain $[0,1]^s \times c $ on the right-hand side. 
In our case, this yields an isomorphism
$$\bigoplus_{|V|=p} H_{q- p}(S \cup V,R)^{\bar{\epsilon}}  \stackrel{\sim}{\rightarrow} H_q(Y_T^p, Y_T^{p-1}; \mathcal{F}_\epsilon)$$
which leads to the conclusion after interpreting the connecting maps as
degeneracy maps.  \qed

It is also useful to take the direct sum over $\epsilon$ in  Theorem~\ref{SarithmeticSS}.

\begin{corollary}  \label{SarithmeticSS2} 
Let $\cF = \bigoplus_{\epsilon} \cF_{\epsilon}$, denote the direct sum\index{$\cF$} 
over all $\epsilon$ (see page \pageref{epsilondef})  and write $\Sigma = S \cup T$. As before set $R = \Z[\frac{1}{ w_F^{(2)}}]$.

There exists an spectral sequence 
 abutting to $H_*(Y(K_{\Sigma}[1/T], \cF)$, where
$$E^{1}_{p,q} = \bigoplus_{V \subset T, |V| = p} H_q(S \cup V, \testring)^{2^{d - p}}.$$
Up to signs, the differential 
$$d_1: H_q(S \cup V \cup \{{\q}\},R)  \rightarrow H_q(S \cup V,R)^{2},$$
is given by the two degeneracy maps.
\end{corollary}

\subsection{  Relationship to level lowering/raising.}

It is now clear that the cohomology of $S$-arithmetic groups is closely tied
to the level raising/lowering complex, and to make significant progress we will need to assume:

\begin{conj} Suppose  \label{conj:Eisenstein}
that $|T| = d$. Then $H_i(Y(K[\frac{1}{T}]),\Z)$ is
Eisenstein for $i \le d$. Equivalently, $H^i(Y(K[\frac{1}{T}]),  \Q/\Z)$ is Eisenstein
for $i \le d$.
\end{conj}

If $i = 1$, this conjecture is a consequence of the congruence subgroup
property (for the associated $S$-arithmetic group $\Gamma$) whenever CSP is known to hold. 
On the other hand, the conjecture for $i = 1$ 
is strictly weaker than asking that $\Gamma$ satisfy the CSP: one (slightly strengthened) version of the
conjecture for $i = 1$ says that any finite index normal subgroup of
a congruence subgroup $\Gamma(N)$ of $\Gamma$ with \emph{solvable} quotient
is congruence. It is an interesting question whether  there are be natural groups which
satisfy this property but \emph{not} the congruence subgroup property, for example,
various compact lattices arising from inner forms of $U(2,1)$.
We make some further remarks concerning this conjecture
in~\S~\ref{section:repeatingmyself}.

\medskip

 It is plausible that this conjecture holds in much greater generally;
see Question~\ref{stable-ash} in~\S~\ref{sec:qns}. 
Conjecture~\ref{conj:Eisenstein} implies that
 
\begin{quote} If $\m$ is an non-Eisenstein maximal
ideal of $\T$, and if $\CM{\constWp}{M}_{\m}$ is finite, the level raising homology complex 
(respectively, the level lowering cohomology complex) 
 is exact. 
\end{quote}

For example, suppose that $H_1(\Sigma,\Z_p)_{\m}$ is finite, $p$ is odd, and that $\m$ is non-Eisenstein.
Then the homology
spectral sequence
completed at $\m$ is empty away from the second row where it consists
of the following complex:
$$0  \leftarrow H_1(S,\Z_p)^{2^d}_{\m}
 \leftarrow \ldots \bigoplus_{\Sigma \setminus S}
H_1(\Sigma/\q,\Z_p)^2_{\m}  \leftarrow
H_1(\Sigma,\Z_p)_{\m}.$$
Thus, if $H_i(Y(K_{\Sigma}[1/T],\cF)_{\m}$ is zero in the appropriate range, this is exact
up until the final term. Since all the groups are finite, we may take the Pontryagin
dual of this sequence to deduce that the level raising complex is exact.

\subsection{Level raising at two primes produces ``genuinely'' new cohomology}
\label{section:twoprimes}

As a sample of what is implied by Conjecture~\ref{conj:Eisenstein}, we revisit the issue raised after Theorem~\ref{theorem:ribet}, namely,
that level raising does not necessarily produce  ``more'' cohomology. We  show that  Conjecture~\ref{conj:Eisenstein} predicts that there  
will \emph{always} be ``more'' cohomology when there are at least two level
raising primes.

Let $\m$ be a non-Eisenstein maximal ideal of $\T_{\Sigma}$  such that $\T/\m = \mathbf{F}_{\ell}$, a field of size $\ell$ not divisible by any orbifold prime; 
 suppose that $\T/\m$ it is generated by Hecke operators
away from $\Sigma, \p,\q$.
Let $\Z_{\ell}$ be the Witt vectors of $\mathbf{F}_{\ell}$ and let $\pi$ be a uniformizer in $\Z_{\ell}$. 

 Suppose that $H_1(\Sigma,\Z)_{\m} \simeq \Z_{\ell}$. Suppose that
$\p$ and $\q$ are primes such that
$$T_{\p} \equiv 1 + N(\p) \mod \pi, \qquad T_{\q} \equiv 1 + N(\q) \mod \pi,$$
and suppose in addition  that $N(\p), N(\q) \not\equiv -1 \mod \pi$.

\begin{theorem}  \label{twoprimes}
Assuming Conjecture~\ref{conj:Eisenstein},    $H_1(\Sigma \cup \{\p,\q\},\Z)_{\m}^{- -} \neq \Z_{\ell}.$ 
 \end{theorem} 
 
 In other words, level raising at {\em two} primes produces ``more'' cohomology. 
Contrast the situation discussed  
in Remark \ref{nonewcohrem}.

\medskip

\begin{proof}
 Suppose not.
  Since $\m$ is assumed not to be Eisenstein, $H_3(\Z)_{\m} = H_0(\Z)_{\m} = 0$.
Theorem~\ref{SarithmeticSS} gives a spectral sequence
abutting to $H_*(\Sigma[\frac{1}{\p\q}], \Z)$. The $E_1$ term of this sequence,
after localizing at $\m$, is
$$
\begin{diagram}
H_2(\Sigma,\Z)_{\m} & \lTo & H_2(\Sigma \cup \{\p\},\Z)^{-}_{\m} \oplus  H_2(\Sigma \cup \{\q\},\Z)^{-}_{\m} & \lTo & 
H_2(\Sigma \cup\{ \p,\q\}),\Z)_{\m}^{--} \\
H_1(\Sigma ,\Z)_{\m} & \lTo & H_1(\Sigma \cup \{\p\},\Z)^{-}_{\m} \oplus  H_1(\Sigma \cup \{\q\},\Z)^{-}_{\m} & \lTo & 
H_1(\Sigma \cup \{\p,\q\},\Z)_{\m}^{--} \\
0 \qquad & \longleftarrow & \qquad 0 \qquad  & \longleftarrow &  \qquad 0 \\
\end{diagram}
$$

\medskip

By Ihara's Lemma, projected to the negative eigenspace of Atkin-Lehner involutions, 
the difference-of-degeneracy maps 
$H_1(\Sigma \cup \p \cup \q, \Z_{\ell})^{--}_{\m} \rightarrow H_1(\Sigma \cup \p, \Z_{\ell})^-_{\m}$
and $H_1(\Sigma \cup \p,\Z_{\ell})^{-}_{\m} \rightarrow H_1(\Sigma,\Z_{\ell})_{\m}$
are surjective, and so must in fact be {\em isomorphisms.}

\medskip
 Put $\eta_{\p} = T_{\p} - 1 - N(\p)$  (considered as an element of $\Z_{\ell}$, 
 since it acts via $\Z_{\ell}$-endomorphisms of $H_1(\Sigma \cup \p \cup \q)$) and define $\eta_{\q}$ similarly. 

\medskip

By our remarks about Ihara's lemma, the second row of the above sequence is exact. 
 However the top row is now easily verified to be isomorphic, as
 a sequence of $\Z_{\ell}$-modules, to the sequence
   $$ \Z_{\ell} \leftarrow \Z_{\ell} \oplus \Z_{\ell} \leftarrow \Z_{\ell},$$
where the first map is $(x,y) \mapsto (\eta_{\p} x + \eta_{\q} y)$, and 
the second map is $x \mapsto (\eta_{\q} x, - \eta_{\p} x)$. 

Write $n$ for the highest power of $\pi$ dividing both $\eta_{\p}$ and $\eta_{\q}$
(in $\Z_{\ell}$).  Then
the $E^2_{1,*}$ and $E^2_{2,*}$ rows of the spectral sequence have the following shape:
$$
\begin{diagram}
\Z_{\ell}/\pi^n &     &  \Z_{\ell}/\pi^n &  & 0  \\
 &   \luTo(4,2)^d  & & &  \\
0 & & 0 & & 0 \\
\end{diagram}
$$

\medskip

We have a contradiction: our conjecture posits that $H_2(Y[1/\p\q],\Z)_{\m}$ is zero,
but the differential $d$ cannot possibly be surjective. 
\end{proof}

It is interesting in this case to consider what the space of new at $\p$ and $\q$ forms will be.
By definition, it will be the kernel of the map:
$$\Z_{\ell} \oplus \Z_{\ell} = H_1(Y_0(\p),\Z)_{\m} \oplus H_1(Y_0(\q),\Z)_{\m}
\longrightarrow H_1(Y_0(\p\q),\Z)_{\m} = \Z_{\ell} \oplus A.$$
We know that the maps send $(x,y) \in \Z_{\ell}^2$ to $\eta_{\q} x - \eta_{\p} y
\in \Z_{\ell}$, but the image in $A$ is inaccessible. If $A = k$, for example, then
the cokernel could either be $\Z_{\ell}/\pi^n$ or $\Z_{\ell}/\pi^n \oplus k$.

 \section{\texorpdfstring{$\zeta(-1)$}{zeta(-1)}  and the homology of \texorpdfstring{$\PGL_2$}{PGL2}}
 
\label{sec:clss}

{\em In this section~\S~\ref{sec:clss}, as in the prior Section,  we continue to  allow $F$ to be a general number field, i.e.,
we relax the assumption that $F$ has a unique complex place. However, we specialize
to the case of $\G = \PGL_2$, rather than an inner form thereof.   Most of the notation established in Chapter~\ref{sec:notn1}  still applies;
we comment on adjustments when necessary.}

We have tried to write this section to be as self-contained as possible, because the main 
Theorem~\ref{theorem:K2popularversion} is of independent interest.

Our main theorem asserts, roughly speaking, that the second $K$-group
$K_2(\OO_F)$ of the ring of integers ``shows up'' in the Hecke-trivial part
of automorphic homology at level $1$ (if $F$ has more than one archimedean place)
or at every prime level ${\p}$ (if $F$ is quadratic imaginary). 
This implies\footnote{The implication is via the ``Birch--Tate'' conjecture for the size of $K_2$; this is a theorem, thanks to Tate's theorem relating $K_2$ and Galois cohomology, and the 
validity of the Main Conjecture of Iwasawa theory.} for example,  that if $F$ is a totally real field of class number $1$
with $[F:\Q] > 2$,   then the odd part of the numerator of $\zeta_F(-1)$  divides the order of (the finite group)
  $H_2(\PGL_2(\OO_F), \Z)$ --- thus the title of this subsection.

As far as this book is concerned, the interest of the theorem is twofold --
it will also be shown to ``match'' with a corresponding
Galois phenomenon in~\S~\ref{section:Eisenstein}; moreover
the methods of the proof shed light on our definition of newform.

As commented, we take $\G =\PGL_2$ but allow
$F$ now to be an arbitrary field. Let $\mathfrak{p}$ be a prime ideal.  Set, as in Chapter~\ref{sec:notn1},  \begin{eqnarray}
Y_0(1) &=&   \PGL_2(F) \backslash \left(    G_{\infty}/K_{\infty} \times \PGL_2(\Afinite) / \PGL_2(\OO_{\adele})) \right)  \\
 Y_0({\p}) &=&   \PGL_2(F) \backslash \left(    G_{\infty}/K_{\infty} \times \PGL_2(\Afinite) / K_{{\p}} \right).  \end{eqnarray}
As discussed in Chapter~\ref{sec:notn1}, if $F$ has odd class number,  we have a canonical identification
\begin{eqnarray*} H_*(Y_0(1), \Z) &\cong& H_*(\PGL_2(\OO_F),   \Z).
  \end{eqnarray*}

Let $\eisideal$ be the ideal of the Hecke algebra
generated by all $T - \deg(T)$, for $T$ relatively prime to the level. 
The ideal $\eisideal$ is (some version of) the 
\emph{Eisenstein ideal}.

\begin{theorem} 
 \label{theorem:K2popularversion} 
Let $w$ be the number of roots of unity in $F$ and set
$\testring = \Z[\frac{1}{6}]$, and other notation as above. 
Then:

\begin{enumerate}

\item[(i)] If  $F$ has at least two archimedean places, 
there is a surjection
$$H_2(Y_0(1), \testring)/ \eisideal \twoheadrightarrow K_2(\OO_F) \otimes R.$$

\item[(ii)] Suppose $F $ is totally imaginary and $\q$ 
is any prime.   Suppose moreover that $H_1(Y_0(1), \C) = 0$. 
Set $$\mathcal{K}_\q =  \ker \left( H_1 (Y_0(\q), \testring)^{-} \rightarrow H_{1, \con}(Y_0(\q), R) \times H_1(Y_0(1), R)\right) .$$ 
Here the $-$ superscript is the negative eigenspace of the Atkin-Lehner involution, and the map $H_1(Y_0(\q))^- \rightarrow H_1(Y_0(1))$ is the difference
of the two degeneracy maps.\footnote{The group $\mathcal{K}_{\q}$ is a dual notion  to the new at $\q$ essential homology; see Theorem~\ref{theorem410} below.}
Then there is a surjection
\begin{equation} \label{nexttolast} \mathcal{K}_\q/\eisideal  \twoheadrightarrow K_2(\OO_F) \otimes R,\end{equation}

The same conclusion holds localized at $p$ (i.e. replacing $R$ by $\Z_p$) so long as
  no cohomological cusp form $f \in H_1(Y_0(1), \C)$ has the property that its Hecke eigenvalues $\lambda_{{\p}}(f)$
satisfy $\lambda_{{\p}}(f) \equiv \Norm(\p)+1$
modulo $p$,  
for almost all prime ideals $\p$. 

\end{enumerate}
 \end{theorem}
 
 The proof of this Theorem takes up the remainder of this section~\S~\ref{sec:clss}. 
 The basic idea is to begin with the known relationship between $K_2(F)$
 and the homology of $\PGL_2(F)$, and descend it to the ring of integers, using the congruence subgroup property many times to effect this. 
 For a numerical example of this theorem for $F = \Q(\sqrt{-491})$, see~\S~\ref{section:K2examples}.

 \begin{remarkable} Some remarks on this theorem and its proof:  \label{remark:ktrivial} {\em
 \begin{itemize}
 \item
Part (ii) of this  theorem actually arose out our attempt to understand some of our numerical observations 
({\em phantom classes}, see~\S~\ref{section:phantomclasses},~\S~\ref{subsection:pathologies}
and Remark~\ref{remark:hangingchad}): speaking very roughly, such a class
shows up at every multiple of a given level ${\frakn}$, but not at level ${\frakn}$ itself (the theorem corresponds to ${\frakn}$ the trivial ideal). 
  We note that the proof of (ii) gives a general result, but for convenience
  we have presented the simplest form here.  
  \item
 Note that the assumption in part (ii) that there are no
cusp forms of level $1$  is not necessarily a rarity --- the non-vanishing theorems
proved of Rohlfs and Zimmert (see~\cite{Rohlfs,Zimmert}) apply only to $\SL_2(\OL_F)$, not
$\PGL_2(\OL_F)$, and indeed one might expect the latter group to often have
little characteristic zero $H_1$.  

More precisely, there are several classes
of representations over $\Q$ which base-change to cohomological representations for $\SL_2(\OO_F)$ (see \cite{FinisGrunewald} for a description); but for $\PGL_2$ such classes exist only when there is a real quadratic field $F'$ such that $F F'$ is unramified over $F$;  then the base change of a weight $2$ holomorphic form of level $\disc(F')$
and Nebentypus the quadratic character attached to $F'$ gives rise to cohomology classes for $\PGL_2(\OO_F)$. See for example \cite[Question 1.8]{FinisGrunewald} for discussion of the possibility
that ``usually'' such base-change classes give all the characteristic zero homology.

\item It is  would be interesting to investigate whether the surjections of the Theorems are isomorphisms. Indeed, more generally, it appears likely that 
``Hecke-trivial'' classes in homology are related to $K$-theory.   It seems very interesting
to investigate this further (see Questions~\ref{stable-ash}  and~\ref{klanglands} in
\S~\ref{sec:qns}, and see also Conjecture~\ref{conj:hecketrivial}). 
 \end{itemize}  }\end{remarkable}

   \begin{remarkable}  \em{
   The theorem implies the existence of torsion in certain Hilbert modular varieties: 
   
   Suppose, for instance, that $F \neq \Q$ is a totally  real field of class number $1$.   Now let $\Sigma$ be the set of embeddings $F \hookrightarrow \R$
and $\mathscr{Y} := \SL_2(\OO_F) \backslash   \left( {\H}^2\right)^{\Sigma}$
the associated Hilbert modular variety. 
Write $Z_F$
   for the $2$-group $\PGL_2(\OO_F)/\SL_2(\OO_F) \cong \OO_F^{\times}/ (\OO_F^{\times})^2$.
It acts on $Y $, and the quotient $\mathscr{Y}/Z_F$ is isomorphic to what we have called $Y_0(1)$. 

We have a decomposition
   $$H_2(\mathscr{Y}, \Z[1/2]) = \bigoplus_{\chi: Z_F \rightarrow \{\pm 1\}} H_2(\mathscr{Y}, \Z[1/2])^{\chi},$$
and moreover --- denoting  by $\chi_0$ the trivial character of $Z_F$ ---
the space  of invariants $H_2(\mathscr{Y}, \Z[1/2])^{\chi_0}$ is naturally isomorphic to $H_2(Y_0(1), \Z[1/2])$,
with notations as above.

Now $H^2(\mathscr{Y}, \C)$ is generated freely by the image of differential forms
   $\omega_i := dz_v \wedge \overline{dz_v}$, where $z_v$ is the standard coordinate on the
   $v$th copy ($v \in \Sigma$) 
   of ${\H}^2$ (see \cite{Freitag}).  
   For $\epsilon \in \OO_F$, the action of the corresponding element of $Z$
   on $\omega_i$ is simply $[\epsilon]^* \omega_i = \mathrm{sign}(v(\epsilon)) \omega_i$. 
   In particular, $H_2(\mathscr{Y},  \C)^{\chi_0} = 0$, since the invariants under $-1 \in \OO_F$ are already trivial. 
   
    This implies that $H_2(\mathscr{Y}, \Z[1/2])^{\chi_0} =
   H_2(Y_0(1), \Z[1/2])$ is pure torsion; thus $H_2(\mathscr{Y}, \Z[1/2])$ has  order divisible by (the odd part of the numerator of) $\zeta_F(-1)$. 
   }\end{remarkable}

 The quotient $\T/\eisideal$ is a quotient of $\Z$. In particular, if $p$
 is any prime, then $\m = (\eisideal,p)$ is a maximal ideal of $\T$ which
 is Eisenstein (of cyclotomic type, see Definition~\ref{df:Eisenstein}).
  Part (ii) of Theorem~\ref{theorem:K2popularversion} implies the  following:
  
  \begin{theorem} \label{theorem410}  Suppose that $F$ is imaginary quadratic.
 Then, for every prime $\q$ with $\dim H_1(Y_0(\q), \C) = 0$, there is a divisibility
  $$\#K_2(\OL_F) \ |  \# H^E_1(Y_0(\q),\Z)$$
  away from the primes $2$ and $3$. 
  Indeed   $H^{E}_1(Y_0(\q),\Z_p)$ has support at~$\m = (p,\eisideal)$ for
  all primes~$\q$.
  \end{theorem}
  
  For examples of this theorem, see~\S~\ref{section:K2examples}.
  This result uses, of course, the definition of essential homology:
 it is simply the kernel of the map from $H_1$ to $H_{1, \con}$,
 and we refer forward to \S  \ref{section:essential} for discussion.
This result is an immediate consequence of the prior Theorem.

\subsection{Background on \texorpdfstring{$K$}{K}-theory}
 \label{section:ktheory}
 
If $F$ is a field, the second $K$-group $K_2 (F)$
is defined to be the universal symbol group $$F^{\times} \wedge F^{\times}/ \langle x \wedge (1-x) : x \in F-\{0,1\} \rangle.$$

If  $\p$ is a finite prime of $F$, there is a {\em tame symbol} $K_2(F) \rightarrow k_{\p}^{\times}$,
defined, as usual, by the rule
$$ x\wedge y \mapsto \frac{x^{v(y)}}{y^{v(x)}} (-1)^{v(x) v(y)} \mbox{ mod } \p.$$
  A nice reference is \cite{CSD}. The $K_2$ of the {\em ring of integers} can be defined as 
$$K_2(\OO_F) := \mathrm{ker}(K_2(F) \longrightarrow \bigoplus_{\p} k_{\p}^{\times} ).$$
This group is known to be finite.  For example, $K_2 (\Z) \cong \{ \pm 1\}$; indeed,
the map that sends $x \wedge y $ to $1$ unless $x,y$ are both negative gives an explicit isomorphism.  On the other hand, if $F = \Q(\sqrt{-303})$, 
then, according to \cite{BeG},  $K_2(\OO_F)$ has order $22$, and is generated by 
$5 \cdot (-(3 \alpha +17) \wedge (\alpha-37))$,
where $\alpha = \frac{1 + \sqrt{-303}}{2}$.

  In the present section we shall use: 
\begin{theorem} [Suslin]
Let $F$ be an infinite field. The map $F^{\times} \rightarrow B(F)$,
given by $x \mapsto \left( \begin{array}{cc} x & 0 \\ 0 & 1 \end{array} \right)$
induces an isomorphism on $H_2( - , \Z[\frac{1}{2}])$. 
  Moreover,  the induced map $$F^{\times} \otimes F^{\times} \rightarrow
H_2(B(F), \Z[\frac{1}{2}]) \rightarrow H_2(\PGL_2(F),\Z[\frac{1}{2}]),$$
where the first map is the cup-product,   is
a universal symbol for $\Z[\frac{1}{2}]$-modules. In particular, for any
$\Z[1/2]$-algebra $R$, 
$$H_2(\PGL_2(F),R) \cong K_2(F) \otimes R.$$ 
\end{theorem}  

Here $B(F)$ denotes the Borel subgroup.
 
 \begin{proof}[sketch]
We indicate a computational proof that follows along very similar
lines to that of Hutchinson~\cite{Hutchinson} and Mazzoleni~\cite{Mazzoleni} for $\GL_2(F)$ and $\SL_2(F)$
respectively (cf. also Dupont and Sah~\cite{DupontSah}.)

Firstly, note that the inclusion $x \mapsto  \left(\begin{array}{cc} x & 0 \\ 0 & 1\end{array} \right)$
induces an isomorphism
$H_i(F^{\times}) \simeq H_i(B(F))$ for $i\in \{1,2\}$ for $F$ infinite; that is a consequence of 
the Hochschild--Serre spectral sequence
applied to
$0 \rightarrow F \rightarrow B(F) \rightarrow F^{\times} \rightarrow 0$
(cf. the proof of Lemma~4 of~\cite{Mazzoleni}). 

\medskip
Denote by $C_k(\mathcal{S})$ denote the free abelian group of
distinct $k+1$ tuples of elements of  $\mathcal{S}$ (for a set $\mathcal{S}$); 
Consider the complex
$$C_0(\PPP^1(F)) \leftarrow C_1(\PPP^1(F)) \leftarrow C_2(\PPP^1(F)) \ldots$$
of $\GL_2(F)$-modules.  Computing its  $G$-hypercohomology (i.e. taking $G$-invariants
on a resolution by injective $G$-modules) gives a spectral
sequence converging to the group homology of $\GL_2(F)$; see 
  p.187 of~\cite{Hutchinson}.

 \medskip

The action of $\GL_2(F)$ on  $\mathbf{P}^1(F)$ factors through $\PGL_2(F)$ ---
and we are interested in the corresponding spectral sequence for $\PGL_2(F)$. The calculations are very similar to those
in~\cite{Hutchinson}, p.185. 
As $\PGL_2(F)$-modules, 
$C_0(\PPP^1(F))$ is the induction of the trivial module $\Z$ from the Borel subgroup $B(F)$, 
$C_1(\PPP^1(F))$ the induction of the trivial module from the torus  $F^{\times} \subset B(F)$, and 
$C_2(\PPP^1(F))$ the induction of the trivial module from the
trivial subgroup (i.e. the regular representation). 
 The first page of the spectral sequence corresponding to the complex
is given as follows (cf.~\cite{Hutchinson}. p.185):

$$\begin{diagram}
H_2(B) &\lTo& H_2(F^{\times}) &\lTo& 0 \\
H_1(B) &\lTo& H_1(F^{\times}) &\lTo& 0 \\
H_0(B)   &  \lTo & H_0(F^{\times})  &  \lTo&
\Z & \lTo &  \bigoplus_{x \notin \{\infty,0,1\}} [x] \Z  & \lTo&  \bigoplus_{x \ne y} [x,y] \Z
\end{diagram}
$$ 
The maps $H_i(F^{\times}) \rightarrow H_i(B)$ 
are given by the composition $w - 1$ with the natural mapping $H_i(F^{\times}) \rightarrow H_i(B)$.  Here $w$
is an element of the normalizer of $F^{\times}$ not belonging to $F^{\times}$.
 Yet the action of $w$ on   $H_1(F^{\times})$ is given by
$$x = \left(\begin{matrix} x & 0 \\0 & 1 \end{matrix} \right) \mapsto
 \left(\begin{matrix} 1 & 0 \\0 & x \end{matrix} \right)  = x^{-1} =(-1) \cdot x,$$
 where the latter indicates the action of $-1$ on $F^{\times} = H_1(F^{\times})$ as a $\Z$-module.
 In particular, the map $H_1(F^{\times}) \rightarrow H_1(B)$ has kernel $\mu_2(F)$, 
 and is an isomorphism after tensoring
 with any $\Z[1/2]$-algebra $R$. 
In particular, the transgression $d_2: E^2_{3,0} \rightarrow E^2_{1,1}$
is trivial after tensoring with such an $R$. 
 Similarly, the map $H_2(F^{\times}) \rightarrow H_2(B)$ is the zero map,
 since $w$ acts trivially on $H_2(F^{\times})$. 
The $\PGL_2$-spectral sequence thus yields the following
exact sequence after tensoring with $R$:
$$P(F) = E^2_{0,3} \rightarrow H_2(B(F)) \rightarrow H_2(\PGL_2(F)) \rightarrow E^{2}_{1,1} \simeq \mu_2(F) \rightarrow 0.$$

In the diagram above, the map from  
$\bigoplus_{x \notin \{\infty, 0, 1\}} [x]\Z$ to $\Z$ is trivial; thus $P(F)$  is generated by classes $[x]$ for $x \in \mathbf{P}^1(F) - \{0,1,\infty\}$
and the image of $[x]$  under
$$ P(F) \rightarrow H_2(B(F)) \stackrel{\sim}{\rightarrow} H_2(F^{\times}) = \wedge^2 F^{\times}$$
is is $2(x \wedge (1-x))$. This is, in effect,
the computation carried out  in the Appendix of~\cite{DupontSah};
the group $P(F)$ corresponds to the group described in (A28) of that reference, 
and the map is their $\varphi$. 

 Hence, tensoring with any $\Z[1/2]$-algebra $R$, we obtain
the desired result.
 \end{proof}

 \subsection{Proof of theorem~\ref{theorem:K2popularversion}}
 We handle first the case when $F$ has more than one archimedean place. In particular,  
 by a result of Serre (\cite{SerreCSP}, Theorem~2, p.498),
 $\SL_2(\OO_F)$ has the congruence subgroup property ``tensored with $R$'', i.e.
 the congruence kernel completed at $R$ is trivial. (The congruence kernel
 is a finite group of order $|\mu_F| = w_F$ which is invertible in $R = \Z[1/6]$.)

We shall apply the spectral sequence of Theorem~\ref{SarithmeticSS}, with
 $\Sigma$ to be the empty set,
 $T$ the set of all finite places, and $\epsilon = +1$ for all primes.
Thus one obtains an $E_1$
sequence converging to $H_*(\PGL_2(F))$: 
 
\begin{diagram}
H_2(1,R)  & \lTo & \bigoplus H_2(\p,\testring)^{-} & \lTo & 
\bigoplus  H_2( \{\p, \q\},\testring)^{--}  \\
 H_1(1,\testring)  & \lTo & \bigoplus  H_1(\p,\testring)^{-}  & \lTo & 
\bigoplus H_1(\{\p,\q\} ,\testring)^{--} \\
H_0(1,R)  & \lTo & \bigoplus H_0(\p,\testring)^{-} & \lTo &  
\bigoplus H_0 ( \{\p, \q\},\testring)^{--} \\
 \end{diagram}

Note that the set of connected components of $Y_0(1), Y_0(\p), Y_0(\p \q)$, etc.
are all identified, via the determinant map, with the quotient $C_F/C_F^2$ of the class group of $F$ by squares.  In what follows we denote this group by $C$, i.e.
$$ C := C_F/C_F^2.$$ 
 Consequently, we may identify $H_0(Y_0(\p), R)$ with
the set of functions  $C \rightarrow R$; now, the action of $w_{\p}$ is then given by multiplication by the class of $\p$, and so we may identify
$$H_0(Y_0(\p), R)^{-} \cong \{ f: C \rightarrow R : \ \ f(I \p) = - f(I)\},$$
and so on.  

\subsubsection{The \texorpdfstring{$H_0$}{H_0}-row.} 
\label{section:H0row}
We claim that the first row of the spectral sequence is exact. 

Indeed if $M$ is any $R[C]$-module,
we may decompose $M$ according to characters $\chi: C \rightarrow \pm 1$ of the elementary abelian $2$-group $C$, i.e.
$$  M = \bigoplus_{\chi} M_{\chi}, \ \ M_{\chi} = \{ m \in M: c \cdot m = \chi(c) m \}.$$
The first row of the spectral sequence splits accordingly; 
if we write $S(\chi)$ for the set of primes $\p$ for which $\chi(\p) = -1$, 
the $\chi$-component is isomorphic to the sequence
\begin{equation} \label{schidef} \mathcal{S}_{\chi}^{\bullet}: R \leftarrow \bigoplus_{\p \in S(\chi)} R \leftarrow \bigoplus_{\p,\q \in S(\chi)} R \leftarrow \dots \end{equation}
This is verified to be  acyclic
as long as $S(\chi)$ is nonempty; by Chebotarev's density, this is so as long as $\chi$ is nontrivial.

In particular, the first row of $E_2$ looks like $R,  \ 0,   \ 0,  \ 0 \dots$.

 \subsubsection{The \texorpdfstring{$H_1$}{H_1}-row.}
 We now analyze similarly the second row of the spectral sequence,
 $$ H_1(1,\testring)  \leftarrow \bigoplus  H_1(\p,\testring)^{-}   \leftarrow 
\bigoplus H_1(\{\p, \q\}  ,\testring)^{--} \leftarrow \cdots $$
We claim it is exact away from the left-most term. 
 
 \medskip
 
The covering $Y_1(\p) \rightarrow Y_0(\p)$ with Galois group $k_{\p}^\times$ (see before~\eqref{A def} for the definition of the space $Y_1$) induces a map
\begin{equation} \label{ra} \theta: \ H_1(Y_0(\p), \Z) \rightarrow k_{\p}^{\times}.\end{equation}

If $X$ is any connected component of $Y_0(\p)$, the congruence subgroup property  implies that the induced map
$H_1(X, R) \rightarrow k_{\p}^{\times} \otimes R$ is an {\em isomorphism}.  (Recall that
we are supposing that $F$ has more than one archimedean place at the moment.)

Therefore, we may canonically identify
$$H_1(Y_0(\p), R) \cong \left( k_{\p}^{\times} \otimes R \right)^{C} =  \{ \mbox{ functions } f: C \rightarrow k_{\p}^{\times}\otimes R\}.$$ 
Explicitly, given a homology class $h \in H_1$, 
we associate to it the function given by $c \in C \mapsto \theta(h|_c)$, 
where $h|_c$ is the homology class that agrees with $h$ on the connected component corresponding to $c$, and is zero on all other components. 
\medskip

Let us compute the action of $w_{\p}$.  Firstly, $w_{\p}$ negates  $\theta$: the automorphism $w_{\p}$ of $Y_0(\p)$
extends to an automorphism $W$ of the covering $Y_1(\p) $
such that$W a W = a^{-1}$ for $a \in k_{\p}^{\times} = \mathrm{Aut}( Y_1(\p) / Y_0(\p) )$. 
Consequently, the action of $w_{\p}$ is given by translation by $\p$ followed by negation.\footnote{Here is the proof: given a homology class $h \in H_1$, with associated function $f_h: C \rightarrow k_{\p}^{\times}$, we have
$$ f_{w_{\p} h}(c) = \theta((w_{\p} h)|c) = \theta( w_{\p}  h|_{w_{\p} c} ) = - h|_{w_{\p} c}
= - f_h(\p c).$$}

Similarly,   
 $$H_1(Y_0(\{\p, \q\}), R) \cong (\left( k_{\p}^{\times} \oplus k_{\q}^{\times} \right)   \otimes R )^{C} =  \{ \mbox{ functions } f: C \rightarrow  \left( k_{\p}^{\times} \oplus k_{\q}^{\times}\right) \otimes R \}.$$
 Here the action of $w_{\p}$ translates by $\p$ on $C$
and by $(-1, 1)$ on $k_{\p}^{\times} \oplus k_{\q}^{\times}$;
the action of $w_{\q}$ translates by $\q$ on $C$ and acts by $(1,-1)$ on 
$k_{\p}^{\times} \oplus k_{\q}^{\times}$.

 Write $V(\p\q)$ for the $--$ part of this space, which we may write as a sum 
 of a $\p$-piece and a $\q$-piece:
\begin{eqnarray*} V(\p\q)_{\p} &=&  \{ f: C \rightarrow  k_{\p}^{\times} \otimes R :
 f(\p I) =  f(I), \ f(\q I) =-  f(I).\} \\ V(\p\q)_{\q} &=&  \{ f: C \rightarrow  k_{\q}^{\times} \otimes R :
 f(\p I) = -  f(I), \ f(\q I) =  f(I).\} 
\end{eqnarray*} 
     The $H_1$-row of the above spectral sequence is thereby identified with 
     \begin{equation} \label{h1m} R \otimes \left(   
     \bigoplus_{\p} V(\p) \leftarrow \bigoplus_{\p, \q} V(\p\q)_{\p} \oplus V(\p\q)_{\q} \leftarrow \dots \right).  \end{equation} 
     
For example, the morphism $V(\p \q)_{\p} \rightarrow  V(\p)$
sends a function $f$ to $f(\q I) - f(I)= 2f$;
this follows from the existence of
a commutative diagram
\begin{small}
$$ \begin{diagram}
Y_1(\{\p, \q\}) &  &  \rTo  &  &  Y_1(\p) \\
\dTo  & &  & & \dTo  \\
Y_0(\{\p, \q\})  &   &  \rTo  & &    Y_0(\p) \\
\end{diagram}$$
\end{small}
that respects the map $k_{\p}^{\times} \oplus k_{\q}^{\times} \rightarrow k_{\p}^{\times}$
on automorphism groups.

     This sequence~\eqref{h1m} decomposes in a natural way indexed by prime ideals $\p$ (the $\p$ piece is ``all those terms involving $k_{\p}^{\times}$; for instance $V(\p\q)_{\p}$
     is the $\p$-part of $V(\p\q)$). 
We further split~\eqref{h1m} up 
  according to characters $\chi: C \rightarrow \pm 1$. 
  The $\chi,\p$-component of the above sequence
is only nonzero if $\chi(\p)= 1$. In that case, it is identified  with 
$$ k_{\p}^{\times} \otimes   \mathcal{S}_{\chi}^{\bullet},$$
at least up to multiplying the morphisms by powers of $2$; here
  $\mathcal{S}_{\chi}^{\bullet}$ is as in~\eqref{schidef}.   As before,
  this is acyclic unless $\chi=1$. In that case, its homology is concentrated
  in one degree and is isomorphic to $R \otimes k_{\p}^{\times}$.

We  conclude
  that the second row of $E_2$ looks like 
  
  \begin{equation} \label{e2abut}
  \begin{diagram}
H_2(Y_0(1),  \testring) / \bigoplus H_2(Y_0(\p),\testring)^{-}    &  &  ? &   & 
\qquad ? \\
0  \qquad & & \bigoplus_{\p} k_{\p}^{\times} \otimes R &  & 
\qquad 0 \\
\testring \qquad &   &  0 \qquad  &   &  \qquad 0  \\
\end{diagram}
 \end{equation}
  
\vspace{0.5cm}

\subsubsection{The edge exact sequence} \label{section:stuffabouttame}
From~\eqref{e2abut} we now obtain an edge exact sequence 
 $$ H_2(Y_0(1), \testring) / \bigoplus_{\p} H_2(Y_0(\p),\testring)^{-}   \longrightarrow H_2(\PGL_2(F), \testring)  \stackrel{\theta}{ \longrightarrow} \bigoplus_{\p}  (k_{\p}^{\times} \otimes \testring) ,$$
  where $\theta$ is the morphism arising from the spectral sequence. 
  We claim that the composite
    $$K_2 F \otimes R \stackrel{\sim}{\leftarrow } H_2(\PGL_2(F), R) \stackrel{\theta}{\rightarrow} 
   \bigoplus  (k_{\p}^{\times} \otimes \testring) $$   
is {\em none other than the tame symbol map}.

This will complete the proof of the theorem
in the case when $F$ has more than two archimedean places: 
we will have exhibited an isomorphism \begin{equation} \label{almostone}  H_2(Y_0(1), \testring) / \bigoplus H_2(Y_0(\p),\testring)^{-}  \cong K_2(\OO) \otimes \testring,\end{equation}
and it remains only to observe that
for any $h \in H_2(Y_0(1), \testring)$
and any prime $\q$, if $\alpha, \beta : Y_0(\q) \rightarrow Y_0(1)$
are the two degeneracy maps, that
\begin{eqnarray} \label{firstline}
 T_{\q} h = \alpha_* \beta^* h &=& \beta_* \beta^* h + (\alpha_*-\beta_*) \beta^* h \\   \nonumber &=&
(\Norm(\q) + 1) h + (\alpha_*-\beta_*) \beta^* h 
\\  \nonumber &\in& (\Norm(\q) + 1) h +  \mathrm{image}(H_2(Y_0(\q), \testring)^{-}).\end{eqnarray}

In other words, each Hecke operator $T_\q - \Norm(\q) - 1$
acts trivially 
on the left-hand of~\eqref{almostone}; therefore $\eisideal$ annihilates this quotient, completing the proof.

    \subsubsection{The edge map is (an invertible multiple of) the tame symbol} \label{edgemapistamesymbol}
    It remains to check that the map $K_2 \rightarrow k_{\p}^{\times}$ is the tame symbol map. 
    To do, compare the above sequence with a corresponding sequence for $F_{\p}$.
    This can be done via the natural maps $\PGL_2(F) \rightarrow \PGL_2(F_{\p})$ and
$ \bar{S} \times \mathscr{B}_T \times \G(\adele^T) / K^T 
\rightarrow \mathscr{B}_{\p}.$  These morphisms preserve the filtrations on the spaces (the filtration on $\mathscr{B}_{\p}$ being two-step: the vertices and the whole space). One obtains then a morphism of the corresponding spectral sequences and thereby  a commuting sequence: 
$$
\begin{diagram}
H_2(Y_0(1),  \testring) / \bigoplus H_2(Y_0(\p),\testring)^{-}    & \rTo &  H_2(\PGL_2(F), \testring)  & \rTo&  \bigoplus  (k_{\p}^{\times} \otimes \testring)  \\
   \dTo &  & \dTo & & \dTo   \\
H_2(\PGL_2(\OO_{\p}), \testring) / H_2(K_{0,\p}, \testring)^{-}    & \rTo &  H_2(\PGL_2(F_v), \testring)  & \rTo&   k_{\p}^{\times} \otimes \testring  \\
  \end{diagram}
$$
 
This reduces us to computing the morphism $H_2(\PGL_2(F_{\p}), R) 
\rightarrow k_{\p}^{\times} \otimes R$.    
This morphism arises from the action of $\PGL_2(F_{\p})$
on the tree $\mathscr{B}_p$; considering the action of $\PSL_2(F_{\p})$ (the image of $\SL_2$ in $\PGL_2$)  on the same gives  a similar map $H_2(\PSL_2(F_{\p}), R) \rightarrow k_{\p}^{\times} \otimes R$. This yields a commutative diagram 

\begin{equation}  \label{pgl2psl2}
\begin{diagram}
(F_p^{\times} \wedge F_{\p}^{\times}) \otimes R & \rTo &  H_2(\PSL_2(F_{\p}), \testring)  & \rTo&     (k_{\p}^{\times} \otimes \testring)  \\
   \dTo^{\times 4} &  & \dTo & & \dTo   \\
(F_p^{\times} \wedge F_{\p}^{\times}) \otimes R    & \rTo &  H_2(\PGL_2(F_{\p}), \testring)  & \rTo&   k_{\p}^{\times} \otimes \testring  \\
  \end{diagram}
\end{equation} 
Here the first top (resp. bottom)  horizontal map is induced by $x\mapsto \left( \begin{array}{cc} x& 0 \\0 &x^{-1} \end{array}\right)$ (resp. $x\mapsto \left( \begin{array}{cc} x& 0 \\0 &1\end{array}\right)$) together with the natural identification $H_2(F_{\p}^{\times}, R) \cong \left( F_{\p}^{\times} \wedge F_{\p}^{\times} \right) \otimes R$.  So it is enough to check that the top row is a multiple of the tame symbol.

We check that by passing to group homology, interpreting the connecting map as a connecting map in the Lyndon sequence, and using the explicit description already given. 
Let $\pi$ be a uniformizer in $F_{\q}$. Take $X = \SL_2(F_{\q}), A = \SL_2(\OL_{\q})$ and $ B =
\left( \begin{array}{cc} \pi& 0 \\0 &1\end{array}\right) A \left( \begin{array}{cc} \pi^{-1} & 0 \\0 &1\end{array}\right) $ and $G = A \cap B$;
then the natural map $A *_G B \rightarrow X$ is an isomorphism. 
Moreover, the connecting homomorphism in group homology 
$$H_2(X, \Z) \rightarrow H_1(G, \Z) \rightarrow k_{\p}^{\times}$$
the latter map induced by $\left( \begin{array}{cc} a & b \\c&d \end{array}\right) \mapsto a$,
is identified with the lower composite of~\eqref{pgl2psl2}. In these terms, we need to evaluate 
\begin{equation} \label{urrghl} \wedge^2 F_{\q}^{\times} = H_2(F_{\q}, \Z) \rightarrow H_2(X, \Z) \rightarrow H_1(G, \Z)\rightarrow k_{\p}^{\times},\end{equation} 
where the first map is again induced by $x \mapsto \left(\begin{array}{cc} x & 0 \\ 0 & x^{-1} \end{array} \right)$.  We will show this is twice the tame symbol.

 Let $\kappa: k_{\p}^{\times} \rightarrow \Q/\Z$ be a homomorphism;  
 we will prove that the composite \begin{equation} \label{urg2} \wedge^2 F_{\q}^{\times} = H_2(F_{\q}, \Z) \rightarrow H_2(X, \Z) \rightarrow H_1(G, \Z)\rightarrow k_{\p}^{\times} \stackrel{\kappa}{\rightarrow} \Q/\Z \end{equation}
 takes $u \wedge \pi \mapsto \kappa(\bar{u})^2$ 
 for any unit $u \in \OO_q$ and any uniformizer $\pi$. This shows
 that~\eqref{urrghl} is twice the tame symbol.

By pullback we regard $\kappa$ as  a homomorphism $\kappa: G \rightarrow \Q/\Z$. By the discussion following~\eqref{moo2}
 we associate to it a central extension $\tilde{X}_{\kappa}$ of $X$ by $\Q/\Z$, equipped
 with splittings over $A$ and $B$. 
  As a shorthand, we denote these splittings as
$ g \in A \mapsto g^{\tilde{A}} \in \tilde{X}_{\kappa},$
and similarly for $B$.  Note that, for $g \in G$, 
$$g^{\tilde{A}} = g^{\tilde{B}} \kappa(g).$$

Set 
$a_x = \left(\begin{array}{cc} x & 0 \\ 0 & x^{-1} \end{array}\right)$.
Let $\widetilde{a_x}$ denote any lift of $a_x$ to $\tilde{X}_{\kappa}$. 
The composite ~\eqref{urg2}
 is given by \begin{equation} \label{commaa} x \wedge y \mapsto  \widetilde{a_x} \widetilde{a_y} \widetilde{a_x}^{-1} \widetilde{a_y}^{-1} \in \mathrm{ker}(\tilde{X}_{\kappa} \rightarrow X) \cong \Q/\Z.\end{equation}    Note this is independent of choice of lift.
  
    Only an ugly computation remains. 
For $u \in \OL_{\q}^{\times}$ a unit,  we may suppose that $\widetilde{a_u} = a_u^{\tilde{A}}$.
On the other hand, we may write $a_\pi = 
     \left( \begin{array} {cc} 0 & 1 \\ -1 & 0 \end{array}\right) 
   \left( \begin{array} {cc} 0 & -\pi^{-1} \\ \pi^{} & 0 \end{array}\right) $, and thus
   may suppose that
   $\widetilde{a_{\pi} } =      \left( \begin{array} {cc} 0 & 1 \\ -1 & 0 \end{array}\right)^{\tilde A} 
   \left( \begin{array} {cc} 0 & -\pi^{-1} \\ \pi^{} & 0 \end{array}\right)^{\tilde B} $.
   
Then
$ \widetilde{a_\pi} \widetilde{a_u} \widetilde{a_\pi}^{-1} \widetilde{a_u}^{-1} 
$ equals 
{\footnotesize
$$ \begin{aligned}
    &  \left( \begin{array} {cc} 0 & 1 \\ -1 & 0 \end{array}\right)^{\tilde A} 
   \left( \begin{array} {cc} 0 & -\pi^{-1} \\ \pi^{} & 0 \end{array}\right)^{\tilde B} 
      \left( \begin{array} {cc} u &  0 \\ \ 0 & u^{-1} \end{array}\right)^{\tilde A} 
         \left( \begin{array} {cc} 0 & \pi^{-1} \\ -\pi^{} & 0 \end{array}\right)^{\tilde B} 
         \\ & \ \times
  \left( \begin{array} {cc} 0 & -1 \\ 1 & 0 \end{array}\right)^{\tilde A}
        \left( \begin{array} {cc} u^{-1} &  0 \\ \ 0 & u \end{array}\right)^{\tilde A} 
        \\& = \kappa(a_u)   
     \left( \begin{array} {cc} 0 & 1 \\ -1 & 0 \end{array}\right)^{\tilde A} 
   \left( \begin{array} {cc} 0 & -\pi^{-1} \\ \pi^{} & 0 \end{array}\right)^{\tilde B} 
      \left( \begin{array} {cc} u &  0 \\ \ 0 & u^{-1} \end{array}\right)^{\tilde B}  
         \left( \begin{array} {cc} 0 & \pi^{-1} \\ -\pi^{} & 0 \end{array}\right)^{\tilde B} 
          \\ & \ \times
  \left( \begin{array} {cc} 0 & -1 \\ 1 & 0 \end{array}\right)^{\tilde A}
        \left( \begin{array} {cc} u^{-1} &  0 \\ \ 0 & u \end{array}\right)^{\tilde A}  \\& =  \kappa(a_u)   
     \left( \begin{array} {cc} 0 & 1 \\ -1 & 0 \end{array}\right)^{\tilde A}  
      \left(   \left( \begin{array} {cc} 0 & -\pi^{-1} \\ \pi^{} & 0 \end{array}\right)       \left( \begin{array} {cc} u &  0 \\ \ 0 & u^{-1} \end{array}\right) 
         \left( \begin{array} {cc} 0 & \pi^{-1} \\ -\pi^{} & 0 \end{array}\right)\right)^{\tilde B}
           \\ & \ \times
      \left( \begin{array} {cc} 0 & -1 \\ 1 & 0 \end{array}\right)^{\tilde A}
        \left( \begin{array} {cc} u^{-1} &  0 \\ \ 0 & u \end{array}\right)^{\tilde A}  \end{aligned}
$$}
 
which becomes
{\small
$$\begin{aligned}
         \\ & =
        \kappa(a_u)      \left( \begin{array} {cc} 0 & 1 \\ -1 & 0 \end{array}\right)^{\tilde A}
 \left( \begin{array} {cc} u^{-1} & 0 \\ 0 & u \end{array}\right)^{\tilde B}
   \left( \begin{array} {cc} 0 & -1 \\ 1 & 0 \end{array}\right)^{\tilde A}
        \left( \begin{array} {cc} u^{-1} &  0 \\ \ 0 & u \end{array}\right)^{\tilde A}  \\ & =
        \kappa(a_u)^2  \left( \begin{array} {cc} 0 & 1 \\ -1 & 0 \end{array}\right)^{\tilde A}
 \left( \begin{array} {cc} u^{-1} & 0 \\ 0 & u \end{array}\right)^{\tilde A}
   \left( \begin{array} {cc} 0 & -1 \\ 1 & 0 \end{array}\right)^{\tilde A}
        \left( \begin{array} {cc} u^{-1} &  0 \\ \ 0 & u \end{array}\right)^{\tilde A}  = \kappa(a_u)^2.
\end{aligned} $$
}

 We have now proven that the composite map  
~\eqref{urrghl}, which is also
  the $\PSL_2$-row of~\eqref{pgl2psl2} 
 is the square of the tame symbol.  By the commutativity of~\eqref{pgl2psl2}, we are done.

\subsubsection{The case when \texorpdfstring{$F$}{F} has one archimedean place; 
proof of (ii) in Theorem~\ref{theorem:K2popularversion} }

In this case, we proceed exactly as for 
~\eqref{almostone} (with the ring of integers
replaced by $\OO_F[\frac{1}{\q}]$ -- i.e., when constructing the spectral sequence as in
\S \ref{sec:sss}, one uses just the filtration on $\mathscr{B}_{T-\{\q\}}$) to construct an isomorphism
\begin{equation} \label{placeholder} H_2(Y[\frac{1}{\q}], R)/\bigoplus_{\p \neq \q} 
H_2(\{\p\} [\frac{1}{\q}], R)^{-}  \stackrel{\sim}{\longrightarrow} K_2(\OO_F[\frac{1}{\q}] ) \otimes R. \end{equation}  
Here $Y(\frac{1}{\q})$ is defined as in~\S~\ref{Sarithmetic}
with $\Sigma = \emptyset$ and $\G=\PGL_2$. 

Note, for later use, that the abstract prime-to-$\q$ Hecke algebra $\mathscr{T}_{\q}$
(see~\S~\ref{CompletionConvention})
acts on $H_2(Y[\frac{1}{\q}], R)$ and -- as discussed around ~\eqref{firstline} -- this action 
descends to the left-hand quotient above; the quotient
action is simply that every $T \in \mathscr{T}_{\q}$
acts by its degree $\deg(T)$.  
  
 Now we have constructed an exact sequence (the $-$ eigenspace of Lemma~\ref{lyndon-lemma}, extended a little to the left)
\begin{small}
\begin{equation}\label{DOHA} H_2(Y_0(\q), R)^{-} \rightarrow H_2(Y_0(1), R) \rightarrow H_2(Y[\frac{1}{\q}], R) 
\twoheadrightarrow  \ker(H_1(Y_0(\q), R)^{-} \rightarrow H_1(Y_0(1), R)). \end{equation} 
\end{small}
Compare~\eqref{DOHA} and its analogue replacing level $1$ by level ${\p}$:

\begin{small} 
\begin{equation} 
\begin{diagram} 
\bigoplus_{\p} H_2(Y_0(\p \q), R)^{-} &\rTo&   \bigoplus_{\p} H_2(Y_0({\p}), R)  &\rTo&   \bigoplus_{\p} H_2(Y\{{\p}\}[\frac{1}{\q}], R)  & \rTo&  \\ 
 \dTo && \dTo && \dTo{\delta} &&  \\ 
 H_2(Y_0(\q), R)^{-} &\rTo&  H_2(Y_0(1), R)  &\rTo&  H_2(Y[\frac{1}{\q}], R) & \rTo&
\end{diagram} 
\end{equation} 
\begin{equation}
\begin{diagram}
&\rOnto&   \bigoplus_{\p} \ker(H_1(Y_0(\p \q), R)^{-} \rightarrow H_1(Y_0(\p), R)) \\ 
&& \dTo \\
&\rOnto&  \bigoplus_{\p} \ker(H_1(Y_0(\q), R)^{-} \rightarrow H_1(Y_0(1), R)). 
\end{diagram}
\end{equation} 
\end{small} 
On the top line, we take summations over all primes $\p$ not equal to $\q$. 
Also,  on the top line, the $-$ subscripts refer to the $\q$-Atkin Lehner involution. 
We take the vertical maps to be differences between the two degeneracy maps. (It would therefore
be possible to  replace the top vertical row by the corresponding $-$ eigenspaces for the $\p$-Atkin Lehner involution.)

\medskip

The cokernel of the map $\delta$, as we have seen, is isomorphic to $K_2(\OO_F[\frac{1}{\q}]) \otimes R$. 
Chasing the above diagram, 
we
obtain a sequence, exact at the middle and final term:
\begin{equation} \label{ryar}  \frac{ H_2(Y_0(1), R)  }{ \sum_{\p} \mbox{ image of } H_2(Y_0({\p}), R) }\rightarrow K_2(\OO_F[\frac{1}{\q}] ) \otimes R  \rightarrow \mathcal{K}^*.\end{equation}
where we set
$$ \mathcal{K}^* =  \frac{ \ker(H_1(Y_0(\q), R)^{-} \rightarrow H_1(Y_0(1), R))}{
\sum_{\p} \mbox{ image of }\ker\left(H_1(Y_0(\p \q ),R)^- \twoheadrightarrow H_1(Y_0(\p),R) \right) }.$$

    The last two groups of~\eqref{ryar} admit compatible maps to $k_{\q}^{\times}$, i.e.,
    there's a commutative square
    \begin{equation} \begin{diagram}
K_2(\OO_F[\frac{1}{\q}]) \otimes R & \rTo & \mathcal{K}^* \\ 
    \dTo && \dTo_{\overline{\theta}} \\
    k_{\q}^{\times} & \rEquals &  k_{\q}^{\times}.
    \end{diagram} \end{equation} 
    where $\overline{\theta}$ is the map induced by the map $\theta$ of~\eqref{ra}.  
    In order to see that $\theta$ descends to $\mathcal{K}^*$, we need to see
    that $\theta$ is trivial on the image, in $H_1(Y_0(\q), R)$, of 
    $H_1(Y_0(\p \q), R)$ (via the difference of the two degeneracy maps).  But $\theta$ pulls back via both degeneracy maps
    to the {\em same} map $H_1(Y_0(\p \q), R) \rightarrow k_{\q}^{\times}$;
    thus $\theta$ descends to $\mathcal{K}^*$, as required.

  That means we get in fact a sequence, exact at middle and final terms:
 \begin{equation} \label{ryar2}   \frac{ H_2(Y_0(1), R) }{ \sum_{\p} \mbox{ image of }H_2(Y_0({\p}), R)} \stackrel{}{\longrightarrow }K_2(\OO_F) \otimes R  \twoheadrightarrow \mbox{a quotient of $\mathcal{K}_{\q}$} ,\end{equation}
 writing, as in the statement of the theorem,  $$\mathcal{K}_{\q} =   \mbox{ the kernel of } H_1(Y_0(\q), R)^{-} \longrightarrow H_1(Y_0(1), R) \times \left(  k_{\q}^{\times} \otimes R\right). $$
 
 \medskip
 
 The assertion~\eqref{nexttolast}  of the theorem  may now be deduced: If there are no cohomological
 cusp forms of level $1$, the first group of~\eqref{ryar2} vanishes. We need to refer ahead to Chapter 6 for the proof of that: In the notations of that Chapter,  in the case at hand, 
 the maps $H_2(\partial Y_0(1), R) \rightarrow H_2( Y_0(1), R)$  
 are surjective (as in~\eqref{forearlierref})
  but, as we check in~\S~\ref{topologyprojectiondown}, the quotient of $H_2(\partial Y_0(1), R)$
 by the image of all $H_2(\partial Y_0({\p}), R)$ is zero. (For this it is important
 that $w_F$ is invertible in $R$). 
  
 \medskip
 
 Now, as in the very final assertion of Theorem~\ref{theorem:K2popularversion} , let us localize at $p$.
Put $R= \Z_p$.  As before let $\mathscr{T}_{\q}$ be the abstract prime-to-$\q$ Hecke algebra. 
Then in fact $\mathscr{T}_{\q}$ acts on all three groups of~\eqref{ryar} where,
according to our discussion after~\eqref{placeholder}, each $T \in \mathscr{T}_{\q}$ 
acts on the middle group by its degree. 
\medskip

Write  $\m = (\eisideal, p) \subset \mathscr{T}_{\q}$, i.e.
it is the kernel of the morphism $\mathscr{T}_{\q} \rightarrow \mathbf{F}_p$
that sends every Hecke operator to its degree (mod $p$). 
Consider, now, the sequence~\eqref{ryar} completed at $\m$. 
\begin{itemize}
\item[-] The first group becomes zero: 
Our assumption (as in the last paragraph of   Theorem~\ref{theorem:K2popularversion}) implies, in particular, that   $H_1(Y_0(1),R)_{\m}$
is finite; that means again that the map $H_2(\partial Y_0(1), R)_{\m} 
\rightarrow H_2(Y_0(1), R)_{\m}$ is surjective, and we can proceed just as before.
\item[-] The middle group equals $K_2(\OO_F) \otimes \Z_p$: 
We already saw that each $T \in \mathscr{T}_{\q}$
acts on the middle group by its degree; thus completing at the ideal $\m$
is the same as completing at $p$. 
\item[-] The third group remains a quotient of $\mathcal{K}_{\q}$:
The group $\mathscr{K}_{\q}$ is finite, and the completion of a finite module
 is isomorphic to a quotient of that module.
\end{itemize}

This concludes the proof of Theorem~\ref{theorem:K2popularversion}.

 \chapter{Analytic torsion and regulators}
 \label{chapter:ch4}
\label{chapter:regulator}

In this Chapter, we review the Cheeger-M{\"u}ller theorem (\S~\ref{sec:CMthm}), 
at least in the case of compact $3$-manifolds.  
It relates the torsion, a Laplacian determinant, and a third quantity: a ``regulator.''
The issue of generalizing~\eqref{mst} to the non-compact case is taken up (and solved partially --- enough for our purposes) --- in
 Chapter~\ref{chapter:ch6}.

The main goal of the rest of the chapter is to better understand the regulator:
\begin{enumerate}
\item We relate  the regulator to central values of $L$-functions (Theorem
\ref{regulatorL}), up to a rational ambiguity,  and discuss potential relationship with
the Faltings height (\S~\ref{subsec:FH}).   We give the proof
in~\S~\ref{regLproof}. Most of the content here is due to Waldspurger. 

\item We explain how the regulator changes under the operation of replacing $Y(K)$
with a cover (Theorem~\ref{regcompare}), at least in the simple situation when there are no newforms;
\item We compare (with less success)
the regulator between two manifolds in a Jacquet--Langlands pair (Lemma~\ref{veryconditional});
we believe that that this is related to level lowering congruences, and give
some conditional evidence for this.

\item There is one further section that perhaps belongs in this chapter:
The methods of~\S~\ref{sec:eisintegrality} 
can be used to compute the Eisenstein part of the regulator in the {\em noncompact} case;
  we postpone this to Chapter 6 so that the necessary definitions in the noncompact
  case are given first. 
\end{enumerate}

\section{The Cheeger-M{\"u}ller theorem}  \label{sec:CMthm}
Here is what the Cheeger--M{\"u}ller theorem  \cite{Cheeger, MEta} says about $Y(K)$ when it is compact (or, indeed,
about any compact hyperbolic $3$-manifold). We assume for the moment that $Y(K)$ has no orbifold points, the modifications in the orbifold case being discussed in~\S~\ref{orbifoldtorsion1}.

\index{regulator (compact case)}
\index{$\reg(H_{\ell})$}
Define the {\em regulator of the $j$th homology} via
\begin{equation} 
\label{reghnsimpledef} \reg(H_\ell(Y(K))) =     \left |\det \int_{\gamma_j} \nu_k \right| \end{equation} 
   where $\gamma_i$ give a basis for $H_\ell(Y(K), \Z) \mbox{ modulo torsion}$, and $\nu_k$ an orthonormal basis of harmonic forms. 
   
   Another way to say this is as follows: the image of $H_{\ell}(Y(K), \Z)$
   inside $H_{\ell}(Y(K), \R)$ is a {\em lattice} in the latter vector space;
   and $\reg(H_{\ell})$ is the {\em covolume} of this lattice, with respect to the measure
   on $H_{\ell}(\R)$ that arises with its identification as the dual to a space of harmonic forms. 
   A more intrinsic way to write the right-hand side would be $  \frac{ \left |\det \int_{\gamma_j} \nu_k \right|} { \left| \det \langle \nu_j, \nu_k \rangle  \right|^{1/2}}$;  in this form it is independent of basis $\nu_j$.

Now we define the {\em regulator of $Y(K)$} via: 
\begin{equation} 
\label{regsimpledef} \reg(Y(K)) =
\frac{\reg(H_1) \reg(H_{3})}{\reg(H_0) \reg(H_{2})}=   \vol  \cdot   \reg(H_1(Y(K))^2, \end{equation} 
where 
$\vol$  is the product of the volumes of all connected components.  The equality
here follows from the easily verified facts that $\reg(H_i) \reg(H_{3-i})= 1$
and $\reg(H_0) = 1/\sqrt{\vol}$.

Denote by $\Delta^{(j)}$ the Laplacian 
on the orthogonal complement of harmonic $j$-forms inside smooth $j$-forms.
Let $H_{1,\tors}$ be the torsion part of $H_1(Y(K), \Z)$. We then define {\em Reidemeister} and {\em analytic} torsion via the formulae: 

\begin{eqnarray} 
\label{rtsimpledef} \RT(Y(K)) &=&  |H_{1,\tors} |^{-1} \cdot \reg(Y(K)) \\
\label{analTsimple} \analT(Y(K)) &=&  \det(\Delta^{(1)})^{-1/2} \det(\Delta^{(0)})^{3/2}. \end{eqnarray}
The determinants appearing in~\eqref{analTsimple} are to be understood via zeta-regularization. Then the {\em equality} of Reidemeister torsion and analytic torsion, 
conjectured by Ray--Singer and proved by Cheeger and M{\"u}ller, asserts simply
\begin{equation} \label{mst} \RT (Y(K)) = \analT(Y(K)). \end{equation}
{\em Intuitively} this formula expresses the fact that the the limit of a chain complex for $M$,
as one triangulates very finely, approaches the de Rham complex.   

\begin{remarkable}
When $Y(K)$ is noncompact, 
the exact analogue of~\eqref{mst} is not known. However, we prove
enough of such an analogue for our purposes in the next Chapter. 
\end{remarkable}

\subsection{Equivariant case/ orbifold case.} \label{orbifoldtorsion1}

  The theorem of Cheeger--M{\"u}ller is for Riemannian manifolds; but we are naturally dealing with {\em orbifolds}.  In short, it continues to hold, but possibly with errors divisible by prime numbers that divide the order of isotropy groups, i.e.
\begin{equation} \label{CMorbifold} \RT (Y(K))= \tau_{an} (Y(K)). u, \ \ u\in \mathbf{Q}^{\times},\end{equation}
where
  the numerator and denominator of $u$ are {\em supported at orbifold primes}.
  By this we mean that $u = a/b$, where $a$ and $b$ are integers divisible
  only by primes  dividing the order of the isotropy group of some point on $Y(K)$.
The  definitions of $\RT, \tau_{an}$ are as previous. (Note that, in the orbifold case, infinitely many homology groups $H_i(Y(K), \Z)$ can be nonzero, but this does not affect our definition, which used only $H_1(Y(K), \Z)_{\tors}$.)

\medskip

The relation~\eqref{CMorbifold} is presumably valid quite generally, but we explain the proof only in the case of interest, where the orbifolds can be expressed as global quotients of manifolds: 
In our situation, where $M = Y(K)$, 
we may choose a sufficiently small subgroup $K'$, normal in $K$, such that $Y(K')$
is a genuine manifold --- that is to say, every conjugate of $K'$ intersects
$\G(F)$ in a torsion-free group. Let $\Delta =K/K'$. 
Then $Y(K)$ is the quotient $Y(K')/\Delta$ in the sense of orbifolds. 

\medskip

Suppose, more generally, that $\tilde{M}$ is an odd-dimensional compact Riemannian manifold with isometric $\Delta$-action,
and let $M=\tilde{M}/\Delta$ (with the induced metric, so that
$\tilde{M} \rightarrow M$ is an isometry).       For any $x \in M$, let
$\Delta_x$ be the (conjugacy class of the) isotropy group 
of a preimage $\tilde{x} \in \tilde{M}$.  The orbifold primes for $M$ are those which divides the order of some 
$\Delta_x$.   

Fix a cell decomposition of $\tilde{M}$ and let $C^*(\tilde{M})$
be the corresponding integral chain complex. Then,
if $\ell$ is not an orbifold prime,  the cohomology groups of
$$C^*(\tilde{M} \otimes\Z_{\ell})^G$$
agree with the orbifold cohomology groups of $M$.\footnote{ Indeed, the orbifold cohomology can be defined as the hypercohomology of $C^*(\tilde{M} \otimes \Z_{\ell}) $
(in the category of complexes of $\Z_{\ell}[G]$-modules) and there is a corresponding spectral sequence, converging to the orbifold cohomology of $M$, 
whose $E_1$ term is $H^i( G, C^j \otimes \Z_{\ell})$.  But for $i > 0$, the group
$H^i(G,C^j(\tilde{M}) \otimes \Z_{\ell})$ vanishes if $\ell$ is not a orbifold prime, by Shapiro's lemma.}
On the other hand, it is shown by Lott and Rothenberg \cite{LottRoth} that the Whitehead torsion of the complex
$$C^*(\tilde{M},\Z)^G$$
coincides with the analytic torsion of $M$,
and our claim~\eqref{CMorbifold} follows easily. 

\medskip

For later applications we will also discuss later (
see~\S~\ref{orbifoldtorsion2})
 the case when $M$  (thus $\tilde{M}$) has a boundary, which has 
been treated by L{\"u}ck \cite{Lueck}.

\section{(Cuspidal) regulators are arithmetic periods}  \label{regcusp}
 
It will follow from the therem of Cheeger and M{\"u}ller that, if $Y,Y'$ are a Jacquet--Langlands pair
(see~\S~\ref{ss:jlc}) then \begin{equation} \label{rrQ} \frac{\reg^{\mathrm{new}}(Y)}{\reg^{\mathrm{new}}(Y')} \in \Q^{\times},\end{equation}
where {\em new} regulators are defined as a certain alternating quotient of regulators.
(See~\S~\ref{sec:newnewnew} and~\S~\ref{sec:coh} for the definition and the deduction of~\eqref{rrQ}, respectively.) 
 In this section, we shall give another proof of~\eqref{rrQ},
{\em independent of the Cheeger--M{\"u}ller theorem:} 
\medskip
We shall show, slightly more precisely, 
that both $\reg(Y)$ and $\reg(Y')$ 
are (up to $\Q^{\times}$) equal to a certain ratio of special values of $L$-functions.

\medskip

In principle, the advantage of  such an independent proof of~\eqref{rrQ}
is that it could lead to a precise understanding of the ratio, and, therefore,
of the effect on torsion. We have nothing unconditional in this direction, but we make a (conjectural) start in this direction in \S
\ref{reglowlevel}.  The examples considered there were motivated by 
numerics described in~\S~\ref{subsec:level lowering}.
\medskip

For our conventions regarding automorphic representations see~\S~\ref{autrepnotns}.

 \subsection{Periods} 
 Let $\pi$ be a cohomological automorphic representation for $\G$. 
 Let $\Q(\pi)$ be the field
generated by all (prime to level) Hecke eigenvalues.  
 We shall attach to $\pi$ a certain ``complex period'' $\reg_{\C}(\pi) \in \C^{\times}/\Q(\pi)^{\times}$. 

\medskip

To define $\reg_{\C}(\pi)$ we use special values of $L$-functions: 

\medskip

We say a pair of quadratic characters $(\chi, \chi')$ of the idele class group of $F$
is {\em admissible} if $\chi \chi'$ is  nontrivial at all real places of $F$. 
 We denote by $\Delta_{\chi}$ the discriminant of the quadratic extension attached to $\chi$ and write $\mathcal{L}_{\pi}(\chi) :=  \Delta_{\chi}^{1/2}  L(\frac{1}{2}, \pi \times \chi)  .$
There exists  a
period $\reg_{\C}(\pi) \in \C^{\times}/\Q(\pi)^{\times}$  such that for every admissible pair $(\chi, \chi')$, we have 
$$\mathcal{L}_{\pi}(\chi) \mathcal{L}_{\pi}(\chi') \in  \Q(\pi)   \cdot  \reg_{\C}(\pi).$$
The existence of such a period follows from the sharper enunciation
on page 174 of \cite{waldsur}.

More generally suppose that $O$ is a set of cohomological automorphic representations $\pi$ that is stable under the action of the Galois group\footnote{Given such $\pi$
and $\sigma \in \mathrm{Gal}(\overline{\Q}/\Q)$, there exists a unique cohomological automorphic representation $\pi^{\sigma}$ whose Hecke eigenvalues are obtained from those of $\pi$ by applying $\sigma$.}.
It is then possible to choose a lifting $\reg_{\C}(O) \in \C^{\times}/\Q^{\times}$
of  $\prod_{\pi \in O} \reg_{\C}(\pi)$ so that, for any admissible pair $(\chi, \chi')$, 
so \begin{equation} \label{galoisorbit} \prod_{\pi \in O}  \mathcal{L}_{\pi}(\chi) \mathcal{L}_{\pi}(\chi')  \in \Q  \cdot \reg_{\C}(O).\end{equation}
This assertion is not explicitly in Waldspurger's paper, but can be derived by the techniques there. If one does not assume it, the result below changes only in that one needs to replace $\mathbf{Q}^{\times}$ by units in the field generated by all Hecke eigenvalues. 

For example, suppose that $\pi$ is associated to an elliptic curve $E$
over $F$; therefore, $\Q(\pi) =\Q$. Then --- if we suppose the validity of the BSD
 conjecture --- $\reg_{\C}(\pi)$ will be the coset $\Q^{\times} \Omega_E$,
where 
\begin{equation} \label{Falt2} \Omega_E = \prod_{\sigma: F \hookrightarrow \C} \int_{E^{\sigma}(\C)} \left| \omega^{\sigma} \wedge \overline{\omega^{\sigma}}\right|, \end{equation}
where $\omega$ is any $F$-rational differential on $E$.  
More generally, if $O$ is the Galois orbit of forms associated to an abelian
variety $A$ of $\GL_2$-type, then $\reg_{\C}(O)$ is again given by a similar formula   with $E$ replaced by $A$,
and with $\omega$ replaced by a basis $\omega_1, \dots, \omega_g$:  
\begin{equation} \label{Falt3} \Omega_A = \prod_{\sigma: F \hookrightarrow \C} \int_{E^{\sigma}(\C)} \left|  \omega_1^{\sigma} \wedge 
\cdots \wedge \omega_g^{\sigma} \wedge \overline{\omega_1^{\sigma}} \wedge \cdots \wedge \overline{\omega_g^{\sigma}}\right|, \end{equation}

\subsection{} 
 Here is the theorem relating regulator and $L$-values.
 \medskip
 
  For $\pi$ a cohomological automorphic
representation for $\G$ write $m(\pi) = \dim \pi_f^{K}$ for the multiplicity of $\pi$ at level $K$.  For $O$ any Galois orbit on cohomological representations, let 
$m(O)$ be the constant value of $m(\pi)$ on $O$. 

\begin{theorem} \label{regulatorL} \label{REGULATORL}
Suppose $D$ is not split. Let $K \subset \G(\Afinite)$ be an open compact, which
equals the image of the units in a maximal order 
at all ramified places for $D$ (i.e., at ramified places, the level structure is the $K_{0,v}$
from \S \ref{ls}). Then 
 \begin{equation} \label{regrat} \reg(H_1(Y(K))^2  \in  c' \Q^{\times} \prod_O \left( 
 \frac{ \reg_{\C}(O)}{ \prod_{\pi \in O}  L(1, \Ad, \pi)} \right)^{m(O)} \end{equation}
where $c' = \left( \pi^{2r_1+1}\cdot \sqrt{\disc_F}\right)^{\dim H_1(Y(K), \C)}$, where $\disc_F$ is the discriminant of $F$,
and $\reg_{\C}(O)$ was defined around \eqref{galoisorbit}.
\end{theorem}

 \begin{remark}
 If $D$ is split, the proof goes through without any change  after replacing $\reg(H_1)$ by $\reg(H_{1!})$ 
(see~\S~\ref{subsec:rtnonsplitdef} for definition), 
taking the product only over cuspidal cohomological representations $\pi$,
and requiring that each $\pi$ 
is Steinberg, up to twist,   at at least one place.   The right-hand side remains unchanged. 
Similar analysis leads to~\eqref{rrQ}   even if one of $Y$ and $Y'$ is split. 
\end{remark}

As we have mentioned, all the key
ideas are already in Waldspurger \cite{waldsur, waldQ}. Nonetheless
we include some detail because we think it is important to investigate
 the precise properties of regulators more carefully, even though this is only a crude first step.

\subsection{Integrality questions and the Faltings height} \label{subsec:FH}

A fundamental question is: 
\begin{quote} Can one pin down $\reg(Y)$ {\em integrally} as opposed to rationally? \end{quote}

It is important to note that our argument to prove Theorem~\ref{regulatorL} gives almost {\em no integral control} on $\reg(Y)$: 
we exhibit cycles in $H_1$ but have {\em no} way of controlling their divisibility.

In some situations it maybe possible to say more:
\begin{enumerate}
\item 
When all forms on $Y$ are base-change from $\Q$. In that case, one can exhibit explicit cycles in $H_2(Y)$ which are Poincar{\'e} dual to the form;  \label{CMform}
this allows one to control the divisibility of cycles in $H_1$.
\item The other case is in the split case, where one can use modular symbols. See
e.g.   our~\S~\ref{sec:eisintegrality}, or work of Sczech~\cite{Sczechelliptic} or Berger~\cite{Berger}.
\end{enumerate}

Suppose, for the sake of our discussion,
that there is a unique orbit $O$ as in Theorem~\ref{regulatorL} with $m(O) = 1$, 
and it is associated to an abelian variety $A$ of $\GL_2$-type over $F$.

It is perhaps reasonable to postulate that~\eqref{regrat}
remains {\em integrally true} if we replace $\reg_{\C}(O)$
by the period $\Omega_A$ defined by ~~\eqref{Falt2}and ~~\eqref{Falt3},
choosing for the $\omega_i$ an integrally normalized basis of differential forms, and 
{\em at least after replacing $A$ by an isogenous $A'$.}
Some evidence for this is the numerical computations of Cremona \cite{Cremona}.
Note that $\Omega_A$ defined thus is (the exponential of the negative of the)
``Faltings height'' of the abelian variety $A$; 
  the work of Faltings shows that $\Omega_A$
and $\Omega_{A'}$ thus defined can only differ at ``small primes,''
so one could perhaps dispense with $A'$ at the cost
of inverting a few explicit, small primes.

Such a relationship would be a striking
relationship of {\em geometric complexity} and {\em arithmetic complexity}: 
it relates $\reg(H_2) = \reg(H_1)^{-1}$ (which is a measure of the geometric complexity
of $2$-cycles representing homology) to the Faltings height of $A$
(which measures the ``arithmetic complexity'' of the variety $A$).

\section{Proof of Theorem~\ref{regulatorL} } \label{regLproof} 
\subsection{} \label{param} 
To prove Theorem~\ref{regulatorL}, we start by parameterizing harmonic forms in the language of automorphic forms.

For any cohomological automorphic representation $\pi$,  denote by $\Omega$ the natural map
\begin{equation} \label{omegadef} \Omega: \Hom_{K_{\infty}}(\mathfrak{g}/\mathfrak{k}, \pi) \rightarrow \lim_{\longrightarrow}  \Omega^1(Y(K))\end{equation}
Indeed, $\Omega^1(Y(K))$ can be considered as functions on
$\G(F) \backslash \left(  \G(\adele) \times \mathfrak{g}/\mathfrak{k} \right) / K_{\infty} K$ that
are linear on each $\mathfrak{g}/\mathfrak{k}$-fiber.

Explicitly, for $X \in \mathfrak{g}/\mathfrak{k}$ and $g \in \G(\adele)$, 
we can regard $(g, X)$ as a tangent vector to $\G(F) g K_{\infty} K \in Y(K)$:
 namely, regard $X$ as a tangent vector at the identity by differentiating
 the $G$-action, and then translate this by $g$. 
 
Denote this vector by $[g, X]$ -- 
and the map $\Omega$ is normalized by the requirement
that, for $f \in  \Hom_{K_{\infty}}(\mathfrak{g}/\mathfrak{k}, \pi) $, we have
\begin{equation} \label{OmegaVeryExplicit}
\Omega(f) ( [g, X])  = f(X) (g).
\end{equation}

 Equipping $\mathfrak{g}/\mathfrak{k}$ with the metric corresponding to the Riemannian structure (see~\S~\ref{H3} for normalizations) 
we obtain a unitary structure on the former space, by requiring  the map $\Omega$ isometric. 
 The map $\Omega$ defines an isomorphism onto a subspace $\mathcal{H}_{\pi} := \varinjlim_K H^1(Y(K), \C)_{\pi}$
of $\mathcal{H} := \varinjlim_{K} H^1(Y(K), \C)$.

\subsection{} \label{sec:data}
From quadratic extensions of $F$ one obtains tori in $\mathbf{G}$
and thereby loops on $Y$; we pin down some details of this.

Recall that we write   $F_{\infty} = F \otimes \R$.  
 Fix an isomorphism $$ \G(F_{\infty}) \stackrel{\sim}{\longrightarrow} \PGL_2(\C) \times \SO_3^r$$
 such that $K_{\infty}$ is carried to $\mathrm{PU}_2 \times \SO_3^r$,
 where $\mathrm{PU}_2$ is the image in $\PGL_2(\C)$ of the  stabilizer of the standard Hermitian form $|z_1|^2+|z_2|^2$.

 Let  $H$ be the subgroup $\left( \mbox{diagonal in $\PGL_2(\C)$} \times \SO_2^r\right)$
 (for an arbitrary choice of $\SO_2$s inside $\SO_3$). There is a natural isomorphism
\begin{equation} \label{rsiso} \bar{H} := H / H \cap K_{\infty} \stackrel{\sim}{\longrightarrow} \mathbf{R}_+ \end{equation}
descending from the 
map $ \left( \begin{array}{cc} z & 0 \\ 0 & 1 \end{array} \right)  \mapsto |z| $.
The isomorphism~~\eqref{rsiso} normalizes an element $X$
in the Lie algebra of $\bar{H}$, namely, the element corresponding to the one-parameter subgroup $t  \in \mathbf{R} \mapsto \exp(t) \in \mathbf{R}_+$. 
Later, we shall regard $X$ as an element of $\mathfrak{g}/\mathfrak{k}$.

We consider data $\mathscr{D}$ consisting of:
\begin{itemize}
\item a maximal subfield $E \subset D$. 

This defines a torus $\T \subset \G$ as the centralizer of $E$ (therefore, $\T(F) = E^{\times}/F^{\times}$); 
\item an element $g \in \G(\adele)$ such that\footnote{We denote by $g_{\infty}$
the image of $g$ under $\G(\adele) \rightarrow \G(F_{\infty})$.} $g_{\infty}^{-1} \T(F_{\infty}) g_{\infty}  = H$;
\item A quadratic character of $\adele_E^{\times}/E^{\times} \adele_F^{\times}$, fixed under the Galois automorphism of $E$ over $\Q$.

This 
induces a character  $\psi:  \T(\adele)/\T(F) =
\adele_E^{\times}/E^{\times} \adele_F^{\times}  \rightarrow \valuefield^{\times}$.
 \end{itemize}

The isomorphism class of the pair $(E, \psi)$ determines, and is determined by,
a pair of quadratic characters $\{\chi_1,\chi_2\}$ of $\adele_F^{\times}/F^{\times}$:
$\psi$ factors through the norm map and therefore there are two characters 
$\chi_1, \chi_2$ such that $\psi = \mathrm{N} \circ \chi_i$; the product $\chi_1 \chi_2$
is the quadratic character $\chi_{E/F}$ attached to $E$ by class field theory.  We often write in what follows $E_v := (E \otimes F_v)^{\times}$.

 For each such datum $\mathscr{D}$, define an element $ \gamma_{\mathscr{D}} \in H_1(Y, \Z)$ as follows: 
 
 Put $K_{\data} = g K_{\infty} K  g^{-1}$.
The quotient $Y_T := \T(F) \backslash \T(\adele) /  K_{\data} \cap \ker(\psi)$ is a compact $1$-manifold (possibly disconnected). 
The map  \begin{equation} \label{ttoY} \ \ \ \ \emb: t \mapsto tg, \ \   \T(F) \backslash \T(\adele)  \rightarrow \G(F)\backslash \G(\adele) \end{equation}  descends to a map
$\emb: Y_T \rightarrow Y(K)$; its image is 
a finite collection of closed  arcs $\{G_i\}_{1 \leq i \leq I}$ on $Y$.
Let $[G_i] \in H_1(Y(K), \Z)$ be the homology class of $G_i$. 
Moreover, $\psi$ factors through to a function from the finite set of $G_i$ to $\pm 1$; 
  we may speak
of $\psi(G_i)$ as the constant value taken by $\psi$ on $G_i$. 
Accordingly, 
put $$\gamma_{\data} = \sum_{1}^{I} \psi(G_i) . [G_i] \in H_1(Y(K), \Z).$$
Note that this will actually be trivial if $K_{\data}$ is not contained in the kernel of $\psi$.

\begin{theorem}\label{thm:rationality}
Let $f \in H^1(Y(K), \overline{\Q})$ be a Hecke eigenform whose Hecke eigenvalues are those of the automorphic representation $\pi$,
and let $\omega_f$ be a harmonic $1$-form representing $f$.  Let $\data$
be a datum as above, with associated class $\gamma_{\data} \in H_1(Y(K), \Q)$
and associated  idele class characters $\chi_1, \chi_2$ of $F$. 

Set
$$S (f) =  \frac{ \left| \int_{\gamma_{\data}} f \right|^2} {\langle \omega_f, \omega_f \rangle_{L^2(Y(K))}} 
  \left(c   \cdot \frac{   \Lambda(1/2, \pi \times \chi_1) \Lambda(1/2, \pi \times \chi_2)  }{ \Lambda(1, \Ad, \pi) } \right)^{-1}$$
where $c = \pi^{ r_1} \sqrt{\disc_F \disc_E}$,  the inner product $\langle \omega_f, \omega_f \rangle$
  is taken with respect to hyperbolic measure,  
    $\Lambda(\dots)$ denotes completed $L$-function, including $\Gamma$-factors.

Then $S(f) \in \overline{\Q}$; moreover, 
$S(f^{\sigma}) = S(f)^{\sigma}$ for $\sigma \in  \Gal(\Qbar/\Q)$.  
Finally, 
we may choose data $\data$ such that $S(f) \neq 0$ (equivalently, $\int_{\gamma_{\data}} f \neq 0$).

\end{theorem}

 Note that this assertion implies~\eqref{galoisorbit}; we'll explicate this later.   

\begin{proof} 
The proof of this theorem is essentially all in  Waldspurger \cite{waldQ, waldsur}, but we shall explicate a few points.   In fact, similar computations have been carried out in a number of contexts (for example, see \cite{Popa})  although usually  the emphasis is on the $L$-function, and thus there
is less  flexibility in the choice of vector $f$.
Our presentation is quite  terse; many points which are briefly treated here  can be found in greater detail in these papers. 

 We shall actually analyze only the case when ``$f$ arises from a pure tensor in the underlying automorphic representation $\pi$'', which is all that is needed for our application; but the general case is not significantly different.

Let $Z$ be the totally real number field generated by the Hecke eigenvalues of $f$,
regarded as a subfield of $\C$. 
This implies that every representation $\pi_v$ for $v$-finite admits a $Z$-structure, 
that is to say, a $Z$-subspace $\pi_v^Z \subset \pi_v$
with the property that the map $\pi_v^Z \otimes_Z \C \rightarrow \pi_v$ is an isomorphism. 
\footnote{More precisely, $\Hom_{K_{\infty}}(\mathfrak{g}/\mathfrak{k}, \pi) =
 \Hom_{K_{\infty}}(\mathfrak{g}/\mathfrak{k}, \pi_{\infty}) \otimes \pi_f  $
 has a $Z$-structure by virtue of $\Omega$.  Thus $\pi_f = \otimes_v \pi_v$, 
 the product taken over finite $v$, has a $Z$-structure $\pi_f^Z \subset \pi_f$. But this easily implies that every factor of the tensor product $\pi_v$ has a $Z$-structure (for example, take fixed vectors under a suitable compact subgroup to construct it). } 
\label{Zstructuredef}

  Later we shall use the following fact about $Z$-structures:
  \begin{itemize}
\item  If $\Pi$ is a representation of $\PGL(2, K)$, with $K$ a nonarchimedean local field,
  admitting a $Z$-structure, and $\sigma: Z \hookrightarrow \C$ an embedding, then $L(1, \Ad, \pi) \in Z$
  and $L(1, \Ad, \Pi^{\sigma}) = L(1, \Ad, \Pi)^{\sigma}$; 
  the same holds for $L(\frac{1}{2}, \Pi)$.      Here $\Pi^{\sigma} = \Pi^{Z} \otimes_{Z,\sigma} \C$. 
\footnote{For example, consider the spherical case: if $\alpha, \beta$ are the Satake parameter of $\Pi$, normalized to be of absolute value $1$ when $\Pi$ is tempered,  then $\alpha, \beta$ need not belong to $Z$, but $\frac{\alpha + \beta}{\sqrt{q}}$ and $\alpha \beta/q$ {\em do} belong to $Z$; 
 thus $L(\frac{1}{2}, \Pi) = (1 - (\alpha + \beta)/\sqrt{q} + \alpha\beta/q)^{-1}$ also belongs to $Z$. 

} 
\end{itemize}

 {\em In this section,} we follow measures normalizations as in  p. 175 of \cite{waldsur}: In particular, this means that we have endowed
  $\G(F) \backslash \G(\adele)$ with Tamagawa measure,
  and $\T(F) \backslash \T(\adele)$ has measure $\mu_T := \Lambda(1, \chi_{E/F})^{-1} \mu_{\Tamagawa}$.    
  We also follow the normalizations of \cite{wald} for measures on $\G(F_v), \T(F_v)$.

 The measure on $\T(F_{\infty})$, identified with $H$ via conjugation by $g_{\infty}$,  pushes down to a measure on $\bar{H}$; this measure on $\bar{H} \cong \mathbb{R}_+$
is computed to be $4^{[F:\Q]} \frac{dx}{x}$, where $dx$ is the usual Lebesgue measure.   The action of $\T(F_{\infty})$
on $Y_T$ factors through $\bar{H}$ -- again, we identify $\T(F_{\infty})$ with $H$ via conjugation by $g_{\infty}$ -- so $Y_T$ carries a natural vector field
$\varchi_0$ obtained from $X \in \mathrm{Lie}(\bar{H})$ by differentiating the $\bar{H}$-action. 
 Let $\omega_0$ be the differential form on $Y_T$ such that $\langle \omega_0, \varchi_0 \rangle = 1$ everywhere and let $|\omega_0|$ be the corresponding measure.  Then\footnote{To see this first note in a similar way $\bar{H}$ carries a invariant differential form $\omega_0^H$ such that 
$\langle \omega_0^H, X \rangle = 1$;  the measure $|\omega_0^H|$
is identified with $\frac{dx}{x}$ under the isomorphism with $\mathbf{R}_+$. Now the measure on $\T(F_{\infty})$
obtained by pullback via $\T(F_{\infty}) \stackrel{\sim}{\rightarrow} H \rightarrow \bar{H}$
is $4^{-[F:\Q]}  \times \left( \mbox{ measure of \cite{waldsur}}\right)$.   On the other hand, the measure
of \cite{waldsur} assigns to the maximal compact subgroup of $\T(\adele_f)$ the measure
$d_{E/F}^{-1/2} 2^?$, where $d_{E/F} = \frac{\disc_E}{\disc_F^2}$, and $?$ is related
to the number of ramified primes in $E/F$. 
 Therefore, the measure of \cite{waldsur} assigns to each connected component $G\subset Y_T$
 the mass $2^{?} \iota^{-1} \cdot d_{E/F}^{-1/2} \cdot \int_{G} |\omega_0|$
 where $\iota$ is 
   the index
of $g K g^{-1} \cap \T(\adele_f)$ inside the maximal compact of $\T(\adele_f)$.
}
 $$\mu_T =  2^{?} \iota^{-1} d_{E/F}^{-1/2} |\omega_0|,$$
for some natural number $?$ and $d_{E/F} = \frac{\disc_E}{\disc_F^2}$, where $\disc_E, \disc_F$
are discriminants of $E, F$ respectively; 
 $\iota$ is 
   the index
of $g K g^{-1} \cap \T(\adele_f)$ inside the maximal compact of $\T(\adele_f)$.

 Take $\varF \in  \Hom_{K_{\infty}}(\mathfrak{g}/\mathfrak{k}, \pi) $, so that
$\Omega(\varF)$ is a $1$-form on $Y(K)$.   

As before let $\emb: Y_T \rightarrow Y(K)$ be the map descending
from $t \mapsto tg$; then
$$\langle \varchi_0,  \emb^* \Omega(\varF) \rangle(t) = \varF(X) (tg),$$
i.e. $\varF(X) \in \pi$ evaluated at $t g \in \G(\adele)$. 
Actually, the left-hand side is simply $\Omega(\varF)$
evaluated at the tangent vector obtained by applying $\emb$
to a tangent vector on $Y_T$ based at $t$ with derivative $X$.
This tangent vector is just $[tg , X]$ in the notation of~\eqref{OmegaVeryExplicit}.

Moreover, $\varF(X) \in \pi$ is a $K_{\infty} \cap H$-fixed vector
and its $K_{\infty}$-type is minimal (i.e., $K_{\infty}. \varF(X)$ is three dimensional);
this follows since $\mathfrak{g}/\mathfrak{k}$ is $3$-dimensional and $X$
is $K_{\infty} \cap H$-fixed. 
It follows from this that 
\begin{eqnarray} \nonumber \int_{\gamma_{\data}} \Omega(\varF)  &=&   \int \langle \varchi_0, \emb^* \Omega(f) \rangle   \psi(t) \cdot  |\omega_0| 
\\  \nonumber&=&
\int_{Y_T} \varF(X)(tg)   \psi(t) \cdot |\omega_0|\\ &=& \label{cow2}
2^{?}\disc_{E/F}^{1/2}    \iota    \int_{t \in T(F) \backslash \T(\adele)}
\varF(X)(tg)  \psi(t) d\mu_T. \end{eqnarray}

Note for later use \footnote{Indeed, 
the norm of the differential form $\Omega(\varF)$
at the point of $Y(K)$ given by $\G(F) g K_{\infty} K$ is precisely
the norm of $\mathrm{ev}_g \circ  \varF \in \Hom(\mathfrak{g}/\mathfrak{k}, \C)$;
here $\mathrm{ev}_g$ denotes ``evaluation at $g$,'' a functional on $\pi$. 
Extend $X \in \mathfrak{g}/\mathfrak{k}$ to an orthonormal basis $X,Y,Z$.
Then the latter norm is $|\mathrm{ev}_g \varF(X)|^2 + |\mathrm{ev}_g \varF(Y)|^2
+ |\mathrm{ev}_g \varF(Z)|^2$.  Integrating we obtain $\|\varF(X)\|^2+ \|\varF(Y)\|^2  + \|\varF(Z)\|^2 
= 3 \|\varF(X)\|^2$; indeed, $\|\varF(X)\| = \|\varF(Y)\|=\|\varF(Z)\|$ because the map 
$\mathfrak{g}/\mathfrak{k} \rightarrow \pi$ must preserve metrics up to a scalar. The factor $\frac{vol}{2}$ arises from transition from Tamagawa to hyperbolic measure.}
$$ \langle \Omega(\varF), \Omega(\varF) \rangle_{Y(K),\hyp} = \frac{3 \  V_{Y(K)}}{2} \| \varF(X)\|^2 ,$$
where the subscript ``hyp'' exists to emphasize the fact that we take the left-hand inner product 
with reference to the hyperbolic metric, whereas on the right-hand side
$\|\varF(X)\|^2$ is defined  with respect  to the Tamagawa measure on $\G$; $V$ is
the volume of $Y(K)$ computed
with respect to the hyperbolic measure (we don't use the notation
$\vol$ to avoid confusion with the product of volumes of connected components).
Therefore, 

\begin{equation} \label{omegaf} \frac{  \left| \int_{\gamma_{\data}} \Omega(\varF)\right|^2 }{\langle \Omega(\varF), \Omega(\varF) \rangle_{\hyp}} \sim \frac{1}{V_{Y(K)}} \frac{\left| \int_{T(F) \backslash T(\adele)} f(tg) \psi(t) d\mu \right|^2 }{\langle f, f \rangle},\end{equation}
where $f = \varF(X) \in \pi$ and $\sim$ denotes equality up to $\Q^{\times}$.

Factorize $\pi = \otimes_v \pi_v$ as unitary representation. Then, for $f = \otimes f_v \in \pi,$   \cite[Prop. 7, page 222]{waldsur} gives  
\begin{equation} \label{WaldF} \frac{\left| \int_{T(F) \backslash T(\adele)} f(t) \psi(t) dt \right|^2 }{\langle f, f \rangle} = \frac{\Lambda(1, \chi_{E/F})  }{2 }  \prod_{v} \frac{H_v(f_v)}{\langle f_v, f_v \rangle} \end{equation} 
 where $H_v(f_v) =  \int_{T(F_v)} \langle t_v \cdot f_v, f_v \rangle  \psi_v(t_v) dt_v$.
Moreover, if we write $$\heart_v :=   \frac{ \zeta_v(2)  L(\frac{1}{2}, \pi \times \chi_{1,v}) L(\frac{1}{2}, \pi \times \chi_{2,v})}{L(1, \Ad, \pi_v) L(1, \chi_{E/F, v}) },$$
then $H_v(f_v) = \heart_v$ for almost all $v$.    One may check that $\heart_v \in Z$ for all $v$ 
(see discussion of $Z$-structures on page \pageref{Zstructuredef}).

In order to proceed, we need to study the rationality properties of $H_v(f_v)$. 
Such matters are already studied in detail by Waldspurger; again, we simply explicate
what we need. 
 Note that $H_v(f_v)/\langle f_v, f_v \rangle$ is independent of choice
of local inner product $\langle \cdot, \cdot \rangle$ 
and also 
 is independent of $f_v$, up to scaling.
Write  \begin{equation} \label{SvE} S_v(f) :=  \frac{ d_{v}^{1/2}}{\heart_v } \frac{H_v(f_v)}{\langle f_v, f_v \rangle}  \end{equation}
Here we put $d_v = d_{E_v}/d_{F_v}^2$, where $d_L$ is the absolute discriminant of the local field $L$,  and understand $d_v = 1$ in the archimedean case.  (For comparison purposes: $S_v$ corresponds to $\alpha$ in \cite[page 222]{waldsur}).

For $v$ archimedean, let  $f_v^{\circ} \in \pi_v$ be  an $H \cap K_{\infty}$-fixed vector in the minimal $K_{\infty}$-type (so that $\varF(X)$ is of the form $\otimes_{v|\infty} f_v^{\circ} \otimes \varF(X)_f$, for some $\varF(X)_f \in \bigotimes_{v  \ \mathrm{finite}}\pi_v$). 
It is proven by Waldspurger
that \cite[Lemme 15, page 230]{wald}
\begin{equation} \label{lemma15wald} \prod_{v | \infty} S_v(f_v^{\circ}) \in \Q^{\times} \pi^{-[F:\Q]}.\end{equation}

As for $v$ finite, we have $S_v(f_v) \in Z$ whenever $f_v \in \pi_v^Z$.
In fact, we may choose an inner product $\langle -,- \rangle$ on $\pi_v$ such that the values of this inner product  on $\pi_v^Z$ lie in $Z$.   This is an exercise in linear algebra (not sesquilinear algebra: $Z$ is totally real so that the restriction of the inner product to $\pi_v^Z$ is a bilinear form). 
This being so, it is easy to see that $d_v^{1/2} H_v(f_v) \in Z$ when $f_v \in \pi_v^Z$, for $d_v^{1/2} H_v(f_v)$ is a (possibly infinite) integral sum of matrix coefficients $\langle t f_v, f_v \rangle$, each of which belong to $Z$;  but it is easy to see that the infinite sum 
has the same algebraic properties as a certain finite truncation of it and in particular belongs to $Z$.\footnote{The sum descends
to the quotient of $\T(F_v)$ by a compact subgroup, i.e. a finitely generated abelian group, say 
$\Z \times A$ with $A$ finite. The restriction of the function to each 
$(n, a) \in \Z \times \{a\}$ is a function that satisfies a finite recurrence for $n $ large
and also a (possibly different) recurrence for $-n$ large. So our result reduces to the following fact: For such a function, $\sum_{-\infty}^{\infty} F(z)$ can be expressed as an algebraic function of finitely many evaluations $F(z)$. } 
Since we have already  noted (without proof) that $\xi_v \in Z$, thus $S_v(f_v) \in Z$
for each finite $v$.

This also handles the effect of Galois automorphisms: let $\sigma$ be an  embedding of $Z$ into $\C$, and set $\pi_v^{\sigma} = \pi_{v} \otimes_{(Z, \sigma)} \C$. 
There's a natural $\sigma$-linear map $\pi_v^Z \rightarrow \pi_v^{\sigma}$,
given by $w \mapsto w \otimes 1$; 
we denote it by $f \mapsto f^\sigma$, and we have
 \begin{equation} \label{sve2} S_v(f_v^\sigma) = (S_v(f_v))^{\sigma} \ \ (f_v \in \pi_v^Z),\end{equation} 
because indeed we may choose an inner product $\langle, \rangle$ on $\pi_v^{\sigma}$
such that $\langle f_1^{\sigma}, f_2^{\sigma} \rangle = \langle f_1, f_2 \rangle^{\sigma}$
whenever $f_1, f_2 \in \pi_v^Z$.  Thus $H_v(f_1, f_2)^{\sigma} = H_v(f_1^{\sigma}, f_2^{\sigma})$, as follows from the same reasoning that shows $H_v(f_1, f_2) \in Z$.

This in fact concludes the proof of all but the last assertion of the Theorem.

 Now we check the final assertion: that there exists $\data$ such that $\int_{\gamma_{\data}} f \neq 0$.  
 If $\sigma$ is a representation of $\PGL_2(k)$, for $k$ a local field,
 let $\varepsilon(\sigma) \in \pm 1$ be the root number;  in the adelic case
 $\varepsilon$ refers to the product of local root numbers.  (In this instance,
 these are independent of the choice of  additive character.)  Recall:

\begin{enumerate}
\item If $\sigma$ is a principal series induced from the character $\beta$, 
then $\varepsilon(\sigma) = \beta(-1)$. 
\item If $\sigma$ is the Steinberg representation (resp. its twist by an unramified quadratic character), then $\varepsilon(\sigma) = -1$ (resp. $1$). 
\end{enumerate}

\begin{lemma} \label{epstwist} Suppose that $\Pi$ is an automorphic representation of $\PGL_2(\adele)$ 
with a Steinberg (up to twist) place $v$ and let $T$
be an arbitrary set of places disjoint from $v$. Then there exists a quadratic character
$\omega$ such that $\varepsilon(\frac{1}{2}, \Pi \times \omega) = - \varepsilon(\Pi)$,
and moreover $\omega = 1$ at all places in $T$.
\end{lemma}
\proof  Indeed, 
enlarge $T$ to contain all archimedean and ramified places of $\Pi$, with the exception of $v$,  and take $\omega$
such that $\omega$ is unramified nontrivial at $v$ and $\omega$ is trivial at all places in $T$.
Then 
$$\frac{\varepsilon(\Pi \times \omega)}{\varepsilon(\Pi)} =   - \prod_{w \notin T \cup \{v\}} \omega(-1) 
= - \prod_{w \in T \cup \{v\}} \omega(-1)   =-1,$$
as claimed.  \qed

Now 
\begin{df} \label{Xvdef} 
For each ramified place $v$ of either $\pi$ or $D$, let $X_v$ be the set of pairs $(\chi_{1,v}, \chi_{2,v})$
of quadratic characters of $F_v^{\times}$ such that:
\begin{enumerate} \item
$ \chi_1 \chi_2(-1) \varepsilon(\frac{1}{2}, \pi_v \times \chi_1) \varepsilon(\frac{1}{2}, \pi_v \times \chi_2)= 
\begin{cases} 1, \mbox{ $D_v$ split} \\ -1, \mbox{ $D_v$ ramified}. \end{cases}$
\item If $D_v$ is not split then $\chi_1 \neq \chi_2$. 
\end{enumerate}   
We say a pair $(\chi_1, \chi_2)$  of quadratic characters of $\Adele^{\times}/F^{\times}$ is $\pi,D$-admissible if $(\chi_1, \chi_2)_v \in X_v$ for all  such $v$. 
\end{df}
 Note that if $D_v$ is split and $\pi_v$ unramified, then $X_v$ consists of {\em all} 
 pairs of quadratic characters.  Therefore the notion of $(\pi, D)$-admissible
 is in fact a finite set of local constraints. 
 
The assumption of
 Theorem~\ref{regulatorL}  implies in particular that
the set 
$X_v$ is then nonempty for all $v$:
 \begin{enumerate}
\item At finite places such that $D_v$ is nonsplit, 
$F_v^{\times}$ has exactly two unramified quadratic characters;
let these be
 $\{\chi_1, \chi_2\}$.
 
 \item At each real place, let $\{\chi_1, \chi_2\}$ be the two quadratic  characters of $F_v^{\times} \cong \mathbf{R}^{\times}$.

\item At every place that $\pi_v$ is ramified and $D_v$ is split, we arrange that $\chi_1, \chi_2$ are both trivial.

 \end{enumerate}

Now there exists a $(\pi, D)$- admissible pair $(\chi_1, \chi_2)$ of quadratic characters such that $$(*) \ \ \ \varepsilon(\frac{1}{2}, \pi \times \chi_1)= \varepsilon(\frac{1}{2}, \pi \times \chi_2) =1 \ \ (j=1, 2).$$
This follows from our prior comments:  Choose any admissible pair $(\chi_1, \chi_2)$; we twist both $\chi_j$ by a character $\omega$ 
as in the Lemma~\ref{epstwist}, taking for $T$ the union of all ramified places for $D$ or $\pi$.

 It    follows from (*) and \cite[Theorem 4, page 288]{waldS}
that we may choose admissible $\chi_1, \chi_2$ such that 
$L(\frac{1}{2}, \pi \times \chi_i) \neq 0$. 
 Let $E$ be the quadratic field associated to $\chi_1 \chi_2$;
 it embeds into $D$ because of condition (2), and fix such an embedding; 
let  $\T \subset \G$ the corresponding torus.  
Let $\psi: \T(\adele)  = \adele_E^{\times}/\adele_F^{\times} \rightarrow \mathbf{C}^{\times}$ be the character
obtained by pulling back $\chi_1$ or $\chi_2$ via the norm. 
For every $v$ 
the space $\Hom_{\T(F_v)} (\pi, \psi)$ is nonzero: 
 in the case where $\pi_v$ or $D_v$ is ramified this is a result of J. Tunnell\footnote{in the case of odd residue characteristic; need to find citation for general case.}; 
 if $D_v$ is split it  is known that any representation that is not a discrete
 series admits such a functional (cf. \cite[Lemme 8]{wald}; it is proved in Tunnell's paper). 
It then follows that there is  $g_v \in \G(F_v)$ for each finite place $v$
so that $H_v(g_v f_v) \neq 0$ (see \cite[Lemme 10]{wald}). 
The claimed nonvanishing follows from~\eqref{WaldF}.  This concludes the proof of Theorem~\ref{thm:rationality}.
\end{proof}

 \begin{remarkable}  \label{simpledatacomputation}\em{
We note some other computations of $S_v(f_v)$ in the
``unramified case.''
These computations will not be necessary in the present proof and can be skipped for the moment:

Data $\data$ as in~\S~\ref{sec:data} is said to be {\em unramified} if:
\begin{enumerate}  
\item[-] For every finite place $v$, 
the image of the maximal compact subgroup of $(E \otimes F_v)^{\times}$
inside $\T(F_v)$ is contained in $K_{\data}= g^{-1} K g$;
\item[-] The ramification of $\pi, D$ and the characters $\chi_1, \chi_2$ associated to $\data$ is disjoint. In other words, at any place $v$, at most one of the following possibilities occur: 
\begin{itemize}
\item[-]$\chi_1$ is ramified;
\item[-] $\chi_2$ is ramified;
\item[-]
$D$ is ramified and $\pi_v$ is trivial;
\item[-] $D_v$ is split and $\pi_v$ is Steinberg. 
\end{itemize}
\end{enumerate} 

The basic fact we will need is this:
\begin{quote}  (*) If $\data$ is unramified and $f_v \in \pi_v$ 
  is fixed under $\T(F_v) \cap K_{\data}$, then $\frac{S_v(f_v)}{L(1, \chi_{E/F,v})}$
  can be expressed in terms of local $L$-factors for $\pi_v$ only. 
\end{quote}

\begin{itemize}
\item[(a)]
In the case where $\pi, D, \chi_i$ are all unramified $S_v=1$; 

\item[(b)] Suppose that either $\chi_1$ or $\chi_2$ is ramified, but not both;
then $S_v = L(1, \chi_{E/F,v})$.

Without loss of generality, $\chi_2$ is ramified. 
Let $\alpha, \beta$ be the Satake parameters of $\pi$, so that
$L(s, \pi) = (1- \alpha q^{-s})^{-1} (1- \beta q^{-s})^{-1}$. 
Let $\varepsilon$ be the value of $\chi_1$ on a uniformizer
for $F_v$, so that $\varepsilon \in \pm 1$. 

Then  our assumption says that $\T(F_v) \subset g K_v g^{-1}$,
and the matrix coefficient $\langle t f_v, f_v \rangle$
in the definition of $f_v$ takes just two values: 
$1$ and $\sqrt{q} \frac{\alpha+\beta}{q+1}$. 
 So
$$ d_v^{1/2} \frac{ H_v(f_v) }{\langle f_v, f_v \rangle} = 1 + \sqrt{q} \frac{\alpha + \beta}{q+1} = \frac{q}{q+1} ( 1 + \frac{\varepsilon \alpha}{\sqrt{q}})
( 1 + \frac{\varepsilon \beta}{\sqrt{q}}),$$
whereas
$$\xi_v = \frac{1 - \varepsilon q^{-1}}{1-q^{-2}} (1-q^{-1}) 
( 1 + \frac{\varepsilon \alpha}{\sqrt{q}})
( 1 + \frac{\varepsilon \beta}{\sqrt{q}}),$$

\item[(c)] If $\chi_1, \chi_2$ are unramified and $D_v$ or $\pi_v$ is ramified
then 
$$S_v(f_v) =   \zeta_v(2)^{-2}  L(1, \Ad, \pi_v) L(1, \chi_{E/F,v}).$$

In fact, if $D_v$ is ramified, 
this follows from $d_v^{1/2} H_v(f_v)/\langle f_v, f_v \rangle = 1$ and $L(\frac{1}{2}, \pi \times \chi_1) L(\frac{1}{2}, \pi \times \chi_2) =  (1-q^{-2})^{-1} = \zeta_v(2)$.

If $\pi_v$ is ramified, then 
  $d_v^{1/2} H_v(f_v) / \langle f_v, f_v \rangle =  L(\frac{1}{2}, \pi_v \times \chi_{1,v}) L(\frac{1}{2}, \pi \times \chi_{2,v}) / \zeta_v(2).$
 \end{itemize}
 
}\end{remarkable}

 \subsection{Proof of Theorem~\ref{regulatorL}}  \label{ratproof}
We follow notation as  in~\S~\ref{param}.

Let now $\pi_1, \dots, \pi_r$ be the set of all cohomological automorphic representations
with multiplicity $m(\pi_i) \neq 0$, that is to say, $\pi_i^K \neq 0$.  
For each such $\pi$, let $Z_{\pi}$ be the field generated 
be the field generated by the Hecke eigenvalues of $\pi$; it is a totally real, finite extension of $\Q$.   

The Galois group $\Gal(\Qbar/\Q)$ acts on $\{\pi_1, \dots, \pi_r\}$. 
Let $O_1, \dots, O_k \subset \{1, \dots, r\}$ be the orbits of $\Gal(\Qbar/\Q)$ acting on $\{\pi_1, \dots, \pi_r\}$.

Let $H^1_{\pi}(Y, \C)$ be the subspace of $H^1(Y, \C)$ corresponding to $\pi$.
This  is a subspace
of dimension equal to $\dim \pi_{f}^K$, described as the image of the map $\Omega$
followed by the natural map from $1$-forms to cohomology. 
For any $O \in \{O_1, \dots, O_k\}$ we denote by $H^1_{O}$ the $O$-isotypical piece of $\Q$-cohomology, i.e.
$$H^1_O := H^1(Y(K), \Q) \cap \bigoplus_{j \in O} H_{\pi_j}^1(Y, \C).$$
Then $H^1(Y, \Q) = \bigoplus_{O} H^1_O$. 
Let $H_{1,O}$ be the corresponding space in homology (in other words, the subspace
of $H_1(Y, \Q)$ that is the annihilator of $\bigoplus_{O' \neq O} H^1_{O'}$).

In what follows, the inner product on $H^1(Y, \C)$ is always that obtained via identification with harmonic forms.  The images of the spaces $H^1_O$ are orthogonal, for distinct $O$, with respect to this inner product.

Choose an orbit $O$, 
and choose $\pi \in O$. 
We fix:
\begin{itemize}
\item A  $\Q$-basis $\Gamma_O $ for $H_{1,O}$. 
\item A   $\C$-basis $\mathcal{B}(\pi) $
for $ \Hom_{K_{\infty}}(\gk, \pi^{K})$  such that
the elements of $\mathcal{B}(\pi)$ are orthogonal and 
for every $f \in \mathcal{B}(\pi)$, the periods of $\Omega(f)$ belong to $Z$. 
 In other words: the cohomology class of $\Omega(f)$ belongs to $H^1(Y, Z) \subset H^1(Y, \C)$. 

 This is indeed possible: 
There exists a $Z$-subspace $H^1_{\pi; Z} \subset H^1(Y, Z)$
whose complexification is $H^1_{\pi}$: the subspace
on which the Hecke operators act the same way as they do on $\pi$. 
Now the inner product induces on $H^1_{\pi; Z}$ a symmetric, 
positive definite, bilinear form; we choose an orthogonal basis with respect to this form
and then lift the resulting basis via $\Omega$. 
  
 \end{itemize}

 For each $f \in \mathcal{B}(\pi)$, we refer to the corresponding cohomology class
 in $H^1(Y(K), Z)$ as $[f]$.  
 If $\sigma \in \Hom(Z, \C)$, then there is a corresponding element 
 $[f]^{\sigma}$ characterized by the property that 
 $ \langle \gamma, [f^{\sigma}] \rangle = \langle \gamma, [f] \rangle^{\sigma}$
 for any $\gamma \in H_1(\Q)$.  
 This cohomology class is represented by a unique element
of $\pi^{\sigma}$ which we denote by $f^{\sigma}$, i.e.
  $f^{\sigma} \in \Hom_{K_{\infty}}(\gk, (\pi^{\sigma})^{K})$
and $\Omega(f^{\sigma})$ represents the class $[f]^{\sigma}$.

  Then 
 $$ \{ [f]^{\sigma}: f \in \mathcal{B}(\pi), \sigma \in \Hom(Z, \C) \}$$ 
gives a basis for $H^1_O$ that is {\em orthogonal} inside $H^1(Y(K), \C)$: 
 Clearly, $[f_1]^{\sigma}$ and $[f_2]^{\tau}$ are orthogonal for distinct $\sigma, \tau$, 
 and any $f_1, f_2 \in \mathcal{B}(\pi)$, since there exists some Hecke operator
with distinct eigenvalues on them.
It remains to check that $[f_1]^{\sigma}$ and $[f_2]^{\sigma}$ are orthogonal, in the same notation. 
This is true for $\sigma=1$ by choice of the basis $\mathcal{B}(\pi)$,
and the rest  follows from linear algebra:
if $\pi_f$ denotes the finite part of an automorphic representation  $\pi$,  
then $\pi_f$ admits a $Z$-structure $\pi_Z$, and 
there's a natural $\sigma$-linear map $\pi_Z \rightarrow \pi_Z^{\sigma}$;
if we denote this map by $v \mapsto v^{\sigma}$, then we may choose
inner products so that $\langle v^{\sigma}, w^{\sigma} \rangle = \langle v, w \rangle^{\sigma}$. 
 
Now define the matrix $A_O$ as the square matrix whose rows are indexed by $\Gamma_O$, 
columns indexed by pairs $\{ (f, \sigma): f \in \mathcal{B}(\pi),  \sigma \in \Hom(Z, \C)\}$
and whose entry, in the position
that is the intersection of the row indexed by $\gamma \in \Gamma_O $
and the column indexed by $(f \in \mathcal{B}(\pi),  \sigma \in \Hom(Z, \C))$, is 

$$ \langle [f^{\sigma}], \gamma \rangle =
  \left( \int_{\gamma}  \Omega(f) \right)^{\sigma}$$
 Each entry is valued in $\overline{\Q}$, and 
automorphisms of $\overline{\Q}$ merely permute rows of $A_O$. Therefore, $\det(A_O)^2$
is invariant under all such automorphisms, so $\det(A_O)^2 \in \Q$. 

Now let $f \in \mathcal{B}(\pi)$ be arbitrary.  Theorem~\ref{thm:rationality} shows that  there exists $\chi_1, \chi_2$ such that 
\begin{equation} \label{waldO}
\prod_{\sigma} \frac{ \langle [f]^{\sigma}, \gamma_{\data} \rangle^2  }{ \langle f^{\sigma}, f^{\sigma} \rangle} \in \Q^{\times} 
\prod_{\pi \in O}  c
 \frac{ \Lambda(1/2, \pi \times \chi_1) \Lambda(\frac{1}{2}, \pi \times \chi_2)}{\Lambda(1, \Ad, \pi)}  
\end{equation}
is nonzero, 
where $c$ is as in Theorem~\ref{thm:rationality}. 
Note that we were able to omit the absolute values, because all the periods $\int_{\gamma} f$
take values in the real field $Z$.  

By direct computation, if $v$ is an
archimedean place, then $ \frac{ L_v(1/2, \pi \times \chi_1)  L_v(\frac{1}{2}, \pi \times \chi_2)}{L_v(1, \Ad, \pi)} $ is a rational multiple of $\pi$.

Write   $B_O = \prod_{\sigma} \prod_{f \in \mathcal{B}(\pi)}  \langle f^{\sigma}, f^{\sigma} \rangle$
for the product of norms of the $f^{\sigma}$. Since $\prod_{\sigma }
\langle \gamma_{\data}, [f]^{\sigma} \rangle \in \Q$, we see
\begin{equation} \label{almost}  \left( \reg_{\C}(O) \prod_{\pi \in O}   \frac{c'}{L(1, \Ad, \pi)} \right)^{m(O)}  B_O \in \Q^{\times}. \end{equation} 
where $c'$ and $m(O)$ are as in statement of Theorem~\ref{regulatorL}.

Thus \begin{equation} \label{Rzform} \prod_{O}   \frac{ \det A_{O}^2}{B_{O} } \in    \Q^{\times} \prod_{O} \reg_{\C}(O) \prod_{\pi \in O}  \frac{c'}{L(1, \pi, \Ad)}   \end{equation} 
which concludes the proof. %

\section{The regulator under change of level}  \label{sec:regchangelevel}

 It is possible to understand precisely the change of regulator when one changes the level
and there are no new forms:
\begin{theorem} \label{regcompare} 
(For $\G$ nonsplit; see~\S~\ref{h2regregcompare}  for $\G$ split). 
Suppose ${\q} \in \Sigma$ and $H_1^{{\q}-\mathrm{new}}(\Sigma, \C) = \{0\}$. 
Write 
$$D := \det \left( T_{\q}^2 - (1 + N(\q))^2 \big|  H_1(\Sigma/\q, \C) \right).$$ 
 Then, up to orbifold primes, 
$$\frac{ \reg(H_1(\Sigma/\q))^2}{  \reg(H_1(\Sigma))} =  \frac{\sqrt{D}}{h_{\lif}(\Sigma/\q;\q)} $$ %
 \end{theorem}
 Recall the definition of $h_{\lif}$ from 
 Lemma~\ref{regulator-compare-1}: It measures the amount of congruence homology
 which lifts to characteristic zero.

\proof    
We write the proof assuming there are no orbifold primes, for simplicity. 

The proof will be an easy consequence of Lemma~\ref{regulator-compare-1}: 
Consider the map
 $\Psi^{\vee}: H_1(\Sigma/{\q},\Z)_{\tf}^2 \rightarrow  H_1(\Sigma,\Z)^2_{\tf}$ and its dual over $\mathbf{R}$: 
 $\Phi^{\vee} _{\R}: H^1(\Sigma, \R) \rightarrow H^1(\Sigma/{\q}, \R)$. 
 Choose an orthonormal basis $\omega_1, \dots, \omega_{2k}$  
 for $H^1(\Sigma,\R)$, and a basis $\gamma_1, \dots, \gamma_{2k}$ for $H_1(\Sigma/{\q}, \Z)_{\tf}^2$.  Then
 
 $$   \det \  \langle \gamma_i, \Phi^{\vee}_{\R} \omega_j \rangle = \det \ \langle \Psi^{\vee}(\gamma_i), \omega_j \rangle  = | \mathrm{coker}(\Psi)|  \cdot \reg(H_1(\Sigma)).$$ 
 
 Now (by Lemma~\ref{regulator-compare-1}) 
 $|\mathrm{coker}(\Psi^{\vee})| = \frac{D}{h_{ \lif}(\Sigma/\q;\q)}$. 
 The $\Psi_{\R}^* \omega_j$ do not form an orthonormal basis; 
 indeed, the map $\Psi_{\R}^*$  multiplies volume by $\sqrt{D}$, i.e.
 $$ \| \Phi_{\R}^{\vee} \omega_1 \wedge \Phi_{\R}^{\vee} \omega_2 \wedge \dots \Phi_{\R}^{\vee}\omega_{2k} \| = \sqrt{D}.$$  Consequently,
 
 $$ \reg(H_1 (\Sigma/{\q}))=  \frac{ \det \  \langle \gamma_i, \Phi^{\vee}_{\R} \omega_j \rangle }{ \| \Phi_{\R}^{\vee} \omega_1 \wedge \Phi_{\R}^{\vee} \omega_2 \wedge \dots \Phi_{\R}^{\vee}\omega_{2k} \| } = \reg(H_1(\Sigma)) \cdot \frac{\sqrt{D}}{h_{\lif}(\Sigma/{\q};\q)},$$ 
which implies the desired conclusion.  \qed

\section{Comparison of regulators and level lowering congruences} \label{reglowlevel}
  
 In this section --- almost entirely conjectural --- we  present a principle ---~\eqref{PP2} below ---
 that seems to govern the relationship between the regulators 
 between the two manifolds in a Jacquet--Langlands pair (see~\S~\ref{section:discusscongruence} for definition).

  A {\em numerical} example where the regulators differ in an interesting way
  is given in
~\S~\ref{subsec:level lowering};  this affects the behavior of the torsion and motivated us to add this section, despite the conditionality of the results.   We advise the reader to glance at this before studying this section. We observe also that
  our later results comparing Jacquet--Langlands pairs (especially Theorem~\ref{theorem:thirdsetting}) also give evidence for~\eqref{PP2}.

\subsection{The case of elliptic curves}

\label{section:prasanna}
Consider, as a warm-up example, the case of an indefinite division algebra $D$ over $\Q$, ramified at a set of places $S$. 

Let $f$ be an {\em integrally normalized} form on the corresponding adelic quotient $X_{S}$,
a certain Shimura curve. 
(The notion of {\em integrally normalized} is defined in~\cite{Prasanna} using algebraic geometry; it is likely that corresponding results hold when $D$ is definite, for the ``naive'' integral normalization.) 
Let $f^{\JL}$ be the Jacquet--Langlands transfer of $f$ to the split group $\GL_2$: it may be considered as a classical
holomorphic modular form for a suitable congruence subgroup and level. 

In  his thesis, Prasanna established the following principle:
\begin{equation} \label{PP} \mbox{Level-lowering congruences for $f^{\JL}$ at $S$
reduces $\langle f, f \rangle$.  }\end{equation}

A heuristic way to remember this principle is the following: Suppose $f^{\JL}$ is of squarefree level $N$,
and that we realize the Jacquet--Langlands correspondence from $X_0(N)$ 
to $X_S$ in some explicit way --- for example, via $\Theta$-correspondence.
This correspondence is insensitive to arithmetic features, so $\Theta(f^{\JL})$
will not be arithmetically normalized. 
A level lowering congruence for $f^{\JL}$ means that it is congruent (mod $\ell$) to a newform $g$ of level $N/p$, for some $p$.  Heuristically, one might expect $\Theta(f^{\JL}) \cong \Theta(g)$;
but $\Theta(g)$ is zero for parity reasons. Therefore, $\Theta(f^{\JL})$ is ``divisible by $\ell$'';
what this means is that to construct the arithmetically normalized form one needs to 
divide by an additional factor of $\ell$, depressing $\langle f, f \rangle$.

\subsection{From \texorpdfstring{$\Q$}{Q} to imaginary quadratic fields}
Now we analyze the case of an imaginary quadratic field;
for this section, then, we suppose $F$ to be imaginary quadratic. 

Let $Y(\Sigma), Y'(\Sigma)$ be a Jacquet--Langlands pair 
(recall~\S~\ref{section:discusscongruence}). 
These therefore correspond to groups $\G$ and $\G'$ ramified
at sets of places $S, S'$ respectively. 

We assume that there is an equality $$\dim H^1(Y(\Sigma), \Q) = \dim H^1(Y'(\Sigma), \Q) = 1$$
and the corresponding motive is an elliptic curve $\ellcurve$.  
Thus,  $\ellcurve$ is a ``modular'' elliptic curve  over $F$ whose conductor
has valuation $1$ at all primes in $\Sigma$; we put modular in quotes
since there is no uniformization of $C$ by $Y(\Sigma)$ in this setting.

\medskip

Our belief (stated approximately) is that: 

\begin{quote}
Level-lowering congruences
for $\ellcurve$ at places of $S' - S$ (respectively $S-S'$,   respectively outside $S \cup S'$)  reduce (respectively increase, respectively have no effect on) \begin{equation}  \label{PP2}   \frac{\mathrm{reg}(Y)}{\mathrm{reg}(Y')}. \end{equation}    \medskip
\end{quote}

Explicitly, we say that $\ellcurve$ admits a level lowering congruence
modulo $\ell$ at a prime ${\q}$ if  
the $\ell$-torsion $\ellcurve[\ell]$
is actually unramified at $\q$. If we working over $\Q$, this would actually mean
there exists a classical modular form at level $\Sigma/\q$ that
is congruent to the form for $\ellcurve$, modulo $\ell$; in our present
case, this need not happen (see again the numerical example of
\S~\ref{subsec:level lowering}).  In this instance, in any case, we expect
that the regulator ratio of~\eqref{PP2}  would contain a factor of $\ell$
either in its numerator or denominator, as appropriate.

\medskip

A little thought with the definition of $\mathrm{reg}(Y)$
shows that this is indeed analogous to~\eqref{PP}: the smaller that $\reg(Y')$ is, 
the smaller the $L^2$-norm of a harmonic representative
for a generator of $H^1(Y', \Z)$. 
We will not formulate a precise conjecture along the lines of~\eqref{PP2}:  at this stage proving it seems
out of reach ---  except in the case where $C$ is the base-change of an elliptic curve from $\Q$, cf. comment on page \pageref{CMform} ---   but we hope it will be clear from the following discussion what
such a conjecture would look like. 

\medskip 

 Note that $\ellcurve$
has multiplicative reduction at $v \in \Sigma$.
For any quadratic character $\chi$ of $\adele_F^{\times}/F^{\times}$, 
 let $\ellcurve_\chi$ be the quadratic twist of $\ellcurve$ by $\chi$,
 and $\Sha(\chi)$ its Tate-Shafarevich group. 
The BSD conjecture predicts that
\begin{equation} \label{bsd} L(\ellcurve_\chi, \frac{1}{2}) = \frac{2}{\sqrt{\Delta_F}} 
\left( \prod_{v} c_v(C_{\chi})  \right) \frac{\Sha(\ellcurve_{\chi}) \Omega(\ellcurve_{\chi})}{\ellcurve_{\chi}(F)_{\tors}^2} \end{equation}
where $\Omega(\ellcurve_\chi)$ is  the complex period of $\ellcurve$, and
$c_v$ is the number of components of the special fiber of the N{\'e}ron model of $C_{\chi}$
at $v$.

The following (very conditional) Lemma, then, provides some evidence
for the belief~\eqref{PP2}.  (We note that our later results with the Cheeger-M{\"u}ller theorem
also gives evidence for the truth of~\eqref{PP2} in some special cases; see Theorem~\ref{theorem:thirdsetting} and the prior discussion. However, it would be much more helpful to have a direct
proof of~\eqref{PP2} independent of the Cheeger-M{\"u}ller theorem.)

 \begin{lemma}  \label{veryconditional}
 Suppose the validity of the conjecture of Birch and Swinnerton-Dyer, in the for 
~\eqref{bsd}. Let $\ell$
 be a prime larger than $3$ that is not Eisenstein for $\ellcurve$ (i.e.,
 the Galois representation of $G_F$ on $\ell$-torsion is irreducible). 
 Suppose moreover that 
 \begin{itemize}
 \item[i.] Cycles $\gamma_{\data}$ associated 
to  unramified data $\data$ (definition as in Remark
\ref{simpledatacomputation}) generate the first homology 
$H_1(Y(\Sigma), \Z/\ell \Z)$ and $H_1(Y'(\Sigma), \Z/\ell \Z)$. 
\item[ii.] Fix a finite set of places $T$ and, for each $w \in T$, 
a quadratic character $\chi_w$ of $F_w^{\times}$; 
then there exist infinitely many quadratic characters $\chi$ such that
$\chi|{F_w} = \chi_w$ for all $w$,  $\ell$ does not divide the order of $\Sha(\ellcurve_{\chi})$,
and $\ellcurve_{\chi}$ has rank $0$. 
\end{itemize} 
  Then\footnote{This ratio of $H_1$-regulators differs by a volume factor 
  from $\reg(Y)/\reg(Y')$, 
  by~\eqref{regsimpledef}.  This arithmetic significance of the volume factors
  is discussed, and completely explained in terms of congruence homology, in
~\S~\ref{section:refined}.}
\begin{equation} \label{geelong} \left( \frac{\reg(H_1(Y,\Z))}{\reg(H_1(Y',\Z))} \right)^2 \sim  \frac{ \prod_{v \in S} v(j) }{\prod_{v \in S'} v(j)} \end{equation} where $\sim$
denotes that the powers of $\ell$ dividing both sides is the same; $j$ is the $j$-invariant of $C$,  and $S$ (resp. $S'$) is the set of places which are   ramified for $D$ (resp. $D'$).
  \end{lemma} 

In order to 
make clear the relation to~\eqref{PP2}, note that prime divisors of $v(j)$ (for $v \in S$)
measure ``level lowering congruences'' at $v \in S$: by Tate's theory, for such $v$, 
the $\ell$-torsion $C[\ell]$ over $F_v$ is isomorphic as a Galois module -- up 
to a possible twist by an unramified character --   to the quotient
$$ \langle q_v^{1/\ell}, \mu_{\ell} \rangle / q_v^{\Z}.$$
In particular, if $v(q_v) = v(j)$ is divisible by $\ell$
and  $\ell$ is co-prime to the residue characteristic
of $F_v$, the mod $\ell$-representation
is unramified at $v$ if the residue characteristic of $v$ is different from $\ell$,
and ``peu ramifi\'{e}e'' (in the sense of
Serre~\cite{Serre2})  if the residue characteristic of $v$ is
$\ell$. This is, conjecturally, exactly the condition that detects
whether $\ell$ should be a level lowering prime.  So the right-hand side indeed compares level lowering congruences on $Y$ and $Y'$.

A word on the suppositions. Item (ii) is essentially Question B(S) of \cite{Prasanna}.
We have absolutely no evidence for (i) but we can see no obvious obstruction to it, 
and it seems likely on random grounds absent an obstruction.

\proof 
 
 In short, the idea is this: To evaluate $\reg(H_1(Y))$ and $\reg(H_1(Y'))$, we integrate
 harmonic forms on $Y$ (resp. $Y'$) against cycles $\gamma_{\data}$ as in
~\S~\ref{regulatorL}.  As in that section, these evaluations are $L$-values; when we use Birch Swinnerton-Dyer,
we will find that the Tamagawa factors $c_v$, as in~\eqref{bsd}, are slightly different on both sides, and this
difference corresponds to the right-hand side of~\eqref{geelong}.  
 
 \medskip

Define $\mathcal{X}$  to be the set of pairs $(\chi_1, \chi_2)$
that satisfy $(\chi_{1,v}, \chi_{2,v}) \in X_v$, where $X_v$ is as on page \pageref{Xvdef},
and moreover $\chi_i$ are unramified at all ramified places for $\pi$ or $D$; 
define $\mathcal{X}'$ similarly, but for $D'$.

Let $\omega$ (resp $\omega'$) be a harmonic form on $Y$ resp. $Y'$ corresponding to $\ellcurve$
and let $\pi$ resp. $\pi'$ be the corresponding automorphic representation. 
Let $R$ be  the generator of the subgroup $\int_{\gamma} \omega \subset \R$, where $\gamma$ ranges over $H_1(Y(\Sigma), \Z)$.  Define similarly $R'$. 
 
 Remark~\ref{simpledatacomputation}, and statement (*) there in particular,  shows that, 
if $\omega$ is a harmonic form arising from a new vector in a representation $\pi$,
all of whose ramification is Steinberg\footnote{that is to say: the transfer of $\pi$
to $\GL_2$ has every local constituent at a finite prime $v$ either unramified or Steinberg.}, and if $\data$ is unramified, then:
\begin{equation} \label{proportionality}  \frac{ |\int_{\gamma_{\data}} \omega |^2}{\langle \omega, \omega \rangle} \propto \Delta_{\chi_1 \chi_2}^{1/2}  L(1/2, \pi \times \chi_1) L(1/2, \pi \times \chi_2)  \end{equation} 
{\em where the constant of proportionality is the same for} $\pi$ on $D$
or its Jacquet--Langlands transfer $\pi'$ on $D'$, and
we wrote $\Delta_{\chi}$ for the discriminant of the quadratic field extension of $F$
attached to the character $\chi$.

For each  unramified datum $\data$ for $Y(\Sigma)$ we may write $\int_{\gamma_{\data}} \omega = m_{\data} R$, where
$m_{\data} \in \Z$; similarly,
for each unramified datum $\data$ for $Y'(\Sigma)$ we may write
$\int_{\gamma_{\data}} \omega' = m_{\data}' R'$. 

Our assumption (i) implies that the greatest common divisor of $\m_{\data}$, 
where $\data$ ranges over unramified data, is relatively prime to $\ell$;
similarly for $m'$.
  
It follows from~\eqref{proportionality} that 
$$\left( \frac{\reg(H_1(Y))}{\reg(H_1(Y'))}\right)^2 \sim \frac{g}{g'},$$ where 
$$ g =  \mathrm{gcd}_{(\chi_1, \chi_2)}  \Delta_{\chi_1 \chi_2}^{1/2} L(1/2, \pi \times \chi_1) L(1/2, \pi \times \chi_2),$$ and $g'$ is defined similarly; again $\sim$ denotes that the powers of $\ell$ dividing both sides coincide.

We shall now apply the conjecture of Birch and Swinnerton-Dyer~\eqref{bsd}:
  $$L(1/2, \pi \times \chi) =2  \frac{ \left( \prod_{v} c_v(\ellcurve_\chi)  \right) \Sha(\ellcurve_{\chi})}{\# \ellcurve_{\chi}(F)_{\tors}^2}  \frac{\Omega_{\ellcurve_{\chi}}}{\Delta_{F}^{1/2}} $$
 We  have
$\Omega(\ellcurve_\chi) = \Omega(\ellcurve) \cdot \Delta_\chi^{-1/2}$ and
$\Delta_{\chi_1} \Delta_{\chi_2} = \Delta_{\chi_1 \chi_2}$;
both of these use the assumption of unramified data $\data$. 
 Since $\ell$ is not Eisenstein,  the $\ell$-part of $\ellcurve_{\chi}(F)_{\tors}$ 
is trivial for any $\chi$.  Also, by~\cite{RubinFudge}, the product
of $c_v(\ellcurve_{\chi})$ over all $v \notin \Sigma$ is a power of $2$.   

Therefore, 
$$ \frac{g}{g'} \sim \frac{  \gcd_{(\chi, \chi') \in \mathcal{X}} \Sha(\ellcurve_{\chi})  
\Sha(\ellcurve_{\chi'}) \prod_{v \in \Sigma} c_v(\ellcurve_\chi)   c_v(\ellcurve_{\chi'}) }{  \gcd_{(\chi, \chi') \in \mathcal{X}'} \Sha(\ellcurve_{\chi}) \Sha(\ellcurve_{\chi'}) \prod_{v \in \Sigma} c_v(\ellcurve_\chi) c_v(\ellcurve_{\chi'})},$$
where, again, $\sim$ means that the powers of $\ell$ dividing both sides coincide. 

\medskip

For $v \in \Sigma$, put 
$A_v:= \gcd_{(\chi_1, \chi_2) \in \varchi} c_v(\ellcurve_{\chi}) c_v(\ellcurve_{\chi'}).$ Now  $c_v(\ellcurve_{\chi})$ is wholly 
determined by the $v$-component $\chi_v$ for $v \in \Sigma$. 
Therefore,   assumption (ii) of Lemma~\ref{veryconditional} implies that
$$g_{\ell} \sim  \prod_{v \in \Sigma} \frac{ A_v}{A_v'},$$
Indeed, requiring that both  $(\chi_1, \chi_2) \in \varchi$ and
$c_v(\ellcurve_{\chi}) c_v(\ellcurve_{\chi'}) \sim A_v$ is a finite set of local constraints, so that 
(ii) is applicable. 
 
Now an examination of the definition of $\mathcal{X}$ and $\mathcal{X}'$ shows that:
$A_v \sim 1$ {\em unless} $D_v$ is ramified, and $A_v' \sim 1$ {\em unless} $D'_v$ is ramified. 
In those cases, $C$ has multiplicative reduction at $v$, and, as we have observed 
$$A_v \sim c_v(\ellcurve) = |v(j)|.$$

 This concludes the proof. 
\qed

\chapter{The split case} 
\label{chapter:ch5}
This chapter is concerned with various issues, particularly of analytic nature, that arise when $\G$ is split. 
 \medskip

One  purpose of this Chapter is to generalize certain
parts of the Chapter 4 to the noncompact case --
these results are concentrated in~\S~\ref{sec:llsplit},~\S~\ref{sec:Eis-series}.
and~\S~\ref{sec:eisintegrality}.  These parts contain some results of independent interest 
(e.g. the upper bound on Eisenstein torsion given in~\S~\ref{sec:eisintegrality}). 

\medskip

However, the main purpose of the chapter is to prove Theorem~\ref{thm:rtatsplit};
it compares  Reidmeister and analytic torsion in the non-compact case.
 More precisely, we compare a {\em ratio} of Reidemeister torsions to a  {\em ratio} of analytic torsions,
 which is a little easier than the direct comparison. 
 This is the  lynchpin of our comparison of torsion homology between a Jacquet--Langlands pair
in the next chapter. 
 The most technical part of the theorem
is~\eqref{rtatsplit}.   

\medskip
For the key analytic theorems we do not restrict to level structure of the form $K_{\Sigma}$. One  technical reason for this is that,
even if one is interested only in those level structures, the Cheeger-M{\"u}ller theorem does not apply to orbifolds. 
Rather we need to apply the equivariant version of this theorem to a suitable covering $Y(K') \rightarrow Y(K_{\Sigma})$. 

 \medskip

 \medskip
 
 Once the definitions are given 
 the flow of the proof goes: 
 \begin{equation} \mbox{~\eqref{rtatsplit} $\impliedby$ Theorem~\ref{theorem:invtrunc} $\impliedby$ Theorem~\ref{prop:SE} } \end{equation}
The proof of Theorem~\ref{prop:SE} is in~\S~\ref{smalleigenvalues};
the left- and right- implication arrows are explained in~\S~\ref{rtatsplitfollows} and~\S~\ref{theorem:invtruncproof}, respectively.
 
Here is a more detailed outline of some of the key parts of this Chapter: 

{\small 
\begin{itemize}
\item~\S~\ref{section:split} introduces the basic vocabulary for dealing with non-compact manifolds
In particular~\S~\ref{subsec:trunc} introduces the notion of  {\em height function} that measures how far one is in a cusp. 

\item~\S~\ref{sec:split2} discusses eigenfunctions of the Laplacian on non-compact manifolds.
In particular,~\S~\ref{subsec:einsteinintro} gives a quick summary of the theory of Eisenstein series, that is to say, the parameterization of the continuous spectrum of the Laplacian on forms and functions. 

\item~\S~\ref{sec:rtatnc} is devoted to a definition of Reidemeister torsion and analytic torsion in the non-compact case.
In particular:
\begin{itemize} 
\item~\S~\ref{ip:polygrowth} 
discusses  harmonic forms of polynomial growth and, in particular, specify an inner product on them. 

\item~\S~\ref{subsec:rtnonsplitdef} gives a definition of Reidemeister torsion. 
 
\item~\S~\ref{subsec:locality} gives a definition of analytic torsion in the non-compact case.
\end{itemize}
\item~\S~\ref{sec:arithmeticmanifolds} analyzes in more detail some of the foregoing definitions
for arithmetic manifolds $Y(K)$. In particular,
\begin{itemize}
\item~\S~\ref{infgeom} specifies a height function on the $Y(K)$ (the same height function that appears for instance in the adelic formulation of the trace formula); 
\item~\S~\ref{subsec:arithhomologyBM} analyzes the difference between homology and Borel--Moore homology from the point of view of Hecke actions. 
\item~\S~\ref{ybybxi} shows that, upon adding level structure at a further prime,  thus
replacing a subgroup $K_{\Sigma}$ by $K_{\Sigma \cup {\q}}$, the cusps
simply double: the cusps of $Y(\Sigma \cup \q)$ are isometric to two copies of the cusps of $Y(\Sigma)$. (This isometry, however, does not preserve the height functions of~\S~\ref{infgeom}). 

\end{itemize} 
\item~\S~\ref{sec:llsplit} extends certain results related to Ihara's lemma to the non-compact case. 

\item~\S~\ref{sec:Eis-series} carries out explicit computations related to Eisenstein series on arithmetic groups, in particular, computation of scattering matrices on functions and forms.
We strongly advise that the reader skip this section at a first reading. 

\item~\S~\ref{sec:eisintegrality} gives an upper bound on modular forms of ``cusp--Eisenstein'' type 

\item~\S~\ref{analysissec} formulates the central result, Theorem~\ref{theorem:invtrunc}:  comparison of Reidemeister and analytic torsion in the non-compact case.  We then explain why~\eqref{rtatsplit} follows from this Theorem. 

However, our result here is adapted solely to our context of interest: we analyze the {\em ratio} of these quantities between two manifolds with isometric cusp structure.
 
 \item 
\S~\ref {smalleigenvalues} is the technical core of the Chapter, and perhaps of independent interest. It computes the ``near to zero'' spectrum of the Laplacian on a truncated hyperbolic manifold. In particular it shows that these near-zero eigenvalues are modelled by the roots
of certain functions $f(s), g(s)$ related to the scattering matrices. 
\begin{itemize}
\item~\S~\ref{subsec:fganal} is purely real analysis: it analyzes zeroes of $f(s), g(s)$. 
\item~\S~\ref{quasimodeargument} shows that {\em every zero of $f$ or $g$
gives rise to an eigenvalue of the Laplacian.}
\item~\S~\ref{U2} goes in the reverse direction: {\em every eigenvalue of the Laplacian
is near to a root of $f$ or $g$.}
\item~\S~\ref{combin} discusses issues related to eigenvalue $0$. 
\end{itemize}

\item~\S~\ref{theorem:invtruncproof} gives the proof of the main theorem (comparison between Reidemeister and analytic), using heavily the results of~\S~\ref{smalleigenvalues}.
\end{itemize} 
}

\section[Noncompact hyperbolic manifolds]{Noncompact hyperbolic manifolds: height functions and homology} \label{section:split}

\subsection{Hyperbolic manifolds and height functions} \label{subsec:trunc}

Let $M$ be a finite volume  hyperbolic $3$-manifold ---  as usual, this may be disconnected-- with cusps.

We think of the {\em cusps} $\mathcal{C}_1, \dots, \mathcal{C}_k \subset M$
as being $3$-manifolds with boundary, so that, firstly, 
$M - \bigcup \mathrm{interior} \  (\mathcal{C}_i)$ is a compact
manifold with boundary; and, moreover, each $\mathcal{C}_j$ is isometric to
\begin{equation} \label{Standard} \{ (x_1, x_2, y) \in {\H}^3:  y \geq 1 \}/\Gamma_{\mathcal{C}},\end{equation}
where a finite index subgroup of $\Gamma_{\mathcal{C}}$ acts as a lattice of translations on $(x_1, x_2)$.   We call a cusp {\em relevant} if all of $\Gamma_{\mathcal{C}}$ acts by translations. \index{relevant cusps}
Cusps which are not relevant have orbifold singularities. 

We may suppose that the cusps are disjoint; for each cusp $\mathcal{C}_i$, fix an isometry $\sigma_i$ onto 
a quotient of the form~\eqref{Standard}. 
For each $x \in M$,  we define the ``height'' of $x$ as a measure
of how high it is within the cusps: 
set \begin{equation} \label{heightdef} \height(x) = \begin{cases} 1,  & x \notin \bigcup \mathcal{C}_i \\ 
\mbox{$y$-coordinate of $\sigma_i(x)$},&  x \in \mathcal{C}_i . \end{cases} \end{equation}

A {\em height function on $M$} is a function $\height $ that arises in the fashion above;
such a function is far from unique, because the choice of the isometry $\sigma_i$ is not unique. 
If, however, we fix a height function on $M$, then
``compatible'' $\sigma_i$s (compatible in that~\eqref{heightdef} is satisfied for sufficiently large heights) 
are unique up to affine mappings  in the $(x_1, x_2)$-coordinates. 
We will often use $(x_1, x_2, y)$ as a system of coordinates on the cusp, with the understanding that it is (after fixing a height function) unique up to affine mappings in $(x_1, x_2)$.

In what follows, we suppose that $M$ is endowed with a fixed choice of height function.
\footnote{It is necessary to be precise about this because we will use the height function to truncate $M$, and then compute the Laplacian spectrum of the resulting manifold; this is certainly dependent on the choice of $\height$. }  
In the case of  $M = Y(K)$ an arithmetic
manifold, we shall specify our choice of height function later. 
We shall  often use the shorthand
$$ y \equiv \height$$
e.g., writing ``the set of points where $y \geq 10$'' rather than ``the set of points
with height $\geq 10$'', 
 since, with respect to a suitable choice of coordinates on $M_B$,
$\height$ coincides with the $y$-coordinate on each copy of $\H^3$.

This being so, we make the following definitions:
\begin{enumerate}
 \item 
For any cusp $\mathcal{C}$, we define $\area(\mathcal{C})$ to be the area of the quotient \index{area of cusp}
$\{(x_1, x_2) \in \R^2\}/\Gamma_{\mathcal{C}}$ ---  i.e., if $\Gamma'$ is a subgroup of $\Gamma_{\mathcal{C}}$ that acts freely, this is by definition $\mathrm{area}(\R^2/\Gamma') \cdot [\Gamma_{\mathcal{C}}:\Gamma']^{-1}$.  The measure on $\R^2$ is the standard Lebesgue measure. 

Once we have fixed a height-function, this is independent of the choice of compatible
isometries $\sigma_i$.  We write $\area(\partial M)$ for the sum $\sum_{\mathcal{C}} \area(\mathcal{C})$, the sum being taken over all cusps of $M$; we refer to this sometimes as the ``boundary area'' of $M$. 

\item For any cusp $\mathcal{C}$, we write $\vol(\mathcal{C})$ for the volume of the connected component of $M$ that contains $\mathcal{C}$. 

\item Put  $$M_B = \coprod_{i} \H^3/\Gamma_{\mathcal{C}},$$ a finite union of hyperbolic cylinders. \index{$M_B$}
Then ``the geometry of $M$ at $\infty$ is modelled by $M_B$.''  Note that $M_B$
is independent of choice of height function. 

The $y$-coordinate on each copy of $\H^3$ descends
to a function on $M_B$. Thus, for example, $y^s$ defines a function on $M_B$
for every $s \in \C$. 

\item 
 The height $Y$ truncation of $M$ denoted $M_{\leq Y}$, is defined as \index{$M_{\leq Y}$}
 $$ M_{\leq Y} =  \{  y \leq Y \} \subset M $$ 
 in words, we have ``chopped off'' the cusps. $M_{\leq Y}$ is a closed manifold with boundary, and its  boundary is isometric to a quotient of a $2$-torus by a finite group ---  or, more precisely,
 a union of such.  We sometimes abbreviate $M_{\leq Y}$ to $M_Y$ when there will be no confusion.  \index{$M_{\leq Y}$}
 \index{$M_{[Y, Y']}$.}
We will similarly use the notation $M_{[Y', Y]}$ to denote elements of $M$ satisfying $Y' \leq y \leq Y$. 
  
{\textbf {\em Convention on boundary conditions:}} When dealing with the Laplacian
 on forms or functions on a manifold with boundary, we shall {\em always suppose
 the boundary conditions to be absolute}:  
 we work on the space of differential forms $\omega$ such that both $\omega$ and $d\omega$, when contracted with a normal vector, give zero.

 \item  We set $\partial M$ to be the boundary of $M_{\leq Y}$ for any sufficiently large $Y$;
 we are only interested in its homotopy class and not its metric structure, and as such this is independent of $Y$.   \index{boundary $\partial M$ of $M$} \label{boundarypage}
 
 We have long  exact sequences
 $$\cdots \rightarrow H^j_c(M,  \Z) \rightarrow H^j(M, \Z) \rightarrow H^j(\partial M, \Z) \rightarrow H^{j+1}_c \rightarrow \cdots $$
and we define as usual {\em cuspidal cohomology}
\begin{equation} \label{cuspcohomologydef} H^j_{!}(M, \Z) = \mathrm{image}(H^j_c(M, \Z) \longrightarrow H^j(M, \Z)) \end{equation} 
  We denote correspondingly
$$H_{j, !}(M, \Z) = \mathrm{image}(H_j(M,\Z) \longrightarrow H_{j, \bm}(M, \Z)),$$
where $\bm$ denotes Borel--Moore homology.

We have the following isomorphisms, where, if $M$ has orbifold points,
one needs to additionally localize $\Z$ away from any orbifold primes:
\begin{equation}  \label{splitduality}
\begin{aligned} 
H^{3-n}_cM,\Z) \simeq  \ & H_n(M,\Z), \\
H^{3-n}(M,\Z) \simeq \ &  H^{\bm}_n(M,\Z), \\
H^{2-n}(\partial M,\Z) \simeq  \ & H_n(\partial M,\Z).
\end{aligned}
\end{equation} 
The first two are Poincar\'{e} duality and Borel--Moore duality respectively.
The last comes from Poincar\'{e} duality for the closed manifold
$\partial Y(K_{\Sigma})$.

 \item We put $b_j(M) = \dim H_j(M, \C)$, the $j$th complex Betti number. In particular, $b_0(M)$ is the number of connected components of $M$. 
 
 \end{enumerate}

As mentioned, in the arithmetic case 
we shall endow $Y(K)$ with a certain, arithmetically defined, height function. 
 In this case we write $Y(K)_T$ and $Y_B(K)$ for $M_{\leq T}$ and $M_B$ respectively.  
We shall later (\S~\ref{infgeom}) describe $Y(K)_T$ and $Y_B(K)$ in adelic terms; the primary
advantage of that description is just notational,  because  it handles the various connected components and
cusps in a fairly compact way. The truncation $Y(K)_T$ is in fact  diffeomorphic to the so-called
 Borel-Serre compactification of $Y(K)$ but all we need is that the inclusion 
 $Y(K)_T \hookrightarrow Y(K)$ is a homotopy equivalence. 
 
\subsection{Linking pairings}  \label{section:linking2}
 We now discuss the analogue of the linking pairing (\S~\ref{section:linking}) in the non-compact case. 
Suppose that $p$ does not divide the order of any orbifold prime for $M$. 

We have a sequence:
{\small
\begin{equation} \label{leftrightseq} H_1(M, \Z_p)  \stackrel{\sim}{\rightarrow} H^1(M, \Q_p/\Z_p)^{\vee}
\stackrel{\sim}{\rightarrow} H^2_c(M, \Z_p) \stackrel{g}{\leftarrow} H^1_c(M, \Q_p/\Z_p)
= H_{1, \bm}(M, \Z_p)^{\vee}\end{equation} 
}
where the morphisms are, respectively, duality between homology and cohomology, 
Poincar{\'e} duality, a connecting map (see below), and duality between
Borel--Moore homology and compactly supported cohomology. 
The map $g$ is a connecting map in the long exact sequence:
 $$ \dots \rightarrow H^1_c(M, \Q_p) \rightarrow H^1_c(M, \Q_p/\Z_p)  \stackrel{g}{ \rightarrow }H^2_c(M, \Z_p) 
\rightarrow H^2_c(M, \Q_p) \rightarrow \dots $$
which we may reformulate as an isomorphism
$$H^1_c(M, \Q_p/\Z_p)/ \left(  H^1_{c,\divi} \right) \stackrel{\sim}{\rightarrow} H^2_c(M, \Z_p)_{\tors}.$$
From~\eqref{leftrightseq} above we now  obtain an isomorphism
$$H_1(M, \Z_p)_{\tors} \stackrel{\sim}{\rightarrow} H_{1, \bm}(M, \Z_p)_{\tors}^{\vee}$$

By means of the map $H_1(M, \Z_p) \hookrightarrow H_{1, \bm} (M, \Z_p)$, 
we obtain a pairing on $H_{1, \tors}(M, \Z_p)$ but it is not perfect.  There are, however,
two situations of interest where one can deduce a perfect pairing from it:

\begin{itemize} \item
Localization: Suppose that $\mathbf{T}$ of correspondences acts on $M$ (e.g. a Hecke algebra)
and $\mathfrak{m}$ is an ideal of $\mathbf{T}$ with the property that
the map $H_1(M, \Z_p) \rightarrow H_{1, \bm}(M, \Z_p)$ is an isomorphism
when completed at $\mathfrak{m}$, then we do obtain a perfect pairing:
$$\left(H_1(M, \Z_p)_{\tors} \right) _{\mathfrak{m}} \times \left( H_1(M, \Z_p)_{\tors}\right)_{\mathfrak{m}} \longrightarrow \Q_p/\Z_p.$$
\item Cusps have no homology: If we suppose that $H_1(\partial M, \Z_p) = 0$
(this will be often true for us, by~\S~\ref{cuspshavenohomology}), then we have a sequence
$$H_1(M, \Z_p) \hookrightarrow H_{1, \bm}(M, \Z_p) \twoheadrightarrow \ker(H_0(\partial M, \Z_p) \rightarrow H_0(M, \Z_p))$$
The final group is torsion free, and so we obtain an
isomorphism 
$$H_1(M, \Z_p)_{\tors} \stackrel{\sim}{\rightarrow} H_{1, \bm}(M, \Z_p)_{\tors}.$$
Thus in this situation we obtain a perfect pairing on $H_{1, \tors}(M, \Z_p)$. 
\end{itemize}

\section{Noncompact hyperbolic manifolds: eigenfunctions and Eisenstein series} \label{sec:split2} 
 
 \subsection{Constant terms}
We continue with our prior notation: $M$ is a non-compact hyperbolic $3$-manifold endowed with a height function $\height$. 

We say that a function (or  differential-form) $f$ on $M$  is of polynomial growth if 
$\| f(x)\| \leq A \cdot  \height(x)^B$ for suitable $A, B$.

The {\em constant term} of the function $f$,  denoted by $f_N$, is the function on $\mathcal{C}$
that  is obtained by ``averaging over toric cross-sections,'' characterized by the property:
\begin{equation} \label{constanttermdef} f_N \circ \sigma_i^{-1}(z) = \frac{1}{\area (\mathcal{C}_i)}  \int_{(x_1, x_2) \in \R^2/\Gamma_{\mathcal{C}_i}} f \circ \sigma_i^{-1}(x_1, x_2, y) dx_1 dx_2,\end{equation}
we shall again give an adelic formulation of this later.  This has a natural analogue  for $f$ a differential form ---  one simply writes $f \circ \sigma_i^{-1}$ 
as a combination $A dx_1 + B dx_2 + C dy$ and replaces each of $A, B, C$
by their averages $A_N, B_N, C_N$.   In fact,
although~\eqref{constanttermdef} defines $f_N$ only on the cusp $\mathcal{C}$,
this definition can be naturally extended to $M_B$ by integrating on the corresponding horoball on $M$.

 Recall that a function (or differential form) is called {\em cuspidal} if $f_N$ is identically zero.

If $f $ is an eigenfunction, then $f$ is asymptotic to $f_N$ ---  denoted  $f \sim f_N$ ---  by which we mean the function $f - f_N$ on $\mathcal{C}$
decays faster than any power of $\height(x)$.  This rapid decay
takes place once the height is sufficiently large; intuitively the ``wave'' corresponding to $f$
cannot penetrate into the cusp once the width of the cusp is shorter than its wavelength.

More precisely:

 \begin{lemma*} If  $f$ is an eigenfunction of eigenvalue $-T^2$, then  there is an absolute constants $b_0, A$  for which
\begin{equation} \label{precisebound}  \left| |f(x) - f_N(x)\right| \leq  \|f\|_{L^2(M_{\leq Y})} \exp(- b_0  \cdot \height(x)), \ \ A (1+T^2) \leq \height(x) \leq Y/2. \end{equation}   
If  $f$ satisfies relative or absolute boundary conditions on $M_{\leq Y}$, then
\begin{equation} \label{precisebound2}  \left|f(x) - f_N(x)\right| \leq  \|f\|_{L^2(M_{\leq Y})} \exp(- b_0  \cdot \height(x)), \ \ \height(x) \geq  A (1+T^2). \end{equation} 
that is to say, the prior estimate \eqref{precisebound} actually holds up to the boundary of $M_{\leq Y}$.

In both cases, similar estimates hold on the right-hand side for any derivative ---  i.e., any monomial in $\partial_{x_1}, \partial_{x_2}, \partial_y$ ---   of fixed order.  
A similar estimate holds if $f$ is an eigenfunction of the Laplacian on $j$-forms.

\end{lemma*}

This Lemma is, in fact, quite important to us. We will later show that eigenfunctions on $M_{\leq Y}$
(for $Y$ large) are well-approximated by restrictions of eigenfunctions on $M$ (at least
for small eigenvalue), and this plays an important role in showing that the eigenfunctions on $M_Y$
are not too big near the boundary. 

The proof is by considering the Fourier expansion; each
Fourier coefficient satisfies a differential equation  that forces it to be an explicit
multiple of a Bessel function, and these can be analyzed by hand. The slightly unusual statement is due to the possibility of the Fourier expansion of $f$ containing terms of exponential increase
in $y$.

\proof   
This is proved by direct analysis of Fourier expansions. We give the proof only for zero-forms, the other cases being similar:

According to \cite[Theorem 3.1]{EGM} for instance (or rather the proof of {\em loc. cit.}), $f-f_N$ admits a Fourier expansion that is a sum of terms expressed, in local coordinates $(x_1, x_2, y)$
on each cusp, as
$$ (x_1, x_2, y) \mapsto \sum a_{\mu} F_{s}(2 \pi |\mu| y) e^{2 \pi i \langle \mu, (x_1, x_2) \rangle}$$
 here $\mu$ varies over the dual lattice to the translation lattice of the cusp,\footnote{  We can identify 
the translation lattice of a cusp with a sublattice of $\mathbf{C}$
via $(x_1, x_2) \mapsto x_1+i x_2$; we can identify then its dual
with a sublattice of $\mathbf{C}$ via the pairing $(z,w ) \mapsto \mathrm{Re}(z \overline{w})$,
and therefore we may think of $\mu \in \mathbf{C}$ also, thus the notation $|\mu|$.}
and $F \in \mathrm{span}\{I_s, K_s\}$ is a linear combination of $I_s$ and $K_s$-Bessel functions, and $1-s^2$ is the eigenvalue of $f$.   
In what follows we are going to assume that $|\mu| \geq 1$ for every such $\mu$. If
this is not true it merely changes the constants $b_0, A$ in the statement of the Lemma.

 The desired relation ~\eqref{precisebound} follows from the asymptotic expansion of $I$ and $K$, namely,  $I(y) \sim \frac{e^y}{\sqrt{2 \pi y}}$
and $K(y) \sim \frac{e^{-y}}{\sqrt{2 y/\pi}}$, this being true   in the range $|y| \gg (1+|s|)^2$.   (See \cite[page 378]{AbramowitzStegun} for more precise asymptotic expansions.)
Now  this asymptotic expansion implies --- under the assumption that $y \geq 10 (1+T^2)$ ---  that

\begin{equation} \label{Yogoanywhereyougo} |F_s(2 \pi | \mu | y)^2| \ll \exp(-\pi  |\mu|  y) \int_{1}^{Y} |F_s(2 \pi | \mu| y)|^2  \frac{dy}{y^3} \ \ (10 (1+T^2)  < y< Y/2)  \end{equation} 
One easily verifies \eqref{Yogoanywhereyougo} for $F = I$ or $F=K$ separately, using their asymptotic expansion. \footnote{
For example, when $F  = K$, the left-hand side is bounded by $e^{-4 \pi |\mu| y} \cdot y^{-1}$, up
to constant factors, whereas the integral on the right-hand side is bounded below
by a constant multiple of $\int_{5 (1+T^2)}^{10 (1+T^2)} e^{-4 \pi |\mu| x}  x^{-4} dx  \gg
e^{- 4 \pi |\mu| (6 (1+T^2))} (1+T^2)^{-4} \gg  e^{-  2.5 \pi |\mu| y} y^{-4} \gg e^{-3  \pi |\mu| y} y^{-1}$. Similarly for $F=I$. }
Now \eqref{Yogoanywhereyougo}  for $F$ an arbitrary linear combination of $K, I$ follows from the
approximate orthogonality of $K, I$, more precisely from the inequality
$ \| a I + b K\|^2_{L^2} \geq  \frac{1}{2}  \left( \| a I \|^2 + \| b K\|^2\right),$
where, on both sides, the norm is on $[|\mu|, |\mu| Y]$ with respect to the measure $dy/y^3$. 
To see this we only need note that the form $\frac{1}{2} \|I\|^2 x^2 +2  \langle I, K \rangle xy 
+ \frac{1}{2} \|K\|^2 y^2$ is positive definite, 
because $\langle  I, K \rangle^2$ is  much smaller than $ \|I\|^2 \cdot \|K\|^2$ so long as $Y$ is large enough.

By Cauchy-Schwarz,  we see that, if $b = \pi/2$ say, 
\begin{multline*} \left|  \sum_{\mu \neq 0} a_{\mu} F_{s}(2 \pi |\mu| y) e^{2 \pi i \langle \mu, (x_1, x_2) \rangle}   \right|^2  \leq    \left| \sum_{\mu \neq 0} a_{\mu}  F_{s}(2 \pi |\mu| y)   e^{-\pi |\mu| y/2} e^{\pi |\mu| y/2}  \right|^2 
\\ \ll  \exp(-by) \sum_{\mu} |a_{\mu}|^2  \int_{A}^Y  F_{s}(2 \pi |\mu| y)|^2  \frac{dy}{y^3}
\\ \leq \exp(-b y) \|f\|_{L^2(M_Y)}.\end{multline*}

To verify~\eqref{precisebound2} we observe that the Fourier expansion of $F$ must consist (for example,
in the case of relative (Dirichlet) boundary conditions) in a linear combination of functions
$G_{s, \mu}(y)$ where
$$G_{s, \mu} := K_s(2 \pi  |\mu|  y) - \frac{K_s(2 \pi   |\mu| Y)}{I_s(2 \pi |\mu| Y)} I_s(2 \pi i  |\mu|  y),$$ this follows
by noting that each Fourier coefficient must individually satisfy Dirichlet boundary conditions,
and so restricts to zero at $y=Y$. 
Similarly, in the case of absolute (Neumann) boundary conditions, 
we obtain a linear combination of similar functions $G_s'(|\mu| y)$ 
with $G_s := K_s(2 \pi i  |\mu| y) - \frac{K_s'(2 \pi i   |\mu|Y)}{I_s'(2 \pi i  |\mu| Y} I_s(2 \pi i   |\mu| y)$.
We proceed in the case of Dirichlet conditions, the Neumann case again being similar: 

Then the asymptotic expansions imply -- again, assuming that  $Y \geq 10 (1+T^2)$, but without the restrictive
assumption that $y < Y/2$ -- that 
\begin{equation} \label{Gsbound} |G_s(|\mu| y)| \leq \exp(-  \pi  |\mu| y) \int_{1}^{Y } |G_s(|\mu|  y)|^2 \frac{dy}{y^3},\end{equation}
hence~~\eqref{precisebound2} follows similarly.  

To check~\eqref{Gsbound},  note that $\left| \frac{K_s(2 \pi i |\mu| Y)}{I_s(2 \pi i |\mu| Y)} I_s(2 \pi i  |\mu| y) \right|
\leq |K_s(2 \pi i  |\mu| Y)| \leq |K_s(2 \pi i |\mu| y)|$, at least for $y$ larger than some absolute constant;
that means, in particular, that $|G_{s, \mu}(y)|$
is bounded by $2 |K_s(2 \pi |\mu| y)|$. But the right-hand side integral exceeds a constant multiple of
the same integral with $G_{s,\mu}$ replaced by
$|K_s(2 \pi i |\mu| y)|^2$, by means of the discussion after \eqref{Yogoanywhereyougo}.  Consequently,~\eqref{Gsbound} results from the fact that,
for $F_s \equiv K_s$,~\eqref{Yogoanywhereyougo} is actually valid for $y < Y$, not merely for $y < Y/2$.

 \qed

We put \index{$\wedge^Y$} \index{truncation} 
$$\wedge^Y f (x)= \begin{cases} f(x), & \height(x) \leq Y \\ f(x) - f_N(x),  & \height(x) > Y.\end{cases} $$

\subsection{Summary of results on Eisenstein series} \label{subsec:einsteinintro} 
We present in summarized form the parameterization of the continuous spectrum
 of the Laplacian on forms and functions on $M$.  This parameterization is in terms
 of the continuous spectrum of the Laplacian on $M_B$. 
  
  We assume familiarity with some of the standard results on Eisenstein series;
  for example, the treatment in Iwaniec \cite{Iwaniec}; although this reference addresses
  hyperbolic two-space and only the case of functions, the generalizations necessary
  for hyperbolic $3$-space and for $i$-forms are quite routine.

\subsubsection{Spaces of functions and forms on \texorpdfstring{$M_B$}{M_B}.}  \label{mbfunction}
We denote by $C^{\infty}(s)$ the space of functions on $M_B$
that are a multiple of $y^{1+s}$ on each component. 
We denote by $\Omega^+(s)$ resp. $\Omega^-(s)$ the space of $1$-forms
on $M_B$ that are a multiple of $y^s (dx_1 +i dx_2)$ resp. $y^{-s} (dx_1 - i dx_2)$
(note the $-s$ in the second definition!) on each component. 

If each cusp has no orbifold points, the dimension of all three spaces equals the number of cusps of $M$.  We denote this number by $h$.  In general the dimensions of $\Omega^{\pm}(s)$
equals the number of {\em relevant} cusps (denoted $\hrel$) and the dimension of $C^{\infty}(s)$ equals the number of cusps:

We sometimes write simply $C^{\infty}$ (resp. $\Omega^+, \Omega^-$)
for these respective spaces when $s=0$. Thus, for $f \in C^{\infty}$, we
have $f \cdot y^s \in C^{\infty}(s)$, and 
our discussion of dimensions says:
$$ \dim \Omega^+= \dim \Omega^- = \hrel, \dim C^{\infty} = h. $$

We introduce inner products on $C^{\infty}, \Omega^+, \Omega^-$ for purely imaginary 
$s$ via the rule 
\begin{equation} \label{cinftyomegaomegaip}  \langle f, g \rangle = \int_{\height(x) = T} \langle f, g \rangle_x , \end{equation} 
where the measure on the set $\height(x) = T$ (this set has an obvious structure of $2$-manifold)  is the measure induced by the hyperbolic metric, and $\langle f, g \rangle_x$ is the inner product on the space of functions/forms at $x$ induced by the hyperbolic metric. Since $s$ is purely imaginary, the resulting inner product is independent of the choice of $T$.

 \subsubsection{Functions}

For $f \in C^{\infty}(0)$ and $s \neq 0$,  the theory of Eisenstein series implies that there exists a unique eigenfunction
$E(s, f)$  of the Laplacian on $M$, with eigenvalue $1-s^2$, and with the property that 
  $$E(s, f) \sim f  \cdot y^{ s} +g \cdot y^{-s},$$
  for some $g \in C^{\infty}(0)$. Indeed, the association $f \mapsto g$ defines  \index{$\Psi(s)$}
  a meromorphic linear operator $\Psi(s): C^{\infty}(0) \rightarrow C^{\infty}(0)$, the so-called {\em scattering matrix}; thus
  $$f \cdot y^s + \left( \Psi(s) f \right)  y^{-s}$$ is the asymptotic part of an eigenfunction. 
 The uniqueness implies that
$\Psi(s) \Psi(-s) = \mathrm{id}.$  \footnote{These statements remain valid in slightly modified form at $s=0$; we leave the formulation to the reader.} We will write $\psi(s) := \det \Psi(s)$ for the determinant of $\Psi$.

  Moreover, 
 choosing an orthonormal basis $\mathcal{B} = \{ f_1, \dots, f_h\}$
for $C^{\infty}(0)$, the functions $E(s, f)$ for $s \in i \R_{\geq 0}$
span the continuous spectrum of the Laplacian ---  that is to say, for any square integrable
$1$-form $F$, we have the relation 
\begin{equation} \label{specexp1}  \|F\|^2 = \sum_{j} \langle F, \psi_j \rangle^2 + \frac{1}{2\pi} \sum_{f \in \mathcal{B}} \int_{t =0}^{\infty} 
|\langle F, E(f,it) \rangle|^2,\end{equation}
where $\{\psi_j\}$ is an orthonormal basis for the discrete spectrum
of the Laplacian on functions.

The so-called Maass-Selberg relations assert that\footnote{The ``short'' mnemonic is that the right-hand side is the regularized
integral $$\int_{\height \geq Y}^{\mathrm{reg}} - E(f,s)_N \cdot \overline{E(f', t)_N}$$
} for $s, t \in i \R$, 
\begin{eqnarray}  \label{MSForms} \langle \wedge^Y E(f, s),  \wedge^Y E(f', t) \rangle &=&   \langle f, f' \rangle \frac{Y^{s-t}}{s-t} + \langle \Psi(s) f, \Psi(t) f' \rangle \frac{Y^{t-s}}{t-s}  \\   \nonumber\nonumber &+& \langle \Psi(s) f, f' \rangle \frac{Y^{-s-t}}{-s-t} + \langle f, \Psi(t) f' \rangle \frac{Y^{s+t}}{s+t}.  \end{eqnarray}

In particular, taking the limit as $s \rightarrow t$ (and still assuming $s$ is purely imaginary), 
\begin{eqnarray}  \nonumber  \| \wedge^Y E(f, s)  \|^2 &=& 2 \log Y \langle f, f \rangle   - \langle \Psi(s)^{-1} \Psi'(s) f, f \rangle \\ &+&  \label{Second} \langle f, \Psi(s) f \rangle \frac{Y^{2s}}{2s} - \langle \Psi(s) f, f \rangle \frac{Y^{-2s}}{2s}
 \end{eqnarray}
 where $\Psi'(s) := \frac{d}{ds} \Psi(s)$; this implies in particular that  
 \begin{equation} \label{psibb} \mbox{ $-\Psi(s)^{-1} \Psi'(s)$ is bounded below, \ \ $s \in i\R$, }\end{equation}  i.e.
 $ - \Psi(s)^{-1} \Psi'(s) + t I$ is positive definite for some $t > 0$ and {\em all $s$.}
 
 The following observation will be of use later: if $h(s)$ is an analytic function
 satisfying $h(-s) = h(s)$, 
 and $J(s) =  \langle f, \Psi(s) f \rangle \frac{Y^{2s}}{2s} - \langle \Psi(s) f, f \rangle \frac{Y^{-2s}}{2s}$ is the second term of~\eqref{Second}, then
 $$ \int_{t \in \R} h(t) J(it) dt- \pi h(0) \langle \Psi(0) f, f \rangle =  \mbox{ decaying as $Y \rightarrow \infty$}$$ 
as one sees   by shifting contours. \footnote{Write $A(s) = \langle \Psi(-s) f, f \rangle Y^{2s}$, so that
$J(it) = \frac{A(it) - A(-it)}{2 it}$.  Shifting contours, 
$$ \int_{t \in \R} h(t) J(it) = \int_{\mathrm{Im}(t) = \delta} h(t) J(it) 
= \int_{\mathrm{Im}=\delta}  \frac{h(t) A(it) }{2 it} dt 
+ \int_{\mathrm{Im}=-\delta} \frac{h(t) A(it)}{2it} dt.$$
But $A(it)$ decays rapidly with $Y$ when $\mathrm{Im}(t)  > 0$. 
Shift the second term to $\mathrm{Im}=\delta$ also, leaving a residue of $\pi h(0) A(0)$.}

 \subsubsection{Forms}   \label{sss:EisForms} 
For each $\omega \in \Omega^{+}(0)$,  \index{$\Phi^+(s)$} \index{$\Phi^-(s)$}
there exists unique $\omega' \in \Omega^{-}(0)$ with the property that
there is an eigenfunction of the Laplacian on $1$-forms, with eigenvalue $-s^2$, and  asymptotic to
\begin{equation} \label{omegaomega} \omega  \cdot y^s + \omega'  \cdot y^{-s}. \end{equation}
 
The association $\omega \mapsto \omega'$ defines 
a meromorphic map $\Phi^{+}(s) : \Omega^+(0) \rightarrow \Omega^{-}(0)$. 
The operator $\Phi^+(s)$, for $s \in i \R$, is {\em unitary} for the unitary structure previously defined.
The inverse of $\Phi^{+}(s)$ is denoted $\Phi^{-}(s)$. We denote the unique differential
form asymptotic to ~\eqref{omegaomega} by $E(\omega, s)$
or $E(\omega', s)$:
\begin{align} E(\omega, s)  &\sim& \omega &\cdot  & y^s & +& \Phi^{+}(s) & \omega \cdot y^{-s} \ \ 
&(\omega \in \Omega^+(0)),  \\
E(\omega', s) & \sim & \omega'& \cdot & y^{-s} &+& \Phi^{-}(s) & \omega' \cdot y^{s} \  \ & (\omega \in \Omega^{-}(0)) .\end{align}
Note the funny normalization of signs! In particular, $\Phi^{-}(s) \Phi^{+}(s)= 1$.

Fixing an orthonormal basis $\omega_1, \dots, \omega_r$ for $\Omega^{+}(s)$, \index{coclosed}
the functions $E(\omega_i, s)$ form an orthonormal basis for {\em co-closed}
$1$-forms on $Y(K)$, i.e., forms satisfying $d^* \omega = 0$;
together, $E(\omega_i, s)$ and $d E(f,s)$ form an orthonormal basis for $1$-forms, that is to say:
$$ \| \nu\|^2 = \sum_{\psi}| \langle \nu, \psi \rangle|^2 + \int_{0}^{\infty} \sum_{f} |\langle \nu, \frac{d E (f,it)}{\sqrt{1+t^2}} \rangle|^2 \frac{ dt}{2 \pi}
+ \sum_{\omega} \frac{1}{2 \pi} \int_{-\infty}^{\infty} | \langle \nu, E(\omega, it) \rangle|^2 dt,$$
where $\psi$ ranges over an orthonormal basis for the spectrum of the form Laplacian
acting on  $1$-forms that belong to $L^2$.

The Maass-Selberg relations assert that
for $\omega_1, \omega_2 \in \Omega^+$ and $s,t \in i \R$, 
\begin{eqnarray} \langle \wedge^Y \mathrm{Eis}(s,\omega_1), \wedge^Y \mathrm{Eis}(t,\omega_2) \rangle \\
= \langle \omega_1 , \omega_2 \rangle \cdot (\frac{Y^{s-t}}{s-t}  )  + \langle \Phi^{+}(t)^{-1} \Phi^{+}(s) \omega_1, \omega_2 \rangle \frac{{Y}^{t-s}}{t-s}.
 \end{eqnarray}
 (one can remember this using the same mnemonic as before, but it is a little bit simpler.) 
 
Taking the limit as $t \rightarrow s$, and a sum over an orthonormal basis
for $\Omega^+)(s)$, we obtain 
\begin{equation} \label{MSR1F} \sum_{\omega} \| \wedge^Y E(s, \omega) \|^2 =  (2 \hrel  \log Y)  -  \mathrm{trace} ( \Phi(s)^{-1} \Phi'(s)).\end{equation}  Note that this reasoning also implies that, as in the previous case, 
\begin{quote} \label{derivbounded} $-\frac{\Phi'(s)}{\Phi(s)}$
is bounded below, \end{quote} i.e. $a I -\frac{\Phi'(s)}{\Phi(s)}$
is positive definite for all $s \in i \R$ and suitable $a > 0$.

\subsubsection{\texorpdfstring{$2$}{2}-forms and \texorpdfstring{$3$}{3}-forms; the 
Hodge \texorpdfstring{$*$}{*}.}

The corresponding spectral decompositions for $2$- and $3$-forms follow by applying the * operator. For the reader's convenience we note that, on $\H^3$:
 
\begin{eqnarray*} *dx_1 = \frac{dx_2 \wedge dy}{y},  \ dx_2  = \frac{dy \wedge  dx_1}{y}, \ * dy = \frac{dx_1 \wedge dx_2}{y},  \\ 
 * (dx_2 \wedge dy) = y dx_1, \  * (dy \wedge dx_1) = y dx_2, \ *(dx_1 \wedge dx_2) = y dy 
 \end{eqnarray*}

 \subsubsection{Residues} \label{subsubsec:residues}
 The function $\Psi(s) : C^{\infty} \rightarrow C^{\infty}$ (and, similarly, $s \mapsto E(s,f)$) has a simple pole at $s=1$.
  We shall show that:
  \begin{quote} The residue $R$ of $\Psi(s)$ at $s=1$  is, in a suitable basis,
  a diagonal matrix, with nonzero entries  
$\area(\partial N)/\vol(N)$, where $N$ varies through connected components of $M$. 
\end{quote}

\begin{proof}   
For any $f \in C^{\infty}$, the residue of $s \mapsto E(s,f)$ at $s=1$
is a function on $M$, constant on each connected component; we denote
it by $c_f$. We also think of it as a function on cusps (namely, its constant value
on the corresponding connected component).   Alternately, 
$R(f)$ is a multiple of $y$ on each connected component of $M_B$, 
and $R(f)/y$ coincides with $c_f$ (considered as a function on cusps). 
 Again, we refer to \cite{Iwaniec} for a treatment of the corresponding facts
in hyperbolic $2$-space, the proofs here being identical.

 A computation with the Maass-Selberg relations shows that, if $f =  y \cdot 1_{\mathcal{C}_i}$ (which is to say: the function which equals $y$ on the component of $M_B$ corresponding
 to $\mathcal{C}_i$, and zero on all other components) then $c_f$ takes the value
\begin{equation} \label{cfvalue} \frac{\mathrm{area}(\mathcal{C}_i)}{\vol(\mathcal{C}_i)} \end{equation} 
on the connected component of $\mathcal{C}_i$, and zero on all other components. 
 
 Indeed, it is enough to check the case of $M$ connected.   The Maass-Selberg relations show that 
 $ \| c_f \|_{L^2(M)}^2 = \langle   R(f) , f \rangle_{C^{\infty}}.$
 In particular, if we write $m_i$ for the value of $c_f$ when $f$ is the
 characterstic function of the cusp $\mathcal{C}_i$, we have, for any scalars $x_i$, 
 $  \vol(M)  \left| \sum x_i m_i  \right|^2 =  \sum x_i m_i  \overline{\sum x_i}  \area(\mathcal{C}_i) $,
 whence \eqref{cfvalue}.   \end{proof} 

Observe for later reference that  the kernel of the map $f \mapsto c_f$ (equivalently $f \mapsto R(f)$) consists 
    of all $f = \sum x_i (y 1_{\mathcal{C}_i})$ such that
\begin{equation} \label{reszero} \sum x_i \area(\mathcal{C}_i) = 0,\end{equation}
    where the $\mathcal{C}_i$ vary through the cusps of any connected component of $M$.
    
  \section{Reidemeister and analytic torsion}  \label{sec:rtatnc}
We define carefully the regulator, Reidemeister torsion, and analytic torsion in the non-compact case.  
 We shall first need to discuss harmonic forms of polynomial growth and how to equip them with an inner product.  It is worth noting that 
for our purposes, when we are comparing two manifolds,  the exact definition is not so important so long as one is consistent.

Let $M$ be any hyperbolic $3$-manifold of finite volume, equipped
with a height function $\height$ as in~\S~\ref{section:split}. Our main interest
is in the case $M = Y(K)$, equipped with $\height$ as in~\S~\ref{infgeom}.

\subsection{The inner product on forms of polynomial growth.} \label{ip:polygrowth}

We will define below a specific space $\mathscr{H}^j$ of harmonic forms $(0 \leq j \leq 2$) uch that: 
 \begin{quote} The  natural map $\mathscr{H}^j \rightarrow H^j(M, \C)$ is an isomorphism. 
\end{quote}

\begin{itemize} \item For $j=0$ we take $\mathscr{H}^j$ to consist of constant functions.
\item For $j=1$, we take $\mathscr{H}^j$ to consist of cuspidal harmonic forms together
with  forms of the type $\Eis(\omega, 0)$, where $\omega \in \Omega^+(0)$. 

\item For $j=2$, we take $\mathscr{H}^j$ to consist of cuspidal harmonic $2$-forms
together with forms of the type $* d E(f,1)$ where $f \in C^{\infty}(0)$
is such that $s \mapsto E(f,s)$ is holomorphic at $s=1$ (equivalently: 
$f$ lies in the kernel of the residue of $\Psi(s)$ at $s=1$). 
\end{itemize} 

One verifies by  direct computation
that the map to $H^j(M, \C)$ from this subspace is an isomorphism.      
For $j=1, 2$ we have a canonical decomposition 
 $$\mathscr{H} =  \Omega^j_{\cusp}(M) \oplus \Omega^j_{\Eis}(M), $$ where 
$\Omega^j_{\cusp}(M)$ is the space of cuspidal $j$-forms (i.e. $\omega_N \equiv 0$),  and $\Omega^j_{\Eis}(M)$ is the
 orthogonal complement\footnote{This makes sense even though
 elements of $\mathscr{H}$ are not square integrable, because elements in $\Omega^j_{\cusp}$ are of rapid decay, and therefore can be integrated against a function of polynomial growth.} 
of $\Omega^j_{\cusp}$.
 Moreover, $\Omega^j_{\cusp}(M)$ maps isomorphically,
 under $\mathscr{H} \stackrel{\sim}{\rightarrow} H^j(M, \C)$, to the cuspidal cohomology
 $H^j_{!}(M, \C)$.

\subsubsection{An inner product on  \texorpdfstring{$\mathscr{H}^1$}{H^1}.}

We introduce on $\Omega^1_{\Eis}(M)$ the inner product 
\begin{equation} \label{constant-term-norm} \| \omega\|^2 =\lim_{Y}   \frac{ \int_{M_{\leq Y}} \langle \omega, \omega \rangle}{\log Y}. \end{equation}
where the inner product $\langle \omega, \omega \rangle$ at each point of $M$
is taken using the Riemannian structure. 
In explicit terms, if $\sigma_i^* \omega = a_i dx_1 + b_i dx_2$, 
then the norm of $\sigma$ equals $\sum_i \area(\mathcal{C}_i) \cdot (|a_i|^2 +| b_i|^2)$. 
 \footnote{An equivalent definition is 
\begin{equation} \label{constant-term-norm2} \| \omega \|^2 = \int_{\partial  M_{\leq Y}}
\omega_N \wedge \overline{\omega_N}, \end{equation}
this holding true for any $Y$ large enough. 
Here, if $\sigma_i^* \omega_N = a dx_1 + b dx_2$, then $\overline{\omega_N}$ is defined
to be the differential form on $M_B$ so 
 $\sigma_i^* \overline{\omega_N} = \bar{b} dx_1 - \bar{a} dx_2$ (this looks more familiar in coordinates $x_1 \pm i x_2$!) . }
This norm has fairly good invariance properties, e.g., the Hecke operators act self-adjointly, etc.  

We now endow 
$$ \mathscr{H}^1 \stackrel{\sim}{\rightarrow} \Omega^1_{\cusp} \oplus \Omega^1_{\Eis}$$
with the direct sum of the two inner products defined. 
\subsubsection{An inner product on \texorpdfstring{$\mathscr{H}^2$}{H^2}.}    For the definition of Reidemeister torsion  later we will use only $H_{2!}$, 
 in fact,  but the general setup will be of use later.

   Unike $\mathscr{H}^1$, there is no ``good'' inner product on $\Omega^2_{\Eis}$ (in a suitable sense,  $\mathscr{H}^1$ grow only logarithmically at $\infty$, but this is not so for $\mathscr{H}^2$).  
 We make the more or less arbitrary definition:
\begin{equation} \label{h2normdef}  \|\omega\|^2 =  \lim_{Y \rightarrow \infty} Y^{-2} \int_{M_{\leq Y}} \langle \omega, \omega \rangle,  \ \ (\omega \in \Omega^2_{\Eis})
\end{equation} 
 which, although analogous to the prior definition, 
has no good invariance properties ---   for exaple, the Hecke operators
 do not act self-adjointly. 
Again we endow 
$$ \mathscr{H}^2 \stackrel{\sim}{\rightarrow} \Omega^2_{\cusp} \oplus \Omega^2_{\Eis}$$
with the direct sum of the standard inner product on the first factor, and the inner product just defined, on the second factor.

The inner product on $\Omega^2_{\Eis}$ is simple to compute explicitly:
 Supposing for simplicity that $M$ is connected (the general case follows component by component). 
Let $\mathcal{C}_1, \dots, \mathcal{C}_h$ be the cusps of $M$.
For $\mathbf{u} = (u_1, \dots, u_h) \in \C^h$ such that $\sum u_i = 0$, 
  let $\omega \in \Omega^2_{\Eis}$ be such that
 $\int_{\mathcal{C}_i} \omega =    u_i$.  Then \begin{equation} \label{badnormdef} \| \omega\|^2 = \frac{1}{2}  \sum_{i} \area(\mathcal{C}_i)^{-1} |u_i|^2 \end{equation}  
Indeed $\sigma_i^* \omega = \frac{u_i}{\mathrm{area}(\mathcal{C}_i)} dx_1 \wedge dx_2 
+ \mbox{ rapidly decaying},$
from where we deduce~\eqref{badnormdef}.

 \subsection{Definition of Reidemeister torsion}
 \label{subsec:rtnonsplitdef}

\index{regulator (noncompact case)}

We {\em define}
\begin{equation} \label{regncdef}  \reg(M) = \frac{\reg(H_1) \reg(H_{3,\bm})}{\reg(H_0) \reg(H_{2!})} \end{equation} 
where
$$\reg(H_{i,?}) =    \left| \det \int_{\gamma_i} \omega_j  \right|  ,$$
and
\begin{itemize}
\item[-] (for $i \in \{0,1\}$) the elements $\gamma_i \in H_i(M, \Z)$ are chosen to project to a basis for $H_i(M, \Z)_{\tf}$, and the $\omega_j$
are an orthonormal basis for $\mathscr{H}^i$;
\item[-] (for $i=2$)   as above, but replace  $H_i(M, \Z)_{\tf}$ by
 the torsion-free quotient of $H_{2,!}(M, \Z)$, 
and the $\omega_j$ are an orthonormal basis for the space of {\em cuspidal} harmonic $2$-forms.   
\item[-] (for $i=3$) 
as above, but replace  $H_i(M, \Z)_{\tf}$ by
 the torsion-free quotient of $H^{\bm}_{3}(M, \Z)$,
 and the $\omega_j$ are an orthonormal basis for the space of  harmonic $3$-forms (= multiples of the volume form on each component). 
\end{itemize}

See also after~\eqref{h2h2shriek} for some further discussion. Our ugly definition here pays off in cleaner theorems later.   

If $M$ is compact this specializes to~\eqref{regsimpledef}. Note that \begin{itemize}
\item[(i)] 
 $\reg(H_0) = \frac{1}{\sqrt{\vol(M)}}  $ and   $\reg(H_3) = \sqrt{\vol(M)}$. 
 We follow our convention that $\vol$ denotes the product of volumes of all connected components. 
 
\item[(ii)]If $M$ were compact, $\reg(H_1) \reg(H_{2,!}) =1$, which recovers our earlier definition 
\end{itemize}

 The definition of Reidemeister torsion given in~\S~\ref{sec:CMthm} -- see equation~\eqref{rtsimpledef} -- now carries over to the non-compact case also, namely:
\begin{eqnarray} \RT(M) :=  |H_1(M, \Z)_{\tors}|^{-1}    \reg(M) \\
=  |H_1(M, \Z)_{\tors}|^{-1} \cdot  \vol(M) \cdot \frac{\reg(H_{1})}{\reg(H_{2!})}.
\end{eqnarray}
As always, $\vol$ here denotes the product of the volumes of all connected components.

 \subsection{Explication of the \texorpdfstring{$H_2$}{H_2}-regulator}
In this section, we relate $\reg(H_2)$ and $\reg(H_{2,!})$, which will be necessary later. 

Begin with the sequence
 $$H^{\bm}_{3}(M, \Z) \rightarrow H_2(\partial M, \Z) \rightarrow H_2(M, \Z) \twoheadrightarrow H_{2, !}(M, \Z),  $$
 where as before $H_{2,!}$ is the image of $H_2$ in $H^{\bm}_{2}$. All of the groups
 above are torsion-free, at least away from any orbifold prime for $M$.  
 In what follows we suppose there are no orbifold primes; if this is not so, 
 our computation must be adjusted by a rational number supported only at those primes. 
  
Fix generators\footnote{For example,
if $M$ is connected without orbifold points and we enumerate the cusps $\mathcal{C}_1, \dots, \mathcal{C}_h$, 
we may take $\delta_i = [\mathcal{C}_{i+1}] - [\mathcal{C}_1]$, where $[\mathcal{C}]$
denotes the fundamental class of a torus cross-section of $\mathcal{C}$.}
 $\delta_1, \dots, \delta_{h}$ for the free group $H_2(\partial M, \Z)/\mathrm{image}(H^{\bm}_{3}(\Z))$.  
 Fix elements $\gamma_1, \dots, \gamma_r \in H_2(M,\Z)$
 whose images span $H_{2, !}$. Then: 
 
 \begin{eqnarray} H_2(M , \Z) &=& \bigoplus_{j} \Z \gamma_j \oplus  \bigoplus_{i} \Z  \delta_i,\end{eqnarray}
where we have identified the $\delta_i$ with their image in $H_2(\partial M, \Z)$. 

 Fix orthonormal bases $\omega_1, \dots, \omega_r$ for $\Omega^2(M)_{\cusp}$
 and $\eta_1, \dots, \eta_h$ for $ \Omega^2(M)_{\Eis}$. 
 Note that the {\em $\omega_i$s are orthogonal to the $\delta_j$}, that is to say
 $$ \int_{\delta_j} \omega_i \equiv 0, $$ 
 because the $\omega_i$s are cuspidal and so decay in each cusp.

 Therefore, the ``matrix of periods'' with respect to these bases is block-diagonal, and so:
\begin{eqnarray} \label{h2h2shriek} \reg(H_2) & =& \left|  \det  \left( \int_{\gamma_i} \omega_j  \right)  \cdot \det \left( \int_{\delta_i} \eta_j \right)  \right| \\ \nonumber
 & =& \reg(H_{2,!}) \cdot  \left| \det \left( \int_{\delta_i} \eta_j \right) \right|.\end{eqnarray} 
where we have used the definition:  $\reg(H_{2, !}) = \det  \left( \int_{\gamma_i} \omega_j  \right) $.

The latter factor $ \det \left( \int_{\delta_i} \eta_j \right)$ is easy to compute from~\eqref{badnormdef}:
if $M$ is connected with cusps $\mathcal{C}_i$,  then there is an equality
$$ \det \left( \int_{\delta_i} \eta_j \right)^{-2} = \prod  (2 \area(\mathcal{C}_i))^{-1} \left(\sum 2 \area(\mathcal{C}_i)  \right).$$
The computation on the right hand side is deduced from~\eqref{badnormdef}: the squared covolume of $\{(x_1, \dots, x_k) \in \Z^k: \sum x_i = 0\}$
in the metric $\sum a_i x_i^2$ is given by $\prod a_i \cdot \sum a_i^{-1}$;
the inverse on the left arises because this is computing a covolume for {\em cohomology}.

More generally, if $N_1, \dots, N_k$ are the connected components of $M$,  we deduce $  \det \left( \int_{\delta_i} \eta_j \right)^2$ equals $  \frac{ \left( \prod_{\mathcal{C}} 2 \area(\mathcal{C}) \right) }{ \prod_{N} 2\area(\partial N) },$
which, in view of the discussion after~\eqref{reszero}, can also be written:
\begin{equation} \label{h2reghelp}   \det \left( \int_{\delta_i} \eta_j \right)^2  =  \frac{ \left( \prod_{\mathcal{C}} 2 \area(\mathcal{C} )\right)   }{\vol(M)}
\cdot \detprime(2 R)^{-1}, \end{equation}
where $R$ is the residue at $s=2$ of $\Psi(s)$, and $\det'$ denotes the product of nonzero eigenvalues.  We have followed the convention that $\vol$ denotes the products of volumes of connected components.

 Although~\eqref{h2reghelp} looks unwieldy, it will be convenient when comparing
the left-hand side between two different manifolds.

\subsection{Regularized analytic torsion} \label{subsec:locality}
We now present the definition of the analytic torsion of a non-compact hyperbolic $3$-manifold.
 
Notation as previous, 
 the operator $e^{-t \Delta}$, and its
analogue for $i$-forms, are not, in general, trace-class. Nonetheless there is
a fairly natural way to regularize their traces, which we now describe: 

Let $K(t; x,y)$ be the integral kernel of $e^{-t \Delta}$ acting on functions, i.e. the heat kernel, and let
$$k_t(x) = K(t; x, x).$$  
(In the case of forms, $K(t; x, x)$ is an endomorphism of the
space of forms at $x$, and we take the trace.) We also set 
$ \wedge^Y k_t(x)$ to be $  \wedge^T_x K(t; x, y)$ (the operator $\wedge^Y$
acts in the first variable $x$) specialized to $x=y$.  By definition,
$$ \wedge^Y k_t(x) = k_t(x) \mbox { if $\height(x) \leq Y$},$$
and also $\wedge^Y k_t(x)$ is of rapid decay: If we fix $t$,  the quantity
$\sup_{\height(x) > Y} \wedge^Y k_t(x)$ decays faster than any polynomial in $Y$. 
In particular,

$$ \int_{M_{\leq Y}} k_t(x) = \int_{M} \wedge^Y k_t(x) + \left( \mbox{  rapidly decaying in $Y$ } \right).$$ 
In fact,  either side   is asymptotic to a linear polynomial in $\log Y$.
This is a consequence of the Selberg trace formula, for example;
see~\cite{Friedman}.
 Assuming this, we may define the regularized trace: 

\begin{df}  \label{RegTraceDef} 
The regularized trace
$\mathrm{tr}^*(e^{-t\Delta})$ is  the constant term of the unique  linear function
of  $\log Y$ asymptotic to  $\int_{M} \wedge^Y k_t$, i.e. it is characterized by the property 
\begin{equation} \label{regdef} \int_M \wedge^Y k_t(x) dx  \sim k_0\log Y + \mathrm{tr}^*(e^{-t \Delta}) ,\end{equation}
Here the notation $A(Y) \sim B(Y)$ means that $A(Y) - B(Y) \rightarrow 0$
as $Y \rightarrow \infty$. 
\end{df}

 Of course this definition of $\mathrm{tr}^*$ depends on the choice of height function. 
We can be more precise about $k_0 = k_0(t)$: e.g., in the case of $0$-forms, 
 $k_0 = \frac{h}{2\pi} \int_{t=0}^{\infty} e^{- t (1+s^2)},$
where $h$ is the number of cusps.\footnote{
This computation of $k_0$ follows easily from the  spectral expansion  
$$k_t(x) = \sum e^{-t \lambda} |\psi_{\lambda}(x)|^2 + \frac{1}{2\pi}  \sum_{f \in \mathcal{B}} \int_{t=0}^{\infty} e^{-t (1+s^2)} |E(f, it)(x)|^2, $$
see \eqref{specexp1}.}

The Selberg trace formula may be used to compute $\mathrm{tr}^*(e^{-t \Delta})$ (see~\cite{Friedman} for the case of functions); in particular, it has an asymptotic expansion near $t=0$:
\begin{equation} \label{asymptoticzero} \mathrm{tr}^*(e^{-t \Delta}) \sim a t^{-3/2} + b t^{-1/2} + c t^{-1/2} \log(t) + d  + O(t^{1/2}), \end{equation}
We do not need the explicit values of the constants $a, b, c, d$, although they are computable. 
We define the regularized determinant of $\Delta$ as usual, but replacing the trace of $e^{-t\Delta}$
by $\mathrm{tr}^*$, i.e.

\begin{eqnarray}  \label{atsecondef} \log \det{} ^* (\Delta) =  - \frac{d}{ds} \Big|_{s=0} \left( \Gamma(s)^{-1} \int_0^{\infty}( H(t)-H(\infty ))  t^s \frac{dt}{t}\right), \\ \nonumber
H(t) := \tr^*(e^{-t \Delta}), H(\infty) = \lim_{t \rightarrow \infty}  H(t)  \end{eqnarray}

 This definition requires a certain amount of interpretation: Split the integral into $I_1(s) = \int_{0}^{1}$ and $I_2(s)= \int_{1}^{\infty}$. The first term is absolutely convergent for $\Re(s) $ large and (by~\eqref{asymptoticzero})
admits a meromorphic continuation to $\Re(s) > -1/2$ with at worst a simple pole at $s=0$; in particular,
$\Gamma(s)^{-1} I_1(s)$ is regular at $s=0$ and so $\frac{d}{ds} \big|_{s=0} \Gamma(s)^{-1} I_1(s)$ is defined.  The second term is absolutely convergent
for $\Re(s) < 1/2$ (as follows, for example, from~\eqref{ASspec} below and elementary estimates) and indeed $\frac{d}{ds} \big|_{s=0}\Gamma(s)^{-1} I_2(s) = \int_{1}^{\infty}
\frac{H(t)-H(\infty)}{t} dt$.   We interpret~\eqref{atsecondef} as {\em defining} $\log \det{}^* \Delta$ as the sum of these quantities.

This leads us to the definition of the regularized analytic torsion:
$$\log \analT(Y) = \frac{1}{2} \sum_{j} (-1)^{j+1} j \log \det{}^* (\Delta_j),$$
where $\Delta_j$ denotes the regularized determinant on $j$-forms.

On the other hand, on the spectral side, we may write (for instance, in the case of functions, that is to say, $0$-forms):

\begin{equation} \label{ASspec} H(t) - H(\infty)=  \sum_{\lambda_j \neq 0}  e^{-t \lambda_j} + \frac{1}{2\pi} \int_{u=0}^{\infty} -\frac{\psi'}{\psi}(iu)  e^{-t (1+u^2) } du + \frac{1}{4} e^{-t} \mathrm{trace} \Psi(0), \end{equation}
where $\psi(s)$ is, as before,  the determinant of the scattering matrix, and $\lambda_j$ are the nonzero eigenvalues of the Laplacian on functions.    On the other hand one has
\begin{equation} \label{ASspec2} H(\infty) = b_0,\end{equation} where $b_0$ is equal to the number of zero eigenvalues of the Laplacian 
in the discrete spectrum, i.e.,  the number of connected components of $M$.
Both~\eqref{ASspec} and~\eqref{ASspec2} follow from~\eqref{Second}.   
They mean, in particular, that $H(t)-H(\infty)$ decays exponentially  at $\infty$
to handle convergence issues in the discussion after \eqref{atsecondef}. 
The case of $1$-forms is slightly different: because $0$ lies in the continuous spectrum, the decay of $H(t)- H(\infty)$ is only as $t^{-1/2}$, but that is still adequate for convergence. 

\subsection{Comparison of regularized trace for two hyperbolic manifolds} \label{subsec:etdcomparison}
 Suppose $M, M'$ are a pair of hyperbolic manifolds with cusps $\mathcal{C}, \mathcal{C}'$. Suppose that there is an  isometry $\sigma: \mathcal{C} \rightarrow \mathcal{C}'$.
 
 Then one can compute the difference 
  $\tr^* (e^{-t \Delta_M}) - \tr^*(e^{-t \Delta_{M'}})$ without taking a limit, which will be useful later. 
 Indeed we show that 
\begin{eqnarray} \label{tdeltadiff} 
 \tr^* (e^{-t \Delta_M}) - \tr^*(e^{-t \Delta_M}) &=& \left( \int_{M-\mathcal{C}} k_t(x) - \int_{M' - \mathcal{C'}} k_t'(x)  \right)  \\  \nonumber &+& 
 \left( \int_{\mathcal{C}} k_t(x) - k_t'(\sigma(x)) \right),\end{eqnarray} 
 where all integrals are absolutely convergent.  This holds for the Laplacian on $j$-forms, for any $j$; we do not include
 $j$ explicitly in our notation below.

To verify~\eqref{tdeltadiff} let $\mathcal{C}_Y = M_{\leq Y} \cap \mathcal{C}$. We may write
  
  $$ \tr^* (e^{-t \Delta_M}) = \lim_{Y \rightarrow \infty}  
  \left( \int_{M-\mathcal{C}} k_t(x) + \int_{\mathcal{C}_Y} k_t(x)  -  k_0(t) \log Y\right), $$
  and similarly for $M'$.   Here $k_0(t)$ is as discussed after  Definition~\ref{RegTraceDef}. The coefficients of $\log Y$, namely, $k_0(t)$ and $k_0'(t)$,
  are identical for $M$ and $M'$; they 
  depend only on the cusps, and were given after~\eqref{regdef}. Subtracting yields  the result once we verify that 
$$\lim_{Y \rightarrow \infty} \left(  \int_{\mathcal{C}_Y} k_t(x)  - \int_{\mathcal{C}_Y} k_t'(x)  \right)=    \int_{\mathcal{C}} k_t(x) - k_t'(\sigma(x)) .$$
But it is easy to verify, by the ``locality'' of the heat kernel, that $k_t - k_t' \circ \sigma$ is of rapid decay on $\mathcal{C}$; in particular it is integrable.  That completes the proof of~\eqref{tdeltadiff}.

\section{Noncompact arithmetic manifolds} \label{sec:arithmeticmanifolds} 
 Thus far our considerations applied to any hyperbolic $3$-manifold.
In this section we specialize to the case of an arithmetic, noncompact $3$-manifold. 

Thus $\G = \PGL_2$ and $F$ is an imaginary quadratic field. 
Recall (\S~\ref{ss:splitnotn}) that $\B$ is the upper triangular Borel subgroup of $\G = \PGL_2$
and $K_{\max}$ is a certain maximal compact subgroup.

 \subsection{The case of adelic quotients} \label{infgeom}
 
 We specify a height function for the manifolds $M = Y(K)$ and then 
  give adelic descriptions of the objects $M_B, M_{\leq T}$, notation as in \S \ref{subsec:trunc}.

 Our presentation is quite terse. We refer to \cite[page 46]{Harder} for a more complete discussion.  \index{$\mathrm{Ht}$}
Define the {\em height function} $\mathrm{Ht}: \G(\adele) \rightarrow \R_{>0}$
 via $$ \mathrm{Ht}(bk)^2= |\alpha(b)|$$ for $b \in \B(\adele), k \in K_{\max}$, 
 where $\alpha: \B \rightarrow \mathbb{G}_m$
 is the positive root,  and $|\cdot|$ is the adelic absolute value. We denote the set $\height^{-1}(1)$ by $\G(\adele)^{(1)}$. {\em Warning ---  } this is nonstandard usage; also, the square on the left-hand side arises because of the difference between the ``usual'' absolute value on $\C$
 and the ``normalized'' absolute value, cf. \ \S \ref{ss:notn:F}.
  \index{$\G(\adele)^{(1)}$}

Recall from~\S~\ref{YMdef} 
the definition of the arithmetic manifold 
\begin{equation} \label{YKredefined} Y(K) := \G(F) \backslash \H^3 \times \G(\Afinite)  /K.\end{equation} 
We define the height $\mathrm{Ht}(x)$ of $x \in Y(K)$ to be the maximal height
of any lift $g_x \in \G(\Adele)$ (where we think of the right-hand side of~\eqref{YKredefined} as the quotient $\G(F) \backslash \G(\adele)/K_{\infty} K$). Reduction theory implies this is indeed a height function in the sense of~\S~\ref{subsec:trunc}; more precisely, there exists a  unique height function $h_0$
in the sense of~\S~\ref{subsec:trunc} such that $\height(x) = h_0(x)$ so long as either side is sufficiently large. 
In the notation of~\S~\ref{subsec:trunc} with $M = Y(K)$, the manifold $M_B$ may be described as
 $$  Y_B(K) := \B(F) \backslash   \H^3 \times \G(\Afinite)/   K $$
 The function $\height$ descends (via the identification $\G(\C)/K_{\infty} \cong \H^3$)
 to a proper function $\height: Y_B(K) \rightarrow \mathbf{R}_{>0}$.  
    We write $Y_B(K)_{\geq T}$ for the preimage
 of $(T, \infty]$ under $\mathrm{Ht}$. Then, for sufficiently large $T$,  the natural projection 
 \begin{equation} \label{projbY} Y_B(K)_{\geq T} \rightarrow Y(K)\end{equation} 
 is a homeomorphism onto its image, and the complement of its image is compact. 
 Fix such a $T$ ---  call it $T_0$ ---  and define $\height: Y(K) \rightarrow [1, \infty]$
 to be $T_0$ off the image of~\eqref{projbY}, and to coincide with $\height$ on the image.

 \subsection{Homology and Borel-Moore homology} \label{subsec:arithhomologyBM}

Let $H_c$ denote compactly supported
cohomology.   As on page \pageref{boundarypage}, 
there is a long exact sequence as follows, for any (constant) coefficients: 
$$H^{i-1}(\partial Y(K),- ) \rightarrow   H^i_c(Y(K),- ) \rightarrow 
H^i(Y(K) ,-) \rightarrow  H^{i}(\partial Y(K),- ) \ldots,$$
where we set $\partial Y = \partial Y_{\leq T}$ for any $T$; the homotopy class although not the isometry class
is independent of $T$. Indeed, the inclusion $\partial Y_{\leq T} \hookrightarrow Y_B$
is a homotopy equivalence.

\medskip

As before, we focus mainly on the case $K = K_{\Sigma}$, when we abbreviate  the cohomology groups  $H^i(Y(K_{\Sigma}), - )$ by $H^i(\Sigma, - )$, 
and similarly for $H^i_c$.
Similarly,
we denote $H^{i}(\partial Y(K_{\Sigma}), - )$ by $H^i(\partial \Sigma, -)$,
and by $H^{i}_!(\Sigma, -)$ the image of
the group $H^i_c(\Sigma, -)$ in $H^i(\Sigma,-)$. 
 
 \medskip

For homology, we have groups $H^{\bm}_{i}$ and corresponding exact sequences
and isomorphisms to the general case:
\begin{equation} \label{forearlierref} H_{2}(\partial Y(K),- ) \rightarrow   H_{1}(\partial Y(K),- ) \rightarrow 
H_1(Y(K) ,-) \rightarrow  H_{1, \bm}(\partial Y(K),- ) \ldots,\end{equation}
 In particular, 
the cuspidal homology $H_{1,!}(\Sigma,\Z)$ is defined to be the quotient of $H_1(\Sigma,\Z)$ by the image of
$H_1(\partial \Sigma,\Z)$, or, equivalently, the image of $H_1(\Sigma, \Z)$ in $H_{1}^{\bm}(\Sigma, \Z)$. 
  
  \medskip

 \subsection{The homology of the cusps}  \label{cuspshavenohomology}

An often useful 
\footnote{We do not make any {\em essential} use of this. For example, we do not assume it is so
for our main analytical computations, and 
in the general case~\S~\ref{sec:eisintegrality} gives a good handle on $H_1(\partial Y(K))$. However,
we often use it for essentially ``cosmetic'' reasons.} fact is that the $H_1$ of cusps vanish for orbifold reasons:
 \begin{lemma} \label{vanishingboundaryhomology}
  \label{lemma:noeis}
$H_1(\partial Y(K_{\Sigma}),\Z_p) = H^1(\partial Y(K_{\Sigma}), \Z_p) = 0$ for $p \ne 2$.
\end{lemma}

This is certainly not true for general level structures, i.e., $H^1(\partial Y(K), \Z_p)$ can certainly be nonzero.
This is simply a special, and convenient, property of the $K_{\Sigma}$-level structures.   Geometrically, the cusps of $Y(K_{\Sigma})$
are all quotients of a torus modulo its hyperelliptic involution (which is a sphere with four cone
points of order $2$), and in particular have trivial $H_1$  in characteristic larger than $2$.

\begin{proof}    Clearly, 
\begin{equation} \label{Ybone} Y_B(K) = \coprod_{z \in \B(F) \backslash \G(\Afinite)/K} \Gamma_z \backslash \H^3, \end{equation} 
where the union is taken over $\B(F)$-orbits on $\G(\Afinite)/K$, and 
for $z \in \G(\Afinite)/K$ representing such an orbit, 
 $\Gamma_z$ is the stabilizer of $z$ in $\B(F)$.
 \medskip
 
  Write $\alpha: \B \rightarrow \mathbb{G}_m$ for the natural map
$ \left( \begin{array}{cc} a & b \\ 0 & c \end{array} \right)\mapsto a c^{-1}$.
 We claim that $\Gamma_z$ always contains an 
 element  (call it $b_0$) such that $\alpha(b_0) =-1$. 
 
 \medskip
 
 In fact, choose any element $b \in \B(F)$ with $\alpha(b) = -1$. 
 Then the set of $x \in F_v$ such that  $b n(x) $ lies in $ g K_v g^{-1}$
 is nonempty open for every $v$. Indeed, by a local computation,  the image  of $\B(F_v) \cap g K g^{-1}$ under $\alpha$ always contains, for any $g \in \G(\Afinite)$, the units $\OO_v^{\times}$. By strong approximation there exists $x_0 \in F$
such that $b n(x_0) \in g K g^{-1}$.  Then $b_0 = b n(x_0)$ will do.

  \medskip

    This is enough for our conclusion: 
The element $b_0$, acting by conjugation on $\Gamma_z^0 := \Gamma_z \cap \ker(\alpha) \cong \mathbf{Z}^2$, necessarily acts by negation;   in particular, its  co-invariants on  $H_1(\Gamma_z^0, \Z_p)$ are trivial. \end{proof}

\medskip

Even in the general case,  i.e. under {\em no} assumptions on the nature of the level structure $K$, 
the Hecke action on the boundary is simple to compute:

\begin{lemma} \label{heckeactiononboundary}
For any Hecke eigenclass  $ h \in H_1(\partial K, \C)$ there is a 
  Grossencharacters $ \psi: \adele_F^{\times}/F^{\times} \rightarrow \C^{\times},$
such that:
\begin{itemize} 
\item The  restriction of $\psi$ to $F_{\infty}^{\times} = \C^{\times}$ is 
given by $z \mapsto \sqrt{z/\bar{z}}$;
\item  $\psi$ is
unramified at any prime where  $K$ is conjugate to $\G(\OO_v)$ or $K_{0,v}$;
\item The eigenvalue $\lambda_{\q}(h)$ of $T_{\q}$ on $h$ is 
 \begin{equation} \label{Heckeactiononboundary} \Norm(\q)^{1/2} \left(  \psi(\q) + \overline{\psi(\q)}\right).  \end{equation}
 \end{itemize}
\end{lemma}

We omit the proof, which can be given by direct computation or via the theory of Eisenstein series. 

\medskip

Note the following consequence:  Let $S$ be the set of places
at which the level structure $K$ is not maximal. For any classes $h_1, h_2$ as in the Lemma, 
 with corresponding characters $\psi_1 \neq \psi_2$, there
 exists $\q$ such that
 $$ \lambda_{\q}(h) \not \equiv \lambda_{\q}(h_2) \mbox{  mod $\mathfrak{l}$.}$$
whenever  $\mathfrak{l}$ is any prime of $\overline{\Q}$
above a prime $\ell$ not dividing $h_F \cdot \prod_{q \in S} q(q-1)$.

\subsection{Relationship between \texorpdfstring{$Y_B(K)$}{Y} and \texorpdfstring{$Y_B(K')$}{Y'}} \label{ybybxi}

The following discussion will be used in the proofs of Section~\S~\ref{analysissec}.
It says, roughly, that the cusps of $Y(\Sigma \cup \q)$ are, metrically speaking,
two copies of the cusps of $Y(\Sigma)$.

\subsubsection{}

 Suppose $\Xi$ is a finite set of finite places disjoint from any level in $K$,
 and let $K' = K \cap K_{\Xi}$, so that $K'$ is obtained from $K$
 by ``adding level at $\Xi$.''

  Then there is an isometry
\begin{equation} \label{levelraiseisometry} Y_B(K') \cong  \{1, 2\}^{\Xi} \times Y_B(K), \end{equation} 
 where $q = |\Xi|$, that is to say; 
 $Y_B(K')$ is  an {\em isometric} union of $2^{|\Xi|}$
copies of $Y_B(K)$.  Moreover, these $2^q$ isometric copies
are permuted simply transitively by the group of Atkin-Lehner involutions $W(\Xi)$ 
(see~\S~\ref{subsec:ali}). 
\medskip

 In view of the definition, 
  $$  Y_B(K) = \B(F) \backslash   \H^3 \times \G(\Afinite)/   K $$
    $$  Y_B(K') = \B(F) \backslash   \H^3 \times \G(\Afinite)/   K' $$

it suffices to construct a $\B(F)$-equivariant bijection 
 \begin{equation} \label{bequivbij} \{1, 2\}^{\Xi}  \times   \G(\Afinite)/K \stackrel{\sim}{\rightarrow} \G(\adele)/K'.\end{equation} 
 
 In fact, the map
\begin{eqnarray*} \Lambda_v:  \{1, 2\} \times \G(F_v)/\G(\OO_v) &  \longrightarrow& \G(F_v)/K_{0, v} \\
(i,  bK_v) &\longmapsto&  \begin{cases} b K_{0,v}, &  i = 1 \\ b w_v K_{0,v}, &  i=2 \end{cases} 
 \end{eqnarray*}
 where $w_v=  \left(  \begin{array}{cc} 0 & 1 \\ \varpi_v & 0 \end{array}\right)$, 
 defines a $\B(F_v)$-equivariant bijection,\footnote{Geometrically speaking, the quotient $\G(F_v)/K_{\Xi, v}$ is identified with arcs on the Bruhat--Tits building; the Borel subgroup $\B(F_v)$ is the stabilizer of a point on the boundary $\partial \mathcal{B}$; the two orbits correspond to arcs that point ``towards'' or ``away from'' this marked point; and the Atkin-Lehner involution reverses arcs. Explicitly, in the action of $\B(F_v)$ on $G/K_{\Xi,v}$, the stabilizer
of both $K_{\Xi,v}$ and $\left(  \begin{array}{cc} 0 & 1 \\ \varpi_v & 0 \end{array}\right)  K_{\Xi, v}$
is the same, namely, elements of $\B(F_v)$ that belong to $\PGL_2(\OO_v)$. }
 and we take the product of $\Lambda_v$ over all $v \in \Xi$ to get  a map~\eqref{bequivbij}
 that is even $\B(\Afinite)$-equivariant.

  \medskip

  \begin{remarkable} \label{rem:heightchange} Note that this bijection {\em does not} preserve the height function.  
  In fact,  since $\height(b w_v) = \height(b) \sqrt{q_v} $, we see that 
  the height function on $Y_B(K')$ pulls back 
  to the function on $\{1, 2\}^{\Xi} \times Y_B(K)$ given by the product
  of the height on $Y_B(K)$ and the functions $1 \mapsto 1, 2 \mapsto \sqrt{q_v}$
  for each $v \in \Xi$.  \end{remarkable}
  
  \medskip
  
 \begin{remarkable} \label{rem:equiv} {\emph Choose an auxiliary place $v$ not contained in $\Xi$;
and suppose $K^*_v \subset \G(F_v)$ is such that $K^*_v$ is normal in $K_v$. 
Then $K_v/K^*_v$ acts on both sides of~\eqref{levelraiseisometry}, where the action on $\{1, 2\}^{\Xi}$ is trivial,  and 
the identification is {\em equivariant} for that action. }
\end{remarkable}

 \subsubsection{} \label{topologyprojectiondown}
The topology of the projection $ \partial Y(\Sigma \cup \q) \rightarrow \partial Y(\Sigma)$ can be deduced similarly: \medskip

 The preimage of each component $P$ of $\partial Y(\Sigma)$ consists of two components  $P_1, P_2$ of $\partial Y(\Sigma \cup \q)$,
 interchanged by the Atkin--Lehner involution. The induced map $H_2(P_i) \rightarrow H_2(P)$
 are given by multiplication by $\Norm(\q)$ resp. trivial. 
 \medskip

 In particular,  if $\mathcal{Q}$ is any set containing all but finitely many prime ideals, 
 the cokernel of 
 $$ \bigoplus_{\q \in \mathcal{Q}} H_2( \Sigma \cup \q, \Z_p)^{w_{\q} -} \longrightarrow H_2(\Sigma, \Z_p),$$
 where the maps are the difference of the two degeneracy maps,
 equals {\em zero} unless $p$ divides $w_F$. Indeed
 $w_F =\mathrm{gcd}_{\q \in \mathcal{Q}}(\Norm(\q)-1)$.

\subsection{Congruence and essential homology, split case} \label{ss:Noncompactcongruence}

The definitions of  congruence  homology of~\S~\ref{subsubsection:reducetogroup} apply equally well to the compact and noncompact cases. 

\section{Some results from Chapter~\ref{chapter:ch3} in the split case} \label{sec:llsplit} 
In this section we discuss the extension of various results  from Chapter~\ref{chapter:ch3} to the split case. 
{\em The assumption on the congruence subgroup property in Lemma~\ref{iharalemma}
is known unconditionally here for $p >2$, making all the corresponding results unconditional.}

\begin{itemize}
 \item Ihara's lemma (Lemma~\ref{iharalemma}) applies:  the original proof made no restriction on $\G$ being nonsplit.
 \item Theorem~\ref{theorem:surjective} still holds; the proof is given in~\S~\ref{theoremsurjectivesplitproof} below.
 
 \item 	 The level-raising theorem (Theorem~\ref{theorem:ribet}) applies, with the same proof, 
 since it relies only on Lemma~\ref{iharalemma} and Theorem~\ref{theorem:surjective}.
\end{itemize}

After we give the proof of Theorem~\ref{theorem:surjective} in the split case, we give
results about comparison of regulators, of exactly the same nature as 
Lemma~\ref{regulator-compare-1} in the nonsplit case.

\subsection{Proof of  Theorem~\ref{theorem:surjective} in the split case}  \label{theoremsurjectivesplitproof}
 
Note that all the natural maps $H^1_c \rightarrow H^1$
 have become isomorphisms, since we have localized
 at a non-Eisenstein maximal ideal.  Indeed, 
as long as the ideal is not cyclotomic-Eisenstein this is true: such a localization kills $H^0(\partial Y(K))$,
and in this setting $H^1(\partial Y(K))$ vanishes by~\S~\ref{sec:llsplit}.
  
Now the proof is exactly as the proof of Theorem~\ref{theorem:surjective}, considering now the diagram: 
$$
\begin{diagram}
0& & \ker(\Phi^{\vee}_{\mathfrak{m}}) &   & \ker(\Psi_{\mathfrak{m}})  & & 0  \\
\dTo & &  \dTo &  & \dTo && \dTo \\ 
H^1_c(\Sigma, \Q_p)_{\mathfrak{m}} & \rTo& H^1_c(\Sigma, \Q_p/\Z_p)_{\m} & \rTo_{\delta_{\Sigma}} & H_1(\Sigma, \Z_p)_{\m}  &\rTo &  H_1(\Sigma, \Q_p)_{\m} \\ 
\dTo_{\Phi^{\vee}_{\Q, \m}} && \dTo_{\Phi^{\vee}_{\m}} & & \dTo_{\Psi_{\m}} &  &  \dTo  \\
H^1_c(\Sigma/\q, \Q_p)_{\mathfrak{m}}^2 &\rTo & H^1_c(\Sigma/\q, \Q_p/\Z_p)_{\m}^2& \rTo_{\delta_{\Sigma/\q}} & H_1(\Sigma/\q, \Z_p)_{\m}^2 & \rTo & H_1(\Sigma/\q, \Q_p)_{\m}^2  \\\dOnto & & \dOnto &  &    \dTo  && \dTo \\
0 && \coker(\Phi^{\vee}_{\m}) &
  & 0&&0   \\
\end{diagram}
$$

  \subsection{Comparison of regulators}

\begin{lemma} \label{h2psiveesplit}

Suppose $\Sigma$ is a level for which $H_1(\Sigma, \C)^{\qnew} = 0$, 
i.e. the level raising map $\Psi^{\vee}_{\C} : H_1(\Sigma/\q, \C)^2 \rightarrow H_1(\Sigma, \C)$ is an isomorphism. 
Then the cokernel of $$\Psi^{\vee}: H_{2,!}(\Sigma/\q, \Z)^2 \rightarrow H_{2, !}(\Sigma, \Z) $$
is of order $h_{\cl}(\Sigma/\q;\q)$, up to orbifold primes. The same conclusion holds with $H_{2!}$ replaced by Borel--Moore homology. 
 \end{lemma}

Here $h_{\cl}$ is as in
\S~\ref{LiftableCongruenceHomology}.

\begin{proof}  
As in the proof of Lemma
\ref{regulator-compare-1} 
the map 
$\Psi^{\vee}_{\tf}: H_1(\Sigma)_{\tf} \rightarrow H_1(\Sigma/\q)_{\tf}^2$ 
 has cokernel  of size  $H_{1}(\Sigma/\q;\q)$, up to orbifold primes.  By dualizing (see~\eqref{splitduality})
$H^{\bm}_{2} (\Sigma/\q, \Z)^2 \stackrel{\Psi^{\vee}}{\rightarrow} H^{\bm}_{2}(\Sigma, \Z)$
also has cokernel of order $h_{\cl}(\Sigma/\q; \q)$, up to orbifold primes. 
The map $H_{2,!} \rightarrow H_{2, \bm}$ is an isomorphism with $\Z[\frac{1}{2}]$ coefficients, because of~\S~\ref{cuspshavenohomology}. 
 \end{proof}
 
 \begin{remarkable} \emph{
 There is a more robust way to pass from $H_{2, \bm}$ to $H_{2,!}$, 
 even if $H_1(\partial Y(K)$ did not vanish for trivial reasons: The Hecke action on
 the cokernel of $\Psi^{\vee}_{\bm}$ is congruence-Eisenstein, by our argument above;
 on the other hand, the Hecke action on $H_1(\partial Y(K))$ has been computed in~\eqref{Heckeactiononboundary}, 
 and it is easy to show these are incompatible for most primes $p$. 
 }
 \end{remarkable}

\begin{theorem}[Comparison of regulators] \label{h2regregcompare} 
Suppose that
 ${\q}$ is a prime in $\Sigma$ and $H_1^{{\q}-\mathrm{new}}(\Sigma, \C) = \{0\}$. 
  Let $b_0$ be the number of connected components of the manifold
  $Y(\Sigma)$. Let $h_{\cl}$ be as in~\S~\ref{LiftableCongruenceHomology} and put
$$D=  \det \left( T_{\q}^2 - (1 + N(\q))^2 \big|  H_1(\Sigma/\q , \C) \right).$$ 
 Let $D_{\cusp}$ be the same  determinant restricted to the cuspidal cohomology. 
Then up to orbifold primes:  
 \begin{equation} \label{h2regcompare} \frac{\reg(H_{1}(Y(\Sigma/\q) \times \{1,2\})}{ \reg(H_{1}Y(\Sigma))} =   \frac{\sqrt{D}}{|h_{\cl}(\Sigma/\q;\q)|}  %
,  \end{equation} 
whereas 
 \begin{equation} \label{h2regcompareprime} \frac{\reg(H_{2!}(Y(\Sigma/\q) \times \{1,2\})}{ \reg(H_{2!}Y(\Sigma))} =   \frac{|h_{\cl}(\Sigma/\q;\q)|}{\sqrt{D_{\cusp}}}  %
 \end{equation}  

 Moreover, in the same notation,  $$  \frac{\reg(H_{0}(Y(\Sigma/\q) \times \{1,2\})}{ \reg(H_{0}Y(\Sigma))} = \frac{\vol(Y(\Sigma))^{1/2}}{\vol(Y(\Sigma/\q))},$$
and the $H_{3,\bm}$ regulators change by the inverse of these. 
 
\end{theorem}

Note in this theorem that, e.g. $\reg(H_1(M \times \{1,2\}))$ is simply $\reg(H_1(M))^2$;
we phrase it in the above way for compatibility with some later statements. 
\begin{proof}
We first give the proof for $H_1$, which is very close to that already given in Theorem~\ref{regcompare}.   
For simplicity we write the proof as if there are no orbifold primes; the general case, of course,
only introduces ``fudge factors'' at these primes. 

 By Ihara's lemma, the cokernel
of  $\Psi: H_1(\Sigma,\Z) \rightarrow  H_1(\Sigma/\q,\Z)^2 $
is   congruence homology (at least at primes away from $2$).  
Note that, {\em just as in the compact case}, 
the map $\Psi \otimes \Q: H_1(\Sigma, \C) \rightarrow H_1(\Sigma/\q, \C)^2$
is still an isomorphism; this is {\em not so for $H_2$}, which is one reason why 
we work with $H_{2!}$ rather than $H_2$.   In particular, 
the map $\Psi_{\tf}:  H_1(\Sigma,\Z)_{\tf} \rightarrow H_1(\Sigma/\q, \Z)_{\tf}^2$ is still injective, and, as before we have the diagram:

\begin{equation}\label{cokerpsitors}
\begin{diagram}
&& (\ker \Psi)_{\tors} & \rTo &  \ker \Psi & \rTo & 0  \\
&&  \dTo & &  \dTo &   & \dTo \\ 
    0 & \rTo &  H_1(\Sigma, \Z)_{\tors} & \rTo_{\Psi} &  H_1(\Sigma , \Z) & \rOnto & H_1(\Sigma, \Z)_{\tf} \\
&& \dTo_{\Psi_{\tors}} & & \dTo_{\Psi} &  & \dTo_{\Psi_{\tf}} \\
 0& \rTo & H_1(\Sigma/\q, \Z)^2_{\tors} & \rTo_{\Psi_{\tf}} &H_1(\Sigma/\q, \Z)^2 & \rOnto & H_1(\Sigma/\q, \Z)^2_{\tf}  \\
&& \dOnto & & \dOnto &  &  \dOnto \\
&&\coker(\Psi_{\tors}) & \rTo
      & \coker(\Psi) & \rTo &  \coker(\Psi_{\tf}) \\
\end{diagram}
\end{equation}

As before, the diagram shows that the order of $\coker(\Psi_{\tf})$ 
equals the order of the cokernel of $\coker(\Psi_{\tors}) \rightarrow \coker(\Psi)$, 
equivalently, the order of the cokernel of $H_1(\Sigma/\q, \Z)_{\tors}^2 \rightarrow
H_1(\Sigma/\q, \Z)_{\mathrm{cong}}$, that is to say
$h_{\cl}(\Sigma/\q;\q)$.

In the sequence $ H_1(\Sigma,\Z)_{\tf} \stackrel{\Psi_{\tf}}{\rightarrow} H_1(\Sigma/\q,\Z)^2_{\tf}
\stackrel{\Psi_{\tf}^{\vee}}{\rightarrow} H_1(\Sigma,\Z)_{\tf}$, 
 the composite
map $\Psi \circ \Psi^{\vee}$ is given explicitly by Lemma~\ref{lemma:composite};  from that we deduce $$| \coker(\Psi_{\tf})|  \cdot  | \coker(\Psi_{\tf}^{\vee}) | =
|\coker(  \Psi_{\tf} \circ \Psi_{\tf}^{\vee} ) | = %
D.$$  Therefore the size of the cokernel of $\Psi^{\vee}_{\tf}$ equals $D/ h_{\cl}(\Sigma/\q;\q)$. 
\medskip

Choose an orthonormal basis $\omega_1, \dots, \omega_{2k}$  
 for $H^1(\Sigma,\R)$, and also choose a
  basis $\gamma_1, \dots, \gamma_{2k}$ for $H_1(\Sigma/{\q}, \Z)_{\tf}^2$.  Then
 $$   \det \  \langle \gamma_i, \Phi^{\vee}_{\R} \omega_j \rangle = \det \ \langle \Psi^{\vee}(\gamma_i), \omega_j \rangle  = \frac{D}{h_{\cl}(\Sigma/\q)}  \cdot  \reg(H_1(\Sigma))$$ 
 Now
 $|\mathrm{coker}(\Psi^{\vee})| = \frac{D}{|H_{1, \lif}(\Sigma/\q)|}$. 
Also 
 $ \| \Phi_{\R}^{\vee} \omega_1 \wedge \Phi_{\R}^{\vee} \omega_2 \wedge \dots \Phi_{\R}^{\vee}\omega_{2k} \| = \sqrt{D}$ and  consequently,
 $$  \reg(H_1(\Sigma/\q))^2=  \frac{ \det \  \langle \gamma_i, \Phi^{\vee}_{\R} \omega_j \rangle }{ \| \Phi_{\R}^{\vee} \omega_1 \wedge \Phi_{\R}^{\vee} \omega_2 \wedge \dots \Phi_{\R}^{\vee}\omega_{2k} \| } = \reg(H_1(\Sigma)) \cdot \frac{\sqrt{D}}{h_{\cl}(\Sigma/{\q})},$$ 
which implies the desired conclusion

We now turn to $H_2$, or rather $H_{2!}$; the proof of~\eqref{h2regcompare} is very similar to that just given. 
 Fix a basis $\gamma_i$ for $H_{2,!}(\Sigma, \Z/\q)^{\oplus 2}$ and 
an orthonormal basis $\omega_1, \dots, \omega_{2k}$  
 for $\Omega^2_{\cusp}(\Sigma)$. Then:
 $$   \det \  \langle \gamma_i, \Phi^{\vee}_{\R} \omega_j \rangle = \det \ \langle \Psi^{\vee}(\gamma_i), \omega_j \rangle  = | \mathrm{coker}(\Psi_!^{\vee})|  \cdot \reg(H_{2!}(\Sigma)) $$

 The $\Psi_{\R}^{\vee} \omega_j$ do not form an orthonormal basis; 
 indeed, the map $\Psi_{\R}^*$  multiplies volume by $\sqrt{D_{\cusp}}$, i.e.
 $ \| \Phi_{\R}^{\vee} \omega_1 \wedge \Phi_{\R}^{\vee} \omega_2 \wedge \dots \Phi_{\R}^{\vee}\omega_{2k} \| = \sqrt{D}.$  Consequently,
 
 $$ \reg(H_{2,!}(\Sigma/\q)^2 =  \frac{ \det \  \langle \gamma_i, \Phi^{\vee}_{\R} \omega_j \rangle }{ \| \Phi_{\R}^{\vee} \omega_1 \wedge \Phi_{\R}^{\vee} \omega_2 \wedge \dots \Phi_{\R}^{\vee}\omega_{2k} \| } = \reg(H_{2!}(\Sigma))  \frac{ h_{\cl}(\Sigma/\q)}{\sqrt{D_{\cusp}}},$$
 where we have used the fact (Lemma~\ref{h2psiveesplit}) that $\# \coker(\Psi^{\vee}_{!}) = h_{\cl}(\Sigma/\q)$. 
This implies ~\eqref{h2regcompare}.

Finally the $H_0$-statement is clear. 
\end{proof}

\section[Scattering matrices]{Eisenstein series for arithmetic manifolds: explicit scattering matrices} \label{sec:Eis-series}  
We briefly recall those aspects of the theory of Eisenstein series. We need this to explicitly evaluate the ``scattering matrices''; {\em if the reader is willing to accept this on faith, 
this section can be skipped without detriment. }
What is needed from this section, for our later Jacquet--Langlands application, is {\em only} a comparison between the scattering matrices at different levels showing that, roughly speaking, ``everything cancels out;'' thus the details of the formulas do not matter, only that they match up at different levels. 

In principle all of the computations of this section can be carried out in classical language;
but the reason we have adelic language
is not mathematical but aesthetic:  it makes the various notational acrobatics involved in comparing the manifolds
at two different levels $Y(K), Y(K')$ much easier.  Recall that these manifolds are both disconnected and have multiple cusps in general!

 The reader who prefers classical language may wish to study the paper of Huxley \cite{Huxley}
 which carries out a similar computation for arithmetic quotients of $\SL_2(\Z)$. 
 
Here is the main result of the section:
\begin{theorem} \label{scatteringmatrixtheorem}
Equip $Y(\Sigma)$ and $Y(\Sigma \cup \q)$ with the height function as in~\S~\ref{infgeom}; 
let $\Psi_{\Sigma}, \Psi_{\Sigma \cup \q}$ be the scattering matrices on functions
and $\Phi_{\Sigma}, \Phi_{\Sigma \cup \q}$ the scattering matrices on $1$-forms. (See~\S~\ref{subsec:einsteinintro} for definitions.)

Also equip $Y(\Sigma) \times \{1,2\}$ with the height function induced by the identification
$Y_B(\Sigma \cup \q) \cong Y_B(\Sigma) \times \{1,2\}$ of~\eqref{levelraiseisometry},
and let $\Psi_{\Sigma \times \{1,2\}}$ and  $\Phi_{\Sigma \times \{1,2\}}$
be the corresponding scattering matrices.
 
For $s \in \C, z \in \C^{\times}$ write $M(s) =  \left(\begin{array}{cc} 1 & 0 \\ 0 & q^{s} \end{array}\right)$ and 
$$ N(z,s) = \frac{1}{(q-z^2 q^{-s})} \left(\begin{array}{cc}  z^2 q^{-s} (q-1) & 1-z^2 q^{-s} \\ q (1- z^2 q^{-s})  & q-1 \end{array} \right)
.$$ 
Here $q = N(\q)$. 
Note that $N(z,s) N(1/z, -s)$ is the identity
and $N(\pm 1, 1)$ is a projection with eigenvalue $1$.

Then there exists a finite set $X$ of Grossencharacters\footnote{That is to say:
characters of $\adele_F^{\times}/F^{\times}$.}, unramified at $\q$,  and matrices
$A_{\chi}(s) \ \ (\chi \in X)$ such that:

\begin{eqnarray}  \label{Psicomparison} 
  \Psi(s) & \sim & \bigoplus_{\chi \in X} A_{\chi}(s) , \\ 
 \Psi_{\Sigma \times \{1,2\} }(s) &\sim& \bigoplus_{\chi \in X}  A_{\chi}(s) \otimes M(s),  \\
 \Psi_{\Sigma \cup \q }(s) &\sim& \bigoplus_{\chi \in X} A_{\chi}(s) \otimes N(\chi(\q), s),
\end{eqnarray} 
where we write $\sim$ for  equivalence of linear transformations.\footnote{Thus, if $A: V \rightarrow V$ and $B: W \rightarrow W$, we write $A \sim B$
if there is an isomorphism $V \rightarrow W$ conjugating one to the other.}

 Exactly the same conclusion holds for $\Phi$, with a different set of Grossencharacters $X$.

  In particular the following hold:
 
 \begin{enumerate}
 \item   If $A_{\chi}$ has size $a_{\chi} \in \mathbb{N}$, then  
\begin{eqnarray} \label{scatmatrix1}
\frac{\det \Psi_{\Sigma \cup \{1,2\}}(s)}{\left( \det \Psi_{\Sigma} \right)^2}   &= & \prod_{\chi \in X} \left(  q^{s}  \right)^{a_{\chi}}  \\
\label{scatmatrix2} 
\frac{\det \Psi_{\Sigma \cup \q}  (s)}{ \left( \det \Psi_{\Sigma}(s)  \right)^2}  &= &  \prod_{\chi \in X}  \left( (\chi(\q)^2 q^{-s}  \cdot \frac{ q -\chi(\q)^{-2} q^{ s}}{q -\chi(\q)^{2} q^{- s}}  \right)^{a_{\chi},}
\end{eqnarray}
with exactly the same conclusion for $\Phi$, allowing a different set of Grossencharacters $X$. 

\item 
At $s=0$, 
  \begin{equation} \label{zeropoint} \trace \Psi_{\Sigma \cup \{1,2\}}(0)  = \trace \Psi_{\Sigma \cup \q}(0) = 2 \trace\Psi(0). \end{equation}

\item  $R_{\Sigma}, R_{\Sigma \times \{1,2\}}, R_{\Sigma \cup \q}$ are the respective residues when $s=1$, 
  then $R_{\Sigma \cup \q} \sim R_{\Sigma} \oplus 0$ whereas $R_{\Sigma \times \{1,2\}} \sim R_{\Sigma} \oplus q R_{\Sigma} \oplus 0$, where $0$ denotes a zero matrix of some size. 
  
  \item   Let $\p$ be an auxiliary prime, not belonging to $\Sigma \cup \{\q\}$. 
  Then the set of Grossencharacters $X$ that arise in the above analysis for $\Psi$
 is the same at levels $\Sigma$ and $\Sigma \cup \{\p\}$.  The same is true for $\Phi$. 
  Moreover, the size of the corresponding matrices satisfies
  $a^{\Sigma \cup \{\p\}}_{\chi} = 2 a^{\Sigma}_{\chi}$. 
  \end{enumerate}
 \end{theorem}

Note that enunciations (1)--(3) are consequences of~\eqref{Psicomparison}, 
and (4) will follow in the course of the proof of \eqref{Psicomparison}, so we will only prove that.

\subsubsection{Measure normalizations}  \label{subsubmeas} 
 
 We have normalized previously measures on $\adele$ and $\adele^{\times}$  (see~\S~\ref{subsec:measures}). 
 We equip $\N(\adele)$ with the corresponding measure arising from the map 
 $y \mapsto \left(\begin{array}{cc} 1 & y \\ 0 & 1 \end{array}\right)$,  and
 we denote by  $d_l b$  the left Haar measure on $\B(\adele)$ given in coordinates $\left(\begin{array}{cc} x & y \\ 0 & 1 \end{array}\right)$ as $|x|^{-1} dx \cdot dy$.
 
 \medskip

 We equip  $\G(\adele)$
with the measure defined to be the push-forward of $\frac{1}{2} d_{l} b \cdot dk$ under $(b,k) \in \B(\adele) \times \Kmax \rightarrow bk \in \G(\adele)$, where $dk$ is the Haar probability measure on $\Kmax$. 

Push forward this measure to $\G(\adele)/\Kmax$. This is a union of copies of $\mathbf{H}^3$:  
for every $g_f \in \G(\Afinite)/\PGL_2(\widehat{\OO})$ we may 
identify  $$ \G(\C) \times g_f \PGL_2(\widehat{\OO}) / \Kmax$$ with $\G(\C)/\mathrm{PU}_2\simeq \H^3$. 
On each such copy of  $\H^3$ the induced measure is the 
same as that induced by the Riemannian structure.

  \medskip
  
We equip $\G(\adele)^{(1)}$ with the measure obtained as the ``fibral'' measure  of $\G(\adele) \stackrel{\height}{\rightarrow} \R_{>0}$, where the measure on $\R_{>0}$ is $dx/x$.

\medskip

Push this measure forward to $\G(\adele)^{(1)}/\Kmax$.   
This  intersects each component $\mathbf{H}^3$ from above
in a level set $y = \mathrm{const}$.  The induced measure on this  is   $\frac{dx_1 dx_2}{y^2}$.

\subsection{The induced space.}  \label{Vxidefsec}

Let $\xi$ be any character of $\C^{\times}$.  \index{$\eta$ (character of $\mathbf{C}^{\times}$)}
In what follows, the character $z \mapsto \frac{z}{|z|^{1/2}} =  \sqrt{z/\bar{z}}$ and its inverse
will play an important role, and consequently we denote it by $\eta$:
$$ \eta: z \in \mathbf{C}^{\times} \mapsto \sqrt{z/\bar{z}}.$$ 
Define $ V_{\xi}(s,K) $ to be: \index{$V_{\xi}(s,K)$} \index{$V_{\xi}(s)$}
\begin{equation} \label{vxidef}  \left\{ f  \in C^{\infty}(\G(\adele)/K): 
 f \left(  \left( \begin{array} {cc} x & \star \\ 0 & 1 \end{array} \right)  g \right) = 
|x|^{1+s} \xi(x)  f(g) ,   \ \  x  \in \C^{\times}. \right\} \end{equation} 

Note that $|x|$ denotes the usual absolute value on $\C$, so that, for example, $|2| = 2$; this is the {\em square root} of the ``canonically normalized'' absolute value that is defined for any local field, i.e.,
the effect of multiplication on Haar measure.   

When the compact $K$ is understood, we write simply $V_{\xi}(s)$. 
This is a representation of $\G(\R)$ of finite length (in particular, any $K_{\infty}$-isotypical component has finite dimension). 

  Multiplication by $\height(g)^{s}$ gives an identification $V_{\xi} \stackrel{\sim}{\rightarrow} V_{\xi}(s)$.
  We shall often use this isomorphism without explicit comment that it is multiplication
  by $\height(g)^{s}$. 
    
We equip $V_{\xi}(s)$ with the unitary structure 
\begin{eqnarray*} \|f\|^2  = \int_{  B(F) \backslash \G(\adele)^{(1)}} |f(g)|^2 dg. \end{eqnarray*}
where  we normalize measures in a moment;
this unitary structure is $\G(\Adele)$-invariant for $s$ purely imaginary.

\subsubsection{}
We relate the spaces just introduced to to the standard presentation of Eisenstein series (see, for example,
\cite{BorelJacquet}):

Note that $V_{\xi}(s,K)$ decomposes according to Grossencharacters
 of $F$:   If we denote by $X(s)$ the set of characters $\adele_F^{\times}/F^{\times}$
 whose restriction to $F_{\infty}^{\times} = \C^{\times}$ is $\xi \cdot | \cdot|^s$ --- 
 we sometimes abbreviate this simply to $X$ when $s$ is understood --- then 
\begin{equation} \label{ind-dec} V_{\xi}(s,K) = \bigoplus_{\chi \in X(s)} \mathcal{I}(\chi) ^K, \end{equation} where, for $\chi$ a character of $\adele_F^{\times}$, 
 we put
 $$ \mathcal{I}(\chi) = \{ f : \G(\adele) \rightarrow \C  \mbox{ smooth}: f(bg) = \chi(b) |b|_{\adele}^{1/2} f(g) \mbox{ for } b \in \B(\adele)\},$$
 Here we write $\chi(b)$ and $|b|_{\adele}$ as shorthand for $\chi (\alpha(b))$ and $|\alpha(b)|_{\adele}$, where $\alpha$
 is the positive root, as in~\eqref{alphadef}.  
 The set $X(s)$ of characters is infinite, but only finitely many $\chi \in X(s)$
 have $\mathcal{I}(\chi)^K \neq 0$. 
 
  We will write $X$ as a shorthand for $X(0)$. 
  Note  that there is a bijection between
 $X(0) =X $ and $X(s)$ given by twisting by $|\cdot|_{\adele}^{s/2}$. 

We make a similar definition if $\chi$ is a character of a localization $F_v^{\times}$, replacing $\adele$ by $F_v$ on the right-hand side.  In that case $\mathcal{I}(\chi)$
is a representation of $\PGL_2(F_v)$, and there is a natural map
$\bigotimes_v \mathcal{I}(\chi_v)  \stackrel{\sim}{\rightarrow} \mathcal{I}(\chi)$, 
where we use the restricted tensor product (see \cite{Flath}). 
  
 The theory of Eisenstein series  gives a canonical intertwiner $\mathrm{Eis} :V_{\xi}(s) \rightarrow C^{\infty}(\G(F) \backslash \G(\adele)/K)$ that is meromorphic --- 
 that is to say, the composition $$V_{\xi}(0,K)   \stackrel{\sim}{\rightarrow} V_{\xi}(s,K) \stackrel{\mathrm{Eis}}{\longrightarrow} C^{\infty}(\G(F) \backslash \G(\adele)/K)$$ is meromorphic, where the first identification is multiplication by $\mathrm{Ht}^{s}$. 
 If $\xi$ is unitary,
this map is holomorphic when $\Re(s) = 0$.    This map is uniquely characterized by the property 
$$ \mathrm{Eis}( v) (e) = \sum_{\gamma \in \B(F) \backslash \G(F)} v(\gamma)$$
whenever $\Re(s) > 2$.

  There is a {\em standard intertwining operator} 
 $V_{\xi}(s) \rightarrow V_{\xi^{-1}}(-s)$ given by (the meromorphic extension of)
\begin{equation} \label{Intertwinerdef}  f \mapsto \left( g \mapsto  \frac{1}{\vol(\N(F) \backslash \N(\adele))} \int_{n \in  \N(\adele)} f(\left( \begin{array}{cc} 0 & -1 \\ 1 & 0 \end{array}\right) n g) dn  \right). \end{equation}
Then
\begin{equation} \label{Eisconstant} \mathrm{Eis}(f)_N = f + M(f). \end{equation}    
 where in this context the constant term $\mathrm{Eis}(f)_N$
 is defined, as usual for automorphic forms, via
$ \frac{1}{\vol(\N(F) \backslash \N(\adele))} \int_{n \in \N(F) \backslash \N(\adele)} 
\mathrm{Eis}(f) (ng ) dn$; this will be compatible with our other usages of the notation.  
 
\subsection{\texorpdfstring{$V_{\xi}$}{V} and the spaces \texorpdfstring{$C^{\infty}(s), \Omega^{\pm}(s)$}{C}.} 

Recall that  ---  for $M$ a hyperbolic manifold ---  we have defined spaces of functions $C^{\infty}(s)$
and of one-forms $\Omega^{\pm}(s)$
in~\S~\ref{mbfunction}. In the case of an arithmetic $M=Y(K)$,
we shall now  relate these to certain $V_{\xi}(s,K)$, defined as in \S \ref{Vxidefsec},  specifically for $\xi = 1, \eta^{-1}, \eta$, where 
   $\eta = \sqrt{z/\bar{z}}$ as before.  
\subsubsection{\texorpdfstring{$\xi$}{xi} trivial.}

When $\xi$ is the trivial character, 
elements of $V_{\xi}(s,K)^{K_{\infty}}$ descend to the quotient
$Y_B(K)$; this gives an isomorphism
$$ V_{\xi}(s,K)^{K_{\infty}} \stackrel{\sim}{\longrightarrow} C^{\infty}(s).$$ 

In this identification, the unitary structure on  $V(s)$ induces the unitary structure on $C^{\infty}(s)$  given by
$$\indexK^{-1} \int_{\height= Y} |f|^2, $$
where the integral is taken with respect to the measure induced by the hyperbolic metric. 
The result is independent of $Y$ when $s \in  i \R$. 
In particular it agrees up to scalars with the unitary structure previously described on $C^{\infty}(s)$. 

\subsubsection{\texorpdfstring{$\xi=\eta$}{xi=eta}.}
Recall the notation ${\p}$ from~\S~\ref{metricmeasurenotn}. 

There is an isomorphism
\begin{equation} \label{6522i} (V_{\eta}(s,K) \otimes {\p}^*)^{K_{\infty}} \stackrel{\sim}{\longrightarrow} 
 \Omega^+(s)   \end{equation} 
whose inverse is described thus:

Firstly an element of $\Omega^+(s)$ can be regarded as a 
$1$-form  $\omega$ on $\H^3 \times \G(\Afinite)/K$, 
 left invariant by $  \B(F)$, and satisfying 
\begin{equation} \label{moogle} \left( \begin{array} {cc} x & \star \\ 0 & 1 \end{array} \right)^* \omega  = |x|^{1+s} \eta(x) \omega 
= |x|^{ s}  x  \cdot \omega. \end{equation} 
Indeed,~\eqref{moogle} is just the same as asking that $\omega$ is a multiple of $y^s (dx_1 + i dx_2)$
on each $\mathbf{H}^3$ component (embedded into $\mathbf{H}^3 \times \G(\Afinite)/K$
via $z \mapsto (z, g_f K)$ for some $g_f \in \G(\Afinite))$.  
 
 \medskip
 
 Now we can construct the inverse to~\eqref{6522i}:
 Given  a  $1$-form $\omega$ on $\H^3 \times \G(\Afinite)/K$, 
 satisfying~\eqref{moogle} and left invariant by $\B(F)$ (so in fact
 by $\N(\A) \B(F)$), we associate to it the ${\p}^*$-valued
 function on $\G(\adele)$
 with the property that $$A(g) = (g^* \omega)|_{e}$$ for $g \in \G(\Adele)$.
 This therefore descends to an element 
$\left( V_{\eta}(s) \otimes {\p}^* \right)^{K K_{\infty}}$.

The inner product on  $\Omega^{\pm}(s)$
 for $s$ purely imaginary that corresponds to the inner product on $V_{\eta}(s)$ is given up to scalars by \begin{equation} \label{omeganormdef} \| \omega \|^2 := \int_{\mathrm{Ht} = T} \omega \wedge \overline{\omega}.\end{equation} 
In fact, the preimage $\mathrm{Ht}^{-1}(T)$ is a $2$-manifold, 
and the result is independent of $T$ because $\omega \wedge \overline{\omega}$
is closed and $[\mathrm{Ht}^{-1}(T)] - [\mathrm{Ht}^{-1}(T')]$ a boundary. 
The space $(V^{+}(s) \otimes {\p}^*)^{K_{\infty}}$ is equipped with a natural inner product, 
derived from that on $V_{\eta}(s)$ and the Riemannian metric on ${\p}$. 
Up to scalars, this coincides with the norm defined in~\eqref{omeganormdef}.\footnote{ Indeed, suppose that $\omega \in \Omega^{+}(s)$ corresponds to $A \in (V_{\eta}(s) \otimes {\p}^*)^{K_{\infty}}$; then 
$$ \| A\|^2 = \int_{ g \in  B(F) \backslash \G(\adele)^{(1)}} \| g^* \omega(1) \|^2 dg $$
But the signed measure $\| g^* \omega(1)\|^2 dg $, pushed down to $\partial Y_B(K) = B(F) \backslash \G(\adele)^{(1)}/K_{\infty} K$,  coincides with the measure $\omega \wedge \overline{\omega}$. 
Indeed, $\|g^* \omega(1)\|^2 = \|\omega(g)\|^2$; 
on the other hand, we have noted earlier that $dg$ pushes down to $\frac{dx_1 dx_2}{y^2}$
on each component of $\partial Y_B(K)$.}

\subsubsection{\texorpdfstring{$\xi = \eta^{-1}$}{xi = eta^-1}.}

Similarly, if we denote by $\Omega^{-}(s)$ the space of one forms on $Y_B(K)$
that are multiples of $y^{-s} (dx_1 - i dx_2)$ on each component, then 
$$ (V_{\eta^{-1}}(-s,K) \otimes {\p}^*)^{ K_{\infty}} \stackrel{\sim}{\longrightarrow}  \Omega^{-}(s).$$ 

The corresponding inner product on $\Omega^{-}(s)$ is given (up to scalars) again by~\eqref{omeganormdef}.

\subsection{Local intertwining operators.}
 
 For a nonarchimedean local field $k$, with ring of integers ${\OO}$, uniformizer $\varpi$, 
 residue field ${\OO}/\pi$ of size $p$, and $z \in \C^{\times}$, 
 let $\chi_z$ be the character $x \in k^{\times} \mapsto z^{\mathrm{val}(x)}$, 
and consider the local induced representation $\mathcal{I}(\chi_z)$:
   $$ I(z) (= I(\chi_z))  :=  \{f :\PGL_2(k) \rightarrow \C: f  \left( \left( \begin{array}{cc} x & * \\ 0 & 1 \end{array}\right) g \right) = |x|^{1/2}  z^{\mathrm{val}(x)} f(g).  \}$$
  Let $M: I(z) \rightarrow I(1/z)$ be  defined by the local analogue  
  of~\eqref{Intertwinerdef}, taking the measure on $N(k)$ to be the one 
  that assigns mass one to $N({\OO})$. 
  
 \medskip

  Let $f_0 \in I(z)$ be the $\PGL_2({\OO})$-fixed vector that takes 
  value $1$ on $\PGL_2({\OO})$.  Set $K_0 \subset \PGL_2({\OO})$
  to be the inverse image of the upper triangular matrices in $\PGL_2({\OO}/\varpi)$, 
  and let $f_1, f_2 \in I(z)$ be so that $f_1| \PGL_2({\OO})$ 
is the  characteristic function of $K_{0}$ and $f_1 + f_2 = f_0$. 
Then \begin{equation} \label{Scalarfactor} M(z) f_0 =   \frac{1-z^2/p}{1-z^2} f_0\end{equation} --- 
the constant of proportionality here may be more familiar as $\frac{L(0, \chi_z)}{L(1, \chi_z)}$ ---  
and the matrix of $M(z)$ in the basis $\{ f_1, f_2\}$ is given by:
\begin{eqnarray*} M(z)  &=&  \frac{1}{p (1-z^2)} \left(\begin{array}{cc}  z^2 (p-1) & 1-z^2 \\ p (1-z^2)  & p-1 \end{array} \right). \\ 
& =& \frac{1-z^2/p}{1-z^2} \frac{1}{p-z^2}  \left(\begin{array}{cc}  z^2 (p-1) & 1-z^2 \\ p (1-z^2)  & p-1 \end{array} \right). \\
\det M(z)  
 &= &  z^2 \left( \frac{1-z^2/p}{1-z^2} \right)^2 \cdot \left( \frac{p-1/z^2}{p-z^2} \right)
\end{eqnarray*}
Moreover the eigenvalues of $\frac{1-z^2}{1-z^2/p} M(z)$ 
are $\frac{pz^2-1}{p-z^2}$ and $1$.

 \subsection{Proof of theorem~\ref{scatteringmatrixtheorem}}
  \label{subsec:compar}   
 Let $\xi \in \{1, \eta, \eta^{-1}\}$.  
  Let $\Xi$ be a finite set of places,  disjoint from any level in $K$. 
Let $K' = K \cap K_{\Xi}$.   
We first analyze the scattering matrix for $Y(K')$, then for $Y(K)$, and then by comparing them
we obtain the statements of the Theorem (the Theorem follows, in fact, from
the case $\Xi= \{\q\}$ and $K = K_{\Sigma}$). 

\subsubsection{Analysis of the scattering matrix for  \texorpdfstring{$Y(K')$}{Y(K')}}
The scattering matrix for $Y(K')$ can be computed by computing the composite: 
 $$ \alpha(s):  V_{\xi}(K', 0)   \stackrel{\height^s}{\longrightarrow} V_{\xi}(K', s)   \stackrel{M}{ \longrightarrow}
V_{\xi^{-1}}(K', -s) \stackrel{\height^{-s}}{\longrightarrow}  V_{\xi^{-1}}(K', 0).$$

 After splitting $V_{\xi}$ into spaces $\mathcal{I}(\chi)$ according to the decomposition
~\eqref{ind-dec},  and   decomposing each $\mathcal{I}(\chi)$ into a tensor products of local representations $\mathcal{I}(\chi_v)$
(defined as in the prior section),  
computation of $\alpha(s)$ above reduces to computing
$$ 
\mathcal{I}(\chi_v)^{K_{0,v}}   \stackrel{\height^s}{\longrightarrow} \mathcal{I}(\chi_{v} |\cdot|_v^{s/2})^{K_{0,v}} 
 \stackrel{M}{ \longrightarrow} 
\mathcal{I}(\chi_{v}^{-1} |\cdot|_v^{-s/2})^{K_{0,v}}   \stackrel{\height^{-s}}{\longrightarrow}   \mathcal{I}(\chi_v^{-1})^{K_{0,v}}  $$ 
in particular for places $v \in \Xi$.  
 Here, $\height_v$ is the local analogue of the function $\height$, so that   $\height = \prod_{v} \height_v$, and $\chi_v$ is the local constituent of a character $\chi \in X$. 
 
 \medskip
 This composite -- with  respect to the standard basis $(f_1, f_2)$ for both
the source and target  prescribed in the prior subsection ---  is given by the matrix 
$M(z)$ of the prior section, where $z$  is the value of $\chi_v$ on a uniformizer of $F_v$, 
and $p = q_v$, the size of the residue field of $F_v$. 
\medskip

This proves (after interpretation)  the assertions of the Theorem concerning $\Sigma \cup \{\q\}$.

\subsubsection{Analysis of the scattering matrix for \texorpdfstring{$Y(K) \times \{1,2\}^{\Xi}$}{Y(K)}} 
Now let $\height'$ be the pullback of $\height$ to $Y_B(K) \times \{1, 2\}^{\Xi} $ under the isometry specified in
\S~\ref{levelraiseisometry}.

\medskip

We analyze the scattering matrix for $Y(K) \times \{1,2\}^{\Xi}$ with height $\height'$
under ~\eqref{levelraiseisometry}. 
   For this write $H = (\C^2)^{\otimes \Xi}$, and regard it as the space of functions on $\{1, 2\}^{\Xi}$;
   in particular  $H$ has an algebra structure by pointwise multiplication. 
   
   \medskip
   
   We may identify $V_{1}(s, K) \otimes H$
  with functions on   $\{ 1, 2 \}^{\Xi} \times Y_B(K)$. In particular, $\height'^s$
  may be considered as an element of $V_1(s,K) \otimes H$, and mutiplication by $\height'^s$
  defines a map $V_{\xi}(K , 0)    \otimes H\stackrel{\height'^s}{\longrightarrow} V_{\xi}(K, s)  \otimes H$.
  
    \medskip
    
 With these identifications, the scattering matrix for $Y(K) \times \{1, 2\}^{\Xi}$
  is given by the composite
 $$ \alpha(s): V_{\xi}(K , 0)    \otimes H\stackrel{\height'^s}{\longrightarrow} V_{\xi}(K, s)  \otimes H  \stackrel{M \otimes \mathrm{id}}{\longrightarrow} 
V_{\xi^{-1}}(K, -s) \otimes H \stackrel{\height'^{s}}{\longrightarrow}  V_{\xi^{-1}}(0) \otimes H.$$

Now we have an identification, as in \eqref{ind-dec},  $$V_{\xi}(K,0) \otimes H = \bigotimes_{\chi \in X(s)}\mathcal{I}(\chi)^K \otimes H
= \bigoplus_{\chi}  \left( \bigotimes_{v} \mathcal{I}(\chi_v)^{K_v} \otimes H_v \right),$$ where $H_v = \C^{\{1,2\}}$ for $v \in \Xi$
and $H_v = \C$ for $v \notin \Xi$.  With respect to this factorization,  $\mathrm{Ht}'$ factors
as the product $\bigotimes \left( \mathrm{Ht}_v \otimes  (1,\sqrt{q_v})\right)$, where $(1,\sqrt{ q_v}) \in H$ is
the vector that has value $1$ on the first $\C$ factor and value $\sqrt{q_v}$ on the second $\C$ factor
(See Remark~\ref{rem:heightchange}). 

Accordingly, computing $\alpha(s)$ reduces to computing the local maps
{\small 
      $$ 
 \begin{diagram}
\alpha_v(s): \mathcal{I}(\chi_v)^{K_v} \otimes \C^2  \stackrel{\height'^s}{\longrightarrow}   \mathcal{I}(\chi_v |\cdot|_v^{s/2})^{K_v} \otimes \C^2   \stackrel{M \otimes \mathrm{id}}{\rightarrow}   \mathcal{I}(\chi_v^{-1} |\cdot|_v^{-s/2})^{K_v} \otimes \C^2 \stackrel{\height'^{s} }{\longrightarrow} \mathcal{I}(\chi_v^{-1})^{K_v} \otimes \C^2
\end{diagram}
$$}
for $\chi_v$ a local constituent of some character $\chi \in X$, and $v \in \Xi$. 

Recall that $K_v = \PGL_2(\OO_v)$ for $v \in \Xi$, for we are supposing that $\Xi$ is disjoint from any level in $K$. 
Write $f_0$ for the element of $\mathcal{I}(\chi_v)$ that takes value $1$ on $K_v$. Let $e_1, e_2$ be the standard basis of $\C^2$. Then
$\alpha_v(s)$,   expressed with respect to the standard basis $f_0 \otimes e_1, f_0 \otimes e_2$
for the source and target --- is expressed by the matrix  $ M
  \left( \begin{array}{cc} 1 & 0 \\ 0 & q_v^{s}\end{array} \right) $, where $M$ is now the scalar by which the intertwining operator
  acts on $f_0$ (see~\eqref{Scalarfactor}). 

This proves (after interpretation)  the assertions of the Theorem concerning $\Sigma \times \{1,2\}$.  

\section{Modular symbols and the Eisenstein regulator} \label{sec:eisintegrality}
In this section, we work with an arbitrary open compact subgroup $K \subset \G(\Afinite)$, 
and discuss the issue of ``integrality'' of Eisenstein series. In other words: the Eisenstein series
gives an explicit harmonic form on $Y(K)$; which multiple of that form is cohomogically integral?

\medskip
This is actually not necessary for our results on Jacquet--Langlands correspondence, and also is vacuous
in the case of $K=K_{\Sigma}$ level structure in view of Lemma~\ref{vanishingboundaryhomology}.
However, it is interesting from two points of view: 
\begin{itemize}
\item[-]  The methods of this section should allow a complete understanding of the Eisenstein part of the regulator
-- in that  sense the present result properly belongs to Chapter~\ref{chapter:ch4};
\item[-] The methods lead to a bound on the cusp-Eisenstein torsion (notation of~\S~\ref{eisdef}), i.e., 
 the torsion in $H_1(Y(K), \Z)$
that lies in the image of $H_1(\partial Y(K), \Z)$.
\end{itemize}

\medskip

{\em We do not strive for the sharpest result}; we include this section largely to illustrate
the method: 
  We use ``modular symbols'' (that is to say: paths between two cusps, which lie a priori in Borel--Moore homology) to generate homology. When one integrates an Eisenstein series over a modular symbol, one obtains (essentially) a linear combination of $L$-values, but this linear combination involves some denominators.  The main idea is to  trade a modular symbol
for a sum of two others to avoid these denominators.   
Since working out the results
we have learned that similar ideas were used by a number of other authors, in particular Sczech
in essentially the same context.  We include it for self-containedness and also give a
clean adelic formulation. 

\medskip

One pleasant feature of the current situation is that the harmonic forms
associated to Eisenstein series have absolutely convergent integral over modular symbols. 
\medskip

  Denote by $\overline{\Z}$ the ring of algebraic integers. 
  In this section we shall define a certain integer $e$ in terms 
  of algebraic parts of certain $L$-values (\S~\ref{edefn}) and prove:
  
    \begin{theorem} \label{theorem:eisintegral}
Let $$ s \in  H^1(\partial Y(K), \overline{\Q})$$
be a Hecke eigenclass which lies in the image of $H^1(Y(K), \overline{\Q})$. 
If $s$ is integral (i.e., lies in the image of $H^1(\partial Y(K), \overline{\Z})$)
then $s$ is in fact in  the image of an  Hecke eigenclass $\tilde{s} \in H^1(Y(K),\overline{\Z}[\frac{1}{e}])$.
  \end{theorem}

The integer $e$ is the product of certainly easily comprehensible primes (as mentioned,
these could probably be removed with a little effort) and a more serious factor: the numerator
of a certain $L$-value.  We have not attempted for the most precise result. In particular, the method
gives a precise bound on the ``denominator'' of $\tilde{s}$, rather than simply a statement about the primes dividing this denominator.

In any case, let us reformulate it as a bound on torsion:  
\begin{quote} The image of $H_1(\partial Y(K), \Z)$ in $H_1(Y(K), \Z)$
contains no $\ell$-torsion unless $\ell$ divides $e$ or is an orbifold prime. 
 \end{quote}
\begin{proof} 
It is equivalent to consider instead the image of 
  $H^1(\partial Y(K), \Z)$ in $H^2_c(Y(K), \Z)$.

Let $\ell$ be a prime not dividing $e$. 
Take  $s \in H^1(\partial Y(K), \Z)$
whose image in $H^2_c(Y(K), \Z)$ is $\ell$-torsion. Let $s_{\Q}$ be the image of $s$ in $H^1(\partial Y(K), \overline{\Q})$. 
The theorem applies to $s_{\Q}$. 
Although $s_{\Q}$ may not be a Hecke eigenclass, we may write 
$s_{\Q} =\sum a_i s_i$, where the $s_i \in H^1(\partial Y(K), \overline{\Q})$ are Hecke eigenclasses and the $a_i$
have denominators only at $e$, as follows from the discussion after
Lemma~\ref{heckeactiononboundary}.
Then, according to the theorem,
there is $N$ such that $e^N s_i$ lifts to $H^1(Y(K), \overline{\Z})$. 
In particular, $e^N s $ maps to zero inside $H^2_c(Y(K), \overline{\Z})$. That proves the claim:
the image of $s$ in $H^2_c$ could not have been $\ell$-torsion. 
 \end{proof}

 \begin{remark*} \emph{One can in some cases obtain a much stronger result towards bounding the image of $H_1(\partial Y(K), \Z)$.
  }
 \emph{
  One method to do so is given by Berger \cite{Berger} who uses the action of involutions, related to the Galois automorphism, on the manifold in question. Another method which is apparently different but 
  seems to have the same range of applicability, is to construct enough cycles
  in $H_1(\partial Y(K), \Z)$ as the boundary of Borel--Moore classes in $H_{2, \bm}(Y(K) , \Z)$,
  by embedding modular curves inside $Y(K)$.   We do not pursue these methods further here; in any case the latter method, and likely also the former, is fundamentally limited to the case
  when all Eisenstein automorphic forms are invariant by the Galois automorphism of $F/\Q$. }
  \end{remark*}

  \subsection{}

    Let $X=X_0$ be the set of characters defined after~\eqref{ind-dec}, corresponding to $s=0$, 
  and let $X(K)$ be the subset with $\mathcal{I}(\chi)^K \neq 0$.  Let $\eta(z) = \sqrt{z/\bar{z}}$ as before.   
  For $\chi \in X$ let $f(\chi)$ be the norm of the conductor of $\chi$, so that \index{$f(\chi)$}
  $f(\chi)$ is a natural number.

  The formulation of the theorem involves an algebraic $L$-value. In order to make sense of this,
  we need to define first the transcendental periods that make it algebraic.

\subsection{Periods and integrality of \texorpdfstring{$L$}{L}-values}  \label{edefn}
    To each imaginary quadratic field $F$ we shall associate
a complex number $\Omega \in \C^{*}$ that is well-defined up to
units in the ring of algebraic integers:

   Let $E = \C/\OO_F$, considered as an elliptic curve over $\C$. 
Then there exists a number field $H \subset \C$ such that $E$ is defined over $H$, and
has everywhere good reduction there.  Let $\pi: \mathscr{E} \rightarrow \OO_H$ be the N{\'e}ron model of $E$ over 
$\OO_H$.   Then $\pi_* \Omega^1_{\mathscr{E}}$ defines a line bundle on $\OO_H$. 
{\em After extending $H$, if necessary, we may suppose this line bundle is trivial.}  
Having done so, we choose a generator $\omega$ for the global sections of that bundle;
this is defined up to $\OO_H^{\times}$.  On the other hand 
  the first homology group of the complex points of $E(\C)$ are  
free of rank one as an $\OO_F$-module. Choose a generator $\gamma$
for $H_1(E(\C), \Z)$ as $\OO_F$-module, and set
$$ \Omega := \int_{\gamma} \omega.$$

Let $\chi \in X$. 
    According to a beautiful result of Damerell \cite{Damerell2}, 
    $$L^{\alg}(1/2, \chi) := \frac{ L(1/2, \chi)}{\Omega}$$
    is algebraic, and, in fact, integral away from primes dividing $6 f(\chi)$. 
    The same is true for 
 $$ L^{\alg}(1/2, \chi^2) := \frac{ L(1/2, \chi^2)}{\Omega^2}.$$
      If $S$ is a finite set of places, we denote by $L^{\alg,S}(1, \chi)$
    the corresponding definition but omitting the factors at $S$ from the $L$-function.

 \subsection{Definition of the integer \texorpdfstring{$e$}{e}}
 If $\alpha$ is algebraic, we denote by $\num(\alpha)$ 
 the product of all primes $p$ ``dividing the numberator of $\alpha$'': that is, 
 such that there exists a prime ${\p}$ of $\Q(\alpha)$
 of residue characteristic $p$ so that the valuation of $\alpha$ at ${\p}$ is positive. 
 
 \medskip
 
 Let $S$ be the set of finite places at which the compact open subgroup $K$
 is not maximal, i.e. the set of $v$ such that $\PGL_2(\OO_v)$ is not contained in $K$. 
 Set $e_1 =  30 h_F  \cdot \disc(F) \cdot  \prod_{v \in S} q_v  (q_v-1)$.   Now put  $$e = e_1 
   \prod \num(  L^{\alg,S}(1, \chi^2));$$
   this involves the much more mysterious (and important) factor of the numerator of an $L$-function. 

As mentioned, we have not been very precise in our subsequent discussion; we anticipate
most of the factors in $e$ could be dropped except the numerator of the $L$-value.

 \subsection{Modular symbols} \label{ss:MS} 
 Let $\alpha, \beta \in \mathbf{P}^1(F)$ and  $g_f \in \G(\Afinite)/K$. Then the geodesic from $\alpha$ to $\beta$
 (considered as elements of $\mathbf{P}^1(\C)$, the boundary of $\mathbf{H}^3$),
 translated by $g_f$, 
 defines a class in 
 $H_{1,\bm}(Y(K))$ (a ``modular symbol'')  that we denote by $\langle \alpha, \beta ; g_f \rangle$. 
 Evidently these satisfy the relation
 $$\langle \alpha, \beta ; g_f \rangle + \langle \beta, \gamma;  g_f \rangle + \langle \gamma, \alpha ; g_f \rangle = 0.$$
 the left-hand side being the (translate by $g_f$ of the) boundary of the Borel--Moore chain defined
 by the ideal triangle by $\alpha, \beta, \gamma$.

If there exists $\epsilon \in \G(F)$ with the property that
$\epsilon . (\alpha, g_f K) = (\beta, g_f K)$  --
that is to say, $\epsilon \alpha = \beta$ and $\epsilon g_f K = g_f K$ -- we term the triple
$\langle \alpha, \beta; g_f \rangle$  {\em admissible}; for such
the image of $\langle \alpha, \beta, K \rangle$ under $H_{1,\bm} \rightarrow H_0(\partial)$
is {\em zero}, and consequently we may lift this to a class
$[\alpha, \beta, K] \in H_1(Y(K), \Z)$; of course we need to make a choice at this point;
the lifted class is unique only up to the image of $H_1(\partial Y(K), \Z)$.

 \begin{lemma*} Classes $[\alpha, \beta, K]$ associated to admissible triples,
together with the group $H_1(\partial Y(K), \Z)$, generates $H_1(Y(K), \Z)$.
\end{lemma*}

\begin{proof} 
Put $\Gamma = \G(F) \cap g_f K g_f^{-1}$. Then 
the connected component of $Y(K)$ containing $1 \times g_f$
is isomorphic to $M := \Gamma \backslash \H^3$.  

Fix $z_0 \in \H^3$. Then
$H_1(M, \Z)$ is generated by the projection to $N$ of geodesic between $z_0$ and $\gamma z_0$.  Taking $z_0$ towards the point $\infty$ yields the result:
the geodesic between $z_0, \gamma z_0$
defines an equivalent class in $H_{1,\bm}$ to the infinite geodesic $\infty, \gamma \infty$. 
 \end{proof}
 
  We now introduce the notion of the ``denominator'' of a triple, which will correspond
  to ``bad primes'' when we compute the integral of an Eisenstein series over it:
  
 We say a finite place $v$ is in the ``denominator'' of the triple $\langle \alpha, \beta; g_f \rangle$ 
 if the geodesic between $\alpha_v , \beta_v \in \mathbf{P}^1(F_v)$ inside
 the Bruhat-Tits building of $\G(F_v)$ does not passes through $g_f \PGL_2(\OO_v)$.  
 This notion is invariant by $\G(F)$, i.e., if we take $\gamma \in \G(F)$, then 
 $v$ is in the denominator of $\langle \alpha, \beta; g_f \rangle$ 
if and only if $v$ is in the denominator of $\gamma \langle \alpha, \beta; g_f \rangle$. 
 
 Note that, because $\G(F)$ acts $2$-transitively on $\mathbf{P}^1(F)$,
 we may always find $\gamma \in \G(F)$ such that
 $\gamma (\alpha, \beta) = (0,\infty)$; 
 thus any modular symbols is in fact equivalent to one of the form $\langle 0, \infty; g_f \rangle$
 for suitable $g_f$.     Now $v$ does not divide
 the denominator of the symbol  $\langle 0, \infty; g_f \rangle$
if and only if $g_f \in \diagA(F_v) \cdot \PGL_2(\OO_v)$.

\medskip

 The key point in our proof is the following simple fact about denominator avoidance.
 Let $p >3$ be a rational prime. 
 \begin{quote}
($\heartsuit$) We may write $\langle \alpha, \beta; g_f\rangle \in H_{1, \bm}$ as a sum of symbols
whose denominators do not contain any place $v$ with either $q_v$ or $q_v-1$
divisible by $p$.
 \end{quote}

 In our application the $q_v-1$ appears to be unnecessary. It may be useful in other contexts
 where one has a poorer {\em a priori} control on denominators.

\begin{proof}  
 By our prior remark it suffices to check the case when $\alpha = 0$ and $\beta=\infty$. 
Now write
$$ \langle 0, \infty; g_f \rangle = \langle 0, x ; g_f \rangle + \langle x, \infty; g_f \rangle.$$
for a suitable $x \in \mathbf{P}^1(F)$.  We claim we can choose $x$ so that, for every $v$ in the denominator of either symbol
 on the right-hand side, neither $q_v$ nor $q_v-1$ is divisible by $p$.
 \medskip

  In fact,  among the places with $g_v  \in   \PGL_2(\OO_v)$,
any place in the denominator of either $\langle 0,  x; g_f \rangle$
 or $\langle x , \infty; g_f \rangle$  satisfies $v(x) \neq 0$. 
 
\medskip

On the other hand, if $v $ is amongst
the finitely many places such that $g_v  \notin  \PGL_2(\OO_v)$,
the requirement that $v$ not bein the denominator of either symbol   amounts
  to requiring that $x$ have a certain (non-zero, non-infinite) reduction in $\mathbf{P}^1(\OO_v/\varpi_v^{n_v} \OO_v)$,
  for suitable $n_v$ (in fact, the distance in the building between $g_v \OO_v^2$
  and the geodesic between $0$ and $\infty$). 
  
  \medskip
  
So, our assertion follows from the following (taking $x=a/b$):  Given an ideal ${\frakn}_0$ and an element
  $\lambda \in (\OO_F/{\frakn}_0)^{\times}$, 
there exists elements $a, b \in \OO_F$ such that:
  \begin{itemize}
  \item $p$ doesn't divide $q(q-1)$,
where $q$ is the residue field size of any prime divisor ${\q}$ of either $a$ or $b$; 
  \item  $a \equiv \lambda b $ modulo ${\frakn}_0$. 
  \end{itemize}

In fact, by the analogue of Dirichlet's theorem, one can  lift any class in $(\OO_F/\mathfrak{m})^{\times}$, for any  ideal $\mathfrak{m}$, 
to a generator of a principal prime ideal. 
Write ${\frakn}_0 = {\frakn}_1 {\frakn}_2$
where ${\frakn}_1$ is prime-to-$p$ and ${\frakn}_2$
is divisible only by primes above $p$. Choose  $\bar{a}, \bar{b}  \in (\OO_F/p {\frakn}_2)^{\times}$ such that 
 $a \equiv \lambda b \ \mbox{mod}  {\frakn}_2$ and the norms of $a,b$
 (under the map $\OO_F/p {\frakn}_2 \rightarrow \OO_F/p \stackrel{\mathrm{N}}{\rightarrow} \Z/p\Z$) are not congruent to $1$ mod $p$. This can be done, for the
 image of the norm map
  $$ (\OO_F/ p)^{\times} \rightarrow (\Z/p)^{\times}$$
has size strictly larger than $2$, as long as $p>5$. 
Now take for $a$ a lift of 
 $$\lambda \times \bar{a} \in (\OO_F/{\frakn}_1) ^{\times} \times  \left( \OO_F/p{\frakn}_2\right)^{\times}
 \simeq (\OO_F/{\frakn}_1 {\frakn}_2 p)^{\times}$$
 to a generator $\pi$ of a principal prime ideal,
 and take $b$ similarly. 
   \end{proof}

 \subsection{Eisenstein integrals over modular symbols} \label{Eismod} 
 
 Let $X, \chi \in X$ and   $\mathcal{I}(\chi)$  as in~\eqref{ind-dec}
 and $\varF \in \Hom_{K_{\infty}}(\mathfrak{g}/\mathfrak{k}, \mathcal{I}(\chi))^K$. 
 We denote by $\bar{\varF} \in  \Hom_{K_{\infty}}(\mathfrak{g}/\mathfrak{k}, \mathcal{I}(\chi^{-1}))$
 the composition of $\varF$ with the standard intertwiner $\mathcal{I}(\chi) \longrightarrow \mathcal{I}(\chi^{-1})$ (defined by the analogue  of formula  ~\eqref{Intertwinerdef}).  Note that, for $\mathcal{I}(\chi)^K$ to be nonempty, \begin{equation} \label{chiunram} 
\mbox{ $\chi$ must be unramified at places $v$ for which $K_v$ is maximal.}\end{equation}

 We denote the automorphic representation which is the image of $\mathcal{I}(\chi)$
 under $\Eis$, by the letter $\Pi$. 
 
We define a $1$-form $\Omega( \Eis \  \varF)$ on $Y(K)$ as in~\S~\ref{param}.

Put $n(t) = \left( \begin{array}{cc} 1 & t \\ 0 & 1 \end{array}\right)$.
  The Lie algebra of $N$ as an $F$-algebraic group is naturally identified with $F$, with trivial Lie bracket (identify $x \in F$ with the one-parameter subgroup $t \mapsto n(xt)$).

  For $x \in F$ we write $\varF_x \in \mathcal{I}(\chi)$ for the evaluation of $\varF$
  at $x$ (considered in $\mathrm{Lie}(N) \rightarrow \mathfrak{g}/\mathfrak{k}$). 
 Therefore, for example, $\varF_1$ makes sense as an element of $\mathcal{I}(\chi)$; 
 it is in particular a function on $\G(\Adele)$.

 \begin{lemma*}
 If the restriction of $\Omega(\Eis \  \varF)$ to $\partial Y(K)$ is integral, i.e. takes
 $\overline{\Z}$-values on $1$-cycles, then $\varF_1(k)$ and $\bar{\varF}_1(k)$ takes 
 values in $N^{-1} \overline{\Z}$ whenever $k \in \PGL_2(\widehat{\OO})$
 where $N$ is divisible only by primes that are nonmaximal\footnote{By this we mean: primes $p$ such that there exists a finite place $w$ of residue characteristic $p$, 
 with $K_w \neq \PGL_2(\OO_{F, w})$.} for $K$ and primes dividing the discriminant of $F$. 
 \end{lemma*}
 \begin{proof}
   Let $ g_f \in \G(\Afinite)$ and $n(t) \in N(F)  \cap g_f^{-1} K g_f$; 
   we identify $t$ with an element of the Lie algebra of $N$.    The pair $(t,g_f)$ defines a $1$-cycle $\mathcal{C}$ in the cusp of $Y(K)$:
 namely the projection of $X \rightarrow n(t_{\infty} X) . g_f : 0 \leq X \leq 1$ to $Y(K)$;
 here $t_{\infty}$ is the image of $t$ under $F \hookrightarrow \C$. 
 
 We compute   $$\int_{\mathcal{C}} \Omega( \Eis \ \varF) = \varF_t(g_f) + \bar{\varF}_t(g_f);$$
 this is a routine computation from~\eqref{Eisconstant}.  Thus 
 $\varF_t(g_f) + \bar{\varF}_t(g_f) \in \overline{\Z}$ whenever $n(t) \in g_f^{-1} K g_f$. 
 
For fixed $g \in \G(\Afinite)$, the map  $x \mapsto \varF_x(g)$, considered as a map $F \rightarrow \C$, 
 is $F$-linear, whereas $x \mapsto \bar{\varF}_x(g)$ is $F$-conjugate-linear. 
  
For each place $v$, choose $m_v > 0$ such that the principal congruence subgroup of level $m_v$ is contained in $K_v$. Then $\{\lambda \in F_v:
\left(  \begin{array}{cc} 1 & \lambda \\ 0 & 1 \end{array} \right) \in g K g^{-1} \}$
 is stable by the order  $\Z_v + \varpi_v^{m_v} \OO_v$ (here $\Z_v$ is the closure of $\Z$ in $\OO_v$).
Let $\mathcal{O} = \{x \in \OO_F: x \in \Z_v + \varpi_v^{m_v} \OO_v \mbox{ for all } v\}$;
it's an order in $\OO_F$, and its discriminant is divisible only by the discriminant of $F$
and by primes at which $K$ is not maximal. 

Then $\{x \in F: n(x) \in g K g^{-1}\}$ is always $\mathcal{O}$-stable;
let $1, t$ be a $\Z$-basis of $\mathcal{O}$.   Then
$\varF + \bar{\varF} (x) \in \overline{\Z}$ and $t \varF + \bar{t} \bar{\varF}(x)\in \overline{\Z}$
implies that $\varF \in (t-\bar{t})^{-1} \overline{\Z} \subset \mathrm{disc}(\mathcal{O})^{-1} \overline{\Z}$. 
\end{proof}

In what follows we suppose $\varF$ to be integral in the sense above.  We shall verify then 
the following statement, which clearly imples Theorem~\ref{theorem:eisintegral}.  in view of~\S~\ref{ss:MS}, especially ($\heartsuit$).

\begin{quote} (*) 
The integral of $\Omega( \Eis \ \varF)$ over any modular symbol $\langle \alpha, \beta; g_f \rangle$
is integral away from primes dividing $e$ {\em and} primes in the denominator of the modular symbol. \end{quote}

We proceed to the proof of (*).  Let $X$ be as in 
defined after~\eqref{rsiso}, taking (in that discussion) the isomorphism
$\G(F_{\infty}) \cong \PGL_2(\C)$ to be the canonical one in this case. 
Thus $X$ is a a generator for the Lie algebra of $\A(\C)/\A(\C) \cap \PU_2$. 
 In fact $X = \left( \begin{array}{cc} 1/2 & 0 \\ 0 & -1/2 \end{array}\right)$
(where we identify $\mathfrak{pgl}_2(\C)$ with trace-free matrices).

\subsubsection{Some measure normalizations} \label{sss:lm} 
Fix an additive character $\psi$ of $\adele_F/F$: for definiteness
we take the composition of the standard character $\adele_{\Q}/\Q$
with the trace. Fix the measure on $\adele_F$ that is self-dual with respect to $\psi$, 
and similarly on each $F_v$. This equips $\N(F_v), \N(\adele)$ with measures via $x \mapsto n(x)$; the quotient measure
on $\N(\adele)/\N(F)$ is $1$. 

  For each $v$ we denote by $d_v$ the absolute discriminant of the local field $F_v$, 
so that the self-dual measure of the ring of integers of $F_v$ equals $d_v^{-1/2}$. 

We denote by $q_v$ the size of the residue field of $F_v$. 

\subsubsection{Normalization of  \texorpdfstring{$F(X)$}{F(X)}} 

Now put $\varF_X := \varF(X) \in \mathcal{I}(\chi)$;
we may suppose without loss of generality that it is a factorizable vector $\bigotimes f_v$.
Here, $f_{\infty}$  lies in the unique $H_{\infty} \cap K_{\infty}$-fixed
line in the unique $K_{\infty}$-subrepresentation isomorphic to $\mathfrak{g}/\mathfrak{k}$. 
We normalize it so that the map $\mathfrak{g}/\mathfrak{k} \rightarrow \mathcal{I}(\chi_{\infty})$
that carries $X$ to $f_{\infty}$ carries the image of $1 \in \mathrm{Lie}(N)$ to a function taking value $1$ at $1$. 

We can be completely explicit: We can decompose
$\mathfrak{g}/\mathfrak{k}$ as the direct sum
$\mathfrak{h} \oplus {\frakn}$, where these are (respectively) the Lie algebras
of $H_{\infty}$ and the image of the Lie algebra of $\N(\C)$. 
Then
$f_{\infty}$ can be taken to be given on $\PU_2$
by the function
$$ k \in \PU_2 \mapsto \mbox{${\frakn}$-component of } \mathrm{Ad}(k) X.$$
 and the image $f_{\infty}^Z$ of any other $Z \in \mathfrak{g}/\mathfrak{k}$
 is given by the corresponding function wherein we replace $X$ by $Z$. 

We may also now suppose that
  $f_v$ takes integral values 
  for all finite $v$.

\subsubsection{The Whittaker function of  \texorpdfstring{$\varF_X$}{varF_X}}
The Whittaker function of $\mathrm{Eis} (\varF_X)$, defined by the rule
$$ \int \mathrm{Eis}(\varF_X)(ng) \psi(n) dn = \prod_{v} W_v(g_v),$$
where  $W_v = \int_{x \in F_v} f_v(w n(x) g ) \psi(x) dx$. This integral is ``essentially convergent,'' i.e. it can be replaced by an integral over a sufficiently large compact set without changing its value. 
Thus $W_v$ is valued in $d_v^{-1} \cdot \overline{\Z}[q_v^{-1}]$,
  as follows from the fact that $\psi$ is $\overline{\Z}$-valued
  and the self-dual measure gives measure $d_v^{-1/2}$ to the maximal compact subring (see~\S~\ref{sss:lm}). 
  If $f_v$ is spherical one may compute more precisely that
   $W_v $ takes values  in $\frac{1}{d_v  L_v(1, \chi_v^2)} \overline{\Z}[q_v^{-1}]$. 
        
        \medskip
        
        Now, let $\omega_v$ be a  unitary character of $F_v^{\times}$ and $g_v \in \PGL_2(F_v)$. 
        
        \medskip
        
      Write $I_v := \int_{y \in F_v^{\times}} W(a(y) g_v) \omega_v(y) d^{\times} y$, where $d^{\times}y = \frac{dy}{|y|_v}$.   We assume that the defining integral is absolutely convergent.  
  Then \begin{equation} \label{Ivdenominator}   I_v \in  d_v^{-1} L_v(\frac{1}{2}, \Pi_v \times \omega_v)  \Z[q_v^{-1}].\end{equation} 
Indeed,   
   $L_v(\frac{1}{2}, \Pi_v \times \omega_v)^{-1} \int_{y \in F_v^{\times}} W_v(a(y)) \omega_v(y) d^{\times} y$
    is in effect a finite sum, with coefficients in $\Z[q^{-1}]$, of values of $W_v \omega_v$, 
    each of which lie in $ d_v^{-1} \overline{\Z}[q_v^{-1}]$. (The multiplication by $L_v$ has the effect of rendering the sum finite.) 
    
     Let $S$ be the set of finite places that lie in the denominator of the modular symbol,
 together with all places $v$ at which $K$ is not maximal. Thus (for $v \notin S$)
 $f_v$ is spherical  and $g_v \in A(F_v) \cdot \PGL_2(\OO_v)$
 From that we deduce that 
\begin{equation} \label{IvvnotinS} 
  I_v =  f_v(e)  L(\frac{1}{2}, \Pi_v \times \psi_v) / L(1, \chi_v^2)\end{equation}  for $v \notin S$.      
  Since $f_v$ is integral, $f_v(e) \in \overline{\Z}$ is an algebraic integer.

  \subsubsection{Integral of \texorpdfstring{$\Omega(\Eis \ \varF)$}{Omega} over a modular symbol} 
We are now ready to attack (*).  Without loss of generality we may suppose $\alpha = 0, \beta =\infty$. 

The integral of $\Omega(\Eis \ \varF)$ over the modular symbol $\langle 0, \infty; g_f \rangle$
  can be expressed as

\begin{eqnarray}  \label{psiunram} \int_{\langle 0, \infty; g_f \rangle} \Omega( \Eis \ \varF) &\stackrel{(a)}{=}& \frac{2^{?}}{h_F} \sum_{\omega} \int_{\A(F) \backslash \A(\adele)} 
  \mathrm{Eis}(\varF_X) (tg) \omega(t) d\mu_T 
  \\  &\stackrel{(b)}{=}&    \frac{2^?}{h_F}  I_{\infty}    \cdot  u\frac{ L(\frac{1}{2}, \Pi \times \omega)}{L^S(1, \chi^2)}  \cdot \prod_{v \in S} \frac{I_v}{L_v(\frac{1}{2}, \Pi \times \omega)}
  \\ \nonumber  &=&  {  \frac{2^?}{h_F}  I_{\infty}   \cdot  u \frac{ L^{\alg}(\frac{1}{2}, \chi \times \omega) L^{\alg}(\frac{1}{2}, \chi^{-1} \times \omega)}{L^{S, \alg}(1, \chi^2)}} \cdot  \prod_{v \in S} \frac{I_v}{L_v(\frac{1}{2}, \Pi \times \omega)} 
  \end{eqnarray} 
  where, at step (a),  we used~\eqref{cow2} in the case where $E = F \oplus F$, and summation is taken over characters $\omega$
       of $\diagA(\adele)/\diagA(F)$
       that are trivial  on $\diagA(\adele) \cap  \left( \diagA(F_{\infty}) \cdot g_f^{-1} K g_f \right)$;       at step (b) we unfolded the integral and used the results of the prior subsection for $v \notin S$, and $u \in \overline{\Z}$ is an algebraic integer. 
  
   Note that the index of the subgroup $\diagA(\Afinite) \cap  g_f^{-1} K g_f $ in the maximal compact of $\diagA(\Afinite)$
  is cancelled by a corresponding term in~\eqref{cow2}. Also note that any $\omega$ as in the prior summation is unramified
   at places $v$ that do not lie
       in the denominator of the modular symbol.

  Now consider the terms: 
  \begin{itemize}
  \item[-] We will see in a moment (\S~\ref{Iinftyeva}) that $I_{\infty}$ is integral; 

  \item[-]  The second  term $\displaystyle{\frac{ L^{\alg}(\frac{1}{2}, \chi \times \omega) L^{\alg}(\frac{1}{2}, \chi^{-1} \times \omega)}{L^{S, \alg}(1, \chi^2)}}$ is integral
  after being multiplied by the numerator of $L^{S, \alg}(1, \chi^2)$
  and a suitable power of $6 f(\chi)  f(\omega)$.
   As we observed  all primes dividing $f(\chi) f(\omega)$
  are primes below elements of $S$.     In particular, all primes dividing $f(\chi) f(\omega)$ also divide $e$. 
  
    \item[-] We saw in~\eqref{Ivdenominator} that the final  term is integral
  after being multiplied by a suitable power of $\disc(F)$ and $\prod_{v \in S} q_v$.  
  And, by choice of $e$, all primes dividing $\disc(F)$ or $\prod_{v \in S} q_v$ also divide $e$. 
  \end{itemize}
  This concludes the proof of (*), and so also of the Theorem.
    
  \subsubsection{Integrality of \texorpdfstring{$I_{\infty}$}{I}} \label{Iinftyeva} 
 
$$ w n(x) =  \left( \begin{array}{cc} 0 & 1 \\ -1 & -x \end{array}\right)  = 
   \left(  \begin{array}{cc} 1 & -\bar{x} \\ 0 & 1 \end{array}\right) 
  \left(  \begin{array}{cc}  \frac{1}{ \sqrt{1+|x|^2}} & 0 \\ 0 &  \sqrt{1+|x|^2}  \end{array}\right)  k$$ 
and we compute that the projection to ${\frakn}$
of $\mathrm{Ad}(k) X$ is given by simply $2\bar{x}$.
Thus $f_{\infty}(w n(x)) =  \bar{x} (1+|x|^2)^{-2} $, and
we wish to compute: 
\begin{align}  I_{\infty} =  \int_{x \in F_{\infty}} f( wn(x) a(y)) \psi(x) \chi(y) d^{\times} y 
       = \int_{x \in F_{\infty}} f(a(y^{-1}) w n(x)) |y|_v \psi(xy)   d^{\times} y 
       \\ = \int_{x}  y  \frac{x}{|x|^2+1} \psi(xy) \chi(y) dx d^{\times} y 
       = \int_{x} y \frac{\bar{x}}{(|x|^2+1)^2}  \psi(xy) dx d^{\times} y\end{align}       
    Now, the Fourier transform of $y/|y|^2$ (by homogeneity) is $\gamma \bar{x}/|x|^2$
    with respect to the self-dual Haar measure; since the measure on $F_{\infty}^{\times}$
    differs from this by a factor $\zeta_{\C}(1) = \pi^{-1},$
    and therefore our integral becomes
    $$\pi^{-1}  \int_{x} \frac{1}{(|x|^2+1)^2} = 1.$$

This concludes the proof of $p$-integrality.

\section{Comparing Reidmeister and analytic torsion:  the main theorems}  \label{analysissec}

 We have previously defined notions of analytic torsion $\analT(Y)$  and  Reidemeister torsion $\RT(Y)$ in the non-compact case. It is an interesting question to relate them,
 but it is not quite necessary for us to completely carry this out.  
 
 We will consider the case where we compare levels  
 $\Sigma, \Sigma \cup \{\p\}, \Sigma \cup \{\q\}$ and $\Sigma \cup \{\p, \q\}$. 
 More generally, we do not restrict to the case of $K = K_{\Sigma}$;
 let $K$ be an arbitrary compact open subgroup of $\G(\Afinite)$
 which is maximal at $\p, \q$, and 
set
\begin{align} \label{ypq}
Y = Y(K) \times \{1,2\}^2,   \ \ & Y_{\q} = Y( K \cap K_{\{\q\}}) \times \{1,2\}; \\ \nonumber Y_{\p} = Y(K \cap K_{\{\p\}}) \times \{1,2\},  \ \ & Y_{\p\q} = Y( K \cap K_{\{\p,\q\}}). \end{align}
Note that all statements and proofs will carry over to the case when $\{\p,\q\}$
is an arbitrary set of primes of size $\geq 2$. 

  These will always be equipped with height functions normalized as follows:
  equip $Y(K_{\q})$ and $Y(K_{\p\q})$ with the heights of~\eqref{infgeom},
  and $Y_{\q} = Y(K_{\q}) \times \{1,2\}$ with the height that restricts to this on each component,
  and similarly for $Y_{\p\q}$; 
  this specifies a height on $Y_{\q}, Y_{\p\q}$. 
Now recall that we have specified (\S~\ref{ybybxi}) an isometry
\begin{eqnarray*} Y(K_{\q})_B  &\cong& Y(K)_B \times \{1, 2\}, 
 Y(K_{\p\q}) \cong Y(K_{\p})_B \times \{1,2\} \end{eqnarray*} 
 and now pull back heights via this isometry to obtain heights on $Y$ and $Y_{\p}$.
 This definition seems somewhat idiosyncratic, but it has been chosen so that:
  
 \begin{quote}
 The cusps of $Y$ and $Y_{\q}$ are isometric in a fashion preserving the height function; the cusps of $Y_{\p}$ and $Y_{\p\q}$
 are isometric in a fashion preserving the height function. 
 \end{quote} 
What we will actually show is:

\begin{theorem}  \label{thm:rtatsplit}
Notation and choices of height function as above, 
 \begin{equation}  \label{rtatsplit}  
 \frac{ \RT (Y_{\p\q}) } {\RT(Y_{\p})}  
 \frac{\RT(Y)}{\RT(Y_{\q}) }  = 
 \frac{ \analT (Y_{\p\q}) } {\analT(Y_{\p})}  
 \frac{\analT(Y)}{\analT(Y_{\q}) },
 \end{equation}
 up to orbifold primes. 
  \end{theorem} 
   ~\eqref{rtatsplit} will prove to be easier than understanding the relation between
each separate term $\analT/\RT$. Indeed, it is well-known 
 that it is much easier to handle {\em ratios} of regularized determinants
than the actual determinants.

The next Corollary properly belongs in the next Chapter, but we give it here
to simplify matters there:
 
\begin{corollary} \label{corollary610} Let $Y'$ be the arithmetic manifold associated to the division algebra $D'$
 ramified precisely at $\p, \q$ and the corresponding open compact subgroup $K'$ of $\G'(\Afinite)$. 
Then, up to orbifold primes, 
\begin{equation} \label{comparisonB}
  \frac{ \RT (Y_{\p\q}) } {\RT(Y_{\p})}  
 \frac{\RT(Y)}{\RT(Y_{\q}) }  = 
\RT(Y').\end{equation}
\end{corollary} 
By ``corresponding open compact subgroup $K'$" we mean that
 $K'_{v}$ is the image of units in a maximal order for $v$ dividing $\p, \q$, and 
 $K'_v$ corresponds to $K_v$ for $v$ not dividing $\p, \q$.  In particular,
 if $K$ is simply a level structure of type $K_{\Sigma}$, then the corresponding $K'$
 is also the level structure $K'_{\Sigma}$, and the $Y, Y'$ are a Jacquet--Langlands pair in the sense
 of~\S~\ref{section:discusscongruence}. 
 
  \begin{proof}
 Fix any $0 \leq j \leq 3$. 
 If $\phi, \phi_{\p}, \phi_{\q}, \phi_{\p\q}$ are the determinants of the scattering matrices for $j$-forms on  $Y ,Y_{\p}, Y_{\q}, Y_{\p\q}$
 then our prior computations (specifically~\eqref{scatmatrix1} and~\eqref{scatmatrix2}, as well as assertion (4) of Theorem \ref{scatteringmatrixtheorem}) show that
\begin{equation}  \label{ar12ef} \left(  \frac{\phi'}{\phi} \right) - \left(  \frac{ \phi_{\p}' }{\phi_{\p}} \right)  = 
  \left(  \frac{\phi'_{\q}}{\phi_{\q}} \right) - \left(  \frac{ \phi_{\p\q}' }{\phi_{\p\q}} \right).
\end{equation}   
  Also, we have seen that, if $\Psi$ is the scattering matrix on functions, then, in similar notation, 
  \begin{equation} \label{ef12ar}  \mathrm{tr}  \ \Psi(0) - \mathrm{tr} \  \Psi_{\p}(0) = 
   \mathrm{tr} \  \Psi_{\q}(0) - \mathrm{tr} \   \Psi_{\p\q}(0).
\end{equation}

   Now consider our  definition~\eqref{atsecondef},~\eqref{ASspec} of analytic torsion in the non-compact case.
~\eqref{ar12ef} and~\eqref{ef12ar} show that all the terms involving scattering matrices cancel out, and we are left with: 
\begin{equation} \label{alternatingratioATone}  \log   \frac{ \analT (Y_{\p\q}) } {\analT(Y_{\p})}  
 \frac{\analT(Y)}{\analT(Y_{\q}) }   =  
  \frac{1}{2}  \sum_{j} (-1)^{j+1} j \log \det{}^* (\Delta_j^{\new}), \end{equation} 
where $\Delta_j$ denote the Laplacian on $j$-forms on $Y(K_{\p\q})$ but restricted
to the {\em new subspace} of $j$-forms on $Y(K_{\p\q})$. 
The Jacquet--Langlands correspondence (skip ahead to 
\S~\ref{ss:JLclassic} for a review in this context) asserts
the spectrum of $\Delta_j^{\new}$    coincides with the spectrum
of the Laplacian $\Delta_j$ on $j$-forms, on $Y'(K)$. Thus the right
hand side equals $\log \analT(Y'(K))$,  and applying the Cheeger-M{\"u}ller theorem
to the compact manifold $Y'(K)$ we are done. 
\end{proof}

We deduce Theorem~\ref{thm:rtatsplit} from the following result. 
   
Suppose we are given two hyperbolic manifolds $M, M'$
such that $M_B$ and $M'_B$ are isometric; we suppose that 
we have an isometry $\sigma: M_B \rightarrow M'_B$ 
and height functions on $M, M'$ that match up with respect to $\sigma$. 
 We suppose that $M, M'$ are of the form $Y(K)$ for some $K$; 
 {\em this  arithmeticity assumption is almost surely unnecessary to the proof}, but it allows us to be a bit lazy at various points in the proof (e.g., we may appeal to standard bounds for the scattering matrix.)

\begin{theorem}  (Invariance under truncation) \label{theorem:invtrunc}    \label{THEOREM:INVTRUNC}As 
 $Y \rightarrow \infty$, 
$$\begin{aligned} \label{desid}    \log  \analT(M)& -&  \log \analT(M_Y) & -&   
 ( \log \RT(M)& -& \log \RT(M_Y)&  ) \\   
- (\log \analT(M') & -&  \log \analT(M'_Y))
&+& (\log \RT(M') &-& \log \RT(M'_Y) & )
\longrightarrow 0,\end{aligned}$$
where we understand the computation of analytic torsion on $M_Y$ to be with respect to absolute boundary conditions. 
 \end{theorem}
 \index{boundary conditions}
When dealing with manifolds with boundary, we {\em always}  compute analytic torsion with respect to absolute boundary conditions in what follows.   
 
Indeed, the following observation was already used in Cheeger's proof: when studying the effect of a geometric operation
 on analytic torsion,  it can be checked more easily that this effect is {\em independent
 of the manifold on which the geometric operation  is performed.}
  In Cheeger's context, the surgery
is usual surgery of manifolds. 
In our context, the geometric operation will be to truncate the cusp of a hyperbolic manifold. 
Cheeger is able to obtain much more precise results by then explicitly computing
the effect of a particular surgery on a simple manifold; we do not carry out the analogue of this step, because we do not need it.
This missing step does not seem too difficult: 
the only missing point is to check the short-time behavior of the heat kernel near the boundary of the truncated manifold; at present we do not see the arithmetic consequence of this computation.

 The actual proof is independent of the Cheeger--M{\"u}ller theorem, although
 it uses ideas that we learned from various papers of both of these authors. 
It will be given in~\S~\ref{theorem:invtruncproof} after
 the necessary preliminaries about small eigenvalues of the Laplacian on $M_T$, which
 are given in~\S~\ref{smalleigenvalues}.  We warn that the corresponding statement is not true for the $\logdetR$ or the $\RT$ terms alone, although we precisely compute both: indeed, both terms individually diverge as $Y \rightarrow \infty$. 
 
 For now we only say
the main idea of the proof:
 \begin{quote} Identify the continuous
parts of the trace formula for $M, M'$ as arising from ``newly created'' eigenvalues on the truncated manifolds $M_Y, M'_Y$. \end{quote}

During the proof we need to work with $j$-forms for $j \in \{0,1,2,3\}$. Our general policy
will be to treat the case $j=1$  in detail. The other cases are usually very similar,
and we comment on differences in their treatment where important. 
 
\subsection{Proof that Theorem~\ref{theorem:invtrunc} implies Theorem~\ref{thm:rtatsplit}.}  \label{rtatsplitfollows}

 For $M = Y(K), Y(K_{\p}), \dots$ or  any compact manifold with boundary, put 
$$ \alpha(M) := \frac{\analT(M)}{\RT(M)};$$
in the cases of $Y(K)$ etc. the definitions of $\analT$ and $\RT$ are those of
\S~\ref{sec:rtatnc} and depend on the choice of height function.  Recall  that we understand $\analT$ to be computed with reference to absolute boundary conditions in the case where $M$ is a manifold with boundary.

For the moment we suppose that all these manifolds are genuine manifolds, i.e., do not have orbifold points. 
We discuss the necessary modifications in the general case in~\S~\ref{orbifoldtorsion2}.

 According to Cheeger \cite[Corollary 3.29]{Cheeger} $\alpha(M)$ depends only on the germ of the metric of $M$ near the boundary. 
The  metric germs at the boundary being identical for $Y, Y_{\q}$ (see comment before~\eqref{rtatsplit}), and similarly for $Y_{\p}, Y_{\p\q}$, 
  Cheeger's result shows that:
 \begin{equation}  \label{rtatsplit2}  
\frac{ \alpha(Y_{ T}) }{ \alpha(Y_{\q,T})}= 1 \mbox{ and } \frac{\alpha(Y_{\p\q,T})}{\alpha(Y(K_{\p,T})  ) } = 1.
 \end{equation} 
 Here, we understand $Y_T$ to be the manifold with boundary obtained by truncating at height $T$, etc. Now Theorem~\ref{theorem:invtrunc} implies, after taking the limit, 
\begin{equation} \label{analC} 
 \frac{\alpha(Y)}{\alpha(Y_{\q})}  \cdot \frac{\alpha(Y_{\p\q}) }{\alpha(Y_{\p})} =1;
\end{equation}
this proves~\eqref{rtatsplit}. 

 \subsection{Modifications in the orbifold case}    \label{orbifoldtorsion2}

Our discussion in the prior section~\S~\ref{rtatsplitfollows}
concerned only the case when $Y(K)$ is strictly a manifold, i.e. has no orbifold points. 
Let us see why  ~\eqref{rtatsplit}
remains true in the general case --- i.e., when $Y(K)$ has orbifold points --
up to a factor in $\mathbf{Q}^{\times}$ ``supported at orbifold primes'', i.e. whose numerator and denominator is divisible by' primes dividing the order of some isotropy group. 

First of all, we note that the proof of Theorem~\ref{theorem:invtrunc} is entirely valid in the orbifold case. Indeed, the proof of this theorem is in fact entirely about properties of Laplacian eigenvalues and harmonic forms,
and is insensitive to whether or not $M$ and $M_Y$ are manifolds or orbifolds.

Now the reasoning of the  prior section~\S~\ref{rtatsplitfollows}
goes through word-for-word in the orbifold case, {\em except} for~\eqref{rtatsplit2}; this remains valid only  when $1$ is replaced by 
a rational number supported only at orbifold primes, according to  $(\diamondsuit)$  below. Note that that indeed the orbifolds $Y(K)$ as well as any orbifolds obtained by truncating it are indeed global quotients, i.e., admit finite covers which are manifolds. 

We now check
\begin{quote}$(\diamondsuit)$
Suppose given two manifolds with boundary $\tilde{M}_1, \tilde{M}_2$,
isometric near the boundary, and whose boundaries have zero Euler characteristic. 
Suppose also given a finite group $\Delta$ acting on $\tilde{M}_i$
preserving this isometry.  Define $M_i$ as the orbifold quotient $ \tilde{M}_i/\Delta$. 

Then  the ratios of analytic and Reidemeister torsion for $M_1, M_2$
differ by a rational number supported only at orbifold primes.
\end{quote}

We proceed similarly to the previous discussion (\S~\ref{orbifoldtorsion1}), but use the work of L{\"u}ck \cite{Lueck}. Note that L{\"u}ck's torsion is the square
of that of \cite{LottRoth}, which causes no problem for us; this is stated in (4.8) 
of \cite{Lueck}. 
In fact, L{\"u}ck studies such equivariant torsion under the following assumptions:
\begin{itemize}
\item[(i)]
the metric near the boundary is a product;
\item[(ii)] A certain assumption of {\em coherence};
the latter is automatically satisfied when the flat bundle under consideration is trivial (or indeed induced by a $G$-representation). 
 \end{itemize}

Indeed  we may first compatibly deform the metrics on the $M_i$ so the metric near the boundary is a product. 
The ratio of Reidemeister and analytic torsion   changes in the same way for $M_1, M_2$
when we do this (this argument is just as in \cite[Corollary 3.29]{Cheeger}). 

This does not change the purely topological statement the $M_i$ can be expressed as global quotients, and we then apply L{\"u}ck's results.

 In this case,  the equivariant Reidemeister and equivariant analytic torsions are not literally {\em equal}; they differ by two terms, one measuring the failure of equivariant Poincar{\'e} duality,
 and the second proportional to the Euler characteristic of the boundary. 
 The latter term vanishes in our case, as we will be dealing with certain truncated hyperbolic $3$-manifolds; their boundaries are unions of $2$-dimensional tori. 
As for the former term: the term measuring the failure of equivariant Poincar{\'e} duality is supported entirely at primes dividing the order of some isotropy group, that is to say,
 orbifold primes for $\tilde{M}/\Delta$.    That follows from \cite[Proposition 3.23(a) and Proposition 5.4]{Lueck}
 and the fact that Poincare duality holds ``away from orbifold primes'' in
 the sense of~\S~\ref{sss:dualityorbifold}.

\section{Small eigenvalues} \label{smalleigenvalues}

We are in the situation of Theorem~\ref{theorem:invtrunc}.  We denote by $\Phi^+, \Phi^-$
the scattering matrices for $1$-forms on $M$, and by $\Psi$ the scattering matrix for functions;
these notations are as in~\S~\ref{subsec:einsteinintro}.

Write $\Delta_{M_Y}$ for the Laplacian on  $1$-forms, on the truncated manifold $M_Y$, 
with absolute boundary conditions. 
The key part of the proof of Theorem~\ref{theorem:invtrunc} is the following result, which describes
the small eigenvalues of $\Delta_{M_Y}$ precisely. 
Through the proof we will make the abbreviation $$ \Tmax := (\log Y)^{100}.$$  \index{$\Tmax$}
This determines the upper threshhold to which we will analyze eigenvalues, i.e.
we will not attempt to compute eigenvalues of $M_Y$ directly that are larger than $\Tmax^2$.

\begin{theorem} \label{prop:SE} [Small eigenvalues.]
Put
$$f(s) = \det( \mathrm{Id} - Y^{-4s} \Phi^{-}(-s) \Phi^+(s)),$$ \index{$f(s)$}

Let $0 \leq a_1 \leq a_2 \leq \dots $ be the union
(taken with multiplicity) of the non-negative  real roots  of $t \mapsto f(it)$ and the set
$$\{t \in \R_{>0}: t^2 \mbox{ is the eigenvalue of a cuspidal co-closed $1$-form on $M$}\},$$
both sets taken themselves with multiplicity. 
Set $\bar{\lambda}_j = a_j^2$.  

Let $0 \leq \lambda_1 \leq \dots$ be the eigenvalues, with multiplicity, of 
$$\Delta_{M_Y} | \ker d^*,$$
i.e. Laplacian eigenvalues on  co-closed $1$-forms on $M_Y$ satisfying absolute boundary conditions. 

There is $a > 0$ such that, with, $\delta = \exp(-a Y)$
and for $Y$ sufficiently large (this depending only on $M$):
\begin{enumerate}
 
\item[(a)]
(This is trivial:) If $b =\dim H^1(M, \C)$, then $\lambda_1 =\dots = \lambda_b = 0$
and $\lambda_{b+1} > 0$.  The same is true for $\bar{\lambda_i}$.

\item[(b)]
 $|\lambda_j  - \bar{\lambda}_j| \leq \delta$ for any $j$ with   $\sqrt{|\lambda_j|} \leq \Tmax$.

 \item[(c)] The same assertion holds for the eigenvalues of $\Delta$ on  $0$-forms (=co-closed $0$-forms) 
 with absolute boundary conditions, replacing $f(s)$ by the function 
  $$g(s) = \det( \mathrm{Id}  + \frac{1-s}{1+s} Y^{-2s} \Psi(s)),$$
  but taking only roots of $g(it)$ where $t$ is strictly positive; 
also replace $b$ by $H^0(M, \C)$ (and $\bar{\lambda_j}$ by $1-a_j^2$.)

  \item[(d)] A similar assertion holds for eigenvalues on co-closed $2$-forms with absolute boundary conditions (equivalently: for eigenvalues on closed $1$-forms with relative
  boundary conditions\footnote{Note that the nonzero such eigenvalues are just the nonzero eigenvalues of the Laplacian on $0$-forms with relative boundary conditions. This is why $\Psi$ intervenes.}) where we now replace the function $f(s)$ by the function
  $$g'(s) = \det(\mathrm{Id} + Y^{-2s} \Psi(s))$$
  we replace $b$ by $\dim H^2(M, \C)$,
  and we  make the  following modification
 to the definition  of $\bar{\lambda}_j$:

  \begin{quote} Let $0 < u_1 \leq \dots \leq u_h < 1$ be the  roots, with multiplicity, of $g'(t)$ for $t \in (0,1]$ together with parameters of cusp forms $\{t \in [0,1]: 1 -t^2 \mbox{ is the eigenvalue
of a cuspidal $3$-form on $M$} \}$,  and let $a_j = i u_j$ for $1 \leq j \leq h$. 

Let $0 < a_{h+1} \leq a_{h+2} \leq \dots$ be the positive real roots of $t \mapsto g'(it)$ 
together with parameters of cusp forms; 
and put $\bar{\lambda}_j  = 1 +a_j^2$. 
  \end{quote}
  {\em Warning: There exist roots $a_j$ of $g' = 0$ very close to $t=1$ related to the residue of $\Psi$;   see  discussion of 
~\S~\ref{sec:modifications}.}

\end{enumerate}
\end{theorem}

What is perhaps not obvious from the statement is {\em that this description ``matches''
in a very beautiful way part of the trace formula for $M$}. 
The proof of the theorem follows along lines that we learned from work of W. M{\"u}ller (see \cite{MEta}) 
and we present it in~\S~\ref{sec-muller-thm}.   We first give the proof of Theorem~\ref{theorem:invtrunc},
using Theorem~\ref{prop:SE}. In fact, we will prove in detail only (a), (b), the other parts being similar; we will however discuss in some detail the ``extra'' eigenvalues of case (e). 
These eigenvalues play an important role in accounting for the failure of Poincar{\'e} duality for the non-compact manifold $M$. 

For the moment we try to only give the main idea by describing:

\subsubsection{Plausibility argument for (b) of the Theorem.} \label{s:plaus} Consider, for $\omega \in \Omega^+(0) $ and $\bar{\omega} \in \Omega^{-}(0)$, 
 the $1$-form $F = E(s, \omega ) +  E(-s, \bar{\omega})$; we shall
attempt to choose $\omega, \bar{\omega}$ so it almost satisfies absolute boundary conditions when restricted to $M_Y$.   Then in any cusp
$$F \sim  y^s \omega + \Phi^{+}(s) y^{-s} \omega +  y^s \baromega + y^{-s} \Phi^{-}(-s) \baromega.$$ 
The sign $\sim$ means, as before, that the difference is exponentially decaying.  
The contraction of $F$ with a boundary normal is automatically of exponential decay;  on the other hand,
up to exponentially decaying factors,  $$s^{-1} \cdot dF  \sim (  y^s  \omega- y^{-s} \Phi^{-}(-s) \bar{\omega})  \wedge \frac{dy}{y} + (y^s \bar{\omega} - y^{-s} \Phi^+(s) \omega) \wedge \frac{dy}{y} .$$

Note that the $s^{-1}$ on the left-hand side suggests that we will have to treat the case when $s$ is near zero separately,
and indeed we do this  (\S~\ref{combin}).

In order that the contraction of this form with a normal vector at the boundary be zero, we require
that both the ``holomorphic part''  $(  y^s  \omega- y^{-s} \Phi^{-}(s) \bar{\omega})  $
and the ``anti-holomorphic part'' $ (y^s \bar{\omega} - y^{-s} \Phi^{-s} \omega) $ be zero at the boundary, that is to say,  
\begin{equation} \label{omegaomegabar} \omega   = Y^{-2s} \Phi^{-}(-s) \bar{\omega} \mbox{ and } Y^s \bar{\omega} = Y^{-s} \Phi^+(s) \omega
, \end{equation} 
which imply that  $\omega$ is fixed by $Y^{-4s} \Phi^{-}(-s) \Phi^+(s)$. 
For such $\omega$ to exist, we must indeed have $f(s) =\det( \mathrm{Id} - Y^{-4s} \Phi^{-}(-s) \Phi^+(s)) = 0$. 

Of course, this analysis is only approximate, and we need to both verify
that a zero of $f(s)$ gives rise to a nearby eigenfunction satisfying absolute boundary conditions,
and, conversely, that any eigenfunction satisfying absolute boundary conditions ``came from''
a zero of $f(s)$.

For later use, it will also be helpful to compute the inner product $\langle \wedge^Y F, \wedge^Y F \rangle$ notations as above. 
We have:

\begin{equation} \label{Fnorm} \| \wedge^Y F\|^2_{L^2(M)} =
\|  \wedge^Y E(s, \omega)\|^2_{L^2(M)}  + \|  \wedge^Y E(-s, \bar{\omega}) \|^2_{L^2(M)} 
\end{equation} 
$$  + \mathrm{Re} \
\langle  \wedge^Y E(s, \omega), \wedge^Y E(-s, \bar{\omega}) \rangle.$$
Assuming that \eqref{omegaomegabar} holds, the second line {\em vanishes identically.}
Indeed, according to the Maass-Selberg relation (see before
\eqref{MSR1F}), the inner product in question equals
$\frac{Y^{2s}}{2s} \langle \omega, \Phi^{-}(s) \bar{\omega} \rangle + 
\frac{Y^{-2s}}{-2s} \langle \Phi^+(s) \omega, \bar{\omega} \rangle,$
and \eqref{omegaomegabar} means that this is purely imaginary . 
Therefore (after some simplifications) one obtains  
\begin{equation} \label{Fnorm2} \| \wedge^Y F\|^2 =  \langle (4 \log Y - A(s)^{-1} \frac{dA}{ds}) \omega, \omega \rangle, \end{equation}
where we have abbreviated $A(s) := 
 \Phi^{-}(-s) \Phi^{+}(s)$. This expression will recur later. 

   \subsection{Analysis of zeroes of \texorpdfstring{$f(s)$}{f(s)} and \texorpdfstring{$g(s)$}{g(s)}.} \label{subsec:fganal}
Before we proceed (and also to help the reader's intuition) we analyze the zeroes of $f(s)$ and $g(s)$.  We discuss $f(s)$, the case of $g$ and $g'$ being similar.

Write $$ A(s) = \Phi^{-}(-s) \Phi^{+}(s): \Omega^+(s) \rightarrow \Omega^+(s).$$  \index{$A(s)$}
Recall 
   $f(s) = \det(1-Y^{-4s} A(s))$.   Let $\hrel$ be the size of the matrix $A(s)$, i.e., the total number of relevant cusps of $M$. 
   
   Recall that a continuously varying family of matrices has continuously varying eigenvalues (in the ``intuitively obvious'' sense; formally, every symmetric function in the eigenvalues is 
   continuous). Consequently, 
we may (non-uniquely) choose continuous functions $\nu_i : \R  \rightarrow \R/\Z$ 
   with the property that $ \{ \exp(  i \nu_i(t)) \}_{1 \leq i \leq h}$
   is the set of eigenvalues of the unitary matrix $A(it)$, where the eigenvalues are
    counted with multiplicity, for every $t \in \R$.
\footnote{   Indeed, it is sufficient to check that this is so in a neighbourhood of every point;
   shrinking the neighbourhood appropriately,  we may suppose (multiplying by a constant if necessary) that
   $-1$ is never an eigenvalue, and  now simply take the $\nu_i$ to be the arguments of the eigenvalues, taken so that $\nu_1 \leq \nu_2 \leq \dots \leq \nu_h$. }
   At $s=0$, we have $A(0) = 1$ so every eigenvalue is $ 1$; we normalize, accordingly,
   so that $\nu_i(0) =0$ for every $i$. 

  In what follows, we write $\frac{A'}{A}$ for $A^{-1} \cdot \frac{d A(s)}{ds}$
  when $A$ is a matrix-valued function of $s$.

The functions $\nu_i$ are real analytic away from a discrete set of points. 
Moreover, $-\nu_i'  = - \frac{d \nu_i}{dt}$  is bounded below whenever differentiable,
as follows from the almost-positive definiteness (see~\eqref{psibb} and its analogue for $1$-forms
after \eqref{MSR1F})  of $-\frac{A'(s)}{A(s)}$:
At any point at which the eigenvalues are distinct,
$-\nu_i'$ is computable by perturbation theory:  
if $v_i$ is an eigenvalue corresponding to $\lambda_i =\exp( i \nu_i t)$, then 
$$\lambda_i' \left( = \frac{d \lambda_i}{dt} \right) =  i \frac{ \langle A' v_i, v_i \rangle }{\langle v_i, v_i \rangle}.$$
Note that $\lambda_i'$ (and also $\nu_i'$)  denote the derivative with respect to the variable $t$,
where $A'$ denotes the derivative with respect to $s=it$. This gives $ i \nu_i' = \frac{\lambda_i'}{\lambda_i} =  i \frac {\langle A'/A  \ v_i, v_i \rangle} {\langle v_i, v_i \rangle}$.  So
$-\nu_i'  = \langle - \frac{A'}{A} v_i, v_i \rangle/\|v_i\|^2$, 
but the Maass-Selberg relations (see~\eqref{psibb} and page \pageref{derivbounded}) imply that this is bounded below.  Note for later use that
\begin{equation} \label{sumnu} \sum -\nu_i' = \mathrm{trace}   \frac{-A'}{A} =  - \frac{\left( \det A\right)'}{\det A} .\end{equation}

 We may also give {\em upper} bounds for $|-\nu_i'|$, although these are more subtle.
 In the arithmetic case, at least, there exists for every $j$ constants $A, B$ such that 
\begin{equation} \label{zerofree} \| (\Phi^{\pm})^{(j)} (s) \| \leq A | \log(2+|s|)|^{B}, \ \ s \in i \R. \end{equation}
See \cite{Sarnak} for a proof of the corresponding fact about $L$-functions. (It is  most unlikely that
we really need this estimate, but it is convenient in writing the proof with more explicit constants.)
   In particular, 
\begin{eqnarray} \label{nuiest} \frac{d}{dt} (- \nu_i + 4 t \log Y )&=& 4  \log Y + O(a (\log \log Y)^b)
\\  &\asymp& (\log Y) \nonumber \end{eqnarray}
(and thus a similar bound for $A$), 
 for suitable constants $a,b$, whenever
   \begin{equation} \label{Tmaxdef} |t| \leq \Tmax := (\log Y)^{100}. \end{equation} 
 We will often use $\Tmax$ as a convenient upper bound for the eigenvalues that
 we consider. Its precise choice is unimportant (although if we only had weaker
 bounds for the size of $\nu_i'$ we would use a smaller size for $\Tmax$). 
 As usual the notation $\asymp$ means
 that the ratio is bounded both above and below. 
 
There are corresponding estimates for the Eisenstein series itself. 
What we need can be deduced from the fact that   the Eisenstein series $E(s,f)$ and $E(s, \nu)$  have no poles
  around $s=it$ in a ball of radius $\geq \frac{C}{\log(1+|t|)}$
(see \cite{Sarnak} for a proof of the corresponding fact about $L$-functions)
and also the   bound there:
\begin{equation} \label{trivesteis}  \| E(\sigma+it,f)\|_{L^{\infty}(M_Y)} \leq a ((1+|t|) Y)^b  \ \ \ \left( |\sigma| \leq \frac{C}{\log(1+|t|)} \right),\end{equation}
for suitable absolute $a,b$.  Similar bounds apply to any $s$-derivative of $E(s,f)$. To obtain, for instance, an estimate
for the derivative of $E$, one uses Cauchy's integral formula and~\eqref{trivesteis}.

 We are interested in the solutions $s = it \ (t \in \R)$ to 
 $$ f(it) := \det(1 - Y^{-4it} A(it) ) = 0,$$
 that is to say,  those $s = it $ for which there exists $j$ with
\begin{equation} \label{nuje} -\nu_j(t) + 4 t \log(Y) \in 2 \pi \Z. \end{equation}
 Now  --- for sufficiently large $Y$ --- the function $-\nu_j(t) + 4 t \log(Y)$
 is monotone increasing, with derivative bounded away from zero. 
 This implies, speaking informally,  that the solutions to~\eqref{nuje} are very regularly spaced (see figure for an example which was actually computed with the scattering matrix for functions
 on $\SL_2(\Z[i])$.)

  \begin{figure}[!ht]
\begin{center}
  \includegraphics[width=70mm]{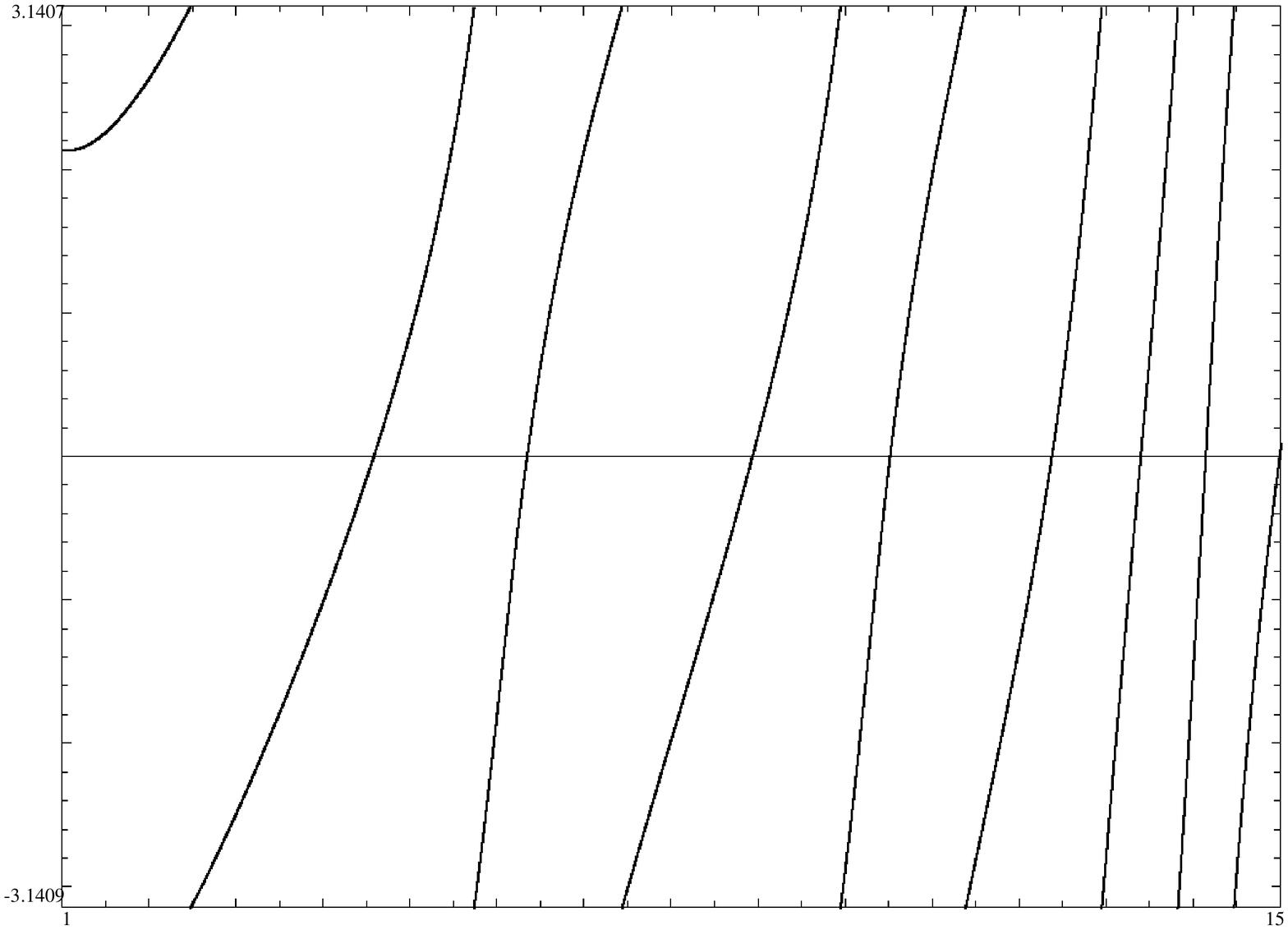}
\end{center}
\caption{What the function $-\nu_j(t) + 4 t \log(Y)$ looks like, modulo $2\pi$}
 \end{figure}

  \vspace{1cm}

\begin{lemma} \label{analysisofzeroes} There exists $Y_0$ such that, for $Y \geq Y_0$:    
  \begin{enumerate}
  
  \item[(i)]   The number of solutions to~\eqref{nuje} in the interval $0 < t \leq T $ is given by
$ \left[ \frac{ 4 T \log(Y) - \nu_j(T)}{2 \pi} \right]$.

  \item[(ii)] Every small value of $f(s)$ is near a root:
  
  If $|f(it) | < \epsilon$ and $|t| \leq \Tmax$,  there exists $t' \in \R$ with $f(it') = 0$ and
  $|t'-t| \ll \epsilon^{1/\hrel}$. 
  
    \item[(iii)] Nearby solutions are almost orthogonal:

   Suppose $t_0 \in  \R$ with $|t_0| \leq \Tmax$;  let $[ x,y]$ be the positive definite
   inner product on $\Omega^+(0)$ defined as $$[x, y]_{it_0} := \langle \left( 4 \log Y-\frac{A'(it_0)}{A(it_0)}  \right) x, y \rangle;$$
   we will sometimes omit the subscript $i t_0$ where clear. 
   
   (As previously remarked,   we interpret
   the fraction as $A(it_0)^{-1} A'(it_0)$. It is an endomorphism of $\Omega^+(0)$. This inner product arises naturally from the Maass-Selberg relations; see \eqref{Fnorm2}. In particular  one verifies that it is positive definite so long as $Y$ is large enough.)

Suppose $t_1 \neq t_2$ satisfy $|t_i-t_0|<\epsilon$, 
and we are given  $v_i$ satisfying the relations $Y^{-4 it_i} A(it_i) v_i =v_i$ for $i=1,2$. 
 Then $$ [v_i, w_i]_{i t_0}= O(\|v_i\| \|w_i\|  (\log Y)^2 \epsilon).$$

    \item[(iv)] Every almost-solution to $Y^{-4s} A(s) v = v$ is near (a linear combination of) exact solutions:
   
    Suppose that $|t| \leq \Tmax, \|v\| = 1$ are such that \begin{equation} \label{vO} \|Y^{-4it} A(it) v - v \| < \epsilon. \end{equation} 
Then there exists pairs $(t_i \in \R, w_i \in \Omega^+(0)): 1 \leq i \leq m$, with $m \leq \hrel = \mathrm{size}(A)$, 
and an absolute constant $M$, 
such that:
\begin{itemize}
\item[-] $Y^{-4 i t_i}  A(i t_i) w_i = w_i$; 
\item[-] $\|t-t_i\| \ll \epsilon^{1/M}$; 
\item[-] $\| v -\sum w_i\| \leq (\log Y)^{C} \epsilon^{1/M}$;
\item[-] $\|w_i\| \leq (\log Y)^C$. 
\end{itemize}

  \end{enumerate}
  \end{lemma}
  
  \begin{proof}  
We note throughout the proof that the assumption that we are considering
  $\nu_i(t)$  when $|t| \leq \Tmax$  implies  $|\nu_i'(t)|  = o(\log Y)$: 
indeed $|A'(s)| + |A''(s)|= o(\log Y)$ for such $s$ by choice of $\Tmax$. 
So for sufficiently large $Y_0$ we may certainly assume: 
\begin{equation} \label{Y0choice1}   (-\nu_i + 4t \log Y)' \geq \log Y \end{equation}
under the same conditions. 
  \medskip 
  
  (i) follows from monotonicity and our assumption that $-\nu_i(0) =0$. 
  \medskip
  
(ii): if $|f(it)| \leq \epsilon$, then there exists an eigenvalue $\lambda$
 of $A(it)$ satisfying $|1 - Y^{-4it} \lambda| \leq \epsilon^{1/h}$. Thus $\lambda = \exp(2 \pi i \nu_i(t))$ for some $i$. The result follows
 from the lower bound ~\eqref{Y0choice1}.
 
 \medskip 
 
(iii): Set $A_s := Y^{-4s} A(s)$ and consider
\begin{eqnarray*} 
0 &=& \\
\langle A_{it_1} v, A_{it_2} w \rangle - \langle v, w \rangle  &=& \langle (A_{it_2}^{-1} A_{it_1}  -1 ) v, w \rangle  \\ &=& i(t_2-t_1) [v,w]_{it_0} + O(|t_1-t_2|^2  (\log Y)^2 \|v\| \|w\|).\end{eqnarray*}
 where the factor of $\log Y$ results from estimating
 the second derivative of $t \mapsto A_{it}$, cf. ~\eqref{zerofree}.  In any case
 the claim follows immediately by rearranging this equation. 
 
 \medskip 
(iv): This proof is a little bit ugly.

 Notice, first of all,
that we may (by adjusting the constants $C, M$)
always suppose that $\epsilon$ is smaller than any fixed power of $(\log Y)^{-1}$. 
Otherwise the statement simply asserts that there exists $t_i$ close
to $t$ such that $Y^{-4 it_i} A(i t_i)$ has a fixed vector --- the assertions
about the position of that vector becomes vacuous ---  and this follows easily
from what we showed in (ii).

Choose $\delta > \sqrt{\epsilon}$ such that $Y^{-4it} A(it)$
has no eigenvalues satisfying $|\lambda -1| \in [\delta, \sqrt{\delta})$. 
Clearly we may do this with $\delta < \epsilon^{2^{-h-1}}$
by considering the intervals $[\epsilon^{1/2}, \epsilon^{1/4}],
[\epsilon^{1/4}, \epsilon^{1/8}], \dots$ in turn, and there are at most $\hrel$ eigenvalues.

Let $W$ be the span of all eigenvectors of $Y^{-4t} A(it)$
whose eigenvalue satisfy $|\lambda - 1| \leq \delta$;
 suppose that these eigenvalues correspond to   $\nu_{i_1}(t), \dots, \nu_{i_m}(t) \in \R/\Z$. Then the orthogonal projection of $v$ to $W$ --- call it $v_W$ --
satisfies
$\|v -v_W\| < \epsilon \delta^{-1} \ll \sqrt{\epsilon}$.  Here $v$ is as in the statement~\eqref{vO}. 

By the monotonicity and~\eqref{Y0choice1}, there exist solutions $t_1, \dots, t_h$ to 
$$ -\nu_j(x) + 4 x \log(Y) \in 2 \pi \Z $$
 satisfying $|t_i - t| \ll \frac{\delta}{\log Y}$.  
Let $v_{i_1}, \dots, v_{i_h}$ be corresponding fixed vectors for
$Y^{-4t_i} A(i t_i)$, normalized so that $\|v_?\| = 1$.  

 We have seen in (iii) that $[v_{i_a}, v_{i_b}] \ll \delta$ for $a \neq b$.  This implies in particular that the $v_?$ are linearly independent. Even more precisely, there is an absolute constant $C$ such that 
\begin{equation} \label{Graminequality} 
\det [v_i, v_j ] \gg (\log Y)^{-C}, \ \  \det \langle v_{i} , v_j \rangle \gg (\log Y)^{-C}, \end{equation}
where the second inner product is the standard one introduced on $C^{\infty}(0)$ 
in~\eqref{cinftyomegaomegaip}.
For the first inequality:  the diagonal entries of the matrix are in size $\gg (\log Y)$,
its off-diagonal entries are $\ll \delta$, and we may suppose, as we stated
at the start of the proof, that $\delta$  smaller
than any fixed power of $(\log Y)^{-1}$. For the second inequality: The ratio
$[v,v]/\langle v, v \rangle$ is bounded above and below by constant multiples of $(\log Y)$.
Therefore
the passage between the
  $\langle -, - \rangle$ metric and $[ - , - ]$ metric distorts  volumes
by a factor that is bounded above and below by constant multiplies of $(\log Y)^h$.

Let us now compute the projection of each $v_{i_a}$ onto the orthogonal complement of $W$. 
Let $u = v_{i_a}$ for some $a$, and write $$u = u_W + u_W' \ \ \ \ ( u_W \in W, u_W' \perp W).$$ 
Since $|t-s_i| \ll \delta/\log Y$, it follows from the remarks after~\eqref{Y0choice1} concerning the choice of $\Tmax$ that   $\|Y^{-4it} A(it) u - u \| \ll \delta$. But, all  the eigenvalues
of $Y^{-4it} A(it) - 1$ on $W^{\perp}$ are $\geq \sqrt{\delta}$, 
and thus it follows that $\|u_W'\| \ll \sqrt{\delta}$. 

Let $v_{i_a, W}$ be the projection of $v_{i_a}$ to $W$ with respect to the
standard metric $\langle -, - \rangle$. 
Then 
\begin{equation} \label{JAns}\det \langle v_{i_a, W} , v_{i_b, W} \rangle \gg (\log Y)^{-C} + C \sqrt{\delta}, \end{equation} 
That follows from the second inequality of~\eqref{Graminequality}
together with the fact we have just proved: $\|v_{i_a, W} - v_{i_a}\| \ll \sqrt{\delta}$. 
Indeed, all the entries of the matrix $\langle v_{i_a}, v_{i_b} \rangle$ 
are at most $1$ in absolute value, and we then are modifying each
entry by at most a constant multiple of $\sqrt{\delta}$, whence the
determinant changes by at most a constant multiple  of $\sqrt{\delta}$. 

We may suppose that 
$\delta$ is smaller than any fixed power of $(\log Y)$, as before;
so we can suppose that the right-hand side of~\eqref{JAns} is $\gg (\log Y)^{-C}$. 
Thus $v_{?, W}$ span $W$. Let $v_{?, W}^* \in W$ denote the dual basis
with respect to the $\langle -, - \rangle$ inner product. 
It also follows from~\eqref{JAns} that $\| v_{i_a}^*\| \ll (\log Y)^{C'}$. 

We may now write:
\begin{eqnarray*} v  & = & (v-v_W) + v_W  \\
& =&  (v-v_W) + \sum_{a} [v_W, v_{i_a,W}^*] v_{i_a,W} 
\\ &=& O(\sqrt{\epsilon}) + \sum_{a} \langle v_W, v_{i_a,W}^*\rangle v_{i_a} -  \sum_{a} \langle v_W, v_{i_a, W}^* \rangle (v_{i_a} -v_{i_a, W}) 
\\ &=& \sum_{a} \langle v_W, v_{i_a,W}^*\rangle v_{i_a}  +  O( (\log Y)^{C'} \sqrt{\delta}).
\end{eqnarray*}
 
 We have shown that there exists 
  absolute constants $C,M$ and
 constants $|c_a| \leq (\log Y)^{C}$ such that:
$$\|v_W - \sum c_a v_{i_a} \| \ll  (\log Y)^{C} \epsilon^{1/M}.$$
 \end{proof}

\subsection{Beginning of the proof} \label{sec-muller-thm}
 
 \medskip 
 {\textbf {\em Important notational warning.}}   \index{Exponentially decaying (Chapter 6)}
In this section, when we write {\em exponentially decaying} we shall   mean:
bounded by $a \exp(-bY)$ where $a,b$ are  constants depending only on $M$. 
For example, when we say ``$F$ is exponentially close to $G$'' for two functions $F, G \in L^2(M)$,
we mean, unless otherwise stated, that they are close in $L^2$-norm:
$$ \| F - G \|_{L^2(M)} \leq a \exp(-b Y),$$
for absolute constants $a,b$. 
The precise values of $a,b$, however, {\em may vary from one instance of this phrase to the next.}
This convention, we hope, improves the readability of the text.   Similarly, it will often be necessary for our argument that a quantity that is exponentially decaying be less than (say) $Y^{-1}$. Such an inequality is automatically valid for large enough $Y$. We will often proceed
implicitly assuming that $Y$ is large enough. This is certainly no loss, since
the statement to be proved allows us to suppose that $Y$ is larger than a fixed constant that we may choose.  We apologize for the implicit imprecision, but again we hope this improves the readability of the text. 
\medskip

We shall carry out the proof in the case of co-closed $1$-forms, proving assertion (b) of the theorem;  the other cases (c)---(e) are handled in the same way and are discussed briefly in~\S~\ref{sec:modifications}.

Denote
by $N(x)$ the number of eigenvalues of the Laplacian on co-closed $1$-forms for $M_Y$ with eigenvalue in $[0 , x^2]$, and $\bar{N}(x) = M_1(x) + M_2(x)$, where \begin{itemize}
\item[$M_1(x)$:] the number of zeroes of $t \mapsto f(it)$ in $[0,x]$;
\item[$M_2(x)$:] the number of eigenvalues in $[0,x^2]$ of the Laplacian on co-closed cuspidal $1$-forms
on $M$.
\end{itemize} 
Then we shall prove that, for certain absolute constants $a, b$  
we have
\begin{equation} \label{thm:coclosed1} \bar{N}(T- a e^{-bY}) \leq N(T) \leq \bar{N}(T+a e^{-bY}), \ \ 
x \leq \Tmax. \end{equation} 
where we interpret $T-a e^{-bY}$ as $0$ if $T< ae^{-bY}$. 
 This proves Theorem~\ref{prop:SE} (b).    
 
 The first inequality will come in  \S \ref{quasimodeargument} (`there is an eigenfunction of $M_Y$ near any root of $f$) and the second equality will be proved  in \S \ref{U2}
(``any eigenvalue of $M_Y$ is near a root of $f(s)$.'') 
And the intervening section  \S \ref{combin} explains what happens when treating eigenvalues very close to zero.

Recall (\S~\ref{sss:EisForms})
that to each $\omega \in \Omega^+(0)$, we have associated
an   
 ``Eisenstein series''  $E(s, \omega)$,
which is a $1$-form on $M$ with eigenvalue $-s^2$ under the form Laplacian
(and eigenvalue $-is$ under the operator $* d$). 
Moreover, 
there exists inverse linear operators
$  \Phi_{\pm}(s):   \Omega^{\pm}   \longrightarrow \Omega^{\mp } $
so that $\Phi_+(s) \Phi_-(s)= \mathrm{id}$, and 
so that the restriction of $E(s)$ to the cusps   looks like:
$$E(s) \sim y^{s}  \omega  +  y^{-s} \Phi^+(s) \omega.$$
 We write, as in the statement of the Theorem, 
$$f(s) = \det( \mathrm{Id} - Y^{-4s} \Phi^{-}(-s) \Phi^+(s)).$$ 
For short we write $\Phi(s) = \Phi^+(s)$. 

In Lemma~\ref{analysisofzeroes} we introduced an inner product on $\Omega^+(0)$, which
is 
\begin{equation} \label{FormSquareIP}   [\omega_1, \omega_2]_{it_0} =
  \langle \left( 4 \log Y- A^{-1}  A'(it_0)  \right) \omega_1, \omega_2 \rangle \end{equation} 
   and the inner product on the right hand side is the standard one introduced in 
~\eqref{cinftyomegaomegaip}. Again, this inner product arises
naturally when computing the norm of the truncated Eisenstein series, 
as in~\eqref{MSForms}.

 \begin{remark*} \em{Some remarks on small eigenvalues will be useful:
 
 At $s=0$,  $Y^{-4s} \Phi^{-}(-s) \Phi^+(s)$ is the identity transformation, and consequently $f(s) = 0$. Indeed, the derivative of $ \mathrm{Id} - Y^{-4s} \Phi^{-}(-s) \Phi^+(s)$ with respect to $s$
 is easily seen to be invertible for large enough $s$, and thus $f(s)$ has a zero
 of order equal to the size of $\Phi$, which also equals the number of zero eigenvalues
 for the Laplacian on co-closed $1$-forms. 
 
It will be useful  to note that smallest nonzero eigenvalue $\lambda_1$ of the Laplacian on zero-forms on $M_Y$, satisfying absolute boundary conditions (i.e. the smallest Neumann eigenvalue) is bounded below 
as $Y$ varies, as we see by an elementary argument; all we need is the weak bound
\begin{equation} \label{smalleigenvaluetrivial} \lambda_1 \geq  \mathrm{const} \cdot (\log Y)^{-1}. \end{equation}
 }  
  \end{remark*}

 \subsection{There is an eigenfunction of \texorpdfstring{$M_Y$}{M_Y}
  near any root of \texorpdfstring{$f$}{f}.} \label{quasimodeargument}
   
We show the first inequality of~\eqref{thm:coclosed1}, based on the idea sketched
in~\S~\ref{s:plaus} of constructing ``approximate eigenfunctions.''

Let $T \leq \Tmax$.
For each $s \in i\R_{\geq 0}$ with $f(s) = 0, |s| \leq T$, we choose
an orthonormal basis $\mathcal{B}_s$ for the set of solutions to \begin{equation} \label{tickytacky} \omega   = Y^{-2s} \Phi^{-}(-s) \bar{\omega} \mbox{ and } Y^s \bar{\omega} = Y^{-s} \Phi^+(s) \omega
\end{equation}
where orthonormality is taken with respect to the inner product $[-, -]_s$ (see~\eqref{FormSquareIP} for definition). 

Let $\mathcal{B}_{\cusp}$ be a basis for the set of coclosed cuspidal $1$-forms on $M$
with eigenvalue $\leq T^2$, and set
$$ \mathcal{B}_{eis}  = \{(s, \omega): f(s) = 0, \omega \in \mathcal{B}_s\}; $$
$$ \mathcal{B} = \mathcal{B}_{eis} \bigcup \mathcal{B}_{\cusp}$$

We shall associate to each element of $\mathcal{B}$ an element of $L^2(M_Y)$
that lies in the subspace spanned by co-closed $1$-forms of eigenvalue $\leq T+a \exp(-bY)$, 
for some absolute positive constants $a,b$, and then show that the resulting set $\mathcal{F}$ is linearly independent.

Let $(s,\omega) \in \mathcal{B}_{eis}$, and let $\bar{\omega}$
be the other component of the solution to~\eqref{omegaomegabar}.
   We normalize so that $\|\omega\| = 1$; the same is true for $\bar{\omega}$
   because $\Phi^{\pm}$ is unitary.   Set $F = E(s, \omega) + E(-s, \overline{\omega})$. We saw in~\S~\ref{s:plaus}
 that $F$ ``almost'' satisfies absolute boundary conditions, in the sense that,
 for example, the contraction of $dF$ with a boundary normal is exponentially small in $Y$. 
 
 \begin{remark*} {\em 
 When $s=0$, then the set of solutions to~\eqref{tickytacky} consists of all of $\Omega^+(0)$; moreover, 
  $\overline{\omega} = \Phi^{+}(0) \omega$
 and $F = 2 E(0, \omega)$. 
 }\end{remark*}

The next step is to modify $F$ slightly so as to make it satisfy absolute boundary conditions, but losing,
in the process, the property of being an exact $\Delta$-eigenfunction.  This can be done in many ways;
for definiteness one can proceed as follows:  fix a smooth monotonically increasing function $h: \R \rightarrow \C$ such that $h(x) = 1$ for $x > 1$
and $h(x) = 0$ for $x < 0$, and let $h_Y :y \mapsto h(Y/2-y)$. Now, near each cusp , considered as isometric to a quotient of $\H^3: y \leq T$,   write $$F =   y^s \omega + \Phi^{+}(s) y^{-s} \omega +  y^s \baromega + y^{-s} \Phi^{-}(s) \baromega + F_0,$$
and replace $F_0$ by $F_0 \cdot h_Y$; thus we replace $F$ by the form $F'$
which is given, in each cusp, by $$F' =   y^s \omega + \Phi^{+}(s) y^{-s} \omega +  y^s \baromega + y^{-s} \Phi^{-}(s) \baromega + F_0 \cdot h_Y,$$
Then $F'$
satisfies absolute boundary conditions. 
On the other hand, it is clear that \begin{equation} \label{F'almost} \|(\Delta + s^2) F'\| \ll \exp(-b Y)\end{equation} for some absolute constant $b$. 
 (Note that $F'$ is no longer co-closed, but it is almost so, as we shall check in a moment.) 

 We now see there is an co-closed eigenfunction close to $F'$: 
 Set $\epsilon = \exp(- b Y)$, where $b$ is the constant of the estimate~\eqref{F'almost}, and
let $W$ be the space spanned by co-closed $1$-forms on $M_Y$ with eigenvalue
 satisfying $\| \lambda + s^2\| \leq \sqrt{\epsilon}$ and let 
 $F = F_{s, \omega}$ be the projection of $F'$ onto $W$.  
 We refer to $F_{s,\omega}$ as {\em the eigenfunction associated to $(s, \omega)$};
 this is a slight abuse of terminology as $F_{s,\omega}$ need not in general
 be a Laplacian eigenfunction, but it behaves enough like one for our purposes. 
 
 Then:
 \begin{itemize}

 \item[(i)] $F_{s,\omega}$ is exponentially close to $E(s, \omega) + E(-s, \bar{\omega})$, where $\bar{\omega}$ is defined by 
~\eqref{omegaomegabar}; in particular,
  $\|F_{s,\omega}\|_{L^2(M_Y)} \asymp (\log Y)$. 
 \item[(ii)]  $F_{s,\omega}$ is a linear combination of eigenfunctions with eigenvalue
 exponentially close to $-s^2$. 
 \end{itemize}
 
 In order to verify these, we need to check that $F'$ is ``almost co-closed'', so that the projection
 onto co-closed forms does not lose too much.
 
 Indeed, first of all, the orthogonal projection of $F'$ onto the orthogonal complement of
 co-closed forms has exponentially small norm. This orthogonal complement
 is spanned by $df$, for $f$ a function satisfying absolute (Neumann) boundary conditions and orthogonal to $1$;  
now  \begin{eqnarray*} \langle F', df \rangle_{M_Y} = \langle d^* F', f \rangle_{M_Y}  
= \langle d^* (F'-F), f \rangle_{M_Y} 
, \end{eqnarray*}
where we used the fact that $d^* F= 0$.   But $F'-F$ is supported entirely ``high in cusps;'' 
we may express the last quantity as a sum, over cusps, of terms 
bounded by $ \| d^* F_0(1-h_Y) \| \|f\|$.
 Now  $\|d f\| \gg \|f\| /\log Y$ because of
~\eqref{smalleigenvaluetrivial}. So, 
  $$\langle F', df \rangle \ll (\log Y) \| d^* F_0(1-h_Y) \| \cdot \|df\|,$$
Since $d^* F_0(1-h_Y)$ is certainly exponentially small, we arrive
at the conclusion that the projection of $F'$ onto   the orthogonal complement of
 co-closed forms has exponentially small norm. 
 Similarly, the orthogonal projection of $F'$ onto the span of eigenvalues
 {\em not in} $[-s^2 - \sqrt{\epsilon}, -s^2 + \sqrt{\epsilon}]$ has, 
 by~\eqref{F'almost}, a norm $\ll \sqrt{\epsilon}$.
 
 This concludes the proof of claim (i) about $F_{s,\omega}$; claim (ii) follows directly from the definition.

 In an exactly similar way, we may associate to every co-closed cusp $1$-form $\psi$
 a co-closed   $F_{\psi}$ on $M_Y$ satisfying corresponding properties:
 
   \begin{itemize}
 \item[(i)] $F_{\psi}$ is exponentially close to $\psi |_{M_Y}$; 
 \item[(ii)] $F_{\psi}$ is a linear combination of eigenfunctions on $M_Y$ with eigenvalue exponentially close to the eigenvalue of $\psi$. 
 \end{itemize}

The maps $(s, \omega) \mapsto F_{s, \omega}, \psi \mapsto F_{\psi}$
give a mapping from $\mathcal{B} = \mathcal{B}_{eis} \cup \mathcal{B}_{\cusp}$ 
to the vector space spanned by all co-closed  forms on $M_Y$
of eigenvalue $\leq T^2+ a \exp(-b Y)$, for some absolute $a,b$. Call the image of this mapping $\mathcal{F}$.   We argue
now that $\mathcal{F}$ is linearly independent:

Indeed,  in any relation
\begin{equation} \label{linearrelation} \sum a_j F_{s_j, \omega_j} + \sum_k  b_k F_{\psi_k}  = 0, \ \ \ (0 \neq a_j \in \C, 0 \neq b_k \in \C).\end{equation}  we may suppose
(by applying an orthogonal projection onto a suitable band of Laplacian eigenvalues) that all the $s_j$ are ``close to one another,'' i.e.,
 there exists $T_0 \leq T+1$ so that $|s_j - i T_0|$ and $|s_{\psi_k} - i T_0|$ 
  are all exponentially small.

By~\S~\ref{subsec:fganal}, the number of solutions
to $f(s_j) = 0$ in the range $|T_0^2 -s_j^2 | < a \exp(-b Y)$
is absolutely bounded for $Y$ large enough. 
As for the number of $\psi_k$ that may appear in~\eqref{linearrelation},
this is (by the trivial upper bound in Weyl's law) polynomial in $T_0$. 
Therefore the number of terms in the putative linear relation is at worst polynomial in $\Tmax$. 

We shall derive a contradiction by showing that the ``Gram matrix'' of inner products
for $F_{s_j, \omega_j}$ and $F_{\psi_j}$ is nondegenerate (here $j,k$
vary over those indices appearing in~\eqref{linearrelation}). 
Clearly $\langle F_{s_j, \omega_j}, F_{\psi_k} \rangle $ is exponentially small
and  $\langle F_{\psi_j}, F_{\psi_k} \rangle$ differs from $\delta_{jk}$ by an exponentially small factor. These statements follow from the property (i) enunciated previously 
for both $F_{s, \omega}$ and $F_{\psi}$.  As for the terms, 
 \begin{equation} \label{innerproductmatrix} \langle F_{s_i, \omega_i} , F_{s_j, \omega_j} \rangle  
 = [\omega_i, \omega_j] + \mbox{exponentially decaying} \end{equation}  
by the Maass-Selberg relations; this is proved as \eqref{Fnorm2}.  

But by the analysis of zeroes, assertion (iii) of Lemma~\ref{analysisofzeroes},
 $[\omega_i, \omega_j]$ is itself exponentially small if $i \neq j$. 
 Consequently, the Gram matrix is nonsingular by ``diagonal dominance;''
 contradiction. 
 
 We have now exhibited a set of linearly independent co-closed $1$-forms 
 with eigenvalue $\leq T^2 + \mbox{ exponentially small}.$ This establishes the first  inequality in~\eqref{thm:coclosed1} (the precise value of the constants $a,b$ in that equation
 are obtained by tracing through the estimates mentioned here).

 \subsection{ Analysis of eigenvalues on \texorpdfstring{$M_Y$}{M_Y}
 near \texorpdfstring{$0$}{0} by passage to combinatorial forms.} \label{combin}
 
 We continue   our proof   Theorem~\ref{prop:SE} (a).  

We will argue that any ``very small'' eigenvalue of the Laplacian on co-closed $1$-forms on $M_Y$ is in fact exactly zero. We prove in fact a stronger statement:
 \begin{enumerate}
 \item[A.]  Any eigenvalue  $\lambda$ on co-closed forms on $M_Y$ satisfying $\lambda \leq A Y^{-B}$
 is actually zero, for suitable constants;
 \item[B.] If $\omega_1, \dots, \omega_r$ form an orthonormal basis for $\Omega^+(0)$,
 then the functions $F_{0, \omega_i}$ constructed in the previous section\footnote{A priori, these functions were not guaranteed to be Laplacian eigenforms; however, the preceding statement (A) guarantees that they are in this case.} are a basis
 for harmonic forms on $M_Y$. 
 \end{enumerate}
  
  For this discussion, it does not, in fact,  matter whether we discuss only co-closed $1$-forms
 or all $1$-forms, because of~\eqref{smalleigenvaluetrivial}: any $1$-form with
 eigenvalue polynomially close to zero must be co-closed. 
  As always, however, we work with respect to absolute boundary conditions.

We will compare the analysis of $j$-forms and of {\em combinatorial} $j$-forms, that is to say, the cochain complex with respect to a fixed triangulation.

Fix a  triangulation of $M_{\leq Y}$ where each triangle has hyperbolic sides of length $\leq 1$. This can be done with $O(\log Y)$ simplices, and for definiteness we do it in the following way:
first  fix $Y_0 \geq 1$ sufficiently large, then triangulate $M_{\leq Y_0}$ for some fixed $Y_0$ (this can be done --- any differentiable manifold may be triangulated). 
Now to triangulate the remaining ``cusp region'' $\mathcal{C}_{Y_0} - \mathcal{C}$,
split it into regions of the type $M \leq \height(x) \leq 2M$ plus an end region
$2^j M \leq \height(x) \leq Y$; we will assume that $2^j M \in [Y/2, Y/4]$. 
Each such region is diffeomorphic to the product $\R^2/L \times [0,1]$
for a suitable lattice $L$, simply via the map $(x_1, x_2, y) \mapsto (x_1, x_2) \times \frac{y-M}{M}$ (with the obvious modification for the end cylinder);   we fix a triangulation for the latter and pull it back. 

We equip the resulting cochain complex with the ``combinatorial'' inner product: The characteristic function of distinct simplices form orthonormal bases. This being done, we define a ``combinatorial'' Laplacian on the cochain complex
just as for the de Rham complex, i.e. $\Delta := d_c d_c^* + d^* d$, where $d_c$
is the differential and $d_c^*$ its adjoint with respect to the fixed inner product.  
 
Suppose given an $N \times N$ symmetric integer matrix, all of whose eigenvalues  
$\mu$ satisfy  the bound $|\mu| \leq A$, where $A$ is real and larger than $1$.
Then also  every nonzero eigenvalue satisfies the lower bound $|\mu| >  A^{-N}$ if $\mu \neq 0$; 
this is so because all eigenvalues are algebraic integers, and in particular
they can be grouped into subsets with integral products.   On the other hand,
every eigenvalue of the combinatorial Laplacian is bounded by the maximal number of simplexes ``adjacent'' to any given one: if we write
out the combinatorial Laplacian with respect to the standard basis, the $L^1$-norm
of each row is thus bounded.   The same goes for $d_c^* d_c$ and $d_c d_c^*$. Therefore, 
\begin{eqnarray} \label{lowerbound} \ 
\mbox{Any nonzero eigenvalue $\lambda$ of the combinatorial Laplacian, }
\\ \nonumber \mbox{or of $d_c^* d_c$, or of $d_c d_c^*$, 
satisfies $|\lambda| \gg Y^{-m}$ }\end{eqnarray}
for some absolute constant $m$.  
This simple bound will play an important role.
\medskip

Any $j$-form $\omega$ on $M_{\leq Y}$ induces (by integration) a combinatorial form, denoted $\omega_c$, 
on the $j$-simplices of the triangulation. 
 If $d \omega  = 0$, then $\omega_c$ is closed, i.e. vanishes on the boundary of any
combinatorial $(j+1)$-cycle.

Let us now consider $V(\epsilon)$, the space spanned by $1$-forms on $M_{\leq Y}$
of eigenvalue $\leq \epsilon$.   For topological reasons, $$\dim V(\epsilon) \geq  \dim V(0) = \dim H^1(M_{\leq Y}),$$
no matter how small $\epsilon$.
 We shall show that $\dim V(\epsilon) = \dim H^1$ holds if $\epsilon <  a Y^{-b}$
for some absolute constants $a,b$. In particular, this shows that all eigenvalues 
less than $aY^{-b}$ are identically zero.

 Suppose to the contrary  that $\dim V(\epsilon) > \dim H^1(M_{\leq Y})$. 
Consider the map from $V(\epsilon)$ to $\mathrm{ker}(d_c)/\mathrm{im}(d_c)$
given by first applying $\omega \mapsto \omega_c$, then taking the orthogonal projection
onto $\mathrm{ker}(d_c)$, and finally projecting to the quotient.
Since we're assuming that $\dim V(\epsilon) > \dim H^1(M_{\leq Y})$ this map must have a kernel. Let $\omega$ be a nonzero element of this kernel.  The idea is to 
construct a combinatorial antiderivative for $\omega$, then to turn it into
an actual ``approximate'' antiderivative; the existence of this will
contradict the fact that $\omega$ is (almost) harmonic.

There is an absolute $m> 0$ such that,  for any $\omega \in V(\epsilon)$, we have:
\begin{equation} \label{CCC} \| d_c   \omega_c \|_{L^2_c} \ll \epsilon^{1/2} Y^{m} \|\omega\|, \end{equation} 
where $L^2_c$ denotes the combinatorial inner product
and $a$ is an absolute constant. 
Indeed, if $K$ is any $2$-simplex, we have
$\langle d_c \omega_c, K \rangle = \langle \omega_c, \partial K \rangle =
\int_{\partial K} \omega = \int_K d\omega$. 
But \begin{equation} \label{domegabound} \|d \omega\|_{L^{\infty}} \ll  \epsilon^{1/2}Y^{m_1} \|\omega\|, \end{equation}
as follows from (say) the Sobolev inequality; 
  the factor $Y^{m_1}$ arises from the fact that the manifold $M_{\leq Y}$
  has injectivity radius $\asymp Y^{-1}$.    That shows~\eqref{CCC}.

Now~\eqref{lowerbound} shows:
$$\mathrm{dist}(\omega_c, \ker(d_c)) \ll  \epsilon^{1/2}  Y^{m_2} \|\omega\|$$
for some absolute constant $m_2$; here the distance $\mathrm{dist}$ is taken with respect
to the $L^2$-structure on the space of combinatorial forms.  
In fact, write $\omega_c = \omega_{c, k} + \omega_{c}'$
where $\omega_{c,k} \in \ker(d_c)$ and $\omega_{c}' \perp \ker(d_c)$. 
Then $\| d_c \omega_c'\| \ll \epsilon^{1/2} Y^m \|\omega\|$; 
in particular $\langle d_c^* d_c \omega_c', \omega_c' \rangle \ll (\epsilon^{1/2} Y^m \|\omega\|^2)$,
and then we apply~\eqref{lowerbound}. 

By our choice of $\omega$, the orthogonal
 projection $\omega_{c,k}$ of $\omega_c$ to $\ker(d_c)$ belongs to $\mathrm{image}(d_c)$. 
So, there exists a combinatorial $0$-form $f_c$ such that:
\begin{equation} \label{bottlenose} \| \omega_c - d_c f_c \| \ll \epsilon^{1/2} Y^{m_2} \|\omega\|.\end{equation}
We now try to promote the combinatorial $0$-form $f_c$ to a function, which
will be an approximate antiderivative for $\omega$. 

Let us now fix a base vertex of the triangulation, $Q$.  For any other three-dimensional simplex $\sigma$
fix a path of minimal combinatorial length $\gamma_{\sigma}$ from $Q$ to a vertex    $Q_{\sigma}$  of the simplex.  For every point in $\sigma$,
fix a linear path\footnote{The notion of linear depends on an identification of the simplex with a standard simplex $\{ x_i \in \mathbf{R}: x_i \geq 0 , \sum x_i = 1\}$. For our purposes, the only requirements on these identifications is that their derivatives with respect to standard metrics should be polynomially bounded in $Y$; it is easy to see that this is achievable.}  $\gamma_P$ from $Q_{\sigma}$ to $P$.  Now define $f: M_{\leq Y} \rightarrow \C$ via 
$$f(P) = \int_{\gamma_{\sigma} + \gamma_P} \omega, $$
where $\gamma_\sigma + \gamma_P$ denotes concatenation of paths.
 
  This defines
$f$ off a set of measure zero; 
$f$ has discontinuities, i.e., doesn't extend to a continuous function on $M$, but they are small because of~\eqref{domegabound} and~\eqref{bottlenose}.
Precisely, for any point $P$ on $M$ and a sufficiently small ball $B$ about $P$, we have
$$\sup_{x,y \in B} |f(x) - f(y)| \ll \ \epsilon^{1/2} Y^{m_3} \|\omega\|$$   for suitable $m_3$. 

On the interior of every simplex $S$,  $df$ is very close to $\omega$: precisely, 
\begin{equation}\label{dfnearomega} \| df - \omega \|_{L^{\infty}(S^{\circ})}  \ll \epsilon^{1/2} Y^{m_4} \|\omega\|. \end{equation}
Again this follows from Stokes' formula and~\eqref{domegabound}. 

Finally, $f$ almost satisfies absolute boundary conditions:
if we denote by $\partial_n$ the normal derivative at the boundary, we have
\begin{equation} \label{partialsmall} 
|\partial_n f |_{L^{\infty}} \leq \epsilon^{1/2} Y^{m_5}  \|\omega\|.\end{equation}
This follows from~\eqref{domegabound} and the fact that $\omega$ itself satisfies absolute boundary conditions.   

Now it is routine to smooth $f$ to 
  get a {\em smooth} function $\tilde{f}$
that satisfies~\eqref{partialsmall} and~\eqref{dfnearomega}, perhaps
with slightly worse $m_3, m_4$ and implicit constants. \footnote{
We explain how to do this in some detail in a coordinate chart near the boundary.
One then splits $f$ into a part supported near and away from the boundary; the part
away from the boundary can be handled by convolving $f$ with 
a $\mathrm{PU)_2}$-bi-invariant smooth kernel on $\PGL_2(\C)$.
Let $(x,y,z) \in  S =\R/\Z \times \R/\Z \times [0,1]$. 
Given a function $f: S \rightarrow \C$ that satisfies
$f_z(x,y,1) =0$, we extend it to a function $f': \R/\Z \times \R/\Z \times [0,2] \rightarrow \C$ by forcing the symmetry
$f(x,y,2-z) = f(x,y,z)$.  Now let $\omega \in \C^{\infty}_c(\R^3)$ 
be such that $\int \omega =1$, 
and let $\omega_{\delta}(\mathbf{x}) = \delta^{-3} \omega(\mathbf{x}/\delta)$. 
Now set $\tilde{f} = f' \star \omega$ (convolution on the abelian group $(\R/\Z)^2 \times \R)$. 
Note in particular that $d\tilde{f} = df \star \omega$. }

Now
\begin{eqnarray} \langle \omega, \omega \rangle &= &\langle d \tilde{f} , \omega \rangle   + \langle d\tilde{f}-\omega, \omega \rangle \\  &=& \langle \tilde{f} , d^* \omega \rangle +  \mbox{boundary term} + \langle d \tilde{f} - \omega, \omega \rangle   .\end{eqnarray}
But $\langle \tilde{f}, d^* \omega \rangle$ vanishes because $\omega$ is co-closed;
the boundary term is bounded by $\epsilon^{1/2} Y^{m_6} \|\omega\|^2$ because $\tilde{f}$ almost satisfies Neumann boundary conditions,
i.e.~\eqref{partialsmall}.   So
$$ | \langle \omega , \omega \rangle| \ll \epsilon^{1/2} Y^{m_7} \|\omega\|^2 $$
for suitable $B$. 
This is a contradiction whenever $\epsilon \leq A Y^{-B}$ for suitable $A,B$, and completes the proof that
$$ \dim V(\epsilon) = \dim V(0) \ \ \ (\epsilon \ll Y^{-B}). $$

\begin{remarkable} \em{ For later usage we compute the asymptotic behavior of the regulator for $M_{\leq Y}$; we'll show
\begin{equation} \label{rtymt} \reg(H_1(M_{\leq Y})) \sim \reg(H_1(M)) (\log Y)^{-\hrel} \end{equation} where $\hrel$ is the number of relevant cusps of $Y$, and $\sim$ denotes that the ratio approaches $1$
as $Y \rightarrow \infty$. 

Let $\psi_1, \dots, \psi_k$
be an orthogonal basis for coclosed harmonic cusp forms on $M$.  (Thus, $k= \dim H^1_c(M, \C)$). 
Let $\nu_i$ be the eigenfunctions associated to $\psi_i$  under the map $\mathcal{B} \rightarrow \mathcal{F}$ of~\S~\ref{quasimodeargument}. (Again, we emphasize that
we only know that these are eigenfunctions because of what we proved in this section \S \ref{combin}.)

Let $\omega_{k+1}, \dots, \omega_{k+\hrel}$ be a basis for $\Omega =\Omega^+(0)$;
note that 
$$1 - Y^{-4s} \Phi^{-}(-s) \Phi^{+}(s) \equiv 0$$
 for $s=0$, because $\Phi^-(s)$ and $\Phi^+(s)$ are inverses; 
thus we may form the eigenfunction $\nu_i := F_{0,\omega_{i}}$ on $M_Y$ associated to $(0, \omega_i)$
for each $k+1 \leq i \leq k+\hrel$.  We set $\psi_i := E(0, \omega_i)$ for $k+1 \leq i \leq k+\hrel$. 
Then $\nu_i$ is exponentially close to $\psi_i$.  (As usual, this means: in $L^2$ norm on $M_Y$).

Again $\nu_1, \dots, \nu_{k+\hrel}$ are linearly independent 
and are harmonic.  Since  $\dim H^1(M_Y, \C) = \dim H^1(M, \C) = k+\hrel$, 
therefore, 
\begin{quote} The $\nu_i$, for $1 \leq i \leq k+\hrel$, form a basis for harmonic forms.
\end{quote}

Let $\gamma_1, \dots, \gamma_{k+\hrel}$ be a basis for $H_1(M_{\leq Y}, \Z)$ modulo torsion. 
Then, to check~\eqref{rtymt}, we note first of all that
$\int_{\gamma_i} \nu_j \stackrel{Y \rightarrow \infty}{\rightarrow} \int_{\gamma_i} \psi_j$ for each $i, j$, i.e. the period matrix for $M_Y$ approaches that for $M$. This requires a little more than
simply the fact that $\nu_j$ is exponentially close to $\psi_j$ which is, a priori, a statement
only in $L^2$-norm. We omit the easy proof. 
Now~\eqref{rtymt} follows since
$$ \det  \left( \langle \nu_i, \nu_j \rangle \right)_{1 \leq i,j \leq k+\hrel}  \sim (\log Y)^{\hrel}  \det \left( \langle \psi_i, \psi_j \rangle \right)_{1 \leq i,j \leq k+\hrel},$$
as follows from the fact that the $\nu_i, \psi_i$ 
are exponentially close, and   the definition ~\eqref{constant-term-norm} of the inner product
on Eisenstein harmonic forms for $M$.
 } \end{remarkable}

 \subsection{Any eigenvalue of \texorpdfstring{$M_Y$}{M_Y} arises from a root  of
 \texorpdfstring{$f(s)$}{f(s)} or a cusp form.}  \label{U2} 
 Again $T \leq \Tmax$. 
 
 Let $\mathcal{F}$ be as constructed in~\S~\ref{quasimodeargument};
 the set of (near)-eigenforms associated to $\mathcal{B}_{cusp} \bigcup \mathcal{B}_{eis}$. 
 Let $\mathcal{F}_0$ be the set of eigenforms of eigenvalue $0$ belonging to $\mathcal{F}$, i.e., 
 the harmonic forms. 
 
 Let \begin{equation} \label{epsilon0def} \epsilon_0 = A Y^{-B}\end{equation}  be as in~\S~\ref{combin}, i.e., such that
 any coclosed eigen-$1$-form of eigenvalue $\lambda \leq \epsilon_0$
 is known to satisfy $\lambda = 0$, and therefore to belong to the span of $\mathcal{F}_0$. 
 As in the statement of the theorem, we write $\delta = a \exp(- b Y)$; the constants
$b$ will be chosen sufficiently small to make various steps in the proof work. 

 We shall show that --- for suitable $b$: \label{686mainassertion} 
 \begin{quote} (*) any co-closed eigen-$1$-form on $M_Y$ with eigenvalue in $(\epsilon_0, (T-\delta)^2]$
 that is orthogonal to all elements of $\mathcal{F} - \mathcal{F}_0$
 is identically zero.  \end{quote} 
  Since we already know that any co-closed $1$-form
 with eigenvalue in $[0, \epsilon_0]$ belongs to the span of $\mathcal{F}_0$,  this gives the second inequality of~\eqref{thm:coclosed1}, 
 and completes the proof of Theorem~\ref{prop:SE} part (b). 
 
The idea of the proof: the Green identity (recalled below)
     shows that any $1$-form $f$ on $M_Y$ that's a Laplacian eigenfunction 
cannot have constant term  which is ``purely'' of the form $y^{-s}  (a dx_1 + b dx_2)$;
its constant term must also\footnote{This corresponds to the  fact that
an ``incoming'' wave reflects and produces an ``outgoing'' wave. In other words there is a symplectic pairing on the possible asymptotics, with respect to which the 
possible realizable asymptotics form a Lagrangian subspace.}
 contain a piece $y^s (a' dx_1 + b' dx_2)$. 
But  by subtracting a suitable Eisenstein series (rather: its restriction to $M_Y$)
from $f$ we can arrange that the constant term of $f$ indeed
looks like $y^{-s} (a dx_1+ bdx_2)$. 
 This leads to a new function $M_Y$ whose constant term is very close to zero;
 we then seek to show it is very close to a restriction of a cusp form from $M$.
 This we do by spectrally expanding $f$ on $M$.  The trickiest point is to control
  the inner product of $f$ with an Eisenstein series
 $E(\omega, t)$ when their eigenvalues are very close; for that, we use 
 a pole-free region for Eisenstein series.

Let us suppose that $\eta = \eta_s$ is an co-closed eigenfunction of the $1$-form Laplacian on $M_{\leq Y}$  with eigenvalue $-s^2 \in (\epsilon_0, (T-\delta)^2]$ and with $\|\eta\|=1$, that is to say,  $\int_{M_Y} \langle \eta, \eta \rangle = 1$.  Note we are always assuming
$T \leq \Tmax$. 
Suppose that $\eta \perp \mathcal{F}$; we shall show $\eta = 0$.   
 
The Green identity gives
\begin{eqnarray}   \label{greenid} \int_{M_{\leq Y}} \langle \omega_1, \Delta \omega_2 \rangle - \langle \Delta \omega_1, \omega_2 \rangle & =&  B_1 + B_2, \\ \nonumber
& B_1 & = \int_{\partial M_{Y}}   \langle d^* \omega_1, \omega_2(X) \rangle -  \langle \omega_1(X), d^* (\omega_2) \rangle, \\ \nonumber & B_2  &  = \int_{\partial M_{Y}} \langle X . d\omega_1, \omega_2 \rangle - \langle  \omega_1,  X.d\omega_2 \rangle.\end{eqnarray} 
where $X$ is a unit normal to the boundary -- in our context it is simply $y \partial_y$ --
and $\omega(X)$ means the evaluation of $\omega$ on $X$, while
$X. d\omega$ is the result of contracting $d\omega$ with $X$.  
  
At first, we do not use the fact that $\eta$ satisfies absolute boundary conditions. 
Consider a fixed cusp of $M$, and write the constant term of $\eta$ (see~\eqref{constanttermdef} for definition) as \begin{eqnarray*} (\eta)_N &=& \omega_{s}y^s + 
\omega_{-s} y^{-s} +\bar{\omega}_{s} y^{s} + \bar{\omega}_{-s} y^{-s} \end{eqnarray*} where $\omega_s , \omega_{-s} \in \Omega^+ = \Omega^+(0), 
\bar{\omega}_s, \bar{\omega}_{-s} \in \Omega^-(0)$.
Note that assumption that $\eta$ is co-closed, i.e., that $d^* \eta \equiv 0$, 
implies that  $\eta_N$ does not contain terms of the form $y^t dy$.

Define $$\omega = \eta - \left( E(s,\omega_{s}) + E(-s, \bar{\omega}_s) )\right) |_{M_Y},$$
a co-closed form on $M_Y$ whose constant terms equals
$$\omega_N :=  (\bar{\omega}_{-s}  + \omega_{-s} - \Phi^+(s) \omega_s - \Phi^{-}(-s) \bar{\omega}_s)  y^{-s} .$$ 
In this way we have killed the $y^s$ part of the constant term of $\eta$. 
\medskip

We apply the Green identity~\eqref{greenid} with $\omega_1 = \omega_2 = \omega$
{\em and also} with $Y$ replaced by $Y' := Y/2$. 
The left-hand side equals zero. By~~\eqref{precisebound},   
 we deduce that
 $$|\omega(X)|_{L^{\infty}(\partial M_{Y/2})}  \ll \exp(-b Y),$$
  on the other hand,  $\| X . d\omega - s   \omega_N \|_{L^{\infty}(\partial M_{Y/2})} \ll \exp(-b Y)$
by the same~\eqref{precisebound}.
The Green identity implies that $(s- \bar{s}) \|\omega_N \|_{L^2(\partial M_{Y/2})}^2 \ll \exp(-bY)$;
to be explicit,  the norm here is simply $\| y^{-s} (a dx_1 + b dx_2)\|_{L^2(\partial M_{Y/2})} 
\propto (|a|^2 + |b|^2)$, the constant of proportionality being the area of the cusp.

This estimate  implies (by assumption, $-s^2 \leq \epsilon_0$, which gives a polynomial-in-$Y$ lower bound on $s-\bar{s}$): 
\begin{equation} \label{harr} \| \omega_N\|^2 \ll \exp(- b Y/2). \end{equation}
Making~\eqref{harr} explicit (by splitting into $dx_1+i dx_2$ and $dx_1 - i dx_2$ parts,
i.e. ``holomorphic'' and ``antiholomorphic'' parts): 
\begin{equation} \label{wk1} \| \omega_{-s} - \Phi^{-}(-s) \bar{\omega}_s\| , \ \|\bar{\omega}_{-s} - \Phi^{+}(s) \omega_s \| \ll \exp(-b Y/2). \end{equation} 

\medskip

In words, \eqref{wk1} says that the constant term of $\eta$ resembles the consant term
of $E(s, \omega_s) + E(-s, \bar{\omega}_s)$. We now
want to check that $s$ needs to be near a root of $f(s)$, which we will
do by seeing what the boundary conditions tell us; we then want to argue
that  $\eta$ itself differs from $E(s, \omega_s) + E(-s, \bar{\omega}_s)$
by (something very close to the) restriction of a cusp form from $M$, which we will do by spectral expansion. 

\medskip

We now  use the absolute boundary conditions on $\eta$: they mean  that $d \eta$ contracted with a boundary vector should be zero, 
which implies that $s  (Y^s \omega_s + Y^s \bar{\omega}_s - Y^{-s} \omega_{-s} - Y^{-s} \bar{\omega}_{-s}) = 0$. Note that we have {\em exact} equality here:   integrating $X . d\eta$
over the boundary picks up solely the constant term of $\eta$. (Thus, the boundary conditions gives one constraint for each Fourier coefficient, and this is the zeroth Fourier coefficient.)

 Thus,  considering ``holomorphic'' and ``antiholomorphic'' components, 
\begin{equation} \label{wk2}  \|Y^s \omega_s - Y^{-s} \omega_{-s}\|= \| Y^s \bar{\omega}_s - Y^{-s} \bar{\omega}_{-s} \|  =0.
\end{equation}

Taken together,~\eqref{wk1} and~\eqref{wk2} imply that
$$ \|Y^s \omega_{s} - Y^{-s} \Phi^{-}(-s) \bar{\omega}_s\|, \|Y^s \bar{\omega}_s - Y^{-s} \Phi^+(s)\omega_s\| \ll \exp(-b Y/2),$$
and since $\Phi^{\pm}$ are unitary, we deduce that  
\begin{equation} \label{wk3} \| \omega_s - Y^{-4s} \Phi^{-}(-s) \Phi^{+}(s) \omega_s\| \ll \exp(-b Y/2).\end{equation} 
This shows (more or less -- details below) that $s$ must be near zero of $f(s)$.

We now show estimates:
\begin{equation}
\label{omegasomegaest} \| \omega_s \| \gg  c_1 Y^{-c_2},     \|\omega\|_{L^2(M_Y)} \ll_N 
Y^{-N} \end{equation}
for suitable constants $c_1, c_2$ and for arbitrary $N$. 

The proofs of both parts of~\eqref{omegasomegaest} have the same flavor:
we need to show that for any co-closed $1$-form  $\xi$ (namely, either $\eta$ or $\omega$),
then $\xi$ is small if its constant term is.  More precisely,  we will prove
\begin{lemma} \label{SmallConstantTermMeansSmall} 
For any $L > 1$, there are constants $a,b,a'$ such that,
for $\xi \in \{\eta, \omega\}$, 
 $$\|\xi\|_{L^2(M_Y)} \ll a Y^b \|\xi_N\|   + a' Y^{-L} \|\eta\|_{L^2(M_Y)},$$
where the norm  on $\xi_N$is defined in \eqref{Nnormdef}.    \end{lemma}

We prove this in~\S~\ref{SCTMSproof}.  Of course $\|\eta\| $ is supposed to be $1$, but
we write it as $\|\eta\|$ to preserve homogeneity of the equation. 

The second assertion of~\eqref{omegasomegaest}
is directly the Lemma applied to $\xi = \omega$, together with~\eqref{harr}; 
  the first assertion of~\eqref{omegasomegaest} by applying this Lemma
to $\xi = \eta$ together with~\eqref{wk1} and~\eqref{wk2}.

But~\eqref{omegasomegaest} concludes the proof of (*) (on page \pageref{686mainassertion}, that is, the main assertion of the present section). In detail: The first estimate implies via Lemma~\ref{analysisofzeroes} part (iv)
 that $\omega_s$ is close to a linear combination of solutions   $\sum \omega'_i$ where $(s_i', \omega_i')$
is an exact solution of $Y^{-4 s_i'} \Phi^-(-s_i') \Phi^+(s_i') \omega_i' = \omega_i'$
and all the $s_i'$ are very close to $s$ (in particular, $s_i' \neq 0$).   The cited lemma gives certain bounds on the norms of the $\omega_i'$, which we use without comment in what follows. The second estimate of ~\eqref{omegasomegaest} implies that $\| \eta - E(s, \omega_s) - E(-s, \bar{\omega}_s)\|_{L^2(M_Y)}$
is $\ll Y^{-N}$, and so also\footnote{Indeed,  put $\bar{\omega}_i '= Y^{-2s} \Phi^{+}(s_i') \omega_i'$, 
so that $\bar{\omega}_s$ is exponentially close to $\sum \bar{\omega_i}'$
then also 
$E(s, \omega_s) + E(-s, \bar{\omega}_s) $ is exponentially close to
$ \sum E(s_i', \omega_i') + 
\sum E(-s_i', \bar{\omega}_i') $. To check the latter assertion
one needs estimates on the derivative of the Eisenstein series in the $s$-variable; see discussion around~\eqref{trivesteis} for this.    In turn, each $E(s_i', \omega_i') + E(-s_i', \bar{\omega}_i')$
is exponentially close to $F_{s_i', \omega_i'}$, and no $s_i'$ is zero.}
$\|\eta - \sum F_{s_i, \omega_i} \|_{L^2(M_Y)} \ll Y^{-N}$. 
In particular, if $\eta$ were orthogonal to all $F_{s, \omega} \in \mathcal{F}- \mathcal{F}_0$
then $\eta  = 0$, as claimed (since the matrix of inner products is nondegenerate -- see discussion around  \eqref{innerproductmatrix}).

\subsubsection{Proof  of Lemma~\ref{SmallConstantTermMeansSmall}} \label{SCTMSproof}
We apologize for not giving a unified treatment of $\xi=\eta$ and $\xi = \omega$, but the two cases are slightly different because $\eta$ satisfies
boundary condtions and $\omega$ does not.

 We regard $\xi$ as a $1$-form on $M$ by extending it by zero;
call the result $\tilde{\xi}$.  Of course $\tilde{\xi}$
is not continuous, a point that will cause some technical trouble in a moment. 

{\em Outline:} 
We shall  expand
 $\tilde{\xi}$ in terms of a basis of $1$-forms on $M$ and take $L^2$-norms:
\begin{equation} \label{const-term}  \tilde{\xi} = \sum_{\psi} |  \langle \tilde{\xi}, \psi \rangle |^2 + \sum_{\nu} \int_{t=-\infty}^{\infty} | \langle \tilde{\xi},  E(\nu,  it) \rangle  |^2  \frac{dt}{2\pi}  +
 \sum_{f} \int_{0}^{\infty} | \langle \tilde\xi, \frac{d E(f, it)}{\sqrt{1+t^2}} \rangle|^2 \frac{dt}{2\pi}\end{equation} 
 where the $\nu$-sum is taken over an orthonormal basis for $\Omega^+$,
 the $f$-sum is taken over an orthonormal basis for $C^{\infty}(0)$ (see
~\S~\ref{subsec:einsteinintro} for discussion) 
 and   the first $\psi$-sum is over {\em cuspidal} forms; we use the fact that arithmetic hyperbolic $3$-manifolds have no discrete spectrum on forms, except for the locally constant $0$- and $3$-forms.

  In fact, we later modify this idea slightly to deal
 with the poor convergence of the right-hand side, which behaves like the Fourier series of a non-differentiable function.

{\em Preliminaries:} 
Note that in either case $\xi - \xi_N$ is ``small'': for $\xi=\eta$
the estimate~\eqref{precisebound2} applies directly:
\begin{equation} \label{preciseboundeta} |\eta(x) - \eta_N(x) | \ll \exp(-b \  \height(x)) \|\eta\|_{L^2}, \ \ \height(x) \leq Y.\end{equation} 
Note that the left-hand side only makes sense when $x$ belongs to one of the cusps, i.e., $\height(x) \geq 1$. 
For $\xi=\omega$ one has  \begin{eqnarray} \nonumber |\omega(x)- \omega_N(x)|  &\leq & |\eta-\eta_N|  + |E(s, \omega)-E(s,\omega)_N| + |E(-s, \bar{\omega}) -
E(-s, \bar{\omega})_N| \\ \label{preciseboundomega}  & \leq &  e^{-b \height(x)} \|\eta\|_{L^2}, \ \ \ \height(x) \leq Y. \end{eqnarray}
where we used  now~\eqref{precisebound} for the latter two terms, taking into account
the Maass-Selberg relation~\eqref{MSR1F}, ~\eqref{zerofree}, and the fact
 that $\|\omega_s\|$ and $\| \bar{\omega}_s\|$ are bounded by constant multiples of $\|\eta\|_{L^2}$.
 
  There are bounds of a similar nature for $d\omega-d\omega_N$
 and $d \eta- d\eta_N$.

 For any form, let $\|F\|_{\bdy} $ be the sum of the $L^{\infty}$-norms
 of $F, dF$ and $d^* F$ along the boundary $\partial M_Y$.  \eqref{precisebound2} and variants imply that for $\xi \in \{\eta, \omega\}$ we have, for suitable $a,b,c$, 
\begin{equation} \label{bdynormbound} \|\xi\|_{\bdy} \ll  Y^c \|\xi\|_N +  a e^{-bY} \|\eta\|, \end{equation}
where we define   \begin{equation} \label{Nnormdef} \|\xi_N\|^2 = |a|^2+|b|^2+|a'|^2+|b'|^2 \end{equation} for $\xi_N  = y^s (a dx_1 + b dx_2) +y^{-s} (a' dx_1+ b' dx_2)$. (This norm is not ``good'' near $s=0$, but we have in any case excluded this.) 

The Green identity implies that if $F$ is an eigenform of eigenvalue $-t^2 \leq T^2$ 
and $\xi \in \{\eta,\omega\}$: 
\begin{equation} \label{Greenconq}  \langle F, \xi \rangle_{M_Y} \ll \frac{ \| F\|_{\bdy} \|\xi\|_{\bdy}} {|s|^2-|t|^2}, \end{equation}

{\em The spectral expansion:} 
We now return to~\eqref{const-term}.  To avoid the difficulties of poor convegence previously mentioned, we replaced $\tilde{\xi}$ by a smoothed version. 
 We choose a ``smoothing kernel'': take a smooth self-adjoint section of $\Omega^1 \boxtimes (\Omega^1)^*$
    on $\H^3 \times \H^3$, invariant under $G_{\infty}$, and supported
    in $\{(P, Q): d(P,Q) \leq C\}$, for some constant $C$. Then (the theory of point-pair invariants) 
    there exist smooth functions $\hat{K}_1, \hat{K}_0$
    such that $K \star \omega = \hat{K}_1(s) \omega$ whenever
    $\omega$ is a co-closed $1$-form of eigenvalue $-s^2$, so that $d^* \omega = 0$, 
    and $K \star \omega =\hat{K}_0(s) \omega$ whenever $\omega$ is of the form $df$, 
    where $f$ has 
  eigenvalue $1-s^2$.    We can arrange matters such that      $\hat{K}_\pm $ takes value in $[0,2]$, $\hat{K}_1(s) = 1$, 
    and so that the decay is rapid:  
    \begin{equation} \label{hatKdecay} \mbox{$\hat{K}_{1,2}(u)$ decays faster than any positive power of
    $(T/u)$.} \end{equation} 
   Now convolving with $K$ and applying spectral expansion:

\begin{eqnarray}      \label{etaparseval}
\begin{aligned}
\langle \tilde{\xi} \star K,  \tilde{\xi} \star K \rangle=  & \ 
\sum_{\psi \ \mathrm{coclosed}}  |K_{1}(s_{\psi})|^2 | \langle \tilde{\xi}, \psi \rangle| ^2 + 
\sum_{\psi \perp \mathrm{coclosed}}  |K_{0}(s_{\psi})|^2 | \langle \tilde{\xi}, \psi \rangle| ^2  \\ & \  + 
 \sum_{\nu} \int_{t=-\infty}^{\infty} \frac{dt}{2 \pi}  |\hat{K}_1(it)|^2 \ | \langle \tilde{\xi}, E(\nu, it) \rangle|^2 
 \\   & \   + 
  \sum_{f} \int_{t}\frac{dt}{2\pi}  |\hat{K}_0(it)|^2 \ | \langle \tilde{\xi}, \frac{dE(f, it)}{\sqrt{1+t^2}} \rangle|^2 
  \end{aligned}
\end{eqnarray}  
Here $-s_{\psi}^2$ is the Laplacian eigenvalue of $\psi$.
The inclusion of the factors
 $|K(\cdots)|^2$ makes the right-hand side rapidly convergent.

We will now show that, for $\xi = \eta$ or $\xi  = \omega$, 
\begin{equation} \label{nowverified} \| \tilde{\xi} \star K\|_{L^2} \ll a Y^b \|\xi\|_{\bdy}  + a' Y^{-N} \|\eta\|_{L^2(M_Y)}\end{equation} 
for arbitrary $N$, and $a,b,a'$ possibly depending on $N$.  This is all we need to finish the proof of the Lemma: Note that $\tilde{\xi} \star K$ agrees with $\xi$ on $M_{\leq c Y}$
(for a suitable constant $c \in (0,1)$ depending on the support of $K$). 
Now  ~\eqref{preciseboundeta} and~\eqref{preciseboundomega} imply that the $L^2$-norm of $\|\xi\|$
 in the region $c Y \leq y \leq Y$ is bounded by $a \exp(-bY) \|\eta\| + a' Y^{b'} \|\xi_N\|$. 
Therefore,
$$
\begin{aligned} 
 \|\xi\|_{L^2}   \ll  & \  \|\tilde{\xi} \star K \|_{L^2} +  \int_{cY \leq \height \leq Y} |\xi|^2  \\
  \ll &\  \|\tilde{\xi} \star K\|_{L^2} + a' Y^{b'} \|\xi_N\| + a \exp(-b Y) \|\eta\| \\
   \ll & \ {a_3} Y^{b_3} \|\xi_N\|  + a_4 Y^{-N} \|\eta\| \end{aligned}
 $$
 where we used, at the last step,~\eqref{nowverified} and~\eqref{bdynormbound}.  
 Again, $a,b, a', \dots$ are suitable constants, which may depend on $N$.

So, a proof of ~\eqref{etaparseval} 
will finish the proof of  Lemma~\ref{SmallConstantTermMeansSmall}.

We discuss each term on the  right-hand side of~\eqref{etaparseval} in turn.
In what follows, when we write ``exponentially small'' we mean
bounded by $e^{-b Y} \|\eta\|_{L^2}$; we shall use without comment
  the easily verified fact that $\frac{\|\omega\|_{L^2(M_Y)}}{\|\eta\|_{L^2}}$
is bounded by a polynomial in $Y$. 
\begin{itemize} 
\item Large eigenvalues:    
Either for $\xi=\eta$ or $\xi=\omega$, the total contribution  to the right-hand side of~\eqref{etaparseval}
of terms with eigenvalue ``large'' (say $\geq Y^{1/10}$) is   $\ll Y^{-N} \|\xi\|$ 
for any $N > 0$; this follows by 
  the assumed decay \eqref{hatKdecay} i.e. that $\hat{K}_{0,1}$ decays faster than any positive power of $(T/u)$.
  
  In what follows, then, we may assume that we are only considering the inner product
  of $\tilde{\xi}$ with eigenfunctions of eigenvalue $\leq Y^{1/10}$. 
  In particular, we may apply~\eqref{precisebound} to such eigenfunctions.

  \item     Coclosed cuspidal $\psi$ with eigenvalue in $[\epsilon_0, T^2]$: 
  
  To analyze $\langle \eta,\psi \rangle$, when $\psi$ is co-closed cuspidal,  we recall that $\eta \perp F_{\psi}$ and $F_{\psi}$ is exponentially close to $\psi$, at least
whenever the eigenvalue of $\psi $ lies in $[\epsilon_0, T^2]$. 
Thus such terms are an exponentially small multiple of $\|\eta\|$. 
 
  To analyze $\langle \omega, \psi \rangle$ we simply note that
$\psi$ is perpendicular to $\eta - \omega = E(s, \omega) + E(-s, \bar{\omega})$ on $M$,
and so (because $\psi$ is cuspidal)  almost perpendicular on $M_Y$: for $\psi$ is exponentially small on $M-M_Y$. 

\item Coclosed cuspidal $\psi$ with eigenvalue not in $[\epsilon_0, T^2]$:

An application of Green's theorem~\eqref{Greenconq} shows that $\langle \tilde \xi, \psi \rangle$
is an exponentially small multiple of $\|\xi\|_{\bdy}$ for such $\psi$. 

Also, the number of such eigenvalues which are $\leq Y^{1/10}$
is at most a polynomial in $Y$.

\item Cuspidal $\psi$ that are perpendicular to co-closed:

To analyze $\langle \eta, \psi\rangle$ when $\psi = d f$ for some
cuspidal eigenfunction $f$, we note that in fact
$\langle \eta, df \rangle = \langle d^* \eta, f \rangle = 0 $;
the boundary term $\int_{\partial M_Y} * \eta \wedge f$ is identically zero
because $\eta$ satisfies absolute boundary conditions. 

As before, $\psi$ is exactly orthogonal to $\eta-\omega$, so almost orthogonal to it on $M_Y$, and so again
  $\langle \omega, \psi \rangle$ is bounded by an exponentially small multiple of $\|\eta\|$. 

\item Eisenstein terms: 
 
From~\eqref{Greenconq} we get the estimate
\begin{equation}\label{cob1}   |\langle \tilde\xi, E(\nu, t) \rangle_{M_Y}|  \ll a Y^b\frac{ \|\xi\|_{\bdy}}{|t|^2-|s|^2}, \end{equation}
and there is a similar expression for the other Eisenstein term.

This estimate is not good when $t$ is very close to $s$ (say, when $|s-t| \leq 1$).
To handle that region, use Cauchy's integral formula and analytic continuation in the $t$-variable  to compute $\langle \tilde{\xi}, E(\nu, t) \rangle$,
this expression being antiholomorphic in $t$; this shows that the estimate
\begin{equation} \label{cob2}   |\langle \tilde\xi, E(\nu, t) \rangle_{M_Y}|  \ll a Y^b \|\xi\|_{\bdy} , \end{equation}
continues to hold when $|t-s|\leq 1$. Here we have  used the fact that there exists a small disc around any $t \in i\R$
such that $E(\nu, t)$ remains holomorphic in that disc, together with bounds on $E(\nu, t)$ there;
see discussion around \eqref{trivesteis}.

Now \eqref{cob1} and \eqref{cob2} suffice to bound the Eisenstein contribution to ~\eqref{etaparseval}.
\end{itemize}
  \qed

\subsection{Modifications for functions, \texorpdfstring{$2$}{2}-forms, and \texorpdfstring{$3$}{3}-forms; conclusion of the proof.}  \label{sec:modifications}

Since the proofs for (c) and (d) of the Theorem are largely the same as the proof of part (b), 
which we have now given in detail, we simply sketch
the differences, mostly related to behavior of $\Psi(s)$ near $1$.

\subsubsection{Informal discussion} \label{sss:informal}   Recall that $\Psi(s)$ acts on an $h$-dimensional space, where $h$ is the number of cusps.
 At $s=1$, it has a pole whose residue is a projection onto a $b_0$-dimensional space,
 where $b_0$ is the number of components. The behavior of $\Psi$
 near $1$ gives rise to two spaces of forms on $M_Y$ of interest,
 one of dimension $h-b_0$ and one of dimension $b_0$. 
 These two spaces are closely related to the kernel and image of
$H^2(\partial M) \rightarrow H^3_c(M)$,
which are $h-b_0$ and $b_0$ dimensional respectively.

\begin{itemize}
\item[(a)] For $f$ that belongs to the kernel of the residue of $\Psi(s)$ at $s=1$,
the function $s \mapsto E(s,f)$ is regular at $s=1$. The $1$-form $d E(s,f)$ 
is then {\em harmonic.}  It (approximately) satisfies relative boundary conditions.

These give  rise (after a small modification)  to an $h-b_0$ dimensional space
of forms spanning 
the image of $H^0(\partial M_Y)$ inside $H^1(M_Y, \partial M_Y)$;
or, more relevant to us, the $2$-forms $* d E(s,f)$ give rise
to an $h-b_0$-dimensional space of forms in $H^2(M_Y)$
that map isomorphically to $\ker( H^2(\partial M_Y) \rightarrow H^3_c(M_Y))$. 

This space was already used  (see \S \ref{ip:polygrowth}).

\item[(b)]   The residue of $\Psi$ at $1$ means that 
one gets $b_0$ solutions to $\det(1 + Y^{-2s} \Psi(s)) = 0$
very close to $s=1$. 

This gives rise (after a small modification, as before) to a $b_0$-dimensional space
of functions  on $M_Y$ with very small eigenvalues (proportional to $Y^{-2}$),
and satisfying relative (Dirichlet) boundary conditions.  Roughly these eigenfunctions
``want to be the constant function,'' but the constant function does not satisfy the correct boundary conditions. 

These give the  ``extra eigenvalues'' in part (d) of the Theorem.

\end{itemize}

\subsubsection{Part (c) of the theorem: co-closed \texorpdfstring{$0$}{0}-forms with absolute conditions}

One way in which  the case of $0$-forms is substantially easier:  the ``passage to combinatorial forms'' is not necessary.
 For $1$-forms, this was necessary only to handle eigenvalues on $M_Y$ very close to zero;
 on the other hand, in the setting of $0$-forms and absolute boundary conditions, the lowest nonzero eigenvalue is actually bounded away from $0$.

More precisely, the analog of statement (*) from \S \ref{U2} now holds in the following form: 
 
 \begin{quote} (*)  Any eigenfunction on $M_Y$ with eigenvalue in 
 $[0, (\Tmax-\delta)^2]$ which is orthogonal to all elements
 of the (analog of the sets) $\mathcal{F}$, is in fact zero.
 \end{quote}
 In the current setting, $\mathcal{F}$ is the set of near-eigenforms 
 constructed from cuspidal eigenfunctions on $M$, and also from
 Eisenstein series at parameters where $g(it ) = 0$ for $t > 0$. \footnote{A priori,
 one needs to consider  here not only $s \in i \R$ but also $s \in [0,1]$; 
 however, for sufficiently large $Y$ there are no roots of $g$ in $[0,1]$, 
 because $\mathrm{Id}$ dominates $\frac{1-s}{1+s}  Y^{-2s} \Psi(s)$. } 
 
 This is proved exactly in the fashion of \S \ref{U2}, with no significant modification:
 Let  $f$ be any eigenfunction on $M_Y$ with eigenvalue $1-s^2$,
and we apply ~\eqref{greenid} with $\omega_1$ replaced by $f$
and $\omega_2$  replaced by the Eisenstein series $E(s)$.    There is no $B_1$-term,
and thus the identity simply shows that the constant terms of $f$ and the constant terms of $E(s)$ almost be proportional.   One subtracts a suitable multiple of $E(s)$ from $f$ and proceeds as before.

Note that in the prior analysis -- e.g. prior to \eqref{harr} --
some of our estimates degnerated when $s$ was close to zero, but this
was just an artifact of using a poor basis; it is better   to replace the roles of $y^s, y^{-s}$ by the nondegenerate basis $e_s = y^s+y^{-s}, f_s =  \frac{y^s-y^{-s}}{s}$,  which eliminates
this issue.   This is related to the issue of why we do not include $s=0$ as a root of $g(s) = 0$ in part (c) of the Theorem, i.e. why we only consider roots of $g(it)$ for $t > 0$.  

\subsubsection{ Part (d) of the Theorem: co-closed \texorpdfstring{$2$}{2}-forms with absolute conditions} \label{697d}

By applying the $*$ operator this is the same as computing eigenvalues on closed $1$-forms
with relative conditions, that is to say, away from the zero eigenvalue,  the same as computing eigenvalues on {\em functions}
with relative boundary conditions. 

The analysis  is now similar to that of the previous subsection, 
again now requiring consideration of $s \in [0,1]$; but now there 
is a significant change related to $s \in [0,1]$.

Namely, $g'$ has roots for large $Y$
in the region $s \in [0,1]$, unlike $g(s)$ or $f(s)$.  The reason for this
is the singularity of $\Psi(s)$ at $s=1$; this didn't show up for $g(s)$ because of the factor
$\frac{1-s}{1+s}$ in front. 

It will transpire that the number of such roots  $s$ is, for large $Y$, precisely the number of connected components of $M$, and they all lie very close to $1$. The corresponding
``monster eigenvalues'' $1-s^2$ are thus very close to zero. 
 
Let $R$ be the residue of $\Psi(s)$ at $s=1$; we have computed it explicitly in 
\S~\ref{subsubsec:residues}: it is a projection with nonzero eigenvalues parameterized by connected components of $M$; the eigenvalue corresponding to the connected component $N$
is $\area(\partial N)/\vol(N)$. 

 Now, 
the solutions to 
$$\det(1 + \Psi(s) Y^{-2s}) = 0$$
very close to $s=1$ are well-approximated, by a routine argument,
by the solutions to
$$\det(1+ \frac{R}{s-1} Y^{-2s})  = 0.$$ 
There is one solution to this for each nonzero eigenvalue $\lambda$ of $R$,
very close to $1 - \frac{\lambda}{Y^2}$ (up to an error $o(1/Y^2)$). 
 In particular the corresponding Laplacian eigenvalue is very close
 to $1-s^2 \approx 2 \frac{\lambda}{Y^2} + o(1/Y^2)$. 
 
 We note for later reference, then, that the product of all nonzero ``monster eigenvalues'' 
 is given by \begin{equation} \label{produit} \det'(2 R/Y^2) (1+o(1)),\end{equation}  
where again $\det'$ denotes the product of all nonzero eigenvalues. 

This concludes our sketch of proof of parts (d), (e) of  Theorem~\ref{prop:SE}. 

\subsubsection{Asymptotic behavior of regulators}
We now consider the relationship between the regulator of $M_Y$ and that of $M$.

Clearly
$$\reg(H_0(M_Y)) \longrightarrow \reg(H_0(M)),$$
as $Y \rightarrow \infty$.

We have seen in \eqref{rtymt} that
\begin{equation}  \reg(H_1(M_Y)) \sim \reg(H_1(M)) (\log Y)^{-\hrel} ,\end{equation} 
where $\hrel$ is the number of relevant cusps.

As for 
the $H^2$-regulator, the analogue of~\eqref{rtymt} is the following, which is proved in a similar way (using the space of forms described in part (a) of \S \ref{sss:informal}): 
\begin{equation} \label{rtymt2} \reg(H_2(M_Y)) \sim \reg(H_2(M)) Y^{-2h'} \end{equation} where $h' = \dim(H_2) - \dim(H_{2!})$, and $\sim$ denotes that the ratio approaches $1$.
Indeed,~\eqref{h2normdef} was chosen in such a way as to make this relation simple.

On the other hand, we have seen that
$$\reg(H_2(M)) =\reg(H_{2!}(M)) \cdot  \left( \left( \prod_{\mathcal{C}} 2 \area(\mathcal{C} )\right)   \prod_{N} \vol(N)^{-1}
\cdot \det{}'(2 R)^{-1} \right)^{1/2} $$
where the $\mathcal{C}$ product is taken over all cusps, the $N$-product is taken over all components $N$ of $M$, and $\vol(N)$ denotes the volume of $N$. 
Finally, because $H_3(M_Y, \C) $ is trivial, 
\begin{equation} \label{rtymt3} \reg(H_3(M_Y)) \sim \reg(H_{3, \bm}(M)) \prod_N \vol(N)^{-1/2} \end{equation}

Consequently --- by~\eqref{h2reghelp} --- the ratio
$$ \frac{\reg(M)}{\reg(M_Y)} \sim  \left( \prod_{\mathcal{C}}2  \area(\mathcal{C} )\right)^{1/2}  (\det {}'(2R)  )^{-1/2}
(\log Y)^{\hrel} Y^{-2 h'} ,$$
where the product is taken over all connected components 
 and taking the ratio for $M$ and $M'$ we arrive at:
 \begin{eqnarray} \nonumber  \frac{\reg(M)}{\reg(M_Y)}  \frac{\reg(M'_Y)}{\reg(M')} & \sim&   \left( \frac{ \det'(2R') Y^{2 h'(M')} }{\det'(2R) Y^{2 h'(M)}} \right)^{1/2} 
 \\ &=& \label{chatham} \left( \frac{ \det{}' (2R'/Y^2)}{\det {}' (2R/Y^2) } \right)^{1/2}. \end{eqnarray}since the terms $\prod \area(\mathcal{C})$ and $\log Y$ terms cancel and also
  $h'(M') - h(M)  =  - b_0(M) + b_0(M')$. 
\section{The proof of Theorem~\ref{theorem:invtrunc}} \label{theorem:invtruncproof}

Notation as in the statement of the Theorem.
 
 In what follows we fix {\em two} truncation parameters $1 < Y' < Y$.

Let $K(t; x, y)$ be the the trace of the heat kernel on $j$-forms on $M$, and $k_t(x) = K(t; x,x)$.
We suppress the $j$ for typographical simplicity, but will refer to it when different $j$s are treated differently. 
 We define similarly $K', k'$ for $M'$, and $K_Y, k_Y$ for $M_Y$, 
$K_Y', k_Y'$ for $M'_Y$.   Recall our convention: in the cases of $M_Y, M'_Y$ we compute with absolute boundary conditions, i.e. 
 we work on the space of differential forms $\omega$ such that both $\omega$ and $d\omega$, when contracted with a normal vector, give zero. 
Set also $k_{\infty}(x)  = \lim_{t \rightarrow \infty} k_t(x)$.
Define similarly $k_{\infty}'$ (for $M'$) and $k_{\infty,Y}, k_{\infty, Y'}$ for $M_Y$ and $M'_Y$ respectively.

Set $$\delta_t(x) = (k_t(x) - k_{Y,t}(x)), (x \in M_Y)$$ where
we identify $M_Y$ with a subset of $M$.  We define $\delta_t'$ similarly.  Set
 \begin{multline}   I_j(Y, t) :=  
  \left(  \tr^* \exp(-t \Delta_M)  - \tr \exp(-t \Delta_{M_Y}) \right)  \\ - \left( \tr^* \exp(-t \Delta_{M'}) - \tr \exp(-t \Delta_{M'_Y}) \right) ,\end{multline}
  where the Laplacians are taken on $j$-forms, and
  $I_j(Y, \infty) := \lim_{t \rightarrow \infty} I_j(Y, t)$;
  for the definition of regularized trace see  Definition \ref{RegTraceDef}.

If we put $I_*(Y,t) = \frac{1}{2} \sum_{j} j (-1)^{j+1} I_j(Y,t),$ then  that part of the expression of Theorem~\ref{theorem:invtrunc} involving analytic torsion equals \begin{equation}\label{deriveddesid} -  \frac{d}{ds}  |_{s=0}\Gamma(s)^{-1} \int_{0}^{\infty}(  I_*(Y,t) - I_*(Y, \infty)) t^s \frac{dt}{t},\end{equation} 
Note that, a priori, this should be understood in a fashion similar to that
described after~\eqref{atsecondef}, that is to say, by dividing into intervals $[0,1]$
and $(1, \infty)$ separately.

In what follows, we fix $j$, and write simply
$I(Y,t)$ for $I_j(Y,t)$; this makes the equations look much less cluttered. 
 Then, using~\S~\ref{subsec:etdcomparison}, we may write \begin{eqnarray} \label{iytdec} I(Y, t) &=& 
 \int_{x \in M_{Y'}} \delta_t(x)  -   \int_{M'_{Y'}} \delta'_t(x)  \\ \nonumber &+&   
 \int_{x \in M_{[Y',Y]}} (k_t(x) - k_t'(x)) - (k_{t,Y}(x) - k'_{t,Y}(x))  \\ \nonumber &+& 
\int_{\mathcal{C}_Y} (k_t(x) - k'_t(x)).  \end{eqnarray} 
and $I(Y, \infty)$ is obtained by substituting $k_{\infty}$ for $k_t$ at each point. 
Note, in the third integral on the left hand side, we have used the natural identification
of $M_{[Y', Y]}$ with $M'_{[Y', Y]}$ to regard $k'$ as a function on $M_{[Y', Y]}$. 
  Now we have computed $\lim_{t \rightarrow \infty} \tr^* e^{-t \Delta_M}$
 in~\eqref{ASspec2}; it equals   $$ \begin{cases} b_0 =\dim H_0(M,\C),  & j=0,3 \\ \dim H_{j, !}(M, \C) , \  & j=1,2.\end{cases},$$
  and similarly for $M'$.
   On the other hand,
 $\lim_{t \rightarrow \infty} \tr  \ e^{-t \Delta_{M_Y}} =  \dim H_j(M_Y, \C)$ which equals
  $\dim H_j(M, \C)$, and
 $$ \dim H_j(M, \C) - \dim H_{j, !}(M, \C) = \begin{cases}   \frac{1}{2} \dim H_1(\partial M), & j=1, \\
 \dim H_2(\partial M)  - b_0(M), & j=2 \end{cases},$$
 so that $I(Y, \infty) =0$ for $j \neq 2,3$ --- since the boundaries $\partial M$ and $\partial M'$ are homotopy equivalent, and we get
\begin{equation} \label{firstgamma} \frac{d}{ds}  |_{s=0}\Gamma(s)^{-1} \int_0^T I(Y, \infty) t^s \frac{dt}{t} =
\begin{cases}  (\log(T) + \gamma) (b_0(M) - b_0(M'))  &  (j=2) \\
(\log T + \gamma) (b_0(M)-b_0(M')), & (j=3). \\
0, &  j \in \{0,1\} \end{cases}
\end{equation}

We claim that for any fixed $T > 0$  \begin{equation} \label{double-limit} \lim_{Y \rightarrow \infty} \int_0^T I(Y, t) = 0. \end{equation} 
 Indeed, reasoning term by term in~\eqref{iytdec} (and, where necessary, first taking $Y $ to $\infty$ and then $Y'$), 

 \begin{itemize}
\item $\delta_t(x), \delta_t'(x)$  both approach zero uniformly for  $ x \in M_{Y'}, Y > 2 Y', t \in [0,T]$. 
This is a consequence of the locality of the heat kernel. 
This disposes of the first two terms. 
\item The terms $  k_t(x) - k_t'(x)  $ is independent of $Y, Y'$; 
 uniformly for $t \in [0,T]$, it approaches zero as $\mathrm{height}(x) \rightarrow \infty$. 
Since the measure of $M_{[Y',Y]}$ is bounded above, this disposes
of the second term.

\item As $Y' \rightarrow \infty$, $k_{t,Y}(x) - k'_{t,Y}(x)  $  approaches zero
uniformly for $x \in M_{[Y',Y]}$, $t \in [0,T]$; this follows from the locality of the heat kernel.

\item Similarly, for $x \in \mathcal{C}_Y$, $k_t(x) -k_t(x')$
approaches zero as $Y \rightarrow \infty$, uniformly for $t \in [0,T]$; this is just as in the case of the second term. This disposes of the final term. 
\end{itemize}

This proves~\eqref{double-limit}.
In combination we have shown
that for any fixed $T$ 
\begin{equation} \label{secondgamma}
\lim_{Y\rightarrow \infty} -\frac{d}{ds} |_{s=0} \Gamma(s)^{-1} \int_0^T \left( I_*(Y, t) - I_*(Y, \infty) \right) t^s \frac{dt}{t} =  \frac{1}{2} (\log T + \gamma) (b_0(M)-b_0(M')). \end{equation}
 Now let us turn to the analysis of $t \in [T, \infty)$; 
 we shall show (using crucially the analysis of small eigenvalues)
{\small
\begin{eqnarray} \label{Istar} \\ \nonumber  
\begin{aligned}
-\frac{d}{ds} |_{s=0} \Gamma(s)^{-1}  & \int_{T}^{\infty}( I_*(Y,t) - I_*(Y,\infty)  t^{s} \frac{dt}{t}   
+ \frac{1}{2}  \left( \log \det{} '(Y^{-2} R) - \log \det{}' (Y^{-2} R') \right)  \\  \nonumber
&   \longrightarrow  - \frac{1}{2} (\log T + \gamma) (b_0(M)-b_0(M')), \mbox{ as $Y \rightarrow \infty$} \end{aligned}
 \end{eqnarray}
 }
 Recall that we are using $\det{}'$ to mean the product of nonzero eigenvalues,
 and $R$ was defined  after \eqref{h2reghelp}
Since we already saw~\eqref{chatham} that
\begin{eqnarray*}
  ( \log \RT(M)- \log \RT(M_Y)    )
-   (\log \RT(M') - \log \RT(M'_Y)   )
\\ + 
\frac{1}{2} \left(  \log \det'(Y^{-2} R) - \log \det' (Y^{-2} R')  \right)  \longrightarrow 0.\end{eqnarray*}
 this, together with~\eqref{firstgamma}, will conclude the proof of the theorem;
 note that the analytic torsion quotient differs from~\eqref{Istar} by a sign,
 owing to the sign difference between $\log \det \Delta$ and $\zeta'(0)$.

Now each term of $I_j(Y,t)$ can be written as a summation
(or integral, in the case of continuous spectrum) 
of $e^{-\lambda t}$ over eigenvalues $\lambda$ of 
$j$-forms on one of $M, M_Y, M', M'_Y$.
We define $I_j^{\coclosed}$ simply be restricting these summations or integrals
to be over co-closed forms, that is to say, the kernel of $d^*$. 
 Then one can write each $I_j$ in terms of $I_j^{\coclosed}$; in fact, 
 $I_j(Y,t) - I_j(Y, \infty) = (I_j^{\coclosed}(Y,t) - I_j^{\coclosed}(Y, \infty))   + (I_{j-1}(Y, t)^{\coclosed} - I_{j-1}(Y,\infty)^{\coclosed})$. 
 In this way we reduce the analysis of \eqref{Istar} to the coclosed case.

Following our general policy of treating $1$-forms in detail, we shall now show
\begin{equation} \label{0formslargetime}\frac{d}{ds} |_{s=0} \Gamma(s)^{-1} \int_{T}^{\infty} \left( I_1^{\coclosed}(Y,t) 
- I_1^{\coclosed}(Y, \infty) \right) t^{s} \frac{dt}{t}  \rightarrow 0
 \end{equation}
(i.e., the contribution of coclosed $1$-forms.)  We shall then explain how the modifications of part (e) of Theorem~\ref{prop:SE} causes a slightly different analysis for $j \neq 1$, 
 giving rise to the terms on the right of~\eqref{Istar}.

Write $N(x)$ for the number of eigenvalues on co-closed $1$-forms on $M_Y$ in $(0,x^2]$; enumerate these eigenvalues $0< \lambda_1\leq \lambda_2 \leq \dots$  
Define correspondingly $N_{\cusp}(x)$ for the number of such eigenvalues on co-closed cuspidal $1$-forms on $M$.  Note, in both cases, we are not counting eigenvalue $0$. 

Set $E(x)$ to be the ``error'' when we approximate $N(x)$ 
by means of Lemma \ref{analysisofzeroes}:
\begin{equation} \label{MoNoFo} E(x) = N(x) - \frac{ \left(4 x \hrel \log(Y)  - \sum_{i=1}^h \nu_i(x) \right)}{2 \pi}  - N_{\cusp}(x), \end{equation} 
where $\hrel$ is the size of the scattering matrix $\Phi^{\pm}$
and $\nu_i$ is as defined as before (cf. Lemma  \ref{analysisofzeroes}, (i)).

Write $\omega(s) = \det (\Phi^{-}(-s) \Phi^+(s))$. Since $\Phi^+(s) \Phi^{-}(s) $ is the identity, it follows that $\omega(s) \omega(-s) = 1$; in particular, as we will use later, the function
$ \omega'/\omega$ is {\em even} (symmetric under $s \mapsto -s$). 
 
According to part (i) of the Lemma  and Theorem~\ref{prop:SE}, 
or more directly  \eqref{thm:coclosed1},
$E(x)$ is bounded when in the range $x \leq \Tmax = \Tmax(Y), Y \geq Y_0$;  even more precisely $E(x)$ is ``well-approximated'' by 
 $- \sum \left\{ \frac{4 x \log Y - \nu_i(x) }{2 \pi} \right\}$ --- in the sense that
\begin{equation} \label{Exapp} \int_{X}^{X+1} \left| E(x)  + \sum_{i}\{\frac{4 x \log(Y) - \nu_i}{2\pi} \} \right|  \leq a \exp(-bY) \end{equation} 
for some absolute constants $a,b$, whenever $X \leq \Tmax$. 
Here $\{ \cdot \}$ denotes, as usual, fractional part.  
 
 Now let us analyze $\sum e^{-\lambda_i t}$, the sum being taken over
 all nonzero igenvalues of $1$-forms on $M_Y$. 
 \begin{eqnarray}  \label{goldstar} \\ \nonumber 
 & &  \sum_{\lambda_i \neq 0} e^{-\lambda_i t} 
      =  \int_{0}^{\infty}  e^{- x^2 t} dN(x)  \\ \nonumber 
 &\stackrel{=}{(a)}&    \int_{0}^{\infty}   N(x) \cdot  2 xt e^{-x^2 t}dx  
\\  \nonumber &\stackrel{=}{(b)}& 
\sum \frac{ - \nu_i(0) }{2 \pi}  +  \frac{1}{2 \pi} \int_{0}^{\infty} e^{-x^2 t} \left(   4 \hrel \log(Y) - \frac{\omega'}{\omega}(ix) +  2 \pi \frac{d N_{\cusp}}{dx} \right) dx  \\ \nonumber &  & +
 \int_{0}^{\infty}  E(x) \cdot 2 xt e^{-x^2 t}  \\  \nonumber & \stackrel{=}{(c)} & 
  \frac{1}{2 \pi} \int_{-\infty}^{\infty} e^{-x^2 t} \left(  2 \hrel \log(Y) - \frac{\omega'}{\omega}(ix)  
  \right) +  \sum_{\lambda \neq 0, \mathrm{ cusp \ eig.}} e^{-\lambda t}
  \\ \nonumber &  & +
   \int_{0}^{\infty}  E(x) \cdot 2 xt e^{-x^2 t}  .
\end{eqnarray}  
 Note when we write, e.g. $\frac{d N_{\cusp}}{dt} \cdot dt$, it should be understood as a distribution (it is, in fact, a measure).

Step (a) is integration by parts.  
At step (b), we first expanded $N(x)$ according to~\eqref{MoNoFo}, then used integration by parts ``in the reverse direction'' and the observation
  that the derivative of $\sum\nu_i$ is $-\frac{\omega'}{\omega}$ by \eqref{sumnu}. 
  In step (c) we used  the fact that all $\nu_i(0) =0$,
  this because 2$\Phi^{-}(-s) \Phi^+(s)$ is the identity when $s=0$ \footnote{Were this not so,  this term would be $\frac{-1}{4} \tr(I - \Phi^- \Phi^+(0))$, which is related to a similar term in the trace formula. } and unfolded the integral from $[0,\infty)$ to $(-\infty, \infty)$, using the fact that $\omega'/\omega$ is an even function.
  
  The prior equation \eqref{goldstar} says that the heat-trace on $M_Y$
  is closely related to the regularized heat-trace on $M$, up to an error term controlled by $E(x)$. 
 
 To proceed further we need to analyze the behavior of $E$,  using:
   \begin{lemma} \label{sillyc}
   Suppose $I$ is an open interval in $\mathbf{R}$ and $m$ a   monotone increasing piecewise differentiable function from $I$
   to $(-1/2,1/2)$, satisfying $m' \in [A, B]$ where $B-A \geq 1$. 
Suppose that $\phi: I \rightarrow \mathbf{R}$ is smooth and $|\phi| + |\phi'| \leq M$. 
Then       $$\int \phi(x) \cdot  m(x)  dx \ll
   M \left( \frac{B-A}{A^2}  \right).$$
   \end{lemma}
   \begin{proof}
Translating $I$, we may suppose that
$I = (a,b)$ and $m(0) =0$, with $a < 0 < b$. 
 Now $Ax \leq m(x) \leq Bx$ for $x > 0$;  therefore, $b \in [\frac{1}{2B}, \frac{1}{2A}]$.
 Similarly, $-a$ belongs to the same interval. 
 
 Also
\begin{eqnarray}  \nonumber \int_{I} \phi(x) m(x) & =&    \int_I A x \phi(0) dx + \int_I Ax (\phi(x) -\phi(0)) +    \int_{I} (m(x) -A x) \phi(x) dx 
\\ \nonumber &=& A \phi(0) (b^2-a^2)/2 + O(A . M. \int_I |x|^2 ) +
O((B-A) M \int_I |x|) 
 \\  \nonumber & = & O ( (B-A) M/A^2). \end{eqnarray}
 where we used the fact that $\int_{I} |x| \ll  \frac{1}{A^2}$. 
\end{proof}

 Return to~\eqref{0formslargetime}. 
 In view of~\eqref{goldstar} and its analogue  on the manifold $M_Y'$ (write $E'$ for the analogue of $E$ on $M'_Y$), together with~\eqref{ASspec} describing the regularized trace, 
we may express   $\int_{T}^{\infty} \left( I_1^{\coclosed}(Y,t) - I_1^{\coclosed} (Y, \infty) \right) \frac{dt}{t} $ as 
 \begin{eqnarray} \int_{T}^{\infty} dt \left( \int_{0}^{\infty} E(x) - E'(x) \right) 2x e^{-x^2 t}  dx
\\  \label{bug2} = 
2 \int_{0}^{\infty} (E(x) - E'(x))  e^{-x^2 T} \frac{dx}{x} \end{eqnarray} 
indeed this is absolutely convergent, so the interchange in order is justified. 

Fix $0 < \epsilon < 1$. We  estimate \eqref{bug2}
by splitting the integral into ranges  of $x$:
$$[0,\frac{1}{100 \log Y}) \cup [\frac{1}{100 \log Y}, \epsilon) \cup
 [ \epsilon, \Tmax) \cup [\Tmax, \infty).$$ Recall that $\Tmax$ was a parameter defined  in \eqref{Tmaxdef}, and our analysis of small eigenvalues 
works ``up to eigenvalue $\Tmax^2$.'' 

\begin{itemize}

 \item[-]
 The integral $\int_{\Tmax}^{\infty}$ is handled by trivial estimates.
 Indeed, $E(x)$ and $E'(x)$ are easily bounded by a polynomial in $x$
  all of whose coefficients are at worst polynomials in $\log Y$. 
  But, given the explicit definition \eqref{Tmaxdef}, the size of $\int_{\Tmax}^{\infty} x^N e^{-x^2 T}$ decays faster than any polynomial in $\log Y$. 

\item[-] 
 The integral $\int_{\epsilon}^{\Tmax}$
is seen to be $O_{\epsilon,T}((\log Y)^{-1/2})$:

We approximate $E$ and $E'$ by means of \eqref{Exapp}.  
Thus the integral in question becomes a sum of integrals of the form
$$ \int_{\epsilon}^{\Tmax}  e^{-x^2 T} \frac{dx}{x}  \sum_{i}\{\frac{4 x \log(Y) - \nu_i}{2\pi}\},$$
up to an exponentially small error.
We handle the integral for each $i$ separately in what follows. 

Now split the range of integration $(\epsilon, T)$ as a union $\bigcup_i (a_j, a_{j+1})$
so that $\{ \frac{4 x \log(Y) - \nu_i}{2 \pi}\}$ goes from $-1/2$ to $1/2$ as $x$ increase from $a_i$ to $a_{i+1}$.
The number of such intervals is $O(\Tmax \cdot \log Y)$. 
(At the endpoints, one can move $\epsilon$ and $\Tmax$ by an amount $\approx 1/\log Y$
to ensure this is true; the integral over these small ``edge intervals''
has in any case size $O_{\epsilon}((\log Y)^{-1}$.)

We now apply Lemma \ref{sillyc} to each $\int_{a_i}^{a_{i+1}}$. 
By virtue of \eqref{nuiest}, we may take $A = 4 \log Y - a (\log \log Y)^b$ and $B = 4  \log Y + a (\log \log Y)^b$
 for appropriate $a,b$, and also $M$ to be the maximum value of $e^{-x^2 T}/x + (e^{-x^2 T}/x)'$
 in the range $(\epsilon,  \Tmax)$, which is bounded by $O_{\epsilon,T}(1)$.   
\item[-]   
 To handle $\int_{0}^{\frac{1}{100 \log Y}}$ note that  \begin{equation} \label{EEbound} |E(x) - E'(x)| \leq \sum_{i} |\nu_i| + |\nu'_i| \end{equation} for $x \leq \frac{1}{100 \log Y}$;  this follows  from the definition~\eqref{MoNoFo} of $E, E'$
 together with the fact that there are no nonzero cuspidal eigenvalues $x^2$ on $M$,
 or any eigenvalues $x^2 $ on $M_Y$, for this range of $x$.  
 Since $|\nu_i(x)|/x$ is  bounded for $x \in (0,  \frac{1}{100 \log Y})$, independent of $Y$, 
  it follows
 that the integral $\int_0^{\frac{1}{100 \log Y}}$ is bounded by $O ( (\log Y)^{-1/2})$.

  \item[-] It remains to consider $\int_{\frac{1}{100 \log Y}}^{\epsilon}$. 
  
  By using \eqref{Exapp} it suffices to bound instead
  
$$ \int_{\frac{1}{100 \log Y}}^{\epsilon} \frac{dx}{x}  \left|   \sum_{i}\{\frac{4 x \log(Y) - \nu_i}{2\pi} \} - 
  \sum_{i}\{\frac{4 x \log(Y) - \nu_i'}{2\pi} \}  \right| ,$$
  the difference between the two integrals going to zero as $Y \rightarrow \infty$. 
  
  To bound this it will suffice to bound, instead, the corresponding difference of integer parts
  $$ \int_{\frac{1}{100 \log Y}}^{\epsilon} \frac{dx}{x} \left|   \sum_{i} \left[ \frac{4 x \log(Y) - \nu_i}{2\pi} \right] - 
  \sum_{i}\left[ \frac{4 x \log(Y) - \nu_i'}{2\pi} \right]  \right| ,$$
  because the difference is bounded, up to constants, by the sum of quantities
  $|\nu_i(x)/x|$; and the integral of $|\nu_i(x)/x|$ over the region in question is bounded by
  $O(\epsilon)$ -- after all, $|\nu_i(x)/x|$ is  again absolutely bounded for $x \in [\frac{1}{100 \log Y},\epsilon)$, independent of $Y$.

  Now the functions $4x \log(Y) - \nu_i$ or $4 x \log(Y) - \nu_i'$, considered
  just in the region $[0,\epsilon)$, are all monotonic,
  and they all cross $2 \pi \Z, 4 \pi \Z$ etc. ``at nearby points.''  More precisely, 
  if we write $a_{i,n}$ for the solution
  to   $4 x \log(Y) - \nu_i(x) = 2 \pi n,$
  then $|a_{i,n} - a_{j,n} | \ll \frac{n}{(\log Y)^2}$;
  and a similar result holds comparing $a_{i,n}$ and the analogously defined $a_{j,n}'$ for $\nu_j'$. 
Each such value of $n$ up to $O(\log Y \cdot \epsilon)$,  then, contributes at most $\frac{n}{(\log Y)^2} \cdot \frac{\log Y}{n} = O(\frac{1}{\log Y})$ to the above integral. In this way,  the integral above is bounded above by $O(\epsilon)$. 

    \end{itemize}
    
These bounds show that $$\limsup_{Y \rightarrow \infty} \left| \int^{\infty}_{T} (I_1^{\coclosed}(Y,t) - I_1^{\coclosed}(Y,\infty)) \frac{dt}{t} \right| = O(\epsilon)$$ but $\epsilon$ is arbitrary; so the limit is zero. 
 We have concluded the proof of~\eqref{0formslargetime}.

 To check~\eqref{Istar}
 we now outline the changes in the prior analysis  
 in the other case (i.e.  co-closed $j$-forms for $j \neq 1$). 
 
  \subsubsection{The analysis for  \texorpdfstring{$I_j$}{I_j} for \texorpdfstring{$j \neq 1$}{j ne 1}.}

  \begin{itemize}
  
  \item The analysis of $0$-forms, i.e. $I_0^{\coclosed}$, is the same, indeed simpler.
  
  The role $\Phi^-(-it) \Phi^+(it)$ is replaced by  $\Psi(it)$.  
The main difference: the lowest eigenvalue on the continuous spectrum  of $M$
is now not zero, and also  the nonzero eigenvalue spectrum of $M_Y$ is bounded away from zero
uniformly in $Y$. 
  This avoids many complications in the above proof. 

 \item The analysis of $2$-forms, i.e. $I_2^{\coclosed}$, 
 is  slightly differently owing  to the different formulation of (d) of Theorem~\ref{prop:SE}.

Indeed, Theorem~\ref{prop:SE} (d) and the discussion of~\S~\ref{sec:modifications} 
especially \S \ref{697d}
describes certain co-closed eigenvalues on $2$-forms for  $M_Y, M'_Y$
that are very close to zero.    
 Denote these ``monster eigenvalues'' by $\lambda_i$ for $M_Y$, and $\mu_j$ for$M'_Y$.

  The same analysis that we have just presented now shows that
 \begin{equation} \label{jardindesplantes} \int_T^{\infty} \frac{dt}{t} (I_2^{\coclosed}(Y, t) - I_2^{\coclosed}(Y,\infty)  =
-  \int_{T}^{\infty} \frac{dt}{t}
 \left( \sum_{i} e^{-\lambda_i t} - \sum_{j} e^{-\mu_j t}   \right) + o(1), \end{equation} 
   Recalling that 
$  \int_T^{\infty} e^{-\lambda t} \frac{dt} {t}  =   \int_{\infty}^{T \lambda} \frac{e^{-u}}{u} du 
\sim - \log(T \lambda) - \gamma + o(1),$  for small $\lambda$, and using the analysis of~\S~\ref{sec:modifications}
(see~\eqref{produit}) 
we obtain  that the left-hand side~\eqref{jardindesplantes} equals
  $$o(1) + \left(  \log \det {} '(2 Y^{-2} R) - \log \det {} '(2 Y^{-2} R') \right)    +  (b_0(M) - b_0(M'))( \log T + \gamma)$$
  as $Y \rightarrow \infty$. 
  This accounts exactly for the extra terms in \eqref{Istar}. 

 \item
 Coclosed $3$-forms are simply multiples of the volume and all have eigenvalue $0$. 
   \end{itemize}

  This finishes the proof of~\eqref{Istar} and so also finishes the proof of   Theorem
~\ref{theorem:invtrunc}.

\chapter{Comparisons  between  Jacquet--Langlands pairs}
\label{chapter:ch6}

Our main goal in this chapter is to compare the torsion homology of two manifolds in a Jacquet--Langlands pair $Y, Y'$.    

We begin with notation (\S \ref{JLChapternotn}), recollections on the classical Jacquet--Langlands correspondence (\S \ref{ss:JLclassic}), and on the notion of newform (\S \ref{sec:newnewnew}). 

After some background, we state and prove the crudest form of a comparison of torsion homology 
(\S \ref{sec:coh}, Theorem~\ref{TJL1}) and we then interpret ``volume factors'' in that Theorem as being related to congruence homology, thus giving a slightly refined statement (\S \ref{section:refined}, Theorem~\ref{TJL2}).   
These statements are quite crude -- they control only certain {\em ratios }of orders of groups, and we
devote the remainder of the Chapter to trying to reinterpret them as relations between the orders of spaces of newforms.

In \S~\ref{section:essential} we introduce the notion of {\em essential homology} and {\em dual-essential homology}. These are two variants of homology which (in different ways)
``cut out'' congruence homology. We observe that Ihara's lemma becomes particularly clean
when phrased in terms of essential homology (Theorem \ref{essentialgood}).

In~\S~\ref{section:relations} we start on the
core business of the chapter: the  matter of trying to interpret Theorem~\ref{TJL2} as an equality between
orders of newforms.  To do so, we need to identify the alternating ratios of orders from Theorem~\ref{TJL2}
with orders of newforms.  
In \S~\ref{section:relations} we consider several special cases
(Theorem~\ref{theorem:firstsetting},~\ref{theorem:newforms},~\ref{theorem:thirdsetting})
where one can obtain somewhat satisfactory results. 

In the final sections (~\S~\ref{section:nonny} and~\S~\ref{sec:nonny2}), we attempt to generalize~\S~\ref{section:relations}, using 
  the spectral  sequence of Chapter~\ref{chapter:ch3}. 
The strength of our results
is limited  because of our poor knowledge about homology of $S$-arithmetic groups; nonetheless, we see that certain mysterious factors from~\S~\ref{section:relations}
are naturally related to the order of $K_2$. 
{\em We advise the reader to skip~\S~\ref{section:nonny} and~\S~\ref{sec:nonny2}  at first reading } -- proceeding directly to~\S~\ref{section:finalsummary}, 
where we summarize the results of these sections.

\section{Notation} \label{JLChapternotn} 

Recall that we have defined the notion of a {\em Jacquet--Langlands pair} in~\S~\ref{ss:jlc}, and we recall here:  given  sets $\Sigma, S, T, S', T'$ of finite places satisfying:
$$\Sigma = S \coprod T = S' \coprod T'$$
 we take $\G, \G'$ 
to be  inner forms of $\PGL(2)/F$ that  are ramified at the set of finite primes $S$ and $S'$ respectively and ramified at all real infinite places of $F$.
 Let $Y(K_{\Sigma})$ be the ``level $\Sigma$'' arithmetic manifold
 associated to $\G$ and $Y'(K_{\Sigma})$ the analogous
construction for $\G'$ (see~\S~\ref{section:ar1} for details).

We refer to a pair of such manifolds as a {\em Jacquet--Langlands pair.}

We recall that $Y(\Sigma)$ and $Y'(\Sigma)$ have the same number of connected components
(see~\eqref{conncompremark}). We refer to this common number as $b_0$
or sometimes as $\# Y = \# Y'$. 

  \section{The classical Jacquet Langlands correspondence} \label{ss:JLclassic}

As mentioned in the introduction to this Chapter, the comparison of torsion orders rests
on the classical Jacquet--Langlands correspondence. 
Accordingly we recall a corollary to it, for the convenience of the reader.
We continue with the notation of~\S~\ref{JLChapternotn}. 

\medskip

For every $j$, let $\Omega^j(Y(K_{\Sigma}))^{\new}$
denote the space of {\em new cuspidal} differential $j$-forms on $Y(K_{\Sigma})$.
By this we shall mean the orthogonal complement, on cuspidal $j$-forms, 
of  the image of all degeneracy maps $Y(K_R) \rightarrow Y(K_{\Sigma})$, 
$S \subset R \subset \Sigma$.  If $Y(K_{\Sigma})$ is noncompact 
(i.e., $S$ is empty) {\em then we restrict only to cuspidal $j$-forms}.
Then the Jacquet--Langlands correspondence implies that
\begin{quote} The Laplacian $\Delta_j$ acting on the two spaces of
forms  $\Omega^j(Y(K_{\Sigma}))^{\new}$
and $\Omega^j(Y'(K_{\Sigma}))^{\new}$ are {\em isospectral}; more precisely, there is,
for every $\lambda \in \R$,  an isomorphism of $\lambda$-eigenspaces:
$$  \Omega^j(Y(K_{\Sigma}))^{\new}_{\lambda} \stackrel{\sim}{\longrightarrow} 
\Omega^j(Y'(K_{\Sigma}))^{\new}_{\lambda}$$
that is equivariant for the action of all Hecke operators outside $\Sigma$. 
\end{quote}

 For the Jacquet--Langlands correspondence in its most precise formulation, 
 we refer to  \cite[Chapter 16]{JL};
 for a partial translation (in a slightly different context) into a language
 closer to the above setting, see \cite{Vigneras}. 

Specializing to the case  $\lambda = 0$, we obtain the Hecke-equivariant isomorphism
 $$H^j( Y(K_{\Sigma}), \C)^{\new} \stackrel{\sim}{\longrightarrow} H^j(Y'(K_{\Sigma}), \C)^{\new}.$$
We recall the notion of new in the (co)homological setting below.

\section{Newforms, new homology, new torsion, new regulator}
\label{sec:newnewnew}

Recall  ( \S \ref{Sec:Newforms}) that the {\em ${\q}$-new} space of the first homology $H_1(Y(K_{\Sigma}), \Z)$ is the quotient of all homology classes coming from level $\Sigma - {\q}$.  

\medskip

So $H_1(Y(K_{\Sigma}), \C)^{\new}$ will be, by definition, the cokernel of the map
$$  \bigoplus_{v \in \Sigma - S}
H_1(\Sigma - \{v\}, \C)^2  \stackrel{\Psi^{\vee}}{ \longrightarrow} H_1(\Sigma, \C) $$
which (since we are working in characteristic zero) is naturally isomorphic to the kernel of the push-forward map 
$H_1(\Sigma, \C)  \stackrel{\Psi}{ \longrightarrow} \bigoplus H_1(\Sigma - \{v\}, \C)^2  $.
However, when we work over $\Z$, the cokernel definition will be the correct one.
 
\medskip

The classical Jacquet--Langlands correspondence (that we have just discussed) proves the existence of an isomorphism
$$H_1(Y(K_{\Sigma}), \C)^{\new} \stackrel{\sim}{\rightarrow} H_1(Y'(K_{\Sigma}), \C)^{\new},$$
which is equivariant for the action of all prime-to-$\Sigma$ Hecke operators. 
\medskip

In the next section we shall formulate certain corresponding theorems for torsion, but only at a numerical level. To motivate the definitions that follow, note (see \S \ref{sec:LRC}) the existence of a sequence, where all maps are degeneracy maps: 
{\tiny 
\begin{equation} \label{charzeroraising} H_1(\Sigma - S, \C)^{2^{|S|}}  \longrightarrow \cdots \longrightarrow   \bigoplus_{\{v, w\} \subset \Sigma - S} H_1(\Sigma-\{v, w\}, \C)^4 \longrightarrow  \bigoplus_{v \in \Sigma - S}
H_1(\Sigma - \{v\}, \C)^2 \longrightarrow H_1(\Sigma, \C) .\end{equation}}
It has homology only at the last term, and this homology is isomorphic to $H_1(\Sigma, \C)^{\new}$. 
Therefore, $H_1(\Sigma, \C)^{\new}$ can be decomposed, as a virtual module for the (prime-to-$\Sigma$) Hecke algebra, as an alternating sum of the left-hand vector spaces. 

We now define: 
 
 \begin{df} 
 If $S$ denotes the set of finite places of $F$ which ramify in $D$, \label{df:newstuff}
we define the {\em new regulator}, the {\em new naive torsion}, and the {\em new volume}   as follows:
\begin{eqnarray*} \reg^{\new}(\Sigma) &=& \prod_{S \subset R \subset \Sigma} \reg(Y(K_R))^{(-2)^{|\Sigma \setminus R|}}\\ 
\torsnew(\Sigma) &=& \prod_{S \subset R \subset \Sigma} |H_{1, \tors}(Y(K_R), \Z)|^{(-2)^{|\Sigma \setminus R|}} \\
\vol^{\new}(\Sigma) &=& \prod_{S \subset R \subset \Sigma} \vol(Y(K(R))^{(-2)^{|\Sigma \backslash R|}} \\ 
 \reg^{\new}_{\one}(H_1(\Sigma)) &=& \frac{\reg^{\new}(\Sigma)}{\vol^{\new}(\Sigma),} \\ 
 \end{eqnarray*}
  \end{df}
  If we wish to emphasize the dependence on $Y$ in this definition (with $\Sigma$ explicit) we write
  $\reg^{\new}(Y)$, etc. etc. Observe that $h^{\new}$ is not, {\em a priori}, the order of any group; 
  indeed much of the later part of this Chapter is concerned with attempting to relate $h^{\new}$
  and the size of a newform space. 
   
     The above definitions are equally valid in the split or non-split case;
      however, the definition of regulator varies slightly
  between the two cases (compare~\S~\ref{sec:CMthm} and~\S~\ref{subsec:rtnonsplitdef}). 
  
The quantity $\reg^{\new}_{\one}$ is meant to signify the ``essential'' part of the regulator, 
i.e. what remains after removing the volume part.  One computes from the definition that:  

      \begin{eqnarray*}  \reg^{\new}_{\one} &=& \prod_{S \subset R \subset \Sigma} \left(   \frac{ \reg(H_{1}(R, \Z))}{ \reg(H_2(R ,\Z) }  \right) ^{(-2)^{|\Sigma \setminus R|}}, \ \mbox{ $\G$ nonsplit},\\
    &=&  \prod_{S \subset R \subset \Sigma} \left(   \frac{ \reg(H_{1}(R, \Z))}{ \reg(H_{2!}(R ,\Z) }  \right) ^{(-2)^{|\Sigma \setminus R|}},
    \   \mbox{$\G$ split}.
   \end{eqnarray*}

\medskip

\section{Torsion Jacquet--Langlands, crudest form} \label{sec:coh}

We continue with the previous notation; in particular, $Y, Y'$ are a Jacquet--Langlands pair.

\begin{theorem}  \label{TJL1}  The quantity 
\begin{eqnarray}  \label{torsionequality}  \frac{ \torsnew (Y)}{  \vol^{\new}(Y) \reg_{\one}^{\new}(Y)} 
 \end{eqnarray}
is the same with $Y$ replaced by $Y'$. 
 \end{theorem}

Our goal in the following section is to refine this statement by trying to interpret this numerical value in terms of
 quantities that relate, conjecturally, to spaces of Galois representations.

\begin{proof} (of Theorem~\ref{TJL1}).

We first give the proof in the compact case, i.e., both $S, S'$ are nonempty. 
First note that there is an equality of analytic torsions
\begin{equation}  \label{analTnewequality} \analT^{\new}(Y) = \analT^{\new}(Y'), \end{equation} 
Because of the sequence \eqref{charzeroraising}, it is  sufficient to show that the spectrum of the Laplacian $\Delta_j$ acting on {\em new}
forms on $Y(K_{\Sigma})$ and $Y'(K_{\Sigma})$ coincide. But that is a direct consequence of the Jacquet--Langlands correspondence as recalled in~\S~\ref{ss:JLclassic}. 

Since we are in the compact case (that is to say, both $S, S'$ are nonempty)
 ~\eqref{analTnewequality}  $\implies$~\eqref{torsionequality}: this follows directly from the Cheeger-M{\"u}ller theorem,
in the form given in~\S~\ref{sec:CMthm}.

Otherwise, we may suppose
that $S =\emptyset$ and $S' \neq \emptyset$, so that $Y$ is noncompact but $Y'$ is compact. 
If we suppose $S' = \{{\p}, {\q}\}$,  then
the statement has already been shown in 
  Corollary~\ref{corollary610}.    For $S=\emptyset$, arbitrary $S'$, 
the proof is essentially the same as for the quoted Corollary: indeed, that corollary
follows in a straightforward way from Theorem~\ref{theorem:invtrunc}, and that deduction
can be carried through in the same fashion for arbitrary $S'$. 

\end{proof}

 \section{Torsion Jacquet--Langlands, crude form: matching volume and congruence homology}
\label{section:refined}

 We now return to the general case and show that {\em the ``volume'' factors in the Theorem~\ref{TJL1} are accounted for by congruence homology.}

Write $\Sigma = S \cup T$, and $R = S \cup V$, where $\emptyset \subseteq V \subseteq T$.
  We now prove an elementary lemma.
 Recall that $w_F$ is the number of roots of unity, and that $w_F$  divides $4$
 unless $F = \Q(\sqrt{-3})$.
\begin{lemma} The groups $H_{1,\con}(R,\Z)$ and \label{lemma:ordercong}
$H^1_{\con}(R,\Q/\Z)$ have order
$$\left( \cdot \prod_{S} (N \q + 1) \cdot  \prod_{V} (N \q - 1)\right)^{\comp}$$
up to  powers of $\ell | w_F$,
 where $\comp$ is the number of connected components of $Y(\Sigma)$.
\end{lemma}

\begin{proof} 
Recall that we have described the connected components of $Y(K_{\Sigma})$
in Remark~\ref{conncompremark}:
$$\YO= \coprod_{\AN}
Y_0(\Sigma,\a) = 
\coprod_{\AN} \Gamma_0(\Sigma,\a) \backslash \H^3,$$
Thus, 
(by definition, as in~\S~\ref{subsubsection:reducetogroup})
it suffices 
 to compute the product of the orders   of the congruence
groups $H_{1,\con}$ of $Y_0(\Sigma,\a)$ for all $\a \in A$.

By strong approximation
this amounts to computing the size of the 
maximal abelian quotients of $K^1_v$ for all places $v$ which
come from abelian quotients of $K_v$ (that is, the inverse image
of the map $K^{1,\ab}_v \rightarrow K^{\ab}_v$).
The independence of the intermediate choice of $K$ follows from
the fact that $K$ is ``$p$-convienient'' (see Definition~\ref{convenient}).
In particular, (following Lemma~\ref{lemma:congcomp}, as
well as equations~\ref{congcompute} and~\ref{congcomputetwo})
this computation is independent of the choice of connected component, and
the answer, up to the exceptional factors above, is
$$\left( \prod_{v}  \left|
\text{inverse image of $K^{\ab}_v$ in $K^{1,\ab}_{v}$}
 \right|
 \right)^{\comp}$$
  For $v \in S$, the level structure $K_v$ arises from the maximal
order in a quaternion algebra over $\OL_v$, and when $v \in V$,
$K_v$ is of  ``$\Gamma_0(v)$-type'' in $\GL_2(\OL_v)$.
If we compute with the norm one element groups $K^1_{v}$, we arrive at the answer above.
In general, we have to ensure that the image of the determinant  (which has $2$-power order)
acts trivially on these quotients (equivalently, the abelian quotients of $K^1_v$ lift to abelian
quotients of $K_v$. This follows by an explicit verification: in the quaternion algebra case,
it follows from the fact that (in the notation of~\S~\ref{section:involutions}) that 
the map $K^1_v \rightarrow l^1$ lifts to a map $K_v \rightarrow l^{\times}$, and in the $\Gamma_0(v)$-case because $K^1_v \rightarrow k^{\times}$ lifts to a map $K_v \rightarrow
(k^{\times})^2$.
\end{proof}

Let us make the following definition (cf. Definition~\ref{df:newstuff}):
 \begin{df}
 If $S$ denotes the set of finite places of $F$ which ramify in $D$,
we define the {\em new essential numerical torsion} as follows: 
\begin{eqnarray*}
\torsnewE(\Sigma) &=& \prod_{S \subset R \subset \Sigma} \frac{ |H_{1, \tors}(R, \Z)|^{(-2)^{|\Sigma \setminus R|}}}
{   |H_{1, \con}(R, \Z)|^{(-2)^{|\Sigma \setminus R|}}  }  \\
\torscong(\Sigma) &=& \prod_{S \subset R \subset \Sigma} |H_{1, \con}(R, \Z)|^{(-2)^{|\Sigma \setminus R|}}
    \end{eqnarray*} 
  \end{df}
As in the case of $h^{\new}$,  the number $h^{E,\new}$ is not, a priori, the order of any group, and again one of our concerns
  later in this Chapter will be to try to interpret it as such.

  Clearly, $\torsnewE(\Sigma) = \torsnew(\Sigma) \cdot \torscong(\Sigma)$.
We have the following refined version of Theorem~\ref{TJL1}.
\begin{theorem}  \label{TJL2}  The quantity 
 $\displaystyle{ \frac{ \torsnewE (Y)}{ \reg_{\one}^{\new}(Y)} }$
is the same with $Y$ replaced by $Y'$. 
 \end{theorem}
 
 \begin{proof}
In light of Theorem~\ref{TJL1}, it suffices to prove that the ratio
$\displaystyle{\frac{\torscong(Y)   }{\vol^{\new}(Y)}}$ is the same when $Y$ is replaced by $Y'$.
This expression can be written as an alternating product, each of whose terms looks like
$$ \frac{  |H_{1, \con}(R, \Z)|} { \vol(Y(K(R))} $$
to the power of $(-2)^{|\Sigma \setminus R|} = (-2)^{|V|}$.
Recall that $\Sigma = S \cup T$, whereas $R = S \cup V$.
By combining Borel's volume formula (Theorem~\ref{theorem:borel})  with Lemma~\ref{lemma:ordercong},
this term (without the exponent) can also be written as the $\comp$th power of:
$$\left( \frac{ \mu 
 \cdot \zeta_{F}(2) \cdot |d_F|^{3/2}}{2^{m} (4 \pi^2)^{[F:\Q] - 1}} \right)^{-1} 
 \cdot   \prod_{S}\frac{(N \q + 1)}{(N \q - 1)} 
 \prod_{T}
\frac{(N \q - 1)}{(N \q + 1)},$$
(recall that $\mu \in \Z$ depends only on $F$). Let us count to what power each term is counted in the alternating product.
\begin{enumerate}
\item The volume (related) term $\displaystyle{\frac{\nu_N (\pi^2)^{[F:\Q] - 1}}{\mu \zeta_{F}(2) |d_F|^{3/2}}}$
occurs to exponent 
$$1  - 2 \binom{d}{1} + 2^2
\binom{d}{2} - 2^3 \binom{d}{3} + \ldots = (1 - 2)^{d} = (-1)^d.$$
\item The factor $\displaystyle{\frac{(N \q + 1)}{(N \q - 1)}}$ for $\q \in S$
occurs to the same exponent, for the same reason.
\item The factor $\displaystyle{\frac{(N \q +1)}{(N \q - 1)}}$ for $\q \in \Sigma \setminus S$,
which occurs to exponent   $-1$ in the expression above and so to exponent
$-(-2)^{|V|}$ in the alternating product, but only appears when $\q \in V$, occurs to exponent
to exponent 
$$-1  + 2 \binom{d-1}{1} -  2^2 \binom{d-1}{2} + \ldots = -(1 - 2)^{d-1} = (-1)^d.$$
\end{enumerate}
Both the volume term and the term
$\displaystyle{\frac{(N \q + 1)}{(N \q - 1)}}$ for $\q \in \Sigma$ do not depend on whether we are considering $Y$
or $Y'$. Thus it suffices to show that in the product above they occur to the same exponent.
 Yet we have
just calculated that they occur to exponent $(-1)^d \comp$ and $(-1)^{d'} \comp'$ respectively,
where $d  = |\Sigma \setminus S|$ and $d' = |\Sigma \setminus S'|$. 
We have already noted (cf. \S \ref{YMdef}) that 
$\comp = \comp'$. On the other hand,
since
$|S| \equiv |S'| \mod 2$, it follows that $d \equiv d' \mod 2$, and so $(-1)^d = (-1)^{d'}$; the result follows.
\end{proof}

\section{Essential homology and the torsion quotient}\label{sec:esshom} 
\label{section:essential}

In this section we discuss notions of ``essential'' homology, that is to say,
we define variant versions of homology which excise the congruence homology.
There are different ways of doing so, adapted to different contexts.

\medskip 

We have defined previously (\S  \ref{s:congess}) a congruence quotient of homology. 
By means of the linking pairing on $H_{1,\tors} (\Sigma,\Z_p)$, we may thereby also
define a {\em congruence subgroup} of homology, namely, the subgroup 
$H_{1, \con*}$ 
orthogonal to the kernel of $H_{1,\tors}(\Sigma, \Z_p) \rightarrow  H_{1,\con}(\Sigma, \Z_p)$. 

Note that the order of $H_{1, \con}$ and $H_{1, \con*}$ need not be the same;
they differ by the $p$-part of {\em liftable} congruence homology (see~\S~\ref{LiftableCongruenceHomology}):
\begin{equation} \label{conconstarlift}  \frac {|H_{1, \con}|} { |H_{1, \con*}|} = h_{\lif}(\Sigma) \mbox{ (equality of $p$-parts.)} \end{equation}

\index{essential homology}

\medskip
\begin{definition}   \label{df:essentialhomology} 
Suppose $p > 2$ is not an orbifold prime. 

We define the {\em essential} and {\em dual-essential} homology;
$$H_{1}^{E}(\Sigma, \Z_p) = \mathrm{ker}( H_1 \rightarrow H_{1, \con}); $$ 
$$ H_{1}^{ E^*}(\Sigma, \Z_p) = H_1 / H_{1, \con*}.$$
 
 We make a similar definition with $\Z_p$ replaced by any ring in which
 $2$ and orbifold primes are invertible. 
 \end{definition} 

The first definition is quite general; the second relies, in the split case,
on the vanishing of boundary homology of the cusps (Lemma~\ref{vanishingboundaryhomology})
and thus is well-defined when $p >2$ and is not an orbifold prime; it has a certain {\em ad hoc} feel to it.  
Although, e.g. $H_1^{E^*}(\Sigma, \Z)$ is not defined for this reason, we will allow ourselves nonetheless to write
\begin{equation} \label{integralcoefficientsexplanation} | H_1^{E^*}(\Sigma, \Z)_{\tors}| \end{equation}
for the   product of orders of $ | H_1^{E^*}(\Sigma, \Z_p)_{\tors} |$  over all non-orbifold primes;
equivalently, for the order of $H_1^{E^*}(\Sigma, \Z')_{\tors}$, where $\Z'$ is obtained by inverting all orbifold primes.
This will be a convenient abbreviation since, in the rest of this chapter, we will usually work only ``up to orbifold primes.''

\medskip

 Note that if the cohomology  $H_1(\Sigma, \Z_p)$ is pure torsion, these  notions are dual, 
 i.e. there is a perfect pairing $H_1^E(\Sigma, \Z_p) \times H_1^{E^*}(\Sigma, \Z_p) \rightarrow \Q_p/\Z_p$.

 \begin{remarkable}\label{SameOrDifferent}
 {\em Even if the congruence homology is nontrivial, the essential homology may or may not differ, as an abstract group, from the actual homology; this is again related to the existence of
 liftable congruence homology (\S~\ref{LiftableCongruenceHomology}). 
 }
 \end{remarkable}

\begin{remarkable}  \label{essnewformremark}
\emph{
The level lowering maps $\Phi^{\vee}$ and $\Psi^{\vee}$ induce
maps $\Phi^{\vee}_E$ and $\Psi^{\vee}_E$ on essential homology, and similarly for dual-essential. By virtue of this remark, one can define correspondingly essential new homology, dual-essential new homology, and so on.
}
\end{remarkable}

\medskip

These notions have different utility in different contexts. In particular, 
Ihara's lemma achieves a very elegant formulation in terms of essential homology,
whereas dual statements seem better adapted to dual-essential homology: 
\medskip

\begin{theorem} \label{essentialgood}
Suppose $p  > 3$ 
and assume that the other conditions for Ihara's lemma (\S~\ref{sec:ihara}) are valid.
Then the level-lowering map on essential homology
$$\Psi: H_{1}^E(\Sigma, \Z_p) \longrightarrow H_{1}^E(\Sigma/\q, \Z_p)^2$$
is surjective. 

Under the same conditions, assume further that $H_1(\Sigma, \Q_p)^{\q-\new} = 0 $. Then  the level-raising map on dual-essential homology
$$ \Psi^{\vee}: H_{1}^{E^*} (\Sigma/\q, \Z_p)^2 \rightarrow H_{1}^{ E^*} (\Sigma, \Z_p)$$
is injective.  More generally the same result holds localized at $\mathfrak{m}$ a maximal ideal
of $\T_{\Sigma}$, assuming that $H_1(\Sigma, \Z_p)_{\mathfrak{m}} $ is torsion. 
\end{theorem}

\subsection{The proof of the first statement in Theorem~\ref{essentialgood}} 

First of all note that the number of roots in the Hilbert class field satisfies $w_H = 2$ if $F$ has a real place; otherwise, the only primes dividing $w_H$ are $2$ and $3$. 

The first statement is an immediate consequence of Ihara's lemma,  
\eqref{disjointness}, and  
\begin{lemma} \label{lemma:ihararefined2}  \label{lemma:refined} Suppose that $p > 3$. 
Under the assumptions of Ihara's Lemma, 
the image of the level raising map on cohomology:
$$\Phi: H^1(\Sigma/\q,\Q_p/\Z_p)^2 \rightarrow H^1(\Sigma,\Q_p/\Z_p)$$
has trivial intersection with the $\q$-congruence homology. Similarly, 
the kernel of the  level lowering
map 
$$\Psi: H_1(\Sigma,\Z_p) \rightarrow H_1(\Sigma/\q,\Z_p)^2$$
on homology surjects onto the congruence homology at $\q$.
\end{lemma}

\begin{proof} Let $\rf$ be an auxiliary prime (different from $\q$, the primes
dividing $p$, and not contained in $\Sigma$) with the property that $p \nmid  N(\rf) - 1$
and that $\rf$ is principal.  Such a prime
exists by the Cebotarev density theorem and the fact that $p \nmid w_H$. 
Let $\Sigmar = \Sigma \cup \{\rf\}$. The prime $\rf$ has the
property that $H^1(\Sigmar, \Q_p/\Z_p)$ does not have any ($p$-torsion) $\rf$-congruence
cohomology, since (\S~\ref{Sarithmeticcongruence}) the $\rf$-congruence homology has order dividing $N(\rf) - 1$ (as a consequence
of our assumptions on $\rf$).

 We have the following commutative diagram:
$$
\begin{diagram}
H^1(\Sigma/\q,\Q_p/\Z_p)^4 &  &  \rTo^{\Phi_{\q} \oplus \Phi_{\q}  } &  &  H^1(\Sigma,\Q_p/\Z_p)^2 \\
\dTo^{\Phi_{\rf} \oplus \Phi_{\rf}} & &  & & \dTo^{\Phi_{\rf}} \\
H^1(\Sigmar /\q,\Q_p/\Z_p)^2  &   &  \rTo^{\Phi_{\q}}  & &   H^1(\Sigmar,\Q_p/\Z_p) \\
\end{diagram}
$$
Here the map $\Phi_{\q} \oplus \Phi_{\q}$ is simply the map $\Phi_{\q}$ applied to the first and second,
and third and fourth copy of $H^1(\Sigma/\q)$, whereas the map
$\Phi^2_{\rf}$ is defined on the first and third, and second and fourth factors respectively.

We first prove that this is commutative.
First note that conjugation by
$\displaystyle{\left(\begin{matrix} \pi_{\q} & 0 \\ 0 & 1 \end{matrix} \right)}$
and
$\displaystyle{\left(\begin{matrix} \pi_{\rf} & 0 \\ 0 & 1 \end{matrix} \right)}$ commute.
It follows that the two nontrivial degeneracy maps $d_{\q}$ and $d_{\rf}$ commute
(the other degeneracys map are the obvious inclusions which certainly commute with everything.)
It follows that the composition of maps on the upper corner is:
$$(\alpha,\beta,\gamma,\delta) \mapsto (\alpha + d_{\q} \beta, \gamma  + d_{\q} \delta)
\mapsto \alpha + d_{\q} \beta + d_{\rf}  \gamma+ d_{\rf} d_{\q} \delta,$$
whereas the lower corner is
$$(\alpha,\beta,\gamma,\delta) \mapsto (\alpha + d_{\rf} \gamma, \beta + d_{\rf} \delta)
\mapsto \alpha + d_{\rf} \gamma + d_{\q} \beta+   d_{\rf} d_{\q} \delta.$$

Assume, now, that $\Phi_{\q}(\alpha,\beta)$ is congruence cohomology at $\q$; to prove the lemma
it suffices to show
that it vanishes. 
We may choose $(\alpha', \beta') $ such that $\Phi_{\q} (\alpha', \beta') 
= - [\rf] \Phi_{\q}(\alpha, \beta)$: since $\rf$ is principal, $\alpha'=-\alpha$ and $\beta'=-\beta$ will do. 
 Then
$$\Phi_{\q} \oplus \Phi_{\q} (\alpha,\beta,\alpha', \beta')
= (\Phi_{\q}(\alpha,\beta),-[\rf] \Phi_{\q}(\alpha,\beta))$$
lies in the kernel of
$\Phi_{\rf}$. By commutivity,  
$$\left( \Phi_{\rf}(\alpha, \alpha'),\Phi_{\rf}(\beta, \beta') \right)$$
lies in the kernel of $\Phi_{\q}$.

By Ihara's lemma, the kernel of $\Phi_{\q}$
consists of classes of the form $([\q]x, - x)$
where $x \in H^1_{\con}(\Sigma /\q, \Q_p/\Z_p)$.

In view of our assumptions on $\rf$,  both of the degeneracy maps
$$H^1_{\con}(\Sigma/\q, \Q_p/\Z_p) \rightarrow H^1_{\con}(\Sigma \rf/\q, \Q_p/\Z_p)$$
are isomorphisms.  In particular, there exists
$a \in H^1_{\con}(\Sigma/\q, \Q_p/\Z_p)$ such that
$$\Phi_{\rf}(\alpha, \alpha') = \Phi_{\rf}(a,0)$$
and similarly there is $b$ such that $\Phi_{\rf}(\beta,\beta') = 
\Phi_{\rf}(b,0)$.

Hence $(\alpha - a, \alpha')$ and $(\beta - b,\beta')$
lie in the kernel of
$\Phi_{\rf}$ -- that is to say,  $\alpha, \beta \in H^1_{\con}(\Sigma/\q)$. 
Then the image of $(\alpha, \beta)$ under $\Phi_{\q}$ cannot be $\q$-congruence
by~\eqref{disjointness}.

The claim about homology follows from an identical (dualized) argument.
\end{proof}

\subsection{The proof of the second statement of Theorem~\ref{essentialgood}}
Note that the second part of Theorem~\ref{essentialgood} is dual to the first statement when $H_1(\Sigma,\Q) = 0$.  Let us now show that it continues to hold when 
 $H_1(\Sigma, \C)^{\qnew} = 0$.  The argument will be the same whether we first tensor with $\T_{\Sigma,\m}$ or not, hence,
we omit $\m$ from the notation.

First of all, to show that $\Psi^{\vee}: H_1^{E^*}(\Sigma/\q,\Z_p)^2
\rightarrow H_1^{E^*}(\Sigma, \Z_p)$ is injective, it is enough to verify the corresponding statement on the torsion subgroups (since
it is clear on torsion-free quotients without any assumptions.)
 
On the other hand,  the duality on $H_{1, \tors}$ shows that 
the injectivity of $\Psi^{\vee}_{\tors}$ is equivalent
to the {\em surjectivity} 
of 
$$ \Psi: H_{1, E, \tors}(\Sigma, \Z_p) \rightarrow H_{1, E, \tors}(\Sigma/\q, \Z_p)^2.$$  But we have shown in the first part of the Theorem that 
$\Psi: H_{1, E}(\Sigma, \Z_p) \rightarrow H_{1,E}(\Sigma/\q, \Z_p)^2$
is surjective; since it is injective on torsion-free parts (this is where the assumption $H_1(\Sigma,\Q_p)^{\new} = 0$ is used), it must also be surjective on torsion.

\section{Torsion Jacquet--Langlands, refined form: spaces of newforms.}
\label{section:relations}

In this section, we outline to what extent we can interpret the \emph{a priori} rational number
$\torsnewE$ as the order of a space of newforms.  In this section,
we will see what we can do only with Ihara's lemma in the case  $S = \emptyset$, $S'= \{\p,\q\}$.  In particular, we will prove ``Theorem~A'' and ``Theorem~B'' from the introduction to the manuscript.  We do not aim for generality -- rather illustrating
some situations where one can obtain interesting results. 

\medskip 
 In the later sections
(\S~\ref{section:nonny} and~\S~\ref{sec:nonny2}) we will attempt to use the spectral sequence
of Chapter 4 to get   results in a general setting.    Although conditional, the results of those sections are valuable
in that they give a much more precise idea of what should be true. 

\medskip

Again, we continue with the same notation as set earlier in the Chapter, so that $Y, Y'$
are a Jacquet--Langlands pair.

\subsection{The case where \texorpdfstring{$H^1(Y(\Sigma), \C) = 0$}{H^1 = 0}. }

Suppose, for example, that $H^1(Y(\Sigma),\C) = 0$. Then we have an equality,
from Theorem~\ref{TJL2}, 
$$\torsnewE(Y) = \torsnewE(Y').$$
Note that
$H^1(Y(R),\C) = 0$ for all $S \subset R \subset \Sigma$, and hence there is an equality:
$$\torsnewE(Y) = \prod_{S \subset R \subset \Sigma} |H^{E}_{1}(Y(R),\Z)|^{(-2)^{\Sigma \setminus R}}.$$
To what extent does this differ from the order of the group
$H^{E}_1(\Sigma,\Z)^{\new}$ defined in Section~\ref{Sec:Newforms} and after Remark
\ref{essnewformremark}?
 \medskip
\begin{theorem}  \label{theorem:firstsetting} Suppose $S = \emptyset, S' = \{{\p}, {\q}\}$. 
 If $H_1(\Sigma,\C) = 0$, then there is an equality 
$$|H^{E^*}_1(Y(\Sigma),\Z)^{\new}| = \chi \cdot  |H^{E^*}_1(Y'(\Sigma),\Z)| .$$
for an integer $\chi  \in \Z$, and away from primes $\leq 3$. 
\label{theorem:simplecase}
\label{theorem:oldforms}
\end{theorem}
Recall (see around \eqref{integralcoefficientsexplanation}) that the order of $H_1^{E^*}(\Sigma, \Z)$
makes sense away from orbifold primes.

Note that even under this quite simple assumption, we still obtain a factor 
$\chi $ which need not be trivial.   We will see that we expect $\chi$ to be related to $K_2$ (see also
\S~\ref{section:phantomclasses} and the references to numerical computations there,
as well as \S  \ref{section:K2examples}).

\begin{proof} 
All the statements that follow are to be understood to be valid only up to the primes $2$ and $3$. 

We have $\torsnewE(Y) = \torsnewE(Y')$, and, in this case, 
$\torsnewE(Y')$ is simply the order of $H^{E^*}_1(Y'(\Sigma), \Z)$.

It remains to analyze the dual-essential new homology for $Y$.  
That group is,  by definition (see Definition~\ref{df:essentialhomology} and subsequent discussion),
the cokernel of the level-raising map:
$$H_1^{E^*}(\Sigma/p)^2 \oplus H_1^{E^*}(\Sigma/\q)^2 \rightarrow H_1^{E^*}(\Sigma).$$

Dualizing (see remark after Definition~\ref{df:essentialhomology})
 this is dual to the kernel of the right-most map in  \begin{equation} \label{chiprelimdef} H_{1}^{E}(\Sigma / \p\q)^4  \leftarrow H_1^{E}(\Sigma/\p)^2 \oplus H_{1}^{E}(\Sigma/\q)^2 \leftarrow H_{1}^{E}(\Sigma). \end{equation} 
Away from primes dividing $w_H$, the left-most map is surjective (Theorem~\ref{essentialgood}) and the alternating ratio of orders of the sizes of the groups
equals $h^{E, \new}_{\tors}(\Sigma)$.  The result follows, where $\chi$
equals the  order of the homology of the sequence at the middle term. 
\end{proof}

\subsection{The nature of $\chi$: preliminary discussion} \label{Adaggerproof} 

 Let us give a sketch of proof of Theorem $A^{\dagger}$ from the introduction, in particular, yielding a relationship between
$\chi$ from the prior Theorem \ref{theorem:simplecase} and $K_2$.  Some of the ideas that follow are further developed in the subsequent sections \S \ref{section:nonny} and \S \ref{sec:nonny2}. 

As in the statement of Theorem $A^{\dagger}$ we suppose that the class number is odd
(so that all the involved manifolds have a single connected component), that 
$\Sigma = \{\p,\q\}$  (so that we are working ``at level $1$'' on the quaternionic side),
and, as before, that $H_1(\Sigma)$ is pure torsion. 
We fix a prime $\ell > 3$, and will compute
$\ell$-parts. Set $g$ to be the $\ell$-part of $\gcd(N\p-1, N\q-1)$.  

We will show that:
\begin{equation} \label{ellpart} \mbox{ the $\ell$-part of } \# K_2(\OO_F) \mbox { divides } \chi \cdot g \end{equation} 
which implies Theorem $A^{\dagger}$.

The integer $\chi$ amounts to the order of the homology at the middle of \eqref{chiprelimdef}.  
Denote that homology group, tensored with $\Z_{\ell}$, by the letter $M$, so that 
$|M|$ is the $\ell$-part of $\chi$. 

Now consider the square, where all homology is taken with coefficients in $\Z_{\ell}$, and we impose Atkin-Lehner signs 
(i.e., all homology at level $\Sigma/\p$ should be in the $-$ eigenspace
for $w_{\q}$, all homology at level $\Sigma/\q$ should be in the $-$ eigenspace
for $w_{\p}$, and at level $\Sigma$ we impose both $-$ conditions):

{\small 
  $$\begin{diagram}
H_{1,\tors}^{E}(\Sigma / \p\q)  &\lTo& H_{1,\tors}^{E}(\Sigma/\p)^- \oplus  H_{1,\tors}^{E}(\Sigma/\q)^- &\lTo& H_{1,\tors}^{E}(\Sigma)^{--}. \\ 
\dTo && \dTo && \dTo  \\
H_{1,\tors}(\Sigma / \p\q)   &\lTo& H_{1,\tors}(\Sigma/\p)^{-} \oplus  H_{1,\tors}(\Sigma/\q)^{-} &\lTo& H_{1,\tors}(\Sigma)^{--}. \\ 
\dTo && \dTo && \dTo  \\
H_{1,\con}(\Sigma / \p\q)  &\lTo& H_{1,\con}(\Sigma/\p)^{-} \oplus  H_{1,\con}(\Sigma/\q)^{-} &\lTo& H_{1,\con}(\Sigma)^{--}. \\ 
 \end{diagram}
$$}

\medskip
Each vertical row is a short exact sequence, by definition. Also, $\Sigma = \{\p,\q\}$,
so that $\Sigma/\p\q$ is the trivial level structure, but we have continued to use $\Sigma$
for compatibility with our discussion elsewhere. 

\medskip

When we speak about ``homology'' of this square in what follows, we always
mean the homology of the horizontal rows. 
For example, the
homology at the middle term of the bottom row is $\left( k_{\p}^{\times} \oplus k_{\q}^{\times} \right) \otimes \Z_{\ell}$ (recall
that there is only one connected component).  
Now let $M$ denote the homology group at the middle of the upper row, 
and let $\tilde{M}$ denote the homology group in the center of the square;
then we have a surjection: 
\begin{equation} \label{tsv2} M \twoheadrightarrow    \ker \left( \tilde{M} \stackrel{f}{ \rightarrow} (k_{\p}^{\times} \oplus k_{\q}^{\times}) \otimes \Z_{\ell} \right). \end{equation}
 But $\tilde{M}$ arose already, in the computations of the homology of $S$-arithmetic groups
(\S \ref{sec:sss}, e.g. \eqref{SSconvenientref}); the spectral sequence of that section shows that we have an exact sequence 
$$\mbox{ a quotient of $\Z_{\ell}/g$} \rightarrow  H_2( Y(K_{\Sigma}[\frac{1}{\p\q}],  \Z_{\ell}) \twoheadrightarrow \tilde{M},$$
 
We have also proven (same proof as
Theorem \ref{theorem:K2popularversion} (i), see \eqref{placeholder}) 
 that there is a surjection
\begin{equation} \label{Fredflintstone} H_2( Y(K_{\Sigma}[\frac{1}{\p\q}]),  \Z_{\ell})  \twoheadrightarrow K_2(\OO_F[\frac{1}{\p\q}]) \otimes \Z_{\ell}. \end{equation} 
Finally, the groups $\tilde{M}, H_2(\dots), K_2(\OO_F[\frac{1}{\p\q}]$
all have maps to $(k_{\p}^{\times} \oplus k_{\q}^{\times}) \otimes \Z_{\ell}$
(the map $f$ from \eqref{tsv2}, the pull-back of $f$ to $H_2(\dots)$, and the tame symbol, respectively), and these three maps are all {\em compatible}. That implies  
that 
$$\mbox{ the $\ell$-part of } \# K_2(\OO_F) \mbox { divides } |M| \cdot g.$$
But $M$ is a direct summand of the group whose order defines $\chi$, so we are done proving
\eqref{ellpart}.

\subsection{No newforms}
 
 We continue under the assumption of $S = \emptyset, S' = \{\p, \q\}$. 
 
Now let us suppose that $H_1(\Sigma,\C)$ is non-zero, but that $H_1(\Sigma,\C)$ is  completely
accounted for by the image of the space of newforms of level $\Sigma/\q$, that is to say,
$H_1(\Sigma,\C) \simeq H_1(\Sigma/\q,\C)^2$. In this case, we see that $\reg^{\new}(Y') = 1$, and we may
also compute $\reg^{\new}(Y)$ as well, namely,
by Theorem~\ref{h2regregcompare},
\begin{equation} \label{regnewoneblah} \reg^{\new}_{\one}(Y) =  \frac{|H_{1,\lif}(\Sigma/\q;\q)|^2}{\Delta},\end{equation}
where $\Delta$ is the order of the cokernel of the map
$$\left(\begin{matrix}  (N(\q) + 1) & T_{\q} \\ T_{\q} & (N(\q) + 1) \end{matrix} \right)$$
on $H^{\tf}_1(\Sigma/\q,\Z)^2$ to itself. Hence in this case, Theorem~\ref{TJL2} manifests itself as an equality
\begin{equation} \label{wugga} \frac{ \torsnewE(Y)}{\reg^{\new}_{\one}(Y)}  = 
\torsnewE(Y'),\end{equation}
or equivalently
  \begin{eqnarray}  \label{moonoom} \torsnewE(Y) \Delta = \torsnewE(Y')  \cdot |H_{1,\lif}(\Sigma/\q;\q)|^2  \end{eqnarray} 
In particular, the factor $\Delta$ manifests itself as torsion on the non-split quaternion algebra
-- for example, if $\ell$ divides $|\Delta|$ but not $|H_{1, \lif}|$,
then $\ell$ must divide the order of a torsion $H_1$ group for $Y'$.

A rather easy consequence is   ``Theorem~B'' (Theorem~\ref{theorem:unnamedtheorem2})
from the introduction:

 \begin{theorem}   \label{theorem:secondtheoremintro} Suppose that 
 $S =\emptyset, S' = \{\p,\q\}$. Suppose that
 $l > 3$ is a prime such that:
 \begin{enumerate}
 \item $H_1(Y(\Sigma/\p),\Z_l) = 0$,
 \item $H_1(Y(\Sigma/\q),\Z_l)$ and  $H_1(Y(\Sigma),\Z_l)$ is torsion free,
 \item $H_1(Y'(\Sigma),\C) = 0$.
 \end{enumerate}
Then $l$ can divide $H_1(Y'(\Sigma),\Z)$ if and only if $l$ divides $\frac{\Delta}{(N(\p)-1)^{\# Y}}$, where $\# Y$ is the number of connected components.   
 \end{theorem}
 Note that this notation is not precisely compatible with that used in Theorem~B, 
 but the content is the same (taking account that Theorem~B is stated
 in the context of odd class number, so that $\# Y = 1$). 
 
 \begin{proof}
 The assumptions imply that there is no mod-$\ell$ 
 torsion at level $\Sigma,\Sigma/\p,\Sigma/\q$ or $\Sigma/\p\q$. 
  Therefore $$ h^{\new}_{\tors}(Y) =1.$$ 

Next we must account for congruence homology. 
Assumption (1) implies that there is no congruence cohomology of characteristic $l$ at level $\Sigma/\p$,
and thus $\ell$ is prime to $H_{1, \con}(\Sigma/\p)$ and $H_{1,\con}(\Sigma/\p\q,\Z)$.
That means that the $l$-part of $H_{1, \con}(\Sigma/\q)$ and $H_{1, \con}(\Sigma)$
arises entirely from $\p$-congruence classes; in particular, they have order equal to (the $l$-part of)  $(N(\p) - 1)^{\# Y}$,
where $\# Y$ is the number of components of $Y$. 
 From that we deduce that
 $$h^{\new}_{\tors, E}(Y) = (N(\p)-1)^{\# Y} \mbox{ (equality of $\ell$-parts)}.$$
 
 As for the quantity $H_{1, \lif}(\Sigma/\q;\q)$ which shows up above -- 
 our assumption (2) implies that it equals simply $(N \p -1)^{\# Y}$. 
 So ~\eqref{moonoom} shows that the size of $H_1^E(Y'(\Sigma), \Z_{\ell})$ 
 (or $H_1^{E^*}$; they both have the same order here)
equals   the $\ell$-part of  $$\frac{\Delta}{(N \p-1)^{\# Y}}.$$
\end{proof}
 
 Now we pass to a more difficult case, to illustrate the type of results that we can still obtain about newforms:

\begin{theorem}  \label{theorem:oldtwo} Suppose $S = \emptyset, S' = \{\p,\q\}$. Suppose that $\ell > 3$ is a prime such that:
\begin{itemize}
\item[(i)]  The natural map $H_1(\Sigma,\Q_{\ell}) \rightarrow H_1(\Sigma/\q,\Q_{\ell})^2$
is an isomorphism, and $H_1(\Sigma/\p, \Q_{\ell}) = 0$; 

\item[(ii)] There are no mod $\ell$ congruences between torsion forms at level $\Sigma/\p \q$
and characteristic zero forms at level $\Sigma/\q$, i.e. 
there is no maximal  ideal for the Hecke algebra $\T_{\Sigma/\q} \otimes \Z_{\ell}$
that has support both in $H_{1, \tors}(\Sigma/\p \q)$ and $H_{1, \tf}(\Sigma/\q)$.

 \end{itemize}

Then, up to $\ell$-units, 
$$|H^{E^*}_1(Y(\Sigma),\Z_{\ell})^{\new}| = |H^{E^*}_1(Y'(\Sigma),\Z_{\ell})| \cdot \chi .$$
for an integer $\chi \in \Z$.
\label{theorem:newforms}
\end{theorem}

Again, the integer $\chi$ is related to $K_2$, as explained in the next section. 

\begin{proof} 
In what follows equalities should be regarded as equalities of $\ell$-parts. 

The proof is a series of diagram chases, which are discussed in a more general context in the next section. 
\medskip
Examine the diagram (recall the definition of $H_{1, \con*}$ from~\S~\ref{sec:esshom})
{\small 
  $$\begin{diagram}
H_{1,\con*}(\Sigma / \p\q)^4  &\lTo& H_{1,\con*}(\Sigma/\p)^2 \oplus  H_{1,\con*}(\Sigma/\q)^2 &\lTo& H_{1,\con*}(\Sigma). \\ 
\dTo && \dTo && \dTo  \\
H_{1,\tors}(\Sigma / \p\q)^4  &\lTo& H_{1,\tors}(\Sigma/\p)^2 \oplus  H_{1,\tors}(\Sigma/\q)^2 &\lTo& H_{1,\tors}(\Sigma). \\ 
\dTo && \dTo && \dTo  \\
H_{1,\tors}^{E^*}(\Sigma / \p\q)^4  &\lTo& H_{1,\tors}^{E^*}(\Sigma/\p)^2 \oplus  H_{1,\tors}^{E^*}(\Sigma/\q)^2 &\lTo& H_{1,\tors}^{E*}(\Sigma). \\ 
 \end{diagram}
$$}
  Inspecting this diagram,  and using \eqref{conconstarlift}, 
 we deduce the equality
\begin{equation} \label{Alph} \torsnewE(Y)  =  \frac{ H_1^{E^*}(\Sigma, \Z)^{\tors} \left( \dots \Sigma/\p\q \cdots \right)^4}{\left( \cdots \Sigma/\p \cdots \right)^2 \left( \dots \Sigma/ \q \cdots \right)^2} \cdot  
 | H_{1,\lif}(\Sigma,\Z)|^{-1} \cdot | H_{1,\lif}(\Sigma/\q,\Z)|^{2} .\end{equation} 
Now consider
{\small 
  $$\begin{diagram}
H_{1,\tors}^{E^*}(\Sigma / \p\q)^4  &\rTo& H_{1,\tors}^{E^*}(\Sigma/\p)^2 \oplus  H_{1,\tors}^{E^*}(\Sigma/\q)^2 &\rTo& H_{1,\tors}^{E^*}(\Sigma). \\ 
\dTo && \dTo && \dTo  \\
H_{1}^{E^*}(\Sigma / \p\q)^4  &\rTo^{A}& H_1^{E^*}(\Sigma/\p)^2 \oplus H_{1}^{E^*}(\Sigma/\q)^2 &\rTo& H_{1}^{E^*}(\Sigma).\\
\dTo && \dTo && \dTo   \\
H_{1,\tf}(\Sigma / \p\q)^4  &\rTo& H_{1,\tf}(\Sigma/\p)^2 \oplus  H_{1,\tf}(\Sigma/\q)^2 &\rTo^{B}& H_{1,\tf}(\Sigma). \\ 
 \end{diagram}
$$ } 
Taking into account Lemma~\ref{regulator-compare-1} to analyze the cokernel of $B$,
and noting that  Theorem~\ref{essentialgood} and assumption (ii) of the theorem statement
prove that  $A$ is injective away from primes dividing $w_H$,  
 the diagram shows that \begin{equation} \label{Balph}  |H^{E^*}_1(Y(\Sigma),\Z)^{\new}|   = \chi \cdot \frac{\Delta}{| H_{1,\lif}(\Sigma/\q;\q)|} \cdot
\frac{   H_1^{E^*}(\Sigma, \Z)^{\tors}  \left( \dots  \Sigma/\p\q \cdots \right)^4}{\left( \cdots \Sigma/\p \cdots \right)^2 \left( \dots \Sigma/\p\q \cdots \right)^2}, \end{equation} 
where $\chi \in \Z$.    As before, $\Delta$ is the determinant of $
( T_{\q}^2 - (1 + N(\q))^2 )$ on $  H_1(\Sigma/\q, \C)$.

Combining~\eqref{Alph} and~\eqref{Balph},  and using the fact (\eqref{regnewoneblah}) that
$\reg^{\new}_{\one}(Y) = h_{\lif}(\Sigma;\q)^2 \Delta^{-1}$, 
\begin{eqnarray} |H^{E^*}_1(Y(\Sigma),\Z)^{\new}|  \nonumber  &=&    \chi \cdot \torsnewE(Y) \cdot \Delta  
\cdot \frac{h_{\lif}(\Sigma)h_{\lif}(\Sigma/\q;\q)^{-1}}{ h_{\lif}(\Sigma/\q)^2} \\  \label{Andr}
&=& \chi \cdot  \frac{ \torsnewE(Y)}{\reg^{\new}_{\one}(Y)}  \cdot \frac{h_{\lif}(\Sigma) h_{\lif}(\Sigma/\q;\q)}{h_{\lif}(\Sigma/\q)^2)}.\end{eqnarray}

We now need to analyze liftable congruence homology.  Consider the  commutative diagram{\small   $$\begin{diagram}
 H_1(\Sigma/\q,\Z)_{\tors}^2 & \lTo^{\Psi_{\tors}} &  H_{1}(\Sigma, \Z)_{\tors} &\lTo& \ker(\Psi_{\tors})   \\
 \dTo & & \dTo & & \dTo \\
 H_1(\Sigma/\q,\Z)_{\con}^2 & \lTo^{id \oplus [\q]} &  H_{1}(\Sigma, \Z)_{ \con} &\lTo & H_1(\Sigma, \Z)_{\q-\con} \\
  \dOnto & & \dOnto \\
H_{1,\lif}(\Sigma/\q)^2 & \lTo & H_{1, \lif}(\Sigma) && \\
 \end{diagram}
$$} 
 By a diagram chase using Lemma 
~\ref{lemma:ihararefined2} -- note that $\ker(\Psi_{\tors}) = \ker(\Psi)$, by assumption --  the induced map $\iota: H_{1,\lif}(\Sigma) \rightarrow H_{1, \lif}(\Sigma/\q)^2$ is injective. 
Write $L$ for  image of $H_{1}(\Sigma/\q, \Z)_{\tors}$
  inside $H_1(\Sigma/\q, \Z)_{\con}$. Then one easily
  sees that the image of  $\iota$ has size equal to the index of $L \cap [\q] L$
  in $H_1(\Sigma/q, \Z)$. So, 
  $$ h_{1,\lif}(\Sigma) = \frac{  |H_1(\Sigma/\q, \Z)_{\con}|}{| L \cap [\q] L|}.$$
On the other hand,  by definition, 
$$ h_{\lif}(\Sigma/\q;\q)  =
\frac{  |H_1(\Sigma/\q, \Z)_{\con}|}{| L  +  [\q] L|} \mbox{ and }
h_{\lif}(\Sigma/\q) =\frac{  |H_1(\Sigma/\q, \Z)_{\con}|}{|L|}.$$
Putting these together, we see that
$ h_{\lif}(\Sigma)  \cdot  h_{\lif}(\Sigma/q;\q)  =   h_{\lif}(\Sigma/\q)^2$;
and, returning to~\eqref{Andr} and recalling \eqref{wugga}, we deduce
\begin{eqnarray*} |H^{E^*}_1(Y(\Sigma),\Z)^{\new}|
&=& \chi \cdot   \torsnewE(Y')  \\ &=& \chi \cdot | H^{E*}_1(Y'(\Sigma),\Z)|.\end{eqnarray*}

Note that the change in essential homology
exactly accounts for the change in regulator, when accounted for properly.
\end{proof}

 \subsection{An example with level-lowering and newforms}

As our last special case, we continue to suppose that $S' = \{\p,\q\}$ and $S= \emptyset$. 

Now suppose that $H_1(\Sigma,\C) \ne 0$, but that we have
$H_1(\Sigma/\q,\C) =0$ for
all $\q \in \Sigma$. Recall that Theorem~\ref{TJL2} in this case says that
$$\frac{\torsnewE(Y)}{\torsnewE(Y')}  = \frac{\reg_{\one}(Y)}{\reg_{\one}(Y')}.$$
We deduce from this the following: Suppose that $H_1(Y(\Sigma),\Z)$ and $H_1(Y'(\Sigma),\Z)$
are $p$-torsion free (or, more generally, that their $p$-torsion subgroups have the same order). 
Let us also assume that $p$ is co-prime to the order of the congruence homology.
Then
\begin{quote} \label{mungle}  $p$ divides $\reg_{\one}(Y)/\reg_{\one}(Y')$ 
 only if there exists a level lowering prime $\m$ of characteristic $p$. \end{quote}
Indeed, if $p$ is to divide $\reg_{\one}(Y)/\reg_{\one}(Y')$, there must be $p$-torsion at level $\Sigma/\p$ or $ \Sigma/\q$; 
Ihara's lemma shows that this must divide a level lowering prime, since there is no $p$-torsion at level $\Sigma$. 
This provides theoretical evidence towards some of our musings in~\S~\ref{reglowlevel}, in particular
the principle ~\ref{PP2}.

More generally, one can prove, similarly to the prior results, that: 
\begin{theorem}  \label{theorem:thirdsetting} Suppose $S=\emptyset, |S'|=2$. 
Suppose that $H_1(\Sigma/\q,\C) = 0$
for all $\q \in \Sigma$. 
Then there is an equality
\begin{equation}\label{mimimi} \frac{|H^{E^*}_1(Y(\Sigma),\Z)^{\new,\tors}|}{|H^{E^*}_1(Y'(\Sigma),\Z)^{\tors}|} \cdot \frac{1}{|\LL(\Sigma)|}
=  \chi \cdot  \frac{\reg_{\one}(Y)}{\reg_{\one}(Y')} ,\end{equation} 
where $p$ divides $\LL(\Sigma)$ if and only if there exists a prime $\m$ of $\T$ of characteristic dividing $p$
 which occurs in the the support of an eigenform of characteristic zero for $Y(\Sigma)$, but also occurs in
 the support of the integral cohomology at some strictly lower level.
 Moreover, $\chi \in \Z$.
\end{theorem}

 Our expectation
is that the two terms on the left-hand side of~\eqref{mimimi}  match termwise with
those on the right-hand side, that is to say:
\begin{enumerate}
\item The spaces $H^{E^*}_1(\Sigma,\Z)^{\new}$ are the same for $Y$ and $Y'$, possibly up to factors arising from
$K_2$-classes which correspond to $\chi$,
\item The ratio of regulators corresponds exactly to level lowering primes $\m$, which are accounted for by the
factor $\LL(\Sigma)$,
\end{enumerate}
and the theorem above shows that the ``product'' of these two statements holds, which can be considered as
a consistency check.

\section{Newforms and the torsion quotient (I)}
\label{section:nonny}

 Theorems~\ref{theorem:simplecase}
and~\ref{theorem:newforms} of the previous section were restricted to
the case when $S = \emptyset$ and $S' =  \{\p,\q\}$.
The goal of the following few sections is to generalize these results
to the general case, where $\Sigma = S \cup T = S' \cup T'$ and the cardinality of
$S$ and $S'$
 can be arbitrary. One case is (in most respects) actually simpler, namely,
when $|T| = |T'|=1$ (so $d=d'=1$). In that case, the ``partial Euler characterstic'' $\chi$ is trivial
 --- provided that  the congruence subgroup property holds for the appropriate
$S$-arithmetic groups (unfortunately the CSP is not known in these cases).
Returning to the general case,
the key to  such a generalization lies is understanding the exactness (or failure to be exact)
of the corresponding level raising and level lowering sequences.
The one tool we have to understand this problem beyond Ihara's Lemma is the spectral
sequence of Corollary~\ref{SarithmeticSS2}, which we recall below.
To a certain extent, however, our understanding of this spectral sequence
is illusory. The sequence abuts to the cohomology of  a $T$-arithmetic group,
which, beyond $H^1$ (where one has access --- in the split case --- to the congruence subgroup property)
is itself difficult to understand. On the other hand, our analysis does give some understanding
of the factors which may arise when comparing the new (dual-) essential cohomology of $Y$ and $Y'$.

Let us begin by recalling the statement of Corollary~\ref{SarithmeticSS2}: 
 
 \medskip
 
 Notation as in that corollary; in particular, $\cF$ is a certain flat local system
 on $Y(K[1/T])$, and we write $\Sigma = S \cup T$,  and $R = \Z[\frac{1}{ w_F^{(2)}}]$. 
 Then there exists an spectral sequence 
 abutting to $H_*(Y(K[1/T], \cF)$, where
$$E^{1}_{p,q} = \bigoplus_{V \subset T, |V| = p} H_q(S \cup V, \testring)^{2^{d - p}}.$$
Up to signs, the differential 
$$d_1: H_q(S \cup V \cup \{{\q}\},R)  \rightarrow H_q(S \cup V,R)^{2},$$
is given by the two degeneracy maps.
 In other words the $E^1$ page looks like this:

 \begin{diagram}
H_3(S,R)^{2^d}  & \lTo & \bigoplus H_3(S \cup \{\p\},\testring)^{2^{d-1}} & \lTo & 
\bigoplus  H_3(S \cup \{\p,\q\},\testring)^{2^{d-2}}  \\
H_2(S,R)^{2^d}  & \lTo & \bigoplus H_2(S \cup \{\p\},\testring)^{2^{d-1}} & \lTo & 
\bigoplus  H_2(S \cup \{\p,\q\},\testring)^{2^{d-2}}  \\
 H_1(S,\testring)^{2^d}  & \lTo & \bigoplus  H_1(S \cup \p,\testring)^{2^{d-1}} & \lTo & 
\bigoplus H_1(S \cup \{\p,\q\} ,\testring)^{2^{d-2}}\\
\testring \qquad & \longleftarrow &  0 \qquad  & \longleftarrow &  \qquad 0  \\
\end{diagram}

We will denote the  $H_1$ row of this by $H_{1, \bullet}$, the
$H_2$ row by $H_{2, \bullet}$; we denote by $H_1^{\bullet}$ and $H_2^{\bullet}$
the same complexes but with arrows reversed, i.e., the
level-{\em raising} complexes that were previously discussed 
in~\S~\ref{sec:LRC}.  In what follows, we will also denote by $\mathcal{F}_p$
the local system $\mathcal{F} \otimes_{\Z} \Z_p$.

\medskip

We wish to discuss what happens when we complete this at a non-Eisenstein prime, more precisely: 
{\em In this section, the Hecke algebra $\T$ will denote the abstract Hecke algebra omitting primes
$T_{\p}$  dividing $\Sigma$. } It acts on all terms (and all pages)  of the spectral sequence
as well as on its abutment.

\medskip

For the rest of~\S~\ref{section:nonny}, we assume that $\m \subset \T$ is a maximal non-Eisenstein\footnote{ Here by Eisenstein we mean in the broadest possible sense,
namely, of type ``D0'' in the classification~\ref{eisdef}. To remind the reader, this is equivalent to 
the (conjectural) Galois representation $\rhobar_{\m}: G_F \rightarrow \GL_2(\Fbar_p)$ being reducible.} prime with finite
residue field, and will study the localization of the spectral sequence at $\m$. 
We return to the case of an Eisenstein prime in the subsequent section.

\begin{lemma} The only non-zero terms in the spectral sequence after localizing at $\m$ occur in 
the first and second row. 
Suppose that $H_{1,!}(\Sigma,\Q) =0$  --- here $H_{1,!}$ is the image of homology in Borel--Moore homology, and so is just $H_1$ in the compact case --- then the only non-zero terms in the spectral
sequence after localizing at $\m$ occur in the first row.
\end{lemma}

\begin{proof}  For the first claim, it suffices to note that $H_3(\Sigma,\Z)$ and $H_0(\Sigma,\Z)$ are always Eisenstein.
For the second claim, this follows automatically in the non-split case. In the split case, it suffices to note that the
action of $\T(g)$ on the boundary term is
--- up to the permutation of components induced by
$\det(g) \in \AN = F^{\times} \backslash \Atimes \slash
\det(K) {\A^{\times 2}}$ ---  either via the degree (in degrees $0$ or $2$) or though
a sum of CM characters in degree $1$ (see Lemma~\ref{heckeactiononboundary}). 
Since both of these actions are Eisenstein (of type D0) the result follows.
\end{proof}

\subsection{A guiding principle} \label{section:hecketrivial}
\label{section:repeatingmyself}

A guiding principle when working with the spectral sequence and
with $S$-arithmetic cohomology is the following conjecture,
which we stated previously (Conjecture~\ref{conj:Eisenstein}).

\begin{conj} \label{conj:hecketrivial} Suppose that $|T| = d$. Then
$H_i(K[1/T],\cF)$ is cyclotomic Eisenstein  for all $i \le d$.
\end{conj}

When $i = 1$, this conjecture follows whenever one knows the congruence
subgroup property for $Y(K[1/T])$. When $i = 2$, this is related
to whether the non-trivial classes in $H_2(K[1/T],\Z)$ \emph{all}
come from $K_2$.
When $i = 3$, the group $H_3(K[1/T],\Z)$ contains characteristic zero
classes which are associated to $K_3(\OL_F)$ --- these are cyclotomic Eisenstein.
The relevance of this conjecture is as follows. It implies that, localizing the
spectral sequence at a \emph{non}-Eisenstein ideal $\m$, that the spectral
sequence converges to zero in the quadrant $p+q \le d$. This puts strong
conditions on the corresponding differentials and often forces exactness.
The conjecture when $i = 1$ exactly corresponds to the claim that the
cokernel of  the level lowering map $H_1(\Sigma,\Z) \rightarrow H_1(\Sigma/\p,\Z)^2$
is congruence Eisenstein, which follows from Ihara's lemma.
When $i \ge 2$, however, one should not expect that
the Eisenstein classes should be easily comprehended, given the relation
to $K_2(\OL_F)$ discussed in Theorem~\ref{theorem:K2popularversion}.

\subsection{Some remarks on complexes}

Fix a positive integer $d$, and
let $\Cat_d = \Cat$ denote the category of chain complexes 
$$C_{\bullet}:=C_0 \leftarrow C_1 \leftarrow \ldots  \leftarrow C_d$$
  such that
   \begin{enumerate}
\item $C_{\bullet}$ denote a complex of finitely generated $\T_{\m}$-modules.
\item All terms of $C_{\bullet}$ are trivial outside the range $[0,d]$. 
\item 
$C_{\bullet} \otimes_{\Z} \Q$ is exact everywhere except possibly at $d$.
\end{enumerate}
The morphisms are maps of chain complexes. 

We assume unless otherwise stated that $C_{\bullet} \in \Cat$.
We may also define a category $\Cat^{\vee}$, consisting of chain complexes
$C^{\bullet}$ satisfying the same conditions but with the maps in the opposite direction.

\medskip

The point of this definition is to model the complexes $(H_{1,\bullet})_{\m}$ and
$(H_{2,\bullet})_{\m}$ as well as $(H^{1,\bullet})_{\m}$ and
$(H^{2,\bullet})_{\m}$, and to establish some basic lemmas relating their
cohomology.  In particular, $(H_{1,\bullet})_{\m}$ and
$(H_{2,\bullet})_{\m}$ are elements of $\Cat$ and $(H^{1,\bullet})_{\m}$ and
$(H^{2,\bullet})_{\m}$ are elements of $\Cat^{\vee}$,  wherer $d$ is taken to be the number of
prime divisors of $T$ = $\Sigma \setminus S$.

Every result of this section is a trivial diagram chase.

\begin{df} Let $\Omega_{d}(C_{\bullet})$, respectively, $\Omega^{d}(C^{\bullet})$  denote the following sums considered
as an element in $K_0$ of the category of $\T_{\m}$-modules of finite length:
$$\Omega_{d}(C_{\bullet}) := \sum_{i=0}^{d-1} (-1)^i [H_i(C_{\bullet})].$$
$$\Omega^{d}(C^{\bullet}) := \sum_{i=0}^{d-1} (-1)^i [H_i(C^{\bullet})].$$
\end{df}
Note that the $i = d$ term is \emph{not} included in the sum.
This definition makes sense, because $H_i(C_{\bullet})$ is finite for
$i \ne d$, by definition.
The following lemma is a diagram chase:
\begin{lemma}  \label{lemma:orderofgroup2} 
  Suppose that $C^{\bullet} \otimes \Q$ is concentrated in degree $d$. 
  Then
    $$[H_d(C^{\bullet})_{\tors}]  =
 (-1)^d \sum_{i=0}^{d} (-1)^i [C^i_{\tors}] - (-1)^d  \Omega^{d}(C^{\bullet}) .$$
  $$=
 (-1)^d \sum_{i=0}^{d} (-1)^i [C^i_{\tors}] - (-1)^d  \sum_{i=0}^{d-1} (-1)^i [H_i(C^{\bullet})] .$$
\end{lemma}
Of course, $C^i_{\tors} = C^i$ unless $i = d$.
This lemma is the reason for the definition of $\Omega$;
it shows that the cokernel of
the final term of sequence is equal to the alternating product of the remaining terms 
--- up to a correction  given by $\Omega$.

\medskip

Given a complex $C_{\bullet}$ such that $C_{\bullet} \otimes \Q = 0$, we may define
the Pontryagin dual complex $(C_{\bullet})^{\vee}:= \Hom(C_{\bullet},\Q/\Z)[d] \in
\Cat^{\vee}$.

\begin{lemma}  \label{lemma:orderofgroupduality} If $C_{\bullet} \otimes \Q = 0$, 
there is an equality
$\Omega_d(C_{\bullet})^{\vee} = \Omega^{d}((C_{\bullet})^{\vee})$.
\end{lemma}

\begin{lemma} Suppose that   \label{lemma:short} $0 \rightarrow A_{\bullet}
\rightarrow B_{\bullet} \rightarrow C_{\bullet} \rightarrow 0$ is an exact sequence
in $\Cat$. 
Then
$$
\begin{aligned} \Omega_d(B_{\bullet}) 
= \ & \  \Omega_d(A_{\bullet}) + \Omega_d(C_{\bullet}) +
(-1)^{d} [\coker(H_d(B_{\bullet}) \rightarrow H_{d}(C_{\bullet}))] \\
= \ & \  \Omega_d(A_{\bullet}) + \Omega_d(C_{\bullet}) +
(-1)^{d} [\im(H_d(C_{\bullet}) \rightarrow H_{d-1}(A_{\bullet}))]\\
= \ & \  \Omega_d(A_{\bullet}) + \Omega_d(C_{\bullet}) +
(-1)^{d} [\ker(H_{d-1}(A_{\bullet}) \rightarrow H_{d-1}(B_{\bullet}))].  \end{aligned}
$$
If $0 \rightarrow A^{\bullet}
\rightarrow B^{\bullet} \rightarrow C^{\bullet} \rightarrow 0$ is an exact sequence in
$\Cat^{\vee}$, then
$$
\begin{aligned} \Omega^d(B^{\bullet}) 
= \ & \  \Omega^d(A^{\bullet}) + \Omega^d(C^{\bullet}) +
(-1)^{d} [\ker(H_d(A^{\bullet}) \rightarrow H_{d}(B^{\bullet}))] \\
= \ & \  \Omega^d(A^{\bullet}) + \Omega^d(C^{\bullet}) +
(-1)^{d} [\im(H_{d-1}(C_{\bullet}) \rightarrow H_{d}(A^{\bullet}))]\\
= \ & \  \Omega^d(A^{\bullet}) + \Omega^d(C^{\bullet}) +
(-1)^{d} [\coker(H_{d-1}(B_{\bullet}) \rightarrow H_{d-1}(C^{\bullet}))]. \end{aligned}
$$
\end{lemma}
Note that the last term either injects into $H_{d-1}(A_{\bullet})$ or is a quotent
of $H_{d-1}(C^{\bullet})$, and so it is finite.
\begin{proof} This is clear from the long exact sequence in homology.
\end{proof}

Given a complex $C_{\bullet}$, we may form the sub-complex $C^{\tors}_{\bullet}$
and quotient complex $C^{\tf}_{\bullet}$ of torsion free quotients, where
clearly
$C^{\tors}_{\bullet}$, $C^{\tf}_{\bullet} \in \Cat$. 
There is a short exact sequence of complexes
$$0 \rightarrow C^{\tors}_{\bullet} \rightarrow C_{\bullet} \rightarrow C^{\tf}_{\bullet}
\rightarrow 0.$$

\begin{lemma} \label{lemma:shortflat} Suppose that $C_{\bullet} \otimes \Q$ is zero
away from degree $d$.
Then
$$\Omega_d(C_{\bullet}) = \Omega_d(C^{\tors}_{\bullet}) 
+ (-1)^d [\im(C^{\tf}_d \rightarrow H_{d-1}(C^{\tors})].$$
Suppose that $C^{\bullet} \otimes \Q$ is zero away from degree $d$.
Then
$$\Omega^d(C^{\bullet}) = \Omega^d(C^{\tors,\bullet}).$$
\end{lemma}

\begin{proof} Because $C^{\tf}_{\bullet}$ is concentrated in degree $d$
in this case, $\Omega_d(C^{\tf}_{\bullet}) = 0$. Since
$C^{\tf}_{\bullet}$ is concentrated in degree $d$, we have
$H_d(C^{\tf}_{\bullet}) = C^{\tf}_d$. Thus we obtain the result by 
Lemma~\ref{lemma:short}. For the second case, the argument is the same, except
now the error term is a quotient of $H_{d-1}(C^{\tf,\bullet}) = 0$.
\end{proof}

\medskip

Suppose that $C_{\bullet}$ is torsion free over $\Z$. \ Then we may define $(C_{\bullet})^* = \Hom(C_{\bullet},\Z)[d]$.

\begin{lemma} \label{lemma:tallflat} Suppose that $C_{\bullet} \otimes \Q$ is exact, and that
$C_{\bullet}$ is flat over $\Z$. Then
$$\Omega_d(C_{\bullet}) + (-1)^d [H_d((C_{\bullet})^*)] + \Omega^d((C_{\bullet})^*) = 0.$$
\end{lemma}

\begin{proof} This is the torsion shifting in the proof of Poincare duality.
More generally, we have isomorphisms
$$H_{i}(C_{\bullet}) \simeq  H_{i+1}((C_{\bullet})^*),$$
and so
$$\sum_{i=0}^{k} (-1)^i [H_{i}(C_{\bullet})]
= (-1) \sum_{i=0}^{k+1} (-1)^i [H_{i}((C_{\bullet})^{*})].$$
\end{proof}

\subsection{Relating newforms to the spectral sequence: The case of finite homology.}
\label{section:caseoffinitehomology}
Let us assume in this section that $H_{1,!}(\Sigma,\Q)_{\m} = 0$. It follows that, completed at $\m$,
the only non-zero terms in the spectral sequence occur in the first row.

Recall (\S~\ref{sec:LRC}) that the space of newforms is defined to be the cokernel of the level
raising complex $H^{E,\bullet}_1$.
Localized at  $\m$, this is the same
as the level raising complex  $(H^{\bullet}_1)_{\m}$.

 The complexes $H_{1, \bullet}$ and $H_1^{\bullet}$ are related (via duality)
by the relation:
$$(H^{\bullet}_1)^{\vee}_{\m} = (H_{1,\bullet})_{\m}.$$
Since there is only one non-zero row in the spectral sequence, it degenerates immediately
on page $2$.
Recall that $H_1(\Sigma,\Z)^{\new}_{\m}$ is defined to be the cokernel of the final map in the level
raising sequence. Since both these sequences consist entirely of finite terms, we deduce that,
as an element in $K_0$ of the category of $\T_{\m}$-modules of finite length, that:
$[H_1(\Sigma,\Z)^{\new}_{\m}]$
is equal to
$$ \sum_{V \subset T}  (-2)^{|\Sigma \setminus S \cup V|}  [H_1(S \cup V,\Z)_{\m}]
 + (-1)^{|T|-1} \sum_{i=0}^{|T| - 1} (-1)^i [H_i((H_{1,\bullet})_{\m})].$$
 This formula may be read as follows: the naive order of the space of newforms is the alternating product of the
 chain complex of which it is the cokernel; the failure of this formula is determined by the cohomology of this chain complex.
 Since the chain complex consists of a sequence of finite groups, this cohomology can also be computed in terms of
 the Pontryagin dual sequence, which is exactly the level lowering sequence that occurs in the first row of the spectral
 sequence. 
 If we examine the last term in this identity, we note, since the spectral sequence has now degenerated, that this
 can be written in terms of the limit. Explicitly, we deduce the following:
 \begin{lemma} In the $K_0$-group of finite length $\T$-modules,
 $[H_1(\Sigma,\Z)^{\new}_{\m}] $ is equal to
$$ \sum_{V \subset T}  (-2)^{|\Sigma \setminus S \cup V|}  [H_1(S \cup V,\Z)_{\m}]
 + (-1)^{|T|} \sum_{i=1}^{|T|} (-1)^i [H_i(K[1/T],\cF)_{\m}].$$
 \end{lemma}
 Note that the last term is shifted from the previous sum
 because  the boundary term $H_n(K[1/T],\cF)$ is related to $E^{\infty}_{1,n-1}$.
 We see in this formula that if $H_i(K[1/T],\cF)$ is Eisenstein for $i$ in the range
 $\{1,\ldots,d\}$ (with $d =|T|$), then the space of newforms in
 the finite case is exactly equal to the alternating products of the orders of the homology groups in lower degree.

\subsection{Relating newforms to the spectral sequence: newforms without oldforms}
\label{section:caseoffinitehomologynewforms}
We extend the analysis of the previous section to the following case:
\begin{enumerate}
\item $H_{1,!}(\Sigma,\Q)_{\m} \ne 0$.
\item $H_{1,!}(\Sigma/\q,\Q)_{\m} = 0$ for all $\q \in \Sigma$.
\end{enumerate}
As in the last section, we complete at a non-Eisenstein ideal $\m$. For ease of notation,
we omit $\m$ as a subscript.
We have short exact sequences as follows:
$$0 \rightarrow  H^{\tors}_{1,\bullet} \rightarrow H_{1,\bullet} \rightarrow H^{\tf}_{1,\bullet} 
\rightarrow 0,$$
$$0 \rightarrow  H^{\tors,\bullet}_{1} \rightarrow H^{\bullet}_{1} \rightarrow H^{\tf,\bullet}_{1} 
\rightarrow 0.$$
By Lemma~\ref{lemma:shortflat}, we deduce that
$\Omega^d(H^{\bullet}_1) = \Omega^d(H^{\tors,\bullet}_1)$.
On the other hand, by Lemma~\ref{lemma:shortflat}, we also deduce
that 
$$ \Omega_d(H^{\tors}_{1,\bullet})
 = \Omega_d(H_{1,\bullet}) + (-1)^d [\im(H_1(\Sigma,\Z)^{\tf}
 \rightarrow H_{d-1}(H^{\tors}_{1,\bullet}))].$$
 Now $(H^{\tors,\bullet}_1)^{\vee} = H^{\tors}_{1,\bullet}$, and so
 $$\Omega^d(H^{\bullet}_1)^{\vee} = \Omega_d(H_{1,\bullet}) + (-1)^d [\im(H_1(\Sigma,\Z)^{\tf}
 \rightarrow H_{d-1}(H^{\tors}_{1,\bullet}))].$$

 \begin{lemma} \label{lemma:support} The support of $[\im(H_1(\Sigma,\Z)^{\tf}
 \rightarrow H_{d-1}(H^{\tors}_{1,\bullet}))]  $  as a $\T_{\m}$-module is nontrivial
 if and only if the following two conditions both hold:
\begin{enumerate}
\item There exists a non-zero characteristic zero form in $H_1(\Sigma,\Z)_{\m}$.
\item There exists a non-zero torsion class in $H_1(\Sigma/\p,\Z)_{\m}$ for some $\p$.
\end{enumerate}
\end{lemma}

\begin{proof}  There is an exact diagram as follows:
$$
\begin{diagram}
H_1(\Sigma,\Z)^{\tors} & \rTo & \bigoplus H_1(\Sigma/\p,\Z)^{2} & \rTo & \ldots \\
\dInto & & \dInto & & \\
H_1(\Sigma,\Z)& \rTo & \bigoplus H_1(\Sigma/\p,\Z)^2 & \rTo & \ldots \\
\dOnto & & \dOnto & & \\
H_1(\Sigma,\Z)^{\tf} & \rTo &   0& \rTo & \ldots \\
\end{diagram}
$$
There is an exact sequence
$H_1(\Sigma,\Z)^{\tf} \rightarrow H_{d-1}(H^{\tors}_{1,\bullet}) \rightarrow H_{d-1}(H_{1,\bullet})$.
A class $[c] \in  H_{d-1}(H^{\tors}_{1,\bullet})$ may be represented by an
element $\gamma \in \bigoplus H_1(\Sigma/\p,\Z)^{2}$ up to the image of
$H_1(\Sigma,\Z)^{\tors}$. To say that it maps to zero in $H_{d-1}(H_{1,\bullet})$ is exactly
to say that $\gamma$ lifts to a class in $H_1(\Sigma,\Z)$, that is, $H_1(\Sigma,\Z) \ne 0$.
To come from a nontrivial class in $H_1(\Sigma,\Z)^{\tf}$ exactly says that 
$\gamma$ is not itself in the image of $H_1(\Sigma,\Z)^{\tors}$. Hence any such class satisfies
the conditions of the Lemma. Conversely, if $\m$ lies in the support of the module in question, then it clearly
lies (nontrivially) in the support of $H_1(\Sigma,\Z)^{\tf}$ and 
$\bigoplus H_1(\Sigma/\q,\Z)$.
\end{proof}

We deduce the following:

 \begin{lemma} Suppose that $H_1(\Sigma/\p,\Q) = 0$ for all $\p \in \Sigma$. 
 Then there is an identity
 $$\begin{aligned}{}
 [H_1(\Sigma,\Z)^{\tors,\new}_{\m}] = & \ \sum_{V \subset T}  (-2)^{|\Sigma \setminus S \cup V|}  [H_1(S \cup V,\Z)^{\tors}_{\m}] \\
& \  + (-1)^{|T|} \sum_{i=1}^{|T|} (-1)^i [H_i(K[1/T],\cF)_{\m}] + [\LL(\Sigma)], \end{aligned}$$
where $\LL(\Sigma)$ denotes a  finite $\T$-module which is nontrivial if and only if the following two conditions both hold:
\begin{enumerate}
\item There exists a non-zero characteristic zero form in $H_1(\Sigma,\Z)_{\m}$.
\item There exists a non-zero torsion class in $H_1(\Sigma/\p,\Z)_{\m}$ for some $\p$.
\end{enumerate}
\end{lemma}

\begin{proof} We apply Lemmas~\ref{lemma:orderofgroup2} 
and~\ref{lemma:short} to the identity relating $\Omega^d(H^{\bullet}_1)$ and $\Omega_d(H_{1,\bullet})$ above.
We use the fact that the only non-zero term in the second row occurs at $(2,d)$, and thus the identification
of the terms in $\Omega_d(H_{1,\bullet})$ with $E^{\infty}_{1,i}$ for $i < d$ allows us to compute
$\Omega_d(H_{1,\bullet})$ in terms of the arithmetic cohomology groups, as in the finite case.
The error term  $[\im(H_1(\Sigma,\Z)^{\tf}
 \rightarrow H_{d-1}(H^{\tors}_{1,\bullet}))]$ is a genuine module (i.e., not virtual), and it is nontrivial
 exactly under the conditions given, as established in Lemma~\ref{lemma:support}.
\end{proof}

\subsection{Splitting at non level-raising primes}
Suppose that $\m$ is a maximal ideal of $\T$ such that $T^2_{\q} - (1 + N(\q))^2 \notin \m$.
Then we claim that the maps
$$
\begin{aligned}
\Psi: \ & \ H_1(S \cup V \cup \{\q\},\Z)_{\m} \rightarrow H_1(S \cup V,\Z)^2_{\m} \\
\Psi^{\vee}: \ & \ H_1(S \cup V,\Z)^2_{\m} \rightarrow H_1(S \cup V \cup \{\q\},\Z)_{\m} 
\end{aligned}$$
are  isomorphisms. This is (essentially) a trivial consequence of the
``only if'' part of level raising 
(Theorem~\ref{theorem:ribet}).
Suppose then that $\m$ is a maximal ideal of $\T$ such that
$T^2_{\q} - (1 + N(\q))^2 \notin \m$ for \emph{some} $\q \in T$, where
$\Sigma = S \cup T$. Then, writing $H_{1,\bullet}(\Sigma)$ and
$H^{\bullet}_{1}(\Sigma)$ for what we have been calling $H_{1,\bullet}$
and $H^{\bullet}_{1}$ respectively, we have isomorphisms:
$$H^{\bullet}_{1}(\Sigma) =H^{\bullet}_{1}(\Sigma/\q) + H^{\bullet}_{1}(\Sigma/\q)[1],$$
$$H_{1,\bullet}(\Sigma) =H_{1,\bullet}(\Sigma/\q) + H_{1,\bullet}(\Sigma/\q)[-1].$$
In particular, we deduce that
$$\Omega_{d}(H_{1,\bullet}) =  \Omega^{d}(H^{\bullet}_1) 
= \Omega^{d}(H^{E,\bullet}_1) =  0,$$
where the last equality comes from the fact that
$T^2_{\q} - (1 + N(\q)) \notin \m$ implies that $H_1(\Sigma,\Z)_{\con,\m} = 0$.
It follows that the only maximal ideals that one need be concerned about
are those for which $T^2_{\q} - (1 + N(\q))^2 \in \m$ for all $\q \in T$.
That is, we have the following:

 \begin{lemma}   \label{lemma:unconcerned} Suppose that $\m$ does not contain $T^2_{\q} - \langle \q \rangle (1 + N(\q))^2$
 for some $\q \in T = \Sigma \setminus S$.
 There is an identity
 $$[H_1(\Sigma,\Z)^{\new,\tors}_{\m}] = \sum_{V \subset T}  (-2)^{|\Sigma \setminus S \cup V|}  [H_1(S \cup V,\Z)^{\tors}_{\m}].$$
 \end{lemma}
 
\begin{proof} If $\m$ does not occur in characteristic zero, then  the finite case argument (\S~\ref{section:caseoffinitehomology}) applies. 
If $\m$ does occur in the support of $H_1(S \cup U,\Z)$ for some $U \subsetneq T$, but
$T^2_{\q} - (1 + N(\q))^2 \notin \m$, then the analysis above shows
that the relevant sequences do not contribute to the space of newforms
at $\m$.
\end{proof}

\subsection{Isolated characteristic zero forms}
We consider the following hypothesis:
\begin{itemize}
\item After localizing at $\m$, there is a single automorphic representation
occuring in cohomology of level dividing $\Sigma/\q$ for some $\q$.
\item The $\m$-torsion vanishes integrally outside levels $\Sigma/\q$ and $\Sigma$.
\end{itemize}
By assumption, the level raising sequence $H^{\bullet}_{1}$ has length $2$ in this case.
This situation was exactly analyzed in 
\S~\ref{sec:lloneprime}, in the more difficult context where $\m$ was not assumed to
be non-Eisenstein. We immediately deduce the following:

 \begin{lemma} Suppose that the only non-zero automorphic newform
 contributing to the cohomology $H_1(\Sigma,\Q)$ after localizing at $\m$ occurs at level $\Sigma/\p$.
 Then there is an identity
 $$\begin{aligned}{}
 [H_1(\Sigma,\Z)^{\tors,\new}_{\m}] = & \ \sum_{V \subset T}  (-2)^{|\Sigma \setminus S \cup V|}  [H_1(S \cup V,\Z)^{\tors}_{\m}] \\
& \  + [H_1(\Sigma/\p,\Z)^{\tf}/(T^2_{\q} - \langle \q \rangle (1 + N(\q))^2)], \end{aligned}$$
\end{lemma}

\section{Newforms and the torsion quotient (II) } \label{sec:nonny2}

This section, which is only for the masochistic reader, analyzes certain ``obvious'' contributions to the spectral sequence.
Namely, the congruence homology contributes to the $H_1$ row (\S~\ref{h1sse})
and the fundamental classes of the relevant $3$-manifold allow the
$H_3$ row to be precisely understood (\S~\ref{section:coco}); in the split case
the fundamental classes of the cusps contribute to the $H_2$ (\S~\ref{section:nococo}).   All of these contributions
were invisible in the previous section, because they vanish after localizing at Eisenstein primes. 

Once these have been analyzed, we define (in a somewhat {\em ad hoc} way)
higher essential homology by excising these.

The final goal of this section is to prove Lemma~\ref{summary}, which relates the orders of various
spaces of newforms to spaces of oldforms and ``correction'' factors coming from the cohomology
of $S$-arithmetic groups. We are able (mostly) to control all terms only when $|T|
= |\Sigma \setminus S|$ has order two, and must otherwise avoid various residue characteristics
which divide certain factors of $N(\p) + 1$ for various $\p$. (Indeed, if one is content to consider the
case $|T| = 2$ then many difficulties of this section disappear, since only a very small portion of the spectral sequence
comes into play. This is why we discussed various special cases concerning $|T| = 2$
in previous sections.)

We continue with the notations of the prior section: {\em the Hecke algebra $\T$ will denote the abstract Hecke algebra omitting primes
$T_{\p}$  dividing $\Sigma$. }

\subsection{The \texorpdfstring{$H_1$}{H_1}-row of the spectral sequence.} \label{h1sse}

 Let us examine more closely the $H_1$ row of the spectral sequence, namely,
  $$H_{1,\bullet} := H_1(S,\Z_p)^{2^d} \leftarrow \bigoplus  H_1(S \p,\Z_p)^{2^{d-1}}   \leftarrow 
\bigoplus H_1(S \p\q  ,\Z_p)^{2^{d-2}} \leftarrow \cdots $$
We have already noted that the first row $H_{0,\bullet}$ is exact away from $i = 0$, and hence
for $i \ge 1$ there exists a boundary map
$$H_i(K[1/T],\cF) \rightarrow H_{i-1}(H_{1,\bullet}).$$
On the other hand, there is a surjective map
$H_{1,\bullet} \rightarrow H_{1,\con,\bullet}$,
where
$H_{1,\con,\bullet}$ is the complex
$$H_{1,\con}(S,\Z)^{2^d}  \leftarrow \bigoplus  H_{1,\con}(S \p,\Z)^{2^{d-1}}   \leftarrow 
\bigoplus H_{1,\con}(S \p\q  ,\Z)^{2^{d-2}}  \cdots $$
The congruence homology is always in the $-$ eigenspace of the
Fricke involution, and hence this sequence is some distance away from being exact.
There is an induced map
$$H_{*} (H_{1,\bullet}) \rightarrow H_{*} (H_{1,\con,\bullet}).$$
We have:
\begin{lemma} \label{lemma:userefined} Assume that the conditions of
Lemma~\ref{lemma:refined} are satisfied; then the map $H_{1,\bullet} \rightarrow H_{1,\con,\bullet}$ induces a
surjective map on homology.
\end{lemma}

\begin{proof} This is an immediate consequence of Lemma~\ref{lemma:refined}.
\end{proof}

We obtain a corresponding boundary map
$$H_i(K[1/T],\cF) \rightarrow H_{i-1}(H_{1,\con,\bullet}).$$
We have the following:
\begin{lemma} We have equalities: \label{lemma:homologycongruence}
$$H_0(H_{1,\con,\bullet}) = H_1(S,\Z)_{\con},$$
$$H_1(H_{1,\con,\bullet}) = \bigoplus_{\p} H_1(S \cup \p,\Z)^{\new}_{\con} \simeq \bigoplus (k^{\times}_{\p}),$$
$$H_i(H_{1,\con,\bullet}) = 0, \quad i \ge 2.$$
The boundary map $H_i(K[1/T],\cF) \rightarrow H_{i-1}(H_{1,\con,\bullet})$ is surjective.
\end{lemma}

\begin{proof}
Note that
the map from $H_1$ is surjective by simply considering the congruence homology in the arithmetic group, and is trivially surjective for $i > 2$ because the target is zero,
so the only subtle case is $i = 2$.
The computation of the $H_0$-row of the 
spectral sequence in~\S~\ref{section:H0row} implies that $E^{\infty}_{0,k}
= E^{2}_{0,k} = 0$ for  any $k \ge 1$. It follows that the boundary maps
$$0 = E^{2}_{3,0} \rightarrow E^{2}_{1,1}, \qquad E^{2}_{1,1} \rightarrow E^{2}_{-1,2} = 0$$
are tautologically trivial, and hence
$E^{\infty}_{1,1} = E^{2}_{1,1}$ and $E^{\infty}_{2,0} = 0$.
The convergence of the spectral sequence implies that
$H_{2}(K[1/T],\cF)$ admits a surjection to $E^{\infty}_{2,0} = 0$
whose kernel surjects onto $E^{\infty}_{1,1} = E^{2}_{1,1}
= H_{1}(H_{1,\bullet})$. 
We deduce that this is a surjective  map
$$H_{2}(K[1/T],\cF) \rightarrow E^{\infty}_{1,1} = H_{1}(H_{1,\bullet}),$$
which then surjects onto $H_{1}(H_{1,\con,\bullet})$ by
Lemma~\ref{lemma:userefined}  above.
\end{proof}

\begin{remark} The map for $i = 2$ may be considered as an
$S$-arithmetic analogue  of the tame
symbol map 
(discussed in~\S~\ref{section:stuffabouttame}).
\end{remark}

\begin{remark} \label{remark:highercong} Note that $H_1(S,\Z)_{\con}$ is
exactly $H_{1,\con}(K[1/T],\Z)$.
\end{remark}

\subsection{The \texorpdfstring{$H_3$}{H_3}-row of the spectral sequence, non-split case} \label{section:coco}
In this section, we assume that the implicit quaternion algebra $D$ is not split.
In particular, the manifold $Y$ and its congruence covers are compact.  Suppose also
that $p$ is not an orbifold prime.

There is an isomorphism $H_3(Y,\Z_p) \simeq \Z_p^C$, where $C$ denotes the number of connected components of $Y(K)$.
In this section, we will be interested in studying the level lowering sequence:
$$H_{3,\bullet}:=H_3(S,\Z_p)^{2^d} \leftarrow  \bigoplus H_3(S\p,\Z_p)^{2^{d-1}}
\leftarrow \bigoplus H_3(S \p \q,\Z_p)^{2^{d-2}} \leftarrow \ldots$$

Recall that, given a subset $\{x_1, \dots, x_n\}$ of $\Z_p$,  
the associated Koszul complex $\mathcal{K}(\{x_1, \dots, x_n\})$  is given by the tensor product
of the length $2$ complex
$\Z_p  \stackrel{\times  x_i}{\longleftarrow}  \Z_p$:
$$ \mathcal{K}(\{x_1, \dots, x_i\}) = \bigotimes_{i=1}^n  [ \Z_p  \stackrel{\times  x_i}{\longleftarrow}  \Z_p].$$

\begin{lemma} \label{Koszul1} The sequence $H_{3,\bullet}(\Z_p)$ is a direct sum of
$2^d$  Koszul complexes of length $0 \le i \le d$, one for each subset
$U \subseteq \Sigma \setminus S$:
$$ \bigoplus_{U} \mathcal{K}(Y_U),$$ where 
$Y_U = \{ N(\p)+1: \p \in U\}$.
 \end{lemma}

\begin{proof} We note that the fundamental class in $H_3(Y_0(S \cup T),\Z)$  is invariant under the Fricke involution
$w_{\p}$ for every $\p \in T$. Let $\eps$ be a choice of sign for the Fricke involution for each
$\p \in \Sigma \setminus S$. Let $U$ denote the subset of elements for which this choice of sign is $1$. Localizing at
$w_{\p} - \eps(\p)$ for each $\p \in T$, the compex $H_{3,\bullet}$ becomes:
$$\Z \leftarrow \bigoplus_{\p \in U} \Z \leftarrow \bigoplus_{\{\p,\q\} \subset U} \Z \leftarrow \ldots$$
where the classes $\Z$ for $V \subseteq U$ correspond to the fundamental class in $H_3(Y_0(S \cup V),\Z)$.
If $V \cup \p \subset U$, we are interested in computing the corresponding map
$$\Z = H_3(S \cup V \cup \p,\Z) \rightarrow H_3(S \cup V,\Z) = \Z.$$
This corresponds to the ``average'' of the two projections coming from the natural map $Y_0(\p)
\rightarrow Y$ and the corresponding map twisted by $w_{\p}$. Both maps are multiplication by the degree, so
the map above is just multiplication by $1 + N(\p)$.
\end{proof}

Using this, we may describe the $E^{1}_{0,3}$ term of the spectral sequence. It contains the following terms:
\begin{enumerate}
\item A copy of $\Z$, coming from the Koszul complex of length $0$.
\item A copy of 
$$\bigoplus_{\q \in T} \Z/(N(\q) + 1)\Z,$$
 coming from the Koszul complexes of length $1$.
\item More generally, a copy of
$$ \Z/\mathrm{gcd}(N(\q_1) + 1,N(\q_2) + 1, \ldots, N(\q_i) + 1)\Z, \qquad \q_i \in U$$
for all $U \subset T$.
\end{enumerate}

There exists a boundary map:
$$H_{i}(H_{\bullet,3}) = E^{2}_{i,3} \rightarrow  H_{3}(Y(K[1/T],\cF) $$
for all $i \ge 0$.
The $\Z$-factor injects.
Its existence corresponds to the fact that $K_3(\OL_F) \otimes \Q \simeq \Q$,
since $F$ is exactly one complex place.

\medskip

It is somewhat of a mystery what the maps $\Z/(N(\p) + 1) \rightarrow H_3(Y(K[1/T],\cF)$ are,
and this provides an obstruction to completely analyzing what happens in the Eisenstein case.
We note, however, the following:

\begin{remarkable}
Suppose that $d = |T| \le 2$, or that the residue characteristic of $\m$ does not not divide
$N(\p) +1$ for any $\p \in T$. Then the $H_3$-row has no effect on the previous calculations
of $H^{E}_1(\Sigma,\Z)^{\tors,\new}$. This is obvious in the second case, and also clear in
the first case since the nontrivial terms are too deep in the spectral sequence to make any difference.
\end{remarkable}

\subsection{The boundary contribution to the \texorpdfstring{$H_2$}{H_2}-row of the spectral sequence, split case} \label{section:nococo}
In this section, we assume that the implicit quaternion algebra $D$ is split.
It follows that $X$ is non-compact, and that the $H_3$-row of the spectral sequence vanishes.
 We shall see here, however, that the same invariants arise (shifted by $1$) via the boundary cohomology.

We have the following sequences:
$$H^{\bm}_{3,\bullet}:=H^{\bm}_3(S,\Z_p)^{2^d} \leftarrow  \bigoplus H^{\bm}_3(S\p,\Z_p)^{2^{d-1}}
\leftarrow \bigoplus H^{\bm}_3(S \p \q,\Z_p)^{2^{d-2}} \leftarrow \ldots$$
$$H(\partial)_{2,\bullet}:=H_2(\partial S,\Z_p)^{2^d} \leftarrow  \bigoplus H_2(\partial S\p,\Z_p)^{2^{d-1}}
\leftarrow \bigoplus H_2(\partial S \p \q,\Z_p)^{2^{d-2}} \leftarrow \ldots$$
$$H_{2,\bullet}:=H_2(S,\Z_p)^{2^d} \leftarrow  \bigoplus H_2(S\p,\Z_p)^{2^{d-1}}
\leftarrow \bigoplus H_2(S \p \q,\Z_p)^{2^{d-2}} \leftarrow \ldots$$
Here, by abuse of notation, we let $H_2(\partial S,Z)$ (etc.) by the cohomology of the boundary of the corresponding manifold.
There is an exact sequence:
$$H^{\bm}_{3,\bullet} \rightarrow H(\partial)_{2,\bullet} \rightarrow H_{2,\bullet}$$

\begin{lemma} 
Suppose $p \neq 2$.
\begin{itemize}

\item[(i)] The $H^{\bm}_{3, \bullet}$ is isomorphic to the complex described in Lemma~\ref{Koszul1}. 

\item[(ii)] 
The complex $H(\partial)_{2, \bullet}$ is isomorphic to:
$$ \bigoplus_{U \subset \Sigma \setminus S}   \mathcal{K}(X_U) ^{h},$$
where $X_U = \{ N(\p) +1: \p \in U\} \bigcup \{ N(\p) - 1: \p \notin  \Sigma \setminus (S \cup U) \}$,
and $h$ is the number of cusps of $Y(K_S)$. 

 \item[(iii)]  The image of the sequence $H(\partial)_{2,\bullet}$ is a direct sum of
complexes
$$ \bigoplus_{U \subset \Sigma \setminus S } \left( \mathcal{K}(X_U) ^{h- b_0(Y)} \oplus  \mathcal{K}( Y_U )[|U|] \right)^{b_0(Y)} ,$$ 
where $Y_U =  \{ N(\p) +1: \p \in U\}$, and $b_0$ is the number of components of $Y(K_S)$. 
 \end{itemize}

\end{lemma}
We do not give a proof since the details of this result are never used;
see, however,~\S~\ref{topologyprojectiondown}

\subsection{\texorpdfstring{$S$}{S}-arithmetic essential homology} \label{section:essentialhomologyone}
In this section we give (somewhat {\em ad hoc}) definitions of the ``essential'' homology of $S$-arithmetic groups, roughly just by cutting out all the contributions which occur
for the reasons already detailed in this section. 
(We use some of the identifications made in
 Lemma~\ref{lemma:homologycongruence} and
 Remark~\ref{remark:highercong}.)

 \begin{df}
 Define $H^E_{i}(K[1/T],\cF)$  as follows: 
 \begin{itemize}
 \item 
 If $i = 1$, then $H^E_{1}(K[1/T],\cF)$ is the kernel of the map
 $$H_1(K[1/T],\cF) \rightarrow H_{0}(H_{1,\con,\bullet}) \simeq H_{1,\con}(K[1/T],\cF).$$
 This agrees with the usual definition when $T = \emptyset$.
  \item If $i = 2$, then $H^{E}_2(K[1/T],\cF)$ is
  the kernel of the map
  $$H_2(K[1/T],\cF) \rightarrow H_{1}(H_{1,\con,\bullet}) 
\simeq \bigoplus (k^{\times}_{\p}) .$$
\item  If $i \ge 3$
and $\G$ is nonsplit, 
let  $H^E_{i}(K[1/T],\cF)$ it is the cokernel of the map
$$H_{i-3}(H_{3,\bullet}) \rightarrow H_i(K[1/T],\cF).$$
\item If $i \geq 3$ and $\G$ is split,  it is the cokernel of the map:
$$H_{i-2}(\mathrm{im}(H(\partial)_{2,\bullet}) \rightarrow H_{i-2} (H_{2,\bullet}) \rightarrow H_i(K[1/T],\cF).$$
\end{itemize}
\end{df}

As we have seen in the prior section, the essential homology can be nonvanishing
for reasons related to $K_2$. 
\begin{lemma} If $i \le d = |T|$, then the groups $H^{E}_i(K[1/T],\cF)$ are finite.
If $\m$ is non-Eisenstein, then $H^{E}_i(K[1/T],\cF)_{\m} \simeq H_i(K[1/T],\cF)_{\m}$.
\end{lemma}

\begin{proof} Considering the spectral sequence over $\Q$, we see that the group
$H_i(K[1/T],\cF)$ is
finite for $i \le d$ with the exception of $i = 3$, which contains a copy of $\Z$. This copy of $\Z$ is
accounted for exactly by the fundamental class (in the co-compact case (see~\S~\ref{section:coco}) and the image of the
boundary map in the non co-compact case (see~\S~\ref{section:nococo}).
\end{proof}

\begin{remark}
For $4 \le i \le d = |T|$, the ratio
$$\frac{H^{E}_i(K[1/T],\cF)}{H_i(K[1/T],\cF)^{\tors}}$$
is only divisible by primes dividing $N(\p) + 1$ for $\p \in T$.
This follows from the arguments 
of~\S~\ref{section:coco} and~\S~\ref{section:nococo}.
For $i = 3$, however, the situation is somewhat more subtle.
 The torsion free part of $H_3(K[1/T],\Z)$ has rank equal
 to the number of connected components of $Y(K)$ providing thatr $|T| \ge 3$.
In particular, after inverting primes dividing 
$H_3(K[1/T],\cF)^{\tors}$,
 the group
$H^E_3(K[1/T],\cF)$ is given by the cokernel of the map
$$H^3(K,\Z)^{\tf} \rightarrow H^3(K[1/T],\Z)^{\tf}.$$
Suppose that $Y(K)$ is connected.
If this map was surjective for all $T$, then by taking the limit, one would also
obtain a surjective map
$$\Z = H^3(Y(K),\Z)^{\tf} \rightarrow K_3(\OL_F)^{\tf} = \Z.$$
However, this is no reason to expect that the fundamental class
of $H^3(Y(K),\Z)$ is not divisible in $K_3(\OL_F)$. 
 \end{remark}

\begin{df} \label{df:pec}
The partial Euler characteristic $\chi(\Sigma)$ is defined as follows. We let
$\chi(\Sigma) := \prod_{p} \chi_p(\Sigma)$, where
$$\chi_p(\Sigma)^{(-1)^d}:= \prod_{i=1}^{d} |H^E_i(K_{\Sigma}[1/T],\cF_p)|^{(-1)^i},$$
and $\mathcal{F}_p$ is the localization at $p$ of the sheaf defined
before Theorem \ref{SarithmeticSS}.
We have
$$\chi(\Sigma) = \left| (-1)^{|T|}  \sum_{i=1}^{d} (-1)^i [H^E_i(K_{\Sigma}[1/T],\cF)] \right| $$ 
where the equality holds up to a rational factor divisible only by orbifold primes.
\end{df}

 The group $H^{E}_1$ vanishes whenever the congruence subgroup
 property holds, and thus we note the following:
 
 \begin{lemma} \label{usefullater}
 If $d = 1$, and the congruence subgroup property holds for
 $K_{\Sigma}[1/T]$, then $\chi(\Sigma) = 1$.
 If $d = 2$, and the congruence subgroup property holds for
 $K_{\Sigma}[1/T]$, then $\chi(\Sigma) \in \Z$.
 \end{lemma}

It turns out that for $d >1$, the quantity $\chi(\Sigma)$ can differ from $1$ even when $H_1(\Sigma, \Q) = 0$;
 see 
~\S~\ref{subsection:hangingchads} for numerical examples; as we have discussed,
 its order is related to $\zeta_F(-1)$ through $K$-theory. 
 
 \medskip
 
 We have the following:
 
  \begin{lemma} \label{summary}
 Suppose that at least one of the following conditions holds:
 \begin{enumerate} 
 \item $|T| \le 2$,
 \item  $\m$ is not Eisenstein,
 \item $\m$ has residue characteristic not dividing $N(\p) + 1$ for any $\p$ in $T$.
 \end{enumerate}
 Then the following hold:
 \begin{enumerate}
 \item If $H_1(\Sigma,\Q)_{\m} = 0$,  there is an identity
 $$\begin{aligned}{} [H^{E}_1(\Sigma,\Z)^{\new}_{\m}] =  & \ \sum_{V \subset T}  (-2)^{|\Sigma \setminus S \cup V|}  [H^{E}_1(S \cup V,\Z)_{\m}] \\  &  \ 
 + (-1)^{|T|} \sum_{i=1}^{|V|} (-1)^i [H^{E}_i(K[1/T],\cF)_{\m}]. \end{aligned}$$
\item   Suppose that $H_1(\Sigma/\p,\Q)_{\m} = 0$ for all $\p \in \Sigma$. 
 Then there is an identity
 $$\begin{aligned}{}
 [H^{E}_1(\Sigma,\Z)^{\tors,\new}_{\m}] = & \ \sum_{V \subset T}  (-2)^{|\Sigma \setminus S \cup V|}  [H^{E}_1(S \cup V,\Z)^{\tors}_{\m}] \\
& \  + (-1)^{|T|} \sum_{i=1}^{|V|} (-1)^i [H^{E}_i(K[1/T],\cF)_{\m}] + [\LL(\Sigma)], \end{aligned}$$
where $\LL(\Sigma)$ denotes a non-negative finite $\T$-module which is nontrivial if and only if
\begin{enumerate}
\item There exists a non-zero characteristic zero form in $H^{E}_1(\Sigma,\Z)_{\m}$.
\item There exists a non-zero torsion class in $H^{E}_1(\Sigma/\p,\Z)_{\m}$ for some $\p$.
\end{enumerate}
\item  Suppose that the only non-zero automorphic newform
 contributing to the cohomology $H_1(\Sigma,\Q)_{\m}$  occurs at level $\Sigma/\p$, and that
 $H_1(S \cup V,\Z)_{\m} = 0$ unless $V = T$ or $T/\q$.
 Then there is an identity
 $$\begin{aligned}{}
 [H^{E}_1(\Sigma,& \Z)^{\tors,\new}_{\m}] =  \ \sum_{V \subset T}  (-2)^{|\Sigma \setminus S \cup V|}  [H^{E}_1(S \cup V,\Z)^{\tors}_{\m}] \\
  + \ & [H^{E}_1(\Sigma/\p,\Z)^{\tf}/(T^2_{\q} - \langle \q \rangle (1 + N(\q))^2)]
- [H^1_{\cl}(\Sigma/\q,\Z)_{\m}] \\
 + \ & (-1)^{|T|} \sum_{i=1}^{|V|} (-1)^i [H^{E}_i(K[1/T],\cF)_{\m}], \end{aligned}$$
\item Suppose that $\m$ does not contain $T^2_{\q} - \langle \q \rangle (1 + N(\q))^2$
 for some $\q \in T = \Sigma \setminus S$.
 There is an identity
 $$[H^{E}_1(\Sigma,\Z)^{\new,\tors}_{\m}] = \sum_{V \subset T}  (-2)^{|\Sigma \setminus S \cup V|}  [H^{E}_1(S \cup V,\Z)^{\tors}_{\m}].$$
 \end{enumerate}
 \end{lemma}
 
 \begin{proof}
 If $\m$ is not Eisenstein, then the lemma has already been proven. In general, we have to worry about the following two issues:
 \begin{enumerate}
 \item Contributions coming from $H_{3,\bullet}$ in the co-compact case, and the image of $H(\partial)_{2,\bullet}$
 in $H_{2,\bullet}$  in
 the split case.
 \item Congruences between cusp forms and Eisenstein series in the split case.
 \item Congruence homology.
 \end{enumerate}
 The congruence homology is dealt with by the modification to the $S$-arithmetic cohomology in degrees $i \le 2$.
 The contributions from the higher rows have no effect when $|T| \le 2$, since they are too far out in the spectral
 sequence.  The only difference occurs when we allow newforms for $|T| \le 2$, and this is accounted for by
  Lemma~\ref{theorem:regcomp}.
  \end{proof}

\section{The general case} \label{section:finalsummary}
In this section, we discuss what can be said in the most general framework without
any assumption, summarizing the prior two sections.

Consider the following ratio:
\begin{equation}\label{ratioR} \frac{|H^{E}_1(Y(\Sigma),\Z)^{\new,\tors}|}{|H^{E}_1(Y'(\Sigma),\Z)^{\tors}|} \cdot \frac{\reg^{\new}_{E}(Y')}{\reg^{\new}_{E}(Y)}. \end{equation} 
What is our expectation for this quantity?

\medskip
 We expect
that  spaces of newforms  should be equal when localized at some
maximal ideal $\m$, with the  exception of a correction factor related to $K$-theoretic classes.

\medskip
Similarly,
we also expect (see the discussion in~\S~\ref{section:prasanna}) that the new regulators
should differ exactly at level lowering primes; that is,  $p$ dividing the residue characteristic of primes $\m$ which
are supported on the space of newforms which lift to characteristic zero, but are \emph{also} supported --- possibly at
the level of torsion --- at lower level.
(Note that, if neither $Y$ nor $Y'$ are split, there can be level lowering primes contributing to both the numerator and
denominator of this quantity.) 

\medskip

On the basis of this discussion,  we expect that this quantity \eqref{ratioR} is always equal to:
\begin{equation} \label{ratioR2}   \frac{\chi_D(Y)}{\chi_D(Y')} \cdot \frac{|\LL(Y)|}{|\LL(Y')|}, \end{equation} 
where $\LL(Y)$ and $\LL(Y')$ are finite $\T$-modules recording the level lowering congruences described above,
and $\chi_D$ is related to $K$-theory. (In particular, if $Y$ is split, then $\LL(Y')$ should be trivial
and so $|\LL(Y)| \in \Z$ is divisible exactly by level lowering primes.)

\medskip

What can we prove?
The first point to note is that we have no exact guess for the order of $\chi_D(Y)$ or $\chi_{D}(Y')$,
or rather, we have a definition given in terms of the $S$-arithmetic cohomology of arithmetic groups (see Definition \ref{df:pec}). 
 The $S$-arithmetic  classes corresponding  are poorly understood, especially
 for $d$ large --- we do not know whether they are
 Eisenstein or not except (under CSP) when $d = 1$.

 \medskip
 
Having said this, the analysis of the prior two sections {\em does} imply that, away from orbifold primes, $$ \frac{|H^{E}_1(Y(\Sigma),\Z)^{\new,\tors}|}{|H^{E}_1(Y'(\Sigma),\Z)^{\tors}|} \cdot \frac{\reg^{\new}_{E}(Y')}{\reg^{\new}_{E}(Y)} 
 =   \frac{\chi_D(Y)}{\chi_D(Y')} \cdot \frac{|\LLL(Y)|}{|\LLL(Y')|}$$
where $\chi_D$ is as in Definition \ref{df:pec}, 
and $\LLL(Y)$ and $\LLL(Y')$
are  virtual $\T$-modules related to $Y$ and $Y'$   with the following property:
 
 \medskip
 
 $\LLL(Y)$ is finite and has support at an ideal $\m$ of $\T$ only if one of the following occures:
 \begin{enumerate}
 \item  $\m$ occurs in the support of
 $H_1(\Sigma,\Z)^{\tf}$, but also  occurs in the support of $H_1(T,\Z)$ for
some $S \subset T \subsetneq \Sigma$ of smaller level, 
\item $d > 2$, $\mathfrak{m}$ has characteristic $p$, and $p$ divides 
$N(\p) - 1$ or $N(\p) + 1$ for some $p \in \Sigma \setminus S \cap S'$.
\end{enumerate}

\medskip

This is very similar to what we predicted above, but it is weaker for the following reasons:
\begin{enumerate}
\item If $d > 2$, we lose control of factors dividing $N(\p) - 1$ or $N(\p) + 1$. This comes from
the difficulty of understanding the interaction of the $H_3$-row of the spectral sequence in the non-split
case (and similarly, the contribution to the $H_2$ row from the fundamental class of the cusps).
\item The level lowering primes $\mathfrak{m}$ we expect to occur in $\LL(Y)$ should actually have support in $H_1(\Sigma,\Z)^{\tf,\new}$;
that is, occur  as characteristic zero \emph{newforms}.
Moreover,  although we know that $\m$ contributing to $\LL(Y)$ above have the listed property, there is no converse;
similarly, we can't control (for big $d$) whether such primes contribute to the numerator or denominator.
\end{enumerate}
It would be good to improve our understanding on any of these  points!

Finally, it would of course be good to show that the terms of \eqref{ratioR} and \eqref{ratioR2}  
match up ``term by term'', and not merely that their product is equal.

\chapter{Eisenstein Deformations, and Even Galois Representations}
\label{chapter:ch7}

In this chapter, we shall discuss certain simple consequences of our general conjectures that can be directly established.   Roughly speaking, in certain situations, 
the generalized multiplicity of a certain Hecke  eigensystem 
is automatically $\geq 2$. This phenomenon has been observed
(by Haluk \Sengun, personal communication; see also
section~~\ref{ss:Serre})
in numerical investigations. We will prove this unconditionally and
then verify that it is indeed a consequence of our Galois-theoretic conjectures.

 Namely, we shall {\em prove} directly (\S~\ref{involutions} for (a), (b); 
 Theorem~\ref{theorem:K2popularversion} for (c)) that:
 
\begin{enumerate}
\item[(a)] Let $\rhobar$ be an even representation $G_{\mathbf{Q}} \rightarrow
\GL_2(k)$, and $\Sigma$ any level such that  $\rhobar|G_F$ contributes to the cohomology of some $Y(\Sigma)$. 

Then the corresponding Hecke eigensystem $H_1(Y(\Sigma), \Z)_{\mathfrak{m}}$ has length
$\ge 2$ as a $\T_{\m}$-module; here $\m$ is the maximal ideal of the Hecke algebra
associated to $\rhobar$.

\item[(b)] (For $\G$ split:)  Suppose ${\q}$ is a prime of $F$, inert over $\mathbf{Q}$, 
such that $\mathrm{Norm}({\q}) \equiv 1$ modulo $p$. 

Then there exists  cylotomic-Eisenstein classes --- of type D1 in the notation~\S~\ref{eisdef} ---
for $H_1(Y_0({\q}), \Z_p)$, {\em beyond those arising from congruence homology}. 

\item[(c)] (For $\G$ split:) Suppose ${\q}$ is {\em any} prime of $F$, and 
the $\eps^{-1}$-part of $\CL(F(\zeta_p))$ is nonzero, where $\eps$ is the mod-$p$ cyclotomic character. 

Then there exists Eisenstein classes  (of trivial plus cyclotomic type, i.e., of type D1) for $H_1(Y_0({\q}), \Z_p)$, {\em beyond those arising from congruence homology}.

\end{enumerate}

On the other hand, we shall show in~\S~\ref{section:Eisenstein} and~\S~\ref{section:even}
 that these results are indeed consistent 
with  the Conjectures of~\S~\ref{sec:gl2inner}. 
More precisely, (a) follows from the conjectures of that section,
whereas (b) and (c) follow from  extensions of those conjectures to the residually reducible case.  
  Thus we regard the validity of (a) --- (c) as providing circumstantial evidence for our conjectures. 
 
 \medskip
 
 In~\S~\ref{section:phantomclasses}, we also speculate
 on extensions of Theorem~\ref{theorem:K2popularversion}(ii)  to higher levels, and their
 connection to so-called phantom classes, which are Eisenstein classes with
 some peculiar properties. The existence of phantom classes is exactly what seems to
 prevent  the spaces $H^1(Y(\Sigma),\Z)^{\new}$ and $H^1(Y'(\Sigma),\Z)^{\new}$ from having the
 same order in general. 
 
 In~\S~\ref{ss:Serre}, we discuss conjectural finiteness theorems for even Galois representations {\em over $\Q$}, as well as the question of finding an even Galois representation over $\Q$ 
 of tame level $1$.

 \section{Involutions} \label{involutions}
 Let us suppose we are given an {\em orientation-reversing isometry} $\sigma$ of some $\YO$
 which is compatible with the Hecke action, in that there exists
a compatible involution (also denoted $\sigma$) of the Hecke algebra. The natural
example of such an involution arises from Galois conjugation; we will discuss it below.

 Let $p \neq 2$ be a non-orbifold prime, and let $A = H_1(\Sigma, \Z)[p^{\infty}]$ be the
 (finite) $p$-torsion subgroup of $H_1(\Sigma,\Z)$. %
 The involution $\sigma$ induces an automorphism of  $A$, and hence
 induces a decomposition
 $A = A^{+} \oplus A^{-}$, where $A^{+}$ is the $\sigma$-fixed
 (respectively, anti-fixed) eigenspace. 
 
 We make use of the linking form, cf.~\S~\ref{section:linking} or~\S~\ref{section:linking2}; note that
even if $\mathbb{G}$ is split we are in the good situation described at the end of~\S~\ref{section:linking2}.
It is a pairing $\tau: A \times A \rightarrow \Q_p/\Z_p$  and
  satisfies $\tau(\sigma x, \sigma y) = - \tau(x, y)$. It follows that
 $A^{+}$ and $A^{-}$ are both isotropic with respect to $\tau$, and are dual to each other. In other terms, 
 $\tau$ induces an isomorphism of $A^{+}$ with the Pontryagin dual $\Hom(A^{-}, \Q_p/\Z_p)$
 of $A^{-}$.  This already means, for example, that $A$ cannot be of order $p$.
 
 A slightly more precise formulation of the phenomenon is:

  \begin{lemma}[Doubling] Notation as in the prior paragraph, suppose that \label{lemma:doubling}
$\m$ is a $\sigma$-stable ideal of $\T_{\Sigma}$ such that $\sigma$ acts trivially on $k:=\T_{\Sigma}/\m$, a field of characteristic $p$.  
   Then either:
  \begin{enumerate}
  \item $\T_{\m} :=\T_{\Sigma,\m} =  W(k)$, the Witt vectors of $k$, and $H_1(\Sigma,\Z)_{\m}$
  is free of rank one as a $\T_{\m}$-module.
  \item  $H_1(\Sigma,\F_p)_{\m}$ has length $\ge 2$ as a $\T_{\m}$-module
  \end{enumerate}
  \end{lemma}
  
  \begin{proof} It suffices to assume that
  $H_1(\Sigma,\F_p)_{\m}$ has length one and prove that there is an isomorphism
  $\T_{\m} = W(k)$.
  Note that $p \in \m$. Hence $H_1(\Sigma,\Z_p)_{\m} \otimes_{\T_{\m}} \T_{\m}/\m$
  is a quotient of 
  $$H_1(\Sigma,\Z_p)_{\m} \otimes_{\T_{\m}} \T_{\m}/p = H_1(\Sigma,\Z_p)_{\m}/p =
  H_1(\Sigma,\F_p)_{\m}.$$
  By assumption, the latter module has length one. Thus, by Nakayama's
  lemma, $H_1(\Sigma,\Z_p)_{\m}$ is cyclic as a $\T_{\m}$-module.
  Because $\T_{\m}$ acts faithfully on $H_1(\Sigma,\Z_p)_{\m}$, it
  follows that $H_1(\Sigma,\Z_p)_{\m}$ is free of rank one as a $\T_{\m}$-module;
  thus 
 $\T_{\m}/p \simeq H_1(\Sigma,\F_p)_{\m} \simeq \T_{\m}/\m$. By Nakayama's
  lemma applied to $\T_{\m}$ as a $W(k)$-module, we deduce that $\T_{\m}$
  is a quotient of $W(k)$. If $\T_{\m} = W(k)$, then we are in case~$(1)$.

  Thus it remains to prove the following: If $\T_{\m} = W(k)/p^n$ for some $n$,
  and $H_1(\Sigma,\Z_p)_{\m} $ is free of rank one over $\T_{\m}$, 
  then a contradiction ensues.

Since  $\sigma$ acts trivially on $k$, it also acts trivially on $\T_{\m}$, 
  and therefore induces a $W(k)$-linear involution of $H_1$;
  necessarily, this must be given by multiplication by an element $x \in (W(k)/p^n)^{\times}$,
  and since $x^2=1$, we must have $x = 1$.  This contradicts the fact that 
  $H_1(\Sigma, \Z_p)_{\m}^+$ and   $H_1(\Sigma, \Z_p)_{\m}^-$ 
are in duality, by the prior discussion. 
  \end{proof}

%

 There is a natural choice of $\sigma$ in the previous lemma, namely, we may take
$\sigma = c$ to be complex conjugation. In order for $c$ to be well defined, we assume
that $F/\Q$ is an imaginary quadratic field, $D^c = D$, and $\Sigma^c = \Sigma$.
For example, one can take $\G = \GL(2)/F$. 
Since it is known  (see~\cite{AC}) that a {\em classical} modular form $G$ over the field $F$ has the property 
that its Hecke eigenvalues are invariant by $c$ if and only if it is a base change, we deduce:

\begin{corollary} Let $\m$ be a maximal ideal of $\T$ of %
odd residue characteristic $p$, and let $c$ be complex conjugation;  suppose that $\m = \m^c, \Sigma = \Sigma^c, D = D^c$
and $c$ acts trivially on $\T/\m$.\label{cor:double}  Then
\begin{equation}  \label{*eqn}
(*)  \  \mbox{$H_1(\Sigma,\mathbf{F}_p)_{\m}$ has length $\ge 2$ as a $\T_{\m}$-module
}
\end{equation}
{\em unless} there exists a classical modular form $f$ of weight two giving rise to
$\m$ via base change from $\Q$ to $F$.
\end{corollary}

\begin{remarkable}
\emph{ One can explicitly describe the level of $f$ --- note, however,
 that
it may have nontrivial  level structure at primes away from $\Sigma$, in particular,
at primes which
ramify in $F/\Q$.
}
\end{remarkable}

\medskip

Another application of the interaction between the involution
induced by complex conjugation and the linking form is as follows.

\begin{theorem}  \label{theorem:square} Suppose that $F$ is an imaginary quadratic field.
Suppose, furthermore, that
$D^c = D$, $\Sigma^c = \Sigma$,  and $H_1(\Sigma,\Q) = 0$.
Then the order of $H_1(\Sigma,\Z)$ is a square up to orbifold primes
and powers of $2$.
More generally, suppose that $\m$ is a maximal ideal of $\T$
so that $H_1(\Sigma,\Z)_{\m}$ is finite, and such that the residue
characteristic $p$ of $\m$ is odd and not an orbifold prime.
Let $\m^* = \m^c$. Then:
\begin{enumerate}
\item If $\m \ne \m^*$, then $|H_1(\Sigma,\Z)_{\m}| = |H_1(\Sigma,\Z)_{\m^*}|$.
\item If $\m = \m^*$, then $|H_1(\Sigma,\Z)_{\m}|$ is a square.
\end{enumerate}
In particular, any odd prime  $p$ which
divides $|H_1(\Sigma,\Z)^{\tors}|$  to odd order is either an orbifold prime
or divides the residue characteristic of a maximal ideal $\m$ of $\T$ such
that $H_1(\Sigma,\Z)^{\tors}$ and $H_1(\Sigma,\Z)^{\tf}$ both have
support at $\m$.
\end{theorem}

\begin{proof}
Suppose that $H^1(\Sigma,\Z)^{\tors}$ has
support at a maximal ideal $\m$  of $\T$. The action of $c$ induces an isomorphism
of $H^1(\Sigma,\Z)$, so if $\m \ne \m^*$, then $H^1(\Sigma,\Z)_{\m}
\simeq H^1(\Sigma,\Z)_{\m^*}$. Suppose that $\m = \m^*$ has
residue characteristic $p$ and
that $H^1(\Sigma,\Z)_{\m}$ is finite. Then the linking form
induces a perfect pairing
$$H^1(\Sigma,\Z)_{\m} \times H^1(\Sigma,\Z)_{\m}
\rightarrow \Q_p/\Z_p.$$
As explained above, this pairing induces a decomposition of
$H^1(\Sigma,\Z)_{\m}$ into $c$-fixed and $c$-anti-fixed isotropic
subspaces which are dual, and hence have the same order.
\end{proof}

\subsection{Restriction of representations from \texorpdfstring{$G_{\Q}$}{G_Q}}
There is a natural class of examples to which the Corollary is applicable, obtained
by restricting representations of $G_{\Q}$: 

Take  $\rhobar^{\Q}: G_{\Q} \rightarrow
\GL_2(k)$ of Serre weight $2$ and squarefree tame level.  According to our conjecture, the restriction $\rhobar$ of $\rhobar^{\Q}$ to $F$
should appear in the cohomology $H_1(\Sigma, \Z)$ for suitable $\Sigma$, without loss of generality $c$-stable; if $\m$ is the corresponding maximal  ideal of $\T_{\Sigma}$, 
then the conditions of the Corollary hold. 

There are two situations of particular interest.

\subsubsection{Restriction of an even Galois representation} \label{sss:reven}
 Suppose characteristic of $k$ is not $2$, and take  $\rhobar^{\Q}$   {\em even}. Then it \emph{cannot} arise from a classical modular form, since such representations
are odd.       Thus we are in case~\eqref{*eqn}.

This has an interesting consequence on the Galois side: 

 If we assume that multiplicity one (Conjecture~\ref{conj:multone}) holds, then it follows
  from~\eqref{*eqn} that $\T_{\m}/p  \ne k$ and hence that
  $\dim_k(\m/(\m^2,p)) \ge 1$. So, once more assuming conjecture~\ref{conj:reciprocity}, 
~\eqref{*eqn} implies that  $\rhobar$ admits a global deformation of fixed determinant that is finite flat at $v|p$.  We will (at least if $p$ splits in $\OL_F$)  prove this unconditionally: see Theorem~\ref{theorem:evendoubling}.
This is a little surprising: odd
representations $\rhobar:G_{\Q} \rightarrow \PGL_2(k)$ need not admit finite flat deformations
of fixed determinant even after restriction to $F$. 

\subsubsection{Restriction of a reducible Galois representation} \label{sss:rreducible}  Take $\rhobar^{\Q}$ of the form $1 \oplus \eps$ where $\eps$
is the cyclotomic character, and $\Sigma = (N)$, where
$N$ is an  prime integer inert in $F$. One sees directly\footnote{
 Consider the group $\Gamma = \PGL_2(\OL_F)$, and suppose
that $p > 3$ (for this discussion). If $p | N^2 - 1$, there is a
surjective (congruence) map:
$$\Gamma_0(N) \rightarrow \Gamma_0(N)/\Gamma_1(N) \rightarrow \Z/p\Z,$$
giving rise to a class in $H_1(\Gamma_0(N),\Z)$. By construction, this class is
congruence, and thus cyclotomic-Eisenstein, and hence $\T_{\Sigma}$ has an cyclotomic-Eisenstein maximal
ideal $\m$ of residual characteristic $p$.
} that $\rhobar = \rhobar^{\Q}|_{G_F}$
is modular at level $Y(\Sigma)$ so long as $N^2 \equiv 1$ modulo $p$. 

Here there are two cases: 
\begin{itemize}
\item[(a)]
 If $N \equiv 1 \mod p$, then  (if $p>3$) there 
\emph{is} a classical (cuspidal) modular form giving rise to $\m$ via base change; this
follows from~\cite{MazEis}.

\item[(b)] If $N \equiv -1 \mod p$, it is no longer necessarily true that
$\m$ arises via base change.  We prove directly that
$\rhobar$ admits interesting deformations, at least under certain assumptions on $F$
(see Theorem~\ref{theorem:inert}).

\end{itemize}

\section{Eisenstein classes}   \label{section:Eisenstein}
Here, motivated by understanding the Galois story corresponding to~\S~\ref{sss:rreducible},
we analyze deformations of residual representations over an imaginary quadratic field. 

   In~\cite{MazEis}, Mazur undertakes a detailed
geometric study of the Eisenstein Ideal over $\Q$, and its ramifications
for the arithmetic of $J_0(N) = \mathrm{Jac}(X_0(N))$ (for $N$ prime). One can also study the
Eisenstein ideal from the perspective of Galois deformations, and this was
was the approach of~\cite{CE1}; we will generalize this analysis to an imaginary quadratic field here.  

For this section, the Eisenstein ideal will be generated
by $T_{\ell} - 1 - N(\ell)$ for all $\ell$ co-prime to the level; these correspond
to cyclotomic-Eisenstein classes. 
Throughout this section, $\q$ will denote a prime ideal of $\OL_F$, and $p$ an odd prime
in $\Z$. We also suppose that $\q$ is of residue characteristic different from $p$.

\medskip

For odd $p$, let
$$\rhobar = \left(\begin{matrix} \eps & 0 \\ 0 & 1 \end{matrix} \right) \mod p,$$
where $\eps$ is the %
 the cyclotomic character.
Since $\rhobar$ admits non-scalar endomorphisms, it does \emph{not} admit
a universal deformation ring. However, the key idea of~\cite{CE1} is the following.
Let $\Vbar$ denote the underlying  $\F_p$-vector space attached to $\rhobar$. Let
$\Lbar$ denote a line in $\Vbar$ that is \emph{not} invariant under $G_F$.
Fix once and for all a choice of inertia group $I_{\q}$ at $\q$. We consider, for
every local  Artinian ring $(A,\m, A/\m = \F_p)$, the set of triples $(V,L,\rho)$ where 
$\rho: G_F \rightarrow \GL(V)$, such that (cf.~\cite{CE1}, p.100):
\begin{enumerate}
\item The triple $(V,L,\rho)$ is a deformation of $(\Vbar,\Lbar,\rhobar)$.
\item The representation $\rho$ is unramified away from $p$ and $\q$, and
is finite flat at $v|p$.
\item The inertia subgroup at $\q$ acts trivially on $L$.
\item The determinant $\rho$ is the cyclotomic character (composed with the 
map $\Z_p \rightarrow A$). 
\end{enumerate}
As in~\cite{CE1}, we have:
\begin{lemma} Let $\Def(A)$ denote the functor that assigns to a local  Artinian
ring $(A,\m,\F_p)$  the collection of such deformations up to strict equivalence. Then
$\Def$ is pro-representable
by a complete Noetherian local $\Z_p$-algebra $R$.
\end{lemma}

Since $\rhobar$ has an obvious deformation to $\Z_p$, it follows that there
is always a map $R \rightarrow \Z_p$. It is then natural to ask: when is $R$
bigger than $\Z_p$? Since $R$ is a $\Z_p$-algebra, this is equivalent to asking
when $\Hom(R,\F_p[x]/x^2)$ is nontrivial, or equivalently when there
exists a nontrivial deformation of $\rhobar$ to $\F_p[x]/x^2$.
Over $\Q$, this is equivalent to asking whether the Eisenstein ideal $\m$ has
cuspidal support (see~\cite{CE1}).
We conjecture that

\begin{conj}
  $\Hom(R, \F_p[x]/x^2)$ is nontrivial iff
$H^E_1(\{\q\},  \Z_p)_{\m} \neq 0$.
\end{conj}

Note that $H_1(\{\q\}, \Z)$ does have congruence cohomology
of degree $N(\q) - 1$. Thus the question is related to whether  the localization
of $H_1(\{\q\},\Z)$ at $\m$ is bigger than the congruence cohomology.

\medskip

Computing $\Def(\F_p[x]/x^2)$ is closely related to computing
$\Ext^1(\Vbar,\Vbar)$ in the category of finite flat group schemes over $\OL_F[1/\q]$.
We make the following assumtions of $F$.
\begin{enumerate}
\item $p$ does not divide $\CL(F)$.
\item $p$ does not divide $w_F d_F$.
\end{enumerate}

The last condition implies that $1 = e_p(F_v/\Q_p) < p -1$ for any $v|p$. 
Under this condition,
that the category of finite flat group schemes over $\OL_F[1/\q]$ is abelian,
and the forgetful functor that takes a finite flat group scheme to its
generic fibre (i.e., considers the   underlying Galois
representation) is fully faithful, by~\cite{Raynaud}.

\medskip
\begin{theorem}  \label{defnonzero} %
 $\Def(\F_p[x]/x^2) \neq 0$
when $N(\q) \equiv 1$ modulo $p$ if $p > 3$ or
$q = N(\q) \equiv 1 \mod 9$ if $p = 3$ under either of the following assumptions: 
\begin{itemize}
\item[(i).] The $\eps^{-1}$-part of the class group\footnote{By the ``$\eps^{-1}$-part'', we mean
those $p$-torsion classes in $\CL(F(\zeta_p))$ for which $x^{\sigma} = \eps(\sigma)^{-1} x$,
for every $\sigma \in \Gal(F(\zeta_p)/F)$.}  of $\CL(F(\zeta_p))$ is nonzero, or,  \item[(ii).]    \label{theorem:inert} $\q$ is inert, generated by some rational prime $N$. 
\end{itemize}
\end{theorem}

{\em Remark.} A similar argument to what we detail below also
 works over $\Q$, and shows that a nontrivial
extension exists (if $p$ is odd) when $p$ divides the numerator of $(N-1)/12$.
When $p = 2$ one needs to consider a slightly different $\rhobar$
(see~\cite{CE1}) in order for the deformation ring to exist.

\subsection{Proof of Theorem~\ref{defnonzero} } \label{subsec:proofdefnonzero}

Writing $\Vbar = \Z/p\Z \oplus \mu_p$, we see that
$\Ext^1(\Vbar,\Vbar)$ (in the category of finite flat group schemes
over $\OL[1/\q]$ that are annihilated by $p$) admits maps to each of the extension groups
$\Ext^1(\Z/p\Z,\Z/p\Z)$,  $\Ext^1(\mu_p,\mu_p)$, 
$\Ext^1(\Z/p\Z,\mu_p)$, and the group  $\Ext^1(\mu_p,\Z/p\Z)$. Let us consider each of these extension
groups   in turn.

\begin{itemize}
\item[(i).] 
$\Ext^1(\Z/p\Z,\mu_p)$. By Kummer theory, the corresponding group
for extensions of {\em Galois modules} is isomorphic to
 $F^{\times}/F^{{\times}p}$. In order for this to arise from a finite flat group scheme
 over $\OL[1/\q]$ the valuation  of a
 class in $F^{\times}/F^{{\times}p}$ must be  zero modulo $p$
 at every prime except $\q$. 
 Thus 
 $$\Ext^1(\Z/p\Z,\mu_p) \simeq \OL_F[1/\q]^{\times}/\OL_F[1/\q]^{{\times}p},$$
where we used the fact that $p \nmid h_F$ to
 make this identification.
This is always non-zero since the group $\OL_F[1/\q]^{\times}$ of $\q$-units has rank $1$.

\item[(ii).]
$ \Ext^1(\mu_p,\mu_p) = \Ext^1(\Z/p\Z, \Z/p\Z)$. 
By the connected-\etale sequence, all
such extensions are \etale and hence unramified at $p$.  This group is therefore non-zero
exactly when
there exists an abelian extension of $F$ of degree $p$ unramified everywhere except for $\q$.
Since we are assuming that $\Cl(F)$ is not divisible by $p$, the $p$-part of
the  ray class group is the maximal $p$-quotient of the cokernel of the map:
$$\OL^{\times}_F \rightarrow (\OL_F/\q)^{\times}.$$
 But $F$ is an imaginary quadratic field and $p$ does not divide $w_F$, and
 hence this is nontrivial exactly when $p$ divides $N(\q) - 1$.
 
 \item[(iii).]  $\Ext^1(\mu_p,\Z/p\Z)$.  By the connected-\etale sequence, 
all such extensions are totally split at $v|p$.  

In particular, the associated extension of Galois modules is unramified
everywhere away from $\q$ after restriction to  $F(\zeta_p)$. Conversely, 
consider any extension $M$ of $\mu_p$ by $\Z/p\Z$  (as Galois modules) with this property.
Let $v$ be a prime of $F$ above $p$, and let $H_v$ be the splitting field of $M$ over $F_v$.
By assumption, the extension $H_v/F_v(\zeta_p)$ is unramified. 
Since $p$ does not divide $d_F$, the extension $F_v(\zeta_p)/F_v$ is
totally ramified. In particular, the inertia subgroup $I$ of $\Gal(H_v/F_v)$
is a normal subgroup which 
maps isomorphically onto $\Gal(F_v(\zeta_p)/F_v)$. This implies that 
$\Gal(H_v/F_v)$ splits as a direct product $\Gal(F_v(\zeta_p)/F_v) \times
\Gal(H_v/F_v(\zeta_p))$. Yet this is incompatible with the structure of $\Gal(H_v/F_v)$
arising from its definition as an extension unless $H_v = F_v(\zeta_p)$. In particular, any such extension
$M$ splits completely at $p$,
and hence  $M$   prolongs to a finite flat
group scheme over $\OO_L[1/\q]$. 
 
 Write $L = F(\zeta_p)$. 
If we are in case (i) --- i.e. $\CL(L)[\eps^{-1}] \neq 0$ --- then 
there exists extensions of $\mu_p$ by $\Z/p\Z$ unramified at
\emph{all} primes, and hence  extensions unramified away 
from $\q$ certainly exist. 
On the other hand, we verify in~\S~\ref{ext1} that, if $\q$ is inert and $N(\q) = N^2$ where
$N \equiv -1$ modulo $p$,  the  cokernel of \begin{equation} \label{oleps} \left( \OL^{\times}_L \otimes_\Z \mathbf{F}_p\right) [\eps^{-1}]  \rightarrow \left( (\OL_L/\q)^{\times}\otimes_\Z \mathbf{F}_p\right)[{\eps^{-1}}]\end{equation} 
is always nontrivial, 
where $\eps$ is the cyclotomic character.  Therefore, in either case (i) or (ii), 
 $\Ext^1(\mu_p,\Z/p\Z) \neq 0$. 
 \end{itemize}

Therefore, in either case (i) or (ii) when $\q = (N)$ with $N \equiv -1$,   all four $\mathrm{Ext}$-groups are  nontrivial, and
one may construct a nontrivial deformation with a suitable line $L$
as in Proposition 5.5 of~\cite{CE1}.  In the case when $p | (N-1)/12$  and  $N \equiv 1$ modulo $p$ is inert, then a nontrivial deformation is constructed by Mazur in~\cite{MazEis}.
Specifically, Mazur constructs a modular abelian variety $A/\Q$  (the Eisenstein quotient)
and a maximal
ideal $\m$ of $\T$ such that $A[\m] = \rhobar$; this implies the existence of a
characteristic zero deformation of $\rhobar$.
This completes the proof
of  the Theorem. 

We record the following Lemma for later use: 
 
 \begin{lemma} \label{lemma:whenexists}
If  $\Cl(L)[\eps^{-1}] = 0$,
then
 $\Def(\F_p[x]/x^2)$ is nontrivial 
{\em exactly}
exactly when $\Ext^1(\mu_p,\Z/p\Z)$ is nontrivial in the category of finite flat group schemes over $\OL_F[1/\q]$, and $N(\q) \equiv 1 \mod p$.
\end{lemma}

\begin{proof}
The computations above show that, under the assumption on $\CL(L)$,  that any nontrivial deformation is ramified at $\q$.
If the deformation is an upper-triangular or lower-triangular definition, then the only
$I_{\q}$-invariant lines will be $G_F$ invariant upon projection to $\Vbar$, contradicting
the assumption on $\Lbar$. Hence, both $\Ext^1(\mu_p,\Z/p\Z)$ and
$\Ext^1(\Z/p\Z,\mu_p)$ must be nontrivial. 
By a similar argument, if the extension group $\Ext^1(\Z/p\Z,\Z/p\Z)$ vanishes,
then the only $I_{\q}$ invariant line will be $G_F$ invariant under
projection (See Lemma~5.2 and the proof of Proposition 5.5
of~\cite{CE1}.)
Hence $N(\q) \equiv 1 \mod p$. 
\end{proof}

 \subsection{Further analysis of  \texorpdfstring{$\Ext^1(\mu_p,\Z/p\Z)$}{Ext}} \label{ext1}
 
We analyze more carefully 
the cokernel of 
~\eqref{oleps}, where $L = F(\zeta_p)$:
$$ \left( \OL^{\times}_L \otimes_\Z \mathbf{F}_p\right) [\eps^{-1}]  \rightarrow \left( (\OL_L/\q)^{\times}\otimes_\Z \mathbf{F}_p\right)[{\eps^{-1}}].$$ 
 In particular, we show that if $\q$ is generated by a rational (inert) prime $N \equiv -1$ modulo $p$, then this cokernel is nonzero. 

Let us first consider the analogous story over $\Q$ (so imagine that $L = \Q(\zeta_p)$), where
the ideal ${\q}$ is generated by a prime number $q$. We will show
that the cokernel is nonzero when $q \equiv \pm 1$ mod $p$ for $p>3$;
for $p=3$ the condition becomes 
  $q \equiv \pm 1$ mod $9$:
  
  All the units of $L$ are products of  real  units and elements of $\mu_p$. The real units
project trivially to the $\eps^{-1}$ eigenspace (if $p \ne 2$) since $\eps$ is odd.
The action of
$\Gal(L/\Q) = (\Z/p\Z)^{\times}$ on $\mu_p$ is via the character $\eps$, which is distinct
from $\eps^{-1}$ if $p$ is neither $2$ nor $3$. Finally, the $\eps^{-1}$-eigenspace on
the right-hand group of~\eqref{oleps} is nontrivial if and only if 
$p$ divides $q^2 - 1$, as follows
from the proof of Lemma~3.9 in~\cite{CE1}.
 When $p = 3$,  $\eps^{-1} = \eps$. Yet
$\eps$-extensions ramified only at $q$ are of the form
$\Q(\zeta_3,\sqrt[3]{q})$. This extension (over $\Q(\zeta_3)$) splits at $3$
only when  $q \equiv \pm 1 \mod 9$.

If $F$ is an imaginary quadratic field, then
 $\OL^{\times}_L \otimes \Q$ is the regular representation
of $(\Z/p\Z)^{\times}$. Hence, determining the cokernel of~\eqref{oleps}
is more subtle.
We begin with a preliminary remark
concerning the case $p = 3$ which is slightly different to
the general case (since $\eps = \eps^{-1}$ in this case). The $\eps^{-1} = \eps$
eigenspace of the unit group is generated by some unit $\beta$ 
\emph{together} with $\zeta_3$. In particular, since $(\zeta_3 - 1)$ is
prime to $\q$ (it has norm a power of $3$),
the image of $\zeta_3$  in $(\OL_L/\q)^{\times} \otimes \F_3[\eps^{-1}]$ is trivial if
$q \equiv 1 \mod 9$, but consists
of the entire group if $q \not\equiv 1 \mod 9$.
Hence, for $p = 3$, we make the additional assumption that
$q \equiv 1 \mod 9$.

 We return to the analysis of $\q = (N)$, where $N \equiv -1$ modulo $p$,
 and $q = N^2 \equiv 1 \mod 9$ if $p = 3$.
  We need to check that the cokernel of 
 $$\OL^{\times}_L/\OL^{\times p}_L(\eps^{-1})
  \rightarrow (\OL_L/N)^{\times}(\eps^{-1}) =  \left(\prod (\OL_L/\PP_L)^{\times}\right)(\eps^{-1})$$
  is divisible by $p$, where the product runs over all primes $\PP_L$ dividing $N$.
Note that the first group is cyclic and generated by some
 unit $\beta$. (When  $p = 3$, it is generated by
 $\beta$ and $\zeta_3$, but the image of $\zeta_3$ is trivial
 since $N^2 \equiv 1 \mod 9$).
  We claim that: 
 \begin{enumerate}
 \item $N_{L/\Q}(\PP_L)  - 1 = N^2 - 1$ is divisible by $p$, for every $\PP_L$ dividing $N$, and the $\eps^{-1}$-invariants
 of $ (\OL_L/N)^{\times}$ have order divisible by $p$.
 \item We may take $\beta$ to lie in a proper subfield $E$, and (thus) the map above factors
 through $\prod (\OL_E/\PP_E)^{\times}$, where $\PP_E = \PP_L \cap \OL_E$.
 \item For every such $\PP_E$, the quantity $N_{E/\Q}(\PP_E) - 1$ equals $N-1$; in particular, this is not divisible by $p$.
 \end{enumerate}
 These suffice to prove the result.
 
 \medskip
 
 Let   $G =\Gal(L/\Q) = (\Z/p\Z)^{\times} \oplus \Z/2\Z$. By assumption ($N$ is inert
 and $N \equiv -1 \mod p$), $ \Frob({\PP_L}) = (-1,1) \in G$.   In particular, $N_{L/\Q}(\PP_L) = N^2$.
 Since the action of $\Gal(L/F)$ on $(\OL_L/\PP_L)^{\times}[p]$ is  via the regular
 representation (as $N$ splits completely from $F$ to $L$), this proves the first claim.
 
Let $\eta$ denote the nontrivial character of $\Z/2\Z = \Gal(K/\Q)$.
This group acts on
$\OL^{\times}_L$ rationally and it decomposes  (at least after inverting $2$) into trivial
and $\eta = -1$ eigenspaces. The trivial eigenspace corresponds to
the units of $\Q(\zeta_p)$. Since these are all the product of a totally real unit
with an element of $\zeta_p$, their projection
to the $\eps^{-1}$-eigenspace in $p$-torsion is trivial. Hence the $\eps^{-1}$-eigenspace is
generated by the projection of the $\eta$-eigenspace.
Hence, we may assume that $\Z/2\Z$ acts on $\beta$ via $\eta$.
Thus $G$ acts on $\beta$ via $(\epsilon^{-1},\eta)$. In particular, if $H = (-1,1)$,
$H$ fixes $\beta$. Let $E = L^{H}$. By construction, $\beta$ lies in $E$, and
we have established part two above.
  Since $\Frob({\PP_L}) = (-1,1) \in H$, it
follows that $\Frob({\PP_E})$ is trivial. Hence $N_{E/\Q}(\PP_E) = N$, and we are done.

\section{Phantom classes} \label{section:phantomclasses}
 Suppose that $\q$ is prime, and that there exists a nontrivial extension class $[c] \in \Ext^1(\mu_p,\Z/p\Z)$ (in the category of finite flat group schemes
 over $\OL_F[\frac{1}{\q}]$).
 If $N(\q) = q \equiv 1 \mod p$, then we have seen -- from analyzing the argument of
~\S~\ref{subsec:proofdefnonzero}  -- that the cylotomic-Eisenstein representation
 $\rhobar = \eps \oplus 1 \mod p$ admits a nontrivial non-reducible deformation (of the kind considered above).
  If $q \not\equiv 1 \mod p$, then such a deformation does not exist.

 At non-prime level, however, the issue of Eisenstein deformations becomes
  quite complicated, even over $\Q$. One can no longer write
down a deformation problem that depends on lifting more than one ``invariant line'',
since the resulting functor will not be representable.

 In light of the numerical evidence of
Chapter~\ref{chapter:ch8}, it is interesting to consider what happens 
whenever there exists such a class when
$q \not\equiv 1 \mod p$.

\begin{conj}  \label{ThePhantom}Let $\q$ be a prime ideal with $N(\q) =  q \not\equiv 1 \mod p$.  Suppose 
that $\Ext^1(\mu_p,\Z/p\Z) \neq 0$, in the category of finite flat group schemes
over $\Spec \ \OL_F[\frac{1}{\q}]$, but the corresponding $\Ext$-group over $\Spec \ \OL_F$ is trivial.
 \label{conj:phantom} 
 
 Then, for every auxiliary prime ideal $\mathfrak{r}$, there exists
 an cyclotomic-Eisenstein maximal ideal $\m$ of the Hecke algebra
 $\T_{\q \mathfrak{r}}$, of residue characteristic $p$, but there is no such
 ideal of $\T_{\q}$. 
\end{conj}
  We call the cohomology classes at level $\q \mathfrak{r}$ annihilated
  by this ideal phantom classes.

In order to motivate this conjecture, let us consider the ``degenerate case'' when $\q = \OL_F$.
(It may seem that this does not fit into the framework of the Conjecture, 
since  $N(\OL_F) = 1 \equiv 1 \mod p$. However, 
that condition was imposed in the conjecture 
in order to ensure 
 that the group $\Ext^1(\Z/p\Z,\mu_p)$ over $\Spec \OL_F[\frac{1}{\q}]$
 was trivial; and the corresponding group of everywhere unramified
 extensions {\em is} trivial  providing that $p \nmid w_F \cdot h_F$.) Then
we might expect that the existence of $[c]$ implies that there exists 
a phantom class of level $\rf$ for any auxiliary prime
$\rf$. On the other hand, the group subgroup of $\Ext^1(\mu_p,\Z/p\Z)$ of extensions unramified
everywhere over $\OL_F$ is (under our assumptions) the same as the $p$-torsion in 
$K_2(\OL_F)$
(see~\cite{BeG}, Theorem~6.2).
Hence, Theorem~\ref{theorem:K2popularversion}(ii) exactly predicts the existence of such 
phantom classes. (In fact, Theorem~\ref{theorem:K2popularversion} was
motivated by the numerical computations that
led to Conjecture~\ref{conj:phantom}.)
It follows that one way to prove the existence of these
more general phantom classes as above is to extend
Theorem~\ref{theorem:K2popularversion}(ii) to nontrivial level. 

\label{section:morephantom}
Here is a rough outline of
what we believe may be happening, although we stress that we are unable to prove
so much at this point.
Suppose that in
 Theorem~\ref{theorem:K2popularversion}(ii), we replace $Y$ by the congruence cover 
 $Y_0(\q)$. Then, for any auxiliary prime $\rf$, we expect that the homology at level $\q \rf$
 will detect classes in a group
 $K_2(\OL_F,\q)$
 which contains $K_2(\OL_F)$. What should the definition of this group be? At the moment, we do
 not have a candidate definition, but let us at least discuss some
 of the properties of this group. The usual $K$-group $K_2(\OL_F)$ admits a second
 \etale Chern class map 
 $$K_2(\OL_F) \rightarrow H^2(F,\Z_p(2)).$$
 A theorem of Tate~\cite{Tate} implies that this map is injective, and
 the image of this map will consist of classes unramified everywhere. 
 For a finite set of primes $S$, let
 $$H^1_S(F,\Q_p/\Z_p(-1)):= \ker H^1(F,\Q_p/\Z_p(-1)) \rightarrow
 \bigoplus_{v \notin S} H^1(F,\Q_p/\Z_p(-1)).$$
 Recall the Poitou--Tate sequence:
$$
\begin{diagram}
 &  & \bigoplus H^0(F,\Q/\Z_p(-1)) & \rTo & H^2(F,\Z_p(2))^{\vee}  & \rTo & 
H^1(F,\Q_p/\Z_p(-1)) &   \\
&  \rTo & \bigoplus H^1(F_v,\Q_p/\Z_p(-1))  & \rTo & 
H^1(F,\Z_p(2))^{\vee} & & \\
\end{diagram}
$$
Hence, if $H^2_S(\Z_p(2)):=H^1_S(F,\Q_p/\Z_p(-1))^{\vee}$, there is an exact sequence:
$$
\begin{diagram}
\bigoplus_S H^1(F_v,\Z_p(2)) & \rTo & H^2_S(F,\Z_p(2)) & \rTo & H^2(F,\Z_p(2)) & \rTo &
\bigoplus H^2(F_v,\Z_p(2))
\end{diagram}
$$
We suspect that $K_2(\OL_F,\q)$ (when correctly defined) 
admits an \etale Chern class map to map to $H^2_S(\Z_p(2))$ with
$S = \{ \p\}$. (More generally, we expect that, for any set of primes $S$,  there will
be a group $K_2(\OL_F, S)$ which admits such a map
to $H^2_S(\Z_p(2))$.)
In particular, there exists a non-zero class in $K_2(\OL_F,\q)$ mapping to zero in $K_2(\OL_F)$
exactly when $H^1_{\q}(F,\Q_p/\Z_p(-1)) \ne H^1_{\emptyset}(F,\Q_p/\Z_p(-1))$, or equivalently,
exactly when there exists a class $[c]$ in $\Ext^1(\mu_p,\Z/p\Z)$ which is unramified everywhere
outside $\q$. In particular, the analogue  of Theorem~\ref{theorem:K2popularversion}(ii)  would 
imply that $[c]$ gives rise to classes at level $\q \rf$ for any auxiliary $\rf$. This analysis
does not obviously predict when such a class arises at level $\q$, or when, alternatively, it is a phantom class;
for this, we use the heuristic of mod-$p$ deformations at level $\q$ to conjecture that it
is a phantom class exactly when $q \not\equiv 1 \mod p$.

\begin{remarkable} \label{remark:failure}
\emph{
We also speculate here on the connection between phantom classes  and the partial
Euler characteristic. Suppose that the conditions of Conjecture~\ref{ThePhantom}
are satisfied, and consider any set $S = \{\p,\rf\}$
of two primes different from $\q$. Then in the level lowering sequence:
$$H_1(\Gamma_0(\q),\Z)^4 \longleftarrow  H_1(\Gamma_0(\p \q),\Z)^2 \oplus
H_1(\Gamma_0(\p \rf),\Z)^2 \leftarrow H_1(\Gamma_0(\p \q \rf),\Z),$$
we find (experimentally) that the middle term of this sequence is ``too large'', 
because the phantom classes at level $\p \q$ and $\p \rf$ must map to 
zero at level $\p$, forcing the homology at the middle term to be non-zero (again, this
is an experimental observation, not an argument).
 It seems, therefore, that
the existence of phantom classes will force the level raising sequence to fail to be exact,
and thus (in the spectral sequence) give rise to terms which contribute to
the partial Euler characteristic $\chi_D(Y)$
via $H^{E}_2(K_0(\q)[1/\p \rf],\Z)$. 
}
\end{remarkable}

For numerical examples of phantom classes, see the
end of~\S~\ref{subsection:pathologies}, 
Remark~\ref{remark:hangingchad},
and~\S~\ref{subsec:chad2}.

 \subsection{Phantom classes over $\Q$}

One may ask whether there is an analogue of phantom classes for $F = \Q$.
In the category of finite flat group schemes over $\Spec(\Z)$, one has  $\Ext^1(\mu_p,\Z/\p\Z) = 0$ for all $p$
(as follows from~\cite{MazEis} --- the point being that the $\eps^{-1}$-part of the class group
of $\Q(\zeta_p)$ is related to the numerator of the Bernoulli number $B_2 = 1/6$).
The extension
group $\Ext^1(\mu_p,\Z/\p\Z)$ of finite flat group schemes over $\Spec(\Z[1/q])$ is non-trivial exactly when
$q \equiv \pm 1 \mod p$ (this follows from the remarks made in~\S~\ref{ext1}). Assume that $q \ge 5$.
If $q \equiv 1 \mod p$,
then one always has modular irreducible characteristic zero deformations
of $1 \oplus \epsilon$  coming from the Eisenstein quotient~\cite{MazEis}.  On the other hand, if
$q \equiv -1 \mod p$, then $1 \oplus \epsilon$ does \emph{not} occur as the semi-simplification of a cuspidal
representation of weight two and level $\Gamma_0(q)$, as also follows from~\cite{MazEis}.
In analogy with the  discussion above,  one expects that for all $q \equiv -1 \mod p$, there exist phantom cyclotomic Eisenstein
classes of  level $qr$ for any auxiliary prime $r$, but not at level $q$ itself.
This is indeed the case, via the following easy lemma:

\begin{lemma} Let $p \ge 5$ be prime, and fix a prime $q \equiv -1 \mod p$. 
If $r$ is any auxiliary prime different form $q$,
then there exists a characteristic zero cuspidal eigenform $\pi$ of level $\Gamma_0(qr)$  such that:
\begin{enumerate}
\item $\pi$ is new at $r$.
\item $\pi \mod p$ is cyclotomic Eisenstein.
\item If $r \ne p$ and $r \not\equiv \pm 1 \mod p$, there exists a lattice for the corresponding $p$-adic Galois
representation $\rho(\pi): G_{\Q} \rightarrow \GL_2(\Qbar_p)$ such that
$\rhobar$ represents the unique class in $\Ext^1(\mu_p,\Z/p\Z)$ which is ramified at $q$ and 
unramified everywhere else.
\end{enumerate}
\end{lemma}

\begin{proof} Let $D$ denote the quaternion algebra ramified exactly at~$q$ and~$r$. 
Let $\Gamma_B$ denotes the
corresponding maximal arithmetic subgroup.
There is a congruence quotient map 
$$H_1(\Gamma_B,\Z)
\rightarrow \F^{\sigma = -1}_{q^2} = \Z/(q+1)\Z.$$
 Since $p \ge 5$ it is not an orbifold prime,  this class lifts to
a characteristic zero eigenform $\pi^{\JL}$ which admits  a corresponding $p$-adic Galois representation $\rho(\pi^{\JL})$ which is cyclotomic Eisenstein modulo~$p$.
 By Jacquet--Langlands, this may be transferred to a form $\pi$ of level $\Gamma_0(qr)$ with the same Galois representation
  which is also new at $q$ and $r$. 
 The $p$-adic Galois representation $r(\pi): G_{\Q} \rightarrow \GL_2(\Qbar_p)$ will be ordinary at $p$ because it is reducible modulo $p$ (this is also true if $r = p$ for a slightly different reason).
 By Ribet's Lemma~\cite{ribetlemma},  there exists a lattice such that
 $$\rhobar = \left( \begin{matrix} 1& * \\ 0 & \epsilon \end{matrix} \right)$$
 where $*$ is non-trivial. Since $p > 3$, comparing this representation
 with the lattice arising from the fact that $\pi$ is ordinary at $p$, we deduce that
 $*$ is unramified at $p$. Hence the representation
 defines a non-trivial extension class in $\Ext^1(\mu_p,\Z/p\Z)$  which is finite flat over $\Spec(\Z[1/qr])$.
  If $r \ne p$ and $r \not\equiv \pm 1 \mod p$, then such extensions are automatically
 unramified at $r$ for local reasons.
 \end{proof}
 
 \begin{remark}One cannot expect the final claim (regarding $\rhobar$) to be true if $r \equiv \pm 1 \mod p$ for symmetry
 reasons (replacing $q$ by $r$ would lead to a contradiction if $r \equiv -1 \mod p$, for example).
 However, in this case, one would expect the corresponding extension at least to occur inside the universal deformation ring, and thus
still be ``modular'' in some sense.
 \end{remark}

\section{Even Representations} \label{section:even}

Here we analyze further the example of~\S~\ref{sss:reven}:

Fix a prime number $p > 2$, and a finite field $k$ of characteristic $p$.
Let 
 $$\rhobar: G_{\Q} \rightarrow \GL_2(k)$$
 be a continuous absolutely irreducible Galois representation.
 We shall assume that $\rhobar$ is either finite flat or ordinary at $p$.
  One of the original motivations of considering torsion homology 
 for $\GL(2)$ over imaginary quadratic fields is precisely that the 
 restricted representation $\rhobar| G_F$ conjecturally contributes to the integral
 cohomology of $Y(K_{\Sigma})$, whether or not $\rhobar$ is even.

  One expects even Galois  representations to have very few deformations
over $\Q$, even without placing any local conditions at $p$. 
(Indeed, by work of the first author, one know nows in many cases that
$\rhobar$ admits no  potentially semi-stable characteristic zero
deformations~\cite{CalEven1, CalEven2}.)

We show here that over $F$ they have \emph{more} deformations with nice
properties than one might \emph{a priori} expect, as is predicted by our conjectures together with
 Corollary~\ref{cor:double} (see discussion of~\S~\ref{sss:reven}).

\subsection{Deformation Functors Revisited}
Let $\rhobar: G_{\Q} \rightarrow \GL_2(k)$ be an absolutely
irreducible continuous representation, and suppose that either:
 \begin{enumerate}
\item $\rhobar$ is finite flat at $p$,
\item $\rhobar$ is ordinary at $p$.
\end{enumerate}
Suppose, moreover,  that $\det(\rho) =  {\eps \chi} \mod p$, where $\eps$  is the $p$-adic cyclotomic
character and $\chi$ is a finite order character that is trivial on $G_{\Q_p}$. 
Let $F/\Q$ be an imaginary quadratic field disjoint from the fixed field of the kernel of
$\rhobar$, and in which $p$ splits. %
Associated to $\rhobar$ one may define the usual classes of deformation problems, namely,
for
local  Artinian rings $A = (A,\m_A,k)$;
\begin{enumerate}
\item $\Dfl(A)$ denotes equivalence classes of  lifts $\rho: G_F \rightarrow \GL_2(A)$ which are finite flat at  both $v|p$,
minimally ramified at all other places, and of determinant $\det(\rho) = \eps \chi$.
\item $\Dord(A)$ denotes equivalence classes of lifts $\rho: G_F \rightarrow \GL_2(A)$ which are ordinary at both $v|p$,
minimally ramified at all other places, and of determinant $\det(\rho) = \eps \chi$.
\end{enumerate}
 
\begin{theorem} Suppose that $\rhobar$ is even.
If $\rhobar$ is finite flat at $p$, then the tangent space
$\Dfl(k[x]/x^2)$ is nontrivial. \label{theorem:evendoubling} If
$\rhobar$ is ordinary at $p$, then $\Dord(k[x]/x^2)$ is nontrivial.
\end{theorem}

\begin{proof}

As in~\cite{CM}, we may identity $\Dord(k[x]/x^2)$  and $\Dfl(k[x]/x^2)$ %
with a Selmer group inside the Galois cohomology group $H^1(F,W)$;
here $W$ is the trace-free elements of $\End(\bar{\rho}, \bar{\rho})$.

Let $c \in \Gal(F/\Q)$ be complex conjugation, and let $\eta$ be the
nontrivial character of $\Gal(F/\Q)$.   Then
$$H^1(F,W) = H^1(\Q,\Ind^{\Q}_F W) \simeq H^1(\Q,W^+) \oplus H^1(\Q,W^-),$$
where $W^+ = \mathrm{End}^0(\bar{\rho}, \bar{\rho})$, and $W^- \cong W^+ \otimes \eta$.

Let $H^1_{\Sigma}(F,W)$ denote the Selmer group which captures deformations that are 
minimally unramified outside $p$. Let $H^1_{\Sigma}(\Q, W^+)$ 
be the preimage of $H^1_{\Sigma}(F, W)$ in $H^1(\Q, W^+)$ under the restriction map, and similarly for $W^+$ replaced by $W^-$.  Then the Euler characteristic formula implies that
$$
\begin{aligned}
\dim H^1_{\Sigma}(\Q,W^+) \ge  0 & \ \  \text{and}  &
\dim H^1_{\Sigma}(\Q,W^-) \ge 3,& \ \text{if $\rhobar$ is \emph{even}},\\
\dim H^1_{\Sigma}(\Q,W^+) \ge  2 & \ \   \text{and}  &
\dim H^1_{\Sigma}(\Q,W^-) \ge 1,& \ \text{if $\rhobar$ is \emph{odd}}.
\end{aligned}
$$
Let $\Hfl$ and $\Hord$ inside $H^1_{\Sigma}$ denote the subspace of classes satisfying the
restricted deformation condition above. Let
$\Hflv$ and $\Hordv$ inside $H^1_{\Sigma}$ denote the subspace of classes satisfying
the restricted deformation condition at \emph{one} prime above $p$ in $F$.

\begin{lemma} Let $* = \mathrm{fl}$ or $\mathrm{ord}$. 
Let $W^+ = W$ and $W^- = W \otimes \eta$.  We have 
an equality 
$$H^1_*(F,W) \cap H^1_{\Sigma}(\Q,W^{\pm}) = H^1_{*,v}(F,W) \cap H^1_{\Sigma}(\Q,W^{\pm}).$$
Equivalently, a class in $H^1_{\Sigma}(F,W)$ that is an eigenvector for $c$ satisfies the necessary local condition
at one $v|p$ if and only if it satisfies the local condition at the other. 
\end{lemma}

\begin{proof}
This can be
seen by considering the action of the group 
$\Gal(F/\Q)$ on $H^1(F,W)$ --- it acts via the
identity on $H^1(\Q,W)$ and by $-1$ on $H^1(\Q,W \otimes \eta)$. On the other hand,
it permutes the local classes at the primes $v|p$.
Thus any class in $H^1_{\Sigma}(\Q,W^{\pm})$ has isomorphic local
classes for each $v|p$.
\end{proof}

Generically, we expect the imposition of a condition $*$ at  $v|p$   imposes two local conditions.
Thus we expect to be able to produce a nontrivial element in  $H^1_{*,v}(F,W) \cap H^1_{\Sigma}(\Q,W^{\pm})$ \emph{if and only if}  the dimension of this space is at least $3$,
that is, if $\rhobar$ is even and we consider $W^{-}$. 
Let us now assume that $\rhobar$ is even.
To complete the proof of the theorem, we require a more accurate computation
of $H^1_{\Sigma}(\Q,W^-)$, and a more accurate count of the local conditions that are imposed.
 Since the Cartier dual of $W^-$ is totally real, by the global
Euler characteristic
formula (Theorem~2.8 of~\cite{Fermat}), we find that
$$\dim H^1_{\Sigma}(\Q,W^-) \ge 3 + \dim H^0(\Q_p,(W^-)^*),$$
where $M*$ denotes Cartier dual of $M$. Since, as a $\Q_p$-representation, $W^- \simeq W$, we may
replace the last term by $H^0(\Q_p,W^*)$. By Tate local duality, we may replace this
by $H^2(\Q_p,W)$. Finally, using the local Euler characteristic formula, we find that
$$\dim H^1_{\Sigma}(\Q,W^-) \ge \dim H^1(\Q_p,W) - \dim H^0(\Q_p,W).$$
(Indeed, we may have derived this directly from the global Euler characteristic formula, but
in the first formulation it was more apparent why this quantity is at least $3$.)
Consider the Cartesian square:
$$
\begin{diagram}
H^1_{*,v}(F,W) & \rTo & H^1_{\Sigma}(F,W) \\
\dTo & & \dTo \\
H^1_{*}(\Q_p,W) & \rTo & H^1(\Q_p,W) 
\end{diagram}$$
where the vertical arrows represent localizing at $v$. 
To prove that $H^1_{*,v}(F,W) \cap H^1_{\Sigma}(\Q,W^{-})$ is non-zero, it suffices to show that
the cokernel of the map to $H^1_{\Sigma}(\Q,W^{-}) \subset H^1_{\Sigma}(F,W)$ has dimension
strictly less than $H^1_{\Sigma}(\Q,W^{-})$ itself. Since the square is
Cartesian, it suffices to show that the kernel of the local map also has dimension
strictly less than $\dim H^1_{\Sigma}(\Q,W^{-})$. Yet, from Proposition~2.27(a) of~\cite{Fermat} (taking into
account Remark~2.26), the codimension of this map is
$$\dim H^1(\Q_p,W) - \dim H^0(\Q_p,W) - 1 < \dim H^1_{\Sigma}(\Q,W^-).$$

\end{proof}

\begin{remarkable}
\emph{
If $p$ is inert (or perhaps even ramified) in $F/\Q$ a similar result can probably be proved using the same techniques.
The key point in that case will be to compute the dimension of the cokernel of the map
$$H^1_{*}(L,W)^{-} \rightarrow H^1(L,W)^{-},$$
where $L/\Q_p$ is a quadratic extension and $-$ denotes the $-1$ eigenspace under the action of $\Gal(L/\Q_p)$.
}
\end{remarkable}

\begin{remarkable}
\emph{
We have excluded $p = 2$ in this theorem. When $p = 2$ the representation $\rhobar$ will actually
be modular in the usual sense. Then one might hope that there exist nontrivial deformations (flat or
ordinary) of $\rhobar$ even over $\Q$. If the Serre conductor of $\rhobar$ is equal to a prime $N$, and $\rhobar$ is
finite flat at $2$, then this is actually a theorem (see~\S~3.2 of~\cite{CE}), although we do not know whether to expect
this to be true in general.
}
\end{remarkable}

\section{Even representations: finiteness and experiments } \label{ss:Serre}

One interesting consequence of Serre's conjecture are finiteness
theorems for the number of odd Galois representations
 $$\rhobar: \Gal(\Qbar/\Q) \rightarrow \GL_2(\Fbar_p)$$
 of any fixed tame level $N$. 
 
There are natural generalizations of Serre's conjecture
to imaginary quadratic fields (see~\cite{Torrey}, as well
 as~\cite{ClozelA,ClozelB}); in this context,
there is no ``odd'' condition. 
 We remark here that these generalizations  have similar, and in fact stronger, implications for even
 Galois representations.
 
 \begin{lemma} Assume Conjecture~\ref{conj:serre}(1)
 and Serre's modularity conjecture over $F$ (in the sense of~\cite[p3]{Torrey}). 
  Fix an integer $N$, and fix an integer $k \ge 2$.
   Then there only exist finitely many
 irreducible even Galois representations:
 $$\rhobar: \Gal(\Qbar/\Q) \rightarrow \GL_2(\Fbar_p)$$
 of tame level $N(\rhobar)$ dividing $N$ and
 Serre weight $k(\rhobar) = k$.
 For a fixed $p$, there are only finitely many $\rhobar$
 of tame level dividing $N$.
 \end{lemma}

\begin{remarkable} 
\emph{
Perhaps surprisingly, if we fix $k(\rhobar) \ge 2$,
this lemma implies that there are only finitely many $p$ for which
(conjecturally) an even $\rhobar$ can exist. This is in marked
contrast to the case of odd Galois representations ---
when $k(\rhobar) = 12$ and $N = 1$, the
representations $\rhobar_{\Delta}$ associated
to  $\Delta = q \prod_{n=1}^{\infty} (1  -q ^n)^{24}$
give examples of such $\rhobar$ for all
$p \notin \{2,3,5,7,691\}$. 
Note that the theorem is not true for $k = 1$, since
(for example) an even  representation with
image in $\GL_2(\C)$ gives rise to an
irreducible $\rhobar$ with $k(\rhobar) = 1$
for all but finitely many $p$.
}
\end{remarkable}

\begin{proof} 
Choose any auxiliary imaginary quadratic field $F$.
The number of $\rhobar$ such that
$\rhobar|G_{F}$ are reducible is easily seen to be finite
by class field theory (we use here that $k \ne 1$).
Assume that there exist infinitely many such
$\rhobar$. Assuming
Conjecture~\ref{conj:serre}, it follows that there exists
infinitely many maximal ideals $\m$ of the ring
of Hecke operators $\T$ acting faithfully on
$H_1(Y_1(N),\Z)$. Since the latter group is finite, only finitely
many of these $\m$ can be supported on torsion classes, and
by a standard pigeonhole principle argument (see~\S~4.6 of~\cite{Serre2}), it
follows that there exists a characteristic zero eigenclass
$[c] \in H_1(Y_1(N),\C)$ which gives rise to infinitely many
of these $\m$ of characteristic $> 2$.

 Let $\pi$ denote the corresponding automorphic 
form; let $\alpha_v$ be the $v$th Hecke eigenvalue of $\pi$,
and let $L$ be the field generated by all $\alpha_v$, with ring of integers
$\OL_L$.  Then for each maximal ideal $\m$ as above,
there is a prime $\lambda_{\m}$ of $\OL_L$
and an isomorphism $\T/\m \stackrel{\sim}{\rightarrow} \OL_L/\lambda_{\m}$
such that $T_v$ is carried to $\alpha_v$ mod $\lambda_{\m}$.

If $v$ is a place of $F$ not dividing $N$, we deduce 
(from the compatibility of $\m$ and  the
correpsonding $\rhobar$) that $\alpha_v \equiv \alpha_{v^c}
\mod \m$ for infinitely many $\m$, and thus that 
$\alpha_v = \alpha_{v^c}$. By multiplicity one for $\GL(2)/F$,
we deduce that $\pi \simeq \pi^c$. By cyclic base change
(\cite{AC}, Theorem~4.2 p.202), it follows
that $\pi$ comes via base change from a $\varpi$ for $\GL(2)/\Q$.
But the Galois representations associated to $\varpi$ are odd, 
contradicting the fact that the $\rhobar$ as above are 
even. 
\end{proof}

In light of this result, it is natural to ask the following
questions concerning the case $N = 1$.

\begin{question} What is the \emph{smallest}
irreducible even Galois representation $\rhobar:\Gal(\Qbar/\Q) \rightarrow \GL_2(\Fbar_p)$
which is unramified away from $p$? More precisely:
\begin{enumerate}
\item What is the smallest $p$ for which such a $\rhobar$ exists?
\item What is the smallest $k$ for which such a $\rhobar$ exists?
\end{enumerate}
\end{question}

Note that, \emph{a priori}, there may not exist any such $\rhobar$ at all!
Here is an approach to studying this question using 
Conjecture~\ref{conj:serre}(1).

\smallskip

It is known that no such representation exists when $p \le 7$ (for $p = 2$ this is
a well known result of Tate~\cite{TateNonExist}, for the other $p$ see
Theorem~1 of~\cite{Moon}).
For any $p$, the Serre weight $k:=k(\rhobar)$ is an integer $2 \le k \le p^2 - 1$. After twisting,
 we may assume that $k$ actually satisfies the inequality
  $2 \le k \le p+1$.
 Given $\rhobar$, we consider the restriction of $\rhobar$ to a small imaginary quadratic
 field of class number $1$
  (for our computations, we used the fields $F = \Q(\sqrt{-2})$ and $\Q(\sqrt{-11})$).
 It is not hard to see that if $N = 1$, that $\rhobar|G_{F}$ is still irreducible,
 since they would otherwise arise
 from inductions of characters of the class group of $F$, which is trivial.
 Hence, assuming Conjecture~\ref{conj:serre}(1), $\rhobar$ gives rise to
  a mod-$p$ cohomology class for $\GL_2(\OL_F)$.
We shall be interested in computing the following quantity:
$$ \dim H_1(\GL_2(\OL_F),\M_k \otimes \F_p) - \dim H_1(\GL_2(\OL_F),\M_k \otimes \C)^{BC},$$
where $\M_k$ is an integral model for the local system corresponding to the representation
 $\Sym^{k-2} \otimes \overline{\Sym^{k-2}}$,
and $BC$ denotes the subspace of forms which are either CM or arise from base change from $\GL(2)/\Q$.
If $\rhobar$ has Serre weight $k$, then the first cohomology group will have larger dimension than the second.
Hence, if we compute the quantity above and it is non-zero, this indicates that there exists a possible $\rhobar$.
Moreover, if we compute this quantity for \emph{different} fields $F$ and still obtain a non-zero value, this gives strong evidence
for the existence of $\rhobar$. (Conversely, if we compute for some field $F$ and find that the quantity above is zero, then $\rhobar$
can not exist).
\medskip

Our computations are summarized in the final table of Chapter~\ref{chapter:ch8}.
 It is the case that in the range of any computation that we do that
$\dim H_1(\M \otimes \C)^{BC} = \dim H_1(\M \otimes \C)$. Hence, when we are able to
compute  $H_1(\M)$ integrally, we can determine when the quantity above is zero for a fixed $k$ and all $p$
simultaneously, namely, the exponents of the torsion subgroup of $H_1(\M)$. 
(If $\dim H_1(\M \otimes \C)^{BC} < \dim H_1(\M \otimes \C)$, we would need to do some extra computations with
Hecke operators to determine the finite set of possible $p$.)
In the range we are able to compute, we find only one irreducible even representaiton.
Specifically, there is a representation
$$G_{\Q} \rightarrow \Gal(K/\Q) \simeq \widetilde{A_4}\rightarrow
\GL_2(\F_{163}),$$
corresponding to a lift of the even $A_4$-extension of
$\Q$ ramified only at $163$.
The determinant of this extension is (one of the) 
characters of $\Q(\zeta_{163})$ of order $3$, which  are all
conjugate to
$\omega^{(163-1)/3} = \omega^{55 - 1}$. Thus this representation corresponds to $(p,k) = (163,55)$.
We can summarize the results of our computation as follows:
\begin{theorem} \label{theorem:comp}
Let
$\rhobar:G_{\Q} \rightarrow \GL_2(\Fbar_p)$ be an absolutely
irreducible continuous even Galois representation
of Serre level $N(\rhobar) = 1$ and of Serre weight $k = k(\rhobar)$.
In particular, $\det(\rhobar) = \omega^{k-1}$, so $k$ is odd.
Then, assuming Conjecture~\ref{conj:serre}(1):
\begin{enumerate}
\item The smallest prime $p$ for which a $\rhobar$ exists
is at least $79$.
\item The smalest $k$ for which a $\rhobar$ exists is at least $33$.
\item If there exists a $\rhobar$ with $k \le 53$, then  $p > 1000$.
\item If there exists a $\rhobar$ with  $k  = 55$, then  $p$ is either $> 200$, or $p = 163$.
Moreover, if $k = 55$ and $p = 163$, then  $\rhobar$ is the
unique representation with projective image $A_4$.
\item If there exists a $\rhobar$ with $k \in \{57, \ldots, 69\}$, then either $p > 200$ or $(p,k)$ corresponds to
a red square 
in Table~\ref{table:mixed}.
\end{enumerate}
\end{theorem}

\chapter{Numerical examples} %
\label{chapter:ch8}
\label{part:examples}

This chapter gives numerical examples of homology for Jacquet--Langlands pairs which illustrate --- and in many cases motivated ---
the results of this book. 
Nathan Dunfield computed much of the data for our main example concerning the
field $F = \Q(\sqrt{-2})$  (using \texttt{Magma})
during the writing of~\cite{CD}. %
The complete data is reproduced as a table in~\S~\ref{section:tables}.
We computed
Hecke actions for a subset of his data using a \texttt{gp-pari} script that is available on request.   Where not otherwise specified, the Hecke algebra refers to the {\em abstract} Hecke algebra of~\S~\ref{CompletionConvention}.

The Chapter considers examples in roughly increasing order of complexity. 
After summarizing some geometric features of the manifolds at hand in~\S~\ref{sec:mans}, 
 we discuss first, in~\S~\ref{sec:nocharzero}, situations similar to that of Theorem   A
 of the introduction, where neither side has nontrivial homology in characteristic zero.  The examples in~\S~\ref{sec:mans}
 illustrate this Theorem  as well as some of the results of Chapter~\ref{chapter:ch3} on level-raising. Next, 
 in~\S~\ref{sec:oldcharzero}, we discuss situations similar to that of Theorem~B in the introduction.
In~\S~\ref{sec:newcharzero} we address the more complicated situation where there exists characteristic zero newforms, 
giving evidence for the conjectural role of level-lowering congruences (cf.~\S~\ref{reglowlevel}) and  illustrating
other results of Chapter 7.  Sections
\S~\ref{sec:Eisenstein7} and~\S~\ref{section:phantomtwo}
 study Eisenstein homology, in particular situations related to ``Theorem C'' (Theorem~\ref{theorem:unamiC})
 of the introduction
and studied also in Chapter 8.  
Section~\S~\ref{section:K2examples} provides numerical examples illustrating Theorems~\ref{theorem:K2popularversion}
and~\ref{theorem410} for the field $F = \Q(\sqrt{-491})$.

\section{The manifolds} \label{sec:mans}

Let $F  = \Q(\theta)$, with $\theta = \sqrt{-2}$. Since $3$ splits in $F$, and $\OL_F$ has class number one, we may
write $3 = \pi \cdot \pibar$, where  $\pi = 1 + \theta$ and $\pi' = 1 - \theta$.
Let    $D/F$ be the unique quaternion algebra
ramified exactly at $\pi$ and $\pibar$, and 
let $B/\OL_F$ be the unique
(up to conjugation) maximal order of $D$.
In this section, we shall compare the homology of arithmetic quotients corresponding
to $\PGL(2)/F$ and $\GL(D^{\times})/\mathbb{G}_m$.

If $\Sigma$ is a set of finite places, the arithmetic orbifold $Y(\Sigma)$
corresponding to $\PGL(2)$ has been defined in~\S~\ref{sec:notn0};
similarly, if $\Sigma$ contains $(\pi)$ and $(\pibar)$
we have an arithmetic orbifold $Y_B(\Sigma)$   corresponding to $D$. 
Note that $w_F = 2$ and $w_F^{(2)} = 8$, so $2$ is the only possible orbifold prime. 

At ``minimal level'', these are explicitly described as follows:

  Let $\Gamma =\PGL_2(\OL_F)$ and
$\Gamma_B$ the image of $B^{\times}$ in $\PGL_2(\C)$. Since $\Cl(\OL_F)$ is trivial,
we have first of all
$$Y := Y(\emptyset) =\Gamma \backslash \H, \ \ Y_B := Y_B(\{ (\pi), (\pibar)\}) =\Gamma_B \backslash \H^3$$
The orbifold $Y$ is a degree two quotient of the Bianchi Orbifold for $F$ corresponding
to $\PSL_2(\OL_F)$ (see~\cite{Bianchi}). 
The orbifold
$Y_B$ was the main object of study in~\cite{CD,BE}, where it was showed that the 
$\pi^n$-adic tower $\Gamma_B(\pi^n)$ had vanishing $\Q$-cohomology for all $n$.
Both orbifolds
have underlying space $S^3$ (as follows from~\cite{Hatcher}  and~\cite{CD},~\S~2.9, respectively),
although this fact is not important for us. Of course, $Y$ and $Y_B$ are not
commensurable as $Y$ is compact whereas 
$Y_B$ has a single cusp whose structure
is given by a torus modulo the hyperelliptic involution.

\medskip

Now take an auxiliary prime $\p \in \OL_F$,
and  take $\Sigma = \{(\pi), (\pibar), \p\}$,
with corresponding quotients $Y(K_{\Sigma})$ (also denoted $Y_0(3\p)$)
and $Y_B(K_{\Sigma})$. These are still quotients of $\H^3$:

 $$Y(3 \p) = \Gamma_0(3\p) \backslash \H^3, \qquad 
Y_B(K_{\Sigma}) = \Gamma_0(B,\p) \backslash \H^3$$
where $\Gamma_0(3\p)$ denotes the level $3\p$ congruence
subgroup of $\Gamma = \PGL_2(\OO_F)$, whereas
$\Gamma_0(B, \p)$ denotes the level $\p$ congruence subgroup
of $B^{\times}$.

According to Conjecture~\ref{conj:reciprocity},
loosely speaking, we expect the following:
$$
\left\{\text{classes in $H_1(\Gamma_0(B,\p),\Z)$} \right\}
\leftrightsquigarrow 
\left\{ \begin{aligned} & \text{\ \  homology classes in $H_1(\Gamma_0(3 \p),\Z)$}  \\ &  \text{that do \emph{not}
arise 
from $H_1(\Gamma_0(\pi \p),\Z)$} \\  & \text{\qquad or $H_1(\Gamma_0(\pi' \p),\Z)$ for   $\pi,\pi' | 3$}
\end{aligned} \right\}.$$

 \subsection{The torsion ratio} Suppose that $H_1(\Gamma_0(3 \p),\Z)$ is finite. By Theorem~\ref{TJL2},
  the ratio: 
$$A_p:= |H_1(\Gamma_0(B,\p),\Z)|  \  \left\slash \  
\frac{|H_1(\Gamma_0(3 \p),\Z)| \cdot  |H_1(\Gamma_0(\p),\Z)|^4}
{|H_1(\Gamma_0(\pi \cdot \p),\Z)|^2 
\cdot |H_1(\Gamma_0(\pi \cdot \p),\Z)|^2} \right.$$
is of $2$-power order. 

 For $p \le 617$, we find
that $A_p = 2^6$ whenever $H_1(\Gamma_0(3 \p),\Z)$ is finite, which happens
exactly  for
$$\begin{aligned}
p =  & \ 19, 41, 43, 59, 97, 137, 163, 179, 227, 251, 281, 283, 379, 
 401, 449, 467,  \\
 & \ 499, 523, 547, 563, 569, 571, 577.\end{aligned}$$ 
The fact that the power of $2$ is constant suggests that our results presumably extend to mod-$2$
classes. 

We address some of these examples in more detail.
Note that by Theorem~\ref{theorem:simplecase} we have in this context an equality
between the orders of $H^{E}_1(\Gamma_0(3 \p),\Z)^{\p\q-\new}$ and
$H^{E}_1(\Gamma_0(B,\p),\Z)$ up to an integral factor $\chi(\Sigma)$.
We shall typically work over the ring $\Z_S:=\Z[1/S]$   \index{$\Z_S$}
where $S$ is a product of small primes, in order to make the
presentation of our groups more  compact.

\section{No characteristic zero forms}
 \label{sec:nocharzero}
\subsection{\texorpdfstring{$N(\p) = 59$}{N(p)=59}} Suppose that $N(\p) = 3 - 5 \theta$. 
The congruence homology
for $\Gamma_0(3\p)$ and $\Gamma_0(B,\p)$ was
effectively computed in Lemma~\ref{lemma:ordercong}; it has order (up to powers of $2$)
$N(\p) - 1 = 29$.
By computing the essential cohomology $H^E_1$ we effectively suprress a factor
of $\Z/29\Z$ in each of the following terms.
Let $S = 2$.
$$
\begin{aligned}
H^E_1(\Gamma_0(\p),\Z_S) =  &  \ 1 \\ 
H^E_1(\Gamma_0(\pi \cdot \p),\Z_S) = & \   (\Z/3)^2 \oplus (\Z/7), \\
H^E_1(\Gamma_0(\pi' \cdot \p),\Z_S) =  &  \ (\Z/3)^2 \oplus
 (\Z/37), \\
H^E_1(\Gamma_0(3 \p),\Z_S) = &  \  (\Z/3)^8 \oplus (\Z/9) \oplus (\Z/7)^2 
\oplus (\Z/ 37 \Z)^2,\\
H^E_1(\Gamma_0(B,\p),\Z_S) =  & \  (\Z/9). \end{aligned}$$
This is an example where no level raising (see~\S~\ref{sec:lr})  occurs. The mod-$7$ and mod-$37$
classes at level $\pi \p$ and $\pi' \p$ give rise to old-forms at level $3 \p$, but admit no
level raising congruences, and so (by Lemma~\ref{summary} $(4)$)
have no effect on the cohomology of 
$\Gamma_0(B, \p)$. %

 Let us explicitly check that the level raising congruences
fail for $\m$ of characteristic $3$. Recall that the level raising congruence
is $T_{\q} =  \pm (1 + N(\q)) \mod \m$. 
\begin{enumerate}
\item 
There is a unique maximal ideal $\m$ of characteristic $3$ and level $\pi \cdot \p$ 
with $\T/\m = \F_9$. We find that $T_{\pi'} \equiv 0 \mod \m$, and hence 
the class $\m$ does not give rise
to a class in $H^1(\Gamma_0(3 \p),\Z_S)^{\new}$
via level raising as in Theorem~\ref{theorem:ribet} 
(and one thus expects that there are no such new classes
--- see Remark~\ref{remark:steinbergnew}.)

\item
There is a unique maximal ideal $\m$ of characteristic $3$ and level $\pi' \cdot \p$
with $\T/\m = \F_9$. We find that the image of $T_{\pi}$ in $\F_9$ has minimal
polynomial $x^2+x-1$, and hence (in the same sense as above)
level raising does not occur.
\end{enumerate}
If we return to the \emph{new} maximal ideal of level $3\p$, we find an $\m$
such that $\T/\m = \F_3$ and $\T_{\m}$ is isomophic to $\Z/9 \Z$.   Our theorems guarantee
that the power of $3$ dividing the order of $H^E_1(\Gamma_0(B,\p),\Z_S)$ is exactly $3^2$, but do not
promise that the Hecke actions are the same
(although, of course, we expect this to be true).

\subsection{\texorpdfstring{$N(\p) = 521$}{N(p)=521}} Suppose that $N(\p) = 521$.
The congruence
 cohomology for $\Gamma_0(3\p)$ and $\Gamma_0(B,\p)$ was
effectively computed in Lemma~\ref{lemma:ordercong}; it has order 
(up to powers of $2$) $N(\p) - 1 = 65$. Let $S$ be the product of all
primes $< 100$.
$$
\begin{aligned}
H_1(\Gamma_0(\p),\Z_S) =  &  \  (\Z/2213) \\ 
H_1(\Gamma_0(\pi \cdot \p),\Z_S) = & \   (\Z/139) \oplus (\Z/169) \oplus (\Z/2213)^2 \\
H_1(\Gamma_0(\pi' \cdot \p),\Z_S) =  &  \ (\Z/2213)^2 \oplus (\Z/6133) \\
H_1(\Gamma_0(3 \p),\Z_S) = &  \  (\Z/139)^2 \oplus (\Z/169)^2 \oplus (\Z/223)
\oplus (\Z/361) \oplus (\Z/401) \ \oplus  \\
&  \ (\Z/2213)^4 \oplus (\Z/6133)^2
\oplus (\Z/1750297)  \oplus (\Z/2613151) \\
H_1(\Gamma_0(B,\p),\Z_S) =  & \  (\Z/223) \oplus (\Z/401) \oplus
(\Z/1750297)  \oplus (\Z/2613151). \end{aligned}$$
This is another example where no level raising occurs, at least for classes mod $p$
of order $> 100$. We present it to illustrate that large primes can divide the orders
of these torsion classes.
(This also provides a compelling check on the numerical data --- for the prime 
 $2613151$ to appear incorrectly in any one computation would be a misfortune, for
 it to incorrectly appear twice would seem like carelessness.)

\subsection{\texorpdfstring{$N(\p) = 97$}{N(p)=97}} Suppose that $\p = 5 + 6 \theta$. 
The congruence cohomology for $\Gamma_0(3\p)$ and $\Gamma_0(B,\p)$ was
effectively computed in Lemma~\ref{lemma:ordercong}; it has order $N(\p) - 1 = 3$.
We are interested in the non-Eisenstein ideal $\m$ of residue characteristic $3$
and level $\pi' \cdot \p$.
We find that
$$
\begin{aligned}
H_1(\Gamma_0(\p),\Z)_{\m} =  &  \ 1 \\ 
H_1(\Gamma_0(\pi \cdot \p),\Z)_{\m} = & \   1 \\ 
H_1(\Gamma_0(\pi' \cdot \p),\Z)_{\m} =  &  \ (\Z/3) \\ 
H_1(\Gamma_0(3 \p),\Z)_{\m} = &  \  (\Z/3) \oplus (\Z/9) \\ 
H_1(\Gamma_0(B,\p),\Z)_{\m} =  & \  (\Z/3). \end{aligned}$$
If $\T$ is the Hecke algebra at level $\pi' \cdot \p$, so $\T_{\m} = \F_3$, we may compute
the small Hecke operators in $\T_{\m}$ explicitly as follows:
 \begin{center}
\begin{tabular}{|c|c|c|c|c|c|}
\hline
$\q$ &  $1+\theta$ & $3 + \theta$ & $3 - \theta$ & $3+2\theta$ & $3 - 2 \theta$ \\
$T_{\q}$ & $1$ & $0$ & $-1$ & $0$ & $0$ \\ 
  \hline
  \end{tabular}
  \end{center}
  From the very first entry on this table, we see that $T_{\pi} = 1 = 1 + N(\pi) \mod \m$, and
  hence $\m$ is a level raising prime.
   And indeed, with $\Sigma = \{\pi,\pi'\}$, the injective map:
$$\Psi^{\vee}: H_1(\Sigma/{\pi},\Z)^2_{\m} \rightarrow H_1(\Sigma,\Z)_{\m},
\qquad \Z/3 \oplus \Z/3 \rightarrow \Z/3 \oplus \Z/9$$
has cokernel $\Z/3$, which is realized on the space of quaternionic homology.

\section{Characteristic zero oldforms} \label{sec:oldcharzero}
Let us now consider the case when there exist oldforms of characteristic zero.
By the classical Jacquet--Langlands theorem, the cohomology group
$H_1(\Gamma_0(B,\p),\Z)$
is finite.
In this context, it still makes sense to consider the na\"{\i}ve ratio:
$$A_p:= |H_1(\Gamma_0(B,\p),\Z)|  \  \left\slash \  
\frac{|H_1(\Gamma_0(3 \p),\Z)^{\tors}| \cdot  |H_1(\Gamma_0(\p),\Z)^{\tors}|^4}
{|H_1(\Gamma_0(\pi \cdot \p),\Z)^{\tors}|^2 
\cdot |H_1(\Gamma_0(\pi \cdot \p),\Z)^{\tors}|^2} \right.$$
where $M^{\tors}$ denotes the torsion subgroup of $M$.
In this context, our theorem (in its crude form Theorem~\ref{TJL1}) does not predict that this quantity is constant, but
rather that it equal to a ratio of regulators on the split side. A particularly
nice case to consider is when $H_1(\Gamma_0(\p),\Q)$ and
$H_1(\Gamma_0(\q \p),\Q) = 0$ for one of $\q \in \{\pi,\pi'\}$.
In this case we have Theorem~\ref{theorem:oldforms}, which guarantees (up to $\chi_D(Y)$) an
numerical equality between the spaces of essential newforms.

\subsection{\texorpdfstring{$N(\p) = 17$}{N(p)=17}} \label{example:data3} Let us take $\p = (2 \theta + 3)$, and let
$S = 10$.
Then:
$$
\begin{aligned}
H_1(\Gamma_0(\p),\Z_S) =  &  \ 1 \\ 
H_1(\Gamma_0(\pi \cdot \p),\Z_S) = & \   \Z_S \\ 
H_1(\Gamma_0(\pi' \cdot \p),\Z_S) =  &  \  1 \\ 
H_1(\Gamma_0(3 \p),\Z_S) = &  \  \Z_S \oplus \Z_S \\ 
H_1(\Gamma_0(B,\p),\Z_S) =  & \  (\Z/3). \end{aligned}$$
In this case, the na\"{\i}ve torsion ratio $A_p$ is equal to $3$. Our results explain this as follows.
Recall that the space of newforms in this case is equal to the cokernel of the map:
$$\Psi^{\vee}:  H_1(\Gamma_0(\pi \cdot \p),\Z_S)^2 \rightarrow H_1(\Gamma_0(3 \p),\Z_S), \qquad
\Z^2_S \rightarrow \Z^2_S.$$
The cokernel of this map is exactly the cokernel  of the map
$$\displaystyle{\left( \begin{matrix} (N(\q) + 1) & T_{\q} \\
T_{\q} &  (N(\q) + 1) \end{matrix} \right)}$$
 with $\q = \pi'$ on $\Z^2_S$.
Explicitly, we find that the value of $T_{\pi'}$ on $H_1(\Gamma_0(\pi \cdot \p),\Z_S)$ is
equal to $-2$, and hence the cokernel of the matrix and
hence the newspace  is isomorphic $\Z/3\Z$.
Thus the computation is in accordance with
Theorem~\ref{theorem:oldtwo}, which predicts
that
$$3 = 
|H^{E^*}_1(\Gamma_0(3 \p),\Z_{3})^{\new}| = |H^{E^*}_1(\Gamma_0(B,\p),\Z_{3})| \cdot \chi 
 = 3 \chi$$
 for some integer $\chi$. Note that, since the maximal ideal $\m$ of residue characteristic
 $3$ is (by computation) not Eisenstein, we would conjecture that $\chi = 1$ in
 this case, which is also borne out by the computation above.

Let us consider this
example in greater detail. Since we have a characteristic zero form with
$\Z$-coefficients in trivial weight,
we have a hope of attaching an elliptic curve via Taylor's theorem~\cite{Taylorimag}.
Indeed, there does exist an elliptic curve $E$ of conductor $\pi \cdot \p$, given as follows:
$$E: y^2 - \theta x y + y = x^3 + (\theta - 1) x^2,$$
of discriminant 
$$\Delta_E  = -14 \theta - 47 = - (1 + \theta)^2 (2 \theta - 3)^2.$$
 Let us study this representation locally at the prime of good reduction
 $\pi' = 1 - \theta$.  
  Reducing $E$ modulo $\pi'$, we obtain the smooth curve
 $$y^2 -  x y + y = x^3,$$
 over $\F_3$, from which we can count points to
 deduce that $a_{\pi'} = 1 + 3 -6 = -2$.
 It follows that $E$ is not only  ordinary at $\pi'$, but that the
 mod-$3$ representation $\rhobar: G_F \rightarrow \GL(E[3])$
 at a decomposition group $D_{\pi'}$ can be written in the form:
 $$
\left( \begin{matrix} \chi & * \\ 0 & 1 \end{matrix} \right) \mod 3.$$
Equivalently, the level raising congruence implies that 
the local representation looks like (a twist of) the Galois representation
attached to a Steinberg representation.
In the classical world of $\GL(2)/\Q$ this would imply that the representation
$\rhobar$ lifts
to a char zero newform of level $3 \p$, but here  it does not, and in fact
the only ``shadow of newness'' appearing is in the nontriviality of
the cokernel of $\Psi^{\vee}$, or equivalently, the torsion occurring in
cohomology on the non-split side.

\subsection{Lifting Torsion Classes}
 It is a consequence of level raising that the mod-$3$ torsion class of~\S~\ref{example:data3}
above %
\emph{does} have a Galois representation, since it lifts to characteristic
zero at level $3\p$. This suggests the very natural question:

\begin{question} Let $\m$ be  a non-Eisenstein maximal ideal of the Hecke algebra $\T$ acting
on $H_1(\Sigma,\Z)$. Does there exist a $\Sigma'$ containing $\Sigma$ such that $\m$ is in
the support of the \emph{torsion free} part of $H_1(\Sigma',\Z)$? Assuming 
Conjecture~\ref{conj:reciprocity}, this amounts to the following. If 
$$\rhobar: G_F \rightarrow \GL_2(\F)$$
is continuous and irreducible, does it lift to a motivic characteristic zero representation
of trivial weight?
\end{question}

We have little to say on either question.
If $\#\F = 2,3$ or $5$, and $\det(\rhobar)$ is the cyclotomic character, then the answer to
the second question is yes, since one can find a lift coming from an elliptic curve. This sheds
little light on the general case, however.
On the other hand, as far as explicitly searching for lifts,
level raising puts a restriction on the places where one should look to
lift any torsion class. There is no point adding $\Gamma_0(\q)$-level structure unless $\q$
is a level raising prime. This provides at least one place to look for any given torsion class.

\subsection{\texorpdfstring{$N(\p) = 409$}{N(p)=409}} Let us take $\p = (11 + 12 \theta)$, and let
$S$ be divisible by all primes $\le 3779$ except  $29$ and $37$.
Then:
$$
\begin{aligned}
H_1(\Gamma_0(\p),\Z_S) =  &  \ 1 \\ 
H_1(\Gamma_0(\pi \cdot \p),\Z_S) = & \   \Z/29 \\ 
H_1(\Gamma_0(\pi' \cdot \p),\Z_S) =  &  \  \Z^6_S \\ 
H_1(\Gamma_0(3 \p),\Z_S) = &  \  \Z^{12}_S \oplus (\Z/29)^3 \\ 
H_1(\Gamma_0(B,\p),\Z_S) =  & \  (\Z/29)^2 \oplus (\Z/37). \end{aligned}$$
There are two maximal non-Eisenstein ideals of residue characteristic $29$ of the abstract Hecke algebra acting on these groups, one supported at level
$\pi \cdot \p$, and another at level $\pi' \cdot \p$. That these ideals are distinct is not surprising, since otherwise one would predict
that it also has support at level $\p$.
Let us study the $6$-dimensional space over $\Q$ at level $\pi' \cdot \p$.
If $\T$ is the Hecke algebra at this level then $H = \T \otimes \Q$ is a
direct sum of two totally
real fields of degrees $2$ and $4$ respectively, with
$T_{\pi}$ eigenvalues $\alpha$ and $\betavar$ which can be determined
explicitly as follows:
$$ \alpha = \frac{-3 + \sqrt{5}}{2}, \qquad \qquad  \alpha^2 + 3 \alpha + 1 = 0 $$
$$\betavar =
\frac{-3 + \sqrt{41}}{4} + \sqrt{\frac{13 + \sqrt{41}}{8}},
\qquad \betavar^4 + 3 \betavar^3 - 5 \betavar^2 - 16 \betavar - 8 = 0.$$
(We obtained these equations by computing the characteristic polynomial
of $T_{\pi}$ modulo the primes $61$ and  $103$ and then comparing the answer.
The resulting polynomial was then checked modulo the prime $32749$, and found to be correct.
Assuming the known bounds towards the Ramanujan conjecture, this can be shown to confirm
our computation.)
We compute that $N(\alpha^2 - 16) =  5 \cdot 29$, and
$N(\betavar^2 - 16) = 2^6 \cdot 5 \cdot 37$.
Thus $H_1(\Gamma_0(3 \p),\Z_S)^{\new} = (\Z/29)^2 \oplus (\Z/37)$, which
--- exactly in accordance with what one expects from Theorem~\ref{theorem:oldtwo} ---
 has the same order as $H_1(\Gamma_0(B,\p),\Z_S)$. 
Note that $A_{409} = 2^{11} \cdot 5^2 \cdot 29 \cdot 37$ in this case.

\section{Characteristic zero newforms and level lowering} \label{sec:newcharzero} In this section, we consider some examples
where there exist \emph{newforms} in characteristic zero. In this context, we may apply Theorem~\ref{theorem:newforms} to
deduce a relationship between primes dividing the ratio of the regulators and level lowering primes.

\subsection{\texorpdfstring{$N(\p) = 11$}{N(p)=11}}  \label{subsec:level lowering}  Take $\p = (3 +  \theta)$. If $S = 10$, we have
$$
\begin{aligned}
H_1(\Gamma_0(\p),\Z_S) =  &  \ 1 \\ 
H_1(\Gamma_0(\pi \cdot \p),\Z_S) = & \   1 \\ 
H_1(\Gamma_0(\pi' \cdot \p),\Z_S) =  &  \ (\Z/3) \\ 
H_1(\Gamma_0(3 \p),\Z_S) = &  \  \Z_S \oplus (\Z/3) \\ 
H_1(\Gamma_0(B,\p),\Z_S) =  & \  \Z_S. \end{aligned}$$
It is a consequence of our results that the ratio of the regulators on the split
side to the non-split side must be divisible by three; according
to our conjectures in~\S~\ref{reglowlevel}, this regulator discrepancy should
be accounted for by level lowering primes.  We confirm this:

In this case,  the elliptic curve
$$E: y^2 + (1 + \theta) x y + (1 + \theta) y = x^3 + (1 - \theta) x^2 - 2 \theta x - \theta$$
has conductor $3 \p$. As we shall see, its level can be lowered {\em modulo $3$};
in other terms, the reduction of this class modulo $3$ comes from level $\pibar \cdot \p$. 

\medskip

Note that $E$ has is multiplicative reduction at $\pi$ and
$\pi'$. We compute that the prime factorization of $\Delta_E$ is
$$\Delta_E = 36 \theta + 9 = (1 + \theta)^3 (1 - \theta)^2 (3 + \theta) = \pi^3 \pi'^2 \p,$$
In particular, since $3$ divides $v_{\pi}(\Delta_E)$, the $3$-torsion is
\emph{finite
flat} at $\pi$ (it would have been equivalent to note that the $\pi$-adic valuation 
of the $j$-invariant $j_E$ of $E$ was also divisible by $3$). It follows from Conjecture~\ref{conj:reciprocity} that the mod-$3$
representation $\rhobar$ should arise at  a lower level, namely, level
$\pi' \p$. In particular, the $\Z/3\Z$ torsion class above at level $\pi' \p$ corresponds the the Galois
representation arising from $E[3]$. 

\subsection{A complicated example,\texorpdfstring{$N(\p) = 89$}{N(p)=89}} In the following example, we observe
 some \label{subsection:pathologies}
of the complications and pathologies that can occur in the general case.
Let $N(\p) = 89$, with $\p = (9 -  2 \theta)$. Let $S = 2$.
We have
$$
\begin{aligned}
H_1(\Gamma_0(\p),\Z_S) =  &  \ (\Z/11)^2 \\ 
H_1(\Gamma_0(\pi \cdot \p),\Z_S) = & \   \Z^2_S \oplus (\Z/3) \oplus (\Z/11)^2  \\
H_1(\Gamma_0(\pi' \cdot \p),\Z_S) =  &  \ \Z_S \oplus (\Z/11)^3 \\ 
H_1(\Gamma_0(3 \p),\Z_S) = &  \  \Z^7_S \oplus (\Z/3) \oplus (\Z/5)^2
\oplus (\Z/7) \oplus (\Z/11)^3  \\ 
H_1(\Gamma_0(B,\p),\Z_S) =  & \  \Z_S \oplus (\Z/5) \oplus (\Z/7) \oplus (\Z/11)
\oplus (\Z/25) \end{aligned}$$
We consider this example in several steps, concentrating
on the maximal ideals $\m$ of $\T$ of characteristic $11$, $3$, and $5$.
First, we identity the characteristic
zero classes. For convenience, we tabulate these now.
The coefficients of the Hecke operators
are $\Z[\phi]$, $\Z$, and $\Z$, where $\phi = \displaystyle{\frac{1 + \sqrt{5}}{2}}$ 
is the golden ratio. Two of these automorphic forms can be attached to elliptic curves.
There is an elliptic curve of conductor $\pi' \p$ given by the following equation:
  $$B: y^2 - \theta yx + (1 - \theta) y = x^3 + 2 \theta x + (3 \theta - 1).$$
  The discriminant of $B$ is $\Delta_B=  (9 - 2 \theta) (1 - \theta)^{12}$.
  The curve $B$ admits a $3$-isogeny  over $F$, and has a
  $3$-torsion point $P = [3-\theta,4]$.
  There is also an elliptic curve of conductor $3 \p$ given as follows:
  $$C: y^2 - x y + y = x^3 + (2\theta - 1) x^2 - (2\theta + 2) x + 2.$$
The curve $C$ has discriminant $\Delta_C=- (9 - 2 \theta) (1 + \theta)^4 (1 - \theta)^2$.
Presumably, there also exists an Abelian variety $A$ of $\GL_2$-type
and conductor $\pi {\p}$ defined 
over $F$ with real multiplication by $\Q(\sqrt{5})$, although we have not tried
to find it explicitly. We compute the first few Hecke eigenvalues of this forms in 
the following table.
 \begin{center}
\begin{tabular}{|c|c|c|c|c|c|c|c|c|}
\hline
$\q$ & $1+\theta$ & $1 - \theta$ & $3 + \theta$ & $3 - \theta$ & $3+2\theta$ & $3 - 2 \theta$ & $1 + 3 \theta$ &
$1 - 3 \theta$  \\
$a_{\q}$ &   & $2 \phi$ & $2 \phi + 2$ & $2 - 2 \phi$ & $-2$ & $4 - 4 \phi$ & $2-6 \phi$ & $4 - 4 \phi$ \\
$b_{\q}$ & $-2$ &  &  $0$ & $0$ & $6$ & $-6$ & $8$ & $2$ \\
$c_{\q}$ &   & &$0$ &$0$ &$-6$ &$-2$ & $0$ & $-4$ \\
\hline
  \end{tabular}
  \end{center}
Here $a_{\q}$, $b_{\q}$, $c_{\q}$ correspond to the characteristic zero newforms
of level $\pi \p$, $\pi' \p$, and $3 \p$ respectively.

\medskip

Let us  consider the various maximal ideals of the Hecke algebras of
residue characteristic $11$.
There is an Eisenstein maximal ideal of characteristic
$11$, arising from congruence homology. 
It turns out that $H_1(\Gamma,\Z_S)_{\m} = \Z/11$ for all $\Gamma$ in the above
table.

 There is a non-Eisenstein maximal ideal $\m$ of level
$\Gamma_0(\p)$, arising via the action of Hecke on the $\Z/11$ summand
of $H_1(\Gamma_0(\p),\Z_S)$.  Let us denote the  image of the Hecke operator
$T_{\q}$ in $\T/\m = \F_{11}$ by  $d_{\q}$.
We may compute the first few values of $d_{\q}$ and record them in the following table.
 \begin{center}
\begin{tabular}{|c|c|c|c|c|c|c|c|c|}
\hline
$\q$ &  $1+\theta$ & $1 - \theta$ & $3 + \theta$ & $3 - \theta$ & $3+2\theta$ & $3 - 2 \theta$ & $1 + 3 \theta$ &
$1 - 3 \theta$  \\
$d_{\q}$ &  $-4$ & $-3$ &  $-1$ & $5$ &  $-2$ & $-1$ & $0$ & $-1$ \\
  \hline
  \end{tabular}
  \end{center}
  The entries of this table lie in $\T/\m = \F_{11}$.

   There are the maximal ideals of characteristic $11$ arising
   from the reduction of characteristic zero classes. The ring $\Z[\phi]$
  has two primes $(11,\phi - 4)$, $(11,\phi + 3)$ above $11$, which gives rise
  to two distinct maximal ideals of $\T$.
 
 \medskip
 
 We record these reductions in the following table:
 {\tiny
   \begin{center}
\begin{tabular}{|c|c|c|c|c|c|c|c|c|}
\hline
$\q$ & $1+\theta$ & $1 - \theta$ & $3 + \theta$ & $3 - \theta$ & $3+2\theta$ & $3 - 2 \theta$ & $1 + 3 \theta$ &
$1 - 3 \theta$  \\
$a_{\q} \mod (11,\phi + 3)$ &   & $5$ & $-4$ & $-3$ & $-2$ & $5$ & $-2$ & $5$ \\
$a_{\q} \mod (11,\phi - 4)$ &   &  $-3$ &  $-1$ & $5$ &  $-2$ & $-1$ & $0$ & $-1$ \\
$b_{\q} \mod 11$ & $-2$ &  &  $0$ & $0$ & $6$ & $-6$ & $8$ & $2$ \\
$c_{\q}$ &   & &$0$ &$0$ &$-6$ &$-2$ & $0$ & $-4$ \\
\hline
  \end{tabular}
  \end{center}
  }
  
  Let $\m$ be the non-Eisenstein maximal ideal of level $\Gamma_0(\p)$ and
characteristic $11$. Then since
$$d_{\pi} = -4 = -1 -3 \mod 11,$$
 $\m$ gives rise to level raising via Theorem~\ref{theorem:ribet}.
By inspection of the first few eigenvalues,
 we may identify $d_{\q}$ with $a_{\q} \mod (11,\phi - 4)$.
 On the other hand, $d_{\pi'} = -3 \not\equiv \pm (1 + 3) \mod 11$, and one finds
 (in accordance with  Theorem~\ref{theorem:ribet} and  Remark~\ref{remark:steinbergnew}.)
that  $\m$ does not give rise to level raising from either level $\p$ to $\pi' \p$ or from
 level $\pi \p$ to $3 \p$. It follows that if we complete at $\m$, we obtain the following:
$$
 \begin{aligned}
H_1(\Gamma_0(\p),\Z)_{\m} =  &  \ \Z/11 \\ 
H_1(\Gamma_0(\pi \cdot \p),\Z)_{\m} = & \   \Z_{11}  \oplus (\Z/11)  \\
H_1(\Gamma_0(\pi' \cdot \p),\Z)_{\m} =  &  \  (\Z/11)^2 \\ 
H_1(\Gamma_0(3 \p),\Z)_{\m}= &  \  \Z^2_{11} \oplus (\Z/11)^2 \\
H_1(\Gamma_0(3 \p),\Z)^{\new}_{\m} = & \ 1  \end{aligned}$$
Together with the congruence homology, this accounts for all the
$11$-torsion. None of the other maximal ideals of residue characteristic
$11$ satisfy level raising congruences, and thus $H_1(\Gamma_0(3\p),\Z)^{\new}$
``explains'' $H_1(\Gamma_0(B,\p),\Z)$.  In particular, the ratio of regulators
(given by $A_p$) is not divisible by $11$.
  
  \medskip
  
Let us now consider the maximal ideals of characteristic $5$. At level $3 \p$ it turns out
  that one of the $\Z/5$-torsion classes gives rise to a new maximal ideal, described by
  the following Hecke eigenvalues:
  \begin{center}
\begin{tabular}{|c|c|c|c|c|c|c|}
\hline
$\q$ &   $13 + \theta$ & $3 - \theta$ & $3+2\theta$ & $3 - 2 \theta$ & $1 + 3 \theta$ & $1 - 3 \theta$ \\
$f_{\q}$ &   $3$ & $0$ &  $0$ & $1$  & $3$ & $3$ \\
  \hline
  \end{tabular}
  \end{center}
  This class is not congruent to the other classes (coming from characteristic zero),
   and
  it (conjecturally) is mirrored by the $\Z/5$-class in the cohomology of $\Gamma_0(B,\p)$. The classes
  coming from characteristic zero are as follows (recalling that $\phi \equiv 3 \mod \sqrt{5}$):
   \begin{center}
\begin{tabular}{|c|c|c|c|c|c|c|c|c|c|c|}
\hline
$\q$ & $1 + \theta$ & $1 - \theta$ & $3 + \theta$ & $3 - \theta$ & $3+2\theta$ & $3 - 2 \theta$ & $1 + 3 \theta$ &
$1 - 3 \theta$  \\
$a_{\q}$   & &  $1$ & $3$ & $1$ & $3$ & $2$ & $4$ & $2$ \\
$b_{\q}$ & $3$ &  &  $0$ & $0$ & $1$ & $4$ & $3$ & $2$ \\
$c_{\q}$ & & &  $0$ &$0$ &$4$ &$3$ & $0$ & $1$ \\
\hline
  \end{tabular}
  \end{center}
  Since $b_{\pi} = 3 \not \equiv \pm (1 + 3) \mod 5$, there is no level raising
  from the maximal ideal of level $\pi' \p$. Since
  $$a_{\pi'} = 1 \equiv -(1+3) \mod \sqrt{5},$$
  we do get a level raising maximal ideal $\m$ in this case. Moreover, at level $\pi \p$, we see
  that $\T_{\m} \simeq \Z_5[x]/(x^2 - 5).$
  Let us consider this maximal ideal at level $3 \p$, where we denote it by
  $\m^{\old}$.
  Over $\F_5$, we compute that $H_1(\Gamma_0(3 \p),\F_5) \otimes \T_{\m^{\old}}$
  has dimension $5$. However, it is useful to note that at level $3 \p$ there is an enriched
  Hecke algebra that acts with the new Hecke operator $U_{v}$ for $v|3$ . The action of $U_{v}$
  on newforms is given by $\pm 1$. On old forms, it acts via the eigenvalues of the
  characteristic polynomial $X^2 - T_{v} X + N(v)$. In particular, if $v = \pi'$, this polynomial
  is given by $X^2 - 2 \phi X + 3 = (X+1)(X+3) \mod 5$. (The fact that $\m$ admits a
  level raising congruence is equivalent to the fact that this polynomial has
  $\pm 1$ as a root.)
 In particular, there exists an enriched maximal ideal
  $\m = (\m^{\old},U_{\pi'} + 1)$ of $\T$ at level $3 \p$.
  With respect to this ideal, we compute that
  $$H_1(\Gamma_0(3 \p),\F_5) \otimes \T_{\m} = (\F_5)^3.$$
  We may compute this Hecke action explicitly. We find that the $\m$-torsion of
  $H_1(\Gamma_0(3 \p),\F_5)$ has rank $1$, from which we deduce that
  $\T_{\m}/5$ has order $5^3$. Surprisingly, perhaps, we find that $\T_{\m}/5 = \F_5[x]/x^3$.
  Explicitly, the first few eigenvalues are given as follows:
    \begin{center}
\begin{tabular}{|c|c|c|c|c|c|c|}
\hline
$\q$ &   $3 + \theta$ & $3 - \theta$ & $3+2\theta$ & $3 - 2 \theta$ & $1 + 3 \theta$ & $1 - 3 \theta$ \\
$T_{\q} \mod \m$ &   $3 + 2x$ & $1 + 3x$ &  $1$ & $2 + x$ & $4  - x + x^2$ & $2 + x + 2 x^2$ \\
  \hline
  \end{tabular}
  \end{center}
  Since there exists a map $\T_{\m} \rightarrow \Z_5[x]/(x^2 - 5)$, which is surjective away
  from a torsion class of order $5$, we deduce that
  $$\T_{\m} \simeq \Z_5[x]/((x^2 - 5)(x,5)).$$
  On the other hand, on the quaternionic side, the only $5$-torsion that is unaccounted
  for is the $\Z/25\Z$-torsion class.
  Explicitly, then, our conjectures in this case amount to the following:
  \begin{conj} Let $\rhobar: G_F \rightarrow \GL_2(\F_5)$ be the Galois
  representation associated to $\phi$-torsion $A[\phi]$ of $A$. 
  Let $R$ denote the universal deformation
  ring of $\rhobar$ which records deformations of $\rhobar$ that have cyclotomic
  determinant, are ordinary at $v|3\p$. Then 
  $$R \simeq \T_{\m}
  =  \Z_5[x]/(x^2 - 5)(x,5)).$$
  Let $R^{\new}$ denote the quotient of $R$ corresponding to deformations $\rho$
  of $\rhobar$ as above such that for the place $v =  \pi'$,
 $$\rho | D_{v} = \left( \begin{matrix} \eps \chi & * \\ 0 & \chi \end{matrix} \right)$$
 for some unramified character $\chi$ of order $2$.
 Then $R^{\new} \simeq \T^{\new}_{\m} = \Z/25$, and the map
 $\T_{\m} \rightarrow \T^{\new}_{\m} = \Z/25$ is given by $x \mapsto 5$. 
 \end{conj}
 An implication for the homology is that the space of $\m$-newforms should be isomorphic to $\Z/25\Z$.
 Let us now try to compute this space.
  Let $A = H_1(\Gamma_0(\pi \p),\Z)_{\m}$. The space of newforms is determined by
  taking the cokernel of the map:
  $\Psi^{\vee}: A^2 \rightarrow A^2 \oplus \Z/5$.
  We know that the cokernel of the composite:
  $$
  \begin{diagram}
  A^2 & \rTo^{\Psi^{\vee}} & A^2 \oplus \Z/5 & \rTo^{\Psi} & A^2 \end{diagram}$$
  is the cokernel of the matrix
    $\displaystyle{\left( \begin{matrix} (N(\q) + 1) & T_{\q} \\
T_{\q} &  (N(\q) + 1) \end{matrix} \right)}$ on $A^2$, which has determinant
$$a^2_{\pi'}
  - (1 + N(\pi))^2 = 4 \phi^2 - 16 = 4 \phi -12 = 2 \sqrt{5}(1 - \sqrt{5}),$$
  of norm $80$.
  Hence the cokernel of the composite map   has order $5$.
    Thus the cokernel of $\Psi^{\vee}$ has order $25$.
    The map on
  one the the $A$ factors will be an isomorphism. 
   Writing $A = \Z^2_5$
 as a $\Z_5$-module, the map $\Psi^{\vee}$ is given explicitly on the other
 $A$-factor by:
 $$\Psi^{\vee}: A = \Z^2_5 \rightarrow A \oplus \Z/5 = \Z^2_5 \oplus \Z/5, \qquad
 (1,0) \mapsto (p,0,\alpha), \qquad (0,1) \mapsto (0,1,\beta).$$
 We check that the projection to the first factor has image $\{(p,0),(0,1)\}$ which has
 index $5$. Let us now consider the quotient.
 If $\alpha \ne 0$, then  the quotient is given by $\Z/25$ via the map
 $$(1,0,0) \mapsto  1, \quad (0,1,0) \mapsto  \beta \alpha^{-1} p, \quad (0,0,1) \mapsto - \alpha^{-1} p.$$
 We expect, but can not (at this point) compute, that $\alpha \ne 0$.
 This would imply that $H_1(\Gamma_0(3 \p),\Z)^{\new}_{\m} =  \Z/25$, which would
 be in agreement with what one obtains on the quaternionic side.
  
  \medskip

Finally, let us  consider the maximal  ideals of residue characteristic $3$.
Since $3$ is inert in $\Z[\phi]$, at level $\pi \cdot \p$ we have two maximal
ideals of residue characteristic $3$, one with $\T/\m = \F_9$ and $T_{\q}$
mapping to the reduction of $a_{\q} \mod 3$, and other given by the
torsion class $\Z/3$ with Hecke eigenvalues as follows:
\begin{center}
\begin{tabular}{|c|c|c|c|c|c|c|c|c|}
\hline
$\q$ &  $1 + \theta$ & $1 - \theta$ & $3 + \theta$ & $3 - \theta$ & $3+2\theta$ & $3 - 2 \theta$ & $1 + 3 \theta$ &
$1 - 3 \theta$  \\
$e_{\q}$ &   & $1$ & $0$ & $0$ &  $0$ & $0$ & $-1$ & $-1$ \\
  \hline
  \end{tabular}
  \end{center}
  The entries of this table lie in $\T/\m = \F_{3}$. We find that $e_{\pi} = 1
  \equiv 1 + 3 \mod 3$ and so $\m$ is a level raising prime.
  However, we also note that $e_{\q} \equiv b_{\q} \mod 3$.
  The elliptic curve $B$ has an $F$-rational
  $3$-torsion point, and hence the maximal ideal $\m$
  is Eisenstein.
  Thus we find that $\m$ is an Eisenstein maximal ideal, but not a congruence
  maximal ideal (and so it satisfies definition $D1$
  of~\S~\ref{section:eiz}, but not $D2$ or $D3$). Moreover, we see that $\m$ is supported
  in level $\pi  \p$ and $\pi'  \p$ but not $\p$. Perhaps surprisingly, at level 
  $\pi  \p$ it occurs as an atomic component $\Z/3$ by itself. Moreover, since 
  $\m$ is supported away from the congruence homology, our level raising  
  theorems still apply, and $\m$ must therefore be supported at any level divisible by
  $\pi \p$ or $\pi' \p$. 
  We have the following:
  $$
 \begin{aligned}
H_1(\Gamma_0(\p),\Z)_{\m} =  &  \ 1 \\ 
H_1(\Gamma_0(\pi \cdot \p),\Z)_{\m} = & \   \Z/3 \\
H_1(\Gamma_0(\pi' \cdot \p),\Z)_{\m} =  &  \  \Z_3 \\ 
H_1(\Gamma_0(3 \p),\Z)_{\m}= &  \  \Z^2_{3} \oplus (\Z/3) \\
H_1(\Gamma_0(\B, \p),\Z)^{\new}_{\m}  = & \ 1 \\
H_1(\Gamma_0(3 \p),\Z)^{\new}_{\m} = & \ ?  \end{aligned}$$
The composite map (mapping from level  $\pi \cdot \p$ to level $3 \p$ and then back again)
$$\Psi \circ \Psi^{\vee}: (\Z/3)^2 \rightarrow \Z^2_3 \oplus \Z/3 \rightarrow (\Z/3)^2$$
has cokernel $\Z/3$, and thus the image of the map $\Psi^{\vee}$ is equal to the torsion subgroup.
On the other hand, the composite map (mapping from level  $\pi' \cdot \p$ to level $3 \p$ and then back again):
$$\Psi \circ \Psi^{\vee}:  \Z^2_3 \rightarrow \Z^2_3 \oplus \Z/3 \rightarrow \Z^2_3$$
also has cokernel $\Z/3$, which can be computed from the determinant of the matrix
$$\det \left( \begin{matrix} 1 + N(\pi) & T_{\pi} \\ T_{\pi} & 1 + N(\pi) \end{matrix} \right)
= \det \left( \begin{matrix} 4 & -2 \\ -2 & 4 \end{matrix} \right) = 12.$$
Since $\Psi$ is surjective (Ihara's Lemma), it follows that the image of $\Psi^{\vee}$ has
index $3$ in the torsion free part of $H_1(\Gamma_0(3 \p),\Z)_{\m}$. In particular, 
if we  consider the level raising map (mapping from old forms of level
$\pi \cdot \p$ and $\pi' \cdot \p$ to level $3 \p$)
$$\Psi^{\vee}:(\Z/3)^2 \oplus \Z^2_3 \rightarrow \Z^2_3 \oplus \Z/3,$$
we deduce that the index of $\Psi^{\vee}$ is at least three (since it has index at least
three on the torsion free quotient). Since the space of
newforms $H_1(\Gamma_0(3 \p),\Z)^{\new}_{\m}$ is the cokernel of $\Psi^{\vee}$,
we deduce:
\begin{lemma} \label{lemma:differentnew}
The following space of newforms are distinct:
$$1 = H_1(\Gamma_0(B, \p),\Z)_{\m} \ne H_1(\Gamma_0(3 \p),\Z)^{\new}_{\m} \ne 1.$$
Moreover, they have different orders. 
\end{lemma}

This lemma provides a counter-example to the most optimistic conjectures regarding
equality between spaces of newforms (although we can still hope for an equality of
orders when $\m$ is not Eisenstein). 
In particular, we see that the ideal $\m$ is not supported on the quaternionic side, even though
it is ramified only at $\p$ and satisfies level raising at $\pi$ and $\pi'$.
Curiously enough, the ideal $\m$, which is not supported at level $\p$ either.
Indeed, we see that $\m$ is exactly a \emph{phantom class}, as discussed
in~\S~\ref{section:phantomclasses}. Namely, $\m$ arises from a nontrivial
extension $\Ext^1(\mu_3,\Z/3\Z)$ unramified outside $\q$ and totally split at $3$, and the fact that
$N(\q) \not\equiv 1 \mod 3$ implies that this class cannot be seen at level $\q$.

\begin{remarkable}  
\emph{
We observe in this case that the level lowering sequence:  \label{remark:hangingchad}
$$0 \longleftarrow (\Z/3)^2 + \Z^2_3 \longleftarrow (\Z/3) + \Z^2_3$$
\emph{cannot} be exact at the middle row (tensor with $\Z/3$). In particular,
it follows in this case that   -- in the notation of Chapter 7 -- $\chi_D(\Sigma)$ is divisible by $3$.
This difference also accounts (in light of
Theorem~\ref{theorem:oldforms}) to the fact that the space of newforms
in Lemma~\ref{lemma:differentnew} have different orders.
See~\S~\ref{section:phantomclasses}  and Remark~\ref{remark:failure}
for some general spectulations concerning this example.
For a related example, see~\S~\ref{subsec:chad2}.
}
\end{remarkable}

\subsection{\texorpdfstring{$N(\p) = 179$}{N(p)=179}}  \label{subsec:chad2}  Take $\p = (9 + 7  \theta)$. 
For this section, $\Z_3$ will denote the $3$-adic integers. Note that $3$ does
not divide $179 -1$, and so there is no congruence homology at this level.
In particular, over $\Z_3$, the essential homology is the same as the usual homology.
$$
\begin{aligned}
H_1(\Gamma_0(\p),\Z_3) =  &  \ 1 \\ 
H_1(\Gamma_0(\pi \cdot \p),\Z_3) = & \   (\Z/3) \\ 
H_1(\Gamma_0(\pi' \cdot \p),\Z_3) =  &  \ (\Z/3) \\ 
H_1(\Gamma_0(3 \p),\Z_3) = &  \  (\Z/3)^6 + (\Z/27) \\ 
H_1(\Gamma_0(B,\p),\Z_3) =  & \  (\Z/3)^5 \end{aligned}$$
We may compute that the classes at level $\pi \cdot \p$ and 
$\pi' \cdot \p$ are both Eisenstein, and hence
are phantom classes (see Definition~\ref{remark:hangingchad}).
Although we were not able to compute
the Hecke action at level $3 \p$, we may still deduce:

\begin{lemma} If $N(\p) = 179$, then
$H^1(\Gamma_0(3 \p),\Z_3)^{\new} \ne H^1(\Gamma_0(B, \p),\Z_3)$.
\end{lemma}

\begin{proof}
The group $H^1(\Gamma_0(B, \p),\Z_3)$ is annihilated by $3$,
 whereas the space of newforms $H^1(\Gamma_0(3 \p),\Z_3)^{\new}$, which is the cokernel of the map
 $$\Psi^{\vee}: (\Z/3)^2 \oplus (\Z/3)^2 \rightarrow (\Z/3)^6 \oplus (\Z/27),$$
is clearly  \emph{not} annihilated by $3$.
\end{proof}

Athough the individual maps  $(\Z/3)^2 \rightarrow H_1(\Gamma_0(3 \p),\Z)$
are injective (by Ihara's Lemma), we believe
(but can not quite compute) that $\Psi^{\vee}$ is not injective, and in particular
that the kernel has order $3^5$ times
$\chi_{D}(\Sigma) =^{?} 3$. This is in accordance with 
Remark~\ref{remark:failure}.
Explicitly, we conjecture that:
$$H^1(\Gamma_0(3 \p),\Z_3)^{\new} = (\Z/3)^4 \oplus (\Z/9),
\qquad H^1(\Gamma_0(B, \p),\Z_3)^{\new} = (\Z/3)^5,$$
the difference in the orders of this group corresponding to
the expected value of  $\chi_D(\Sigma)$.

  \section{Deformations of Eisenstein classes} \label{sec:Eisenstein7} 
  Let us study the (cyclotomic-) Eisenstein
   classes at prime level $\p$ for $F = \Q(\theta)$ with
  $\theta = \sqrt{-2}$,
  using the results of~\S~\ref{section:Eisenstein}.
  Note that $w_F = 2$ and $\Cl(F) = 1$. Moreover, the appropriate
Bernoulli--Hurwitz number in this case has no denominator, so the
$\eps^{-1}$-part of the class group of $F(\zeta_p)$ is trivial for all odd
primes $p$.

\subsection{\texorpdfstring{$p = 3$}{p=3}} We may write $\q = a + b  \theta$.
From~\S~\ref{ext1}, an
Eisenstein
deformation exists if and only if $a + b \theta \equiv \pm 1 \mod 9 \OL_{v}$ for
both primes above $3$, and if $q = N(a + b \theta) \equiv 1 \mod 3$.
Since 
$\theta \equiv \pm 4 \mod 9$ (where the sign depends on the choice of $v|3$), 
we are reduced to the equations
$$a + 4b \equiv \pm 1 \mod 9, \qquad a - 4 b \equiv \pm 1 \mod 9.$$
Writing $q = N(\q)$, 
Respectively adding or subtracting these equations, we arrive at the
congruence relations $a \equiv \pm 1 \mod 9$ and $b \equiv 0 \mod 9$,
or  $a \equiv 0 \mod 9$ and $b \equiv \pm 2 \mod 9$. In the latter case, $q \equiv -1 \mod 9$,
contradicting the assumption that $q \equiv 1 \mod 3$.
Hence a deformation exists whenever $q \equiv 1 \mod 9$ and
$q = a^2 + 162 b^2$. The only primes
$q \le 617$ (that is, within reach of our tables) of this form are $q = 163$ and  $523$.

When $N(\q) = 163$, and $\m$ is Eisenstein, we find that
$$H_1(\Gamma_0(\p),\Z)_{\m} = \Z/3\Z \oplus \Z/729\Z,$$
whereas the congruence cohomology 
is cyclic of order $81$. Moreover, we compute that the $\m$-torsion has
rank $1$. From the self-duality induced from the linking form, we deduce that
$H_1(\Gamma_0(\p),\Z)_{\m}$ is free of rank one as a $\T_{\m}$-module.

\medskip

\subsection{\texorpdfstring{$p = 5$}{p=5}}
Suppose that $p = 5$, so $L = F(\zeta_p) = \Q(\sqrt{-2},\zeta_5)$. We find that
$\Cl(L) = 1$.
If $N(\q) = q$, then a necessary condition
for Eisenstein deformations is that  $q \equiv 1 \mod 5$, and so
$q$ splits completely in $L$.
We explicitly compute the $\eps^{-1}$-invariants of $\OL^{\times}_L$. The invariants
 are cyclic and generated
by some  unit $\beta$. In particular:
$$L = \Q[x]/(x^8 - 2 x^6 + 4 x^4 - 8 x^2 + 16), \quad \beta = 
\frac{3}{4} x^7 + \frac{5}{4} x^6 + \frac{3}{4} x^5 - \frac{5}{2} x^4 - 6 x^3 + 6 x + 13.$$
Even more explicitly, 
$\displaystyle{\beta = 18 - 5 \sqrt{5} - 3 \sqrt{50 - 20 \sqrt{5}}}$
is a root of the polynomial  $z^4 - 72 z^3 + 794 z^2 + 72 z + 1 = 0$.
By Lemma~\ref{lemma:whenexists},
there exists a nontrivial Eisenstein deformation if and only if $\beta$ is  a $p$-th power modulo 
 one prime (equivalently, any prime) above $\q$. The prime $q$
 splits completely in $L$ if   $q \equiv 1,11 \mod 40$. For such primes $q \le 617$,
$\beta$ is a $5$th power only when
 $q = 251$ and $331$.
 Thus, we predict that the Eisenstein part of
 $H_1(\Gamma_0(\p),\Z)$ will be \emph{exactly equal} to the congruence homology
 unless $N(\q) = 251$ or $331$.
 For $N(\q) = 331$, we find that
 $$H_1(\Gamma_0(\q),\Z)_{\m} = \Z/25\Z,$$
 which is clearly larger than the congruence homology (which has order $5$), whereas
 for $N(\q) = 251$, 
 $$H_1(\Gamma_0(\q),\Z)_{\m} = \Z/25\Z \oplus \Z/125\Z,$$
 which is also larger than the congruence homology (the $5$-part of which is cyclic of order $125$). In both cases, the homology is
 free of rank one as a $\T_{\m}$-module.
 
 \section{Non-congruence Eisenstein classes}   \label{section:phantomtwo}
 We expect (\S~\ref{section:phantomclasses})
 that phantom classes exist exactly
 when $\Ext^1(\mu_p,\Z/p\Z)$ contains
 extensions which are unramified everywhere outside $\q$ and
 $N(\q) = q$ is not $1 \mod p$. When $p = 3$, this is equivalent
 to asking that \label{subsection:hangingchads}
$N(\q) = q$ is of the form $81 a^2 +  2 b^2$ where $b \equiv \pm 2 \mod 9$.
The first two primes of this form are $N(\q) = 81  \cdot 1^2 + 2 \cdot 2^2 = 89$ and 
$N(\q) = q = 81 \cdot 1^2 + 2 \cdot 7^2 = 179$, where we \emph{do}
indeed find phantom classes, and which were discussed above.

\subsection{\texorpdfstring{$p = 5$}{p=5}}
In a similar vein, we expect that non-congruence phantom classes
occur for $p = 5$ when $\Ext^1(\mu_5,\Z/5 \Z)$ is non-zero but $q = N(\q) \not\equiv 1 \mod 5$.
Such classes will exist when $q \equiv 9,19 \mod 40$ and when
$\beta \subset L = \Q(\beta)$ is a $5$th power in $\OL_L/{\q}$, for a prime ${\q}$
of norm $q^2$ in $L$. (Note that the corresponding finite field has order $q^2 - 1$, which is divisible
by $5$). The smallest prime with this property is $q = 419$. Indeed, in this case, the
cohomology groups have a similar flavour to what happens for the non-congruence Eisenstein
classes for $p = 3$ and $N(\q) = 89$, specifically;
  $$
 \begin{aligned}
H_1(\Gamma_0(\p),\Z)_{\m} =  &  \ 1 \\ 
H_1(\Gamma_0(\pi \cdot \p),\Z)_{\m} = & \   \Z/5 \\
H_1(\Gamma_0(\pi' \cdot \p),\Z)_{\m} =  &  \  \Z^2_5 \\ 
H_1(\Gamma_0(3 \p),\Z)_{\m}= &  \  \Z^4_{5} \oplus (\Z/5) \\
H_1(\Gamma_0(B, \p),\Z)_{\m} = & \ 1,  \end{aligned}$$
In particular, we deduce as in the proof of Lemma~\ref{lemma:differentnew} 
and subsequent remarks that:
\begin{enumerate}
\item $H_1(\Gamma_0(3 \p),\Z)^{\new}_{\m}$ is nontrivial.
\item $H_1(\Gamma_0(3 \p),\Z)^{\new}_{\m}$ and $H_1(\Gamma_0(B, \p),\Z)_{\m}$ are distinct, and indeed
have different orders.
\item $\chi_D(\Sigma)$ is divisible by $5$.
\end{enumerate}
Although we have not attempted to compute similar examples for primes $p > 5$, presumably
they exist, and conjecturally  correspond exactly to primes $\q$ for which
there exist nontrivial extensions $\Ext^1(\mu_p,\Z/p \Z)$ 
which are unramified away from $\q$ and such that $N(\q) \not\equiv 1 \mod p$ 
(this will force $N(\q) \equiv -1 \mod p$).

\medskip

\section{\texorpdfstring{$K_2(\OL_F)$}{K2} and \texorpdfstring{$F =\Q(\sqrt{-491})$}{F = Q(sqrt(-491))}}
\label{section:K2examples}

Let $F = \Q(\sqrt{-491})$. If $\Gamma = \PGL_2(\OL_F)$, then
$$H_1(\Gamma,\Z) \simeq (\Z/2\Z)^{26}.$$
Theorems~\ref{theorem:K2popularversion} and~\ref{theorem410}
 thus predict that, for each prime $\p$ in $\OL_F$,
 there  exists a non-trivial class in $H^1_{E}(\Gamma_0(\p),\Z/p\Z)$ whenever
 $p > 2$ divides $K_2(\OL_F)$. According to computations of Belabas and Gangl~\cite{BeG},
 there should be an isomorphism 
 $K_2(\OL_F) \simeq \Z/13\Z$, and moreover $|K_2(\OL_F)|$ divides $13$.
 Following a principle espoused by Doron Zeilberger, we produce evidence for this result by computing some
 numerical examples.  
 
 \medskip
 
  The groups $K_2$ for quadratic imaginary fields of discriminant $-D$ for
 (relatively) small $D$ tend to have very small order (mainly $1$ or $2$), in which case our
 Theorems~\ref{theorem:K2popularversion} and~\ref{theorem410} are not so interesting.
  The particular $F$ considered in this section is exceptional in this regard, which
 is the reason for considering it.
  While the data for $F = \Q(\sqrt{-2})$ was computed in $2005$ by Nathan Dunfield using
  a presentation for the corresponding groups found by hand with the help of {\texttt{SnapPea}}, the 
   data for $F = \Q(\sqrt{-491})$ was computed in September $2012$ using Dunfield's  {\texttt{Magma}} script and
   a  presentation for $\PGL_2(\OL_F)$ found for us by Aurel Page using his
   {\texttt{Magma}} package {\texttt{KleinianGroups}}.
   
 \medskip
 
   The table below does indeed confirm that $|K_2(\OL_F)| =^{?} 13$
 divides  the order of
 $$|H^E_1(\Gamma_0(\p),\Z)|$$
  in the indicated range.
 The congruence quotient of $H_1(\Gamma_0(\p),\Z)$ has order $N(\p) - 1$. The only $N(\p)$
 in the table below which are $1 \mod 13$ are $N(\p) = 79$ and $N(\p) = 131$. In both these cases one has
$H_1(\Gamma_0(\p),\Z_{13}) = \Z/13^2 \Z$ and
 $H^{E}_1(\Gamma_0(\p),\Z_{13}) = \Z/13\Z$.
 In the table below, we only give the prime factorization of the order for primes $p < 10000$, in
 part because the other prime divisors were very large, and also because it there was some difficulty
 factorizing the numbers. (Even for $N(\p) = 3$ the order of $H_1$ is divisible by
 primes as large as $30829$.)
 We note that it  follows from the computation of $H^{E}_1(\Gamma_0(\p),\Z)$ for $N(\p) = 3$
 and $N(\p) = 4$
 that the odd part of the order of $K_2(\OL_F)$ divides $13$. However, computing the cohomology
 of arithmetic groups does not seem to be an optimal way to compute (upper bounds) for $K_2(\OL_F)$.
   Note also when $\p = (2)$ and $(7)$, the order of $H_1(\Gamma_0(\p),\Z)$ is a perfect square
   up to the orbifold primes  $2$ and $3$, in accordance with
   Theorem~\ref{theorem:square}. 
   
   \medskip
   
For each prime $p$ which splits in $F$, there are two primes $\p$ of norm $p$, but
because of complex conjugation, the structure of $H_1(\Gamma_0(\p),\Z)$  does not depend on the choice.
On the other hand, if $p$ and $q$ are two such primes, then there are four ideals of norm $pq$, and there
are two distinct isomorphism classes.
  Let $\alpha = \frac{1 + \sqrt{-491}}{2}$.
Suppose that $\p = (3,2 \alpha)$ and $\q = (11, 2 \alpha + 1)$ are ideals of norm $3$ and $11$
respectively. Then the ideal $\p \q$ has norm $33$ and the  small primes dividing the
 homology of $\Gamma_0(\p \q)$ are given in the table below under $33A$. At the prime $13$,
 one more explicitly finds that
$$H_1(Y_0(\p \q),\Z_{13}) \simeq (\Z/13)^2 \oplus (\Z/13^2).$$
It follows from
Theorem~\ref{TJL2}  that, if $\Gamma_B$ denotes the arithmetic
group arising from the non-split quaternion algebra ramified at
$\p$ and $\q$ then $H_1(\Gamma_B,\Z_{13})$ is trivial. 
On the other hand, the space of newforms $H^{E}_1(\Gamma_0(\p \q),\Z_{13})^{\p\q-\new}$
 is the cokernel of the map:
$$
\begin{diagram}
 H^E_1(\Gamma_0(\p),\Z_{13})^2
\oplus H^E_1(\Gamma_0(\q,\Z_{13})^2  & \rTo & 
H^E_1(\Gamma_0(\p \q),\Z_{13}) \\
\dEquals & & \dEquals \\
(\Z/13\Z)^4 & \rTo &  (\Z/13\Z)^2 \oplus (\Z/13^2 \Z),\end{diagram}$$
which, by virtue of the structure of the groups in question, has order divisible by $13$. In particular, 
$H^{E}_1(\Gamma_0(\p \q),\Z_{13})^{\p\q-\new}$ is divisible by $13$, and the level raising sequence is not exact and
has homology of order divisible by $13$, providing  an example
of Theorem $A^{\dagger}$ (Theorem~\ref{theorem:theoremD}).

 \section{Tables} \label{section:tables}
We reproduce here three tables of data.

\medskip

The first table records the computations obtained for $F = \Q(\sqrt{-491})$  made in September 2012 using
Aurel Page's \texttt{KleinianGroups} package to compute
a presentation of $\PGL_2(\OL_F)$ and a  \texttt{magma} script of Nathan Dunfield to compute
the homology of congruence subgroups.
 The second set of data   was computed for us by Nathan Dunfield in November of 2005. It records the 
first  homology of the four manifolds $Y_0(\p)$, $Y_0(\pi \cdot \p)$, $Y_0(\pibar \cdot \p)$, $Y_0(3 \cdot \p)$ corresponding
 to the congruence subgroups of $\GL_2(\OL_F)$ with $F = \Q(\sqrt{-2})$ and $\p$ a prime of degree one in $\OL_F$, as well as the first homology of
 the congruence subgroup $Y_{B,0}(\p)$ of the manifold $Y_B$ corresponding to the quaternion algebra ramified
 over $\Q(\sqrt{-2})$ exactly at $\pi$ and $\pibar$. The groups are recorded by their abelian invariants, so 
 $(81,2)$ corresponds to $(\Z/81\Z)^2$ and $(0,3)$ to $\Z^3$.
 
\newpage

 The homology of congruence subgroups of $\PGL_2(\OL_F)$ for $F = \Q(\sqrt{-491})$.
   
   \begin{center}
\begin{tabular}{|c|rl|}
\hline
$N(\p)$ &  odd prime power divisors & of $\mathrm{gcd}(10000!,|H_1(\Gamma_0(\p),\Z)|)$ \\
\hline
$3$ &  $5$, $7^2$,&$\mathbf{13}$,  \\
$4$ & $3$,  &$\mathbf{13}^4$, \\ 
$5$ & $7^2$, & $\mathbf{13}^3$, $23$, $29$, $31$, $293$ \\
$11$ & $5^2$, $11$, & $\mathbf{13}$, $17^3$, $101$ \\
$13$ & 
 $11$,
& $\mathbf{13}$, $61$, $607$, $3529$ \\
$17$ & $7$, & $\mathbf{13}^4$, $89$, $103$, $173$, $1069$  \\
$31$ & $3$, $5^2$, $7$, & $\mathbf{13}$, $1657$, $2887$ \\
$37$ & $3^2$, $5^2$, & $\mathbf{13}$, $19^2$, $97$, $131$, $2539$, $4639$, $7331$ \\
$41$ & $5^2$, & $\mathbf{13}$, $17$, $71$ \\
$43$ & $3$, $7^2$, & $\mathbf{13}$, $19$, $113$, $191$, $347$, $9719$ \\
$49$ & $3^3$, $5^2$, $11^2$, & $\mathbf{13}^2$,
$31^2$, $53^2$, $59^2$,  $547^2$,  $4337^2$\\
$61$ & $3^4$, $5^4$, & $\mathbf{13}$ \\
$71$ & $5^2$, $7$, & $\mathbf{13}$, $29$, $197$ \\
$79$ & $3$, $5$, $7$, & $\mathbf{13^2}$, $2543$ \\
$83$ & & $\mathbf{13}$, $41^3$, $73$ \\
$97$ & $3$, $5^3$, $11$, & $\mathbf{13}$, $53$, $229$, $383$, $607$, $2039$, $5441$ \\
$101$ & $3^6$, $5^4$, $7$, &  $\mathbf{13}$ \\
$107$ & $5$, $7^2$,  & $\mathbf{13}$, $31$, $53$ \\
$127$ & $3^7$, $7$,  & $\mathbf{13}^2$ \\
$131$ & $5$, $7$,  & $\mathbf{13^2}$, $23$ \\
$139$ & $3$, $5$, & $\mathbf{13}^2$, $23$, $97$, $193$, $2341$ \\
$163$ & $3^6$, $5$, &  $\mathbf{13}$, $47$, $4241$ \\
$179$ & $3^5$, &  $\mathbf{13}^2$, $71$, $79$, $89$, $139$ \\
$181$ & $3^6$, $5$, $7^2$,  &  $\mathbf{13}$, $41$, $59$ \\
$197$ & $7^2$  &  $\mathbf{13}$, $47$, $61$, $1187$ \\
$199$ & $3^2$, $5$, $11^2$,  &  $\mathbf{13}$, $23$, $29^2$, $2269$ \\
$223$ & $3$, $5$, $7$, &  $\mathbf{13}$, $17$, $37^2$, $41$, $223$,
$449$, $701$ \\
$227$ & $5^2$, $11$, &  $\mathbf{13}$, $113$, $607$ \\
$229$ & $3$, $5$, $7$, &  $\mathbf{13}^3$, $19$,  $47$ \\
$233$ & $3^6$, $5^2$,  &  $\mathbf{13}^2$, $29$, $307$ \\
$239$ & $3^6$, $5$, $7^3$, &  $\mathbf{13}$, $17$,  $19$, $197$, $241$, $6569$ \\
$241$ &  $3^3$, $5^2$, $11$,  &  $\mathbf{13}$ \\
$33A$ & $3$, $5^6$, $7^7$, $11^2$,   &  $\mathbf{13^4}$ , $17^6$, $19$, $97$, $101^2$, $239$, \\
& & $311$,
$3167$, $3209$, $4447$ \\
$33B$ &  & $\infty$   \\
\hline
\end{tabular}
\end{center}

\newpage

    \begin{center}
\begin{tabular}{|l|l|}
\hline
$N(\p) = 11$ &  \\
\hline
$ H_1(Y_0(\p))$ &  $ (0, 0), (2, 3), (5, 1) $ \\
$ H_1(Y_0(\pi \cdot \p))$ & $ (0, 0), (2, 6), (5, 1) $ \\
$H_1(Y_0(\pibar \cdot \p))$ & $  (0, 0), (2, 6), (3, 1), (5, 1) $ \\
$H_1(Y_0(3 \cdot \p))$ & $  (0, 1), (2, 9), (3, 1), (4, 2), (5, 1) $ \\
$H_1(Y_{B,0}(\p))$ & $  (0, 1), (2, 3), (4, 1), (5, 1) $ \\
  \hline
  $N(\p) = 17$ &  \\
\hline
$H_1(Y_0(\p))$ & $  (0, 0), (2, 4), (8, 1) $ \\
$H_1(Y_0(\pi \cdot \p))$ & $  (0, 1), (2, 3), (4, 1), (16, 1) $ \\
$H_1(Y_0(\pibar \cdot \p))$ & $  (0, 0), (2, 3), (4, 2), (16, 1) $ \\
$H_1(Y_0(3 \cdot \p))$ & $  (0, 2), (2, 9), (4, 1), (8, 1), (64, 1) $ \\
$H_1(Y_{B,0}(\p))$ & $  (0, 0), (2, 4), (3, 1), (4, 3), (32, 1) $ \\
\hline  
$N(\p) =  19 $ & \\
\hline 
$ H_1(Y_0(\p))$  & $  (0, 0), (2, 3), (9, 1) $ \\
$H_1(Y_0(\pi \cdot \p))$ & $  (0, 0), (2, 6), (9, 1) $ \\
$H_1(Y_0(\pibar \cdot \p))$ & $  (0, 0), (2, 6), (9, 1) $ \\
$H_1(Y_0(3 \cdot \p))$& $  (0, 0), (2, 9), (4, 2), (8, 1), (9, 1) $ \\
$H_1(Y_{B,0}(\p))$ & $  (0, 0), (2, 2), (4, 2), (9, 1), (16, 1) $ \\
\hline  
$N(\p) =  41 $ & \\
 \hline 
 $ H_1(Y_0(\p))$  & $  (0, 0), (2, 5), (4, 1), (5, 1) $ \\
$H_1(Y_0(\pi \cdot \p))$ & $  (0, 0), (2, 7), (4, 2), (5, 1), (8, 1) $ \\
$H_1(Y_0(\pibar \cdot \p))$ & $  (0, 0), (2, 9), (4, 1), (5, 1), (8, 1) $ \\
$H_1(Y_0(3 \cdot \p))$& $  (0, 0), (2, 19), (4, 6), (5, 1), (19, 1), (32, 1) $ \\
$H_1(Y_{B,0}(\p))$ & $  (0, 0), (2, 4), (4, 3), (5, 1), (16, 1), (19, 1) $ \\
\hline  
$N(\p) =  43 $ & \\
 \hline 
 $ H_1(Y_0(\p))$  & $  (0, 0), (2, 3), (3, 1), (7, 1) $ \\
$H_1(Y_0(\pi \cdot \p))$ & $  (0, 0), (2, 6), (3, 1), (7, 1), (13, 1) $ \\
$H_1(Y_0(\pibar \cdot \p))$ & $  (0, 0), (2, 7), (3, 1), (7, 1), (8, 2) $ \\
$H_1(Y_0(3 \cdot \p))$& $  (0, 0), (2, 11), (3, 1), (4, 2), (7, 1), (8, 5), (13, 2) $ \\
$H_1(Y_{B,0}(\p))$ & $  (0, 0), (2, 2), (3, 1), (4, 1), (7, 1), (8, 2) $ \\
\hline 
$N(\p) =  59 $ &  \\
\hline
 $ H_1(Y_0(\p))$  & $  (0, 0), (2, 3), (29, 1) $ \\
$H_1(Y_0(\pi \cdot \p))$ & $  (0, 0), (2, 6), (3, 2), (7, 1), (29, 1) $ \\
$H_1(Y_0(\pibar \cdot \p))$ & $  (0, 0), (2, 6), (3, 2), (29, 1), (37, 1) $ \\
$H_1(Y_0(3 \cdot \p))$& $  (0, 0), (2, 10), (3, 8), (4, 3), (7, 2), (9, 1), (16, 1), (29, 1), (37, 2) $ \\
$H_1(Y_{B,0}(\p))$ & $  (0, 0), (2, 4), (4, 1), (9, 1), (16, 2), (29, 1) $ \\

\hline  $N(\p) =  67 
$ & \\ \hline $ H_1(Y_0(\p))$  & $  (0, 0), (2, 4), (3, 2), (11, 1) $ \\
$H_1(Y_0(\pi \cdot \p))$ & $  (0, 0), (2, 8), (3, 3), (7, 1), (11, 2) $ \\
$H_1(Y_0(\pibar \cdot \p))$ & $  (0, 0), (2, 8), (3, 3), (9, 1), (11, 1), (17, 1) $ \\
$H_1(Y_0(3 \cdot \p))$& $  (0, 1), (2, 16), (3, 4), (4, 2), (7, 1), (9, 2), (11, 3), (17, 2) $ \\
$H_1(Y_{B,0}(\p))$ & $  (0, 1), (2, 4), (3, 1), (4, 1), (8, 1), (11, 1) $ \\
\hline
\end{tabular}

\begin{tabular}{|l|l|}

\hline  $N(\p) =  73 
$ & \\ \hline $ H_1(Y_0(\p))$  & $  (0, 0), (2, 4), (4, 1), (9, 1), (19, 1) $ \\
$H_1(Y_0(\pi \cdot \p))$ & $  (0, 0), (2, 4), (4, 1), (8, 1), (9, 1), (19, 2), (32, 1) $ \\
$H_1(Y_0(\pibar \cdot \p))$ & $  (0, 1), (2, 4), (3, 1), (4, 1), (9, 1), (16, 1), (19, 2) $ \\
$H_1(Y_0(3 \cdot \p))$& $  (0, 3), (2, 8), (3, 1), (4, 3), (9, 1), (11, 1), (16, 2), (19, 4), (32, 1) $ \\
 & $ (256, 1) $ \\
$H_1(Y_{B,0}(\p))$ & $  (0, 1), (2, 4), (4, 4), (8, 1), (9, 1), (11, 1), (32, 1) $ \\

\hline  $N(\p) =  83 
$ & \\ \hline $ H_1(Y_0(\p))$  & $  (0, 0), (2, 3), (4, 1), (41, 1) $ \\
$H_1(Y_0(\pi \cdot \p))$ & $  (0, 1), (2, 6), (3, 1), (7, 1), (11, 1), (16, 1), (41, 1) $ \\
$H_1(Y_0(\pibar \cdot \p))$ & $  (0, 0), (2, 6), (4, 2), (41, 1) $ \\
$H_1(Y_0(3 \cdot \p))$& $  (0, 3), (2, 9), (3, 4), (4, 2), (7, 1), (11, 2), (16, 2), (41, 1), (49, 1) $ \\
$H_1(Y_{B,0}(\p))$ & $  (0, 1), (2, 3), (3, 2), (4, 1), (5, 1), (7, 1), (41, 1) $ \\
\hline  $N(\p) =  89 
$ & \\ \hline $ H_1(Y_0(\p))$  & $  (0, 0), (2, 4), (4, 1), (11, 2) $ \\
$H_1(Y_0(\pi \cdot \p))$ & $  (0, 2), (2, 3), (3, 1), (4, 3), (11, 2) $ \\
$H_1(Y_0(\pibar \cdot \p))$ & $  (0, 1), (2, 4), (4, 1), (11, 3), (16, 1) $ \\
$H_1(Y_0(3 \cdot \p))$& $  (0, 7), (2, 10), (3, 1), (4, 1), (5, 2), (7, 1), (8, 1), (11, 3), (32, 1) $ \\
$H_1(Y_{B,0}(\p))$ & $  (0, 1), (2, 4), (4, 4), (5, 1), (7, 1), (11, 1), (16, 1), (25, 1) $ \\

\hline  $N(\p) =  97 
$ & \\ \hline $ H_1(Y_0(\p))$  & $  (0, 0), (2, 4), (3, 1), (16, 1), (41, 1) $ \\
$H_1(Y_0(\pi \cdot \p))$ & $  (0, 0), (2, 3), (3, 1), (4, 2), (17, 1), (19, 1), (32, 1), (41, 2) $ \\
$H_1(Y_0(\pibar \cdot \p))$ & $  (0, 0), (2, 3), (3, 2), (4, 2), (7, 1), (32, 1), (41, 2) $ \\
$H_1(Y_0(3 \cdot \p))$& $  (0, 0), (2, 13), (3, 2), (4, 2), (7, 3), (8, 2), (9, 1), (16, 2), (17, 2) $ \\
  &  $  (19, 2), (41, 4), (128, 1) $ \\
$H_1(Y_{B,0}(\p))$ & $  (0, 0), (2, 8), (3, 2), (4, 3), (7, 1), (16, 2), (64, 1) $ \\

\hline  $N(\p) =  107 
$ & \\ \hline $ H_1(Y_0(\p))$  & $  (0, 0), (2, 3), (53, 1) $ \\
$H_1(Y_0(\pi \cdot \p))$ & $  (0, 0), (2, 8), (4, 1), (53, 1), (61, 1) $ \\
$H_1(Y_0(\pibar \cdot \p))$ & $  (0, 0), (2, 6), (53, 1), (67, 1), (163, 1) $ \\
$H_1(Y_0(3 \cdot \p))$& $  (0, 1), (2, 14), (4, 5), (17, 1), (53, 1), (61, 2), (67, 2), (163, 2) $ \\
$H_1(Y_{B,0}(\p))$ & $  (0, 1), (2, 4), (4, 1), (8, 1), (17, 1), (53, 1) $ \\

\hline  $N(\p) =  113 
$ & \\ \hline $ H_1(Y_0(\p))$  & $  (0, 0), (2, 4), (4, 2), (7, 1), (8, 1) $ \\
$H_1(Y_0(\pi \cdot \p))$ & $  (0, 1), (2, 4), (4, 1), (7, 1), (8, 5), (11, 1) $ \\
$H_1(Y_0(\pibar \cdot \p))$ & $  (0, 0), (2, 4), (3, 2), (4, 3), (7, 1), (8, 2), (64, 1) $ \\
$H_1(Y_0(3 \cdot \p))$& $  (0, 3), (2, 9), (3, 6), (7, 1), (8, 7), (9, 1), (11, 2), (16, 3), (32, 1) $ \\
 & $ (64, 1) $ \\
$H_1(Y_{B,0}(\p))$ & $  (0, 1), (2, 4), (3, 3), (4, 5), (7, 1), (9, 1), (64, 1) $ \\
\hline
\end{tabular}

\begin{tabular}{|l|l|}

\hline  $N(\p) =  131 
$ & \\ \hline $ H_1(Y_0(\p))$  & $  (0, 0), (2, 3), (5, 1), (13, 1) $ \\
$H_1(Y_0(\pi \cdot \p))$ & $  (0, 0), (2, 6), (5, 2), (13, 1), (19, 1), (71, 1) $ \\
$H_1(Y_0(\pibar \cdot \p))$ & $  (0, 0), (2, 6), (3, 1), (5, 1), (7, 1), (13, 1), (179, 1) $ \\
$H_1(Y_0(3 \cdot \p))$& $  (0, 1), (2, 9), (3, 1), (4, 3), (5, 4), (7, 3), (8, 1), (11, 2), (13, 1), (19, 2) $ \\
 & $ (31, 1), (71, 2), (179, 2) $ \\
$H_1(Y_{B,0}(\p))$ & $  (0, 1), (2, 2), (4, 3), (5, 2), (7, 1), (8, 1), (11, 2), (13, 1), (31, 1) $ \\

\hline  $N(\p) =  137 
$ & \\ \hline $ H_1(Y_0(\p))$  & $  (0, 0), (2, 7), (4, 2), (17, 1) $ \\
$H_1(Y_0(\pi \cdot \p))$ & $  (0, 0), (2, 11), (4, 4), (8, 1), (17, 1), (83, 1) $ \\
$H_1(Y_0(\pibar \cdot \p))$ & $  (0, 0), (2, 7), (4, 3), (7, 1), (8, 5), (17, 1) $ \\
$H_1(Y_0(3 \cdot \p))$& $  (0, 0), (2, 15), (4, 10), (5, 2), (7, 3), (8, 8), (11, 1), (17, 1), (32, 1) $ \\
 & $ (83, 2), (443, 1) $ \\
$H_1(Y_{B,0}(\p))$ & $  (0, 0), (2, 4), (4, 3), (5, 2), (7, 1), (11, 1), (16, 1), (17, 1), (443, 1) $ \\

\hline  $N(\p) =  139 
$ & \\ \hline $ H_1(Y_0(\p))$  & $  (0, 0), (2, 3), (3, 1), (23, 1), (25, 1) $ \\
$H_1(Y_0(\pi \cdot \p))$ & $  (0, 0), (2, 6), (3, 1), (13, 1), (23, 1), (25, 1), (125, 1) $ \\
$H_1(Y_0(\pibar \cdot \p))$ & $  (0, 0), (2, 6), (3, 2), (23, 2), (25, 1), (125, 1) $ \\
$H_1(Y_0(3 \cdot \p))$& $  (0, 1), (2, 10), (3, 4), (4, 2), (5, 1), (7, 1), (8, 1), (9, 1), (11, 1), (13, 2) $ \\
 & $ (23, 3), (25, 2), (41, 1), (125, 2) $ \\
$H_1(Y_{B,0}(\p))$ & $  (0, 1), (2, 4), (3, 2), (4, 1), (5, 1), (7, 1), (9, 1), (11, 1), (16, 1), (23, 1) $ \\
 & $ (41, 1) $ \\

\hline  $N(\p) =  163 
$ & \\ \hline $ H_1(Y_0(\p))$  & $  (0, 0), (2, 3), (3, 1), (729, 1) $ \\
$H_1(Y_0(\pi \cdot \p))$ & $  (0, 0), (2, 6), (3, 3), (13, 1), (49, 1), (59, 1), (81, 1), (729, 1) $ \\
$H_1(Y_0(\pibar \cdot \p))$ & $  (0, 0), (2, 6), (3, 1), (27, 1), (31, 1), (41, 1), (229, 1), (729, 1) $ \\
$H_1(Y_0(3 \cdot \p))$& $  (0, 0), (2, 9), (3, 5), (4, 2), (8, 1), (9, 2), (13, 2), (27, 1), (29, 1), (31, 2) $ \\
 & $ (41, 2), (49, 2), (59, 2), (81, 1), (229, 2), (729, 1), (6561, 1) $ \\
$H_1(Y_{B,0}(\p))$ & $  (0, 0), (2, 2), (3, 5), (4, 1), (8, 2), (29, 1), (2187, 1) $ \\

\hline  $N(\p) =  179 
$ & \\ \hline $ H_1(Y_0(\p))$  & $  (0, 0), (2, 3), (31, 1), (89, 1) $ \\
$H_1(Y_0(\pi \cdot \p))$ & $  (0, 0), (2, 6), (3, 1), (13, 1), (31, 2), (89, 1), (769, 1) $ \\
$H_1(Y_0(\pibar \cdot \p))$ & $  (0, 0), (2, 6), (3, 1), (5, 1), (13, 1), (31, 2), (49, 1), (89, 1), (107, 1) $ \\
$H_1(Y_0(3 \cdot \p))$& $  (0, 0), (2, 12), (3, 6), (4, 2), (5, 3), (13, 4), (19, 1), (27, 1), (31, 4) $ \\
 & $ (32, 1), (49, 2), (89, 1), (97, 1), (107, 2), (769, 2) $ \\
$H_1(Y_{B,0}(\p))$ & $  (0, 0), (2, 4), (3, 5), (4, 1), (5, 1), (8, 1), (19, 1), (64, 1), (89, 1), (97, 1) $ \\

\hline  $N(\p) =  193 
$ & \\ \hline $ H_1(Y_0(\p))$  & $  (0, 0), (2, 4), (3, 1), (32, 1), (251, 1) $ \\
$H_1(Y_0(\pi \cdot \p))$ & $  (0, 0), (2, 3), (3, 1), (4, 1), (7, 1), (9, 1), (61, 1), (64, 1), (128, 1) $ \\
 & $ (251, 2) $ \\
$H_1(Y_0(\pibar \cdot \p))$ & $  (0, 0), (2, 3), (3, 1), (4, 1), (5, 1), (64, 1), (67, 1), (251, 2), (512, 1) $ \\
$H_1(Y_0(3 \cdot \p))$& $  (0, 5), (2, 12), (3, 2), (4, 3), (5, 2), (7, 1), (8, 1), (9, 1), (32, 1), (61, 2) $ \\
 & $ (67, 1), (251, 4), (256, 2), (293, 1), (2048, 1) $ \\
$H_1(Y_{B,0}(\p))$ & $  (0, 5), (2, 8), (3, 1), (4, 2), (8, 3), (128, 1), (293, 1) $ \\
\hline
\end{tabular}

\begin{tabular}{|l|l|}

\hline  $N(\p) =  211 
$ & \\ \hline $ H_1(Y_0(\p))$  & $  (0, 0), (2, 3), (3, 1), (5, 2), (7, 1) $ \\
$H_1(Y_0(\pi \cdot \p))$ & $  (0, 0), (2, 6), (3, 1), (5, 2), (7, 1), (25, 1), (1759, 1) $ \\
$H_1(Y_0(\pibar \cdot \p))$ & $  (0, 1), (2, 9), (3, 1), (5, 3), (7, 1), (73, 1), (3691, 1) $ \\
$H_1(Y_0(3 \cdot \p))$& $  (0, 2), (2, 10), (3, 1), (4, 7), (5, 3), (7, 1), (13, 1), (16, 2), (23, 1), (25, 2) $ \\
 & $ (32, 1), (73, 2), (1759, 2), (3691, 2) $ \\
$H_1(Y_{B,0}(\p))$ & $  (0, 0), (2, 6), (3, 1), (4, 1), (5, 1), (7, 1), (8, 5), (13, 1), (23, 1), (64, 1) $ \\

\hline  $N(\p) =  227 
$ & \\ \hline $ H_1(Y_0(\p))$  & $  (0, 0), (2, 3), (37, 1), (113, 1) $ \\
$H_1(Y_0(\pi \cdot \p))$ & $  (0, 0), (2, 8), (4, 2), (7, 2), (37, 2), (41, 1), (113, 1) $ \\
$H_1(Y_0(\pibar \cdot \p))$ & $  (0, 0), (2, 6), (37, 2), (113, 1), (157, 1), (10487, 1) $ \\
$H_1(Y_0(3 \cdot \p))$& $  (0, 0), (2, 13), (4, 6), (7, 5), (8, 1), (9, 2), (37, 4), (41, 2), (113, 1) $ \\
 & $  (157, 2), (257, 1), (397, 1), (10487, 2) $ \\
$H_1(Y_{B,0}(\p))$ & $  (0, 0), (2, 3), (4, 1), (7, 1), (9, 2), (32, 1), (113, 1), (257, 1), (397, 1) $ \\

\hline  $N(\p) =  233 
$ & \\ \hline $ H_1(Y_0(\p))$  & $  (0, 0), (2, 4), (4, 1), (29, 1), (431, 1) $ \\
$H_1(Y_0(\pi \cdot \p))$ & $  (0, 0), (2, 4), (3, 3), (4, 1), (5, 1), (8, 1), (19, 1), (29, 1), (32, 1), (431, 2) $ \\
$H_1(Y_0(\pibar \cdot \p))$ & $  (0, 1), (2, 4), (4, 3), (7, 1), (29, 1), (83, 1), (347, 1), (431, 2) $ \\
$H_1(Y_0(3 \cdot \p))$& $  (0, 3), (2, 10), (3, 7), (4, 4), (5, 1), (7, 2), (8, 1), (11, 1), (19, 2), (29, 2) $ \\
 & $  (32, 1), (41, 1), (61, 1), (64, 1), (83, 2), (347, 2), (431, 4), (6427, 1) $ \\
$H_1(Y_{B,0}(\p))$ & $  (0, 1), (2, 4), (3, 1), (4, 4), (11, 1), (16, 1), (29, 2), (41, 1), (61, 1) $ \\
 & $ (6427, 1) $ \\
\hline  $N(\p) =  241 
$ & \\ \hline $ H_1(Y_0(\p))$  & $  (0, 0), (2, 4), (3, 1), (5, 1), (8, 1), (1861, 1) $ \\
$H_1(Y_0(\pi \cdot \p))$ & $  (0, 3), (2, 5), (3, 3), (4, 1), (5, 3), (16, 1), (1861, 2) $ \\
$H_1(Y_0(\pibar \cdot \p))$ & $  (0, 0), (2, 3), (3, 2), (4, 1), (5, 1), (11, 1), (16, 1), (64, 1), (373, 1) $ \\
 & $ (1861, 2) $ \\
$H_1(Y_0(3 \cdot \p))$& $  (0, 8), (2, 11), (3, 5), (4, 2), (5, 4), (8, 1), (11, 1), (13, 2), (16, 1), (32, 1) $ \\
 & $ (64, 1), (107, 1), (373, 2), (839, 1), (1153, 1), (1861, 4) $ \\
$H_1(Y_{B,0}(\p))$ & $  (0, 2), (2, 3), (3, 1), (4, 4), (5, 1), (8, 1), (11, 1), (13, 2), (16, 1), (107, 1) $ \\
 & $ (839, 1), (1153, 1) $ \\

\hline  $N(\p) =  251 
$ & \\ \hline $ H_1(Y_0(\p))$  & $  (0, 0), (2, 3), (17, 1), (25, 1), (625, 1) $ \\
$H_1(Y_0(\pi \cdot \p))$ & $  (0, 0), (2, 8), (8, 2), (13, 1), (17, 2), (23, 1), (25, 1), (125, 1), (151, 1) $ \\
 & $ (625, 1) $ \\
$H_1(Y_0(\pibar \cdot \p))$ & $  (0, 0), (2, 6), (3, 2), (5, 1), (9, 1), (17, 2), (25, 1), (125, 1), (251, 1) $ \\
 & $ (625, 1) $ \\
$H_1(Y_0(3 \cdot \p))$& $  (0, 0), (2, 19), (3, 6), (4, 2), (5, 2), (8, 4), (9, 2), (13, 2), (16, 1), (17, 5) $ \\
 & $ (23, 2), (25, 1), (64, 2), (125, 2), (151, 2), (163, 1), (251, 2), (625, 2) $ \\
$H_1(Y_{B,0}(\p))$ & $  (0, 0), (2, 7), (3, 2), (4, 3), (17, 1), (32, 2), (64, 1), (163, 1), (625, 1) $ \\
\hline
\end{tabular}

\begin{tabular}{|l|l|}
\hline  $N(\p) =  257 
$ & \\ \hline $ H_1(Y_0(\p))$  & $  (0, 0), (2, 8), (49, 1), (128, 1) $ \\
$H_1(Y_0(\pi \cdot \p))$ & $  (0, 2), (2, 6), (4, 6), (7, 1), (8, 2), (47, 1), (49, 1), (128, 1) $ \\
$H_1(Y_0(\pibar \cdot \p))$ & $  (0, 0), (2, 14), (3, 1), (4, 3), (49, 2), (256, 1), (5821, 1) $ \\
$H_1(Y_0(3 \cdot \p))$& $  (0, 4), (2, 25), (3, 1), (4, 12), (7, 2), (8, 8), (9, 2), (47, 2), (49, 2), (67, 1) $ \\
 & $ (71, 1), (181, 1), (512, 1), (5821, 2) $ \\
$H_1(Y_{B,0}(\p))$ & $  (0, 0), (2, 8), (3, 2), (4, 5), (8, 1), (9, 1), (16, 1), (67, 1), (71, 1), (181, 1) $ \\
 & $ (512, 1) $ \\

\hline  $N(\p) =  281 
$ & \\ \hline $ H_1(Y_0(\p))$  & $  (0, 0), (2, 8), (4, 1), (5, 1), (49, 1) $ \\
$H_1(Y_0(\pi \cdot \p))$ & $  (0, 0), (2, 12), (4, 3), (5, 1), (7, 1), (29, 1), (49, 1), (128, 1), (251, 1) $ \\
$H_1(Y_0(\pibar \cdot \p))$ & $  (0, 0), (2, 4), (4, 5), (5, 1), (7, 1), (8, 5), (49, 1), (128, 1) $ \\
$H_1(Y_0(3 \cdot \p))$& $  (0, 0), (2, 11), (4, 10), (5, 2), (7, 2), (8, 10), (16, 1), (29, 2), (41, 1) $ \\
 & $ (49, 2), (71, 1), (128, 1), (149, 1), (197, 1), (251, 2), (256, 3), (347, 1) $ \\
$H_1(Y_{B,0}(\p))$ & $  (0, 0), (2, 4), (4, 4), (5, 2), (16, 2), (41, 1), (49, 1), (71, 1), (149, 1) $ \\
 & $ (197, 1), (347, 1) $ \\

\hline  $N(\p) =  283 
$ & \\ \hline $ H_1(Y_0(\p))$  & $  (0, 0), (2, 6), (3, 1), (47, 2) $ \\
$H_1(Y_0(\pi \cdot \p))$ & $  (0, 0), (2, 16), (3, 1), (4, 4), (47, 3), (121, 1) $ \\
$H_1(Y_0(\pibar \cdot \p))$ & $  (0, 0), (2, 12), (3, 1), (7, 2), (23, 1), (47, 3), (271, 1), (1409, 1) $ \\
$H_1(Y_0(3 \cdot \p))$& $  (0, 0), (2, 29), (3, 1), (4, 11), (7, 4), (11, 1), (16, 1), (17, 2), (23, 2) $ \\
 & $ (32, 1), (47, 5), (89, 1), (121, 2), (271, 2), (1409, 2) $ \\
$H_1(Y_{B,0}(\p))$ & $  (0, 0), (2, 2), (3, 1), (4, 3), (11, 1), (16, 1), (17, 2), (47, 1), (64, 1) $ \\
 & $ (89, 1) $ \\

\hline  $N(\p) =  307 
$ & \\ \hline $ H_1(Y_0(\p))$  & $  (0, 0), (2, 5), (9, 1), (17, 1), (97, 1) $ \\
$H_1(Y_0(\pi \cdot \p))$ & $  (0, 0), (2, 12), (3, 2), (4, 1), (7, 3), (9, 1), (11, 1), (17, 1), (97, 2) $ \\
 & $ (107, 1) $ \\
$H_1(Y_0(\pibar \cdot \p))$ & $  (0, 0), (2, 10), (3, 1), (5, 3), (7, 2), (9, 1), (17, 1), (25, 1), (31, 1), (97, 2) $ \\
$H_1(Y_0(3 \cdot \p))$& $  (0, 3), (2, 21), (3, 7), (4, 5), (5, 6), (7, 7), (9, 2), (11, 2), (16, 2), (17, 1) $ \\
 & $ (19, 1), (25, 1), (31, 2), (49, 1), (97, 4), (103, 1), (107, 2), (125, 1) $ \\
$H_1(Y_{B,0}(\p))$ & $  (0, 3), (2, 6), (3, 3), (4, 1), (5, 1), (8, 2), (9, 1), (16, 1), (17, 1), (19, 1) $ \\
 & $ (103, 1) $ \\

\hline  $N(\p) =  313 
$ & \\ \hline $ H_1(Y_0(\p))$  & $  (0, 0), (2, 5), (3, 1), (4, 1), (13, 1), (79, 1), (101, 1) $ \\
$H_1(Y_0(\pi \cdot \p))$ & $  (0, 0), (2, 9), (3, 1), (4, 1), (5, 1), (8, 1), (11, 1), (13, 1), (25, 1), (43, 1) $ \\
 & $ (79, 2), (101, 2) $ \\
$H_1(Y_0(\pibar \cdot \p))$ & $  (0, 0), (2, 7), (3, 3), (4, 1), (8, 3), (13, 1), (79, 2), (101, 2), (103, 1) $ \\
 & $ (1033, 1) $ \\
$H_1(Y_0(3 \cdot \p))$& $  (0, 1), (2, 14), (3, 6), (4, 6), (5, 2), (8, 9), (11, 3), (13, 1), (16, 1), (25, 2) $ \\
 & $ (37, 1), (43, 2), (47, 2), (71, 1), (79, 4), (101, 4), (103, 2), (128, 1) $ \\
  & $ (1033, 2) $ \\
$H_1(Y_{B,0}(\p))$ & $  (0, 1), (2, 6), (3, 2), (4, 5), (5, 1), (8, 4), (11, 1), (13, 1), (37, 1), (47, 2) $ \\
 & $ (64, 1), (71, 1) $ \\
\hline
\end{tabular}

\begin{tabular}{|l|l|}

\hline  $N(\p) =  331 
$ & \\ \hline $ H_1(Y_0(\p))$  & $  (0, 0), (2, 3), (3, 1), (11, 2), (25, 1), (199, 1) $ \\
$H_1(Y_0(\pi \cdot \p))$ & $  (0, 1), (2, 9), (3, 2), (11, 3), (23, 1), (25, 1), (149, 1), (197, 1), (199, 2) $ \\
$H_1(Y_0(\pibar \cdot \p))$ & $  (0, 0), (2, 6), (3, 1), (5, 2), (11, 3), (25, 1), (199, 2), (317, 1) $ \\
 & $ (27277, 1) $ \\
$H_1(Y_0(3 \cdot \p))$& $  (0, 2), (2, 16), (3, 4), (4, 3), (5, 3), (11, 6), (19, 1), (23, 2), (25, 1) $ \\
 & $ (32, 1), (47, 1), (113, 1), (149, 2), (197, 2), (199, 4), (317, 2), (433, 1) $ \\
  & $  (991, 1), (27277, 2) $ \\
$H_1(Y_{B,0}(\p))$ & $  (0, 0), (2, 4), (3, 1), (4, 1), (8, 1), (9, 1), (11, 2), (19, 1), (25, 1), (47, 1) $ \\
 & $  (64, 1), (113, 1), (433, 1), (991, 1) $ \\

\hline  $N(\p) =  337 
$ & \\ \hline $ H_1(Y_0(\p))$  & $  (0, 1), (2, 3), (3, 1), (4, 3), (7, 1), (127, 1) $ \\
$H_1(Y_0(\pi \cdot \p))$ & $  (0, 2), (2, 8), (3, 1), (4, 2), (7, 1), (8, 3), (23, 1), (25, 1), (127, 2) $ \\
& $  (137, 1), (641, 1) $ \\
$H_1(Y_0(\pibar \cdot \p))$ & $  (0, 2), (2, 4), (3, 1), (4, 2), (7, 2), (8, 2), (11, 1), (16, 1), (25, 1), (127, 2) $ \\
 & $ (331, 1) $ \\
$H_1(Y_0(3 \cdot \p))$& $  (0, 5), (2, 23), (3, 3), (4, 7), (5, 1), (7, 2), (8, 6), (11, 2), (16, 2), (23, 2) $ \\
 & $ (25, 3), (49, 1), (127, 4), (137, 2), (256, 1), (331, 2), (641, 2) $ \\
$H_1(Y_{B,0}(\p))$ & $  (0, 1), (2, 10), (3, 3), (4, 7), (7, 2), (16, 2), (32, 1) $ \\

\hline  $N(\p) =  347 
$ & \\ \hline $ H_1(Y_0(\p))$  & $  (0, 0), (2, 3), (4, 3), (173, 1) $ \\
$H_1(Y_0(\pi \cdot \p))$ & $  (0, 0), (2, 6), (4, 6), (173, 1), (1129, 1), (1392367, 1) $ \\
$H_1(Y_0(\pibar \cdot \p))$ & $  (0, 0), (2, 8), (3, 1), (4, 7), (173, 1), (5197, 1), (5441, 1) $ \\
$H_1(Y_0(3 \cdot \p))$& $  (0, 1), (2, 13), (3, 5), (4, 20), (16, 1), (23, 1), (64, 1), (173, 1), (263, 1) $ \\
 & $ (1129, 2), (5197, 2), (5441, 2), (1392367, 2) $ \\
$H_1(Y_{B,0}(\p))$ & $  (0, 1), (2, 13), (3, 4), (4, 1), (8, 1), (23, 1), (173, 1), (256, 1), (263, 1) $ \\

\hline  $N(\p) =  353 
$ & \\ \hline $ H_1(Y_0(\p))$  & $  (0, 0), (2, 7), (4, 1), (11, 1), (16, 1) $ \\
$H_1(Y_0(\pi \cdot \p))$ & $  (0, 0), (2, 11), (3, 1), (4, 2), (8, 1), (11, 2), (17, 1), (128, 1), (24421, 1) $ \\
$H_1(Y_0(\pibar \cdot \p))$ & $  (0, 1), (2, 15), (4, 5), (11, 1), (16, 1), (32, 1), (41, 1), (173, 1) $ \\
$H_1(Y_0(3 \cdot \p))$& $  (0, 4), (2, 35), (3, 2), (4, 10), (11, 4), (16, 4), (17, 2), (23, 1), (32, 1) $ \\
& $ (41, 1), (59, 2), (67, 1), (173, 2), (256, 1), (257, 1), (512, 1), (24421, 2) $ \\
$H_1(Y_{B,0}(\p))$ & $  (0, 2), (2, 9), (3, 1), (4, 3), (8, 1), (11, 2), (16, 2), (23, 1), (59, 2), (67, 1) $ \\
 & $ (257, 1), (512, 1) $ \\

\hline  $N(\p) =  379 
$ & \\ \hline $ H_1(Y_0(\p))$  & $  (0, 0), (2, 6), (4, 4), (7, 1), (13, 1), (27, 1) $ \\
$H_1(Y_0(\pi \cdot \p))$ & $  (0, 0), (2, 12), (4, 8), (7, 2), (9, 2), (13, 2), (27, 1), (29, 1), (43, 1) $ \\
 & $ (20359, 1) $ \\
$H_1(Y_0(\pibar \cdot \p))$ & $  (0, 0), (2, 14), (4, 8), (7, 1), (8, 2), (13, 2), (27, 1), (113, 1), (2063, 1) $ \\
$H_1(Y_0(3 \cdot \p))$& $  (0, 0), (2, 25), (4, 18), (5, 1), (7, 2), (8, 5), (9, 2), (11, 1), (13, 4) $ \\
 & $ (27, 3), (29, 2), (31, 1), (43, 2), (49, 1), (59, 1), (67, 1), (113, 2) $ \\
  & $  (127, 1), (337, 1),(569, 1), (2063, 2), (20359, 2) $ \\
$H_1(Y_{B,0}(\p))$ & $  (0, 0), (2, 2), (3, 2), (4, 2), (5, 1), (7, 2), (11, 1), (16, 1), (27, 1), (31, 1) $ \\
 & $ (59, 1), (67, 1), (127, 1), (337, 1), (569, 1) $ \\
\hline
\end{tabular}

\begin{tabular}{|l|l|}
\hline  $N(\p) =  401 
$ & \\ \hline $ H_1(Y_0(\p))$  & $  (0, 0), (2, 4), (8, 1), (25, 1), (37, 1), (71, 1) $ \\
$H_1(Y_0(\pi \cdot \p))$ & $  (0, 0), (2, 3), (4, 1), (16, 1), (17, 1), (25, 1), (29, 1), (32, 1), (37, 2) $ \\
 & $ (71, 2), (199, 1), (4357, 1) $ \\
$H_1(Y_0(\pibar \cdot \p))$ & $  (0, 0), (2, 5), (4, 2), (5, 1), (7, 2), (8, 1), (16, 2), (23, 1), (25, 1), (37, 2) $ \\
 & $ (71, 2), (1879, 1) $ \\
$H_1(Y_0(3 \cdot \p))$& $  (0, 0), (2, 17), (3, 2), (4, 5), (5, 1), (7, 4), (8, 3), (11, 1), (16, 2), (17, 2) $ \\
 & $ (23, 2), (25, 2), (29, 2), (32, 1), (37, 4), (49, 1), (53, 1), (64, 2), (71, 4) $ \\
  & $ (103, 1), (199, 2), (419, 1), (1879, 2), (2269, 1), (4357, 2) $ \\
$H_1(Y_{B,0}(\p))$ & $  (0, 0), (2, 8), (3, 2), (4, 7), (5, 1), (7, 2), (11, 1), (25, 1), (32, 1), (53, 1) $ \\
 & $ (103, 1), (419, 1), (2269, 1) $ \\

\hline  $N(\p) =  409 
$ & \\ \hline $ H_1(Y_0(\p))$  & $  (0, 0), (2, 5), (3, 1), (4, 1), (5, 1), (17, 1), (49, 1), (101, 1) $ \\
$H_1(Y_0(\pi \cdot \p))$ & $  (0, 0), (2, 3), (3, 1), (4, 2), (5, 2), (8, 3), (17, 1), (29, 1), (49, 2) $ \\
 & $  (101, 2),(1693, 1) $ \\
$H_1(Y_0(\pibar \cdot \p))$ & $  (0, 6), (2, 9), (3, 1), (4, 4), (5, 1), (8, 1), (17, 1), (49, 2), (89, 1) $ \\
 & $  (101, 2) $ \\
$H_1(Y_0(3 \cdot \p))$& $  (0, 12), (2, 22), (3, 1), (4, 8), (5, 3), (8, 4), (16, 3), (17, 1), (23, 1) $ \\
 & $ (29, 3), (32, 1), (49, 4), (89, 2), (101, 4), (241, 1), (929, 1), (1693, 2) $ \\
  & $ (3779, 1) $ \\
$H_1(Y_{B,0}(\p))$ & $  (0, 0), (2, 5), (3, 1), (4, 7), (5, 1), (8, 1), (16, 3), (17, 1), (23, 1) $ \\
 & $  (25, 1),(29, 2), (37, 1), (241, 1), (929, 1), (3779, 1) $ \\

\hline  $N(\p) =  419 
$ & \\ \hline $ H_1(Y_0(\p))$  & $  (0, 0), (2, 3), (7, 1), (9, 1), (11, 1), (19, 1) $ \\
$H_1(Y_0(\pi \cdot \p))$ & $  (0, 0), (2, 6), (3, 1), (5, 1), (7, 1), (9, 1), (11, 2), (19, 1), (27, 1) $ \\
 & $ (49, 1), (59, 1), (79, 1), (1481, 1) $ \\
$H_1(Y_0(\pibar \cdot \p))$ & $  (0, 2), (2, 6), (3, 5), (4, 2), (7, 1), (9, 1), (11, 2), (19, 1), (243, 1) $ \\
 & $ (1427, 1) $ \\
$H_1(Y_0(3 \cdot \p))$& $  (0, 4), (2, 9), (3, 9), (4, 6), (5, 1), (7, 2), (9, 4), (11, 4), (13, 1), (16, 1) $ \\
 & $ (17, 1), (19, 1), (27, 1), (49, 1), (59, 2), (79, 2), (121, 1), (243, 1) $ \\
  & $ (729, 1), (1427, 2), (1481, 2), (10303, 1), (26251, 1) $ \\
$H_1(Y_{B,0}(\p))$ & $  (0, 0), (2, 3), (3, 2), (4, 1), (7, 1), (9, 1), (11, 2), (13, 1), (17, 1), (19, 1) $ \\
 & $ (64, 1), (10303, 1), (26251, 1) $ \\

\hline  $N(\p) =  433 
$ & \\ \hline $ H_1(Y_0(\p))$  & $  (0, 0), (2, 4), (7, 1), (8, 1), (19, 1), (27, 1), (2437, 1) $ \\
$H_1(Y_0(\pi \cdot \p))$ & $  (0, 0), (2, 3), (3, 2), (4, 2), (7, 2), (13, 2), (16, 1), (19, 2), (27, 1) $ \\
 & $ (307, 1), (2437, 2), (52561, 1) $ \\
$H_1(Y_0(\pibar \cdot \p))$ & $  (0, 1), (2, 6), (3, 2), (4, 2), (7, 1), (8, 1), (13, 1), (16, 1), (19, 2), (27, 1) $ \\
 & $ (157, 1), (512, 1), (2437, 2) $ \\
$H_1(Y_0(3 \cdot \p))$& $  (0, 3), (2, 15), (3, 6), (4, 7), (7, 2), (8, 6), (9, 2), (13, 7), (16, 6), (19, 4) $ \\
 & $ (27, 1), (157, 2), (307, 2), (512, 2), (641, 1), (1979, 1), (2437, 4) $ \\
  & $ (52561, 2) $ \\
$H_1(Y_{B,0}(\p))$ & $  (0, 1), (2, 3), (3, 2), (4, 5), (8, 10), (13, 1), (16, 2), (27, 1), (641, 1) $ \\
 & $ (1979, 1) $ \\
\hline
\end{tabular}

\begin{tabular}{|l|l|}

\hline  $N(\p) =  443 
$ & \\ \hline $ H_1(Y_0(\p))$  & $  (0, 0), (2, 3), (13, 1), (49, 1), (289, 1), (353, 1) $ \\
$H_1(Y_0(\pi \cdot \p))$ & $  (0, 1), (2, 9), (7, 1), (13, 1), (17, 1), (49, 1), (289, 1), (353, 2) $ \\
 & $ (33391, 1), (200159, 1) $ \\
$H_1(Y_0(\pibar \cdot \p))$ & $  (0, 0), (2, 6), (5, 1), (13, 1), (17, 1), (49, 2), (289, 1), (353, 2) $ \\
 & $  (39877, 1), (57920143, 1) $ \\
$H_1(Y_0(3 \cdot \p))$& $  (0, 2), (2, 9), (3, 2), (4, 6), (5, 2), (7, 3), (8, 1), (11, 1), (13, 1), (16, 2) $ \\
 & $ (17, 3), (29, 1), (49, 2), (61, 1), (89, 1), (137, 1), (289, 1), (313, 1) $ \\
  & $ (353, 4), (4913, 1), (33391, 2), (39877, 2), (200159, 2), (57920143, 2) $ \\
$H_1(Y_{B,0}(\p))$ & $  (0, 0), (2, 5), (3, 2), (4, 1), (7, 1), (8, 4), (11, 1), (13, 1), (17, 1), (29, 1) $ \\
 & $ (32, 1), (61, 1), (89, 1), (137, 1), (313, 1), (4913, 1) $ \\

\hline  $N(\p) =  449 
$ & \\ \hline $ H_1(Y_0(\p))$  & $  (0, 0), (2, 4), (7, 1), (32, 1), (1425821, 1) $ \\
$H_1(Y_0(\pi \cdot \p))$ & $  (0, 0), (2, 3), (4, 2), (7, 1), (37, 1), (64, 1), (4639, 1), (9151, 1) $ \\
 & $ (1425821, 2) $ \\
$H_1(Y_0(\pibar \cdot \p))$ & $  (0, 0), (2, 3), (3, 3), (4, 1), (7, 1), (29, 1), (64, 2), (167, 1), (227, 1) $ \\
 & $ (1425821, 2) $ \\
$H_1(Y_0(3 \cdot \p))$& $  (0, 0), (2, 9), (3, 7), (4, 1), (5, 1), (7, 1), (8, 1), (17, 1), (25, 1) $ \\
& $ (29, 1), (31, 1), (37, 2), (43, 1), (64, 1), (128, 1), (167, 2), (227, 2) $  \\
 & $ (256, 1) $ \\
 & $ (289, 1), (841, 1), (937, 1), (4639, 2), (4657, 1), (9151, 2) $ \\
  & $ (1425821, 4) $ \\
$H_1(Y_{B,0}(\p))$ & $  (0, 0), (2, 4), (3, 1), (4, 3), (5, 1), (7, 1), (17, 1), (25, 1), (29, 1) $ \\
 & $ (31, 1), (43, 1), (128, 1), (289, 1), (937, 1), (4657, 1) $ \\

\hline  $N(\p) =  457 
$ & \\ \hline $ H_1(Y_0(\p))$  & $  (0, 0), (2, 7), (3, 2), (4, 2), (19, 1), (343, 1) $ \\
$H_1(Y_0(\pi \cdot \p))$ & $  (0, 1), (2, 14), (3, 6), (4, 4), (5, 1), (8, 1), (11, 1), (19, 1), (49, 1) $ \\
 & $ (83, 1), (89, 1), (343, 1) $ \\
$H_1(Y_0(\pibar \cdot \p))$ & $  (0, 0), (2, 11), (3, 3), (4, 4), (8, 1), (9, 1), (19, 1), (25, 1), (41, 1) $ \\
 & $  (61, 1), (343, 2), (397, 1), (1637, 1) $ \\
$H_1(Y_0(3 \cdot \p))$& $  (0, 2), (2, 26), (3, 9), (4, 14), (5, 2), (8, 4), (9, 4), (11, 2), (13, 1) $ \\
 & $  (16, 3), (19, 1), (25, 1), (32, 1), (41, 2), (49, 2), (53, 1), (61, 2) $ \\
  & $ (83, 2), (89, 2), (107, 1), (125, 1), (343, 2), (397, 2), (983, 1) $ \\
  & $ (1637, 2), (6241, 1) $ \\
$H_1(Y_{B,0}(\p))$ & $  (0, 0), (2, 8), (3, 2), (4, 5), (5, 1), (8, 5), (9, 1), (13, 1), (16, 2) $ \\
 & $  (19, 1), (53, 1), (107, 1), (983, 1), (6241, 1) $ \\
\hline
\end{tabular}

\begin{tabular}{|l|l|}

\hline  $N(\p) =  467 
$ & \\ \hline $ H_1(Y_0(\p))$  & $  (0, 0), (2, 3), (3, 3), (17, 1), (49, 1), (233, 1) $ \\
$H_1(Y_0(\pi \cdot \p))$ & $  (0, 0), (2, 6), (3, 6), (5, 1), (17, 2), (37, 1), (49, 2), (233, 1) $ \\
 & $ (503483, 1) $ \\
$H_1(Y_0(\pibar \cdot \p))$ & $  (0, 0), (2, 6), (3, 6), (17, 2), (23, 2), (31, 1), (37, 1), (49, 3), (233, 1) $ \\
 & $ (8861, 1) $ \\
$H_1(Y_0(3 \cdot \p))$& $  (0, 0), (2, 9), (3, 15), (4, 2), (5, 4), (8, 1), (17, 4), (23, 4), (25, 2) $ \\
 & $ (31, 3), (37, 4), (41, 1), (49, 6), (61, 1), (73, 2), (83, 1), (233, 1) $ \\
  & $  (8861, 2), (321467, 1), (503483, 2) $ \\
$H_1(Y_{B,0}(\p))$ & $  (0, 0), (2, 3), (3, 3), (4, 1), (5, 2), (25, 2), (31, 1), (32, 1), (41, 1) $ \\
& $ (61, 1), (73, 2), (83, 1), (233, 1), (321467, 1) $ \\

\hline  $N(\p) =  491 
$ & \\ \hline $ H_1(Y_0(\p))$  & $  (0, 0), (2, 3), (5, 1), (32, 1), (49, 1),  (421, 1) $ \\
$H_1(Y_0(\pi \cdot \p))$ & $  (0, 0), (2, 6), (5, 1), (17, 1), (19, 1), (49, 1), (128, 2), (199, 1) $ \\
 & $  (421, 2), (129763, 1) $ \\
$H_1(Y_0(\pibar \cdot \p))$ & $  (0, 0), (2, 6), (3, 2), (5, 1), (7, 1), (11, 1), (25, 1), (32, 2), (49, 1) $ \\
 & $  (73, 1), (421, 2), (787, 1), (5779, 1) $ \\
$H_1(Y_0(3 \cdot \p))$& $  (0, 1), (2, 9), (3, 3), (4, 9), (5, 1), (7, 3), (8, 1), (9, 1), (11, 3), (16, 1) $ \\
 & $ (17, 2), (19, 3), (25, 2), (32, 1), (49, 1), (73, 2), (128, 4), (199, 2) $ \\
 & $  (227, 1), (421, 4), (547, 1), (787, 2), (809, 1), (5779, 2), (129763, 2) $ \\
$H_1(Y_{B,0}(\p))$ & $  (0, 1), (2, 2), (3, 1), (4, 9), (5, 2), (7, 1), (8, 1), (11, 1), (16, 1), (19, 1) $ \\
 & $ (32, 1), (49, 1), (227, 1), (547, 1), (809, 1) $ \\

\hline  $N(\p) =  499 
$ & \\ \hline $ H_1(Y_0(\p))$  & $  (0, 0), (2, 3), (3, 1), (23, 1), (43, 1), (83, 1), (1291, 1) $ \\
$H_1(Y_0(\pi \cdot \p))$ & $  (0, 0), (2, 6), (3, 7), (19, 1), (23, 2), (43, 2), (81, 1), (83, 1), (619, 1) $ \\
& $ (1291, 2), (1423, 1) $ \\
$H_1(Y_0(\pibar \cdot \p))$ & $  (0, 0), (2, 6), (3, 2), (5, 1), (23, 2), (37, 1), (43, 2), (83, 1), (1291, 2) $ \\
& $ (14243, 1), (129529, 1) $ \\
$H_1(Y_0(3 \cdot \p))$& $  (0, 0), (2, 9), (3, 16), (4, 3), (5, 3), (9, 4), (16, 1), (19, 3), (23, 4) $ \\
 & $ (25, 1), (27, 1), (32, 1), (37, 2), (43, 4), (81, 1), (83, 1), (243, 1) $ \\
& $  (619, 2),  (797, 1), (1291, 4), (1423, 2), (4969, 1), (14243, 2) $ \\
 & $ (129529, 2) $ \\
$H_1(Y_{B,0}(\p))$ & $  (0, 0), (2, 2), (3, 10), (4, 3), (5, 3), (9, 2), (19, 1), (32, 2), (83, 1) $ \\
 & $  (797, 1), (4969, 1) $ \\
\hline
\end{tabular}

\begin{tabular}{|l|l|}

\hline  $N(\p) =  521 
$ & \\ \hline $ H_1(Y_0(\p))$  & $  (0, 0), (2, 5), (3, 1), (4, 1), (5, 1), (7, 1), (13, 1), (2213, 1) $ \\
$H_1(Y_0(\pi \cdot \p))$ & $  (0, 2), (2, 3), (3, 3), (4, 2), (5, 1), (7, 3), (8, 1), (13, 1), (19, 1), (31, 1) $ \\
 & $ (139, 1), (169, 1), (2213, 2) $ \\
$H_1(Y_0(\pibar \cdot \p))$ & $  (0, 0), (2, 5), (3, 1), (4, 2), (5, 2), (7, 2), (8, 1), (9, 1), (13, 1), (17, 1) $ \\
 & $ (19, 1), (53, 1), (2213, 2), (6133, 1) $ \\
$H_1(Y_0(3 \cdot \p))$& $  (0, 4), (2, 9), (3, 8), (4, 3), (5, 3), (7, 7), (8, 1), (9, 1), (13, 2), (17, 2) $ \\
 & $ (19, 3), (31, 2), (32, 1), (53, 2), (97, 1), (139, 2), (169, 2), (223, 1) $ \\
  & $ (361, 1), (401, 1), (2213, 4), (6133, 2), (1750297, 1), (2613151, 1) $ \\
$H_1(Y_{B,0}(\p))$ & $  (0, 0), (2, 4), (3, 3), (4, 3), (5, 1), (7, 1), (11, 1), (13, 2), (16, 1) $ \\
 & $  (19, 1) , (97, 1), (223, 1), (401, 1), (1750297, 1), (2613151, 1) $ \\

\hline  $N(\p) =  523 
$ & \\ \hline $ H_1(Y_0(\p))$  & $  (0, 0), (2, 3), (3, 1), (8, 2), (9, 1), (29, 1), (31, 1), (289, 1) $ \\
$H_1(Y_0(\pi \cdot \p))$ & $  (0, 0), (2, 6), (3, 2), (8, 4), (9, 1), (17, 1), (29, 1), (31, 2), (53, 1) $ \\
& $ (289, 2), (6343, 1), (413143, 1) $ \\
$H_1(Y_0(\pibar \cdot \p))$ & $  (0, 0), (2, 6), (3, 5), (5, 1), (8, 4), (9, 1), (11, 2), (19, 1), (23, 1) $ \\
 & $ (29, 1), (31, 2), (43, 1), (289, 2), (2585447, 1) $ \\
$H_1(Y_0(3 \cdot \p))$& $  (0, 0), (2, 9), (3, 12), (4, 14), (5, 5), (8, 9), (9, 3), (11, 4), (17, 2) $ \\
 & $ (19, 2), (23, 2), (29, 1), (31, 4), (37, 1), (43, 3), (53, 2), (109, 1) $ \\
  & $ (229, 1), (289, 4), (2657, 1), (6343, 2), (413143, 2), (2585447, 2) $ \\
$H_1(Y_{B,0}(\p))$ & $  (0, 0), (2, 2), (3, 4), (4, 14), (5, 3), (9, 2), (16, 1), (29, 1), (37, 1) $ \\
 & $  (43, 1), (109, 1), (229, 1), (2657, 1) $ \\
\hline  $N(\p) =  547 
$ & \\ \hline $ H_1(Y_0(\p))$  & $  (0, 0), (2, 3), (3, 1), (7, 1), (13, 1), (1376237, 1) $ \\
$H_1(Y_0(\pi \cdot \p))$ & $  (0, 0), (2, 6), (3, 1), (7, 3), (13, 1), (24943, 1), (670619, 1) $ \\
 & $ (1376237, 2) $ \\
$H_1(Y_0(\pibar \cdot \p))$ & $  (0, 0), (2, 6), (3, 4), (5, 1), (7, 2), (11, 1), (13, 1), (23, 1), (27, 1) $ \\
 & $ (983, 1), (1376237, 2) $ \\
$H_1(Y_0(3 \cdot \p))$& $  (0, 0), (2, 9), (3, 9), (4, 2), (5, 4), (7, 7), (8, 1), (11, 3), (13, 2), (23, 2) $ \\
 & $ (25, 1), (27, 1), (83, 1), (181, 1), (243, 1), (409, 1), (631, 1), (821, 1) $ \\
  & $ (983, 2), (15581, 1), (24943, 2), (670619, 2), (1376237, 4) $ \\
$H_1(Y_{B,0}(\p))$ & $  (0, 0), (2, 2), (3, 3), (4, 1), (5, 4), (7, 1), (8, 2), (9, 1), (11, 1), (13, 2) $ \\
 & $ (83, 1), (181, 1), (409, 1), (631, 1), (821, 1), (15581, 1) $ \\
\hline
\end{tabular}

\begin{tabular}{|l|l|}

\hline  $N(\p) =  563 
$ & \\ \hline $ H_1(Y_0(\p))$  & $  (0, 0), (2, 3), (5, 1), (19, 1), (281, 1), (5507, 1) $ \\
$H_1(Y_0(\pi \cdot \p))$ & $  (0, 0), (2, 6), (5, 2), (13, 1), (19, 2), (89, 1), (281, 1), (383, 1), (647, 1) $ \\
 & $ (1663, 1), (5507, 2) $ \\
$H_1(Y_0(\pibar \cdot \p))$ & $  (0, 0), (2, 6), (5, 3), (19, 2), (121, 1), (281, 1), (5393, 1), (5507, 2) $ \\ 
 & $ (323379209, 1) $ \\
$H_1(Y_0(3 \cdot \p))$& $  (0, 0), (2, 9), (3, 9), (4, 2), (5, 6), (7, 1), (8, 1), (9, 1), (11, 2), (13, 2) $ \\
 & $ (17, 1), (19, 4), (41, 1), (89, 2), (101, 1), (121, 3), (233, 1), (281, 2) $ \\
  & $ (383, 2), (647, 2), (1601, 1), (1663, 2), (5393, 2), (5507, 4) $ \\
   & $ (323379209, 2) $ \\
$H_1(Y_{B,0}(\p))$ & $  (0, 0), (2, 3), (3, 9), (4, 1), (7, 1), (9, 1), (11, 2), (17, 1), (32, 1), (41, 1) $ \\
 & $ (101, 1), (121, 1), (233, 1), (281, 2), (1601, 1) $ \\

\hline  $N(\p) =  569 
$ & \\ \hline $ H_1(Y_0(\p))$  & $  (0, 0), (2, 5), (4, 1), (9, 1), (71, 1), (25169, 1) $ \\
$H_1(Y_0(\pi \cdot \p))$ & $  (0, 0), (2, 3), (3, 9), (4, 2), (5, 1), (7, 1), (8, 3), (9, 1), (27, 1), (71, 1) $ \\
 & $ (1669, 1), (25169, 2) $ \\
$H_1(Y_0(\pibar \cdot \p))$ & $  (0, 0), (2, 5), (3, 2), (4, 2), (7, 1), (8, 1), (9, 2), (11, 1), (71, 1), (191, 1) $ \\
 & $ (467, 1), (5477, 1), (25169, 2) $ \\
$H_1(Y_0(3 \cdot \p))$& $  (0, 0), (2, 9), (3, 21), (4, 3), (5, 2), (7, 4), (8, 4), (9, 4), (11, 2), (13, 1) $ \\
 & $ (16, 1), (27, 3), (31, 1), (32, 1), (67, 1), (71, 1), (107, 1), (173, 1) $ \\
  & $ (191, 2), (281, 1), (467, 2), (1489, 1), (1669, 2), (5477, 2), (25169, 4) $ \\
  & $ (554209, 1) $ \\
$H_1(Y_{B,0}(\p))$ & $  (0, 0), (2, 4), (3, 4), (4, 3), (9, 1), (13, 1), (16, 1), (31, 1), (67, 1), (71, 1) $ \\
 & $ (107, 1), (173, 1), (281, 1), (1489, 1), (554209, 1) $ \\

\hline  $N(\p) =  571 
$ & \\ \hline $ H_1(Y_0(\p))$  & $  (0, 0), (2, 3), (3, 1), (5, 2), (19, 2), (15889, 1) $ \\
$H_1(Y_0(\pi \cdot \p))$ & $  (0, 0), (2, 6), (3, 1), (5, 3), (17, 1), (19, 3), (25, 1), (29, 1), (15889, 2) $ \\ 
& $ (19709, 1), (865043, 1) $ \\
$H_1(Y_0(\pibar \cdot \p))$ & $  (0, 0), (2, 6), (3, 1), (5, 3), (7, 1), (13, 1), (19, 3), (139, 1), (1901, 1) $ \\
 & $ (15889, 2), (922069, 1) $ \\
$H_1(Y_0(3 \cdot \p))$& $  (0, 0), (2, 9), (3, 4), (4, 3), (5, 6), (7, 3), (13, 3) $ \\
 & $ (16, 1), (17, 2), (19, 5) $ \\ 
  & $ (25, 2), (29, 2), (31, 1), (64, 1), (139, 2), (167, 1), (307, 1), (313, 1) $ \\
   & $ (353, 1), (1373, 1), (1607, 1), (1901, 2), (1933, 1), (15889, 4) $ \\
    & $  (19709, 2), (865043, 2), (922069, 2) $ \\
$H_1(Y_{B,0}(\p))$ & $  (0, 0), (2, 2), (3, 4), (4, 3), (5, 2), (7, 1), (13, 1), (16, 1), (19, 1), (31, 1) $ \\
 &$  (128, 1), (167, 1), (307, 1), (313, 1), (353, 1), (1373, 1), (1607, 1) $ \\
  & $ (1933, 1) $ \\
\hline
\end{tabular}

\begin{tabular}{|l|l|}

\hline  $N(\p) =  577 
$ & \\ \hline $ H_1(Y_0(\p))$  & $  (0, 0), (2, 4), (4, 2), (9, 1), (32, 1), (283, 1), (821, 1) $ \\
$H_1(Y_0(\pi \cdot \p))$ & $  (0, 0), (2, 6), (3, 1), (4, 2), (7, 1), (8, 2), (9, 1), (23, 1), (37, 1) $ \\
 & $ (64, 1), (227, 1), (283, 2), (821, 2), (37199, 1) $ \\
$H_1(Y_0(\pibar \cdot \p))$ & $  (0, 0), (2, 8), (3, 3), (4, 2), (8, 2), (9, 1), (11, 1), (47, 1), (64, 1) $ \\
 & $ (239, 1), (283, 2), (821, 2), (35837, 1) $ \\
$H_1(Y_0(3 \cdot \p))$& $  (0, 0), (2, 19), (3, 9), (4, 6), (5, 3), (7, 2), (8, 8), (9, 1), (11, 3), (13, 3) $ \\
 & $ (16, 1), (17, 1), (23, 3), (32, 1), (37, 2), (47, 2), (121, 1), (227, 2) $ \\
 & $ (239, 2), (256, 1), (283, 4), (293, 1), (337, 1), (607, 1), (821, 4) $ \\
  & $ (853, 1), (35837, 2), (37199, 2) $ \\
$H_1(Y_{B,0}(\p))$ & $  (0, 0), (2, 7), (3, 1), (4, 4), (5, 3), (8, 4), (9, 1), (11, 3), (13, 3), (16, 1) $ \\
 & $ (17, 1), (23, 1), (128, 1), (293, 1), (337, 1), (607, 1), (853, 1) $ \\

\hline  $N(\p) =  587 
$ & \\ \hline $ H_1(Y_0(\p))$  & $  (0, 0), (2, 3), (3, 1), (5, 1), (8, 1), (149, 1), (293, 1) $ \\
$H_1(Y_0(\pi \cdot \p))$ & $  (0, 0), (2, 8), (3, 1), (5, 3), (7, 1), (9, 1), (16, 2), (25, 1), (37, 1), (59, 1) $ \\
 & $ (149, 2), (293, 1), (1187, 1), (431891, 1) $ \\
$H_1(Y_0(\pibar \cdot \p))$ & $  (0, 0), (2, 10), (3, 1), (4, 4), (5, 2), (8, 2), (13, 1), (27, 1), (49, 1) $ \\
 & $ (149, 2), (293, 1), (409, 1), (2711, 1) $ \\
$H_1(Y_0(3 \cdot \p))$& $  (0, 1), (2, 21), (3, 2), (4, 11), (5, 6), (7, 2), (8, 1), (9, 1), (11, 1), (13, 2) $ \\
 & $ (16, 4), (17, 2), (25, 2), (27, 1), (37, 2), (41, 1), (49, 2), (53, 1), (59, 2) $ \\ 
 & $ (73, 1), (81, 1), (149, 4), (293, 1), (409, 2), (887, 1), (1187, 2) $ \\
 & $ (2711, 2), (6977, 1), (188143, 1), (431891, 2) $ \\
$H_1(Y_{B,0}(\p))$ & $  (0, 1), (2, 2), (3, 1), (4, 3), (5, 1), (8, 1), (11, 1), (13, 1), (17, 2), (41, 1) $ \\
 & $ (53, 1), (73, 1), (293, 1), (887, 1), (6977, 1), (188143, 1) $ \\

\hline  $N(\p) =  593 
$ & \\ \hline $ H_1(Y_0(\p))$  & $  (0, 0), (2, 9), (4, 5), (8, 1), (17, 1), (37, 1) $ \\
$H_1(Y_0(\pi \cdot \p))$ & $  (0, 0), (2, 7), (3, 1), (4, 12), (8, 8), (16, 1), (17, 2), (32, 1), (37, 1) $ \\
 & $ (179, 1), (443, 1) $ \\
$H_1(Y_0(\pibar \cdot \p))$ & $  (0, 0), (2, 18), (3, 2), (4, 10), (5, 1), (8, 2), (16, 1), (17, 2), (37, 1) $ \\
& $ (53, 1), (439, 1) $ \\
$H_1(Y_0(3 \cdot \p))$& $  (0, 1), (2, 20), (3, 3), (4, 25), (5, 2), (8, 23), (9, 1), (16, 1), (17, 4) $ \\
 & $ (23, 1), (31, 1), (37, 1), (53, 2), (64, 3), (79, 1), (97, 1), (128, 1) $ \\
  & $ (179, 2), (229, 1), (379, 1), (439, 2), (443, 2), (1009, 1), (1087, 1) $ \\
   & $ (18523, 1) $ \\
$H_1(Y_{B,0}(\p))$ & $  (0, 1), (2, 11), (4, 5), (8, 1), (16, 1), (23, 1), (31, 1), (32, 1), (37, 1) $ \\
& $ (79, 1), (97, 1), (128, 1), (229, 1), (379, 1), (1009, 1), (1087, 1) $ \\
 & $ (18523, 1) $ \\
\hline
\end{tabular}

\begin{tabular}{|l|l|}

\hline  $N(\p) =  601 
$ & \\ \hline $ H_1(Y_0(\p))$  & $  (0, 0), (2, 4), (3, 1), (4, 1), (25, 1), (11964811, 1) $ \\
$H_1(Y_0(\pi \cdot \p))$ & $  (0, 0), (2, 4), (3, 1), (4, 1), (11, 1), (17, 1), (25, 1), (32, 1), (41, 1) $ \\
 & $ (47, 1), (64, 1), (191, 1), (282349, 1), (11964811, 2) $ \\
$H_1(Y_0(\pibar \cdot \p))$ & $  (0, 0), (2, 6), (3, 1), (4, 1), (7, 2), (8, 1), (25, 1), (32, 1), (47, 1), (97, 1) $ \\
 & $ (317, 1), (78203, 1), (11964811, 2) $ \\
$H_1(Y_0(3 \cdot \p))$& $  (0, 1), (2, 15), (3, 5), (4, 8), (7, 5), (8, 3), (11, 2), (13, 1), (16, 1), (17, 4) $ \\
& $ (25, 1), (32, 1), (41, 2), (47, 4), (64, 1), (97, 2), (163, 1), (191, 2) $ \\
 & $  (223, 1), (256, 2), (277, 1), (317, 2), (523, 1), (919, 1), (1249, 1) $ \\
 & $ (78203, 2), (282349, 2), (11964811, 4) $ \\
$H_1(Y_{B,0}(\p))$ & $  (0, 1), (2, 7), (3, 5), (4, 7), (7, 1), (8, 5), (13, 1), (16, 1), (17, 2), (25, 1) $ \\
& $ (163, 1), (223, 1), (277, 1), (523, 1), (919, 1), (1249, 1) $ \\

\hline  $N(\p) =  617 
$ & \\ \hline $ H_1(Y_0(\p))$  & $  (0, 0), (2, 4), (4, 1), (7, 2), (11, 3), (59, 1), (103, 1) $ \\
$H_1(Y_0(\pi \cdot \p))$ & $  (0, 1), (2, 4), (4, 1), (7, 3), (9, 1), (11, 5), (16, 1), (19, 1), (59, 4) $ \\
 & $  (103, 2), (565583, 1) $ \\
$H_1(Y_0(\pibar \cdot \p))$ & $  (0, 0), (2, 4), (3, 1), (4, 3), (7, 3), (8, 1), (9, 1), (11, 5), (59, 2), (79, 1) $ \\
 & $ (103, 2), (107, 1), (2017, 1), (11161, 1) $ \\
$H_1(Y_0(3 \cdot \p))$& $  (0, 4), (2, 15), (3, 2), (4, 4), (5, 1), (7, 6), (8, 1), (9, 2), (11, 8), (13, 1) $ \\
 & $ (19, 2), (23, 1), (49, 1), (59, 8), (64, 1), (79, 2), (81, 3), (83, 1) $ \\
  & $ (103, 4), (107, 2), (121, 1), (256, 1), (439, 1), (2017, 2), (11161, 2) $ \\
   & $ (28559, 1), (565583, 2) $ \\
$H_1(Y_{B,0}(\p))$ & $  (0, 2), (2, 8), (4, 6), (5, 1), (7, 2), (13, 1), (23, 1), (27, 2), (49, 1) $ \\
 & $ (64, 1), (81, 1), (83, 1), (121, 1), (439, 1), (28559, 1) $ \\
 \hline
  \end{tabular}
  \end{center}
  
  \newpage

The third table represents computations of the cohomology of
$\PGL_2(\OL_F)$ (with $F = \Q(\sqrt{-2})$ ) of the local system
$\M_k$ corresponding to $\Sym^{k-2} \otimes \overline{\Sym^{k-2}}$,
The entry in column $(p,k)$ denotes the quantity:
$$ \dim H_1(\GL_2(\OL_F),\M_k \otimes \F_p) - \dim H_1(\GL_2(\OL_F),\M_k \otimes \C)^{BC}.$$
 Complex conjugation acts on both spaces,
which explains why every entry in the table is even. (If the representation was
invariant under complex conjugation, then the indicated number would still
be even by Lemma~\ref{lemma:doubling}.)
The computations for $F = \Q(\sqrt{-2})$ and $\Q(\sqrt{-11})$ were
carried out using the methods indicated by Theorem~8.10 of~\cite{CM},
with the geometry of the Bianchi manifolds for these fields
given by~\cite{Cremona}.
We make some specific remarks on the tables.
\begin{enumerate}
\item
If $k > p+1$, the quantity is greyed out, because any $\rhobar$ is twist equivalent to a representation
with $k(\rhobar) \le p+1$.
\item We only list what happens for odd $k$. If $\rhobar$ has $k(\rhobar) = k$, then
$\det(\rhobar) = \omega^{k-1}$ which is only even if $k$ is odd (at least when $p$ is odd. When $p = 2$,
one can appeal to~\cite{TateNonExist}.)
\item
The green squares correspond to weights for which the non-existence of $\rhobar$ is known, whereas magenta squares correspond
to weights in which  a similar conclusion is known modulo GRH. (See~\cite{Moon}).
\item The yellow squares denote pairs $(p,k)$ for which the computation above yields a non-zero quantity
for $F = \Q(\sqrt{-2})$, but such that the same computation for the
field $F = \Q(\sqrt{-11})$ yields the result zero. In particular, for yellow squares,  there are no even representations
$\rhobar$ (assuming Conjecture~\ref{conj:serre}.)
\item The blue squares are the squares for which our quantity is non-zero for
\emph{both} $F = \Q(\sqrt{-2})$ and $F = \Q(\sqrt{-11})$. These are the only
squares which could plausibly correspond to even
Galois representations. The only blue square corresponds to $(p,k) = (163,55)$,
which is accounted for by an $A_4$-representation.
\item  If $k \le 21$, the results are in agreement with computations
of \Sengun  \  for $F = \Q(\sqrt{-2})$, which were computed using another method (see~\cite{Sengun}).
\item If $k \le 32$, the integral cohomology for $F = \Q(\sqrt{-1})$ was computed by \Sengun~\cite{Sengun}.
(There is a conflict of notation here --- our $\M_k$ corresponds to \Sengun's $M_{k-2,k-2}$.)
If follows from those computations that the only possible even $\rhobar$ of weight
$k(\rhobar) = k \le 32$ must occur with $(p,k) = (89,31)$. Yet the corresponding
computation with $F = \Q(\sqrt{-2})$ yields zero.
\end{enumerate}

\newpage 
  {\small 
\begin{table}[h]
\definecolor{agray}{gray}{0.8}
\begin{center}
\begin{tabular}{rr|*{13}{c}}
\multicolumn{2}{c}{} & \multicolumn{13}{c}{Serre weight $k$} \\
 &    &  $3$ & $5$ & $7$     &$9$ &  $11$ & $13$ & $15$ & $17$ & $19$ & $21$ & $23$ & $25$ & $27$ \\
\cline{2-15}\\[-12.9pt]
& $3$    &   \cellcolor{green} $0$
   &\cellcolor{agray}   & \cellcolor{agray}   &\cellcolor{agray}   & \cellcolor{agray}  &\cellcolor{agray}   & \cellcolor{agray}    & \cellcolor{agray}   & \cellcolor{agray}    & \cellcolor{agray}   &  \cellcolor{agray} &\cellcolor{agray} &
 \cellcolor{agray} \\
& $5$   & \cellcolor{green} $0$   &\cellcolor{green}  $0$  & \cellcolor{agray}      &\cellcolor{agray}   & \cellcolor{agray}   &\cellcolor{agray}   & \cellcolor{agray}    & \cellcolor{agray}   & \cellcolor{agray}    & \cellcolor{agray}   &  \cellcolor{agray} &\cellcolor{agray} &
 \cellcolor{agray} \\
& $7$   & \cellcolor{green} $0$  &\cellcolor{green}  $0$   &\cellcolor{green}  $0$     &\cellcolor{agray}   & \cellcolor{agray}   &\cellcolor{agray}   & \cellcolor{agray}    & \cellcolor{agray}   & \cellcolor{agray}    & \cellcolor{agray}   &  \cellcolor{agray} &\cellcolor{agray} &
 \cellcolor{agray} \\
& $11$   &\cellcolor{green}  $0$  &\cellcolor{green}  $0$   &\cellcolor{green}   $0$    & \cellcolor[cmyk]{0,0.5,0,0}$0$
  & \cellcolor[cmyk]{0,0.5,0,0}$0$
  &\cellcolor{agray}   & \cellcolor{agray}    & \cellcolor{agray}   & \cellcolor{agray}    & \cellcolor{agray}   &  \cellcolor{agray} &\cellcolor{agray} &
 \cellcolor{agray} \\
& $13$     & \cellcolor{green} $0$  &\cellcolor{green}  $0$   & \cellcolor{green}  $0$    & \cellcolor[cmyk]{0,0.5,0,0}$0$
  &  \cellcolor[cmyk]{0,0.5,0,0}$0$
   & $0$  & \cellcolor{agray}    & \cellcolor{agray}   & \cellcolor{agray}    & \cellcolor{agray}   &  \cellcolor{agray} &\cellcolor{agray} &
 \cellcolor{agray} \\
Prime $p$     &  $17$  &\cellcolor{green}  $0$   &\cellcolor{green}  $0$   & \cellcolor{green}  $0$  & \cellcolor[cmyk]{0,0.5,0,0}$0$
   & \cellcolor[cmyk]{0,0.5,0,0}$0$
   & \cellcolor[cmyk]{0,0.5,0,0}$0$
  &$0$    &$0$   & \cellcolor{agray}    & \cellcolor{agray}   &  \cellcolor{agray} &\cellcolor{agray} &
 \cellcolor{agray} \\
   & $19$    & \cellcolor{green} $0$  & \cellcolor{green}  $0$   & \cellcolor{green}  $0$    & \cellcolor[cmyk]{0,0.5,0,0}$0$
 & \cellcolor[cmyk]{0,0.5,0,0}$0$
  & \cellcolor[cmyk]{0,0.5,0,0}$0$
   &$0$    &$0$   & $0$    & \cellcolor{agray}  & \cellcolor{agray}  &\cellcolor{agray} &
 \cellcolor{agray} \\
  & $23$     & $0$  & $0$   & $0$    & $0$  & $0$   & $0$   &$0$    &$0$   &$0$    & $0$   & $0$ & \cellcolor{agray} & \cellcolor{agray} \\
    & $29$     & $0$  & $0$   & $0$    & $0$  & $0$   & $0$   &$0$    &$0$   &$0$    & $0$   & $0$ & $0$ & $0$ \\
    & $\vdots$     & $\vdots$  & $\vdots$   & $\vdots$    & $\vdots$  & $\vdots$   & $\vdots$   &$\vdots$    &$\vdots$   &$\vdots$    & $\vdots$   & $\vdots$ & $\vdots$ & $\vdots$ \\
    & $67$    & $0$  & $0$   & $0$     & $0$  & $0$   & $0$   &$0$    &$0$   &$0$    & $0$   & $0$ & $0$ & $0$ \\
    & $71$    & $0$  & $0$   & $0$     & $0$  & $0$   & $0$   &$0$    &$0$   &$0$    & $0$   & $0$ & $0$ & \cellcolor{yellow}  $2$ \\
       & $\vdots$     & $\vdots$  & $\vdots$   & $\vdots$    & $\vdots$  & $\vdots$   & $\vdots$   &$\vdots$    &$\vdots$   &$\vdots$    & $\vdots$   & $\vdots$ & $\vdots$ & $\vdots$ \\
           & $103$    & $0$  & $0$   & $0$     & $0$  & $0$   & $0$   &$0$   & $0$  &$\cellcolor{yellow} 2$   &$0$     & $0$ & $0$ & $0$ \\
           & $\vdots$     & $\vdots$  & $\vdots$   & $\vdots$    & $\vdots$  & $\vdots$   & $\vdots$   &$\vdots$    &$\vdots$   &$\vdots$    & $\vdots$   & $\vdots$ & $\vdots$ & $\vdots$ \\
           \end{tabular}
\end{center}
\end{table}}

\newpage

\begin{table}[ht]
\definecolor{agray}{gray}{0.8}
\begin{center}
\begin{tabular}{rr|*{13}{c}}
\multicolumn{2}{c}{} & \multicolumn{13}{c}{Serre weight $k$} \\
 &          &$29$ &  $31$ & $33$ & $35$ & $37$ & $39$ & $41$ & $43$ & $45$ & $47$ & $49$ & $51$ & $53$ \\
\cline{2-15}\\[-12.9pt]
& $29$ & $0$  &\cellcolor{agray}  &\cellcolor{agray}   & \cellcolor{agray}    & \cellcolor{agray}   & \cellcolor{agray}    & \cellcolor{agray}   &  \cellcolor{agray} &\cellcolor{agray} & 
 \cellcolor{agray}  \cellcolor{agray}   &  \cellcolor{agray} &\cellcolor{agray} & \cellcolor{agray}   \\
& $31$   & $0$ & $0$   &\cellcolor{agray}   & \cellcolor{agray}    & \cellcolor{agray}   & \cellcolor{agray}    & \cellcolor{agray}   &  \cellcolor{agray} &\cellcolor{agray} &
 \cellcolor{agray}  \cellcolor{agray}   &  \cellcolor{agray} &\cellcolor{agray} &  \cellcolor{agray} \\
 & $37$  & $0$   & $0$   & $0$   &  \cellcolor{yellow} $2$    &$0$   & \cellcolor{agray}    & \cellcolor{agray}   &  \cellcolor{agray} &\cellcolor{agray} &
 \cellcolor{agray}  \cellcolor{agray}   &  \cellcolor{agray} &\cellcolor{agray} &  \cellcolor{agray} \\
Prime $p$   & $41$  & $0$  & $0$   & $0$   &$0$    &$0$  & $0$    &$0$   &  \cellcolor{agray} &\cellcolor{agray} &
 \cellcolor{agray}  \cellcolor{agray}   &  \cellcolor{agray} &\cellcolor{agray} &  \cellcolor{agray}  \\
   & $43$  & $0$  & $0$   & $0$   &$0$    &$0$   & $0$    & $0$  & $0$  &\cellcolor{agray} &
 \cellcolor{agray}  \cellcolor{agray}   &  \cellcolor{agray} &\cellcolor{agray} &  \cellcolor{agray}  \\
  & $47$   & $0$  & $0$   & $0$   &$0$    &$0$   &$0$    & $0$   & $0$ & $0$ & $0$  & \cellcolor{agray} & \cellcolor{agray} & \cellcolor{agray}  \\
   & $53$   & $0$  & $0$   & $0$   &$0$    &$0$   &$0$    & $0$   & $0$ & $0$ & $0$  & \cellcolor{yellow} $2$  & $0$ &  $0$ \\
      & $\vdots$     & $\vdots$  & $\vdots$   & $\vdots$    & $\vdots$  & $\vdots$   & $\vdots$   &$\vdots$    &$\vdots$   &$\vdots$    & $\vdots$   & $\vdots$ & $\vdots$ & $\vdots$ \\
      & $71$    & $0$  & $0$   & $0$     & $0$  & \cellcolor{yellow}  $6$   & $0$   &$0$    &$0$   &$0$    & $0$   & $0$ & $0$ & $0$ \\
     & $\vdots$     & $\vdots$  & $\vdots$   & $\vdots$    & $\vdots$  & $\vdots$   & $\vdots$   &$\vdots$    &$\vdots$   &$\vdots$    & $\vdots$   & $\vdots$ & $\vdots$ & $\vdots$ \\
      & $103$     &
     $0 $   & $0$ & $0$   &\cellcolor{yellow}$2$    &  $0$  & $0$ & $0$ & $0$
    & $0$ &
    $0$ & $0$ & $0$ & $0$  \\
   & $\vdots$     & $\vdots$  & $\vdots$   & $\vdots$    & $\vdots$  & $\vdots$   & $\vdots$   &$\vdots$    &$\vdots$   &$\vdots$    & $\vdots$   & $\vdots$ & $\vdots$ & $\vdots$ \\
       & $113$    & $0$  & \cellcolor{yellow}$2$  & $0$     & $0$  & $0$   & $0$   &$0$    &$0$   &$0$    & $0$   & $0$ & $0$ & $0$ \\
       & $\vdots$     & $\vdots$  & $\vdots$   & $\vdots$    & $\vdots$  & $\vdots$   & $\vdots$   &$\vdots$    &$\vdots$   &$\vdots$    & $\vdots$   & $\vdots$ & $\vdots$ & $\vdots$ \\
     \end{tabular}
\end{center}
\end{table}

\begin{table}[ht]
\definecolor{agray}{gray}{0.8}
\begin{center}
\begin{tabular}{rr|*{13}{c}}
\multicolumn{2}{c}{} & \multicolumn{13}{c}{Serre weight $k$} \\
 &          &$29$ &  $31$ & $33$ & $35$ & $37$ & $39$ & $41$ & $43$ & $45$ & $47$ & $49$ & $51$ & $53$ \\
\cline{2-15}\\[-12.9pt]
     & $211$    & $0$  & $0$   & $0$     & $0$  & $0$   & $0$   &$0$    &$0$   &$0$    & $0$   & $0$ & $0$ & $0$ \\
     & $223$    & $0$  & $0$   & $0$     & $0$  & $0$   & \cellcolor{yellow}$2$   &$0$    &$0$   &$0$    & $0$   & $0$ & $0$ & $0$ \\
       & $\vdots$     & $\vdots$  & $\vdots$   & $\vdots$    & $\vdots$  & $\vdots$   & $\vdots$   &$\vdots$    &$\vdots$   &$\vdots$    & $\vdots$   & $\vdots$ & $\vdots$ & $\vdots$ \\
       & $353$    & $0$  & $0$   & $0$     & $0$  & $0$   & $0$   &$0$
  & $0$  & $0$   & $0$  & $0$ & \cellcolor{red}  &
       \cellcolor{red} \\
         &  $\ddots$ &   &    &      &   &    &    &    &   &    &  &   &    & \\
     Prime $p$      &  $523$    & $0$  & $0$   & $0$     & $0$  & $0$   & $0$   &$0$  
& \cellcolor{yellow}$2$  & $0$   & $0$  &$0$ & \cellcolor{red}  &
              \cellcolor{red}  \\
         &  $\ddots$ &   &    &      &   &    &    &    &   &    &  &   &    & \\
         & $599$    & $0$  & $0$   & \cellcolor{yellow}$2$     & $0$  & $0$   & $0$   &$0$  & $0$  & $0$   & $0$  & $0$ & \cellcolor{red}  &
              \cellcolor{red}  \\
         &  $\ddots$ &   &    &      &   &    &    &    &   &    &  &   &    & \\
         & $643$    & $0$  & $0$   & $0$     & $0$  & $0$   & $0$   &$0$  & $0$  & $0$   & $0$  & \cellcolor{yellow}$2$ & \cellcolor{red}  &
              \cellcolor{red}  \\
         &  $\ddots$ &   &    &      &   &    &    &    &   &    &  &   &    & \\

         & $701$    & $0$  & $0$   & $0$     & $0$  & $0$   & $0$  & 
\cellcolor{yellow}$2$ 
 & $0$  & $0$  & $0$  & $0$ & \cellcolor{red}  &
              \cellcolor{red}  \\
         &  $\ddots$ &   &    &      &   &    &    &    &   &    &  &   &    & \\
         & $743$    & $0$  & $0$   & $0$     & $0$  & $0$   & $0$  & $0$
 & $0$  & $0$   & 
\cellcolor{yellow}$2$  & $0$ & \cellcolor{red}  &
              \cellcolor{red}  \\
         &  $\ddots$ &   &    &      &   &    &    &    &   &    &  &   &    & \\
          & $983$    & $0$  & $0$   & $0$     & $0$  & $0$   & $0$   & $0$  
& $0$  & $0$   & $0$  & $0$ & \cellcolor{red}  &
        \cellcolor{red} \\
     \end{tabular}
\end{center}
\end{table}

\begin{table}[ht]
\definecolor{agray}{gray}{0.8}
\begin{center}
\begin{tabular}{rr|*{13}{c}}
\multicolumn{2}{c}{} & \multicolumn{13}{c}{Serre weight $k$} \\
 &          &$55$ &  $57$ & $59$ & $61$ & $63$ & $65$ & $67$ & $69$ & $71$ & $73$ & $75$ & $77$ & $79$ \\
\cline{2-15}\\[-12.9pt]
& $59$ & $0$  & $0$  & $0$   & \cellcolor{agray}    & \cellcolor{agray}   & \cellcolor{agray}    & \cellcolor{agray}   &  \cellcolor{agray} &\cellcolor{agray} & 
 \cellcolor{agray}  \cellcolor{agray}   &  \cellcolor{agray} &\cellcolor{agray} & \cellcolor{agray}   \\
 & $61$  & $ $0$ $ &  $0$  & $0$   &$0$    & \cellcolor{agray}    & \cellcolor{agray}    & \cellcolor{agray}   &  \cellcolor{agray} &\cellcolor{agray} &
 \cellcolor{agray}  \cellcolor{agray}   &  \cellcolor{agray} &\cellcolor{agray} &  \cellcolor{agray} \\
Prime $p$   & $67$  & $0$  & $0$   & $0$   &$0$    &$0$  & $0$    & $0$  &  \cellcolor{agray} &\cellcolor{agray} &
 \cellcolor{agray}  \cellcolor{agray}   &  \cellcolor{agray} &\cellcolor{agray} &  \cellcolor{agray}  \\
   & $71$  &  $0$  & $0$  & $0$   & $0$      & $0$   & $0$  & $0$ & $0$  
    & $0$ &
 \cellcolor{agray}  \cellcolor{agray}   &  \cellcolor{agray} &\cellcolor{agray} &  \cellcolor{agray}  \\
  & $73$   & $0$  & $0$   & $0$   &$0$    &$0$   &$0$    &  $0$   &$0$  & $0$  &  $0$    & \cellcolor{agray} & \cellcolor{agray} & \cellcolor{agray}  \\
    & $79$    & $0$   & $0$   &$0$    &$0$   &$0$    &  $0$
      &$0$ & $0$  
    & \cellcolor{red} &
 \cellcolor{red}  \cellcolor{red}   &  \cellcolor{red} &\cellcolor{red} &  \cellcolor{red}  \\
    & $83$    & 
   $0$   &$0$   & \cellcolor{yellow}$2$   & $0$      & $0 $   & $0$ & \cellcolor{red}   & \cellcolor{red}  
    & \cellcolor{red} &
 \cellcolor{red}  \cellcolor{red}   &  \cellcolor{red} &\cellcolor{red} &  \cellcolor{red}  \\
     & $89$    & 
    $0$ & $0$ & $0$ & $0$  & $0$   & $0$  & \cellcolor{red}   & \cellcolor{red}  
    & \cellcolor{red} &
 \cellcolor{red}  \cellcolor{red}   &  \cellcolor{red} &\cellcolor{red} &  \cellcolor{red}  \\
     & $97$    & 
     $0$ & $0$ & $0$ & $0$ & $0$ & \cellcolor{yellow}$2$   &$0$   & \cellcolor{red}  
    & \cellcolor{red} &
 \cellcolor{red}  \cellcolor{red}   &  \cellcolor{red} &\cellcolor{red} &  \cellcolor{red}  \\
   & $101$    & 
     $0$   & $0$ & $0$   & $0$     & $0$   & $0$   & $0$  & \cellcolor{red}  
    & \cellcolor{red} &
 \cellcolor{red}  \cellcolor{red}   &  \cellcolor{red} &\cellcolor{red} &  \cellcolor{red}  \\
  & $103$    &  $0$ &   $0$   & $0$   &$0$    & $0$   & $0$   & \cellcolor{red}   & \cellcolor{red}  
    & \cellcolor{red} &  \cellcolor{red}   & \cellcolor{red}    & \cellcolor{red}   & \cellcolor{red}   \\
  & $107$    &  $0$ &   $0$   &
 $0$   & $0$    & \cellcolor{red}   & \cellcolor{red}   & \cellcolor{red}   & \cellcolor{red}
    & \cellcolor{red} &  \cellcolor{red}   & \cellcolor{red}    & \cellcolor{red}   & \cellcolor{red}   \\
 & $109$    & $0$ &   $0$   &
 $0$   & $0$    & \cellcolor{red}   & \cellcolor{red}   & \cellcolor{red}   & \cellcolor{red}
    & \cellcolor{red} &  \cellcolor{red}   & \cellcolor{red}    & \cellcolor{red}   & \cellcolor{red}   \\
 & $113$    &  $0$ &   $0$   &
 $0$   & $0$    & \cellcolor{red}   & \cellcolor{red}   & \cellcolor{red}   & \cellcolor{red}
    & \cellcolor{red} &  \cellcolor{red}   & \cellcolor{red}    & \cellcolor{red}   & \cellcolor{red}   \\
 & $127$    &  $0$ &   \cellcolor{red}   &
 \cellcolor{red}   &\cellcolor{red}    & \cellcolor{red}   & \cellcolor{red}   & \cellcolor{red}   & \cellcolor{red}
    & \cellcolor{red} &  \cellcolor{red}   & \cellcolor{red}    & \cellcolor{red}   & \cellcolor{red}   \\
 & $131$    &  $0$ &  $0$  &
 $0$   & $0$    & \cellcolor{red}   & \cellcolor{red}   & \cellcolor{red}   & \cellcolor{red}
    & \cellcolor{red} &  \cellcolor{red}   & \cellcolor{red}    & \cellcolor{red}   & \cellcolor{red}   \\
 & $137$    &  $0$  &   \cellcolor{red}   &
 \cellcolor{red}   &\cellcolor{red}    & \cellcolor{red}   & \cellcolor{red}   & \cellcolor{red}   & \cellcolor{red}
    & \cellcolor{red} &  \cellcolor{red}   & \cellcolor{red}    & \cellcolor{red}   & \cellcolor{red}   \\
 & $139$    &  $0$ &   \cellcolor{red}   &
 \cellcolor{red}   &\cellcolor{red}    & \cellcolor{red}   & \cellcolor{red}   & \cellcolor{red}   & \cellcolor{red}
    & \cellcolor{red} &  \cellcolor{red}   & \cellcolor{red}    & \cellcolor{red}   & \cellcolor{red}   \\
 & $149$    &  $0$ &   \cellcolor{red}   &
 \cellcolor{red}   &\cellcolor{red}    & \cellcolor{red}   & \cellcolor{red}   & \cellcolor{red}   & \cellcolor{red}
    & \cellcolor{red} &  \cellcolor{red}   & \cellcolor{red}    & \cellcolor{red}   & \cellcolor{red}   \\
 & $151$    & $0$ & $0$ & $0$ & $0$ 
    & \cellcolor{red}   & \cellcolor{red}   & \cellcolor{red}   & \cellcolor{red}
    & \cellcolor{red} &  \cellcolor{red}   & \cellcolor{red}    & \cellcolor{red}   & \cellcolor{red}   \\
 & $157$    &  $0$ & $0$ & $0$ & $0$ 
    & \cellcolor{red}   & \cellcolor{red}   & \cellcolor{red}   & \cellcolor{red}
    & \cellcolor{red} &  \cellcolor{red}   & \cellcolor{red}    & \cellcolor{red}   & \cellcolor{red}   \\
      & $163$    &   \cellcolor[cmyk]{0.5,0.2,0,0}$2$  & 
$0$ & $0$ & $0$   
    & \cellcolor{red} &  \cellcolor{red}   & \cellcolor{red}    & \cellcolor{red}   & \cellcolor{red}  &
 \cellcolor{red}  \cellcolor{red}   &  \cellcolor{red} &\cellcolor{red} &  \cellcolor{red}  \\
 & $167$    &  $0$ &   \cellcolor{red}   &
 \cellcolor{red}   &\cellcolor{red}    & \cellcolor{red}   & \cellcolor{red}   & \cellcolor{red}   & \cellcolor{red}
    & \cellcolor{red} &  \cellcolor{red}   & \cellcolor{red}    & \cellcolor{red}   & \cellcolor{red}   \\
 & $173$    &  $0$ &   \cellcolor{red}   &
 \cellcolor{red}   &\cellcolor{red}    & \cellcolor{red}   & \cellcolor{red}   & \cellcolor{red}   & \cellcolor{red}
    & \cellcolor{red} &  \cellcolor{red}   & \cellcolor{red}    & \cellcolor{red}   & \cellcolor{red}   \\
 & $179$    &  $0$  &   $0$   &
 $0$   & $0$    & \cellcolor{red}   & \cellcolor{red}   & \cellcolor{red}   & \cellcolor{red}
    & \cellcolor{red} &  \cellcolor{red}   & \cellcolor{red}    & \cellcolor{red}   & \cellcolor{red}   \\
     & $181$    &  $0$  &   $0$   &
$0$  & \cellcolor{red}    & \cellcolor{red}   & \cellcolor{red}   & \cellcolor{red}   & \cellcolor{red}
    & \cellcolor{red} &  \cellcolor{red}   & \cellcolor{red}    & \cellcolor{red}   & \cellcolor{red}   \\
     & $191$    &  $0$  &   $0$   &
$0$  & \cellcolor{red}     & \cellcolor{red}   & \cellcolor{red}   & \cellcolor{red}   & \cellcolor{red}
    & \cellcolor{red} &  \cellcolor{red}   & \cellcolor{red}    & \cellcolor{red}   & \cellcolor{red}   \\
     & $193$    &  $0$  &   $0$   &
 $0$   & $0$    & \cellcolor{red}   & \cellcolor{red}   & \cellcolor{red}   & \cellcolor{red}
    & \cellcolor{red} &  \cellcolor{red}   & \cellcolor{red}    & \cellcolor{red}   & \cellcolor{red}   \\
     & $197$    &  $0$  &   $0$   &
 $0$   & \cellcolor{red}    & \cellcolor{red}   & \cellcolor{red}   & \cellcolor{red}   & \cellcolor{red}
    & \cellcolor{red} &  \cellcolor{red}   & \cellcolor{red}    & \cellcolor{red}   & \cellcolor{red}   \\
       & $199$    &  $0$ &   \cellcolor{red}   &
  \cellcolor{red} & \cellcolor{red}    & \cellcolor{red}   & \cellcolor{red}   & \cellcolor{red}   & \cellcolor{red}
    & \cellcolor{red} &  \cellcolor{red}   & \cellcolor{red}    & \cellcolor{red}   & \cellcolor{red}   \\
      \end{tabular}
\end{center}
\label{table:mixed}
\end{table}

 \chapter{Some Miscellaneous Remarks}
 \label{chapter:ch9} 
\section{Higher weights}  \label{higherweights}

 It is natural to also consider the ``higher weights,'' that is to say, cohomology with coefficients
in a nontrivial local system:

 Fixing $E$ that contains the Galois closure of $F$ over $\Q$, 
 and let $\Sigmainf$ be the set of embeddings $F \hookrightarrow E$. 
We may parameterize local systems  by integer-valued functions $\weightk: \Sigmainf
\rightarrow \{2,4,6,\dots\}$ as follows:

If $S$ is any quaternion algebra (or the split matrix algebra) over $E$, and $m$ a  positive integer,  we write $[S]_m$ for the $2m+1$-dimensional $E$-vector space 
$ \Sym^m S^0 /  \Delta \cdot \Sym^{m-2} S^0 ,$
where $S^0$ is the set of trace-zero elements in $S$, and
  $\Delta$ is any nonzero invariant element of $\Sym^2 S^0$.
Then the conjugation action yields an action of $S^{\times}/E^{\times}$
on $[S]_m$.

Now $\sigma \in \Sigmainf$ gives rise to a map $\iota_{\sigma}: D \hookrightarrow D_{\sigma} :=  (D \otimes_{(F, \sigma)}  E)$.     Fixing $\weightk$ as above, we write
 $k_{\sigma}$ for the value $\weightk(\sigma)$ on $\sigma \in \Sigma_{\infty}$. 
 Then $D^{\times}/F^{\times}$ acts (by means of the embeddings $\iota_{\sigma}$)
 on the vector space
 $$V_{\weightk} := \bigotimes [D_{\sigma}]_{\weightk(\sigma)/2 -1 }.$$

 This gives rise to a local system $\localsystem_{\weightk}$ of $E$-vector spaces over $Y(K)$: 
$$   \G(F) \backslash (\mathbf{H}^3 \times \G(\Afinite) \times V_{\weightk}) / K \rightarrow Y(K),  $$
where the action of $\G(F) = D^{\times}/F^{\times}$ on $V_{\weightk}$ is as just described. 

  This local system admits an ``integral structure,''
that is to say, if $R$ is the ring of integers of $E$, there is a local system
$\localsystem_{\weightk}^R$ such that the induced map
$\localsystem_{\weightk}^R \otimes_R E \longrightarrow \localsystem_{\weightk}$
is an isomorphism.   We omit the easy proof. 

 \medskip
 
All the questions analyzed in this book have natural analogues   in this setting, e.g.
compare the cohomology $H_1(Y(K_{\Sigma}), \localsystem^R_{\weightk})$
and $H_1(Y'(K_\Sigma), \localsystem^R_{\weightk})$. 
We anticipate that all of our results in trivial weight extend to this case without substantial change,
although the analysis of congruence homology and cohomology is somewhat more complicated
at ``small primes,'' i.e. primes that do not satisfy all of the following conditions:
  $p - 1 > \min\{k_{\sigma}\}$; for all $v|p$, $K_v = \GL_2(\OL_v)$;  $p \nmid w_F$.

There is one way in which this general case is {\em easier} than the case 
of trivial coefficients.  
For an embedding $\lambda: E \hookrightarrow \C$, say that
$\weightk$ is $\lambda$-parallel if the value of $\weightk$ on the two complex places
coincide.  
Then (for any $\lambda$) the cohomology $H_i(\YO,\localsystem_{\weightk})  = 0$ 
for all $i$ if $\weightk$ is not $\lambda$-parallel.  Of course, $H_i(\YO, \localsystem_{\weightk}^R)$ need not be zero, even
 when $\weightk$ is not parallel.

In the case of non-parallel $\weightk$, so that $\localsystem_{\weightk} $ 
is acyclic,  our results become substantially stronger.
 {\em The absence of characteristic zero cohomology means that regulators vanish, and therefore
 we obtain precise  numerical Jacquet--Langlands theorems in all cases. }

\subsection{Completed cohomology and \texorpdfstring{$p$}{p}-adic Jacquet--Langlands}
 
What might one expect be the ultimate version of an integral Jacquet--Langlands theorem \label{section:padic}
for homology to be?  A plausible candidate for a formulation (completed at  a prime $p$) would use the notion of completed cohomology:

Consider the following completed homology groups
(see~\cite{E,CE2})
$$\Ht_i(\C):= \lim_{\underset{K}{\leftarrow}} H_i(Y(K),\C), \qquad
\Ht_i(\Z_p):= \lim_{\underset{K}{\leftarrow}} H_i(Y(K),\Z_p).$$
The first group has a discrete action of $\G(\Afinite)$, and is the usual space on which
one constructs adelic representations attached to cohomological automorphic forms.
The second group, however, carries a continuous action of $\G(\OL_p) = \prod_{v|p} \G(\OL_v)$, which records
all the mod-$p$ congruences between automorphic forms.

A local--global compatibility conjecture would link the representation module
$\Ht_i(\Z_p) \otimes \Q_p$ (with the discrete action  of $\GL(\OL_v)$ for $v$ not dividing
$p$ and the continuous action at $p$) with an explicit representation coming
from local Langlands at each $v \nmid  p$ and $p$-adic local Langlands for $v|p$.
A local $p$-adic Jacquet--Langlands correspondence would then compare this
with the analogous construction for $\G'$. As many of the ingredients to define such
a conjecture do not yet exist, however, we shall not try to formulate any coherent conjecture.
 
 One feature to expect of such any
 properly formulated
 conjecture  is {\em independence of the weight:}
 the conjecture at any two weights
 would be equivalent.
 Indeed,  fixing $\weightk$, put $\localsystem_p = \localsystem^R \otimes_{\Z} \Z_p$.
One may define  the corresponding completed homology groups
$\Ht_i(\Mp)$ by the corresponding formula:
$$\Ht_i(\Mp):= \lim_{\underset{K}{\leftarrow}} H^i(Y(K),\Mp).$$
 The ``independence of the weight'' is the following assertion:   %
There is an isomorphism
$$\Ht_i(\Mp) \simeq \Ht_i(\Z_p) \otimes \localsystem_p,$$
which ``pulls out'' the coefficients.

 We have decided, however, to eliminate a detailed discussion of the
 higher weight case from this book, essentially because it would add
 substantial length without contributing any fundamentally
 new ideas.

\section{\texorpdfstring{$\GL(3)/\Q$}{GL(3)} and other groups}

\medskip

An interesting question is to extend our analysis to inner forms of $\PGL_3$ (and, indeed, to other groups; but, as we explain below, there is no shortage of new difficulties even in the case of $\PGL_3$).

Let $D/\Q$  be a division algebra over $\Q$ of dimension $9$, 
and let $\G = \GL_1(D)/\mathbb{G}_m$.
Assume that $D$ is indefinite, so $D \otimes \R = \GL_3(\R)$.
The corresponding symmetric  space  $Y$ is $5$-dimensional, and so
the arithmetic quotients are not hermitian, as for quaternion algebras over $F$.
We expect (see \cite{AB}) that $H_2(Y, \Z)$ will often have a large torsion part,
and it is natural, then, to relate it to $H_2(Y', \Z)$, where $Y'$ is an arithmetic  quotient
for the {\em split} group $\G = \PGL_3$. 
 
Here there are several new difficulties: 

\begin{itemize}
\item[(i)]    
 In the case when $\G$ is nonsplit, we cannot {\em prove}
 that all the interesting torsion occurs in $H_2$. 
 
(If $\G$ is split, all interesting torsion in homology must occur in $H_2$: 
 because of the congruence subgroup property 
$$H_1(Y(K),\Z) = H_{1,\con}(Y(K),\Z).$$ Also by Poincar{\'e} duality, $H_3(Y(K), \Z)$
is torsion-free  if $Y(K)$ is a manifold.)

\item[(ii)] Even if $\G$ is  split, we have no analogue of Ihara's Lemma: 
Since the interesting homological torsion is in degree $2$, 
an analogue of Ihara's lemma would require an analysis of the {\em second} homology
of $S$-arithmetic groups such as $\PGL_3(\Z[\frac{1}{p}])$; we do not know how to do this. 
    
\item[(iii)]  In the nonsplit case,  there is no natural analogue of the 
analysis of~\S~\ref{regcusp} --- that is to say, we cannot establish a direct relationship
between the regulator and an $L$-function. 

    Although we can still construct two-dimensional cycles
    on $Y(K)$ starting with tori in $\G$, 
    and it is quite plausible that such cycles span $H_2$, 
the integral of a (cohomological) cusp form over such a torus has no known
direct relationship to an $L$-function.

\item[(iv)] The analysis in the {\em split} case appears substantially harder than what
we carried out in Chapter~\ref{chapter:ch5}; in particular, it is not clear in the present
case either what the most natural ``truncation'' of the manifold $Y(K)$ is, 
nor how one can compute all the low-lying eigenvalues of this truncation. 
\end{itemize}

Taking these together, it seems to be far from routine
to extend the results of this manuscript to the $\PGL_3$ setting.  Such an extension
seems both interesting and desirable.

\section{Further Questions}  \label{sec:qns}

\begin{enumerate}

\item  \textbf{Do periods of torsion classes detect classes in Galois cohomology?}

 It is well-known that periods of automorphic forms are related to $L$-functions,
and so (under Birch--Swinnerton-Dyer type conjectures) to the orders of suitable Selmer groups. 
In our setting, $L$-functions are unavailable; nonetheless, is still makes sense to ask
whether mod $p$ periods vanish when certain Selmer groups are nonzero. 
Our preliminary experiments
suggest that this may be so, and we will return to this issue in a sequel.

\medskip

\item \textbf{Do torsion classes lift to characteristic zero after raising the level?}

Our experiments showed this was often so for torsion classes of relatively low level
and modulo small primes. It seems difficult to imagine that this is always true 
for torsion classes modulo large primes; if so, one would have to raise the level {\em greatly}.
Nonetheless we cannot rule it out; is there an argument showing that it cannot happen?

\medskip

\item \label{stable-ash}  \textbf{Are low degree cohomology classes always Eisenstein?}

Given 
an arithmetic group $\Gamma$ in a group of real rank $r$, is it true that 
all classes in $H_j(\Gamma, \Z)$ for $j \leq r$ are acted upon by every Hecke operator $T$
through the degree $\deg(T)$?
For characteristic zero classes  this is a consequence of the 
Kumaresan vanishing theorem (see~\cite{VZ}, Table 8.2).
After a first draft of this book appeared, a version of this conjecture
has been verified more generally by the first-named author (F.C.) and M. Emerton. 

\medskip
\item \label{regulators} \textbf{Study integral properties of regulators.}

In our context, we were able to show that regulators coincide up to $\mathbf{Q}^{\times}$
with products of $L$-values. It would be most desirable to show an integral version of this, i.e. that regulators coincide
up to sign with periods of the corresponding motives; the precise statement is not yet clear. 
In an analytic direction, it would be desirable to show that regulators do not grow or decay too fast with the complexity of the underlying manifold.  In yet a third direction, it is of interest
to study regulators in a context such as an inner form of $\PGL_4$ over $\Q$, where there are not enough sub-locally symmetric spaces
to represent interesting homological cycles. 

\medskip

\item \label{klanglands} \textbf{Comparing $K$-theory and a torsion Langlands correspondence.}

Both give (different) insights into the structure of low degree  homology of general linear groups
over integer rings. Is it possible to ``reconcile'' the two descriptions? Our discussion
is a very small example of this : We constructed torsion classes from $K_2$ and
discussed the associated Galois representations. 
\medskip 
\item \textbf{How far can one ``enrich'' Jacquet--Langlands correspondences?}

In this book, we have seen that the classical Jacqute--Langlands correspondence, in the case of hyperbolic 
$3$-manifolds, reflects deeper connections; the corresponding manifolds have related 
volumes and torsion homology. 
It may be that the relationship between Jacquet--Langlands pairs goes deeper than we yet understand. For instance, is there any direct relation between the profinite completions of the two discrete groups?

\medskip

\item \textbf{Higher analytic Torsion.}

 Is there some notion of higher analytic torsion that would
allow us to compare interesting torsion cohomology classes on arithmetic quotients of
general symmetric spaces? On a quotient of $G_{\infty}/K_{\infty}$ the analytic torsion vanishes
unless $\mathrm{rank}(G_{\infty}) - \mathrm{rank}(K_{\infty}) = 1$. 

It may be worth noting that if, e.g., $\mathrm{rank}(G_{\infty}) - \mathrm{rank}(K_{\infty})= 3$,
the vanishing of analytic torsion is not ``topologically obvious,'' and carries nontrivial 
information.
 
 \medskip

\item \textbf{Weight one forms over $\Q$.}  There are certain parallels between
torsion classes for $\PGL_2$ over an imaginary quadratic field  and mod $p$ weight one forms over $\Q$.  Although the latter case is far better understood because of the availability of congruences to higher weight, it still may be of interest to develop a corresponding picture to this book, with the role of analytic torsion and the Cheeger--M{\"u}ller theorem replaced by the Arakelov--Riemann-Roch formula. 
\end{enumerate}

\bibliographystyle{plain}
\bibliography{Cave}

\printindex

\end{document}